\documentclass[12pt,a4paper,oneside,onecolumn,openright,final]{memoir}

\usepackage[english]{babel}
\usepackage[T1]{fontenc}
\usepackage[utf8]{inputenc}

\usepackage{fourier}
\usepackage{fontawesome}
\usepackage{miama}

\DeclareFontFamily{U}{skulls}{}
\DeclareFontShape{U}{skulls}{m}{n}{ <-> skull }{}
\newcommand{\Skull}{%
\text{\usefont{U}{skulls}{m}{n}\symbol{'101}}%
}

\settrimmedsize{297mm}{210mm}{*}
\settypeblocksize{660pt}{460pt}{*}
\setulmargins{4cm}{*}{*}
\setlrmargins{*}{*}{1}
\setmarginnotes{17pt}{51pt}{\onelineskip}
\setheadfoot{\onelineskip}{2\onelineskip}
\setheaderspaces{*}{2\onelineskip}{*}
\checkandfixthelayout

\usepackage{amsmath}
\usepackage{amssymb}
\usepackage{amsthm}
\usepackage{mathtools}

\newtheorem{theorem}{Theorem}[section]
\newtheorem{lemma}[theorem]{Lemma}

\newtheorem{proposition}[theorem]{Proposition}
\newtheorem{corollary}[theorem]{Corollary}
\theoremstyle{remark}

\newtheorem{exercise}[theorem]{Exercise}
\newtheorem{exercisesk}[theorem]{$\Skull\,$Exercise}

\numberwithin{equation}{section}

\usepackage{graphicx}
\graphicspath{{figures/}}
\usepackage[dvipsnames]{xcolor}
\usepackage[textfont={small,sl},labelfont={small,bf}]{caption}

\usepackage{authblk}
\makeindex

\def \R {\mathbb{R}}
\def \C {\mathbb{C}}
\def \N {\mathbb{N}}
\def \Z {\mathbb{Z}}
\def \Q {\mathbb{Q}}
\newcommand{\code}[1]{\texttt{#1}}

\renewcommand{\Re}{\mathfrak{Re}}
\renewcommand{\Im}{\mathfrak{Im}}
\renewcommand{\leq}{\leqslant}
\renewcommand{\le}{\leqslant}
\renewcommand{\geq}{\geqslant}
\renewcommand{\ge}{\geqslant}

\DeclareMathOperator{\arccosh}{arccosh}
\DeclareMathOperator{\arccot}{arccot}

\usepackage{hyperref}

\title{A primer on Fourier Series}

\author{%
Serena Dipierro\thanks{SD \& EV have removed their affiliations from this work as a form of civil protest.},
David Pfefferl\'e
and 
Enrico Valdinoci\thanksmark{1}%
}

\begin{document}

\frontmatter
\begin{titlingpage}
\maketitle
\end{titlingpage}


\clearpage
\thispagestyle{empty}
{\fmmfamily
  
  \vspace*{4em}
  \begin{quotation}
  \Huge
      ``Le savant n'\'etudie pas la nature parce que cela est utile; il l'\'etudie parce qu'il y prend plaisir''. 
      \sourceatright{\Large Henri Poincar\'e, \emph{Science et M\'ethode}, Part~I, Chapter~1, 1908.}
  \end{quotation}
  
\vfill

\begin{quotation}
 \Huge{``Physics is like sex: sure, it may give some practical results, but that's not why we do it''.
 
 \sourceatright{\Large Unfoundedly attributed to Richard P. Feynman, appearing in Robyn Williams, \emph{Scary Monsters and Bright Ideas}, page~44, 2000.}}
\end{quotation}

\vfill
}
\clearpage


\tableofcontents


\chapter*{Preface}
This is a book on Fourier Series, suitable for both undergraduate and graduate levels. 
The presentation aims to be precise and culturally informed, yet accessible and not overly technical,
making the book well-suited for a high-quality undergraduate course in Fourier Series and its applications, typically offered in mathematics, physics, and engineering departments.

The book also offers insights and ideas that occasionally venture beyond standard approaches, making it valuable for graduate students as well.

It is also well-suited for independent study, as detailed solutions are provided for many exercises. 
We hope that these exercises are expediting the development, 
encouraging participation on the part of the reader.

The worked examples and applications may also be of interest to professionals and applied scientists.

In short, this is a book for anyone interested in exploring the beauty and utility of Fourier Series, appealing to students, instructors, practitioners, and curious readers alike.

This book cannot compare with all the excellent books already available on the topic in terms of extensiveness and depth of the material treated, but we hope it can be useful and pleasant for readers looking at arguments that are both rigorous and close to the intuition, as well as willing to challenge themselves with interesting exercises.

Because the topic of Fourier Series is already monumental, we sadly do not treat Fourier Transforms in this book (we will address these in a separate project~\cite{NEXT}).

The prior knowledge assumed from the readers are rather minimal and are essentially limited to familiarity with
one-dimensional and multi-dimensional calculus, basic measure theory, Lebesgue spaces, and approximation by smooth functions (e.g. mollification). In any case, when some specific famous result is used, a precise classical reference is given for it.

In fact, we emphasize that this book is not intended solely for readers already familiar with measure theory. We do make use of measure-theoretic concepts whenever necessary: after all, any rigorous mathematical argument will, at some point, rely\footnote{According to~\cite[pages~xi--xii]{MR2380238},
``I find it both surprising and immensely satisfying that the search for understanding of Fourier series continued to be one of the principal driving forces behind the development of analysis well into the twentieth century. The tools that these mathematicians had at hand were not adequate to the
task. In particular, the Riemann integral was poorly adapted to their needs''. See also~\cite{MR1857634}
and~\cite[Chapter~9]{MR3676533}.}
on them. However, we are fully aware that many undergraduate programs introduce measure theory only in the final term.

Similarly, some notions from complex analysis,
such as holomorphic functions,
may appear occasionally. Yet, we do not assume the reader has taken a course in complex variables. On the contrary, we make no attempt to hide technical details or obscure the inherent difficulties: mathematics is deeply interconnected, and the development of any theory rarely follows a straightforward, linear path. The notion that mathematical learning can be neatly divided into sequential courses with rigid ``prerequisites'' is an oversimplification. In reality, high-quality mathematical understanding evolves through a dynamic process, often involving the anticipation of advanced ideas and their later revisitation. Scientific maturity is gradual and requires continuous rethinking.

In line with these principles, our approach is to clearly highlight whenever more advanced material arises and to provide precise references. Readers are then free to adopt one of two learning strategies: either they may consult the references and learn the advanced material immediately to bridge the gap, or they may temporarily accept it as given and continue reading, with the intention of returning to the referenced topics later in their long-term learning journey. Both approaches are valid and viable. This book offers the flexibility for readers to choose the method that best suits their individual needs, preferences, and learning styles.

In any case, some sections are perhaps a bit advanced, or ``fancy'', and maybe they can be skipped on a first pass.
These sections are marked with a~\faBomb.

The book contains the basic theory, but also a number of specific applications.
The applications are essentially independent from each other and the reader can freely pick their favourite ones and disregard the others (and so can the lecturer to select the ones to present in class).

The text is complemented by many exercises. Some of the exercises are useful for the development of the theory itself; some exercises constitute parts of a proof that a lecturer leaves for the audience to check. Such exercises could be given as assignments before lectures, so to encourage attendees to engage with the content and to improve their readiness for class. 

The hardest (in our humble opinion) exercises are marked with a~$\Skull$.
The solutions of some of the exercises are presented in Chapter~\ref{SOLUTIONS}
(if the reader finds some exercises too hard, it is an excellent training to take a little peep at the solution, try again, and possibly iterate the process).

The book is full of digressions. It is our opinion that digressions should be part of a fruitful learning experience. They are intended to expand horizons and help the reader seize the ``bigger picture''.

We claim no expertise in Fourier Analysis and this work may well contain typographic and conceptual errors: if you spot any of them, please let us know and we will correct them for future versions of this product.

We hope you enjoy the rest of your reading.


\mainmatter

\chapter{Where to start}

The first line of a book about Fourier analysis is quite problematic.
One would like to convince the reader about the beauty and the importance of the topic, and the first line is supposed to be catchy and engaging. There are several classical arguments to do this.

For example, Chapter~1 in~\cite{MR1970295}, 
entitled \emph{The Genesis of Fourier Analysis}, relies on historical and physical motivations; it is inspiring and captivating. It begins with a quotation by Fourier himself, stressing the innovation of his vision with respect to the predecessors' (namely, d'Alembert and Euler, and we will consider this quotation later, see page~\pageref{CITA-F}).
Then, it wisely explains how to use Fourier methods to deal with
wave and heat equations (and we will discuss this in Section~\ref{PDESE}).
This indeed would provide a very solid motivation for Fourier Series (provided that the reader has already some familiarity with partial differential equations): the reader who would like to start from there can indeed profitably exploit~\cite{MR1970295}.

The strategy of~\cite{MR1145236} shares, in a sense, a common tradition, but, possibly to make
the above introductory topics more accessible with less prerequisites, this great textbook
starts with an ``Overture'', to introduce the relevant ``equations of mathematical physics'' and explain the notion of ``linear differential operators'' and the method of ``separation of variables''.

The introductions in~\cite{zbMATH03004117, MR4380761} also assume a historical point of view, and the reader willing to receive a thorough account of the original references is encouraged to look into them.

The classical textbook~\cite{MR545506} begins with a brief, though inspired, historical introduction, to quickly move to a presentation from a bird's eye perspective of the major results of the theory.
Similarly, the incipit of~\cite[Chapter~3]{MR3243734} presents an intriguing half-page historical presentation,
then moving quickly to the mathematical theory of Fourier Series in dimension~$n$.
Instead, the starting point of~\cite{MR44660} is directly from the technical side.

The fantastic textbook~\cite{MR442564} also begins with a historical introduction, explaining how Fourier introduced series and integrals to deal with partial differential equations, comparing with the previous state of the art and the subsequent developments. 
Interestingly, it is also briefly mentioned that some precursors of Fourier Series had been utilised in the antiquity to study astronomical problems, dating this trend at least back to the Babylonians (in Section~\ref{ASTRP} of this set of notes we will indeed describe in some detail some classical features of ancient astronomy in the light of Fourier methods).
Again, a reader wishing to start from this perspective is warmly encouraged to look at~\cite{MR442564}.

Another nice example of ``incipit'' is given by Chapter~1 in~\cite{MR2382058}, which motivates Fourier Series from the notion of ``synthesis'' and ``analysis''.
In a nutshell, it is often desirable to ``synthesize'' a rather arbitrary function from a set of certain ``elementary'' functions. The typical example is that of the Taylor Series, which aims at reproducing sufficiently nice (i.e. analytic) functions from a power series representation; that is by an infinite superposition of monomials (we will retake this analogy on pages~\pageref{ANL} and~\pageref{LAFOOTTF} and in Exercises~\ref{NONSepCFGDV2} and~\ref{NICELI}).

Other excellent books, such as~\cite{MR1800316, MR2039503, MR3331143, MR4404761},
do not indulge much on the history\footnote{For further information about the story of Fourier and of his times, see e.g.~\cite{MR419139}, \cite{MR605488},
\cite[Chapters~92 and~93]{MR4404761}, 
\cite{MR636934}, \cite[Chapter~1]{MR3616140}, and the references therein.} of the problem, trying to address as soon as possible the mathematical aspects of the theory, and this is also an effective and impactful way to start.

Thus, we now face a difficult choice: where do {\em we} start? Since we cannot compare with the Masters, we have decided to break all the traditions and to start {\em from nowhere}.
We completely skate over any historical introduction and we postpone all the discussions
about the application of Fourier methods to Section~\ref{APPL-SEZ1} (see also~\cite{MR3616140} for accurate descriptions
of many models and phenomena deeply relying on Fourier methods).

This radical choice of ours provides some important simplifications. First of all, the history related to Fourier methods is broad and complex, and the applications are very deep and far from intuitive, therefore
it is not easy to start from there (it may be easier instead to come back to these subtle points
after having familiarised with the basic concepts).

Moreover, no historical introduction can really make justice of the incredible power of Fourier methods, which go well beyond the original scopes for which they were invented. Fourier theory is one of the few core background of all scientists (including, but not limited to mathematicians, physicists, chemists, engineers,  statisticians, computer scientists, data scientists, etc.) and is a paradigmatic example of a field in which the applied and theoretical points of view are intimately intertwined. So, in a sense, we have deliberately decided to avoid any ``incipit'' that could have been misinterpreted or could have misrepresented
the enormous power of Fourier theory: saying that Fourier methods are good to solve the heat equation or to describe the vibration of the string of a violin, or that predecessors of Fourier Series were used thousands of years ago to describe the motion of celestial bodies, is a bit like saying that complex numbers were invented to find roots of cubic polynomials -- it may be true, but it may also sound reductive.

A final comment on these notes and on Fourier analysis in general. A natural question is {\em Why does all this have to be so freaking hard?} Well, Fourier analysis is difficult, and probably it is intrinsically so, because going from a function to its Fourier coefficients, or vice versa, is a bit like going through a wormhole. That is, over the years, scientists have looked into the possibility that black holes could be doors to other universes, and going through tunnels connecting disparate universes is iconic in the sci-fi genre. Fourier analysis is, in a sense, like this: anytime we use Fourier methods, we travel instantaneously from a function in ``normal'' variables to its counterpart in a ``frequency'' universe. These two universes may well be somewhat related, but, as it happens in sci-fi, the connection may be knotty and require some thoughtful inventiveness to be discovered.

So, let us start our adventure.


\chapter{Fourier Series: development of the theory}

\section{Periodic functions (and periodic extensions)}\label{SEC:PER:FUNC:913uo}

The object of interest in this chapter are periodic functions\footnote{The case \label{casnoco}
of periodic functions from~$\R$ to~$\C$ can also be taken into account, just separating the real and imaginary parts of the values attained by the function. Here, we stick to the case of a real range of the function not to overburden the notation.} from~$\R$ to~$\R$.
For concreteness, we will restrict ourselves to the case of functions with period~$1$
(i.e., functions~$f:\R\to\R$ such that~$f(x+1)=f(x)$ for every~$x\in\R$).
The case of periodic functions of different period can be easily reduced to this situation, see Section~\ref{ANY} for full details.

The functions of this chapter will be assumed to be in~$L^1((0,1))$, namely they are Lebesgue measurable and the integral of their absolute value on the interval~$(0,1)$ is finite (see e.g., \cite{MR3381284} for a solid introduction to measure theory and Lebesgue spaces). When further regularity is needed, we will mention it explicitly.

It is tempting to develop Fourier methods only for ``regular'' functions, and many books do so, but this approach presents severe downsides on the long run. Indeed, measure theory and Lebesgue spaces were not invented by 
dummy pure mathematicians for the sake of generalisation and to make the life of future generations of students miserable: they have been invented because everybody badly needs them.

Indeed: one of the main applications \label{OSK99f-32ut494hg}
of Fourier Series is partial differential equations and, on many occasions, solutions are known to be irregular (for example, we know that waves in the ocean can break on our surfs and kayaks, sometimes with unpleasant practical consequences).  Many partial differential equations also lie in an uncharted territory and the regularity of their solutions is, still today, mostly unknown (for example, fully understanding the regularity of the solutions of the \index{Navier-Stokes equation} Navier-Stokes equation is worth the award of a million US dollars). And also for partial differential equations for which we expect to always have regular solutions, it is almost unavoidable to look for solutions in a broader sense (thus allowing, at least in principle, irregular solutions) and only later to check that the solutions found in this way are actually regular enough (this trick of ``enlarging the space of possible solutions'' being a fundamental step to gain the necessary ``compactness'' that ensures that a solution indeed exists). So, no way to avoid irregular functions
if one deals with partial differential equations.

And how about numerics? Well, the issue is that numerical errors need to be controlled (otherwise the computer always gives an answer, and the answer is~42) and for this suitable quantitative estimates are needed. Very often, Lebesgue spaces (and even more sophisticated spaces) are of paramount importance to obtain viable bounds.

What about physics? Well, once again, physical processes are often dictated by the notion of energy, which typically requires a sound knowledge of Lebesgue spaces (and, again, most of the time, even more sophisticated spaces).

Same in computer science, digital electronics and signal processing, in which the use of binary devices often makes it convenient to use ``square waves'' to alternate steady frequencies between minimum and maximum values, and the use of sawtooth waveforms is also quite common\footnote{For example, sawtooth waves were used in the early developments of radars, see~\cite[page~66]{MR4404761}.} for chips and power-supplies, hence, even in this case, restricting to continuous/smooth functions would be too limiting. \label{GADIS}

But, in a sense even more importantly, spaces of irregular functions are very often needed from the conceptual point of view, to deal with a set of objects which is ``closed'' under reasonable (topological) assumptions: that is, if we have a sequence of functions satisfying certain properties, we want that the limit satisfies that property as well, otherwise, when using approximation arguments, we risk to fall out of the category of objects we are looking at (in this, Lebesgue measure is of paramount importance, because it assures that the pointwise limit of measurable functions is still a measurable function, which is not true, e.g., for the Riemann measure, compare with Exercise~\ref{RIEM}).

\begin{figure}[h]
\includegraphics[height=3.8cm]{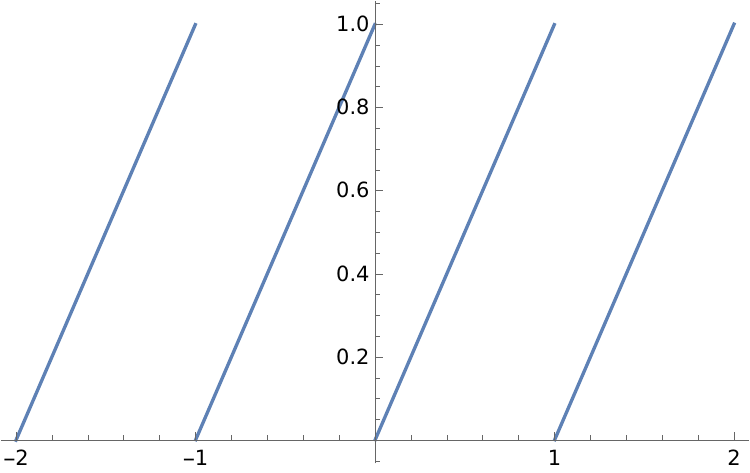}$\,\;\qquad\quad$\includegraphics[height=3.8cm]{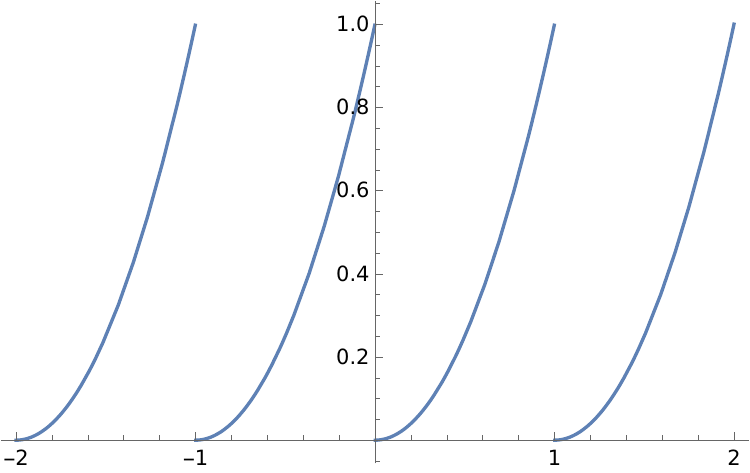}\\ $\,$
\\
\includegraphics[height=3.8cm]{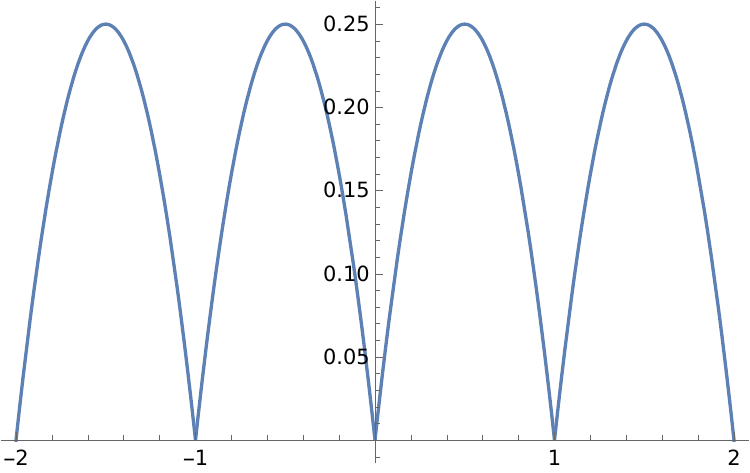}$\,\;\qquad\quad$\includegraphics[height=3.8cm]{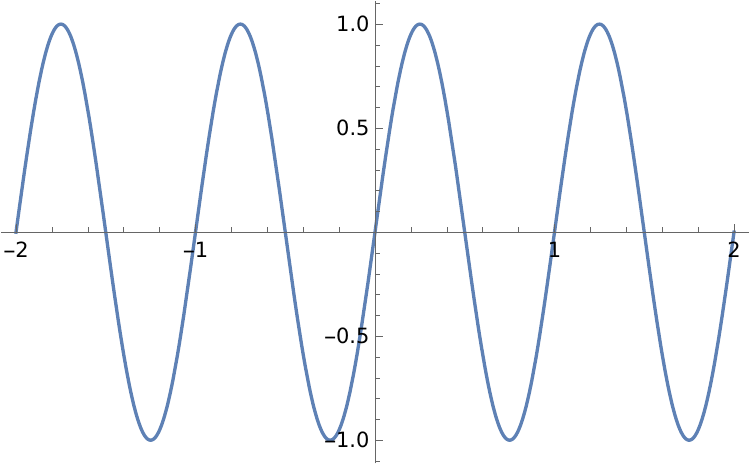}
\centering
\caption{Periodic extensions of the functions~$f(x)=x$, $f(x)=x^2$, $f(x)=x(1-x)$, and $f(x)=\sin(2\pi x)$.}\label{RET12P}
\end{figure}

If we accept to work with irregular functions, we also have a ``natural'' way to extend any function~$f:[0,1)\to\R$ to a periodic function defined in the whole of~$\R$.
Namely, given~$x\in\R$, one considers the \emph{floor function} \index{floor function} $\lfloor x\rfloor$, defined as the largest integer less than or equal to~$x$ (in this way, $\lfloor \pi\rfloor=3$, $\lfloor -e\rfloor=-3$, $\lfloor 0\rfloor=0$, etc.), and\footnote{Warning: \code{Mathematica} has the deplorable habit of using a different notation. For the fractional part, to be consistent with the notation here, do not use \code{FractionalPart[x]}, use \code{SawtoothWave[x]} instead.}
 the \emph{fractional part} \index{fractional part} $\{x\}$, defined as~$x-\lfloor x\rfloor$ \label{FRaPARTADE}
(so that~$\{1.23\}=0.23$, $\{-0.1\}=0.9$, $\{1\}=0$), and defines the \index{periodic extension}
\emph{periodic extension} of a function~$f:[0,1)\to\R$, for every~$x\in\R$, by
\begin{equation}\label{FPER} f_{\text{per}}(x):= f\big(\{x\}\big).\end{equation}
Some examples of this construction are given in Figure~\ref{RET12P}.

Notice that~$f$ can be as smooth as one wants in~$[0,1)$ (and even in~$[0,1]$), but~$f_{\text{per}}$
may lose continuity and smoothness at all the integers (one more reason for not reducing Fourier analysis to the case of regular functions).

\begin{exercise} Prove that the fractional part $x\mapsto \{x\}$ is a periodic function of period~$1$.
\end{exercise}

\begin{exercise} Prove that~$f_{\text{per}}$, as defined in~\eqref{FPER}, is a periodic function of period~$1$.
\end{exercise}

\begin{exercise} Prove that if~$f$ is periodic of period~$1$, then~$f_{\text{per}}$, as defined in~\eqref{FPER}, coincides with~$f$.
\end{exercise}

\begin{exercise} Prove that if~$f\in L^1((0,1))$, then~$f_{\text{per}}\in L^1_{\text{loc}}(\R)$.
\end{exercise}

\begin{exercise} \label{fr12}
Let~$f:\R\to\R$ be periodic of period~$1$ and $f\in L^1((0,1))$. Prove that, for all~$r\in\R$,
$$ \int_0^1 f(x)\,dx=\int_r^{r+1} f(x)\,dx.$$
\end{exercise}

\begin{exercise} Let~$f:\R\to\R$ be periodic of period~$1$. Prove that the following conditions are equivalent:
\begin{itemize}
\item $f\in L^1((0,1))$,
\item for all~$r\in\R$, one has that~$f\in L^1((r,r+1))$,
\item there exists~$r\in\R$ for which~$f\in L^1((r,r+1))$.\end{itemize}
\end{exercise}

\begin{exercise}
    Let $f:\R\to\R$ be periodic of period $1$. If $f\in C^m(\R)$, prove that $D^mf$ (the $m$th derivative) is periodic of period~$1$. What can be said of $D^mf$ with the assumption $f\in C^m((0,1))$ ?
\end{exercise}

\begin{exercise}\label{LINC} Prove that the linear combinations of periodic functions of period~$1$
is also periodic of period~$1$.

That is, assume that~$f_1,\dots,f_k:\R\to\R$ are periodic of period~$1$, let~$a_1,\dots,a_k\in\R$,
and prove that the function~$a_1f_1+\dots+a_kf_k$ is also periodic of period~$1$.
\end{exercise}

\begin{exercise}\label{PRD1}
For every~$\omega\in\left[-\frac1{2\pi},\frac1{2\pi}\right]$ and~$n\in\N$, let
$$ I_n(\omega):=\int_0^{1/4} \cos^n(2\pi t)\,\cos(2\pi\omega t)\,dt.$$
Prove that, for every~$n\ge2$,
$$ (n^2-\omega^2)\,I_{n}(\omega)=n(n-1)\,I_{n-2}(\omega).$$
\end{exercise}

\begin{exercise}\label{PRD2}
Let~$I_n$ be as in Exercise~\ref{PRD1}. Prove that, for every~$\omega\in\left[-\frac1{2\pi},\frac1{2\pi}\right]$ and~$m\in\N$,
$$ \frac{I_0(\omega)}{I_0(0)}=\frac{I_{2m}(\omega)}{I_{2m}(0)}\,\prod_{j=1}^m\frac{(2j)^2-\omega^2}{(2j)^2}.$$
\end{exercise}

\begin{exercise}\label{PRD3}
Let~$I_n$ be as in Exercise~\ref{PRD1}. Prove that, for every~$\omega\in\left[-\frac1{2\pi},\frac1{2\pi}\right]$,
$$\lim_{m\to+\infty}\frac{I_m(\omega)}{I_m(0)}=1.$$\end{exercise}

\begin{exercise}\label{PROSI} Euler's sine product formula states that, for any~$z\in\C$, \index{Euler's sine product formula}
\begin{equation}\label{EUJNC-P}
\frac{\sin (\pi z)}{\pi z}=\prod_{j=1}^{+\infty}\left(1-\frac{z^2}{ j^2}\right).\end{equation}
Prove this formula.
\end{exercise}

\begin{exercise}\label{PROSIW}
In 1656, John Wallis stated that \index{Wallis's formula}
$${\frac{\pi }{2}}=\prod_{j=1}^{+\infty }{\frac{4j^{2}}{4j^{2}-1}}.$$
Prove this identity.\end{exercise}

\begin{exercise}\label{STIRL-va}
Prove that
$$ \lim_{n\to+\infty}\frac{2^{2n}\,(n!)^2}{\sqrt{n}\,(2n)!}=\sqrt\pi.$$
\end{exercise}

\begin{exercise}\label{STIRL}
Stirling's approximation for the factorial states that \index{Stirling's approximation for the factorial}
$$\lim_{n\to+\infty} \frac{\sqrt {2\pi n}}{n!}{\left({\frac{n}{e}}\right)}^{n}=1.
$$
Prove this formula.
\end{exercise}

\begin{exercise}\label{RIEM} Let
$$ f(x):=\begin{cases}
1 & {\mbox{ if }}x\in\Q,\\
0& {\mbox{ if }}x\in\R\setminus\Q.
\end{cases}$$
Prove that~$f$ is not Riemann measurable, but
there exists a sequence of Riemann measurable functions~$f_k$ such that, for all~$x\in\R$, $f_k(x)\to f(x)$
as~$k\to+\infty$.
\end{exercise}

\begin{exercisesk}\label{SnsdFrhndiYUJS-2}
The Fundamental Theorem of Algebra \index{Fundamental Theorem of Algebra}
states that every non-constant polynomial with complex coefficients has a zero.

Prove the Fundamental Theorem of Algebra by looking at the derivative of a periodic function.
\end{exercisesk}

\begin{exercisesk}\label{JSLD.023oler-NOMPO}
Prove that there exist continuous functions that are nowhere monotonic.

More specifically, for every~$x\in\R$, let
\begin{equation}\label{JSLD.023oler-NOMPOn} f(x):=\sum_{k=0}^{+\infty}\frac{ \cos(100^k\pi x)}{2^k}.\end{equation}
Prove that~$f$ is a continuous function but there exists no interval in which~$f$ is monotonic.
\end{exercisesk}

\begin{exercisesk}\label{DDEF8AMPMvS1}
Prove that the function in~\eqref{JSLD.023oler-NOMPOn} is nowhere differentiable.
\end{exercisesk}

\section{How to ``optimally'' approximate periodic functions}\label{L2pe}

Natural examples of periodic functions of period~$1$ are~$\sin (2\pi x)$ and~$\cos(2\pi x)$.

Actually, for every~$k\in\Z$, the functions~$\sin (2\pi kx)$ and~$\cos(2\pi kx)$ are also nice
examples of periodic functions of period~$1$.

And obviously constant functions are periodic.

Thus, since the linear combination of periodic functions of period~$1$
is also periodic of period~$1$ (see Exercise~\ref{LINC}), we have that, for every~$a_0,\dots,a_N\in\R$
and~$b_1,\dots,b_N\in\R$, the function
\begin{equation}\label{TRIG-FO} \frac{a_0}2+\sum_{k=1}^N\Big(a_k\cos(2\pi kx) + b_k\sin(2\pi kx)\Big)\end{equation}
is also periodic of period~$1$.

Now, trigonometric functions are nice, but exponential are nicer (indeed, the algebraic rules for trigonometric functions are equivalent to, but much more complicated than, the ones for exponentials), therefore we have the strong preference of rewriting trigonometric functions in terms of exponentials by using the classical formulas
$$ \cos \vartheta=\frac{e^{i\vartheta}+e^{-i\vartheta}}2\qquad{\mbox{and}}\qquad
\sin \vartheta=\frac{e^{i\vartheta}-e^{-i\vartheta}}{2i}.$$
This reduces the expression in~\eqref{TRIG-FO} to an equivalent one in the form
\begin{equation}\label{EXP-FO}
\sum_{{k\in\Z}\atop{|k|\le N}} c_k \,e^{2\pi i k x},
\end{equation}
for suitable complex numbers~$c_k$,
see Exercise~\ref{smc203e}.

Expressions of the type~\eqref{TRIG-FO}, or equivalently~\eqref{EXP-FO}, are called \index{trigonometric polynomials}
``trigonometric polynomials'' (specifically, trigonometric polynomials of degree~$N$).

The core of Fourier methods is to approximate (whenever possible) a given periodic function
with a linear combination of complex exponentials. That is, one would like to approximate
a periodic function~$f$ with expressions as in~\eqref{EXP-FO}.
The issue is how to choose~$c_k$ such that such an approximation is ``good'' ?

To answer this question, given a periodic function~$f:\R\to\R$ of period~$1$, a possible strategy
is to choose the coefficients~$c_k$ in order to minimise the distance between the function in~\eqref{EXP-FO} and $f$.
A delicate matter however is what we mean by distance, or, in general, what is the most convenient notion of distance in order to implement such a strategy.

It is far from being obvious, but the most suitable notion of distance to make this idea work is provided by a Lebesgue space, namely~$L^2((0,1))$ (and, once again, this highlights the immense power of Lebesgue spaces). This is a very elegant result, sometimes falling under the name of \index{best approximation}
\emph{best approximation},
and goes as follows:

\begin{theorem}\label{BEST}
Let~$c_k$ be a sequence of complex numbers.

Let $f\in L^2((0,1))$ and, for every~$k\in\Z$, define
\begin{equation}\label{FOUCO}
\widehat f_k:=\int_0^1 f(x)\,e^{-2\pi i kx}\,dx.
\end{equation}

Then, for all~$N\in\N$,
\begin{equation}\label{bestapp}
\left\| f-\sum_{{k\in\Z}\atop{|k|\le N}} \widehat f_k \,e^{2\pi i k x}\right\|_{L^2((0,1))}\le
\left\| f-\sum_{{k\in\Z}\atop{|k|\le N}} c_k\,e^{2\pi i k x}\right\|_{L^2((0,1))}.
\end{equation}

Additionally, the inequality above is strict unless~$c_k=\widehat f_k$ for all~$|k|\le N$.
\end{theorem}

The quantities in~\eqref{FOUCO} are usually\footnote{Some authors use the nomenclature of \index{Fourier constants} \emph{Fourier constants} rather than Fourier coefficients, see e.g.~\cite{MR44660}.}
called \emph{Fourier coefficients} \index{Fourier coefficients}
of the function~$f$
(and we stress that these coefficients are complex numbers). Another motivation will be given in Exercise~\ref{LADEBD}
for why the Fourier coefficients must be defined as in~\eqref{FOUCO} if one wants to reconstruct~$f$
out of a superposition of complex exponentials. Notice that the expression~\eqref{FOUCO} makes sense under the weaker assumption of $f\in L^1((0,1))$, and will be used under this generality throughout the book.

The notation
\begin{equation}\label{FOSE}
\sum_{k\in\Z} \widehat f_k \,e^{2\pi i kx}
\end{equation}
is often\footnote{Looking at~\eqref{FOSE}, in analogy with acoustics, one can think that the $k$th Fourier coefficient
represents the ``amplitude'' of the elementary wave with ``frequency''~$k$.} called \index{Fourier Series}
\emph{Fourier Series of~$f$}. There is a caveat, however: this is in principle just a ``formal'' expression
(say, a picturesque way of listing all the Fourier coefficients of~$f$)
and the actual convergence of this series to the original function~$f$, under the appropriate assumptions, is one of the important topics that will be discussed below.

We stress that since~$f$ is periodic of period~$1$ (hence so is the function~$f(x)\, e^{-2\pi i k x}$), the Fourier coefficients of~$f$ can be computed by an integral over any interval of unit length (see Exercise~\ref{fr12}), i.e., for all~$r\in\R$,
\begin{equation}\label{ERMCBVSUEFASDG}
    \widehat f_k = \int_{r}^{r+1} f(x)\, e^{-2\pi i k x}\,dx.
\end{equation}

We also remark that~$f$ takes values in~$\R$ and therefore, for all~$k\in\Z$,
\begin{equation}\label{fasv}
{\mbox{$\widehat f_{-k}$ is equal to the complex conjugate of~$\widehat f_k$,}}
\end{equation}
see Exercise~\ref{fasv2}.

In this context, \eqref{bestapp} states that among all the possible linear combinations of complex exponentials, the one that best approximates~$f$ (at least in the sense of~$L^2((0,1))$) is the one obtained by utilising the Fourier coefficients (hence, it will be natural to try to reconstruct~$f$ from a finite sum obtained in this way, and this will be the topic of Section~\ref{CONVES0}).

To prove Theorem~\ref{BEST}, it is useful to consider, given~$g, h\in L^2((0,1),\C)$, the (complex, also known as \emph{Hermitian}) inner product \index{Hermitian product} defined by
$$ \langle g,h\rangle_{L^2((0,1))}:=\int_0^1 g(x)\,\overline{h(x)}\,dx,$$
where, as usual, the ``bar'' denotes the complex conjugation,
and notice that
\begin{equation}\label{GL2} \|g\|_{L^2((0,1))}^2=\int_0^1 |g(x)|^2\,dx=\langle g,g\rangle_{L^2((0,1))}.
\end{equation}

A tremendous advantage of the complex exponentials for computations\footnote{The list of complex exponentials $\{e^{2\pi i k x}\}_{k\in \Z}$ is said to form an orthonormal system. That the complex exponentials is actually a complete basis of this Hilbert space is equivalent to showing uniqueness and convergence of Fourier Series in $L^2$, see Section \ref{COL2} and references therein.} is that they are orthonormal with respect to this inner product.
\begin{lemma}\label{orthogonality-fourier-basis}
  Let $j,k\in \Z$. Then,
\begin{equation}\label{GL4} \langle e^{2\pi ikx},e^{2\pi ijx}\rangle_{L^2((0,1))}=\delta_{k,j}:=\begin{cases}
    1 & {\mbox{ if }}k=j,\\0&{\mbox{ if }}k\ne j.\end{cases}
\end{equation}
\end{lemma}

\begin{proof}
  see Exercise~\ref{ORTHO}).
\end{proof}

The Fourier coefficients of any function, defined in~\eqref{FOUCO}, can be rewritten with this notation as
\begin{equation}
  \label{GL3} \widehat f_k=\langle f,e^{2\pi ikx}\rangle_{L^2((0,1))}.
\end{equation}
The Fourier coefficients can thus be thought of as the $L^2$-projection of the function onto the complex exponentials.

Now, we make the following nice and instructive calculation:

\begin{lemma}\label{INS}
Let~$c_k$ be a sequence of complex numbers and~$f\in L^2((0,1))$.

Then,
\begin{equation}\label{jqdsc0fjn}
\left\| f-\sum_{{k\in\Z}\atop{|k|\le N}} c_k \,e^{2\pi i k x}\right\|_{L^2((0,1))}^2=\|f\|_{L^2((0,1))}^2+
\sum_{{k\in\Z}\atop{|k|\le N}} |\widehat f_k-c_k|^2-\sum_{{k\in\Z}\atop{|k|\le N}}|\widehat f_k|^2.
\end{equation}
\end{lemma}

\begin{proof} The linearity of the inner product (see Exercise~\ref{SCALA}) gives that the left-hand side of~\eqref{jqdsc0fjn} equals
\begin{eqnarray*}&& \langle f,f\rangle_{L^2((0,1))} -\sum_{{k\in\Z}\atop{|k|\le N}}\overline{c_k}
\langle f, e^{2\pi i k x}\rangle_{L^2((0,1))}\\&&\qquad-\sum_{{k\in\Z}\atop{|k|\le N}} c_k
\langle e^{2\pi i k x},f\rangle_{L^2((0,1))}+\sum_{{-N\le k\le N}\atop{-N\le j\le N}}c_k \,\overline{c_j}
\langle e^{2\pi i k x}, e^{2\pi i j x}\rangle_{L^2((0,1))}.
\end{eqnarray*}
This (recall~\eqref{fasv}, \eqref{GL2}, \eqref{GL3}) is equal to
\begin{equation}\label{0oiwejfP} \| f\|^2_{L^2((0,1))} -\sum_{{k\in\Z}\atop{|k|\le N}}\overline{c_k}\,\widehat f_k-\sum_{{k\in\Z}\atop{|k|\le N}} c_k\,\overline{\widehat f_k}+\sum_{{k\in\Z}\atop{|k|\le N}} |c_k|^2.
\end{equation}
Since
$$  |\widehat f_k-c_k|^2=\big(\widehat f_k-c_k\big)\big(\overline{\widehat f_k}-\overline{c_k}\big)=
|\widehat f_k|^2+|c_k|^2-\overline{c_k}\,\widehat f_k-c_k\,\overline{\widehat f_k},$$
the expression in~\eqref{0oiwejfP} equals the right-hand side of~\eqref{jqdsc0fjn}, as desired.
\end{proof}

As a consequence of Lemma~\ref{INS}, we have:
\begin{corollary}
Let~$f\in L^2((0,1))$.

Then,
\begin{equation} \label{Bty}
\left\| f-\sum_{{k\in\Z}\atop{|k|\le N}} \widehat f_k \,e^{2\pi i k x}\right\|_{L^2((0,1))}^2=\|f\|_{L^2((0,1))}^2-\sum_{{k\in\Z}\atop{|k|\le N}}|\widehat f_k|^2.
\end{equation}
\end{corollary}

\begin{proof}
It suffices to apply~\eqref{jqdsc0fjn} for the special choice of coefficients equal to~$\widehat f_k$.\end{proof}

\begin{corollary}
Let~$c_k$ be a sequence of complex numbers and~$f\in L^2((0,1))$.

Then,
\begin{equation}\label{B1RAS}
\left\| f-\sum_{{k\in\Z}\atop{|k|\le N}} c_k \,e^{2\pi i k x}\right\|_{L^2((0,1))}^2-\left\| f-\sum_{{k\in\Z}\atop{|k|\le N}} \widehat f_k \,e^{2\pi i k x}\right\|_{L^2((0,1))}^2=
\sum_{{k\in\Z}\atop{|k|\le N}} |\widehat f_k-c_k|^2.\end{equation}
\end{corollary}

\begin{proof} The identity in~\eqref{B1RAS} follows
by subtracting~\eqref{Bty} to~\eqref{jqdsc0fjn}.
\end{proof}

With this, we can provide the proof of Theorem~\ref{BEST}:

\begin{proof}[Proof of Theorem~\ref{BEST}] 
The quantity on the right-hand side of~\eqref{B1RAS} is non-negative, and in fact strictly positive unless~$c_k=\widehat f_k$ for all~$|k|\le N$. This yields the desired result.
\end{proof}

\begin{exercise}\label{LINC-bisco} Prove that the Fourier coefficients of a
linear combination of functions is the linear combination of the Fourier coefficients.

That is, assume that the functions~$f_1,\dots,f_\ell$ belong to~$L^1((0,1))$ and
are periodic of period~$1$, and let~$\alpha_1,\dots,\alpha_\ell\in\R$.

Define~$f:=\alpha_1f_1+\dots+\alpha_\ell f_\ell$ and prove that, for every~$k\in\N$,
$$ \widehat f_k=\alpha_1\widehat f_{1,k}+\dots+\alpha_\ell\widehat f_{\ell,k},$$
where, for each~$m\in\{1,\dots,\ell\}$, we have denoted by~$\widehat f_{m,k}$ the $k$th Fourier coefficient of~$f_m$.
\end{exercise}

\begin{exercise}\label{PKS0-3-21} The Fourier Series of a trigonometric polynomial is the trigonometric polynomial itself. More precisely, if
$$ f(x)=\sum_{{k\in\Z}\atop{|k|\le n}} c_k\,e^{2\pi ikx},$$
prove that, for all~$N\ge n$ and~$x\in\R$,
$$ \sum_{{k\in\Z}\atop{|k|\le N}} \widehat f_k\,e^{2\pi ikx}=f(x),$$
with $\widehat f_k=0$ whenever $k>n$.
\end{exercise}

\begin{exercise} \label{FBA} Prove that,
if~$f\in L^1((0,1))$ is a periodic function, for all~$k\in\Z$ we have that
$$ |\widehat f_k|\le\int_0^1 |f(x)|\,dx.$$\end{exercise}

\begin{exercise}\label{NONSepCFGDV}
Let~$f\in L^1((0,1))$ be a periodic function of period~$1$.

Show that if~$f$ is a non-negative function, then,
for all~$k\in\Z$,
\begin{equation}\label{12wewiatf}|\widehat f_k|\le\widehat f_0.\end{equation}
\end{exercise}

\begin{exercise}\label{NICELI}
One of the versions of Cauchy Integral Formula in complex analysis (see e.g.~\cite[Corollary~4.2 on page~47]{MR1976398})
states that if~$F$ is holomorphic
in some complex open disk of radius~$R$ centred at some point~$z_0\in\C$ and if~${\mathcal{C}}_r$
is a circle of radius~$r\in(0,R)$ centred at~$z_0$ and travelled anticlockwise, then, for all~$k\in\N$,
$$ D^k F(z_0)=\frac{k!}{2\pi i}\oint_{{\mathcal{C}}_r}\frac{F(z)}{(z-z_0)^{k+1}}\,dz.$$

Prove this formula by computing the Fourier coefficients of a suitable function.
\end{exercise}

\begin{exercisesk} \label{NONSepCFGDV2}
Let~$F$ be a holomorphic function in the complex unit disk and suppose that~$\Re F$ is a non-negative function. 

Prove that, for all~$k\in\N\cap[1,+\infty)$,
$$ |D^k F(0)|\le 2k!\; \Re F(0),$$
where~$D^k$ denotes the~$k$th derivative in the complex variable~$z$.\end{exercisesk}

\begin{exercise} \label{NONSepCFGDV2-cont}
Does the claim in Exercise~\ref{NONSepCFGDV2} hold true for~$k=0$?
\end{exercise}

\begin{exercise} \label{09iulappcoemmdba}
Let~$f_j$ be a sequence of functions, periodic of period~$1$, and belonging to~$L^1((0,1))$.
Let~$\widehat f_{j,k}$ be the $k$th Fourier coefficient of~$f_j$.

Suppose that~$f_j$ converges in~$L^1((0,1))$ to some function~$f$, periodic of period~$1$.

Prove that
$$ \widehat f_k=\lim_{j\to+\infty}\widehat f_{j,k}.$$
\end{exercise}

\begin{exercise}\label{ILTRI-FAA}
Let~$P\in C^\infty(\R)$. Prove that, for all~$\ell\in\N$ and~$x\in\R$,
$$ D^\ell_x\big(P(e^{2\pi ix})\big)=e^{2\pi i\ell x}\,D^\ell P(e^{2\pi ix})+\sum_{j=1}^{\ell-1} c_{\ell,j}\,e^{2\pi ijx}\,D^j P(e^{2\pi ix}),$$
for suitable coefficients~$c_{\ell,1},\dots,c_{\ell,\ell-1}$.
\end{exercise}

\begin{exercise} \label{ILTRI}
Let~$r\in\R$ and~$I=[r,r+1)$.
Prove that a trigonometric polynomial of degree~$N$ admits at most~$2N$ zeros 
(counted\footnote{For all~$k\in\N$, one says that~$x_\star\in\R$ is a zero of the trigonometric polynomial~$f$ of multiplicity~$k+1$ if all the derivatives of~$f$ up to the order~$k$ vanish at~$x_\star$.} with multiplicity) in the interval~$I$.
\end{exercise}

\begin{exercise}\label{ILTRI-22}
Is it true that a trigonometric polynomial of degree~$N$ admits exactly~$2N$ zeros (counted with multiplicity)?
\end{exercise}

\begin{exercise}\label{smc203e}
Prove that the expressions in~\eqref{TRIG-FO} and~\eqref{EXP-FO} are equivalent through the relation among coefficients given by
\begin{equation}\label{jasmx23er}
\begin{split} &
{\mbox{for all~$k\in\N$,}}\\&
a_k=\widehat f_k+\widehat f_{-k}=2\,\Re(\widehat f_k)=
2\int_0^1 f(x)\,\cos(2\pi kx)\,dx\\
&{\mbox{and,
for all~$k\in\N$ with~$k\ge1$,}}\\ &b_k=i(\widehat f_k-\widehat f_{-k})=-2\,\Im(\widehat f_k)=
2\int_0^1 f(x)\,\sin(2\pi kx)\,dx.\end{split}
\end{equation}
\end{exercise}

\begin{exercise}\label{smc203e2rf436b.2900rjm4on} Let~$f:\R\to\R$ be periodic of period~$1$ and $f$ belongs to $L^1((0,1))$.

Assuming that~$f$ is even (i.e., that~$f(x)=f(-x)$ for all~$x\in\R$),
prove that~$\widehat f_k=\widehat f_{-k}$ and that, in the notation of~\eqref{jasmx23er}, $b_k=0$.

Similarly, assuming that~$f$ is odd (i.e., that~$f(x)=-f(-x)$ for all~$x\in\R$), prove that~$\widehat f_k=-\widehat f_{-k}$ and that, in the notation of~\eqref{jasmx23er}, $a_k=0$.
\end{exercise}

\begin{exercise} \label{920-334PKSXu9o2fg}
Check that, in the notation of Exercise~\ref{smc203e},
the Fourier Series in~\eqref{FOSE} assumes the trigonometric form
\begin{equation}\label{TRISI} \frac{a_0}2+\sum_{k=1}^{+\infty} a_k\cos(2\pi kx)+b_k\sin(2\pi kx).\end{equation}
\end{exercise}

\begin{exercise} \label{920-334PKSXu9o2fgfbsmos}
Let~$f\in L^1((0,1))$ be periodic of period~$1$. Prove that if~$f$ is an even function then its
Fourier Series in trigonometric form~\eqref{TRISI} contains only cosines,
while if~$f$ is an odd function then its
Fourier Series in trigonometric form~\eqref{TRISI} contains only sines.
\end{exercise}

\begin{exercise}\label{SQ:W} Calculate the Fourier coefficients of the \index{square wave}
square wave. That is, for every~$x\in[0,1)$, let
$$ f(x):=\begin{dcases}
1&{\mbox{ if }}x\in\displaystyle\left[0,\frac12\right),\\
-1&{\mbox{ if }}x\in\displaystyle\left[\frac12,1\right).
\end{dcases}$$ Let~$f_{\text{per}}$ be the periodic extension of period~$1$ of~$f$, as defined in~\eqref{FPER}.
Prove that the~$k$th Fourier coefficient of~$f_{\text{per}}$ equals~$-\frac{2i}{\pi k}$ when~$k$ is odd, and vanishes otherwise.
\end{exercise}

\begin{exercise}\label{SA:W} Calculate the Fourier coefficients of the \index{sawtooth waveform}
sawtooth waveform. That is, for every~$x\in[0,1)$, let~$f(x):=x-\frac12$. Let~$f_{\text{per}}$ be the periodic extension of period~$1$ of~$f$, as defined in~\eqref{FPER}.
Prove that the~$k$th Fourier coefficient of~$f_{\text{per}}$ equals~$\frac{i}{2\pi k }$ when~$k\ne0$ and~$0$ when~$k=0$.
\end{exercise}

\begin{figure}[h]
\includegraphics[height=2.8cm]{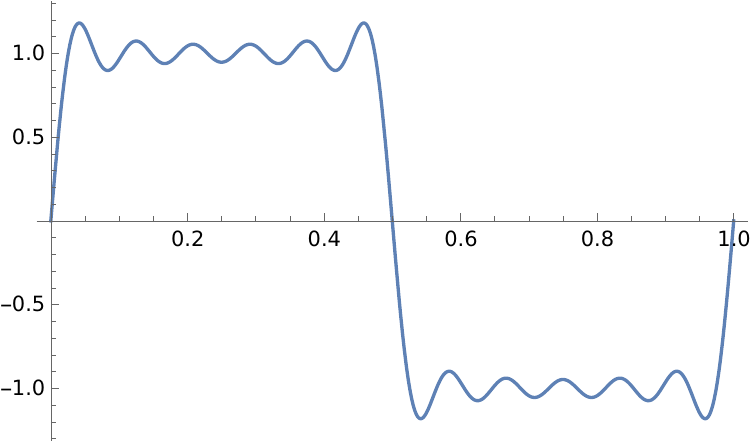}$\quad$\includegraphics[height=2.8cm]{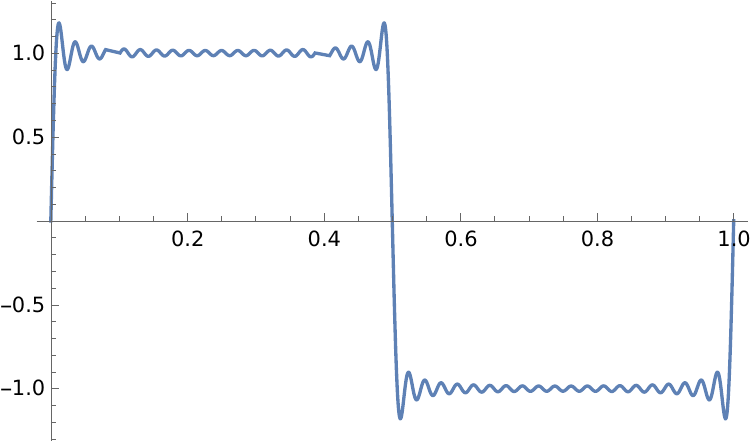}$\quad$\includegraphics[height=2.8cm]{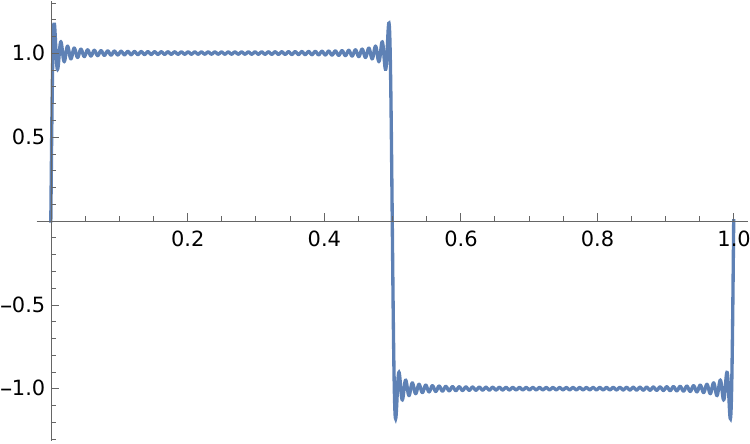}
\centering
\caption{Plot of a square wave approximation~$\displaystyle
\sum_{j=0}^{N} \frac4{\pi(2j+1)}\sin(2\pi(2j+1)x)$, with~$N\in\{5,20,50\}$.}\label{09iuygfy876trfdfghjjnn2meI20}
\end{figure}

\begin{exercise}\label{ojld03-12d} Write the Fourier Series of the square wave and of the sawtooth wave in trigonometric form. Specifically, check that the square wave in Exercise~\ref{SQ:W} can be written in the trigonometric form~\eqref{TRISI} as
$$\sum_{j=0}^{+\infty} \frac4{\pi(2j+1)}\sin(2\pi(2j+1)x)$$
and that the sawtooth wave\footnote{It is also interesting to compare the trigonometric
Fourier Series of the sawtooth wave with~\eqref{SBPF2-E2}.} 
in Exercise~\ref{SA:W} takes the form
$$-\sum_{k=1}^{+\infty}\frac1{\pi k}\sin(2\pi kx).$$
\end{exercise}

\begin{figure}[h]
\includegraphics[height=2.8cm]{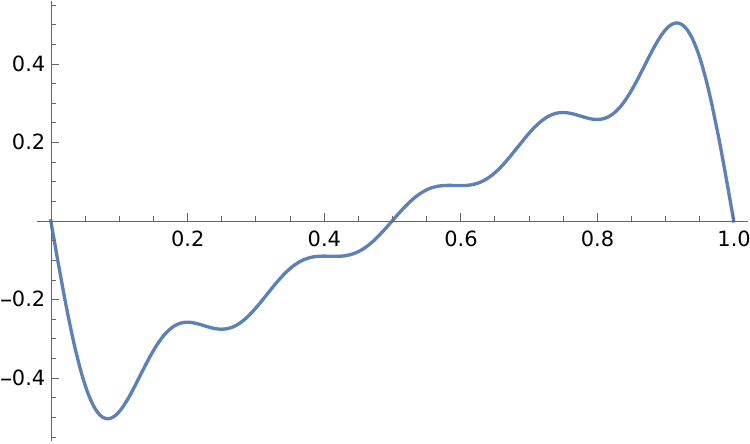}$\quad$\includegraphics[height=2.8cm]{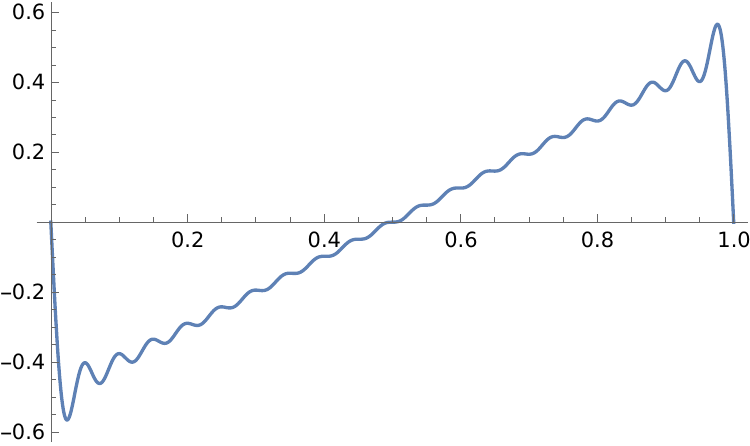}$\quad$\includegraphics[height=2.8cm]{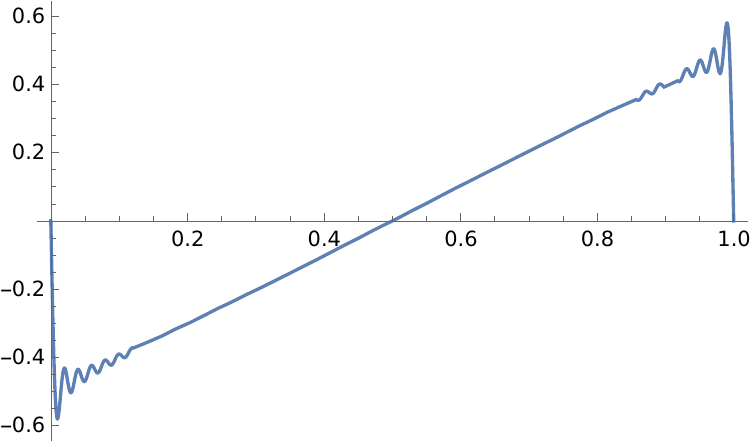}
\centering
\caption{Plot of a sawtooth wave approximation~$\displaystyle-\sum_{k=1}^{N}\frac1{\pi k}\sin(2\pi kx)$, with~$N\in\{5,20,50\}$.}
\label{09iuygfy876trfdfghjjnn2meI202}
\end{figure}

A visualisation of the situation described in Exercise~\ref{ojld03-12d} is given in Figures~\ref{09iuygfy876trfdfghjjnn2meI20} and~\ref{09iuygfy876trfdfghjjnn2meI202}.
These diagrams suggest that Fourier Series have some capacity of reconstructing the original waves,
up to suitable errors which seem to be more significant at the jumps of the original function.
A more solid understanding about how Fourier Series indeed provide an approximation procedure will come from the forthcoming Sections~\ref{CONVES0}, \ref{CONVESob}, \ref{CONVES}, \ref{UNIFORMCO:SECTI}, \ref{COL2}, \ref{COLP},
\ref{FEJERKESE}, \ref{SEC:GIBBS-PH} and~\ref{EXCE} (this will be indeed a rather subtle topic,
but of paramount importance).

\begin{exercise}\label{FO:DE:MAGIB} Let~$w$ be the square wave in Exercise~\ref{SQ:W}, considered as a periodic function of period~$1$. 

Calculate the derivative of the function$$\sum_{{k\in\Z}\atop{|k|\le N}}\widehat w_k\,e^{2\pi ikx}-w(x)$$ for all~$N\in\N$ and~$x\in\left(0,\frac12\right)$.
\end{exercise}

\begin{exercise} \label{ka-1} Calculate the Fourier coefficients of the \index{triangular wave}
triangular wave. That is,
for every~$x\in[0,1)$, let~$f(x):=\min\{x,1-x\}-\frac14$. Let~$f_{\text{per}}$ be the periodic extension of period~$1$ of~$f$, as defined in~\eqref{FPER}.
Prove that the~$k$th Fourier coefficient of~$f_{\text{per}}$ equals~$-\frac{1}{\pi^2 k^2}$ when~$k$ is odd and~$0$ otherwise.
\end{exercise}

\begin{exercise}\label{ka-2}
For every~$x\in[0,1)$, let~$f(x):=x(1-x)$. Let~$f_{\text{per}}$ be the periodic extension of period~$1$ of~$f$, as defined in~\eqref{FPER}.
Prove that the~$k$th Fourier coefficient of~$f_{\text{per}}$ equals~$-\frac{1}{2\pi^2 k^2}$ when~$k\ne0$ and~$\frac16$ when~$k=0$.
\end{exercise}

\begin{exercise} Write the Fourier Series of the functions in Exercises~\ref{ka-1} and~\ref{ka-2} in both exponential and trigonometric forms. 
\end{exercise}

\begin{exercise}\label{NU90i3orjf:12oeihfnvZMAS} For all~$x\in\left[-\frac12,\frac12\right)$ let
$$ f(x):=\begin{dcases}
x^2 & {\mbox{ if }}\displaystyle x\in\left(-\frac13,\frac13\right),\\ 0&{\mbox{ otherwise.}}\end{dcases}$$
Calculate the Fourier Series of~$f$ in trigonometric form.\end{exercise}

\begin{exercise}\label{OJSNILCESNpfa.sdwpqoed-23wedf}
Let~$f:\R\to\R$ be a continuous and periodic function of period~$1$.

Assume that~$f$ is differentiable in~$(0,1)$, with~$f'\in L^1((0,1))$. Prove that the Fourier Series of the function~$f'$ (periodically extended outside~$(0,1)$) coincides with the ``termwise'' derivative of the Fourier Series of~$f$, namely that
the Fourier Series of~$f'$ is
\begin{equation}\label{pqjdwlfmaonsveaN} \sum_{{k\in\Z}} 2\pi ik\widehat f_k\,e^{2\pi ikx}.\end{equation}
\end{exercise}

\begin{exercise}\label{OJSNILCESNpfa.sdwpqoed-23wedf.1}
In Exercise~\ref{OJSNILCESNpfa.sdwpqoed-23wedf}, is the assumption that~$f$ is continuous in~$\R$ really necessary
or can it be removed (replacing it with~$f\in L^1((0,1))$)?
\end{exercise}

\begin{exercise}\label{G0ilMajsx912e} Let~$f:\R\to\R$ be continuous and periodic of period~$1$.

Assume that~$f$ is differentiable in~$(0,1)$, with~$f'\in L^2((0,1))$.
Prove that
$$\sum_{k\in\Z}|\widehat f_k|<+\infty.$$
\end{exercise}

\begin{exercise}\label{G0ilMajsx912e-14} Let~$M\in\N$ and~$a_0,\dots,a_M\in[0,1]$ with~$0=a_0<a_1<\dots<a_{M-1}<a_M=1$.

Let~$f:\R\to\R$ be continuous and periodic of period~$1$.

Assume that~$f$ is differentiable in~$(a_{j-1},a_j)$ and~$f'\in L^2((a_{j-1},a_j))$ for each~$j\in\{1,\dots,M\}$.

Prove that
$$\sum_{k\in\Z}|\widehat f_k|<+\infty.$$
\end{exercise}

\begin{exercise}\label{LIFCOA}
Prove that the operation of taking Fourier coefficients is linear, in the sense that, if~$f_1,\dots,f_M\in L^1((0,1))$ and~$\alpha_1,\dots,\alpha_M\in\R$, then, for every~$k\in\Z$, the $k$th Fourier coefficient of the function~$\alpha_1 \,f_1+\dots+\alpha_M\,f_M$ is equal to~$\alpha_1 \,\widehat f_1+\dots+\alpha_M\,\widehat f_M$.\end{exercise}

\begin{exercise}\label{fasv2}
Prove formula~\eqref{fasv}.\end{exercise}

\begin{exercise}\label{ORTHO} Prove formula~\eqref{GL4}.\end{exercise}

\begin{exercise}\label{LADEBD}
Consider a sequence of complex numbers~$\{c_k\}_{k\in\Z}$ with~$c_{-k}=\overline{c_k}$ such\footnote{The condition that~$c_{-k}$ equals the complex conjugate of~$c_k$ is only needed to ensure that the function~$f$ takes values in the reals. If one is willing to consider complex-valued functions, such condition can be dropped.} that
\begin{equation}\label{AMs-ccg} \sum_{k\in\Z}|c_k|<+\infty.\end{equation}
Let
\begin{equation}\label{AMs-ccf} f(x):=\sum_{k\in\Z} c_k\,e^{2\pi ikx}.\end{equation}
Prove that necessarily~$c_k=\widehat f_k$.\end{exercise}

\begin{exercise}\label{FO:DE:MA}
Let
$$ \phi(x):=\max\{ 0,\,1-|x|\}.$$

For every~$a\in\left(0,\frac12\right)$ and~$x\in\left[-\frac12,\frac12\right)$, let also
$$ \phi_a(x):=\frac1a\,\phi\left(\frac{x}a\right).$$
Extend~$\phi_a$ to a function periodic of period~$1$ and calculate its Fourier Series.
\end{exercise}

\begin{exercise}\label{RESID} Consider the function
$$ \frac{\cos(2\pi x)+1}{\cos(2\pi x)+5\pi}$$
and write its Fourier Series in trigonometric form.
\end{exercise}

\begin{exercise}\label{AKxz3eMS-90394585865}
Write the Fourier series of the function
$$ \R\ni x\mapsto f(x):=\frac{1}{5-4\cos(2\pi x)}.$$
\end{exercise}

\begin{exercise}\label{ASKSMc-1}
Write the Fourier series of~$\cos^2 x$ in trigonometric form.
\end{exercise}

\begin{exercise}\label{JSOLTRA}
Prove that translating a function results in the Fourier coefficients being multiplied by a phase: namely,
given a real number~$a$ and a function~$f\in L^1((0,1))$ periodic of period~$1$, and setting
$$ T_a f(x):=f(x+a) ,$$
check that, for any~$k\in\Z$, the $k$th Fourier coefficient of~$T_a f$ is~$  e^{2\pi i k a}\,\widehat f_k$.\end{exercise}

\begin{exercise} Prove that, if~$g$, $h\in L^2((0,1),\C)$, then
$$ \langle h,g\rangle_{L^2((0,1))}=\overline{ \langle g,h\rangle_{L^2((0,1))}}.$$
\end{exercise}

\begin{exercise} \label{SCALA}
Prove that, if~$g_1,\dots,g_M$, $h_1,\dots,h_L\in L^2((0,1),\C)$ and~$\alpha_1,\dots,\alpha_M$, $\beta_1,\dots,\beta_L\in\C$, then
$$ \left\langle \sum_{k=1}^M \alpha_k\,g_k,\sum_{j=1}^L\beta_j\, h_j\right\rangle_{L^2((0,1))}=
\sum_{{1\le k\le M}\atop{1\le j\le L}} \alpha_k\,\overline{\beta_j}\,\langle g_k,h_j\rangle_{L^2((0,1))}.$$
\end{exercise}

\begin{exercise} \label{pecoOSJHN34jfgj}
Given $f$, $g\in L^1((0,1))$ two periodic functions of period~$1$, we define\footnote{The periodic convolution
is sometimes called also ``resultant'' \index{resultant}
or \index{Faltung}
(with a German word) ``Faltung'', see e.g.~\cite[page~10]{MR44660}.} their \index{periodic convolution}
periodic convolution as
\begin{equation}\label{pecoOSJHN34jfgj-1}
f\star g(x):=\int_0^1 f(x-y)\,g(y)\,dy.
\end{equation}
Prove that:
\begin{itemize}
\item $f\star g=g\star f$,
\item $f\star g\in L^1((0,1))$,
\item $f\star g$ is periodic of period~$1$,
\item for every~$k\in\Z$, the $k$th Fourier coefficient of~$f\star g$ is equal to the product~$\widehat f_k\,\widehat g_k$.\end{itemize}
\end{exercise}

\begin{exercise}\label{8ytfd23ertghn6rdsfghj0oi-P}
The classical convolution \index{convolution}
between two (non-periodic) functions~$f$ and~$g$ is defined by
\begin{equation}\label{pecoOSJHN34jfgj-2}
f* g(x):=\int_\R f(x-y)\,g(y)\,dy,
\end{equation}
provided the integral exists.

Suppose that~$f\in L^1((0,1))$ is a periodic function of period~$1$.
Let~$r\in\R$ and assume that~$g\in L^1((r,r+1))$ and that~$g$ vanishes outside~$(r,r+1)$.

Let~$\underline{g}$ be the periodic extension of~$g$ of period~$1$ outside~$(r,r+1)$, i.e. let~$\underline{g}:\R\to\R$ be such that~$\underline{g}=g$ in~$(r,r+1)$ and~$\underline{g}(x+1)=\underline{g}(x)$ for all~$x\in\R$.

Prove that~$f*g=f\star \underline{g}$, i.e. the classical convolution in~\eqref{pecoOSJHN34jfgj-2}
coincides with the \index{convolution, periodic}
periodic convolution in~\eqref{pecoOSJHN34jfgj-1}, up to identifying~$g$ with its periodic extension~$\underline{g}$.
\end{exercise}

\begin{exercise}\label{34rtghjd23ertghn6rdsfghj0oi-P}
Suppose that~$f\in L^1((0,1))$ is a periodic function of period~$1$.
Let~$r\in\R$ and assume that~$g\in L^1((r,r+1))$ and that~$g$ vanishes outside~$(r,r+1)$.
Let~$\underline{g}$ be the periodic extension of~$g$ outside~$(r,r+1)$.

Prove that the $k$th Fourier coefficient of~$f*g$ is equal to~$\widehat f_k\,\widehat{\underline{g}}_k$.
\end{exercise}

\begin{exercisesk} \label{FERR-2}
The\footnote{This result was actually first conjectured by Lip\'ot Fej\'er and later proven by Frigyes Riesz.
See~\cite[page~13]{MR545506} for further details about this.} Fej\'er-Riesz Theorem \index{Fej\'er-Riesz Theorem}
states that if~$f:\R\to[0,+\infty)$ is a trigonometric polynomial
of degree~$N$, then there exists a trigonometric polynomial~$g$ of degree~$N$ such that~$f=|g|^2$.

Prove this result.\end{exercisesk}

\section{Some further observations in~$L^2((0,1))$}\label{SFOL2}

Some berries are delicious, but they end up always at the bottom of the bowl.
Mathematical results are sometimes like berries, one has to find them at the bottom of the bowl. The following two results are perhaps of this sort, since they follow by scraping the bottom of the bowl provided in Section~\ref{L2pe}.

The first result is known by the name of \emph{Bessel's Inequality} \index{Bessel's Inequality}:

\begin{theorem} Let~$f\in L^2((0,1))$.

Then,
\begin{equation}\label{AJSa}
\sum_{k\in\Z} |\widehat f_k|^2 \leq \|f\|^2_{L^2((0,1))}. 
\end{equation}
\end{theorem}

\begin{proof} It suffices to observe that the left-hand side of~\eqref{Bty} is non-negative
and then send~$N\to+\infty$.
\end{proof}

The second result is a useful equivalence:

\begin{proposition}\label{AJSaa}
Let~$f\in L^2((0,1))$.

Then, the identity
\begin{equation}\label{PARS}
    \sum_{k\in\Z} |\widehat f_k|^2 = \|f\|^2_{L^2((0,1))}
\end{equation}
holds true if and only if
\begin{equation}\label{1PCc}
\lim_{N\to+\infty} \left\| f-\sum_{{k\in\Z}\atop{|k|\le N}} \widehat f_k \,e^{2\pi i k x}\right\|_{L^2((0,1))}=0.
\end{equation}\end{proposition}

Formula~\eqref{PARS} is known by the name of \emph{Parseval's Identity} \index{Parseval's Identity}
(basically, it is the equality case in Bessel's Inequality~\eqref{AJSa})
and it is sometimes considered as an ``energy identity'' (the energy of the function in the right-hand side of~\eqref{PARS} being equal to the energy of its coefficients in the left-hand side).

The interest of Proposition~\ref{AJSaa} is that it says that the validity of
Parseval's Identity~\eqref{PARS} is equivalent to the convergence
of the Fourier Sum \index{Fourier Sum}
\begin{equation}\label{FourierSum}
S_{N,f}(x):=\sum_{{k\in\Z}\atop{|k|\le N}} \widehat f_k \,e^{2\pi i k x}
\end{equation}
as~$N\to+\infty$ to the original function~$f$ in the space~$L^2((0,1))$, as stated in~\eqref{1PCc}.

The notation~$S_{N,f}$ for the Fourier Sum will be utilised throughout this chapter
(and its limit as~$N\to+\infty$ will be called Fourier Series). When no ambiguity arises, we will use the short notation
$$S_N:=S_{N,f}.$$

We will come back to results in the spirit of Proposition~\ref{AJSaa} in the forthcoming Section~\ref{COL2}.

\begin{proof}[Proof of Proposition~\ref{AJSaa}] By~\eqref{Bty},$$
\lim_{N\to+\infty}\left\| f-\sum_{{k\in\Z}\atop{|k|\le N}} \widehat f_k \,e^{2\pi i k x}\right\|_{L^2((0,1))}^2=\|f\|_{L^2((0,1))}^2-\sum_{k\in\Z}|\widehat f_k|^2,$$
from which the desired result follows.
\end{proof}

\begin{exercise}\label{A6KMSD345678ifgb}
Prove that, for all~$k\in\Z$ with~$|k|\le N$, the $k$th Fourier coefficient of~$S_{N,f}$ equals that of~$f$.
\end{exercise}

\begin{exercise}\label{G:PIVO}
Show that horizontal and vertical translations \index{translations}
do not affect the pointwise convergence of a Fourier Series
(thus allowing one to reduce the pointwise convergence problem at the origin).

More precisely, prove the following statement.

Let~$f\in L^1((0,1))$ be periodic of period~$1$.

Given~$x_0\in\mathbb{R}$, let~$g(x):=f(x+x_0)-f(x_0)$.

Prove that
\begin{itemize}
\item $g(0)=0$,
\item we have that
$$ \widehat g_k=\begin{cases} \widehat f_0-f(x_0) & {\mbox{ if }}k=0,\\
e^{2\pi i k x_0}\,\widehat f_k  & {\mbox{ if }}k\ne0,
\end{cases}$$
\item we have that, for every~$x\in\R$, $$ S_{N,g}(x)=S_{N,f}(x+x_0)-f(x_0),$$
\item we have that, for every~$x\in\R$, 
$$ S_{N,g}(x)-g(x)=S_{N,f}(x+x_0)-f(x+x_0),$$
\item we have that~$S_{N,f}$ converges to~$f$ as~$N\to+\infty$
at the point~$x_0$ if and only if~$S_{N,g}$ converges to~$g$ as~$N\to+\infty$ at the point~$0$.\end{itemize}\end{exercise}

\begin{exercise}\label{SBPF} Prove the \index{summation by parts}
following ``summation by parts'' formula, valid for all sequences~$\{\alpha_k\}_{k\in\N}$, $\{\beta_k\}_{k\in\N}$ and any~$n$, $m\in\N$ with~$m\le n$:
\begin{equation}\label{BYPA} \sum_{k=m}^{n}\alpha_{k}(\beta_{k+1}-\beta_{k})=\left(\alpha_{n}\beta_{n+1}-\alpha_{m}\beta_{m}\right)-\sum_{k=m+1}^{n}\beta_{k}(\alpha_{k}-\alpha_{k-1}).\end{equation}
\end{exercise}

\begin{exercisesk}\label{SBPF2} Prove that, for every~$x\in(0,1)$,
\begin{equation}\label{SBPF2-E1} \sum_{k=1}^{+\infty}\frac{\cos(2\pi kx)}k=-\frac12\Big(\ln 2+\ln(1-\cos(2\pi x))\Big)\end{equation}
and
\begin{equation}\label{SBPF2-E2}\sum_{k=1}^{+\infty}\frac{\sin(2\pi kx)}k=\pi\left(\frac12- x\right).\end{equation}
\end{exercisesk}

\begin{exercise} \label{NOCAP}
  Prove that the closed unit ball in~$L^2((0,1))$ is not compact.
\end{exercise}

\begin{exercise} \label{LEGEPO}
For any~$k\in\N$, the \index{Legendre polynomials} Legendre polynomials are defined by
$$ P_k(x):=\frac{\sqrt{2k+1}}{2^{k+\frac12}\,k!}\left(\frac{d}{dx}\right)^k(x^2-1)^k.$$
Prove that these polynomials form an orthonormal system in~$(-1,1)$, that is, for every~$k$, $m\in\N$,
$$ \int_{-1}^1 P_k(x)\,P_m(x)\,dx=\delta_{k,m}.$$
\end{exercise}

\begin{exercise} \label{HAARBA} For every~$n\in\N$ and~$k\in\N$ with~$k\in[1,2^n]$ let
$$ h_{k,n}(x):=
\begin{dcases}
2^{n/2} & {\mbox{ if }}x\in\displaystyle\left(\frac{k-1}{2^n},\frac{k-\frac12}{2^n}\right),\\
-2^{n/2} & {\mbox{ if }}x\in\displaystyle\left(\frac{k-\frac12}{2^n},\frac{k}{2^n}\right),\\
0&{\mbox{ otherwise.}}
\end{dcases}$$
These functions are called \index{Haar functions} Haar functions.

Prove that they form an orthonormal system in~$(0,1)$, that is, for every~$k_1$, $k_2$, $n_1$, $n_2\in\N$,
with~$k_1\in[1,2^{n_1}]$ and~$k_2\in[1,2^{n_2}]$,
\begin{equation}\label{po09iuyhglsqd-2.1} \int_{0}^1 h_{k_1,n_1}(x)\,h_{k_2,n_2}(x)\,dx=\delta_{k_1,k_2}\,\delta_{n_1,n_2}.\end{equation}
\end{exercise}

\begin{exercise}\label{CLA:AUNIFO} Let~$a\in\left(0,\frac12\right)$. Prove that the series $$ \sum_{k=1}^{+\infty}\frac{\cos(2\pi kx)}{k}\qquad{\mbox{and}}\qquad
\sum_{k=1}^{+\infty}\frac{\sin(2\pi kx)}{k}$$ converge uniformly for~$x\in[a,1-a]$.\end{exercise}

\begin{exercise}\label{1-CLA:AUNIFO}
Prove that the series
$$ \sum_{k=2}^{+\infty}\frac{\sin(2\pi kx)}{k\,\ln k}\qquad{\mbox{and}}\qquad
\sum_{k=2}^{+\infty}\frac{\cos(2\pi kx)}{k\,\ln k}$$
converge uniformly in~$\R$.
\end{exercise}

\begin{exercise}\label{MASCHE}
The \index{Euler-Mascheroni constant}
Euler-Mascheroni constant\footnote{A major unsolved problem is to decide whether the Euler-Mascheroni constant is irrational (and, if so, whether it is transcendental).} is defined by
$$ \gamma :=\int_{1}^{+\infty }\left({\frac{1}{\lfloor x\rfloor }}-{\frac{1}{x}}\right)\,dx.$$
Prove that
$$ \lim_{x\to+\infty} \ln x+\gamma-\sum_{k=1}^{\lfloor x\rfloor}\frac1k=0.$$
\end{exercise}

\begin{exercise}\label{ABEDI1} The following criterion is known as \index{Abel's Test}
Abel's Test.
Let~$\{\sigma_k\}_{k\in\N}$ be a sequence of complex numbers and~$\{\tau_k\}_{k\in\N}$ be a sequence of real numbers.

Assume that the series
$$\sum_{k=0}^{+\infty} \sigma_k$$ is convergent.

Assume also that~$\tau_k$ is monotone and bounded.

Then, the series
$$ \sum_{k=1}^{+\infty }\sigma_k\,\tau_k$$is convergent.

Prove this statement.\end{exercise}

\begin{exercise}\label{ABEDI2} The following criterion is known as Dirichlet's Test. \index{Dirichlet's Test}
Let~$\{\sigma_k\}_{k\in\N}$ be a sequence of complex numbers and~$\{\tau_k\}_{k\in\N}$ be a sequence of real numbers.

Assume that there exists~$M\ge0$ such that, for each~$N\in\N$,$$\left|\sum_{k=0}^{N}\sigma_k\right|\leq M.$$

Assume also that~$\tau_k$ is monotone and infinitesimal as~$k\to+\infty$.

Then, the series$$ \sum_{k=1}^{+\infty }\sigma_k\,\tau_k$$is convergent.

Prove this statement.\end{exercise}

\section{Decay of Fourier coefficients}\label{DECAY:e:sFOL}

In this section, we discuss the decay of the Fourier coefficients\footnote{Of course, using~\eqref{jasmx23er},
the decay of~$\widehat f_k$ is equivalent to that of the corresponding trigonometric coefficients~$a_k$ and~$b_k$.}
as defined in~\eqref{FOUCO}. 

Understanding the decay of Fourier coefficients is a very difficult art, but it is of pivotal importance: after all, one of the main objectives of Fourier methods is to approximate functions with a Fourier Series. With this respect, if we want the Fourier Series in~\eqref{FOSE} to converge, say, at the point~$x=0$, we need
$$ \sum_{k\in\Z} \widehat f_k $$
to converge, for which a necessary condition is that
$$\lim_{k\to\pm\infty} \widehat f_k=0.$$

The first result in this direction is the following one,
which is often\footnote{When this result is used in further generality,
  it is also sometimes \index{Mercer's Theorem}
  called Mercer's Theorem, see e.g.~\cite[Section~2.4]{MR44660}.} called the \index{Riemann-Lebesgue Lemma}
\emph{Riemann-Lebesgue Lemma}:

\begin{theorem}\label{RLjoqwskcdc}
Let~$f\in L^1((0,1))$. Then,
\begin{equation} \label{SMXC}\lim_{k\to\pm\infty} \widehat f_k=0.\end{equation}
\end{theorem}

The intuition\footnote{The Riemann-Lebesgue Lemma already suffices to deduce in one shot the convergence of Fourier Series for sufficiently regular functions, which in turn is one of the cornerstones of the theory, ensuring that Fourier Series can indeed reconstruct any arbitrary (sufficiently nice) function. The reader who wishes to go straight to this feature is very warmly {\em encouraged to jump directly to Section~\ref{CONVESob} and then come back here}, where we construct a suitable architecture for a full-bodied system of ideas.
Going first to Section~\ref{CONVESob} and then come back is, in our opinion, pedagogically very appropriate, because it allows one to fix the convergence issue (at least in a provisional way) as soon as possible, for peace of mind, and the rest of the theory will fall into its proper position later on.}
behind this result is that~$\widehat f_k$ is obtained by integrating~$f$ against a very oscillatory function (as~$|k|$ gets larger and larger). This oscillatory function, of period~$\frac1k$, takes opposite values, say~$1$ and~$-1$, at very close-by points, say~$x_0$ and~$x_0+\frac{1}{2k}$.
So, when computing~$\widehat f_k$, the value of~$f$ at~$x_0$ and the value of~$f$ at~$x_0+\frac1{2k}$ are multiplied by opposite values of the highly oscillating function.
If~$f$ is continuous, the value of~$f$ at~$x_0$ and that at~$x_0+\frac1{2k}$ are ``almost'' the same, hence these two contributions ``almost'' cancel out, and this cancellation effect becomes more and more relevant as $|k|$ gets larger and larger.

Interestingly, not only this argument can be made rigorous, but also it can be extended to discontinuous functions as well, thanks to an approximation argument.
The details\footnote{See Exercises~\ref{S91uejfnv10moOrEo023.4fSIAL1}, \ref{S91uejfnv10moOrEo023.4fSIAL2}, and~\ref{S91uejfnv10moOrEo023.4fSIAL3} for alternative proofs of Theorem~\ref{RLjoqwskcdc}. } go as follows:

\begin{proof}[Proof of Theorem~\ref{RLjoqwskcdc}] First of all, we observe that
\begin{equation}\label{SMXC2}
{\mbox{if~$f\in L^2((0,1))$, then~\eqref{SMXC} holds true.}}\end{equation}
Indeed, if~$f\in L^2((0,1))$ we can use
Bessel's Inequality~\eqref{AJSa} to see that $$
\sum_{k\in\Z} |\widehat f_k|^2<+\infty.$$
{F}rom this, we obtain~\eqref{SMXC2}, as desired.

The result in~\eqref{SMXC2} proves Theorem~\ref{RLjoqwskcdc} when~$f\in L^2((0,1))$,
which is a stronger assumption than~$f\in L^1((0,1))$ (due to H\"older's Inequality), so we need now to consider the case in which~$f$ only belongs to~$L^1((0,1))$.

For this, we decompose~$f$ in the form
\begin{equation}\label{8id2}
  f=g_j+h_j,
\end{equation}
with sequences of functions $\{g_j\}_{j\in \N}$ and $\{h_j\}_{j\in \N}$ such that ~$g_j\in L^\infty((0,1))$ (so in particular~$g_j\in L^2((0,1))$)
and~$h_j$ such that
\begin{equation}
  \label{ladjja}
  \lim_{j\to+\infty}\int_0^1 |h_j(x)|\,dx=0.
\end{equation}
See e.g. Exercise~\ref{LASESQ} for such a decomposition.

Thus, by~\eqref{SMXC2},
if we denote by~$\widehat g_{j,k}$ the $k$th Fourier coefficient of~$g_j$,
\begin{equation}\label{kasdwe} \lim_{k\to\pm\infty} \widehat g_{j,k}=0.\end{equation}

Besides, using~\eqref{ladjja},
\begin{equation}\label{ksmdmdmf}\begin{split}& \lim_{j\to+\infty} \lim_{k\to\pm\infty}|\widehat h_{j,k}|=
\lim_{j\to+\infty} \lim_{k\to\pm\infty}
\left|\int_0^1 h_j(x)\,e^{-2\pi i kx}\,dx\right|\\&\qquad\le\lim_{j\to+\infty} \lim_{k\to\pm\infty}
\int_0^1 |h_j(x)|\,dx=\lim_{j\to+\infty}
\int_0^1 |h_j(x)|\,dx=0.\end{split}\end{equation}

Also, by~\eqref{8id2} and the linearity of the operation of taking Fourier coefficients (recall Exercise~\ref{LIFCOA}) we conclude that
\begin{equation}\label{kasdwe-b} \widehat f_k=\widehat g_{j,k}+\widehat h_{j,k}.\end{equation}

From~\eqref{kasdwe}, \eqref{ksmdmdmf} and~\eqref{kasdwe-b}, we arrive at
\begin{eqnarray*}&& 
\lim_{k\to\pm\infty}|\widehat f_k|=
\lim_{j\to+\infty} \lim_{k\to\pm\infty}|\widehat f_k|\le
\lim_{j\to+\infty} \lim_{k\to\pm\infty}\big( |\widehat g_{j,k}|+|\widehat h_{j,k}|\big)=0,\end{eqnarray*}
as desired.
\end{proof}

We stress that the converse of Theorem~\ref{RLjoqwskcdc} is false: there are infinitesimal sequences which are not the Fourier coefficients of any locally integrable periodic function, see Exercise~\ref{ESIPNo}
(but the matter is delicate, since ``small'' modifications can produce rather different results,
compare with Exercise~\ref{ESIPNo-2}).

Theorem~\ref{RLjoqwskcdc} is pivotal, as we will appreciate in the convergence proof that will be presented for
Theorem~\ref{ONEGO},
but it does not give any estimate on the rate of decay of the Fourier coefficients: and it cannot do so, since, without further assumptions, the decay to zero of the Fourier coefficients can be arbitrarily slow,
see Exercise~\ref{FTCPASNo2}.

For more regular functions, however, one can estimate the decay of the Fourier coefficients, as pointed out in the following result.

\begin{theorem}\label{SMXC22} Let~$m\in\N$.
Let~$f\in C^m(\R)$ be periodic of period~$1$.

Then, for all~$k\in\Z$,
\begin{equation} \label{CARGi0} |\widehat f_k|\le \frac{\displaystyle\sup_{x\in\R}|D^m f(x)|}{|2k\pi |^m}.\end{equation}
\end{theorem}

We remark that the order of decay provided by Theorem~\ref{SMXC22} is somewhat optimal (see e.g. Exercise~\ref{0-1-303kuHNAsq1-1}; see however Exercises~\ref{BEDD-pre} and~\ref{BEDD}
for suitable refinements and modifications of Theorem~\ref{SMXC22}).

Theorem~\ref{SMXC22} is actually part of an interesting observation linking Fourier coefficients and derivatives, which we present here:

\begin{theorem}\label{SMXC22b}
Let~$m\in\N$. Let~$f\in C^m(\R)$ be periodic of period~$1$.

Then, for all~$k\in\Z$,
\begin{equation}\label{CARGi}  \widehat{D^m f}_k=(2\pi i k)^m\widehat f_k.\end{equation}
\end{theorem}

\begin{proof}[Proof of Theorems~\ref{SMXC22} and~\ref{SMXC22b}] We start by proving~\eqref{CARGi}.
For this, we argue by induction over~$m$. When~$m=0$, the claim in~\eqref{CARGi} boils down to the definition of Fourier coefficient in~\eqref{FOUCO}. Suppose now that~\eqref{CARGi} holds true for some~$j\in\{0,1,\dots,m-1\}$. Then, we define~$g:=D^jf$,
we use the inductive hypothesis and integrate by parts to see that, for all~$k\in\Z$,
\begin{eqnarray*}&& (2\pi i k)^j\widehat f_k=\widehat{D^jf}_k=\widehat g_k=\int_0^1 g(x)\,e^{-2\pi i kx}\,dx=
-\frac{1}{2\pi ik} \int_0^1 g(x)\,D(e^{-2\pi i kx})\,dx\\&&\qquad=
-\frac{1}{2\pi ik}\left(  
  g(1)\,e^{-2\pi i k}-g(0)\,e^{0}-
\int_0^1 D g(x)\,e^{-2\pi i kx}\,dx\right)\\&&\qquad=
-\frac{1}{2\pi ik}\left(  
D^jf(1)-D^jf(0)-
\int_0^1 D^{j+1}f(x)\,e^{-2\pi i kx}\,dx\right)\\&&\qquad=0+
\frac{1}{2\pi ik}\int_0^1 D^{j+1}f(x)\,e^{-2\pi i kx}\,dx=\frac{\widehat{D^{j+1}f}_k}{2\pi ik}.
\end{eqnarray*}
This establishes~\eqref{CARGi} for the index~$j+1$ and thereby completes the inductive step.
Thus, we have proved~\eqref{CARGi}.

Now, we prove~\eqref{CARGi0}. To this end, we use~\eqref{CARGi} and we see that
$$ |\widehat f_k|=\frac{|\widehat{D^m f}_k|}{|2k\pi |^m}.$$
Then, (see Exercise~\ref{FBA}, applied here to the function~$D^mf$),
$$ |\widehat f_k|\le\frac{1}{|2k\pi |^m}\int_0^1 |{D^m f}(x)|\,dx,$$
from which we obtain~\eqref{CARGi0}.
\end{proof}

Results such as Theorem~\ref{SMXC22} highlight the intimate interplay between 
local features (such as smoothness of a function) and global features (such as \label{OjndFNSmgA}
rapid decay of its Fourier coefficients). This ``local-global duality'' occurs very often in Fourier analysis.

The relation between the decay of Fourier coefficients and regularity of the function will be revisited in Section \ref{DEC2}.

\begin{exercise}
Give an alternative proof of~\eqref{SMXC2} by using Theorem~\ref{SMXC22}.
\end{exercise}

\begin{exercise}\label{LASESQ}
Let~$f\in L^1((0,1))$. Construct sequences of functions~$\{g_j\}_{j\in\N}$ and~$\{h_j\}_{j\in \N}$ satisfying
$$ f=g_j+h_j$$
such that~$g_j\in L^\infty((0,1))$ and~$
h_j\in L^1((0,1))$ with
$$ \lim_{j\to+\infty}\int_0^1 |h_j(x)|\,dx=0.$$
\end{exercise}

\begin{exercise}\label{BEDD-pre} Let~$m\in\N$ and~$f\in C^m(\R)$ be periodic of period~$1$. Prove that, for every~$k\in\Z$,
$$ |\widehat f_k|\le\frac{\|D^mf\|_{L^1((0,1))}}{|2\pi k|^m}.$$\end{exercise}

\begin{exercise}\label{BEDD} Let~$m\in\N$ and~$f\in C^m(\R)$ be periodic of period~$1$. Prove that
$$ \sum_{k\in\Z}|2k\pi |^{2m}\,|\widehat f_k|^2\le\int_0^1 |D^m f(x)|^2\,dx.$$\end{exercise}

\begin{exercise}\label{0-1-303kuHNAsq1-1} Let~$m\in\N$, $\alpha\in(0,1]$ and $$f(x):=
\sum_{\ell=0}^{+\infty} \frac{\cos(2^{\ell+1}\pi x)}{2^{\ell(m+\alpha)}}.$$
Prove that~$f\in C^m(\R)$ and periodic of period~$1$,
and that the decay of the Fourier coefficients of~$f$ is no faster than~$\frac{1}{|k|^{m+\alpha}}$, namely
$$ \limsup_{k\to\pm\infty} |k|^{m+\alpha}\,|\widehat f_k|>0.$$\end{exercise}

\begin{exercise}\label{S91uejfnv10moOrEo023.4fSIAL1}
Give an alternative proof of the Riemann-Lebesgue Lemma in Theorem~\ref{RLjoqwskcdc} by using the density of smooth functions in Lebesgue spaces
(see e.g.~\cite[Theorem~9.6]{MR3381284}).
\end{exercise}

\begin{exercise}\label{S91uejfnv10moOrEo023.4fSIAL2}
Give an alternative proof of the Riemann-Lebesgue Lemma in Theorem~\ref{RLjoqwskcdc} by using the density of step functions in Lebesgue spaces
(see e.g.~\cite[Theorems~2.26 and~2.41]{MR1681462}).
\end{exercise}

\begin{exercise}\label{S91uejfnv10moOrEo023.4fSIAL3}
Give an alternative proof of the Riemann-Lebesgue Lemma in Theorem~\ref{RLjoqwskcdc} by using
the continuity of the translations in Lebesgue spaces,
i.e. the fact that
$$ \lim_{\epsilon\searrow0}\int_0^1 \big|f(x-\epsilon)-f(x)\big|^p\,dx=0,$$
see e.g.~\cite[Theorem~8.19]{MR3381284}.
\end{exercise}

\section{The quest for convergence}\label{CONVES0}

The convergence of Fourier Series is not obvious at all.
To stress its conceptual difficulty, let us recall Fourier's citation, as reported in~\cite[Chapter~1]{MR1970295}:
\begin{quote}
  {\em ``Regarding the researches of d'Alembert and Euler could one not add that if they knew this expansion, they \label{CITA-F} made but a very imperfect use of it. They were both persuaded that an arbitrary and discontinuous function could never be resolved in series of this kind, and it does not even seem that anyone had developed a constant in cosines of multiple arcs, the first problem which I had to solve in the theory of heat''.}
\end{quote}

That is, while it may be intuitive that Fourier Series can reproduce well {\em some} periodic functions,
it is, in our opinion, far less obvious, or not obvious at all, that they can in fact reproduce
{\em any, rather arbitrary, periodic function} (with some relevant
exceptions that will be discussed in Section~\ref{EXCE}),
and one has all the rights of being sceptical\footnote{The possibility, offered by Fourier Series, of recovering functions with jumps and corners from infinite sums of complex exponentials is indeed remarkable and surprising, because the sum of {\em finitely many} complex exponentials (i.e., of finitely many sines and cosines) is a smooth (and even analytic) function. Therefore, the great flexibility of the Fourier method lies precisely in its limit process and at that time contributed also to highlight the conceptual difference between a finite sum and an {\em infinite} series.}
about the possibility, say, of approximating
a piecewise constant function, or a sawtooth function, with complex exponentials (i.e., with sines and cosines).

In a sense, the difficulty of reproducing a function from an approximate sum surfaces already when considering Taylor Series, \label{ANL} since it is not obvious that one can obtain functions such as~$e^x$,
$\sin x$, $\sin(\cos(e^x))$, etc. out of an infinite sum of monomials. As we will see, the situation regarding Fourier Series is even more subtle, since for power series a standard method is to use the Triangle Inequality to ``place the absolute values inside the summation'', while the convergence of Fourier Series critically depends on a finer analysis of oscillations and cancellations.

To highlight the conceptual difficulty of Fourier Series, let us remark that when~$f$ is smooth enough
(e.g., at least~$C^2(\R)$)
we already know that, for every~$x\in\R$, the Fourier Sum~$S_N(x)$, as defined in~\eqref{FourierSum}, converges as~$N\to+\infty$
(simply because, by Theorem~\ref{SMXC22}, its $k$th term is bounded by~$\frac{C}{|k|^2}$, for some~$C\ge0$, which is the $k$th term of a convergent series).
But so what? The real issue is to decide {\em whether or not~$S_N(x)$ converges to the right value~$f(x)$} (or, more generally, whether or not~$S_N$ approaches~$f$ in some sense as~$N\to+\infty$).

After all, it is very easy, given a function~$f$, to construct a converging series out of it (for instance, one could associate to any~$f$ a sequence of vanishing coefficients, in this way the associated series would simply be zero): but the interesting feature is whether the series can recover the original function~$f$
(and, if so, in which sense).

The general theory\footnote{As a spoiler, and to reassure the reader about the fact that
mathematicians do not completely waste their time,
let us mention that Fourier Series are appropriately designed
to ``converge to the appropriate object when they converge'' (e.g. almost everywhere), see Theorem~\ref{ILGIOB}. 
Also, a sufficiently fast decay of the Fourier coefficients does ensure the convergence (in fact, the uniform convergence)
of the Fourier Series to the right function, see Theorem~\ref{BASw}, so Fourier did invent a decent theory after all.
These facts are not completely trivial though, and in some cases crazy things can happen, see Section~\ref{EXCE}
(and there are scary examples of Fourier Series diverging everywhere, see~\cite{zbMATH02586338}).} of the convergence of Fourier Series is very broad and complex, so we will not aim at being exhaustive, but rather at presenting, hopefully in an accessible way, some results of general use.

For starters, in Section~\ref{CONVESob}, we show a sufficiently straightforward proof of the pointwise convergence of Fourier Series for ``sufficiently nice'' functions. A full-bodied structure of this theory will be presented later in Section~\ref{CONVES}.
In spite of its ``provisional'' character, the gist of Section~\ref{CONVESob} is to show that, while convergence issues are deep and very delicate, some powerful ideas can be implemented in a rather direct way, without major prerequisites. Also, these ideas can be suitably expanded and rethought for a deeper understanding and the construction of a comprehensive system of great impact (but sometimes at the cost of some additional technicalities: the goal of Section~\ref{CONVESob} is instead to avoid additional technicalities as much as possible, without losing coherence and identity).

The discussion about the convergence of Fourier Series will be further expanded in Section~\ref{UNIFORMCO:SECTI},
to present sufficient conditions ensuring uniform convergence, and in Section~\ref{COL2},
to deal with convergence in~$L^2((0,1))$ (the case of~$L^p((0,1))$ will be sketched in Section~\ref{COLP}). The intriguing case of jump discontinuities will be dealt with in Section~\ref{SEC:GIBBS-PH}.

As a useful digression, we will also highlight in Section~\ref{FEJERKESE} how to use ``averaging procedures'' to obtain
an approximation of a given function by trigonometric polynomials obtained as a mean of Fourier Sums.

Also, since no theory is really complete in the absence of negative results, we will present some
convergence issues in Section~\ref{EXCE}.

\section{Pointwise convergence of Fourier Series ``all in one breath''}\label{CONVESob}

Now we show, in one go, that:

\begin{theorem}\label{ONEGO}
Let~$f\in C^1(\R)$ be periodic of period~$1$.

Then, for all~$x\in\R$,
$$ \lim_{N\to+\infty}S_N(x)=f(x).$$
\end{theorem}

Finer results will be provided in the forthcoming sections (and the impatient reader can completely jump this section and proceed with the rest). The goal here is to present a proof which is, possibly, as simple as it can be and highlight some of the main ideas.

\begin{proof}[Proof of Theorem~\ref{ONEGO}]
Without loss of generality, we can focus on the convergence of~$S_{N}$ to $f$ at the point~$0$
under the additional assumption that~$f(0)=0$ (see Exercise~\ref{G:PIVO}: indeed, one could replace~$f$ by the function~$g$ constructed there).

Hence,
\begin{equation}\label{ojqlwd2o-9939}\begin{split}&
S_N(0)-f(0) = S_N(0)=
\sum_{{k\in\Z}\atop{|k|\le N}} \widehat f_k = \sum_{{k\in\Z}\atop{|k|\le N}}\int_0^1 f(x) \,e^{-2\pi i k x} \,dx\\&\qquad
= \int_0^1 f(x) \sum_{{k\in\Z}\atop{|k|\le N}} e^{-2\pi i kx} \,dx.
\end{split}\end{equation}

Moreover, exploiting the geometric sum, for every~$r\in\R\setminus\{1\}$, 
\begin{equation*}\begin{split}&
\sum_{{k\in\Z}\atop{|k|\le N}} r^k=
\sum_{k=0}^N r^k+\sum_{k=0}^N \left(\frac1r\right)^k-1=
\frac{1-r^{N+1}}{1-r}+\frac{1-(1/r)^{N+1}}{1-(1/r)}-1\\&\qquad=
\frac{1-r^{N+1}}{1-r}+\frac{r-(1/r)^{N}}{r-1}-1=
\frac{r^{-N}-r^{N+1}}{1-r},\end{split}
\end{equation*}
and in fact this holds true for~$r=1$ too, provided that we consider the latter expression in a limit sense.

Substituting for~$r:=e^{-2\pi ix}$, we find that
\begin{equation}\label{KASMqwdfed123er}
\sum_{{k\in\Z}\atop{|k|\le N}} e^{-2\pi i kx}= \frac{e^{2\pi iNx}-e^{-2\pi i(N+1)x}}{1-e^{-2\pi ix}}.
\end{equation}
Plugging this information into~\eqref{ojqlwd2o-9939} and defining
\begin{equation}\label{hcde} h(x):=\frac{f(x)}{1-e^{-2\pi ix}},\end{equation}
we arrive at
\begin{equation}\label{RieRLjoqwskcdc}
S_N(0)-f(0) = \int_0^1 h(x)\,\Big(e^{2\pi iNx}-e^{-2\pi i(N+1)x}\Big)\, dx=
\widehat h_{-N}-\widehat h_{N+1}.
\end{equation}

Now, using de l'H\^{o}pital's Rule,
$$ \lim_{x\to0}h(x)=\frac{f'(0)}{2\pi i}$$
and moreover~$h$ is periodic of period~$1$.

As a result, $h$ is continuous\footnote{We stress that the continuity of~$h$ was used here merely to deduce that~$h\in L^1((0,1))$.
Also the smoothness condition on~$f$ was only used to deduce the continuity of~$h$.\label{LAFOOD}

Hence, if one could find any condition (alternative to the smoothness of~$f$) that guarantees that~$h\in L^1((0,1))$, the rest of the proof would go through. This is indeed the strategy that we will implement in the forthcoming proof of Theorem~\ref{DINITS}.} and, in particular, in~$L^1((0,1))$. Hence, by the Riemann-Lebesgue Lemma (see Theorem~\ref{RLjoqwskcdc}),
$$ \lim_{N\to\pm\infty}\widehat h_N=0.$$
We combine this with~\eqref{RieRLjoqwskcdc} and we deduce that
$$ \lim_{N\to\pm\infty} \Big(S_N(0)-f(0) \Big)=\lim_{N\to\pm\infty}\Big(
\widehat h_{-N}-\widehat h_{N+1}\Big)=0-0=0,$$
showing the desired convergence.
\end{proof}

The careful reader may have already spotted that perhaps the assumptions on~$f$ in Theorem~\ref{ONEGO} could be weakened:
after all, what we need in its proof is only that~$h$ belongs to~$L^1((0,1))$
(the fact that we prove that~$h$ is even continuous may be superfluous, recall footnote~\ref{LAFOOD}).
We will show that indeed one can relax the regularity requirement, replacing it with
a sharper one that will be presented in~\eqref{DINI} below.

Indeed, roughly speaking, the function~$h$ in~\eqref{hcde}, for small~$x$, captures the behaviour of~$\frac{f(x)}{x}$ or, better to say (since here~$f(0)=0$), of~$\frac{f(x)-f(0)}x$. Again, roughly speaking,
this suggests that a suitable control on the derivative of~$f$ (or, more specifically, of its incremental quotient) guarantees the convergence of its Fourier Series.
This will indeed be discussed in detail in the forthcoming Theorem~\ref{DINITS} and Exercises~\ref{DINI1}, \ref{DINI2} and~\ref{DINI3}.

Also, the main structural assumption of Theorem~\ref{ONEGO} is of ``global'' nature, namely that~$f$ is {\em everywhere} continuously differentiable, but the thesis of Theorem~\ref{ONEGO} is somewhat of ``local'' nature, namely convergence at a given point (though, at any given point), so one could ambitiously wonder whether it is possible to obtain local theses out of local hypotheses. For instance, whether one could obtain the convergence of the Fourier Series in a given interval by only knowing that~$f$ is, say,
continuously differentiable there. This is indeed the case
(a positive answer to this question following from Theorem~\ref{DINITS}, see Exercise~\ref{DINITS-EXEC1}),
but we stress that these kinds of ``localisation properties'' are sometimes very far from being intuitive, because the Fourier coefficients involve the function in its entirety, being an integral on the whole period, and it is quite remarkable that one can actually ``average out''
far away contributions and obtain local convergence statements of Fourier Series (in fact, the interplay between local and non-local phenomena is a very intriguing aspect of Fourier methods).

Moreover, though very useful, formula~\eqref{KASMqwdfed123er} is aesthetically questionable, \label{QSTNB}
since its left-hand side is symmetric if one exchanges~$N$ with~$-N$ but the right-hand side is not.
Of course this is only a cosmetic observation, but one could guess that this loss of symmetry in the formula is due to the fact that the number of indexes in the left-hand side (summing from~$-N$ to~$N$) equals~$2N+1$, which gets split in the right-hand side into two indices~$N$ and~$N+1$. It is thereby tempting to conjecture that maybe splitting into two pairs, both equal to~$N+\frac12$, would have led to a more symmetric, more pleasant to look at, formula.
This is indeed the case, and we will give credit to aesthetic in the forthcoming Lemma~\ref{KASMqwdfed123erDKLI}.

\begin{exercise}\label{SPDCD-0.01}
Let~$\alpha$, $\beta\in\R$.

Prove that
$$ \cos(\alpha-\beta)-\cos(\alpha+\beta)=2\sin\alpha\sin\beta.$$
\end{exercise}

\begin{exercise}\label{SPDCD-0.01anc}
Let~$\alpha$, $\beta\in\R$.

Prove that
$$ \sin(\alpha-\beta)+\sin(\alpha+\beta)=2 \sin\alpha \cos\beta.$$
\end{exercise}

\begin{exercise}\label{SPDCD-0.02}
Let~$N$, $\ell\in\N$.

Prove that, for all~$x\in\R$,
$$ \sum_{{j\in\Z\setminus\{0\}}\atop{|j|\le\ell}}\frac{\cos(2\pi(N-j)x)}{j}=
2\sin(2\pi Nx)\sum_{j=1}^{\ell}\frac{\sin(2\pi jx)}j .$$\end{exercise}

\section{Pointwise convergence of Fourier Series: a refined version}\label{CONVES}

We now give a more general version of Theorem~\ref{ONEGO}, in which the assumption that the function~$f$ is continuously differentiable is replaced by a weaker one (this can be useful to treat non-differentiable functions):

\begin{theorem}\label{DINITS}
Let~$f:\R\to\R$ be periodic of period~$1$. Assume that~$f\in L^1((0,1))$.

Let~$x\in\R$ and suppose that there exists~$\delta>0$ such that \begin{equation}\label{DINI}
\int_{-\delta}^{\delta}\left|\frac{f (x + t) - f(x)}{t}\right|\,dt<+\infty.\end{equation}

Then,
$$ \lim_{N\to+\infty}S_N(x)=f(x).$$
\end{theorem}

Assumption~\eqref{DINI} (as well as some variations of it)
is known in the jargon\footnote{In the literature, there are other fine criteria for the convergence of
Fourier Series, such as 
\emph{de la Vall\'ee-Poussin's Criterion},
\emph{Dirichlet Criterion},
\emph{Jordan Criterion},
\emph{Lebesgue Criterion},
\emph{Young Criterion}, etc.: we will not dive here into all these conditions and on their mutual relations, but we refer the interested reader to~\cite[Sections~4.4--4.7]{MR44660},
\cite[Chapter~III]{MR171116}, \cite[Section 10]{MR545506}
and the references therein for further details on these criteria.

We also point out that Dini's Condition is sharp, since when this assumption is dropped
Fourier Series can diverge, see~\cite[Theorem~2.4]{MR1963498} and~\cite{MR612724}. \label{EXCE-foodf}

We will also discuss convergence problems of Fourier Series in Section~\ref{EXCE}.} as \index{Dini's Condition} \emph{Dini's Condition} (or \emph{Dini's Criterion}, or \emph{Dini's Test}) at the point~$x$.
We stress that ``sufficiently regular'' functions automatically satisfy this condition
(see Exercises~\ref{DINI1}, \ref{DINI2}, \ref{DINI3}, \ref{DINI-addC-1}, and~\ref{DINI-addC-2}).

\begin{proof}[Proof of Theorem~\ref{DINITS}]
As discussed at the beginning of the proof of Theorem~\ref{ONEGO},
without loss of generality, we can suppose that the point~$x$ at which we are analysing the convergence of the Fourier Series is~$0$
and that~$f(0)=0$. Thus, condition~\eqref{DINI} becomes
\begin{equation}\label{DINIZE}
\int_{-\delta}^{\delta}\left|\frac{f (t) }{t}\right|\,dt<+\infty.\end{equation}
Thus, as observed in footnote~\ref{LAFOOD},
to prove the desired convergence in Theorem~\ref{DINITS},
it would suffice to consider the function~$h$ in~\eqref{hcde}
and prove that~\eqref{DINIZE} implies that~$h\in L^1((0,1))$.
Actually, since~$h$ is periodic, the latter would be equivalent to
\begin{equation}\label{DINITS-ca}
h\in L^1\left(\left(-\frac12,\frac12\right)\right).
\end{equation}
To check this, we observe that, for all~$x\in\left(-\frac12,\frac12\right)$,
\begin{eqnarray*}&& |1-e^{-2\pi ix}|^2 = (1-\cos(2\pi x))^2+(\sin(2\pi x))^2=2\big(1-\cos(2\pi x)\big)\\&&\qquad=4\sin^2(\pi x)=
4\sin^2(\pi |x|)\end{eqnarray*}
and therefore, since
\begin{eqnarray*}&& 2\sin (\pi |x|)=2\int_0^{\pi|x|}\cos\vartheta\,d\vartheta\ge
\begin{dcases}\displaystyle
2\int_0^{\pi/4}\cos\vartheta\,d\vartheta=\sqrt{2} & {\mbox{ if }}|x|\in\left[\frac14,\frac12\right),\\
\displaystyle 2\int_0^{\pi|x|}\cos\frac\pi4\,d\vartheta =\sqrt{2}\pi|x|& {\mbox{ if }}|x|\in\left[0,\frac14\right),
\end{dcases}
\end{eqnarray*}
we deduce that$$|1-e^{-2\pi ix}|\ge c\,|x|,$$
for some~$c>0$, e.g. $c=2\sqrt{2}$.

Hence, recalling the definition of~$h$ in~\eqref{hcde}, for all~$x\in\left(-\frac12,\frac12\right)$,
$$ |h(x)|=\frac{|f(x)|}{|1-e^{-2\pi ix}|}\le \frac{|f(x)|}{c\,|x|}.$$

As a consequence, defining~$I:=\left(-\frac12,\frac12\right)\setminus(-\delta,\delta)$,
$$ \int_{-1/2}^{1/2}|h(x)|\,dx\le
\int_{-\delta}^{\delta}\frac{|f(x)|}{c\,|x|}\,dx
+\int_{I}\frac{|f(x)|}{c\,\delta}\,dx
\le\int_{-\delta}^{\delta}\frac{|f(x)|}{c\,|x|}\,dx
+\frac{\|f\|_{L^1((0,1))}}{c\,\delta}<+\infty,
$$
thanks to~\eqref{DINIZE}, which establishes~\eqref{DINITS-ca}, as desired.
\end{proof}

Regarding results such as Theorem~\ref{DINITS}, it may be surprising that these statements have a ``local nature'', \label{LAMS}
namely they make structural assumptions, such as~\eqref{DINI}, only at a given point~$x$ at which one checks the desired convergence: indeed, as mentioned above, the Fourier coefficients depend ``globally'' on the function
(due to their integral definition in~\eqref{FOUCO}) and
any slight modification of the function (even if away from the given point~$x$)
does change the Fourier coefficients and therefore could heavily impact the convergence of the Fourier Series at the point~$x$.
Instead, quite remarkably, the convergence of Fourier Series is a local property, as clarified in Theorem~\ref{C1uni-c} below.

\begin{exercise}\label{DINI1} Show that if $f$ is locally H\"older continuous \index{H\"older continuous}
at $x$, i.e. if there exist~$\alpha\in(0,1)$,
$\delta>0$ and~$M\ge0$ such that, for all~$h\in(-\delta,\delta)$,
$$ |f (x + h) - f(h)|\le M|h|^\alpha,$$ then~$f$ satisfies Dini's Condition~\eqref{DINI} at the point~$x$.
\end{exercise}

\begin{exercise}\label{DINI2} Show that if $f$ is locally \index{Lipschitz continuous}
Lipschitz continuous at $x$, i.e. if there exist~$\delta>0$ and~$M\ge0$ such that, for all~$h\in(-\delta,\delta)$,
$$ |f (x + h) - f(h)|\le M|h|,$$ then~$f$ is locally H\"older continuous at $x$.
\end{exercise}

\begin{exercise}\label{DINI3} Show that if $f$ is of class~$C^1$ (i.e., continuous and with continuous first derivative) in a neighbourhood of a point~$x\in\R$, then~$f$ is locally Lipschitz continuous at $x$.
\end{exercise}

\begin{exercise}\label{DINI-addC-1}
Let~$f:\R\to\R$ be periodic of period~$1$. Assume that~$f\in L^1((0,1))$ and that it satisfies Dini's Condition~\eqref{DINI} at some point~$x\in\R$, for some~$\delta>0$. 

Can~$f$ present a jump discontinuity at the point~$x$?\end{exercise}

\begin{exercise}\label{DINI-addC-2}
Let~$f:\R\to\R$ be periodic of period~$1$. Assume that~$f\in L^1((0,1))$ and that it satisfies Dini's Condition~\eqref{DINI} at some point~$x\in\R$, for some~$\delta>0$. 

Does this mean that~$f$ is necessarily continuous\footnote{It will be interesting to compare Exercise~\ref{DINI-addC-2} with the forthcoming Exercise~\ref{DINI-addC-2-maunif}.} at the point~$x$?\end{exercise}

\begin{exercise}\label{DINITS-EXEC1} Let~$f\in L^1((0,1))$ be periodic of period~$1$. Let~$x\in\R$ and~$\delta>0$.  Suppose that~$f$ is continuously differentiable in~$(x-\delta,x+\delta)$. Prove that the Fourier Sum of~$f$ at~$x$ converges to~$f(x)$. 
\end{exercise}

\section{Uniqueness results}\label{SEC:UNIQ:3}

A useful, and not completely trivial, observation is that Fourier coefficients
unambiguously identify a function, namely that if two functions possess the same Fourier coefficients then the two functions must necessarily coincide:

\begin{theorem}\label{UNIQ}
Let~$f$, $g$ be periodic functions of period~$1$, with~$f$, $g\in L^1((0,1))$.

Assume that~$\widehat f_k=\widehat g_k$ for all~$k\in\Z$.

Then, $f(x)=g(x)$ for all~$x\in\R$, up to a set of null Lebesgue measure.
\end{theorem}

\begin{proof} We let~$h:=f-g$ and our goal is thus to show that
\begin{equation}\label{FGOA}
{\mbox{$h(x)=0$ for all~$x\in\R$, up to a set of null Lebesgue measure.}}\end{equation}
To this end, we use the linearity of the Fourier coefficients (see Exercise~\ref{LIFCOA}) to see that, for all~$k\in\Z$,
\begin{equation}\label{EZED} \widehat h_k=\widehat f_k-\widehat g_k=0.\end{equation}
Now, if~$h$ were regular enough to satisfy Dini's condition~\eqref{DINI}, then we would be done, because
Theorem~\ref{DINITS} would have guaranteed that, for all~$x\in\R$,
$$ h(x)=\lim_{N\to+\infty}S_{N,h}(x)=\lim_{N\to+\infty}\sum_{{k\in\Z}\atop{|k|\le N}} \widehat h_k\,e^{2\pi ikx}=0.$$

The catch is that we do not know that~$h$ is regular enough, so we cannot apply this argument directly.
However, we can reduce to it by a convolution argument (recall~\eqref{pecoOSJHN34jfgj-2} for the definition of convolution). Namely, for every~$\epsilon\in(0,1)$, to be taken as small as we wish in what follows,
we consider a non-negative compactly supported\footnote{The support of a function is the subset of the domain over which the function takes a nonzero value. A function is compactly supported when the closure of the support is a compact subset of the domain. For our purposes and in the situation where a function is defined on an open interval of the real line such as $(-\frac{1}{2},\frac{1}{2})$, requesting that the function is compactly supported implies that the function is zero outside a closed (bounded) subset of the domain.} function~$\psi\in C^\infty_0\left(\left(-\frac12,\frac12\right)\right)$ with
$$ \int_{-1/2}^{1/2}\psi(x)\,dx=1$$
and define
$$ \psi_\epsilon(x):=\frac1\epsilon\, \psi\left(\frac{x}\epsilon\right).$$
Notice that~$\psi_\epsilon\in C^\infty_0\left(\left(-\frac12,\frac12\right)\right)$ and
$$ \int_{-1/2}^{1/2}\psi_\epsilon(x)\,dx=1.$$
Hence (see e.g.~\cite[Theorem~9.6]{MR3381284}) we know that the function~$h_\epsilon:=h*\psi_\epsilon$ converges to~$h$ in~$L^1((0,1))$ as~$\epsilon\searrow0$.
Therefore (see e.g. \cite[Theorem 4.9(a)]{MR2759829}), up to a sub-sequence,
\begin{equation}\label{Amsl022345}
  {\mbox{$h_\epsilon\to h$ a.e. in~$\R$.}}
\end{equation}

We denote by~$\underline\psi_\epsilon$ the periodic extension of~$\psi_\epsilon$ outside~$\left(-\frac12,\frac12\right)$.
Let also~$\widehat{ \underline\psi}_{\epsilon,k}$ the $k$th Fourier coefficient of~$\underline\psi_\epsilon$
and~$\widehat h_{\epsilon,k}$ the $k$th Fourier coefficient of~$h_\epsilon$. 
Then (see Exercise~\ref{34rtghjd23ertghn6rdsfghj0oi-P}) we have that
$$ \widehat h_{\epsilon,k}=\widehat h_k\,\widehat{ \underline\psi}_{\epsilon,k}.$$
This and~\eqref{EZED} give that, for all~$k\in\Z$,
\begin{equation}\label{Amsl02} \widehat h_{\epsilon,k}=0.\end{equation}

We also remark that~$h_\epsilon$ is a smooth function (see e.g.~\cite[Theorem~9.3]{MR3381284}) hence, as pointed out above, in this situation, we can apply Theorem~\ref{DINITS} and deduce from~\eqref{Amsl02} that~$h_\epsilon(x)=0$ for every~$x\in\R$.

From this and~\eqref{Amsl022345}, we obtain~\eqref{FGOA}, as desired.
\end{proof}

In a nutshell, Theorem~\ref{UNIQ} says that if the
Fourier coefficients of an integrable function are all equal to zero, then the function itself is equal to zero (up to sets of null measure): this property is sometimes stated using the sentence \index{completeness of the trigonometric system}
``the trigonometric system is complete''.

Another proof of Theorem~\ref{UNIQ} will be proposed in Exercise~\ref{UNIQ:ALAmowifj30othgore0-1}.

See e.g.~\cite[Theorems~97--100]{MR44660} and~\cite[Section~7, Chapter~I]{MR2039503} for more general uniqueness results.

\begin{exercise}\label{ojqdwn23EEE} Let~$m\in\N$, $m\ge1$ and
$$ f(x):=\sum_{k=1}^{+\infty} \frac1{k^{m+1}}\cos(2\pi kx).$$
Prove that~$f$ is periodic of period~$1$ and~$f\in C^{m-1}(\R)$, but~$f\not\in C^m(\R)$.
\end{exercise}

\begin{exercise} \label{anjqnxn}
Let~$f\in L^1((0,1))$. Assume that, for every~$n\in\N$,
$$ \int_0^1 f(x)\,x^n\,dx=0.$$
Prove that~$f(x)=0$ for a.e.~$x\in\R$.
\end{exercise}

\begin{exercise}\label{PKS0-3-21bismo} 
If the Fourier Series of a function is a trigonometric polynomial
(i.e., it has only finitely many terms), the original function coincides with this trigonometric polynomial.
More precisely, given~$f\in L^1((0,1))$, periodic of period~$1$, if~$ \widehat f_k=0$ for all~$k\in\Z$ with~$|k|\ge N+1$,
prove that
$$ f=S_{N,f}.$$
\end{exercise}

\begin{exercise}\label{PKS0-3-21bismo2} 
Prove that the Fourier Series of a function~$f\in L^1((0,1))$ contains only finitely many terms
if and only if the function~$f$ is a trigonometric polynomial.\end{exercise}

\begin{exercise}\label{PKS0-3-21bismo2.NECE} Consider a sequence of functions~$f_j\in L^1((0,1))$, periodic of period~$1$.

Suppose that there exists~$M\in\N$ such that the Fourier coefficients~$\widehat f_{j,k}$ of~$f_j$
satisfy, for all~$k\in\Z$ with~$|k|\ge M$, that
\begin{equation}\label{ICDCMMANHSMCGIEVIMN} \lim_{j\to+\infty}\widehat f_{j,k}=0.\end{equation}
Assume that~$f_j$ converges uniformly in~$\R$ to some function~$f$. 

Prove that~$f$ is necessarily a trigonometric polynomial.
\end{exercise}

\begin{exercise}\label{KPSLM.01e2ourfhvb0odc.2ewrfevgbx3xdnfv0-1}
Let~$m\in\N$ and
$$ f_m(x):=\prod_{j=0}^m \cos(2\pi jx).$$
Prove that, for every~$x\in\R$,
$$ f_m(x)=\sum_{{\sigma_0,\dots,\sigma_m\in\{-1,1\}}}
\frac{e^{2\pi i\phi_m(\sigma) x}}{2^m},$$
where, for all~$\sigma=(\sigma_0,\dots,\sigma_m)\in\{-1,1\}^{m+1}$,
\begin{equation}\label{phimwlxc} \phi_m(\sigma):=\sum_{k=0}^m k\sigma_k.\end{equation}
\end{exercise}

\begin{exercise}\label{phimwlxc2}
Let~$m\in\N$, $\sigma=(\sigma_0,\dots,\sigma_m)\in\{-1,1\}^{m+1}$,
and~$\phi_m(\sigma)$ be as in~\eqref{phimwlxc}.
Assume that~$\phi_m(\sigma)=0$. Prove that either~$\frac{m}4\in\N$ or~$\frac{m-3}4\in\N$.
\end{exercise}

\begin{exercise}\label{phimwlxc3} Let~$m\in\N$, $\sigma=(\sigma_0,\dots,\sigma_m)\in\{-1,1\}^{m+1}$, and~$\phi_m(\sigma)$ be as in~\eqref{phimwlxc}.

Assume that~$\frac{m}4\in\N$ and that, for all~$j\in\{0,\dots,\ell-1\}$,
\begin{equation}\label{SIKMSDgmwiedj1} \sigma_{4j+1}=1,\qquad\sigma_{4j+2}=-1, \qquad\sigma_{4j+3}=-1\qquad{\mbox{ and }}\qquad\sigma_{4j+4}=1.\end{equation}
Prove that~$\phi_m(\sigma)=0$.

On the same note, assume that~$\frac{m-3}4\in\N$ and that,
for all~$j\in\{0,\dots,\ell\}$,
\begin{equation}\label{SIKMSDgmwiedj12} \sigma_{4j}=1,\qquad\sigma_{4j+1}=-1, \qquad\sigma_{4j+2}=-1\qquad{\mbox{ and }}\qquad\sigma_{4j+3}=1.\end{equation}
Prove that~$\phi_m(\sigma)=0$.
\end{exercise}

\begin{exercise} \label{LAVppnUNOMASMAJONC0989okj744124IJN}
This is a variation of a question posed in the 1985 William Lowell Putnam Mathematical Competition
(see~\cite[pages~58--59]{MR1933844}).

Let~$m\in\N$ and
\begin{equation}\label{LAVppnUNOMASMAJONC0989okj744124IJN.0}
\R\ni x\longmapsto f_m(x):=\prod_{j=0}^m \cos(2\pi jx).\end{equation}
Check that the function~$f_m$ is continuous and periodic of period~$1$ and write the Fourier Series of~$f_m$ in trigonometric form as
\begin{equation}\label{LAVppnUNOMASMAJONC0989okj744124IJN.7}\frac{a_{0}}2+\sum_{k=1}^{+\infty}\Big(a_{k}\cos(2\pi kx) + b_{k}\sin(2\pi kx)\Big).\end{equation}

Prove that:
\begin{eqnarray}\label{LAVppnUNOMASMAJONC0989okj744124IJN.1}
&&{\mbox{$b_{k}=0$ for every~$k\in\N$, $k\ge1$,}}\\ 
\label{LAVppnUNOMASMAJONC0989okj744124IJN.2}&&{\mbox{only finitely many~$a_{k}$'s are nonzero,}}\\ {\mbox{and }}
&& \label{LAVppnUNOMASMAJONC0989okj744124IJN.3} \frac{a_{0}}2+\sum_{k=1}^{+\infty}a_{k}=1.
\end{eqnarray}

Moreover, for which values of~$m$ is~$a_0\ne0$?
\end{exercise}

\section{More on the decay of Fourier coefficients}
\label{DEC2}
We return to the decay of Fourier coefficients now that we have uniqueness results about Fourier Series.

As a further example of the duality between decay of Fourier coefficients and regularity of the function, we show that if the Fourier coefficients of~$f$ decay fast enough, then the function~$f$ is regular:

\begin{theorem}\label{ojqdwn23E}
  Let~$m\in\N$ and~$\epsilon>0$. Let~$f:\R\to\R$ be periodic of period~$1$ and $f\in L^1((0,1))$.

  Assume that there exists~$C\ge0$ such that,
for all~$k\in\Z\setminus\{0\}$,
\begin{equation}\label{0909b}|\widehat f_k|\le\frac{C}{|k|^{m+1+\epsilon}}.
\end{equation}

Then, $f\in C^m(\R)$.
\end{theorem}

We point out that Theorem~\ref{ojqdwn23E} is optimal, in the sense that one cannot take~$\epsilon=0$
(see Exercises~\ref{ojqdwn23EEE},~\ref{ojqdwn23EEE0921ef}, and \ref{ojqdwn23EEE0921ef-LN}).
To prove Theorem~\ref{ojqdwn23E}, we point out the following auxiliary result:

\begin{lemma}\label{Le-ojqdwn23E}
Let~$m\in\N$ and let~$\{c_k\}_{k\in\Z}$ be a sequence of complex numbers such that~$c_{-k}=\overline{c_k}$ and
\begin{equation}\label{0909}\sum_{k\in\Z}|k|^m |c_k|<+\infty.\end{equation}

Then, there exists a function~$g\in C^m(\R)$ which is periodic of period~$1$ such that
for all $x\in \R$
\begin{equation}
  \label{geqseries}
\lim_{N\to+\infty}\sum_{{k\in\Z}\atop{|k|\le N}} c_k\,e^{2\pi ikx}=g(x).  
\end{equation}
Moreover, the convergence in \eqref{geqseries} is uniform in $\R$.

Also, $c_k=\widehat g_k$ and, for all~$j\in\{0,1,\dots,m\}$,
\begin{equation}\label{CELLOD}
D^jg(x)=\lim_{N\to+\infty}\sum_{{k\in\Z}\atop{|k|\le N}} (2\pi ik)^j c_k\,e^{2\pi ikx}.
\end{equation}
\end{lemma}

\begin{proof} Let
$$ \phi_{0,N}(x):=\sum_{{k\in\Z}\atop{|k|\le N}} c_k\,e^{2\pi ikx}.$$
For~$j\in\{0,1,\dots,m\}$, let
$$ \phi_{j,N}(x):=D^j\phi_{0,N}(x)=\sum_{{k\in\Z}\atop{|k|\le N}} (2\pi ik)^j c_k\,e^{2\pi ikx}.$$
By~\eqref{0909} and the Limit Comparison Test, for each~$j\in\{0,1,\dots,m\}$ and $\epsilon>0$, there exists $N_\epsilon\in\N$ such that
\begin{equation*}
  \sum_{k\in\Z\atop |k|\geq N_\epsilon} |k|^j|c_k| < \epsilon.
\end{equation*}
As a result, for every $N>M\geq N_\epsilon$, we have
\begin{align*}
   \sup_{x\in\R} |\phi_{j,N}-\phi_{j,M}|
  = \sup_{x\in\R}\left|\sum_{{k\in\Z}\atop{M+1\le |k|\le N}} (2\pi ik)^j c_k e^{2\pi i k x}\right|
  \leq  (2\pi)^j \sum_{k\in\Z\atop |k|\geq N_\epsilon} |k|^j |c_k| < (2\pi)^j\epsilon,
\end{align*}
meaning that $\phi_{j,N}$ converges uniformly in $\R$ as~$N\to+\infty$.

In fact, we have that~$\phi_{0,N}$ converges uniformly to some~$g\in C^m(\R)$ (see e.g.~\cite[Theorem 7.17]{MR385023}),
with~$\phi_{j,N}$ converging to~$D^j g$ as~$N\to+\infty$, for all~$j\in\{0,1,\dots,m\}$.

As a result (see Exercise~\ref{LADEBD}) we have that~$c_k=\widehat g_k$.
\end{proof}

\begin{proof}[Proof of Theorem~\ref{ojqdwn23E}] Let~$c_k:=\widehat f_k$. Then, condition~\eqref{0909} holds true, thanks to~\eqref{0909b}. From Lemma~\ref{Le-ojqdwn23E}, there exists $g\in C^m(\R)$ such that $\widehat{g}_k=c_k=\widehat f_k$. By Theorem \ref{UNIQ}, $f$ and $g$ coincide up to a set of null Lebesgue measure, thus $g\in C^m(\R)$ is a representative of the class of $f\in L^1((0,1))$.
\end{proof}

In view of Theorems~\ref{SMXC22} and~\ref{ojqdwn23E}, it would be desirable to have a perfect characterisation of spaces such as~$C^m(\R)$ (say, for periodic functions) in terms of the Fourier coefficients, but this, unfortunately, is not possible, since there is a discrepancy between the decays in~\eqref{CARGi0} and~\eqref{0909b}. Actually, for a perfect correspondence between decay and regularity in a suitable sense, one has to take into account slightly more sophisticated functional spaces (see e.g.~\cite[Section~5.8.4]{MR2597943}) and we do not address this point here. However, Fourier coefficients characterise well~$C^\infty(\R)$, as noticed by the following result:

\begin{theorem}\label{0CINFITH1}
Let~$f\in L^1((0,1))$ be periodic of period~$1$.

Then, $f\in C^\infty(\R)$ if and only if for every~$m\in\N$ there exists~$C_m\ge0$ such that, for every~$k\in\Z\setminus\{0\}$,
\begin{equation}\label{m3edfgyui02woeiujfhnceu}
|\widehat f_k|\le\frac{C_m}{|k|^m}.
\end{equation}
\end{theorem}

\begin{proof} Suppose that~$f\in C^\infty(\R)$. Then, for all~$m\in\N$, we have that~$f\in C^m(\R)$, whence, in view of Theorem~\ref{SMXC22}, we know that~\eqref{m3edfgyui02woeiujfhnceu} holds true, with$$C_m:=\frac{\sup\limits_{x\in\R}|D^m f(x)|}{(2\pi)^m}.$$

Conversely, assume that~\eqref{m3edfgyui02woeiujfhnceu} holds true and let~$m_0\in\N$.
Taking~$m:=m_0+2$ in~\eqref{m3edfgyui02woeiujfhnceu}, we see that
$$ |\widehat f_k|\le\frac{C_{m_0+2}}{|k|^{m_0+2}}.$$
In particular, condition~\eqref{0909b} is fulfilled and therefore~$f\in C^{m_0}(\R)$.
Since~$m_0$ is arbitrary, we have that~$f\in C^\infty(\R)$.
\end{proof}

The duality between regularity of a function and decay of its Fourier coefficients carries over to the analytic case, namely:

\begin{theorem}\label{7un.mSIjkaja21dcA.13th}
Let~$f\in L^1((0,1))$ be periodic of period~$1$.

Then, $f$ is real analytic if and only if there exist~$C\in[0,+\infty)$ and~$\sigma\in(0,+\infty)$ such that, for every~$k\in\Z$,
\begin{equation}\label{AN2rfJSDAm3edfgyui02woeiujfhnceu}
|\widehat f_k|\le Ce^{-\sigma|k|}.
\end{equation}
\end{theorem}

\begin{proof} Assume that~$f$ is real analytic. Then (see e.g.~\cite[Corollary~1.2.8]{MR1916029}) for every~$p_0\in[0,1]$
there exist~$C_{p_0}$, $r_{p_0}\in(0,+\infty)$ such that, for every~$x\in[p_0-r_{p_0},p_0+r_{p_0}]$ and every~$m\in\N$,
\begin{equation} \label{MiGmr6defr-1}|D^m f(x)|\le\frac{C_{p_0}\;m!}{r_{p_0}^m}.\end{equation}
We now cover~$[0,1]$ with open intervals~$\left(p_0-\frac{r_{p_0}}2,p_0+\frac{r_{p_0}}2\right)$, varying~$p_0\in[0,1]$.
By compactness, we can extract a finite subcover and write that
\begin{equation}\label{MiGmr6defr-3} 
[0,1]\;\subseteq\;\bigcup_{\ell=1}^K \left(p_\ell-\frac{r_{p_\ell}}2,p_\ell+\frac{r_{p_\ell}}2\right),
\end{equation}
for suitable points~$p_\ell\in[0,1]$ and some~$K\in\N$.

We thereby define
\begin{equation}\label{MiGmr6defr-4} C_\star:=\max_{1\le\ell\le K}C_{p_\ell}\qquad{\mbox{and}}\qquad r_\star:=\min_{1\le\ell\le K} r_{p_\ell}
\end{equation}
and we claim that, for all~$x\in[0,1]$ and~$m\in\N$,
\begin{equation}\label{MiGmr6defr-2}
|D^m f(x)|\le\frac{C_\star\;m!}{r_\star^m}.
\end{equation}
To prove this, pick~$x\in[0,1]$. Then, by~\eqref{MiGmr6defr-3}, there exists~$\ell_x\in\{1,\dots,K\}$ for which
$$x\in
\left(p_{\ell_x}-\frac{r_{p_{\ell_x}}}2,p_{\ell_x}+\frac{r_{p_{\ell_x}}}2\right).$$

In particular, $x\in[p_{\ell_x}- r_{p_{\ell_x}},p_{\ell_x}+r_{p_{\ell_x}}]$ and therefore, by~\eqref{MiGmr6defr-1},
$$ |D^m f(x)|\le\frac{C_{p_{\ell_x}}\;m!}{r_{p_{\ell_x}}^m}.$$
From this and~\eqref{MiGmr6defr-4}, one obtains~\eqref{MiGmr6defr-2}, as desired.

As a consequence, in light of Theorem~\ref{SMXC22b}, for all~$k\in\Z$ and~$m\in\N$,
\begin{equation}\label{TBCAS}  |\widehat f_k|=\frac1{(2\pi| k|)^m}|\widehat{D^m f}_k|\le\frac1{(2\pi| k|)^m}\sup_{x\in[0,1]}|D^mf(x)|\le\frac{C_\star\;m!}{(2\pi r_\star |k|)^m}.
\end{equation}

It comes in handy now to observe that~$m$ here above is still for us to choose. Thus, we define
$$ \sigma:=\frac{\pi r_\star\,\ln 2}2\qquad{\mbox{and}}\qquad C:=\max\left\{C_\star, 2\max_{{j\in\Z}\atop{\pi r_\star \,|j|<2}}|\widehat f_j|\right\}.$$

In this way, we see that, if~$\pi r_\star \,|k|\ge2$ then the interval~$\left[\frac{\pi r_\star \,|k|}2,\pi r_\star \,|k|\right]$
has length at least~$1$ and thus it contains at least one integer~$m_k$.

Accordingly, choosing~$m:=m_k$ in~\eqref{TBCAS},
\begin{equation}\label{masmgfrmAqPtA}
|\widehat f_k|\le\frac{C_\star\;{m_k}!}{(2\pi r_\star |k|)^{m_k}}\le\frac{C_\star\;m_k^{m_k}}{(2\pi r_\star| k|)^{m_k}}
\le\frac{C_\star}{2^{m_k}}\le\frac{C_\star}{2^{\frac{\pi r_\star \,|k|}2}}
\le C e^{-\sigma|k|}
.\end{equation}

Moreover, if~$\pi r_\star \,|k|<2$,
\begin{eqnarray*}
|\widehat f_k|\le\max_{{j\in\Z}\atop{\pi r_\star \,|j|<2}}|\widehat f_j|\le \frac{C}2\le C e^{-\sigma|k|}.
\end{eqnarray*}
This observation and~\eqref{masmgfrmAqPtA} show that~\eqref{AN2rfJSDAm3edfgyui02woeiujfhnceu} holds true.

Suppose now that~\eqref{AN2rfJSDAm3edfgyui02woeiujfhnceu} is satisfied. Then, by Theorem~\ref{0CINFITH1} we know that~$f\in C^\infty(\R)$. Also, in light of~\eqref{CELLOD}, for all~$j\in\N$ and~$x\in\R$,
$$ D^jf(x)=\lim_{N\to+\infty}\sum_{{k\in\Z}\atop{|k|\le N}} (2\pi ik)^j \widehat f_k\,e^{2\pi ikx}.$$
This and~\eqref{AN2rfJSDAm3edfgyui02woeiujfhnceu} yield that
$$ |D^jf(x)|\le\sum_{k\in\Z} (2\pi |k|)^j |\widehat f_k|\le C\sum_{k\in\Z} (2\pi |k|)^j e^{-\sigma|k|}\le
2C\sum_{k=0}^{+\infty} (2k\pi )^j e^{-\sigma k}.$$
Hence (see Exercise~\ref{ui80-0})
$$ |D^jf(x)|\le \frac{L\, j!}{\rho^j},$$
for some~$L$, $\rho>0$.

This (see e.g.~\cite[Lemma~1.2.9]{MR1916029}) gives that~$f$ is real analytic.
\end{proof}

For general results about the decay of Fourier coefficients, see e.g.~\cite[Chapter~3]{MR3243734},
\cite{2018arXiv180502445N}, and the literature referred to there.

\begin{exercise}\label{ojqdwn23EEE0921ef}
Observe that the square and sawtooth waveforms considered in Exercises~\ref{SQ:W} and~\ref{SA:W},
are discontinuous functions whose $k$th Fourier coefficient are bounded in norm by~$\frac1{k}$.
\end{exercise}

\begin{exercisesk}\label{ojqdwn23EEE0921ef-LN}
Prove that there exists a function~$\phi\in L^1((0,1))$, periodic of period~$1$, such that
\begin{equation}\label{OSK-1-013-1}\widehat\phi_k=\begin{dcases}\frac{1}{2|k|}&{\mbox{ if }}k\ne0,\\0&{\mbox{ if }}k=0,\\
\end{dcases}\end{equation}
whose Fourier Series in trigonometric form has the expression
$$\sum_{k=1}^{+\infty}\frac1k\cos(2\pi kx)$$
and such that
\begin{equation}\label{GIJMPPrvU}\lim_{x\to0}\phi(x)=+\infty.\end{equation}
\end{exercisesk}

\begin{exercise} \label{ui80-0}
Let~$a\ge 0$ and~$\sigma>0$. Prove that there exist~$M$, $\rho>0$, depending only on~$a$ and~$\sigma$, such that, for every~$j\in\N$,
$$ \sum_{k=0}^{+\infty} (a k)^j e^{-\sigma k}\le \frac{M\, j!}{\rho^j}.$$
\end{exercise}

\begin{exercise}\label{EPPmdcer0}
Let~$E\subseteq(0,1)$ be a measurable set.
Let~$\{\xi_k\}_{k\in\N}$ be a sequence of real numbers such that
$$ \lim_{k\to+\infty}\xi_k=+\infty.$$
Prove that
\begin{equation}\label{EPPmdcer0.3ns} \lim_{k\to+\infty}\int_E \cos^2(2\pi kx+\xi_k)\,dx=\frac{|E|}2,\end{equation}
where~$|E|$ denotes the Lebesgue measure of~$E$. 
\end{exercise}

\begin{exercisesk}\label{EPPmdcer0.1}
Consider sequences of real numbers~$\{a_k\}_{k\in\N}$ and~$\{b_k\}_{k\in\N}$.
The \index{Cantor-Lebesgue Theorem} Cantor-Lebesgue Theorem states that if
\begin{equation}\label{EPPmdcer0.3} \lim_{k\to+\infty}a_k\cos(2\pi kx)+b_k\sin(2\pi kx)=0\end{equation}
for all~$x$ belonging to a measurable set~$E\subseteq(0,1)$ with positive measure, then
\begin{equation}\label{EPPmdcer0.4} \lim_{k\to+\infty}a_k=0\qquad{\mbox{and}}\qquad\lim_{k\to+\infty}b_k=0.
\end{equation}

Prove this result.
\end{exercisesk}

\begin{exercisesk}\label{EPPmdcer0.1-029.1}
Let~$f:\R\to\R$ be continuously differentiable and periodic of period~$1$ and, for all~$x$, $t\in\R$, define
\begin{equation}\label{EPPmdcer0.1-029.231.4.23.1}F(x,t):=\int_0^1 f(\xi)\,e^{i(\xi^2t-\xi x)}\,d\xi.\end{equation}
Let~$x(t)$ describe a curve such that
\begin{equation}\label{EPPmdcer0.1-029.231.74.2489} \lim_{t\to+\infty}\frac{x(t)}t=\rho,\end{equation}
for some~$\rho\in\R$.

Prove that, as~$t\to+\infty$,
\begin{equation}\label{EPPmdcer0.1-029.231.4.23}
F(x(t),t)=\frac{1+i}{2}\sqrt{\frac\pi{2t}}\,e^{-\frac{ix^2(t)}{2t}} f\left(\frac{\rho}{2}\right)+o\left(\frac1{\sqrt{t}}\right)
.\end{equation}
\end{exercisesk}

\begin{exercise}\label{WRO-eq00}
Learning from mistakes is important. Sometimes, in mathematics, a fruitful way to learn new concepts relies on understanding why an argument is {\em wrong}. Here we will provide an erroneous statement, with an incorrect proof, and the reader will be asked to spot what's wrong with it.

The fallacious claim is that
\begin{equation}\label{WRO-eq1}\begin{split}&
{\mbox{if a Fourier Series of a function~$f$ converges}}\\&{\mbox{and~$\widehat f_k=0$ for all~$k\in\Z\cap(-\infty,0)$,}}\\&{\mbox{then~$f$ is real analytic.}}\end{split}
\end{equation}

The invalid proof of~\eqref{WRO-eq1} goes as follows. Given~$x\in\R$, one uses the notation~$z:=e^{2\pi ix}$ and
\begin{equation}\label{WRO-eq20} g(z):=\sum_{k=0}^{+\infty} \widehat f_k\,z^k.\end{equation}
Then,
\begin{equation}\label{WRO-eq2}\begin{split}&
{\mbox{since the Fourier Series of the function~$f$ converges,}}\\&{\mbox{the
power series defining~$g$ converges and therefore~$g$ is analytic.}}
\end{split}\end{equation}
Also (see e.g.~\cite[Proposition~1.3.3]{MR1916029}), one has that
\begin{equation}\label{WRO-eq3}\begin{split}
{\mbox{the composition of analytic functions is analytic.}}
\end{split}\end{equation}
Since the complex exponential map~$x\mapsto e^{2\pi ix}=z$ is analytic, the claim in~\eqref{WRO-eq1} follows from~\eqref{WRO-eq2} and~\eqref{WRO-eq3}!

What's wrong\footnote{This exercise is indeed technically demanding
(though perhaps not overly difficult)
but conceptually quite subtle, and, historically, it
highlights one of the main contributions of Fourier to the theory that is
nowadays named after him, see~\cite[Chapter~I]{MR605488}.

Indeed, it appears that infinite series of sines and cosines were already adopted by
Daniel Bernoulli and retaken by Euler, who seemed to believe, however,
that ``such a series of analytic functions could
only represent an analytic function.
 His opinion was shared
 (with some variations) by almost all other mathematicians of
his time, and no progress was made on this question until the
beginning of Fourier's work. [...]
Fourier was able to show on many
examples of non analytic functions that the corresponding
Fourier series converged to~$\frac12(f(x_+) + f(x_-))$, and 
expressed his conviction that this was true for {\em arbitrary} 
functions...''. 

It is a pleasure to thank Xavier Cabr\'e for narrating this historical anecdote and for clarifying the importance of Fourier's pioneering idea.} with this?
\end{exercise}

\section{Uniform convergence results}\label{UNIFORMCO:SECTI}

Now, we refine the results in Section~\ref{CONVES} in order to obtain uniform convergence
of the Fourier Series to the original function under suitable additional assumptions.
These assumptions can be taken either on the Fourier coefficients (a prototypical result in this context
being that presented in Theorem~\ref{BASw}) or on the regularity of the function (as exemplified in Theorem~\ref{C1uni}).

Concerning the first of these two perspectives, we remark that the absolute convergence of the series of the Fourier coefficients is sufficient to guarantee the uniform convergence of the Fourier Series to the original function:

\begin{theorem}\label{BASw}
Let~$f:\R\to\R$ be periodic of period~$1$, with~$f\in L^1((0,1))$.

Assume that
\begin{equation}\label{AKSMaLS} \sum_{k\in\Z}|\widehat f_k|<+\infty.\end{equation}

Then, the Fourier Series of~$f$ converges uniformly\footnote{More specifically, the Fourier Series converges uniformly in $\R$ to some continuous function $g$, which is equal to $f$ up to a set of null Lebesgue measure. In cases like this, it is customary to identify $f$ with its continuous representative $g$.} to~$f$.

Moreover, $f$ is necessarily a continuous function.
\end{theorem}

\begin{proof}
The proof is a variation of that of Theorem \ref{ojqdwn23E} and the details are as follows. By Lemma \ref{Le-ojqdwn23E} used here with $m:=0$ and $c_k:=\widehat{f}_k$, we know that $S_{N,f}$ converges uniformly to some continuous function~$g$ with $\widehat{g}_k=c_k=\widehat{f}_k$.

This and Theorem~\ref{UNIQ} yield that~$f$ and~$g$ coincide, as desired.
\end{proof}

For other results ensuring uniform convergence, see e.g.~\cite[Theorem~15.3]{MR4404761}, \cite[Chapter~II]{MR1963498}, and~\cite[Section~6, Chapter~I]{MR2039503}.

As an immediate consequence of Theorems~\ref{SMXC22} and~\ref{BASw}, one obtains that if~$f\in C^2(\R)$
and periodic of period~$1$, then its Fourier Series converges uniformly to~$f$.
But one can do better than this, both in terms of the regularity required on the function and on the localisation of such an assumption:
as a paradigmatic example of uniform convergence of the Fourier Series due to a regularity assumption on the original function, we present the following result.

\begin{theorem}\label{C1uni}
Let~$f\in L^1((0,1))$ be periodic of period~$1$.

Let~$I$ be a bounded interval in~$\R$. Assume that for every~$\epsilon>0$
there exists~$\delta_\epsilon>0$ such that
\begin{equation}\label{DINI-UNO}
\sup_{x\in I}\int_{-\delta_\epsilon}^{\delta_\epsilon}\left|\frac{f (x + t) - f(x)}{t}\right|\,dt\le \epsilon.\end{equation}

Then, the Fourier Series of~$f$ converges to~$f$ uniformly in~$I$.
\end{theorem}

Comparing with~\eqref{DINI}, one can consider~\eqref{DINI-UNO} a sort \index{Dini's Condition, uniform}
of \emph{uniform Dini's Condition}.

To prove Theorem~\ref{C1uni} we split the argument into two separate parts, plus an addendum:\begin{itemize}
\item first, we will revisit formula~\eqref{KASMqwdfed123er} by translating it at every given point
(and, to make justice of the ``aesthetic'' principles discussed on page~\pageref{QSTNB},
we will obtain an expression whose right-hand side depends more symmetrically on~$N$): obtaining the counterpart
of formula~\eqref{KASMqwdfed123er} centred at any point is technically advantageous, because while in
Sections~\ref{CONVESob} and~\ref{CONVES} we were dealing with convergence at a point (which we could conveniently reduce to be the origin) and therefore we only cared about formula~\eqref{KASMqwdfed123er} at a specific point (actually, the origin),
here we will need to deal with uniform convergence and accordingly we will have to express the sum of finitely many complex exponentials in a given interval, not just at a point, and deduce suitable uniform estimates on that expression;
for this, we will need an expression for~$f-S_N$ at every point,
\item then, we complete the proof of Theorem~\ref{C1uni}, \item as an addendum, in Theorem~\ref{C1uni-c}, we will stress (quite surprisingly, given the ``global'' character of the Fourier coefficients) that
the type of local convergence of a Fourier Series
depends only on the local behaviour of the function.\end{itemize}

Regarding the ``polished'' version of formula~\eqref{KASMqwdfed123er}, we have:

\begin{lemma}\label{KASMqwdfed123erDKLI} We have that
\begin{equation} \label{PAKSw-L4}
\sum_{{k\in\Z}\atop{|k|\le N}} e^{-2\pi i kx}=
\frac{\sin\big((2N+1)\pi x\big)}{\sin(\pi x)}=:D_N(x),
\end{equation}
with continuous extension of this identity when the denominator vanishes.
\end{lemma}

Formula~\eqref{PAKSw-L4} is often called the \emph{Lagrange Trigonometric Identity} \index{Lagrange Trigonometric Identity}
and $D_N$ is called \index{Dirichlet Kernel}
the \emph{Dirichlet Kernel}.

\begin{figure}[h]
\includegraphics[height=4.3cm]{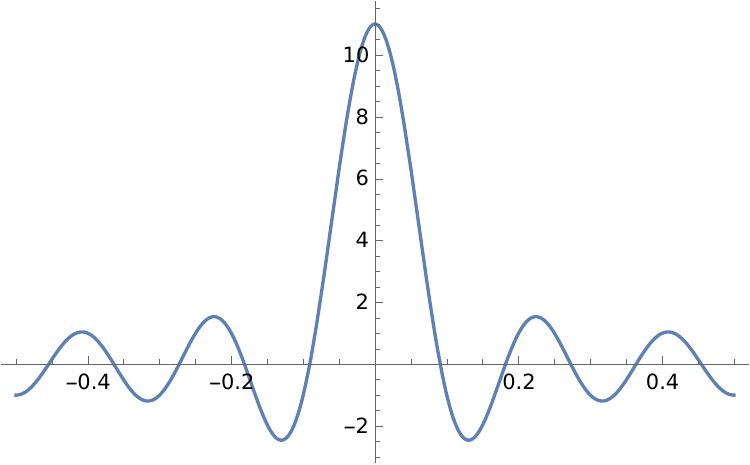}$\quad$\includegraphics[height=4.3cm]{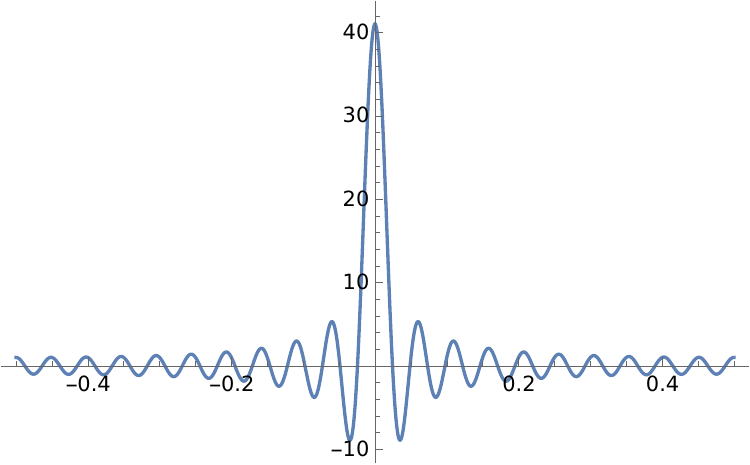}
\centering
\caption{The Dirichlet Kernel (left, with~$N=5$; right, with~$N=20$).}\label{amdi}
\end{figure}

See Figure~\ref{amdi} for a diagram showcasing
the ``concentration'' property of the Dirichlet Kernel. See also~\cite{MR4509245, MR4596335} for an ingenious ``proof without words'' of the Lagrange Trigonometric Identity.

In our context, the utility of Lemma~\ref{KASMqwdfed123erDKLI} is to allow for an explicit expression, valid at every point, that describes the discrepancy between a given function and its Fourier Sum, which goes as follows:

\begin{corollary}\label{9.0km2df0.2wnedfo.1-2x}
Let~$f\in L^1((0,1))$ be periodic of period~$1$.

Then, for all~$N\in\N$,
\begin{equation}\label{023wefv4567uygfdfgyhuizo0iwhfg0eb5627e01.pre}
S_N(x)=\int_0^1 f(y)\,D_N(x-y)\,dy
\end{equation}
and\footnote{Using the notation of the periodic convolution introduced in Exercise~\ref{pecoOSJHN34jfgj}, one can also
rewrite~\eqref{023wefv4567uygfdfgyhuizo0iwhfg0eb5627e01.pre} in the form
$$ S_N=f\star D_N.$$}
\begin{equation}\label{023wefv4567uygfdfgyhuizo0iwhfg0eb5627e01} f(x)-S_N(x)=
\int_0^1 \big(f(x)-f(y)\big)\,\frac{\sin\big((2N+1)\pi (y-x)\big)}{\sin(\pi(y- x))}\,dy.
\end{equation}
\end{corollary}

Now we prove Lemma~\ref{KASMqwdfed123erDKLI}, which will promptly imply Corollary~\ref{9.0km2df0.2wnedfo.1-2x}.

\begin{proof}[Proof of Lemma~\ref{KASMqwdfed123erDKLI}] Retaking~\eqref{KASMqwdfed123er},
\begin{eqnarray*}&&
\sum_{{k\in\Z}\atop{|k|\le N}} e^{-2\pi i kx}= \frac{e^{-\pi ix} \Big(e^{2\pi i\left(N+\frac12\right)x}-e^{-2\pi i\left(N+\frac12\right)x}\Big)}{1-e^{-2\pi ix}}\\&&\qquad=\frac{2ie^{-\pi ix} \sin\left(2\pi \left(N+\frac12\right)x\right)}{1-e^{-2\pi ix}}
=\frac{2i\sin\left(2\pi \left(N+\frac12\right)x\right)}{e^{\pi ix}-e^{-\pi ix}}=\frac{\sin\left( \left(2N+1\right)\pi x\right)}{
\sin(\pi x)}.\qedhere
\end{eqnarray*}
\end{proof}

\begin{proof}[Proof of Corollary~\ref{9.0km2df0.2wnedfo.1-2x}]
In light of Lemma~\ref{KASMqwdfed123erDKLI},
\begin{eqnarray*}S_N(x)&=&\sum_{{k\in\Z}\atop{|k|\le N}}\int_0^1 f(y)\,e^{-2\pi ik(y-x)}\,dy
\\&=&\int_0^1 f(y)\,\sum_{{k\in\Z}\atop{|k|\le N}}e^{-2\pi ik(y-x)}\,dy
\\&=&\int_0^1 f(y)\,\frac{\sin\big((2N+1)\pi (y-x)\big)}{\sin(\pi(y- x))}\,dy.
\end{eqnarray*}

This gives~\eqref{023wefv4567uygfdfgyhuizo0iwhfg0eb5627e01.pre} and also that
\begin{eqnarray*} f(x)-S_N(x)&=&f(x)-\int_0^1 f(y)\,\frac{\sin\big((2N+1)\pi (y-x)\big)}{\sin(\pi(y- x))}\,dy.
\end{eqnarray*}
Since (see Exercises~\ref{fr12} and~\ref{K-1PIO})
\begin{eqnarray*}&& \int_0^1 \frac{\sin\big((2N+1)\pi (y-x)\big)}{\sin(\pi(y- x))}\,dy=
\int_{-x}^{1-x} \frac{\sin\big((2N+1)\pi t\big)}{\sin(\pi t)}\,dt\\&&\qquad=
\int_{0}^{1} \frac{\sin\big((2N+1)\pi t\big)}{\sin(\pi t)}\,dt=
\int_0^1 D_N(t)\,dt
=1,\end{eqnarray*}
we obtain~\eqref{023wefv4567uygfdfgyhuizo0iwhfg0eb5627e01}.
\end{proof}

Now we complete the proof of Theorem~\ref{C1uni}:

\begin{proof}[Proof of Theorem~\ref{C1uni}] First of all, we can suppose that
\begin{equation}\label{men34}
{\mbox{the length of the interval~$I$ is less than or equal to~$\frac34$.}}
\end{equation}
Indeed, if~$I$ has length larger than~$\frac34$, we can write~$I=I_1\cup\dots\cup I_L$,
where, for all~$j\in\{1,\dots,L\}$, $I_j$ is an interval with length less than or equal to~$\frac34$.
Since~\eqref{DINI-UNO} holds true with~$I$ replaced by~$I_j$,
we could then apply Theorem~\ref{C1uni} to~$I_j$, obtaining that for any~$\eta>0$ there exists~$N_{\eta,j}\in\N$ such that for all~$N\ge N_{\eta,j}$ one has that
$$ \sup_{x\in I_j}|f(x)-S_{N}(x)|\le\eta.$$
Hence, we define
$$N_\eta:=\max_{j\in\{1,\dots,L\}}N_{\eta,j}$$
and conclude that, when~$N\ge N_\eta$,
$$ \sup_{x\in I}|f(x)-S_{N}(x)|\le\max_{j\in\{1,\dots,L\}}
\sup_{x\in I_j}|f(x)-S_{N}(x)|\le\eta$$
and therefore~$S_N$ converges to~$f$ uniformly in~$I$.

This shows that we can assume~\eqref{men34} for the rest of this proof.

Additionally, up to a translation, one can take the centre of~$I$ to be at~$\frac12$.
From this and~\eqref{men34}, it follows that
\begin{equation}\label{men345}
I\subseteq\left[\frac18,\frac78\right].
\end{equation}

Now, given~$\epsilon>0$, we deduce from condition~\eqref{DINI-UNO} (where we can take~$\delta_\epsilon\in\left(0,\frac1{10}\right)$ without loss of generality) that, for every~$x\in I$,
\begin{eqnarray*}
&&\int_{(0,1)\cap(x-\delta_\epsilon,x+\delta_\epsilon) } \big|f(x)-f(y)\big|\,\frac{\big|\sin\big((2N+1)\pi (y-x)\big)\big|}{\big|\sin(\pi(y- x))\big|}\,dy\\
&&\quad\le\sup_{\theta\in(-\delta_\epsilon,\delta_\epsilon)}\left|\frac{\theta}{\sin(\pi\theta)}\right|
\int_{(0,1)\cap(x-\delta_\epsilon,x+\delta_\epsilon) } \left|\frac{f(y)-f(x)}{y-x}\right|\,
\big|\sin\big((2N+1)\pi (y-x)\big)\big|\,dy\\
&&\quad\le\sup_{\theta\in(-1/10,1/10)}\left|\frac{\theta}{\sin(\pi\theta)}\right|
\int_{-\delta_\epsilon}^{\delta_\epsilon} \left|\frac{f(x+t)-f(x)}{t}\right|\,
\,dt\\&&\quad\le C\epsilon,
\end{eqnarray*}
for some~$C>0$, and therefore we infer from~\eqref{023wefv4567uygfdfgyhuizo0iwhfg0eb5627e01} that,
for every~$x\in I$,
\begin{equation}\label{OPK-oDer58j-09ij8765tyhhji9205132021f41}\begin{split}&| f(x)-S_N(x)|\\
&\qquad\le C\epsilon+\left|
\int_{(0,1)\setminus(x-\delta_\epsilon,x+\delta_\epsilon) } \big(f(x)-f(y)\big)\,\frac{\sin\big((2N+1)\pi (y-x)\big)}{\sin(\pi(y- x))}\,dy\right|.
\end{split}\end{equation}

Now we observe that, in light of~\eqref{men345},
if~$y\in(0,1)\setminus(x-\delta_\epsilon,x+\delta_\epsilon)$ then~$|\pi(y-x)|\ge\pi\delta_\epsilon$ and
$$ |\pi(y-x)|=\pi\max\big\{y-x,x-y\big\}\le\pi\max\left\{1-\frac18,\frac78-0\right\}=\frac{7\pi}8\le\pi(1-\delta_\epsilon).$$
As a consequence, if~$y\in(0,1)\setminus(x-\delta_\epsilon,x+\delta_\epsilon)$
\begin{equation*} |\sin(\pi(y- x))|\ge\sin(\pi\delta_\epsilon).\end{equation*}
Hence, we take~$\varphi\in C^\infty_0((-1,1),[0,1])$ with unit mass and, for~$\eta>0$ to be taken conveniently small, we define
$$ f_\eta(x):=\frac1\eta\int_{-\infty}^{+\infty} f(x-y)\varphi\left(\frac{y}\eta\right)\,dy.$$
We know (see e.g.~\cite[Theorem~9.6]{MR3381284}) that if~$\eta=\eta_\epsilon$ is sufficiently small then
$$\|f-f_{\eta_\epsilon}\|_{L^1((0,1))}\le {\epsilon}\,{\sin(\pi\delta_\epsilon)}.$$
Moreover, for all~$x\in I$,
\begin{eqnarray*}&&
|f(x)-f_{\eta_\epsilon}(x)|\le 
\frac1{\eta_\epsilon}\int_{-\eta_\epsilon}^{\eta_\epsilon} |f(x)-f(x-y)|\,dy\\&&\qquad\qquad\le
\int_{-\eta_\epsilon}^{\eta_\epsilon} \left|\frac{f(x)-f(x-y)}{y}\right|\,dy\le {\epsilon}\,{\sin(\pi\delta_\epsilon)},
\end{eqnarray*}
thanks to condition~\eqref{DINI-UNO}, as long as~$\eta_\epsilon$ is sufficiently small.

As a result, for all~$x\in I$,
\begin{eqnarray*}&&
\int_{(0,1)\setminus(x-\delta_\epsilon,x+\delta_\epsilon) } \big|f(x)-f_{\eta_\epsilon}(x)-f(y)+f_{\eta_\epsilon}(y)\big|\,\frac{\big|\sin\big((2N+1)\pi (y-x)\big)\big|}{\big|\sin(\pi(y- x))\big|}\,dy
\\&&\quad\le\frac{1}{\sin(\pi\delta_\epsilon)}\left(
\big|f(x)-f_{\eta_\epsilon}(x)\big|+
\int_0^1\big|f(y)+f_{\eta_\epsilon}(y)\big|\,dy
\right)\\&&\quad\le2\epsilon.
\end{eqnarray*}
Combining this information and~\eqref{OPK-oDer58j-09ij8765tyhhji9205132021f41} we conclude that,
for every~$x\in I$,
\begin{equation}\label{087bgdvit65Rynsg1ilonjmtian}\begin{split}&| f(x)-S_N(x)|\\&\qquad\le (C+2)\epsilon+\left|
\int_{(0,1)\setminus(x-\delta_\epsilon,x+\delta_\epsilon) } \big(f_{\eta_\epsilon}(x)-f_{\eta_\epsilon}(y)\big)\,\frac{\sin\big((2N+1)\pi (y-x)\big)}{\sin(\pi(y- x))}\,dy\right|\\&\qquad= (C+2)\epsilon+\left|
\int_{(0,1)\setminus(x-\delta_\epsilon,x+\delta_\epsilon)} g_\epsilon(x,y)\,\sin\big((2N+1)\pi (y-x)\big)\,dy\right|,\end{split}
\end{equation}
where
$$ g_\epsilon(x,y):=\frac{\big(f_{\eta_\epsilon}(x)-f_{\eta_\epsilon}(y)\big)}{\sin(\pi(y- x))}.$$

Furthermore, integrating by parts,
\begin{eqnarray*}&&
\left|\int_{(0,1)\setminus(x-\delta_\epsilon,x+\delta_\epsilon)} g_\epsilon(x,y)\,\sin\big((2N+1)\pi (y-x)\big)\,dy\right|\\&&\qquad=
\frac1{(2N+1)\pi}\left|\int_{(0,1)\setminus(x-\delta_\epsilon,x+\delta_\epsilon)} g_\epsilon(x,y)\,\frac{d}{dy}\Big(\cos\big((2N+1)\pi (y-x)\big)\Big)\,dy\right|\\&&\qquad\le\frac1{(2N+1)\pi}\left(\big| g_\epsilon(x,x+\delta_\epsilon)-g_\epsilon(x,x-\delta_\epsilon)\big|
+\int_0^1 |\partial_y g_\epsilon(x,y)|\,dy
\right)\\&&\qquad\le\frac{C_\epsilon}{2N+1},
\end{eqnarray*}
for some~$C_\epsilon>0$.

On this account and~\eqref{087bgdvit65Rynsg1ilonjmtian}, we have that
$$ | f(x)-S_N(x)|\le (C+2)\epsilon+\frac{C_\epsilon}{2N+1}$$
and accordingly
$$ \lim_{N\to+\infty}\sup_{x\in I}| f(x)-S_N(x)|\le \lim_{N\to+\infty}\left((C+2)\epsilon+\frac{C_\epsilon}{2N+1}\right)=(C+2)\epsilon.$$
Hence, sending now~$\epsilon\searrow0$ we conclude that
$$ \lim_{N\to+\infty}\sup_{x\in I}| f(x)-S_N(x)|=0,$$
as desired.
\end{proof}

We now return to the ``local'' nature of the convergence results for Fourier Series, already mentioned on page~\pageref{LAMS}, and we present a result known under the name of \index{Riemann Localisation Principle}
\emph{Riemann Localisation Principle}, stating that
if two functions agree in a neighbourhood of a given point, then their Fourier Series behave in the same
way at that point:

\begin{theorem}\label{C1uni-c}
Let~$f$, $g\in L^1((0,1))$ be periodic of period~$1$.

Let~$I$ be an interval and assume that~$f$ and~$g$ coincide a.e. in~$I$.

Let~$J$ be an interval such that the closure of~$J$ is contained in the interior of~$I$.

Then,
\begin{equation}\label{eeAJXMSCX} \lim_{N\to+\infty} S_{N,f}-S_{N,g}=0\end{equation}
uniformly in~$J$.

In particular, for every~$x$ in the interior of~$ I$, the Fourier Series of~$f$ at~$x$ converges if and only if so does that of~$g$,
in which case the two Fourier Series converge at the same value.
\end{theorem}

We observe that the domain of uniform convergence of Theorem~\ref{C1uni-c} is optimal, see Exercise~\ref{C1uni-c-EXEOP}.

\begin{proof}[Proof of Theorem~\ref{C1uni-c}] Up to replacing~$f$ with~$f-g$,
without loss of generality we can suppose that~$g=0$. In this setting, the desired claim in~\eqref{eeAJXMSCX} becomes
\begin{equation}\label{eeAJXMSCX0-n} \lim_{N\to+\infty} S_{N,f}=0\end{equation}
uniformly in~$J$, as long as~$f$ vanishes a.e. in~$I$.

Hence, we let~$\delta_\star>0$ be the minimal distance between points in~$J$ and~$\partial I$:
in this way, if~$x\in J$, then~$(x-\delta_\star,x+\delta_\star)\in I$ and consequently
\begin{equation*}
\sup_{x\in I}\int_{-\delta_\star}^{\delta_\star}\left|\frac{f (x + t) - f(x)}{t}\right|\,dt=
\sup_{x\in I}\int_{-\delta_\star}^{\delta_\star}\left|\frac{0 -0}{t}\right|\,dt=0.\end{equation*}
This gives that condition~\eqref{DINI-UNO} is satisfied and therefore the claim in~\eqref{eeAJXMSCX0-n}
follows from Theorem~\ref{C1uni}.
\end{proof}

The topic of pointwise and uniform convergence of Fourier Series is deep and rich of results.
See e.g.~\cite[Chapter~15]{MR4404761} for additional information.

\begin{exercise}\label{NEWLA.pre}
Let~$M\in\N$ and~$a_0,\dots,a_M\in[0,1]$ with~$0=a_0<a_1<\dots<a_{M-1}<a_M=1$.

Let~$f:\R\to\R$ be continuous and periodic of period~$1$.

Assume that~$f$ is differentiable in~$(a_{j-1},a_j)$ and~$f'\in L^2((a_{j-1},a_j))$ for each~$j\in\{1,\dots,M\}$.

Prove that the Fourier Series of~$f$ converges to~$f$ uniformly in~$\R$.
\end{exercise}

\begin{exercise}\label{NEWLA}
Prove that the Lagrange Trigonometric Identity~\eqref{PAKSw-L4}
can be rewritten in the form
$$\sum_{k=0}^N \cos(2\pi kx)=\frac12+
\frac{\sin\big((2N+1)\pi x\big)}{2\sin(\pi x)}.$$\end{exercise}

\begin{exercise}\label{LGEGMCTLFDCD5D325E28A9N}
Prove that
$$ \max_{x\in\R}\big|D_N(x)\big|=2N+1=D_N(0).$$\end{exercise}

\begin{exercise}\label{K-0PIO}
Prove that $$
D_N\left(\frac12-x\right)-D_N\left(\frac12+x\right)=0.$$\end{exercise}

\begin{exercise}\label{K-1PIO}
Calculate
$$ \int_0^1 D_N(x)\,dx.$$
\end{exercise}

\begin{exercise}\label{K-1PIO-bis}
Calculate
$$ \int_0^{1/2} D_N(x)\,dx.$$
\end{exercise}

\begin{exercise}\label{K-2PIO}
Calculate
$$ \int_0^1 \big(D_N(x)\big)^2\,dx.$$
\end{exercise}

\begin{exercise}\label{K-3PIO}
Prove that\footnote{The quantities in the left-hand side of~\eqref{makm.swela0oerf.ZZ}
are sometimes called \index{Lebesgue's numbers} \emph{Lebesgue's numbers}.} there exists~$c>0$ such that, for all~$N\in\N$,
\begin{equation}\label{makm.swela0oerf.ZZ} \int_{-1/2}^{1/2} |D_N(x)|\,dx\ge c\ln N.\end{equation}
\end{exercise}

\begin{exercise}\label{K-3PIO.eDhnZmasKvcTYhFA.1} Is the estimate in Exercise~\ref{K-3PIO} optimal?
\end{exercise}

\begin{exercise}\label{IFEJ} Let
$$ F_{N}(x):=\frac{1}{N}\sum_{k=0}^{N-1}D_{k}(x).$$
This is called the \index{Fej\'er Kernel}
\emph{Fej\'er Kernel} and will play a major role\footnote{It seems that Fej\'er was aged only 19 when he introduced his kernel and made a number of wonderful discoveries about Fourier Series: however (see~\cite[page~5]{MR4404761})
``any reader discouraged by Fej\'er's precocity should note that a few years earlier his school considered him so weak in mathematics as to require extra tuition''.}
in Section~\ref{FEJERKESE}.

Prove that
$$ F_N(x)=\frac1{N}\left( \frac{\sin(N\pi x)}{\sin(\pi x)}\right)^2.$$
\end{exercise}

\begin{exercise}\label{FEJPO}
Prove that~$F_N(x)\ge0$ for all~$x\in\R$.
\end{exercise}

\begin{exercise}\label{VjweWRFVSHEgKFB8ndVIAPMAQ}
Let~$a\in\left(0,\frac12\right)$ and~$\phi_a$ be the function defined in
Exercise~\ref{FO:DE:MA}, periodically extended as a function of period~$1$
outside~$\left[-\frac12,\frac12\right)$.

Let~$g_a:=\phi_a-1$.

Let also~$\{p_k\}_{k\in\N}$ be a sequence of non-negative real numbers such that
$$ \sum_{k=0}^{+\infty}p_k<+\infty.$$

Define
$$ f(x):=\sum_{k=3}^{+\infty} p_k\,g_{1/k}(x).$$

Prove that~$f\in L^1((0,1))$, $f$ is even and periodic of period~$1$, and
\begin{equation*} \int_{0}^{1}f(x)\,dx=0.\end{equation*}
\end{exercise}

\begin{exercise}\label{VjweWRFVSHEgKFB8ndVIAPMAQ-bis}
Write the Fourier Series in trigonometric form of the function~$f$ in 
Exercise~\ref{VjweWRFVSHEgKFB8ndVIAPMAQ}.
\end{exercise}

\begin{exercise}\label{K-4PIO}
Use the Dirichlet Kernel in~\eqref{PAKSw-L4} to prove that\footnote{The integrand in~\eqref{CASIEQIN} is a very popular function in Fourier analysis and it is sometimes called the cardinal sine. \index{cardinal sine}
For instance, it will appear again in Section~\ref{SEC:GIBBS-PH}.}
\begin{equation}\label{CASIEQIN}\int_0^{+\infty}\frac{\sin x}{x}\,dx=\frac\pi2.\end{equation}
\end{exercise}

\begin{exercise}\label{CO213421FSBCA.02-ejNS12c2g2n4}
Prove that
\begin{equation} \label{CO213421FSBCA.02-ejNS12c2g2n4.E010}
\sup_{m\in\N\setminus\{0\}}\int_0^{\pi m}\frac{\sin x}{x}\,dx=\int_0^{\pi}\frac{\sin x}{x}\,dx>0\end{equation}
and that
\begin{equation} \label{CO213421FSBCA.02-ejNS12c2g2n4.E011}\sup_{m\in\N\setminus\{0\}}\int_0^{2\pi m}\frac{\sin x}{x}\,dx=\frac\pi2.\end{equation}
\end{exercise}

\begin{exercise}\label{DINI1UNI} Show that if $f$ is globally H\"older continuous, \index{H\"older continuous}
i.e. if it is bounded and there exist~$\alpha\in(0,1)$
and~$M\ge0$ such that, for all~$x$, $y\in\R$,
$$ |f (x)-f(y)|\le M|x-y|^\alpha,$$ then~$f$ satisfies
the uniform version of Dini's Condition~\eqref{DINI-UNO} for any bounded interval~$I$ (and, in fact, even replacing~$I$ with the whole of~$\R$).
\end{exercise}

\begin{exercise}\label{DINI2UNI} Show that if $f$ is globally \index{Lipschitz continuous}
Lipschitz continuous, i.e. if it is bounded and there exists~$M\ge0$ such that, for all~$x$, $y\in\R$,
$$ |f (x)-f(y)|\le M|x-y|,$$ then~$f$ is globally H\"older continuous.
\end{exercise}

\begin{exercise}\label{DINI3UNI} Show that if $f$ is bounded and differentiable, and~$f'$ is bounded, then~$f$ is globally Lipschitz continuous.\end{exercise}

\begin{exercise}\label{DINI-addC-2-maunif}
Let~$f:\R\to\R$ be periodic of period~$1$. Assume that~$f\in L^1((0,1))$ and that it satisfies the uniform version of Dini's Condition~\eqref{DINI-UNO} in some interval~$I$. 

Does this mean that~$f$ is necessarily continuous\footnote{It is interesting to compare Exercise~\ref{DINI-addC-2-maunif} with Exercise~\ref{DINI-addC-2}.} at every point of~$I$?\end{exercise}

\begin{exercise}\label{C1uni-c-EXEOP}
In the setting of Theorem~\ref{C1uni-c}, show that it is not possible to replace the sentence ``uniformly in~$J$'' with ``uniformly in~$I$''.

That is, prove that there exist functions~$f$, $g\in L^1((0,1))$ that are periodic of period~$1$ and coincide on an interval~$I$ but are such that
$$ \liminf_{N\to+\infty} \sup_{x\in I}|S_{N,f}(x)-S_{N,g}(x)|>0.$$
\end{exercise}

\begin{exercise}\label{7un.mSIjkaja21dcA.13}
Prove that if a function~$f:\R\to\R$ is real analytic and periodic,
then its Fourier Series converges exponentially fast, that is there exist~$C$, $\sigma>0$ such that
$$ \sup_{x\in\R}|f(x)-S_N(x)|\le Ce^{-\sigma N}.$$
\end{exercise}

\begin{exercise}\label{7un.mSIjkaja21dcA.13BIS}
The Fourier Series of smooth functions converge polynomially fast.
More specifically: let~$m\in\N$, with~$m\ge2$, and~$f\in C^m(\R)$ be periodic of period~$1$.

Prove that there exists~$C>0$ such that
$$ \sup_{x\in\R}|f(x)-S_N(x)|\le \frac{C}{N^{m-1}}.$$
\end{exercise}

\begin{exercise}\label{PIPIO-010}
Prove that, for all~$m$, $k\in\Z$ and~$x\in\R$,
$$\left|\sum_{{j\in\Z}\atop{m\le j\le k-1}}e^{2\pi ijx}\right|\le\frac1{|\sin(\pi x)|}.$$\end{exercise}

\begin{exercise}\label{PIPIO-0}
Consider a sequence~$\{\gamma_k\}_{k\in\N}$ such that~$\gamma_k\ge \gamma_{k+1}\ge0$.

Prove that, for all~$m$, $n\in\Z$ and~$x\in\R$,
$$\left|\sum_{{k\in\Z}\atop{m\le k\le n}}\gamma_k e^{2\pi ikx}\right|\le \frac{\gamma_m}{|\sin(\pi x)|}.$$\end{exercise}

\begin{exercise}\label{PIPIOmaq}
Consider a sequence~$\{\gamma_k\}_{k\in\N}$ such that~$\gamma_k\ge \gamma_{k+1}\ge0$ and~$\gamma_k\searrow0$ as~$k\to+\infty$.

Prove that, for every~$\epsilon\in\left(0,\frac12\right)$, the series
$$\sum_{k=1}^{+\infty}\gamma_k\,e^{2\pi ikx}$$
converges uniformly for~$x\in(\epsilon,1-\epsilon)$.
\end{exercise}

\begin{exercisesk}\label{PIPIO}
Consider a sequence~$\{\gamma_k\}_{k\in\N}$ such that~$\gamma_k\ge \gamma_{k+1}\ge0$ and~$k\gamma_k\le1$.

Prove that, for all~$N\in\N$ and~$x\in\R$,
$$\left|\sum_{k=1}^N\gamma_k\sin(2\pi kx)\right|\leq\pi+1.$$
\end{exercisesk}

\begin{exercise}\label{RENOLLME}
Does the result in Exercise~\ref{PIPIO} remain valid if one replaces the sine function with the cosine function?
\end{exercise}

\begin{exercise}\label{PIPIO-00-123-321}
Consider a sequence~$\{\gamma_k\}_{k\in\N}$ such that~$\gamma_k\ge \gamma_{k+1}\ge0$ and~$k\gamma_k\searrow0$ as~$k\to+\infty$.

Prove that
$$ \lim_{\delta\to0}\sum_{k=0}^{+\infty} \gamma_k\,\sin(2\pi k\delta)=0.$$
\end{exercise}

\begin{exercisesk}\label{ESIPNo}
Does there exist~$f\in L^1((0,1))$, periodic of period~$1$, whose Fourier Series in trigonometric form is
\begin{equation}\label{muFVPFbaGtpf} \sum_{k=2}^{+\infty} \frac{\sin(2\pi kx)}{\ln k}\; ?\end{equation}
\end{exercisesk}
`
\begin{exercisesk}\label{ESIPNo-2}
Does there exist~$f\in L^1((0,1))$, periodic of period~$1$, whose Fourier Series in trigonometric form is
\begin{equation*} \sum_{k=2}^{+\infty} \frac{\cos(2\pi kx)}{\ln k}\; ?\end{equation*}
\end{exercisesk}

\begin{exercisesk}\label{LESENUEDGE}
Let~$\{\gamma_k\}_{k\in\N}$ be a sequence of real numbers such that~$\gamma_k\ge\gamma_{k+1}\ge0$
and~$\gamma_k\searrow0$ as~$k\to+\infty$.

Suppose in addition that~$\gamma_k$ is a convex sequence, namely that
$$ \gamma_{k+1}+\gamma_{k-1}-2\gamma_k\ge0.$$

Under these assumptions, does there exist~$f\in L^1((0,1))$ whose Fourier Series in trigonometric form is
$$\sum_{k=0}^{+\infty} \gamma_k\,\cos(2\pi kx)\;?$$\end{exercisesk}

\begin{exercise}\label{ESIPNo-3}
Does there exist~$f\in L^1((0,1))$, periodic of period~$1$, whose Fourier Series in trigonometric form is
\begin{equation} \label{ITFS8901P13cde}\sum_{k=2}^{+\infty} \frac{\sin(2\pi kx)}{k\,\ln k}\; ?\end{equation}
\end{exercise}

\begin{exercise}\label{ESIPNo-2-ap2}
Does there exist~$f\in L^2((0,1))$, periodic of period~$1$, whose Fourier Series in trigonometric form is
\begin{equation*} \sum_{k=2}^{+\infty} \frac{\cos(2\pi kx)}{\ln k}\; ?\end{equation*}
\end{exercise}

\begin{exercise}\label{FTCPASNo1}
Prove that there exists a positive, monotone, and convex sequence which tends to zero arbitrarily slowly.

More specifically, let~$\{\sigma_k\}_{k\in\N}$ be a sequence of real numbers such that~$\sigma_k\to0$ as~$k\to+\infty$.

Prove that there exists a sequence~$\{\gamma_k\}_{k\in\N}$ of real numbers
such that~$\gamma_k\ge\gamma_{k+1}\ge0$, $\gamma_k\to0$ as~$k\to+\infty$, 
$$ \gamma_{k+1}+\gamma_{k-1}-2\gamma_k\ge0$$
and~$|\sigma_k|\le\gamma_k$.\end{exercise}

\begin{exercise}\label{FTCPASNo2}
Prove that there exist Fourier Series of integrable functions of period~$1$ whose Fourier coefficients decay to zero arbitrarily slowly.

More specifically, let~$\{\sigma_k\}_{k\in\N}$ be a sequence of real numbers such that~$\sigma_k\to0$ as~$k\to+\infty$.

Prove that there exists~$f\in L^1((0,1))$, periodic of period~1, whose Fourier Series in trigonometric form is
\begin{equation}\label{090-0987yhn-q0woieuxn83gfakowrynt0-w.MAQ-bis}
\sum_{k=0}^{+\infty} \gamma_k\,\cos(2\pi kx),\end{equation}
with~$\gamma_k\ge|\sigma_k|$.
\end{exercise}

\begin{exercise}\label{FTCPASNo2CONYI}
Give another example for Exercise~\ref{FTCPASNo2} relying on Exercises~\ref{VjweWRFVSHEgKFB8ndVIAPMAQ} and~\ref{VjweWRFVSHEgKFB8ndVIAPMAQ-bis}.
\end{exercise}

\begin{exercise}\label{ojdkfnvioewyr098765rewe67890iuhgvgyuuygfr43wdfgh8bft5xs}
Let~$f:\R\to\R$ be periodic of period~$1$ and continuous. Prove that for every~$\epsilon>0$ there exists a trigonometric polynomial~$S_\epsilon$ such that
$$ \sup_{x\in\R}|f(x)-S_\epsilon(x)|<\epsilon.$$
\end{exercise}

\begin{exercise}\label{ojdkfnvioewyr098765rewe67890iuhgvgyuuygfr43wdfgh8bft5xs-EAN}
Prove that trigonometric polynomials are dense in~$L^p((0,1))$ for all~$p\in[1,+\infty)$.
\end{exercise}

\begin{exercise}\label{0ei02r3290234pfp42} Prove that for all~$x\in\left[-\frac12,\frac12\right]$
$$ x^4=\frac1{80}+\sum_{k=1}^{+\infty}\frac{(-1)^k(\pi^2 k^2-6)}{2\pi^4 k^4}\,\cos(2\pi kx).$$\end{exercise}

\begin{exercise} \label{OAKSPiojlfS023a}
For every~$k\in\Z$, let
$$ \mu_k:=\begin{dcases}
\displaystyle\frac{(-1)^{{\operatorname{sign}}(k)}}{i\pi^{\ell}} & {\mbox{ if either~$k=2^\ell$ or $k=-2^\ell$ for some~$\ell\in\N$,}}\\
\\
0& {\mbox{ otherwise,}}
\end{dcases}$$
where
$$ {\operatorname{sign}}(k):=\begin{dcases}-1 &{\mbox{ if }}k<0,\\
0 &{\mbox{ if }}k=0,\\1 &{\mbox{ if }}k>0.
\end{dcases}$$
Does there exist a function~$f:\R\to\R$ which is continuous, periodic of period~$1$, and such that~$\widehat f_k=\mu_k$?

If so, calculate
$$ \lim_{x\to0}\frac{f(x)}{x}.$$
\end{exercise}

\begin{exercise}\label{ZIMadmATnZ-1}
Prove that, for all~$x\in\R$, the series
\begin{equation}\label{LA:zkI1}
\sum_{k=0}^{+\infty}\cos(2\pi kx)
\end{equation}
does not converge, but, for all~$x\in(0,1)$,
\begin{equation}\label{ZIMadmATnZ-1.01}
\sup_{N\in\N}\left|\sum_{k=0}^{N}\cos(2\pi kx)\right|<+\infty.
\end{equation}
\end{exercise}

\begin{exercise}\label{ZIMadmATnZ-2}
Prove that, for all~$x\in\left(0,\frac12\right)$, the series
\begin{equation}\label{LA:zkI2}
\sum_{k=0}^{+\infty}\sin(2\pi kx)
\end{equation}
does not converge, but, for all~$x\in\R$,
\begin{equation}\label{ZIMadmATnZ-1.02}
\sup_{N\in\N}\left|\sum_{k=1}^{N}\sin(2\pi kx)\right|<+\infty.
\end{equation}
\end{exercise}

\begin{exercise}\label{ZIMadmATnZ-3}
Does there exist a periodic function in~$L^1((0,1))$ such that~\eqref{LA:zkI1}
is its Fourier Series in trigonometric form?

And does there exist a periodic function in~$L^1((0,1))$ such that~\eqref{LA:zkI2}
is its Fourier Series in trigonometric form?
\end{exercise}

\section{Convergence in~$L^2((0,1))$}\label{COL2}

We now discuss the reconstruction of a function in~$L^2((0,1))$ from its Fourier Series.
After the previous sections, one might think that this is just a matter of generalising the previous methods
and that this generalisation has been performed merely for academic reasons, such as writing one more paper, getting promoted, making students' life miserable. While these (and especially the latter) are arguably some of the favourite activities of some academics,
the convergence of Fourier Series in~$L^2((0,1))$ is a topic of paramount importance, and it is technically and conceptually different from the pointwise and uniform convergence (though the convergence in~$L^2((0,1))$ implies, up to a sub-sequence,
the pointwise convergence up to a set of null measure, see e.g.~\cite[Theorem~4.9]{MR2759829}, and the
uniform convergence obviously implies the one in~$L^2((0,1))$).

The importance of the development of the theory in class~$L^2((0,1))$ is due to many facts. One is that, as already discussed on page~\pageref{OSK99f-32ut494hg}, it would be too cheap to always assume that we can work with smooth functions, since in many situations of concrete interest this is not the case.

Even more significantly, the development of the Fourier Series in~$L^2((0,1))$ allowed a major progress, since~$L^2((0,1))$ is a Hilbert Space
(see e.g.~\cite[Chapter~5]{MR2759829}), thus providing a universal structure for this kind of spaces
(see e.g.~\cite[Theorems~8.33 and~8.35]{MR3381284}), allowing a construction of Fourier Series that does not depend
on the complex exponentials (see e.g.~\cite[Chapter~8]{MR3381284}), which in turn is of exceptional importance in applications
(for instance, it allows a broad use of Fourier methods in partial differential equations based on eigenfunctions and spectral analysis,
a coherent formulation of many principles of quantum physics, an efficient use of more appropriate objects than exponentials, such as ``wavelets'', to treat problems in edge detection, graphics, data compression, etc., \label{DEDESECTDET:F}
see e.g. Section~\ref{DEDESECTDET}). Since this book
is of introductory type, we will restrict our attention to Fourier Series constructed via the standard basis
of complex exponentials, but the theory of Fourier Series, and especially its Hilbert space analysis, is much richer\footnote{Roughly speaking, the issue with Fourier methods in several practical applications
is that they become ``expensive'' when the signal is ``highly non-smooth'', in the sense that too much information
is required to recapture the minutiae of a signal. Irregular signals, unfortunately, often occur in concrete situations,
for instance due to transient phenomena, and the use of plain Fourier methods in these cases would be too slow in transmitting
and manipulating the signal.

To get around these problems and efficiently capture the information of an irregular signal over a large range of scales, among various techniques, the ``wavelet'' analysis has been introduced. The Fourier Series, in this framework, is replaced by
a series of wavelets depending on two parameters, one accounting for fine scale, the other for spatial translations.
This allows one, in a sense, to ``zoom in'' wherever the irregular behaviour takes place. For more details on this topic, see~\cite{MR1162107, MR2479996, MR2742531}, \cite[Chapter~VI]{PICARDELLO},
and the literature referred to there.} than that
(see e.g.~\cite[Section~2.8]{MR44660}, \cite[Chapters~8 and~12]{MR3381284}, and the references therein to broaden one's horizons on this topic).

For our humble purposes, the main convergence and representation result in~$L^2((0,1))$ goes as follows.

\begin{theorem}\label{THCOL2FB} We have the following results:

  \begin{itemize}
\item[{(i).}] Let~$f\in L^2((0,1))$ be periodic of period~$1$.

Then,
\begin{equation}\label{L2THM.0-01}
{\mbox{the Fourier Series of~$f$ converges to~$f$ in~$L^2((0,1))$.}}\end{equation}

Moreover, \begin{equation}\label{L2THM.0-02}
{\mbox{Parseval's Identity in~\eqref{PARS} holds true.}}\end{equation}

\item[{(ii).}] Conversely, given a sequence of complex numbers~$\{c_k\}_{k\in\Z}$ such that~$c_{-k}=\overline{c_k}$,
if
\begin{equation}\label{NTK5rfct68jn} \sum_{k\in\Z}|c_k|^2<+\infty,\end{equation}
then the series
$$ \sum_{k\in\Z}c_k \,e^{2\pi ikx}$$
converges in~$L^2_{\text{loc}}(\R)$ to a function~$f$ taking real values,
periodic of period~$1$, and such that, for all~$k\in\Z$, $\widehat f_k=c_k$.\end{itemize}
\end{theorem}

\begin{proof} Let~$f\in L^2((0,1))$. Then, by~\eqref{GL4}, for all~$N$, $M\in\N$,
$$ \left\| \sum_{{k\in\Z}\atop{N\le|k|\le M}}\widehat f_k\,e^{2\pi ikx}\right\|_{L^2((0,1))}^2
=\sum_{{k,j\in\Z}\atop{N\le|k|,|j|\le M}}\widehat f_k\,\overline{\widehat f_j}\,\langle e^{2\pi ikx},e^{2\pi ijx}\rangle_{L^2((0,1))}
=\sum_{{k\in\Z}\atop{N\le|k|\le M}}|\widehat f_k|^2,$$
which, by Bessel's Inequality~\eqref{AJSa}, is the tail of a convergent series.

This says that the series
$$ \sum_{{k\in\Z}}\widehat f_k\,e^{2\pi ikx}$$
converges to some function~$g$ in~$L^2((0,1))$, and we extend~$g$ to a periodic function of period~$1$.

Hence, given~$\epsilon>0$, we take~$N_\epsilon\in\N$ sufficiently large such that, for all~$N\ge N_\epsilon$,
$$ \left\|g-\sum_{{k\in\Z}\atop{|k|\le N}}\widehat f_k\,e^{2\pi ikx}\right\|_{L^2((0,1))}\le\epsilon.$$
In this way, for all~$j\in\Z$, we take~$N\ge\max\{j,N_\epsilon\}$ and we see that
\begin{eqnarray*}&&
|\widehat g_j-\widehat f_j|=
\left| \int_0^1 g(x)\,e^{-2\pi ijx}\,dx-\widehat f_j
\right|\\&&\qquad\le\left| \int_0^1 
\sum_{{k\in\Z}\atop{|k|\le N}}\widehat f_k\,e^{2\pi i(k-j)x}\,dx-\widehat f_j
\right|+\left\|g-\sum_{{k\in\Z}\atop{|k|\le N}}\widehat f_k\,e^{2\pi ikx}\right\|_{L^1((0,1))}\\&&\qquad\le\left| \sum_{{k\in\Z}\atop{|k|\le N}}\int_0^1 
\widehat f_k\,e^{2\pi i(k-j)x}\,dx-\widehat f_j
\right|+\left\|g-\sum_{{k\in\Z}\atop{|k|\le N}}\widehat f_k\,e^{2\pi ikx}\right\|_{L^2((0,1))}
\\&&\qquad\le|\widehat f_j-\widehat f_j|+\epsilon\\&&\qquad=\epsilon.
\end{eqnarray*}
As a result, taking~$\epsilon$ as small as we wish, we infer that all the Fourier coefficients of~$g$ coincide with those of~$f$.

Therefore, in light of Theorem~\ref{UNIQ}, we have that~$f=g$, up to a set of null Lebesgue measure, and this establishes~\eqref{L2THM.0-01}.

This and Proposition~\ref{AJSaa} also entail~\eqref{L2THM.0-02}, as desired.

Let us now take a sequence~$\{c_k\}_{k\in\Z}$ as in the statement of Theorem~\ref{THCOL2FB}
and define
$$ F_N(x):=\sum_{{k\in\Z}\atop{|k|\le N}}c_k \,e^{2\pi ikx}.$$
We notice that~$F_N$ is periodic of period~$1$.

Also, $F_N$ takes values in the reals because, for every~$x\in\R$,
\begin{eqnarray*}
\Im\big(F_N(x)\big)&=&\Im\left(c_0+
\sum_{k=1}^N \big(c_k \,e^{2\pi ikx}+
c_{-k} \,e^{-2\pi ikx}\big)\right)\\&=&\Im\left(
\sum_{k=1}^N \big(c_k \,e^{2\pi ikx}+
\overline{c_k} \,e^{-2\pi ikx}\big)\right)\\&=&\Im\left(
\sum_{k=1}^N \big(c_k \,e^{2\pi ikx}+
\overline{c_k \,e^{2\pi ikx}}\big)\right)\\&=&
\sum_{k=1}^N \Im\big(c_k \,e^{2\pi ikx}+
\overline{c_k \,e^{2\pi ikx}}\big)\\&=&0.
\end{eqnarray*}

Finally, by~\eqref{GL4}, for all~$N$, $M\in\N$,
$$ \left\| \sum_{{k\in\Z}\atop{N\le|k|\le M}}c_k\,e^{2\pi ikx}\right\|_{L^2((0,1))}^2
=\sum_{{k,j\in\Z}\atop{N\le|k|,|j|\le M}}c_k\,\overline{c_j}\,\langle e^{2\pi ikx},e^{2\pi ijx}\rangle_{L^2((0,1))}
=\sum_{{k\in\Z}\atop{N\le|k|\le M}}|c_k|^2,$$
which, by~\eqref{NTK5rfct68jn}, is the tail of a convergent series.

This ensures the desired convergence of~$F_N$ in~$L^2((0,1))$, and thus, by periodicity, in~$L^2_{\text{loc}}(\R)$.

Calling~$f$ this $L^2_{\text{loc}}(\R)$-limit (and thus~$L^1_{\text{loc}}(\R)$-limit), we have (see Exercise~\ref {09iulappcoemmdba}) that,
for all~$k\in\Z$,
$$ \widehat f_k=\lim_{N\to+\infty}\widehat F_{N,k}=c_k.$$
The desired result is thereby complete.
\end{proof}

Theorem~\ref{THCOL2FB} (or some form of it) is sometimes called \index{Riesz-Fischer Theorem}
the \emph{Riesz-Fischer Theorem}.

The $L^2$-theory of Fourier Series is very rich and possesses many results that go well beyond the scopes of this book.
For example, while a straightforward consequence of Theorem~\ref{THCOL2FB}(i) 
and of a classical result of Lebesgue spaces (see e.g.~\cite[Theorem 4.9(i)]{MR2759829})
is that the Fourier Series
of a function~$f\in L^2((0,1))$ converges almost everywhere to~$f$ {\em up to a sub-sequence},
a classical, and very deep, result by Carleson \label{Carleson}
is that the Fourier Series
of a function~$f\in L^2((0,1))$ converges almost everywhere to~$f$, {\em without any need of extracting a sub-sequence},
see~\cite{MR199631, MR2091007}
(and also~\cite{MR1783613} for a difficult\footnote{Quite impressively, the MathSciNet review of~\cite{MR1783613} by Loukas Grafakos starts with the sentence ``This is one of the greatest papers written in Fourier analysis''. Then elaborates: ``this remarkable proof of Carleson's almost everywhere convergence theorem is the shortest available today''. And thus finishes with ``This is a truly great paper which brings up a new framework of study in the perplexing and mysterious world of almost everywhere convergence of Fourier series''.}
but relatively short proof).

This result about almost everywhere convergence of Fourier Series has been also extended to the~$L^p((0,1))$ case, for all~$p\in(1,+\infty]$, by Hunt, see~\cite{MR238019}
(and the result is known to be false in~$L^1((0,1))$, since
there exist functions in~$L^1((0,1))$, periodic of period~$1$, \label{MEND:pale}
whose Fourier Series diverges everywhere, see~\cite{zbMATH02586338},
\cite{zbMATH00775031},
\cite[Chapter~VIII, Section~4]{MR1963498},
and~\cite[Chapter~II, Section~3.6]{MR2039503}; see also Section~\ref{EXCE}
\label{0o2jrf98-l.1} below for other convergence issues for Fourier Series).

\begin{exercise}\label{PARSCA}
Let~$f$, $g\in L^2((0,1))$ and periodic of period~$1$.

Prove that\footnote{Equation~\eqref{YOTHCLCHINDMAKSTABS} is also sometimes considered as a version of
\index{Parseval's Identity} Parseval's Identity.}
\begin{equation}\label{YOTHCLCHINDMAKSTABS} \int_0^1 f(x)\,g(x)\,dx=\sum_{k\in\Z}\widehat f_k\,\overline{\widehat g_k}.\end{equation}
\end{exercise}

\begin{exercise}\label{ISOMel2}
Let~$\ell^2$ be the set of all the sequence of complex numbers~$\{c_k\}_{k\in\Z}$ such that~$c_{-k}=\overline{c_k}$
and
\begin{equation*} \sum_{k\in\Z}|c_k|^2<+\infty.\end{equation*}
Using the notation~$c=\{c_k\}_{k\in\Z}$, we endow the space~$\ell^2$ with the norm
$$ \|c\|_{\ell^2}:=\sqrt{\sum_{k\in\Z}|c_k|^2}.$$
Prove that there exists a linear and invertible map~$T: L^2((0,1))\to\ell^2$ such that, for every~$f\in L^2((0,1))$,
\begin{equation}\label{L2THM.0-02-makdefg} \|T(f)\|_{\ell^2}=\|f\|_{L^2((0,1))}.\end{equation}
\end{exercise}

\begin{exercise}\label{OJSNILCESNpfa.sdwpqoed-23wedf-87a}
In analogy with Exercise~\ref{OJSNILCESNpfa.sdwpqoed-23wedf}, one may wonder whether Fourier Series can be integrated ``termwise'' (of course, with some attention to linear terms which are not periodic functions and which can still pop up in the integration process).

In this sense, consider a periodic function~$f\in L^2((0,1))$ and define
$$ F(x):=\int_0^x f(t)\,dt.$$
Prove that, for all~$x\in(0,1)$,
$$ F(x)=\widehat f_0\,x+\sum_{k\in\Z\setminus\{0\}} \frac{\widehat f_k}{2\pi ik}\,e^{2\pi ikx},$$ with the latter series converging uniformly in~$(0,1)$.
\end{exercise}

\section{\faBomb Convergence in~$L^p((0,1))$}\label{COLP}

As a (rather nontrivial) variant of the material exposed in Section~\ref{COL2},
let us mention that the convergence in norm of Fourier Series has also been studied in the class~$L^p((0,1))$.
In particular, the claim in~\eqref{L2THM.0-01} can be carried through the class~$L^p((0,1))$ for all~$p\in(1,+\infty)$.
More precisely,  we have that:

\begin{theorem}\label{0pXe-244kf.1} If~$f$ is a periodic function of period~$1$ belonging to~$L^p((0,1))$ with~$p\in(1,+\infty)$, then the Fourier series of~$f$ converges to~$f$ in~$L^p((0,1))$.\end{theorem}

See Exercise~\ref{0pXe-244kf.1EX} for a proof of this statement (see also~\cite[Lemma~1.8 and Corollary~3.7]{MR1800316},
\cite[Chapter~II, Section~1]{MR2039503}, \cite{MIAOJING}, \cite{LZHANG},
and~\cite[Section~7.6.7]{PICARDELLO}
for more details).

We remark that
\begin{itemize}
\item When~$p\in[1,+\infty)$, a weaker statement will be independently proved in Theorem~\ref{FeJ-th.3}, where we will show that
a suitable average of the Fourier Series converge to the original function in~$L^p((0,1))$,
\item When~$p=+\infty$, the claim in Theorem~\ref{0pXe-244kf.1} does not hold true, see Theorem~\ref{BaireCategoryTh},
\item When~$p=1$, the claim in Theorem~\ref{0pXe-244kf.1} does not hold true, see Exercise~\ref{KDKAMSDW:SDCMILCSKMc},
\item When~$p=1$, Theorem~\ref{FeJ-th.3} will ensure that
a suitable average of the Fourier Series converge to the original function in~$L^1((0,1))$.
\end{itemize}

For other convergence results related to Fourier Series in~$L^p((0,1))$ see e.g.~\cite[Part~1, Chapter~III]{MR1408905}
and the references therein.

\begin{exercise}\label{0pXe-244kf.1EX-PRE.pe}
If~$S$ is a trigonometric polynomial of the form
\begin{equation}\label{cjamsdperrea-23} S(x)=\sum_{{k\in\Z}\atop{|k|\le M}} c_k \,e^{2\pi i k x},\end{equation}
with
\begin{equation}\label{cjamsdperrea}
c_{-k}=\overline{c_k},
\end{equation}
one defines the conjugate trigonometric polynomial as
\begin{equation} \label{0pXe-244kf.1EX-PRE.pe-djenfEF}
{S}^\Diamond(x)=-i\sum_{{k\in\Z}\atop{|k|\le M}} {\operatorname{sign}}(k)\, c_k \,e^{2\pi i k x},\end{equation}
where\footnote{A similar notion of conjugacy will be introduced for the periodic Poisson Kernel in~\eqref{POISSONKERN-ex2-form1},
\eqref{POISSONKERN-ex2-form1-ILCO}, and~\eqref{LAMCIKCONJA0-2}. See also Exercise~\ref{LPmPOIKE01ijdqsp21rj90jtm}
and compare~\eqref{0pXe-244kf.1EX-PRE.pe-djenfEF}
with footnote~\ref{MKDFOOTOBECMSCoojwlf395utjhgnhbdfvgwherh5} on page~\pageref{MKDFOOTOBECMSCoojwlf395utjhgnhbdfvgwherh5}.}
$$ {\operatorname{sign}}(k):=\begin{dcases}
1 &{\mbox{ if }}k>0,\\
-1 &{\mbox{ if }}k<0,\\
0 &{\mbox{ if }}k=0.
\end{dcases}$$

We also consider the \index{Riesz Projections} Riesz Projections\footnote{The next exercises are ancillary to the proof of the convergence of Fourier Series in~$L^p((0,1))$, as stated in Theorem~\ref{0pXe-244kf.1}. The proof of this deep result will be finally accomplished in Exercise~\ref{0pXe-244kf.1EX}. In this, the usefulness of the Riesz Projections is greater than what it may look at a first glance. 

Indeed, thanks to the density of trigonometric polynomials (see Exercise~\ref{ojdkfnvioewyr098765rewe67890iuhgvgyuuygfr43wdfgh8bft5xs-EAN}), the gist of the proof of Theorem~\ref{0pXe-244kf.1} could be first to approximate the given periodic function~$f\in L^p((0,1))$ with a trigonometric polynomial~$S_\epsilon$, say with~$\|S_\epsilon-f\|_{L^p((0,1))}\le\epsilon$, where~$\epsilon>0$ can be taken as small as we wish.

Now, dealing with trigonometric polynomials may seem as easy as pie, but perhaps it is not so, let's see why. On the one hand, to establish Theorem~\ref{0pXe-244kf.1}, we need to show that
$$ \lim_{N\to+\infty}\|S_{N,f}-f\|_{L^p((0,1))}=0.$$
It is therefore tempting to ``swap'' $f$ with~$S_\epsilon$ and~$S_{N,f}$ with the corresponding Fourier Sum of~$S_\epsilon$, that is~$S_{N,S_\epsilon}$, i.e. to engage with the inequality
$$ \|S_{N,f}-f\|_{L^p((0,1))}\le \|S_{N,f}-S_{N,S_\epsilon}\|_{L^p((0,1))}+\|S_{N,S_\epsilon}-S_\epsilon\|_{L^p((0,1))}+\|S_\epsilon-f\|_{L^p((0,1))}.$$
The last term is fine, being bounded by~$\epsilon$. The second last term is also fine, actually being null for~$N$ large enough, since the Fourier Series of a trigonometric polynomial is the trigonometric polynomial itself (and, more precisely, $S_{N,S_\epsilon}=S_\epsilon$ whenever~$N$ gets larger than the degree of~$S_\epsilon$ (see Exercise~\ref{PKS0-3-21}).

Hence, the real issues come from the term~$\|S_{N,f}-S_{N,S_\epsilon}\|_{L^p((0,1))}$, or equivalently~$\|S_{N,S_\epsilon-f}\|_{L^p((0,1))}$. It would be highly desirable to bound by above this term by something like~$\|S_\epsilon-f\|_{L^p((0,1))}$, which we know is controlled by~$\epsilon$.

But this essentially reduces the proof of Theorem~\ref{0pXe-244kf.1} to an estimate of the type~$\|S_{N,g}\|_{L^p((0,1))}\le C\,\|g\|_{L^p((0,1))}$, for some~$C>0$ (in fact, for Theorem~\ref{0pXe-244kf.1} we would only need this estimate with~$g:=S_\epsilon-f$, but proving an estimate of this type for this specific~$g$ looks as hard as proving it for any~$g$ whatsoever).

Hence, the essence of the proof of Theorem~\ref{0pXe-244kf.1} will be to estimate the norm of any Fourier Sum with the norm of the original function (and this will be accomplished in Exercise~\ref{0pXe-244kf.1EX-SK0-1o23lk-0923}). To this end, one may try to focus on the case in which~$g$ is just a trigonometric polynomial, again by some density argument. The catch however is now that the degree of this trigonometric polynomial may be very high, and in particular higher than the Fourier Sum~$S_{N,g}$. Therefore, it would be very convenient to have a ``cutting device'' which, given a trigonometric polynomial, only preserves some band of frequencies (this mathematical idea will be also retaken in the application highlighted in Section~\ref{SE:FIL:023er-1}).

This cutting device is provided precisely by the Riesz Projections (strictly speaking, these objects cut at the level of the zero frequency, but, as we will see in Exercise~\ref{0pXe-244kf.1EX-PRE.0p23rl.oqjsdlwnfvMSAS01maspe}, different frequencies can be easily selected by multiplications with complex exponentials).

Though very handy, the Riesz Projections cutting devices have a drawback, namely, as expressed in~\eqref{EQARProjections}, they require the introduction of a ``conjugate'' function. That is, if we use Riesz Projections to write Fourier Sums to estimate their norms, it will be necessary to estimate also the norm of this conjugate object, and this will be taken care of in Exercise~\ref{0pXe-244kf.1EX-SK}.

The bottom line is that the proof of Theorem~\ref{0pXe-244kf.1} is challenging, but also fun and culturally very elevated.
}
\begin{equation*}\begin{split}
& P^+_{S}(x):=\sum_{{k\in\Z}\atop{1\le k\le M}}c_k \,e^{2\pi i k x}
\\{\mbox{and }}\qquad&
P^-_{S}(x):=\sum_{{k\in\Z}\atop{-M\le k\le -1}}c_k \,e^{2\pi i k x}.
\end{split}\end{equation*}

Prove that
\begin{equation}\label{EQARProjections}\begin{split}
& S=P^+_{S}+P^-_{S}+c_0\,,\\&
{S}^\Diamond=-iP^+_{S}+iP^-_{S}\,,\\{\mbox{and }}\quad&
P^+_{S}=\frac{S+i{S}^\Diamond-c_0}2.
\end{split}\end{equation}
\end{exercise}

\begin{exercise}\label{0pXe-244kf.1EX-PRE.0p23rl} In the notation of Exercise~\ref{0pXe-244kf.1EX-PRE.pe},
prove that~$S$ and~${S}^\Diamond$ take values in the reals.
\end{exercise}

\begin{exercise}\label{0pXe-244kf.1EX-PRE.0p23rlmanyOnt} In the notation of Exercise~\ref{0pXe-244kf.1EX-PRE.pe},
prove that~$P^+_{S}$ and~$P^-_{S}$ do not necessarily take values in the reals.\end{exercise}

\begin{exercise}\label{0pXe-244kf.1EX-PRE.0p23rl.oqjsdlwnfvMSAS01} Let~$S_1$ and~$S_2$ be two trigonometric polynomials of the form
\begin{equation*} S_j(x)=\sum_{{k\in\Z}\atop{|k|\le M}} c_{j,k} \,e^{2\pi i k x},\end{equation*}
satisfying the reality condition~$c_{j,-k}=\overline{c_{j,k}}$, for~$j\in\{1,2\}$.

Let~$S:=S_1+iS_2$ and~$P^+_S(x):=P^+_{S_1}+iP^+_{S_2}$.

Prove that
$$ S(x)=\sum_{{k\in\Z}\atop{|k|\le M}} c_{k} \,e^{2\pi i k x}$$
and that
\begin{equation}\label{gbsdfnigbfa.S01maspe}
P^+_S(x)=\sum_{{k\in\Z}\atop{1\le k\le M}}c_k \,e^{2\pi i k x},\end{equation}
with~$c_k:=c_{1,k}+ic_{2,k}$.
\end{exercise}

\begin{exercise}\label{0pXe-244kf.1EX-PRE.0p23rl.oqjsdlwnfvMSAS01maspe} In the notation of Exercise~\ref{0pXe-244kf.1EX-PRE.0p23rl.oqjsdlwnfvMSAS01}, let
$$ S^{(n)}(x):=e^{2\pi inx} S(x)$$
and prove that, for all~$n\in\Z$ with~$n+ M\ge1$,
$$ e^{-2\pi inx}P^+_{S^{(n)}}(x)=\sum_{{k\in\Z}\atop{1-n\le k\le M}} c_{k} \,e^{2\pi i k x}.$$
\end{exercise}

\begin{exercise}\label{0pXe-244kf.1EX-PRE.0p23rl.oqjsdlwnfvMSAS01maspe.bismdc} In the notation of Exercise~\ref{0pXe-244kf.1EX-PRE.0p23rl.oqjsdlwnfvMSAS01maspe}, 
prove that, for all~$N\in\N$ with~$N\le M-1$,
$$ e^{-2\pi i(N+1)x}P^+_{S^{(N+1)}}(x)-e^{2\pi iNx}P^+_{S^{(-N)}}(x)
=\sum_{{k\in\Z}\atop{-N\le k\le N}} c_{k} \,e^{2\pi i k x}.$$
\end{exercise}

\begin{exercise}\label{0pXe-244kf.1EX-PRE} In the notation of Exercise~\ref{0pXe-244kf.1EX-PRE.pe}, assume that~$c_0=0$.

Given~$\ell\in\N\cap[1,+\infty)$, prove that
$$ \int_0^1 \big( S(x)+i{S}^\Diamond(x)\big)^{\ell}\,dx=0.$$
\end{exercise}

\begin{exercise}\label{0pXe-244kf.1EX-PRE-012} In the notation of Exercise~\ref{0pXe-244kf.1EX-PRE.pe}, assume that~$c_0=0$.

Given~$\ell\in\N\cap[1,+\infty)$ and~$M\in\N$, prove that
$$ \|{S}^\Diamond\|_{L^{2\ell}((0,1))}^{2\ell}
\le\sum_{j=1}^{2\ell} \binom{2\ell}{j}\,\|S\|_{L^{2\ell}((0,1))}^{j}\,\|{S}^\Diamond\|_{L^{2\ell}((0,1))}^{2\ell-j}.$$
\end{exercise}

\begin{exercise}\label{0pXe-244kf.1EX-PRE-012.b} Let~$\ell\in\N\cap[1,+\infty)$ and~$X\ge0$.
Assume that
$$ X^{2\ell}\le\sum_{j=1}^{2\ell}\binom{2\ell}{j} X^{2\ell-j}.$$
Prove that there exists~$C>0$, depending only on~$\ell$, such that~$X\le C$.
\end{exercise}

\begin{exercise}\label{0pXe-244kf.1EX-PRE-012.c} 
In the notation of Exercise~\ref{0pXe-244kf.1EX-PRE.pe}, let~$\ell\in\N\cap[1,+\infty)$.
Prove that there exists~$C>0$, depending only on~$\ell$, such that for all trigonometric polynomials~$S$ one has that
\begin{equation}\label{0pXe-244kf.1EX-PRE-012.d} \|{S}^\Diamond\|_{L^{2\ell}((0,1))}\le C\|S\|_{L^{2\ell}((0,1))}.\end{equation}
\end{exercise}

\begin{exercisesk}\label{0pXe-244kf.1EX-SK} Let~$p\in(1,+\infty)$. Prove that there exists~$C>0$, depending only on~$p$,
such that for all trigonometric polynomials~$S$ we have that
\begin{equation}\label{0pXe-244kf.1EX-PRE-012.d-96} \| {S}^\Diamond\|_{L^p((0,1))}\le C\,\|S\|_{L^p((0,1))}.\end{equation}
\end{exercisesk}

\begin{exercise}\label{0pXe-244kf.1EX-SK-ilcoM} In the setting of Exercise~\ref{0pXe-244kf.1EX-SK},
show that~\eqref{0pXe-244kf.1EX-PRE-012.d-96} holds true for complex-valued
trigonometric polynomials as well, namely if~$S_1$ and~$S_2$ are real-valued trigonometric polynomials and~$S:=S_1+iS_2$, then, setting~${S}^\Diamond:={S}^\Diamond_1+i{S}^\Diamond_2$ and up to renaming~$C$,
\begin{equation*}
\| {S}^\Diamond\|_{L^p((0,1),\C)}\le C\,\|S\|_{L^p((0,1),\C)}.\end{equation*}
\end{exercise}

\begin{exercisesk}\label{0pXe-244kf.1EX-SK0-1o23lk-0923} Let~$p\in(1,+\infty)$. Prove that there exists~$C>0$, depending only on~$p$,
such that for all~$f\in L^p((0,1))$ periodic of period~$1$ and all~$N\in\N$ we have that
\begin{equation}\label{0pXe-244kf.1EX-PRE-012.d-96-0923k} \| S_{N,f}\|_{L^p((0,1))}\le C\,\|f\|_{L^p((0,1))}.\end{equation}
\end{exercisesk}

\begin{exercise}\label{0pXe-244kf.1EX}
Prove Theorem~\ref{0pXe-244kf.1}.\end{exercise}

\section{\faBomb Averaging procedures for Fourier Series}\label{FEJERKESE}

So far, we have seen that the very desirable property that the Fourier Series reconstruct a given periodic function holds true on a multitude of occasions, but always under suitable hypotheses (such as a regularity assumption on the function, or a decay assumption on its Fourier coefficients). When these assumptions are not met, the Fourier Series may fail to reconstruct the given function: this rather unpleasant situation was mentioned already on page~\pageref{MEND:pale} and will be further discussed in Section~\ref{EXCE}.

Such is life. However, even in adverse situations, one does not lose hope and courage, and can still rely on the rule of thumb that convergence ``in average'' could still hold when stronger notions of convergence fail (see Exercise~\ref{CESA}).

In this sense, it would be worth to look at the average of~$S_N$ and check if this converges to~$f$: to this end, one could recall Corollary~\ref{9.0km2df0.2wnedfo.1-2x} and write the average of~$S_N$ in the concise form
\begin{equation}
  \label{eq4}
  \frac1N\sum_{k=0}^{N-1}S_k(x)=
  \frac1N\sum_{k=0}^{N-1} f\star D_k(x)=f\star\left(
    \frac1N\sum_{k=0}^{N-1}  D_k(x)\right),
\end{equation}
where
\begin{equation*}
  D_N(x) = \frac{\sin\big((2N+1)\pi x\big)}{\sin(\pi x)},
\end{equation*}
as introduced in Lemma \ref{KASMqwdfed123erDKLI}.

The last term in \eqref{eq4} suggests that it is worth investigating
the ``average of the Dirichlet Kernel'', which is called\footnote{We observe that, to start familiarising with it, we had already
introduced the Fej\'er Kernel in Exercise~\ref{IFEJ} and we have used it also in
Exercise~\ref{ESIPNo-2}. This confirms the fact that the Fej\'er Kernel is indeed useful, especially when convergence issues become delicate.}%
the \index{Fej\'er Kernel} \emph{Fej\'er Kernel}, namely
\begin{equation}\label{01foqld03rkfg.GnwedRA.a}
  F_{N}(x):=\frac{1}{N}\sum_{k=0}^{N-1}D_{k}(x).
\end{equation}

We will see here below that indeed an analysis of this kernel will lead to a convenient notion of convergence in average for the Fourier Series, which also presents several useful consequences.

To start this analysis, here is a summary of some basic properties of the Fej\'er Kernel:

\begin{theorem}
We have that
\begin{equation}\label{IFEJ-FORM} F_N(x)=\frac1{N}\left( \frac{\sin(N\pi x)}{\sin(\pi x)}\right)^2\end{equation}
and
\begin{equation}\label{IFEJ-FORM.bi} F_N(x)=\sum_{{k\in\Z}\atop{|k|\le N-1}}\left(1-\frac{|k|}{N}\right)\,e^{2\pi ikx}.\end{equation}

Moreover, for all~$x\in\R$,
\begin{equation}\label{IFEJ-FORM2}
F_N(x)\ge0.\end{equation}

In addition,
\begin{equation}\label{IFEJ-FORM2.bn} \int_0^1 F_N(x)\,dx=1.\end{equation}
\end{theorem}

\begin{proof} See Exercise~\ref{IFEJ} for the proof of the identity in~\eqref{IFEJ-FORM}
and Exercise~\ref{FEJPO} for the simple proof of the inequality in~\eqref{IFEJ-FORM2}.

See also Exercise~\ref{IFEJ-FORM2.bn-EX} for a proof of~\eqref{IFEJ-FORM2.bn}.

In addition, in view of~\eqref{PAKSw-L4}, exchanging the order of summation we have that
\begin{eqnarray*}
&&\frac{1}{N}\sum_{k=0}^{N-1}D_{k}(x)=
\frac{1}{N}\sum_{k=0}^{N-1}\sum_{j=-k}^k e^{-2\pi i jx}=
\frac{1}{N}\sum_{{k\in\Z}\atop{0\le k\le N-1}}\sum_{{j\in\Z}\atop{|j|\le k}} e^{-2\pi i jx}\\
&&\qquad=\frac{1}{N}\sum_{{j\in\Z}\atop{|j|\le N-1}}\sum_{{k\in\Z}\atop{|j|\le k\le N-1}} e^{-2\pi i jx}=
\frac{1}{N}\sum_{{j\in\Z}\atop{|j|\le N-1}}\big(N-|j|\big) e^{-2\pi i jx},
\end{eqnarray*}
which establishes~\eqref{IFEJ-FORM.bi}. 
\end{proof}

\begin{figure}[h]
\includegraphics[height=4.3cm]{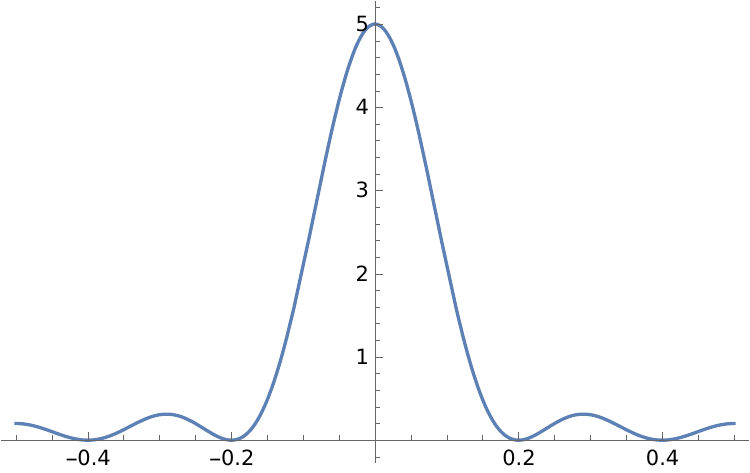}$\quad$\includegraphics[height=4.3cm]{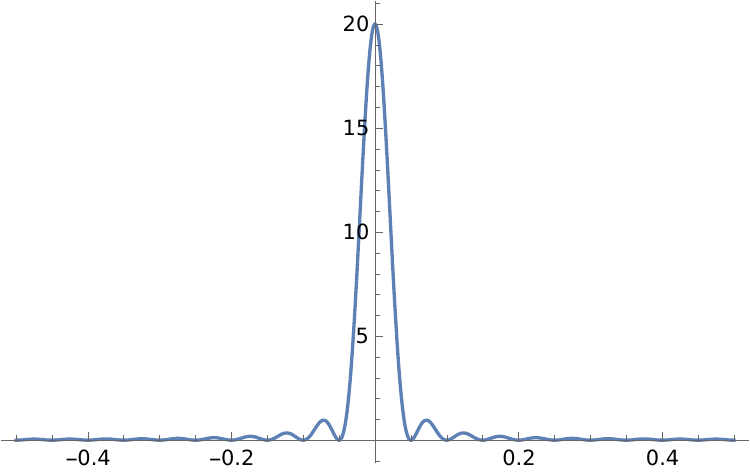}
\centering
\caption{The Fej\'er Kernel (left, with~$N=5$; right, with~$N=20$).}
\label{amdi.fe}
\end{figure}

See Figure~\ref{amdi.fe} for a sketch of the Fej\'er Kernel. The pattern represented there can be compared with that of Figure~\ref{amdi}. In a sense, the Fej\'er Kernel is shown to enjoy ``similar'' concentration properties to the Dirichlet Kernel, with two major advantages that pop up from the comparison of these figures:
firstly, in contrast with the sign-changing Dirichlet Kernel, the Fej\'er Kernel is non-negative; secondly, the oscillations of the Fej\'er Kernel are ``milder'' than those\footnote{Sometimes, it is said that the Fej\'er Kernel is a ``good kernel'', while the Dirichlet Kernel is not. Though a sentence of this sort is perhaps a bit ungenerous to the Dirichlet Kernel, one can formalise and fruitfully use the notion of good kernel, see~\cite{MR1970295} for further details on this aspect of the general theory.} of the Dirichlet Kernel,
as a consequence of the ``averaging procedure'' in~\eqref{01foqld03rkfg.GnwedRA.a} (see also Exercise~\ref{GKFND}).

For further observation comparing the Dirichlet and Fej\'er Kernels see~\cite[Volume~1, Section~52]{MR171116}.

The average of Fourier Sums can be conveniently represented via the Fej\'er Kernel. Indeed, we have the following result, which can be considered as an ``averaged version'' of Corollary~\ref{9.0km2df0.2wnedfo.1-2x}

\begin{corollary}
Let~$f\in L^1((0,1))$ be periodic of period~$1$.

Then, for all~$N\in\N$,
\begin{equation}\label{023wefv4567uygfdfgyhuizo0iwhfg0eb5627e01.p1rebi}
\frac1N\sum_{k=0}^{N-1}S_k(x)=\int_0^1 f(y)\,F_N(x-y)\,dy
\end{equation}
and\footnote{Using the notation of the periodic convolution introduced in Exercise~\ref{pecoOSJHN34jfgj}, one can also
rewrite~\eqref{023wefv4567uygfdfgyhuizo0iwhfg0eb5627e01.p1rebi} in the form
$$ \frac1N\sum_{k=0}^{N-1}S_k(x)=f\star F_N(x).$$}
\begin{equation}\label{023wefv4567uygfdfgyhuizo0iwhfg0eb5627e01.p1rebipi}f(x)-
\frac1N\sum_{k=0}^{N-1}S_k(x)=
\frac1{N}\int_0^1 \big(f(x)-f(y)\big)\,\left( \frac{\sin(N\pi (x-y))}{\sin(\pi( x-y))}\right)^2\,dy.
\end{equation}
\end{corollary}

\begin{proof} Notice that~\eqref{023wefv4567uygfdfgyhuizo0iwhfg0eb5627e01.p1rebi} follows by exchanging the finite sum with the integral after~\eqref{023wefv4567uygfdfgyhuizo0iwhfg0eb5627e01.pre} and
recalling~\eqref{01foqld03rkfg.GnwedRA.a}.

Recalling~\eqref{IFEJ-FORM2.bn},
we see that the claim in~\eqref{023wefv4567uygfdfgyhuizo0iwhfg0eb5627e01.p1rebipi}
is a consequence of~\eqref{IFEJ-FORM} and~\eqref{023wefv4567uygfdfgyhuizo0iwhfg0eb5627e01.p1rebi}.\end{proof}

We now discuss the convergence properties of the averages of the Fourier Sums.
We deal with three specific cases, namely, pointwise convergence at continuity points, uniform convergence to continuous functions, and convergence in~$L^p((0,1))$ (see also Exercise~\ref{FEJPOAE1edc4.1} for the almost everywhere convergence).

\begin{theorem}\label{FeJ-th.1}
Let~$f\in L^1((0,1))$ be periodic of period~$1$.

Assume that~$f$ is continuous at some point~$x\in\R$.

Then,
$$ \lim_{N\to+\infty}
\frac1N\sum_{k=0}^{N-1}S_k(x)=f(x).$$
\end{theorem}

\begin{theorem}\label{FeJ-th.2}
Let~$f$ be a continuous function, periodic of period~$1$.

Then,
$$ \lim_{N\to+\infty}\sup_{x\in\R}\left|
\frac1N\sum_{k=0}^{N-1}S_k(x)-f(x)
\right|=0.$$
\end{theorem}

\begin{theorem}\label{FeJ-th.3}
Let~$p\in[1,+\infty)$. Let~$f\in L^p((0,1))$ be periodic of period~$1$.

Then,
\begin{equation}\label{FeJ-th.3.NOIJMS.021324v4c25} \lim_{N\to+\infty}\left\|
\frac1N\sum_{k=0}^{N-1}S_k-f
\right\|_{L^p((0,1))}=0.\end{equation}
\end{theorem}

\begin{proof}[Proof of Theorem~\ref{FeJ-th.1}]
Let~$\delta\in\left(0,\frac14\right)$. Then, 
we define
$$ {\mathcal{D}}_\delta:=\bigcup_{\ell\in\Z}(\ell+x-\delta,\,\ell+x+\delta)$$
and we notice that if~$y\in(0,1)\setminus{\mathcal{D}}_\delta$ then~$\big|\sin(\pi(x-y))\big|\ge\sin(\pi\delta)$.

As a result,
\begin{equation}\label{FAGBS:0-0.1}
\begin{split}&
\left|\int_{(0,1)\setminus{\mathcal{D}}_\delta} \big(f(x)-f(y)\big)\,\left( \frac{\sin(N\pi (x-y))}{\sin(\pi (x-y))}\right)^2\,dy
\right|\\
&\qquad\le
\frac1{\sin^2(\pi\delta)}
\int_{0}^1 \big(|f(x)|+|f(y)|\big)\,dy\\&\qquad\le
\frac{|f(x)|+\|f\|_{L^1((0,1))}}{\sin^2(\pi\delta)}.
\end{split}
\end{equation}

Furthermore,
\begin{equation}\label{FAGBS:0-0.2}
\begin{split}&
\left|\int_{(0,1)\cap{\mathcal{D}}_\delta} \big(f(x)-f(y)\big)\,\left( \frac{\sin(N\pi (x-y))}{\sin(\pi (x-y))}\right)^2\,dy
\right|\\
&\qquad\le2\int_{\{|t|\le\delta\}} \big|f(x)-f(x+t)\big|\,\left( \frac{\sin(N\pi t)}{\sin(\pi t)}\right)^2\,dt.
\end{split}
\end{equation}

Now, since~$f$ is continuous at~$x$, given~$\epsilon>0$, we can choose~$\delta_\epsilon\in\left(0,\frac14\right)$ such that
whenever~$|t|\le\delta_\epsilon$ we have that~$\big|f(x)-f(x+t)\big|\le\epsilon$. In this way, we infer from~\eqref{023wefv4567uygfdfgyhuizo0iwhfg0eb5627e01.p1rebipi}, \eqref{FAGBS:0-0.1}, and~\eqref{FAGBS:0-0.2} that, for every~$\epsilon>0$,
\begin{eqnarray*}&&\left|
f(x)-\frac1N\sum_{k=0}^{N-1}S_k(x)\right|\le
\frac1{N}\left[
\frac{|f(x)|+\|f\|_{L^1((0,1))}}{\sin^2(\pi\delta_\epsilon)}
+2\epsilon\int_{\{|t|\le\delta_\epsilon\}} \left( \frac{\sin(N\pi t)}{\sin(\pi t)}\right)^2\,dt
\right].
\end{eqnarray*}
Hence, taking the limit as~$N\to+\infty$ and recalling~\eqref{IFEJ-FORM} and~\eqref{IFEJ-FORM2.bn},
\begin{eqnarray*}
\lim_{N\to+\infty}\left|
f(x)-\frac1N\sum_{k=0}^{N-1}S_k(x)\right|\le
\lim_{N\to+\infty}\frac{2\epsilon}{N}\int_{0}^1 \left( \frac{\sin(N\pi t)}{\sin(\pi t)}\right)^2\,dt
=2\epsilon\lim_{N\to+\infty}\int_{0}^1 F_N(t)\,dt=2\epsilon.
\end{eqnarray*}
Since we can now take~$\epsilon$ arbitrarily small, we have established the desired result.
\end{proof}

\begin{proof}[Proof of Theorem~\ref{FeJ-th.2}] This is a careful variation of the proof of Theorem~\ref{FeJ-th.1}.
Namely, in this case~$f$ is uniformly continuous and therefore, for every~$\epsilon>0$ there exists~$\delta_\epsilon>0$ such that
$$ \sup_{{x\in\R}\atop{|t|\le\delta_\epsilon}}\big|f(x)-f(x+t)\big|\le\epsilon.$$
Hence, we deduce from~\eqref{FAGBS:0-0.1} and~\eqref{FAGBS:0-0.2} that
\begin{eqnarray*}\sup_{x\in\R}\left|\int_0^1\big(f(x)-f(y)\big)\,\left( \frac{\sin(N\pi (x-y))}{\sin(\pi (x-y))}\right)^2\,dy\right|\le
\frac{2\|f\|_{L^\infty((0,1))}}{\sin^2(\pi\delta_\epsilon)}+2\epsilon\int_{\{|t|\le\delta_\epsilon\}} \left( \frac{\sin(N\pi t)}{\sin(\pi t)}\right)^2\,dt.\end{eqnarray*}

Combining this with~\eqref{IFEJ-FORM}, \eqref{IFEJ-FORM2.bn}, and~\eqref{023wefv4567uygfdfgyhuizo0iwhfg0eb5627e01.p1rebipi}, we obtain that, for every~$\epsilon>0$,
\begin{eqnarray*}\lim_{N\to+\infty}
\sup_{x\in\R}\left|f(x)-\frac1N\sum_{k=0}^{N-1}S_k(x)\right|\le
\lim_{N\to+\infty}\left(
\frac{2\|f\|_{L^\infty((0,1))}}{N\sin^2(\pi\delta_\epsilon)}+2\epsilon\int_{0}^1 F_N(t)\,dt\right)=2\epsilon,
\end{eqnarray*}from which we obtain the desired result.
\end{proof}

\begin{proof}[Proof of Theorem~\ref{FeJ-th.3}] Let~$\epsilon>0$ and~$f_\epsilon$
be a continuous function, periodic of period~$1$, such that~$\|f-f_\epsilon\|_{L^p((0,1))}\le\epsilon$ (see e.g.~\cite[Theorem~9.6]{MR3381284}).

Thus, we use Young's Convolution Inequality (see e.g.~\cite[Theorem~9.2]{MR3381284}) and we deduce that, for every~$g\in L^p((0,1))$,
$$ \|g\star F_N\|_{L^p((0,1))}\le \|F_N\|_{L^1((0,1))}\, \|g\|_{L^p((0,1))}.$$

Hence, since (see~\eqref{IFEJ-FORM2} and Exercise~\ref{IFEJ-FORM2.bn-EX})
$$ \|F_N\|_{L^1((0,1))}=\int_0^1 F_N(x)\,dx=1,$$
we arrive at \begin{equation*}
\|g\star F_N\|_{L^p((0,1))}\le \|g\|_{L^p((0,1))}.\end{equation*}

Moreover, after Theorem~\ref{FeJ-th.2} we find~$N_\epsilon\in\N$ such that, for all~$N\ge N_\epsilon$,
$$ \sup_{x\in\R}\left|
\frac1N\sum_{k=0}^{N-1}S_{k,f_\epsilon}(x)-f_\epsilon(x)
\right|\le\epsilon.$$

In particular, for all~$N\ge N_\epsilon$,
\begin{eqnarray*}&&
\left\|
\frac1N\sum_{k=0}^{N-1}S_{k,f}-f
\right\|_{L^p((0,1))}\\&\le&
\left\|
\frac1N\sum_{k=0}^{N-1}S_{k,f}-\frac1N\sum_{k=0}^{N-1}S_{k,f_\epsilon}
\right\|_{L^p((0,1))}
+
\left\|
\frac1N\sum_{k=0}^{N-1}S_{k,\epsilon}-f_\epsilon
\right\|_{L^p((0,1))}+\|f-f_\epsilon\|_{L^p((0,1))}
\\&\le&\left\|
\frac1N\sum_{k=0}^{N-1}S_{k,f-f_\epsilon}
\right\|_{L^p((0,1))}
+\epsilon+\epsilon
\\&=&\left\|(f-f_\epsilon)\star F_N
\right\|_{L^p((0,1))}
+2\epsilon\\&\le&\|f-f_\epsilon\|_{L^p((0,1))}
+2\epsilon\\&\le&3\epsilon,
\end{eqnarray*}
yielding the desired result.
\end{proof}

With the results developed so far, we can now check that the convergence of Fourier Series is appropriately
designed to ``converge to the right function'' whenever they converge almost everywhere:

\begin{theorem}\label{ILGIOB}
Let~$f\in L^1((0,1))$ be periodic of period~$1$ and assume that for almost every~$x\in(0,1)$
we have that~$S_N(x)$ converges to some quantity~$g(x)$ as~$N\to+\infty$.

Then, for almost every~$x\in(0,1)$, it holds that~$g(x)=f(x)$.
\end{theorem}

\begin{proof} Let
$$A_N(x):=\frac1N\sum_{k=0}^{N-1}S_k(x).$$
By Theorem~\ref{FeJ-th.3}, we know that, as~$N\to+\infty$, 
$A_N$ converges to~$f\in L^1((0,1))$.

As a consequence (see e.g.~\cite[Theorem 4.9]{MR2759829}), for a certain\footnote{Actually,
one could also argue that the full sequence~$A_N$ converges almost everywhere, see Exercise~\ref{FEJPOAE1edc4.1},
but this fact is not needed here.} sub-sequence~$N_j$
we have that
\begin{equation}\label{pqosjdl9.01}
{\mbox{$A_{N_j}$ converges to~$f$ almost everywhere in~$(0,1)$.}}\end{equation}

By hypothesis, we also know that~$S_N$ converges to~$g$ almost everywhere in~$(0,1)$, therefore (see Exercise~\ref{CESA})
also~$A_N$ converges to~$g$ almost everywhere in~$(0,1)$.

Comparing with~\eqref{pqosjdl9.01}, we conclude that~$f$ and~$g$ coincide almost everywhere in~$(0,1)$.
\end{proof}

For other fine notions of convergence of Fourier Series, see e.g.~\cite[Chapters~5 and~6]{MR545506},
\cite[Part1, Section~2.2 and Chapter~IV]{MR1408905},
and~\cite[Section~3.4.1]{MR3243734}.
See also~\cite[page~60]{MR3601106} for interesting discussions comparing the
``standard'' summation method of Fourier Series and the Fej\'er averaging procedure.
The literature presents several ingenious techniques to ``regularize'' divergent series to somewhat produce a convergent
result: in the setting of Fourier Series, a rather popular one is the notion of \index{Abel-Poisson summability}
Abel-Poisson summability. We do not enter into this field here and we refer the interested reader to~\cite[Chapter~VIII, Section~6]{MR171116},
\cite[Sections~6--7 in Chapter~III and Section~6, Chapter~IV]{MR1963498}, and~\cite[Section~12.5]{MR3381284}
(however, the reader will catch sight of some of these methods in Exercises~\ref{POISSONKERN-ex1}, \ref{POISSONKERN-ex2},
and~\ref{POISSONKERN-ex3}, as well as in the applications given in Section~\ref{POISSONKERN-ex4}).

\begin{exercise}\label{CESA} The Ces\`aro mean of a sequence is the asymptotic limit of the sequence of arithmetic means.
That is, let~$\{\gamma_k\}_{k\in\N}$ be a sequence of real numbers and define its Ces\`aro mean the quantity
$$ \lim_{N\to+\infty}\frac1N\sum_{k=0}^{N-1}\gamma_k,$$
whenever this limit exists. \index{Ces\`aro mean}

Prove that if~$\gamma_k$ converges to some~$\gamma\in\R\cup\{+\infty\}\cup\{-\infty\}$ as~$k\to+\infty$, then
its Ces\`aro mean is also equal to~$\gamma$.

Exhibit also an example of positive sequence whose Ces\`aro mean is equal to zero but such that
$$ \sup_{k\in\N}\gamma_k=+\infty.$$
\end{exercise}

\begin{exercise} What is the Ces\`aro mean of the sequence~$\gamma_k:=(-1)^k$?
\end{exercise}

\begin{exercise}\label{NOINF-FeJ-th.3}
Proof that Theorem~\ref{FeJ-th.3} does not hold when~$p=+\infty$.
\end{exercise}

\begin{exercise}\label{IFEJ-FORM2.bn-EX}
Prove that, for all~$N\in\N$,
$$\int_{0}^{1} F_N(x)\,dx=1.$$
\end{exercise}

\begin{exercise}\label{GKFND}
Prove that
$$ \sup_{N\in\N}\int_0^1 |D_N(x)|\,dx=+\infty$$
but
$$ \sup_{N\in\N}\int_0^1 |F_N(x)|\,dx<+\infty.$$
\end{exercise}

\begin{exercise}\label{UTSJIFPOAX2e921rhnr} Prove that if~$f$ is bounded and periodic of period~$1$, then, for all~$x\in\R$,
$$\left|\frac1N\sum_{k=0}^{N-1}S_{k}(x)\right|\le\|f\|_{L^\infty(\R)}.$$\end{exercise}

\begin{exercise}\label{UTSJIFPOA} Let~$f\in L^1((0,1))$ be periodic of period~$1$.

Prove that $$\frac1N\sum_{k=0}^{N-1}S_{k}(x)=
\sum_{{k\in\Z}\atop{|k|\le N-1}}\left(1-\frac{|k|}{N}\right)\,\widehat f_k\,e^{2\pi ik x}.$$\end{exercise}

\begin{exercise}\label{UNIQ:ALAmowifj30othgore0-1} Give another proof of
Theorem~\ref{UNIQ} by using~\eqref{IFEJ-FORM.bi} and~\eqref{023wefv4567uygfdfgyhuizo0iwhfg0eb5627e01.p1rebi}.\end{exercise}

\begin{exercise}\label{UNIQ:ALAmowifj30othgore0-1ILBVIS} Let~$f\in L^1((0,1))$ be periodic of period~$1$ and
assume that
$$ \sum_{k\in\Z}|\widehat f_k|<+\infty.$$
Prove that~$\displaystyle\frac1N\sum_{k=0}^{N-1}S_{k}(x)$ converges uniformly for~$x\in\R$ as~$N\to+\infty$.
\end{exercise}

\begin{exercise}\label{FEJPOAE1edc4.10} Prove that there exists~$C>0$ such that, for all~$x\in\left[-\frac12,\frac12\right]$ and~$N\in\N\cap[1,+\infty)$,
$$ F_N(x)\le\frac{CN}{1+N^2x^2}.$$
\end{exercise}

\begin{exercise}\label{ojdkfnvioewyr098765rewe67890iuhgvgyuuygfr43wdfgh8bft5xsBIS} Provide an alternative solution to Exercise~\ref{ojdkfnvioewyr098765rewe67890iuhgvgyuuygfr43wdfgh8bft5xs} by using Theorem~\ref{FeJ-th.2}.
\end{exercise}

\begin{exercisesk}\label{PKSMDu0ojf65bv94.11} Let~$f$ be an even function,
bounded and periodic of period~$1$.
Suppose that the Fourier Series of~$f$ in trigonometric form is
$$ \frac{a_0}2+\sum_{k=1}^{+\infty} a_k\cos(2\pi kx),$$
with
\begin{equation}\label{PKSMDu0ojf65bv94.112}
{\mbox{$a_k\ge0\,$ for all~$\,k\in\N$.}}\end{equation}
Prove that there exists a constant~$C>0$ such that, for all~$x\in\R$ and~$N\in\N$,
$$ |S_N(x)|\le C\|f\|_{L^\infty(\R)}$$
and that~$\displaystyle\frac1N\sum_{k=0}^{N-1}S_{k}(x)$ converges uniformly for~$x\in\R$ as~$N\to+\infty$.
\end{exercisesk}

\begin{exercisesk}\label{PKSMDu0ojf65bv94.11124erf54y} Let~$f$ be an odd function, bounded and periodic of period~$1$.
Suppose that the Fourier Series of~$f$ in trigonometric form is
$$ \sum_{k=1}^{+\infty} b_k\sin(2\pi kx),$$
with
\begin{equation}\label{PKSMDu0ojf65bv94.112.1123}
{\mbox{$b_k\ge0\,$ for all~$\,k\in\N\cap[1,+\infty)$.}}\end{equation}
Prove that there exists a constant~$C>0$ such that, for all~$x\in\R$ and~$N\in\N$,
$$ |S_N(x)|\le C\|f\|_{L^\infty(\R)}$$
and that~$\displaystyle\frac1N\sum_{k=0}^{N-1}S_{k}(x)$ converges uniformly for~$x\in\R$ as~$N\to+\infty$.
\end{exercisesk}

\begin{exercisesk}\label{PKSMDu0ojf65bv94.11124erf54y.bi} In the setting of Exercises~\ref{PKSMDu0ojf65bv94.11} and~\ref{PKSMDu0ojf65bv94.11124erf54y}, suppose additionally that~$f$ is continuous
and prove that the Fourier Series of~$f$ converges to~$f$ uniformly in~$\R$.
\end{exercisesk}

\begin{exercise}\label{PKSMDu0ojf65bv94.11124erf54y.bi.PY} The so-called Paley's Theorem \index{Paley's Theorem}
states that if~$f$ is continuous, periodic of period~$1$, and the Fourier coefficients of~$f$ in trigonometric form
are all non-negative, then the Fourier Series of~$f$ converges uniformly to~$f$.
Prove this result.
\end{exercise}

\begin{exercisesk}\label{FEJPOAE1edc4.1}
Let~$f\in L^1((0,1))$ be periodic of period~$1$.

Then, for almost every~$x\in\R$,
$$ \lim_{N\to+\infty}
\frac1N\sum_{k=0}^{N-1}S_k(x)=f(x).$$
\end{exercisesk}

\begin{exercisesk}\label{FEJPO-010}
The following question was posed in the~1980 Mikl\'os Schweitzer Mathematical Contest (see~\cite[page~30 and~507]{MR1416130}).

Let~$f$ be a non-negative, integrable function on~$(0, 1)$ whose Fourier series in trigonometric form is
\begin{equation}\label{A0kmalseb.a}
\frac{\lambda_0}2 + \sum_{k=1}^{+\infty} \lambda_k \cos(2\nu_k\pi x),\end{equation}
for suitable~$\lambda_k\in\R$ and~$\nu_k\in\N\setminus\{0\}$.

Assume that none of the~$\nu_k$ divides another. 

Prove that, for all~$k\in\N\setminus\{0\}$,
\begin{equation}\label{A0kmalseb.a0-ln}|\lambda_k | \le \frac{\lambda_0}2.\end{equation}\end{exercisesk} 

\begin{exercise}\label{kpoqsld0O3.1P2P37jsd1}
Are the bounds obtained in Exercises~\ref{NONSepCFGDV} and~\ref{FEJPO-010} optimal?
\end{exercise}

\begin{exercisesk}\label{MSCOCOS}
Let~$\{\alpha_k\}_{k\in\N}$ and~$\{n_k\}_{k\in\N}$ be sequences, with~$\alpha_k\in(0,+\infty)$ and~$n_k\in\N\cap[2,+\infty)$,
such that\footnote{For example, one can take~$\alpha_k:=\frac1{k^2+1}$ and~$n_k:=2^{k^3+1}$ to satisfy all the conditions~\eqref{ZLAM-01}, \eqref{ZLAM-01a}, and~\eqref{ZLAM-02}.

As a notational remark, we stress that, as customary, $2^{k^3+1}$ has to be interpreted as~$2^{ ( k^3 +1) }$, here and in the rest of this \label{PERGBDSE0232345021} work.}
\begin{equation}\label{ZLAM-01}
\sum_{k=0}^{+\infty}\alpha_k<+\infty,
\end{equation}
\begin{equation}\label{ZLAM-01a}
n_{k+1}>n_k,
\end{equation}
and
\begin{equation}\label{ZLAM-02}
\limsup_{k\to+\infty} \alpha_k\,\ln n_k>0.
\end{equation}
Let
\begin{equation}\label{91wueihdncos2pimx.01oewjdflnv} Q(x,\ell):=\sum_{{j\in\Z\setminus\{0\}}\atop{|j|\le\ell}}\frac{\cos(2\pi(2\ell-j)x)}{j}\end{equation}
and consider the series
\begin{equation}\label{CSPPCSjKS} \sum_{k=1}^{+\infty} \alpha_k \,Q(x,n_k).\end{equation}

Does this series represent a continuous function~$f$ periodic of period~$1$?

And, if so, does the Fourier Series of~$f$ converge uniformly?
\end{exercisesk}

\begin{exercise}\label{POISSONKERN-ex1} Let~$f\in L^1((0,1))$ be periodic of period~$1$. Given~$r\in(0,1)$ and~$x\in\R$, the
\index{Abel mean} $r$-Abel mean of~$f$ at the point~$x$ is given by
$$ \sum_{{k\in\Z}} r^{|k|}\,\widehat f_k\,e^{2\pi ikx}.$$
Prove that this series\footnote{The method sketched in Exercise~\ref{POISSONKERN-ex1}
highlights the interesting possibility of somewhat relating Fourier and Taylor Series. \label{LAFOOTTF}
Let us stress however, that Fourier and Taylor Series are conceptually quite different. For instance,
Taylor Series require a much higher degree of smoothness of the given function and they can only
provide approximation in a small neighbourhood of a given point, while the Fourier Series can be written
under minimal assumptions on the function and can provide a global, not just local approximation.

However, nothing comes for free. The conceptual reason for which Fourier Series seem to have the above
advantages with respect to Taylor Series is that Fourier Series may be more ``expensive'':
while the computation of the Taylor coefficients only require the knowledge of a function in the vicinity of the expansion point,
the Fourier coefficients require the knowledge of the function on all its domain.

In the same vein, Taylor coefficients remain unchanged if the given function is modified away from the expansion point,
while Fourier coefficients are affected by any modification of the function.

See
Exercises~\ref{NICELI} and~\ref{NONSepCFGDV2}, as well as~\cite[pages~58--59]{MR1145236},
\cite{MR1435557},
and~\cite[pages~94--98]{MR3601106}, for more insight about similarities and differences regarding Fourier and Taylor Series.}
is convergent, and that the convergence is uniform in~$x$.
\end{exercise}

\begin{exercise}\label{POISSONKERN-ex2} The function
\begin{equation}\label{POISSONKERN-ex2-form1} [0,1)\times\R\ni (r,\theta)\longmapsto {\mathcal{P}}(r,\theta):=1+2\sum_{k=1}^{+\infty}r^k\cos(2\pi k\theta)\end{equation}
is called the \index{periodic Poisson Kernel} periodic Poisson Kernel.

Prove that, for all~$r_0\in(0,1)$, the series above is uniformly convergent
for~$(r,\theta)\in[0,r_0]\times\R$, that
\begin{equation}\label{POISSONKERN-ex2-form2}
{\mathcal{P}}(r,\theta)=\frac{1-r^2}{ (1-r)^2+4r\sin^2(\pi\theta) },
\end{equation}
and that, in the notation of Exercise~\ref{POISSONKERN-ex1}, the $r$-Abel mean of~$f$ can be written in the form
\begin{equation}\label{XIUJHGV9iuhg-2.1pjolsnENASMDP}
\Re\left(\widehat f_0+2\sum_{k=1}^{+\infty} r^k\,\widehat f_k\,e^{2\pi ikx}\right)
\end{equation}
and in the form
\begin{equation}\label{POISSONKERN-ex2-form3}
\int_0^1 f(x-y)\,{\mathcal{P}}(r,y)\,dy.
\end{equation}
\end{exercise}

\begin{exercise}\label{POISSONKERN-ex3} In the setting of Exercise~\ref{POISSONKERN-ex1}, suppose additionally that~$f$ is continuous and prove that the 
limit of the $r$-Abel mean of~$f$ converges to~$f$ everywhere in~$\R$ as~$r\nearrow1$.
\end{exercise}

\begin{exercise}\label{LPmPOIKE01} In the setting of Exercise~\ref{POISSONKERN-ex1}, suppose
that~$f\in L^p((0,1))$ for some~$p\in[1,+\infty)$
and prove that the 
limit of the $r$-Abel mean of~$f$ converges to~$f$ in~$L^p((0,1))$ as~$r\nearrow1$.
\end{exercise}

\begin{exercise}\label{LPmPOIKE01.002l42er21fygnb9mnvf} 
The function
\begin{equation}\label{POISSONKERN-ex2-form1-ILCO} [0,1)\times\R\ni (r,\theta)\longmapsto {\mathcal{Q}}(r,\theta):=2
\sum_{k=1}^{+\infty}r^k\sin(2\pi k\theta)\end{equation}
is called the \index{periodic Poisson Kernel, conjugate} conjugate periodic Poisson Kernel.

Prove that the series above is uniformly convergent and that
\begin{equation}\label{POISSONKERN-ex2-form2medl213mrpo21t}
{\mathcal{Q}}(r,\theta)=\frac{2r\sin(2\pi\theta)}{ (1-r)^2+4r\sin^2(\pi\theta)}.
\end{equation}

Moreover, in the notation of Exercise~\ref{POISSONKERN-ex2}, prove that \begin{equation}\label{LAMCIKCONJA0-2}
{\mathcal{P}}(r,\theta)=\Re\left(\frac{1+re^{2\pi i\theta}}{1-re^{2\pi i\theta}}\right)\qquad{\mbox{ and }}\qquad{\mathcal{Q}}(r,\theta)=\Im\left(\frac{1+re^{2\pi i\theta}}{1-re^{2\pi i\theta}}\right).\end{equation}
\end{exercise}

\begin{exercise}\label{POISSONKERN-exsadsdghj3Xs23954868.0-1} In the setting of Exercise~\ref{LPmPOIKE01.002l42er21fygnb9mnvf},
Prove that if~$f\in L^1((0,1))$ is periodic of period~$1$, then, for all~$x\in\R$,
\begin{equation}\label{lmq-21iDgbnmpeTArhmewxi} \int_0^1 f(y)\,{\mathcal{Q}}(r,x-y)\,dy=-i\sum_{k\in\Z}r^{|k|}\,{\operatorname{sign}}(k)\,\widehat f_{k}\,e^{2\pi ikx},\end{equation} where
\begin{equation}\label{lmq-21iDgbnmpeTArhmewxibis} {\operatorname{sign}}(k):=\begin{dcases}
1 &{\mbox{ if }}k>0,\\
-1 &{\mbox{ if }}k<0,\\
0 &{\mbox{ if }}k=0.
\end{dcases}\end{equation}
\end{exercise}

\begin{exercise}\label{LPmPOIKE01ijdqsp21rj90jtm} In the setting of Exercise~\ref{POISSONKERN-exsadsdghj3Xs23954868.0-1},
suppose
that~$f\in L^2((0,1))$ and we denote by~$f^\Diamond_r(x)$ the function in~\eqref{lmq-21iDgbnmpeTArhmewxi}.

Prove that~$f^\Diamond_r$ converges to some function~$f^\Diamond$ in~$L^2((0,1))$.

This function is called\footnote{Notice that, formally, by passing to the limit the expression in~\eqref{lmq-21iDgbnmpeTArhmewxi},
one could be tempted to use the suggestive expression \label{MKDFOOTOBECMSCoojwlf395utjhgnhbdfvgwherh5}
$$ f^\Diamond(x)=-i\sum_{k\in\Z} {\operatorname{sign}}(k)\,\widehat f_{k}\,e^{2\pi ikx}.$$
Compare this expression with the one in~\eqref{0pXe-244kf.1EX-PRE.pe-djenfEF}.}
the \index{Conjugate Function}
Conjugate Function of~$f$.
\end{exercise}

\section{Persistent overshooting phenomena}\label{SEC:GIBBS-PH}

This section deals with discontinuous functions.
A natural question is whether it is worth spending our time trying to understand such functions, or should we rather confine ourselves to functions that are as smooth as we want and spend more time at the beach.

This is a very deep dilemma, and for this we need to seek help from the Masters.

According to~\cite[page~479]{MR698947},
\begin{quote}
    ``physicists rightly believe that all functions (with, possibly, some exceptions) are analytic''.
\end{quote}

Also, a postulate attributed to Aristotle and often retaken by many scholars is that
``natura abhorret vacuum'' (that is, ``nature abhors a vacuum''), thus suggesting that discontinuities do not really model ``natural'' phenomena.

But according to~\cite[page~59]{MR4404761}, 
\begin{quote}
    ``not all the functions which occur in classical physics are continuous''.
\end{quote}

Well, now one is perhaps even more confused, but some comments about how wise it could be to include discontinuous functions in the analysis of an important phenomenon was proposed already on page~\pageref{GADIS}.
The bottom line is of course that the Masters are always right, and in mathematics, as much as in life, one has to do what one can (and, when one can, do what one likes).

That is, on the one hand, it is totally fine to try to squeeze the best of a theory under additional regularity assumptions,
since in many cases of interest these assumptions are met and then one possesses a strong theory that can be used as a nutcracker in plenty of situations. 

On the other hand, many complex problems are more suitably modelled by irregular functions,
since the model in itself can be regular but one is interested in obtaining quantitative information regardless of some of the parameters involved (e.g. the modulus of continuity of the data).

Moreover, the study of hypothetically irregular cases can finally allow one to establish that the solution found
is in fact regular (yet, to allow the solution to exist, very often it comes in very handy not to postulate its regularity to start with, so to gain the necessary compactness properties in the appropriate mathematical spaces).

To cut a long story short, we proceed with our plan and we discuss what happens to Fourier Series at the jump discontinuities of a given function (but readers more interested in spending time at the beach are free to skip this section, or to read this section at the beach).

Actually, the study of the convergence of Fourier Series at jump discontinuities arose from a concrete experiment, which we now briefly recall.
Albert Michelson retook the machines called harmonic analysers formerly invented by Lord Kelvin to predict tides (see Section~\ref{TIDES}): these devices were capable of adding together elementary oscillations to produce a complex graph
(a process called in the jargon \index{synthesis} ``synthesis''), and, conversely, to approximately reconstruct the Fourier modes from data.
One of the driving reasons to use these machines in Michelson's work is to study various physical phenomena, including the light emitted by flames. Michelson wanted to determine via  Fourier methods the different electromagnetic frequencies that compose it. 

\begin{figure}[h]
\includegraphics[height=2.8cm]{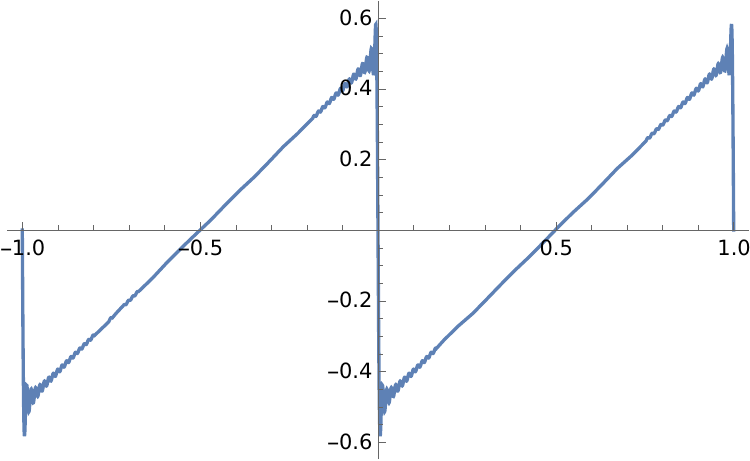}$\quad$\includegraphics[height=2.8cm]{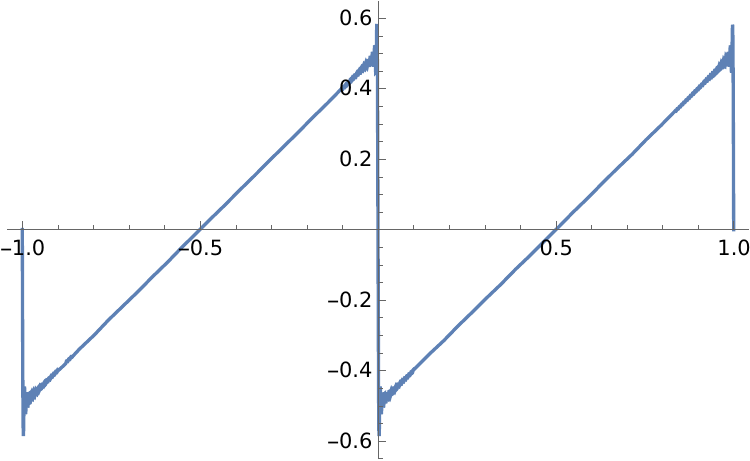}$\quad$\includegraphics[height=2.8cm]{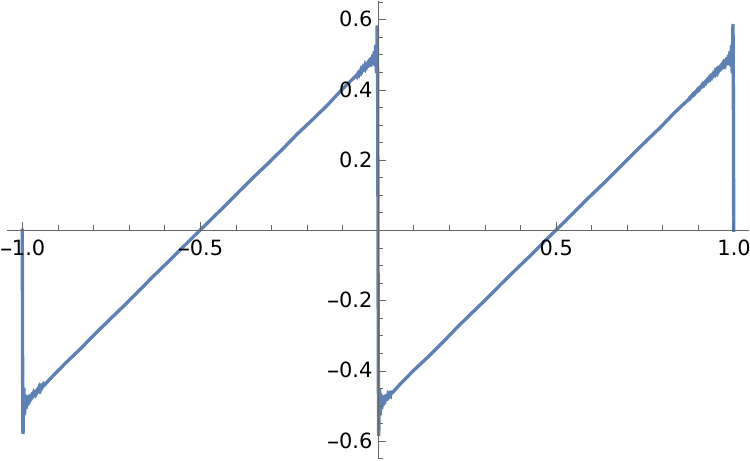}\\
\includegraphics[height=2.8cm]{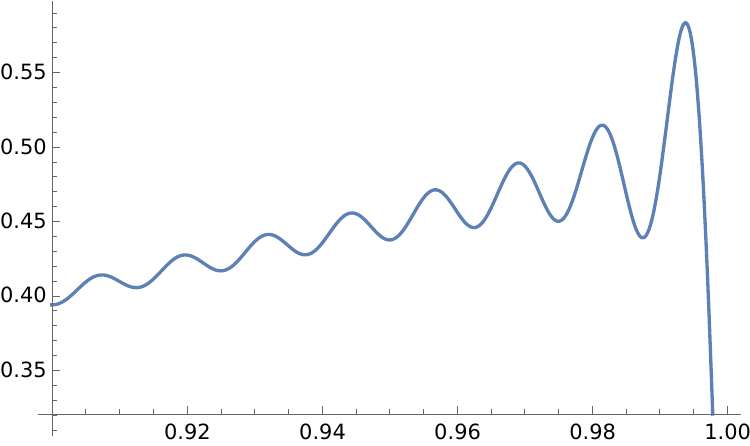}$\quad$\includegraphics[height=2.8cm]{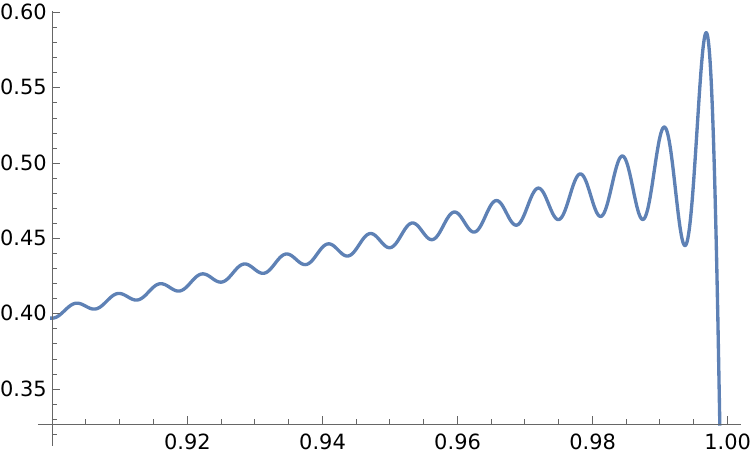}$\quad$\includegraphics[height=2.8cm]{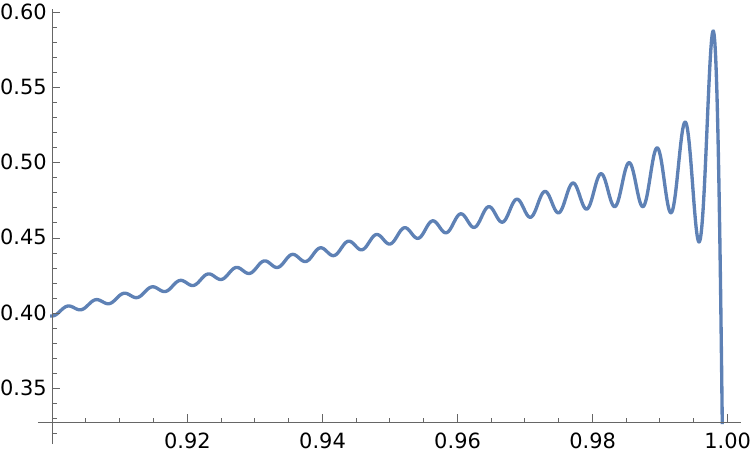}
\centering
\caption{Plot of a sawtooth wave approximation~$\displaystyle-\sum_{k=1}^{+\infty}\frac1{\pi k}\sin(2\pi kx)$, with~$k\in\{80,160,240\}$, and a zoom for~$x\in\left(\frac9{10},1\right)$.}\label{0FSSIS.K9iuygfy876trfdfghjjnn2meI202}
\end{figure}

This would have been a pivotal step in Michelson's work, eventually leading him to setting up the so-called \index{Michelson-Morley experiment} Michelson-Morley experiment to compare the speed of light in perpendicular directions aiming at detecting the relative motion of matter through the ``luminiferous aether'' (the negative result obtained in this experiment was in turn crucial to Einstein's invention of special relativity).

In his study of flames, Michelson performed at first the necessary Fourier analysis by hand, a task which however turned out to be rather laborious: in his own words,
\begin{quote}
``every one who has had occasion to calculate or to construct graphically the resultant of a large number of simple harmonic motions has felt the need of some simple and fairly accurate machine which would save the considerable time and labour involved in such computations''.  
\end{quote}
 The available machine for this at the time was the one introduced by Lord Kelvin (see Figure~\ref{o90pjHAhanaokdfR134ETPx.22} in Section~\ref{TIDES}) and relying on complex systems of strings and ropes around pulleys. For a very precise experimentalist such as Michelson, these ropes were an unforgivable imperfection, given their stretch and imperfect flexibility. The attempt of improving their performances by increasing the number of elements in the machine was also producing accumulated errors, thus preventing viable solutions.

To eliminate these experimental problems, Michelson proposed several solutions, including the addition of spiral springs, see~\cite{STRATT}. %
See also~\cite{MI2CHEBO}
(freely available to view and share for non-commercial purposes at

\url{https://engineerguy.com/fourier/pdfs/albert-michelsons-harmonic-analyzer.pdf}) for a very detailed and visual account
of Michelson's machines.

It is often reported (see e.g.~\cite[page~62]{MR4404761})
that, to test his machine with challenging, but somewhat explicit and verifiable cases, Michelson fed in the first 80 Fourier coefficients of the sawtooth waveform (see Exercise~\ref{SA:W})
and, rather surprisingly, the resulting plot seemed rather inaccurate at the discontinuity points, adding two little ``blips'' on either sides, and these blips seemed to be persistent even when all Michelson's technical experimental skills were put in place to reduce any possible source of errors: in fact, even after increasing the number of Fourier coefficients used in the synthesis, the blips would have just moved closer and closer to the jump of the wave but their size would have remained pretty much the same, with a significant overshoot of about 18\% the correct value (see Figures~\ref{09iuygfy876trfdfghjjnn2meI202}
and~\ref{0FSSIS.K9iuygfy876trfdfghjjnn2meI202} for a description of this phenomenon about the sawtooth wave).

\begin{figure}[h]
\includegraphics[height=5.8cm]{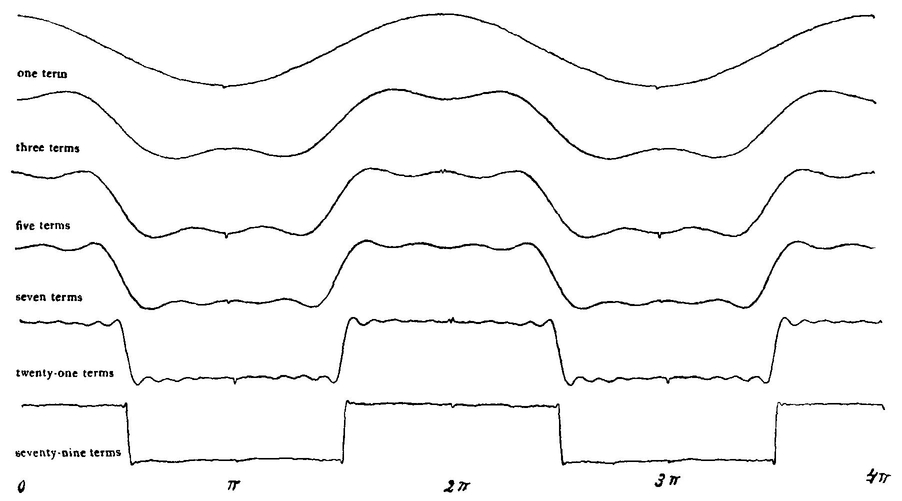}
\centering
\caption{Plot of a square wave approximation
according to~\cite{STRATT}
(to be compared with Figure~\ref{09iuygfy876trfdfghjjnn2meI20}).}\label{01FS.234.45fSIS.K9iuygfy876trfdfghjjnn2meI202}
\end{figure}

It is however controversial whether or not Michelson appreciated this overshoot in its full significance or still considered it as a byproduct of the manufacturing flaws of the physical device used for the synthesis.
It is also unclear whether Michelson's machine at the time was accurate enough to detect clearly and quantitatively the overshooting effect. At any rate, no mention of this overshooting is given in~\cite{STRATT} (though some images there did show some overshooting, see Figure~\ref{01FS.234.45fSIS.K9iuygfy876trfdfghjjnn2meI202}) and the link between Michelson and the jump discontinuity question only arose in a less-than-one-column letter to Nature~\cite{SMC12}, with the rather vague title of ``Fourier's Series'', in which no hint to the blip phenomenon is given, mostly discussing that ``the idea that a real discontinuity can replace the sum of continuous curve is so utterly at variance with the physicists' notions of quantity'' (so, apparently, at the time physicists were surprised that a continuous function, such as the sum of sines and cosines, could converge pointwise to a discontinuous function: fortunately, nowadays all physicists are fully at ease with the notions of pointwise, almost everywhere, and uniform convergence, as well as with the convergence in the more subtle functional spaces, given their ubiquitous importance in physical phenomena).

Michelson's brief digression did not remain unnoticed and Augustus Love replied to it in the very next issue of Nature~\cite{SMC13}, using the same unspecific title of the original letter. Love's article was also very short, again about one column, and the main focus of those pages was about marvellous natural phenomena (aurora borealis and lightning flashes during a thunderstorm). Love criticises Michelson's paper, considering the statements incorrect and the processes employed invalid, but the overshooting phenomenon of Fourier Series at jump discontinuities is not explicitly mentioned.

The dispute about convergence of Fourier Series and jump discontinuities was retaken by J. Willard Gibbs in another letter to Nature~\cite{GIBBSNA1} (appearing right after a very quick reply by Michelson to Love's observations). In this one-page note, Gibbs tries to distinguish between the ``limiting line'' (likely, the sawtooth waveform) and the ``limiting positions'' of the approximating Fourier Sum ``with fixed vertical transversals'' (arguably some kind of maximum limit).

The reader who does not find the statements in~\cite{GIBBSNA1} completely transparent
should appreciate that so did Gibbs, who published another half-a-column letter in Nature,
still under the same plain title of ``Fourier's Series'' ``to correct a careless error'' and\footnote{Given the lack of unambiguous terminology in the analysis carried out at that time, it may be difficult to point the finger to the specific error in~\cite{GIBBSNA1}. According to~\cite[page~648]{MR1491051}, the erroneous statement in~\cite{GIBBSNA1}
should be related to the fact that Gibbs ``seemed to imply that the oscillations decayed with~$N$'' (the core of the phenomenon is instead that the overshoot is persistent regardless of the number of the terms included in the Fourier sum).}
an ``unfortunate blunder'', see~\cite{GIBBSNA2}. In this additional note, Gibbs tries to make precise the distinction between the ``limiting form of the graphs of the functions expressed by the sum'' obtained via Fourier methods and ``the graph of the function expressed by the limit of that sum'' (arguably, the original discontinuous waveform). Interestingly, Gibbs complains about the lack of rigour in the notation used at time, mentioning that ``a misunderstanding on this point is a natural consequence of the usage which allows us to omit the word limit in certain connections'' which leads to the confusion between ``the limit of the graphs'' and ``the graph of the limit'', but no explicit formula to distinguish these concepts of limit is given. However, Gibbs correctly quantifies
the overshoot (see the last formula in display in~\cite{GIBBSNA2}, which corresponds to~\eqref{GIB:FO:CO} in our notation), but he gives no hint whatsoever about the proof of his claim.

\begin{figure}[h]
\includegraphics[height=11.7cm]{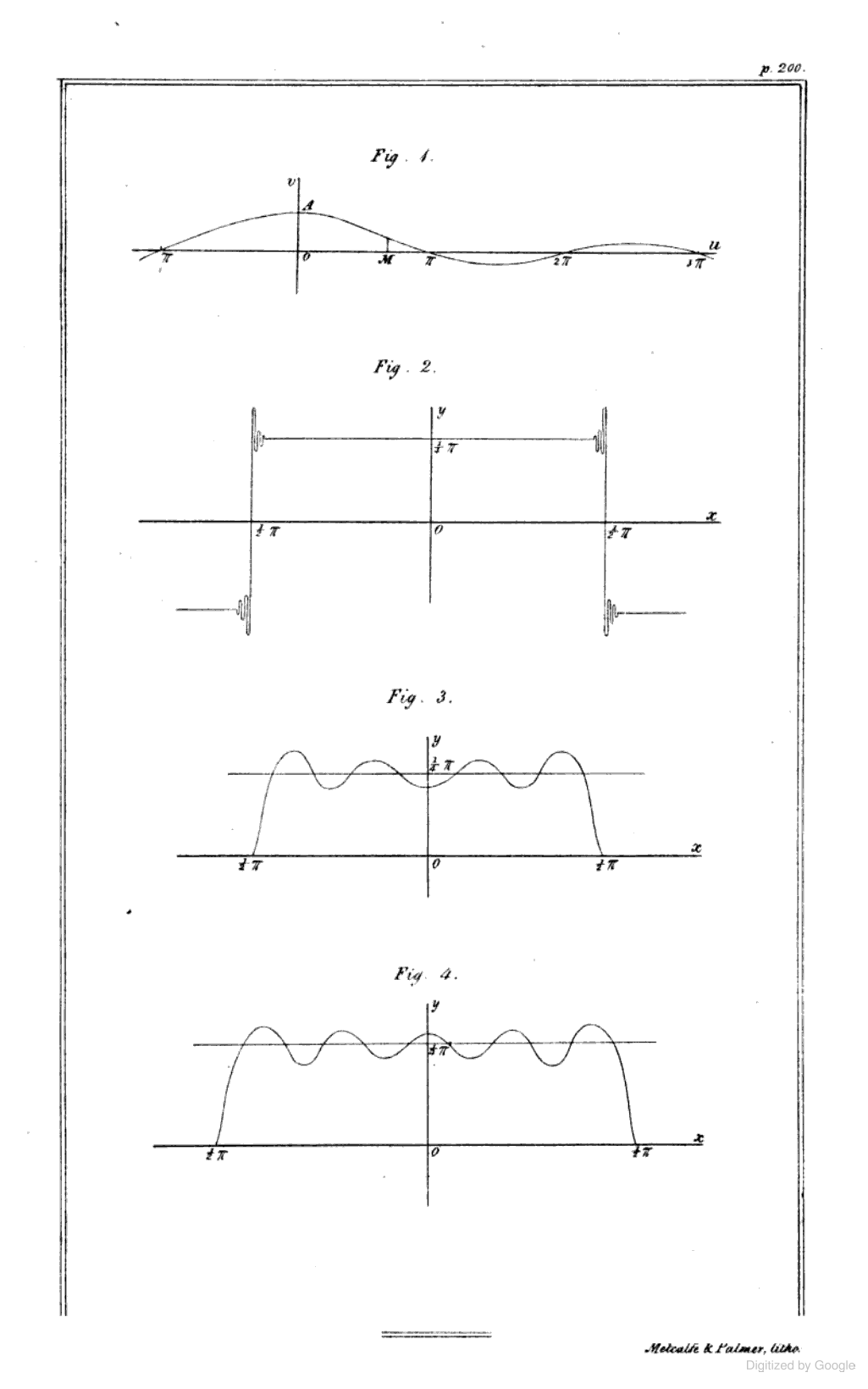}
\centering
\caption{A page of~\cite{Wilbraham}, depicting the approximation of a square wave
(to be compared with Figure~\ref{09iuygfy876trfdfghjjnn2meI20}).
The reader who wishes to complain for the lack of accuracy of this image should remember that
these graphs were computed manually.}
\label{n0FSSIS.K9iuygfy876trfdfghjjnn2meI202SQW}
\end{figure}

It appears however that the overshooting phenomenon (known nowadays as \index{Gibbs phenomenon} Gibbs phenomenon) was in fact first pointed out and analyzed in detail fifty years before by Henry Wilbraham~\cite{Wilbraham} (who also credited Francis Newman for part of the analysis, though Newman did not spot the overshoot). See figure~\ref{n0FSSIS.K9iuygfy876trfdfghjjnn2meI202SQW} for a page of Wilbraham's work, which had also correctly quantified the overshoot.

The rather unspecific title of Wilbraham's article, namely ``On a certain periodic function'', may not have helped popularising the specific result. Still and all, it is very likely that the other researchers (Michelson, Love, Gibbs, as well as nearly everyone else)
were simply unaware of the previous contribution by Wilbraham.

As the Scottish-Australian mathematician Horatio Carslaw pointed out~\cite{zbMATH02589958},
``we may still call this property of Fourier's series (and certain other series) Gibbs's phenomenon; but we must no longer claim that the property was first discovered by Gibbs''. Even so, the list of misnamed theorems and formulas in mathematics is lengthy and in fact, according to the so-called Boyer's Law~\cite{zbMATH04211975},
``mathematical formulas and theorems are usually not named after their original discoverers''
(or even more radically, by an adage attributed to Alfred Whitehead,
``everything of importance has been said before by somebody who did not discover it'').

But that's not too bad: after all, no mathematician really cares about who really invented calculus, or who really discovered some theorem. Mathematicians are fortunately antisocial, they don't care much about loud people fighting for priority, they study and create mathematics because they enjoy doing so. Mathematicians appreciate mathematics, not giving much importance to the names that are more or less artificially attached to it.

See~\cite{MR1502321} for a classical, very detailed, rigorous description of the persistent overshooting
and~\cite[pages~51--54]{MR3601106} for several interesting observations.
See also~\cite{zbMATH03660184} and the references therein for more information about the history of the Gibbs's phenomenon and several technically interesting observations.

We now dive into some mathematical details about Gibbs's phenomenon. Without aiming at the broadest possible generality,
we present the following result here:

\begin{theorem}\label{IGIDE}
Let~$f\in L^1((0,1))$ be periodic of period~$1$. Let~$x_0\in\R$ and suppose that there exists~$\mu\in\left(0,\frac12\right)$
such that~$f$ is of class~$C^1$ in the intervals~$(x_0-\mu,x_0)$ and~$(x_0,x_0+\mu)$, with
\begin{equation} \label{ILBA}\sup_{x\in(x_0-\mu,x_0)}|f'(x)|+\sup_{x\in(x_0,x_0+\mu)}|f'(x)|<+\infty.\end{equation}

Assume also that the following limits exist and are finite
$$ \ell_-:=\lim_{x\nearrow x_0}f(x),\qquad L_-:=\lim_{x\nearrow x_0}f'(x),\qquad
\ell_+:=\lim_{x\searrow x_0}f(x)
\qquad{\mbox{and}}\qquad L_+:=\lim_{x\searrow x_0}f'(x),$$
with~$\ell_+\ge\ell_-$.

Then, for all~$\delta\in(0,\mu)$,
\begin{equation}\label{GIB-a1} \lim_{N\to+\infty}\sup_{x\in\left(x_0-\delta,x_0\right)\cup\left(x_0,x_0+\delta\right)}\big(S_{N,f}(x)-f(x)\big)=
\frac{\ell_+-\ell_-}\pi\int_0^{\pi } \frac{\sin \tau }{\tau}\,d\tau-\frac{\ell_+-\ell_-}2.\end{equation}

Moreover,
\begin{equation}\label{GIB-a1.las} \lim_{N\to+\infty}
S_{N,f}\left(x_0+\frac1{2N}\right)=\frac{\ell_++\ell_-}2+ \frac{\ell_+-\ell_-}\pi\int_0^{\pi } \frac{\sin \tau }{\tau}\,d\tau.\end{equation}
\end{theorem}

We will prove this result by reducing to a ``model case'', namely\footnote{The technical advantage of
reducing to a specific situation is that this case can have simpler algebraic properties. For instance, by Exercise~\ref{FO:DE:MAGIB}, the derivative
of the Fourier Sum of the square wave can be reduced to a geometric sum which is explicitly computable. This comes in handy to detect the critical points of the corresponding Fourier sum.
Also, the Fourier Sum itself can be nicely reconstructed from the derivative, via
the Fundamental Theorem of Calculus.
The proof of Lemma~\ref{MODELEM} will leverage these observations.} the one of the square wave:

\begin{lemma} \label{MODELEM}
Let~$w$ be the square wave in Exercise~\ref{SQ:W}.

Then, for each~$\delta\in\left(0,\frac12\right)$,
\begin{equation}\label{GIB:FO:CO} \lim_{N\to+\infty}\sup_{x\in(-\delta,0)\cup(0,\delta)}\big( S_{N,w}(x)-w(x)\big)=
\frac2\pi\int_0^{\pi } \frac{\sin \tau }{\tau}\,d\tau-1.\end{equation}

Moreover,
\begin{equation}\label{GIB:FO:CO.002} \lim_{N\to+\infty} S_{N,w}\left(\frac1{2N}\right)=
\frac2\pi\int_0^{\pi } \frac{\sin \tau }{\tau}\,d\tau.\end{equation}
\end{lemma}

For concreteness, let us mention\footnote{Interestingly, the value of the integral in~\eqref{GIB:FO:CO} has been oftentimes miscalculated, also by the
Masters of Mathematics, see e.g.~\cite[pages~153 and~158]{zbMATH03660184}:
as Latin poet Horace mentioned in his
``Ars Poetica'', ``quandoque bonus dormitat Homerus''. Maybe, we have also miscalculated here, in our goofy attempt to imitate the Masters.}
 that the value of the right-hand side of~\eqref{GIB:FO:CO}
is approximately~$0.17898$. That is, given that the jump of the square wave in Exercise~\ref{SQ:W} is equal to~$2$ (and its maximal value equal to~$1$), the overshooting in~\eqref{GIB:FO:CO} is about~$9\%$ of the size of the jump (hence, not negligible in concrete applications).

For the sake of clarity, let us also stress that,
despite the Gibbs phenomenon, the Fourier Series of the square wave
converges pointwise to the original function except at the discontinuity point,
and in fact the convergence is uniform in all sets of the form
$$ \R\setminus\left(\bigcup_{j\in\Z}\left(\frac{j}2-a,\frac{j}2+a\right)\right),$$
for each~$a>0$, owing to Theorem~\ref{C1uni}.

The gist for it is that the overshoot does maintain its vertical size, but it gets closer and closer to the discontinuities when one incorporates more terms in the Fourier Sum.

\begin{proof}[Proof of Lemma~\ref{MODELEM}] By Exercise~\ref{ojld03-12d},
$$ S_{N,w}(x)=\sum_{j=0}^{N} \frac4{\pi(2j+1)}\sin(2\pi(2j+1)x).$$
Hence, for all~$x\in\left(0,\frac12\right)$,
$$ F_N(x):=S_{N,w}(x)-w(x)=\sum_{j=0}^{N} \frac4{\pi(2j+1)}\sin(2\pi(2j+1)x)-1.$$

We now use Exercise~\ref{FO:DE:MAGIB}, and specifically~\eqref{EFDEN}, to see that
$$ F_N'(x)=\frac{4\sin(4\pi (\widetilde N+1)x)}{\sin(2\pi x)},$$
with
$$ \widetilde N:=\begin{dcases}\displaystyle\frac{N-1}2 &{\mbox{ if $N$ is odd,}}\\
\displaystyle\frac{N-2}2 &{\mbox{ if $N$ is even,}}
\end{dcases}$$
therefore the critical points of~$F_N$ are of the form~$x_m=\frac{m}{4(\widetilde N+1)}$, with~$m\in\N\setminus\{0\}$,
and notice that~$x_m\in(0,\delta)$ as long as~$\widetilde N>\frac{m}{4\delta}-1$.

We have that\footnote{It is useful to observe that the supremum in~\eqref{kijnmlp0-9IKmGB.12343}
will not be attained for~$x\searrow0$, because
$$ \lim_{x\searrow0} S_{N,w}(x)-w(x) =0-1=-1,$$
nor for~$x\nearrow\delta\in\left(0,\frac12\right)$, since, by the uniform convergence result in Theorem~\ref{C1uni},
for all~$\delta'\in\left(\delta,\frac12\right)$,
$$ \lim_{N\to+\infty}\sup_{x\in(\delta/2,\delta')}\big|S_{N,w}(x)-w(x)\big|=0.$$}
\begin{equation}\label{kijnmlp0-9IKmGB.12343}\begin{split}\lim_{N\to+\infty}\sup_{x\in(0,\delta)}\big(S_{N,w}(x)-w(x)\big)&=
\lim_{N\to+\infty}\sup_{x\in(0,\delta)}F_N(x)\\&=
\lim_{N\to+\infty}\sup_{{m\in\N}\atop{0<m<4\delta (\widetilde N+1)}}F_N(x_m)\\&=
\lim_{N\to+\infty}\sup_{{m\in\N}\atop{0<m<4\delta (\widetilde N+1)}}S_{N,w}(x_m)-w(x_m)\\&=
\lim_{N\to+\infty}\sup_{{m\in\N}\atop{0<m<4\delta (\widetilde N+1)}}S_{N,w}(x_m)-1.
\end{split}\end{equation}

Now we observe that, for all~$x\in\left(0,\frac12\right)$,
\begin{equation*}\begin{split}&
S_{N,w}(x)=F_N(x)+w(x)=F_N(x)+1=F_N(x)-F_N(0^+)\\&\qquad=
\int_0^x F'_N(y)\,dy=\int_0^x \frac{4\sin(4\pi (\widetilde N+1)y)}{\sin(2\pi y)}\,dy\end{split}
\end{equation*}
and thus, using the substitution~$\tau:=4\pi (\widetilde N+1)y$,
\begin{equation}\label{msmNCSNqdk.1}
S_{N,w}(x)=\int_0^{4\pi(\widetilde N+1)x} \frac{\sin \tau }{
\pi(\widetilde N+1)\sin\left(\frac{\tau}{2(\widetilde N+1)}\right)}\,d\tau.
\end{equation}

This, combined with~\eqref{kijnmlp0-9IKmGB.12343}, returns that
\begin{equation*}\begin{split}&\!\!\!\!\!\lim_{N\to+\infty}\sup_{x\in(0,\delta)}\big(S_{N,w}(x)-w(x)\big)\\&=
\lim_{\widetilde N\to+\infty}\sup_{{m\in\N}\atop{0<m<4\delta (\widetilde N+1)}} \int_0^{4\pi(\widetilde N+1)x_m} \frac{\sin \tau }{
\pi(\widetilde N+1)\sin\left(\frac{\tau}{2(\widetilde N+1)}\right)}\,d\tau-1\\&=
\lim_{\widetilde N\to+\infty}\sup_{{m\in\N}\atop{0<m<4\delta (\widetilde N+1)}} \int_0^{\pi m} \frac{\sin \tau }{
\pi(\widetilde N+1)\sin\left(\frac{\tau}{2(\widetilde N+1)}\right)}\,d\tau-1.
\end{split}\end{equation*}

On this account (see Exercise~\ref{ALLSFAMNLS}) we have that
\begin{equation*}\begin{split}\lim_{N\to+\infty}\sup_{x\in(0,\delta)}\big(S_{N,w}(x)-w(x)\big)&=
\frac2\pi\lim_{\widetilde N\to+\infty}\sup_{{m\in\N}\atop{0<m<4\delta (\widetilde N+1)}} \int_0^{\pi m} \frac{\sin \tau }{\tau}\,d\tau -1\\&=
\frac2\pi\sup_{{m\in\N\setminus\{0\}}} \int_0^{\pi m} \frac{\sin \tau }{\tau}\,d\tau-1.
\end{split}\end{equation*}
The desired result in~\eqref{GIB:FO:CO} thus follows (see Exercise~\ref{CO213421FSBCA.02-ejNS12c2g2n4},
and notice that the functions under consideration are odd symmetric).

To prove~\eqref{GIB:FO:CO.002}, we retake~\eqref{msmNCSNqdk.1} to see that
\begin{equation}\label{MSasdX-kdmcS.a}
S_{N,w}\left(\frac1{4(\widetilde N+1)}\right)=
\int_0^{\pi} \frac{\sin \tau }{
\pi(\widetilde N+1)\sin\left(\frac{\tau}{2(\widetilde N+1)}\right)}\,d\tau\\
=\int_0^{\pi} \frac{\sin \tau }{\tau}\cdot\frac{1}{
\frac{\pi(\widetilde N+1)}\tau\;\sin\left(\frac{\tau}{2(\widetilde N+1)}\right)}\,d\tau.
\end{equation}

Additionally, since~$4(\widetilde N+1)\in\{ 2N,\,2(N+1)\}$, we find that
\begin{eqnarray*}&&
\limsup_{N\to+\infty}\left|S_{N,w}\left(\frac1{4(\widetilde N+1)}\right)-S_{N,w}\left(\frac1{2N}\right)\right|\le
\limsup_{N\to+\infty}\left|S_{N,w}\left(\frac1{2(N+1)}\right)-S_{N,w}\left(\frac1{2N}\right)\right|\\
&&\qquad=
\limsup_{N\to+\infty}\left|
\sum_{j=0}^{N} \frac4{\pi(2j+1)}\left(\sin\left(\frac{\pi(2j+1)}{N+1}\right)-\sin\left(\frac{\pi(2j+1)}{N}\right)\right)
\right|\\&&\qquad\le4\limsup_{N\to+\infty}
\sum_{j=0}^{N} \left(\frac{1}{N}-\frac{1}{N+1}\right)=\limsup_{N\to+\infty}\frac4N=0.
\end{eqnarray*}
Combining this and~\eqref{MSasdX-kdmcS.a} we gather that
\begin{eqnarray*}&&\limsup_{N\to+\infty}\left| S_{N,w}\left(\frac1{2N}\right)-
\frac2\pi\int_0^{\pi } \frac{\sin \tau }{\tau}\,d\tau\right|\le
\limsup_{N\to+\infty}\left| S_{N,w}\left(\frac1{4(\widetilde N+1)}\right)-
\frac2\pi\int_0^{\pi } \frac{\sin \tau }{\tau}\,d\tau\right|\\&&\qquad=\limsup_{\widetilde N\to+\infty}\left| \int_0^{\pi} \frac{\sin \tau }{\tau}\cdot\frac{1}{
\frac{\pi(\widetilde N+1)}\tau\;\sin\left(\frac{\tau}{2(\widetilde N+1)}\right)}\,d\tau-
\frac2\pi\int_0^{\pi } \frac{\sin \tau }{\tau}\,d\tau\right|
\end{eqnarray*}

From this, taking notice of the positivity intervals of the sine function and using the Dominated Convergence Theorem, we obtain~\eqref{GIB:FO:CO.002}, as desired.
\end{proof}

\begin{proof}[Proof of Theorem~\ref{IGIDE}] Up to a horizontal translation, we can assume that~$x_0:=0$. Let also~$\delta\in(0,\mu)$, $\mu'\in(\delta,\mu)$, and~$\xi\in C^\infty_0\left( (-\mu,\mu)\right)$ with~$\xi(x)=1$ for all~$x\in\left[-\mu',\mu'\right]$ and define
$$g(x):=f(x)\xi(x).$$
With a slight abuse of notation, we identify~$g$ with its periodic extension of period~$1$ outside~$\left[-\frac12,\frac12\right)$.

We observe that
\begin{equation}\label{AMs-RLMAP-0}{\mbox{$f=g$ over~$\left(-\mu',\mu'\right)$}}\end{equation} and accordingly, by the Riemann Localisation Principle in
Theorem~\ref{C1uni-c}, 
\begin{equation}\label{AMs-RLMAP-1}
\lim_{N\to+\infty}\sup_{x\in\left(-\delta,\delta\right)}| S_{N,f}(x)-S_{N,g}(x)|=0.
\end{equation}

Now we take a globally Lipschitz function~$\varphi\in C^1\left(
\left(-\frac12,0\right)\cup\left(0,\frac12\right)
\right)$, periodic of period~$1$ and such that
\begin{equation}\begin{split}\label{AMs-RLMAP-1.192eihdfn}
&\sup_{x\in\left(-\frac12,0\right)}|\varphi'(x)|+\sup_{x\in\left(0,\frac12\right)}|
\varphi'(x)|<+\infty,\qquad\varphi(0)=0,\\ &\lim_{x\searrow0}\varphi'(x)=L_+,\qquad
{\mbox{and}}\qquad\lim_{x\nearrow0}\varphi'(x)=L_-.\end{split}\end{equation}
We define
$$ h:=g-\varphi$$
and we remark that
$$ \lim_{x\nearrow 0}h(x)=\ell_-\,,
\qquad\lim_{x\searrow 0}h(x)=\ell_+\,,
\qquad{\mbox{and}}\qquad
\lim_{x\nearrow 0}h'(x)=
\lim_{x\searrow 0}h'(x)=0.$$

Also, by Theorem~\ref{C1uni} (see also Exercises~\ref{DINI1UNI} and~\ref{DINI2UNI}),
\begin{equation}\label{AMs-RLMAP-1.192eihdfn2} \lim_{N\to+\infty}\sup_{x\in\R}| \varphi(x)-S_{N,\varphi}(x)|=0\end{equation}
and, as a result,
\begin{equation}\label{ojls.d5fgokmerft56hZX1} \begin{split}&
\lim_{N\to+\infty}\sup_{x\in\left(-\delta,\delta\right)}
\left||g(x)-S_{N,g}(x)|-|h(x)-S_{N,h}(x)|\right|\\&\qquad\le
\lim_{N\to+\infty}\sup_{x\in\left(-\delta,\delta\right)}
|g(x)-h(x)-S_{N,g}(x)+S_{N,h}(x)|\\&\qquad=\lim_{N\to+\infty}\sup_{x\in\left(-\delta,\delta\right)}
|\varphi(x)-S_{N,\varphi}(x)|\\&\qquad=0.
\end{split}\end{equation}

Now, let~$w$ be the square wave in Exercise~\ref{SQ:W} and define
$$ W(x):=h(x)-\frac{\ell_++\ell_-}2-\frac{\ell_+-\ell_-}2w(x).$$
In this way,
\begin{equation}\label{ojls.d5fgokmerft56hZX1.owfnIAN} \lim_{x\nearrow0}W(x)=\lim_{x\searrow0}W(x)=0\qquad{\mbox{and}}\qquad
\lim_{x\nearrow0}W'(x)=\lim_{x\searrow0}W'(x)=0.\end{equation}

Now we check that~$W$ satisfies a suitable uniform version of the Dini's Condition
(see~\eqref{DINI-UNO}), namely
that for every~$\epsilon>0$ there exists~$\delta_\epsilon>0$ such that
\begin{equation}\label{OJHSD.hy769ijhk76GD98ewygievn7b7598unv6pojhgM10}
\sup_{x\in (-\delta,\delta)}\int_{-\delta_\epsilon}^{\delta_\epsilon}\left|\frac{W (x + t) - W(x)}{t}\right|\,dt\le \epsilon.
\end{equation}
For this purpose, we first remark that~$W$ is of class~$C^1$ in~$(-\delta,0)\cup(0,\delta)$, with
$$ W'=h'=g'-\varphi'=f'\xi+f\xi'-\varphi'$$
and therefore, in virtue of~\eqref{ILBA}
and~\eqref{AMs-RLMAP-1.192eihdfn},
\begin{equation}\label{OJHSD.hy769ijhk76GD98ewygievn7b7598unv6pojhgM10.2p} \sup_{x\in(-\mu',0)}|W'(x)|+\sup_{x\in(0,\mu')}|W'(x)|\leq C,\end{equation}
for some~$C>0$.

As a consequence, when~$|t|<\mu'-\delta$,
\begin{equation}\label{OJHSD.hy769ijhk76GD98ewygievn7b7598unv6pojhgM10.2}
\sup_{x\in (-\delta,\delta)}|W(x+t)-W(x)|\le C|t|.
\end{equation}
Indeed, if the origin does not lie between~$x+t$ and~$x$, then~\eqref{OJHSD.hy769ijhk76GD98ewygievn7b7598unv6pojhgM10.2} follows from~\eqref{OJHSD.hy769ijhk76GD98ewygievn7b7598unv6pojhgM10.2p}; if instead
the origin lies between~$x+t$ and~$x$, then we denote by~$\alpha\ge0$ the distance
of~$x+t$ to the origin and by~$\beta\ge0$ the distance
of~$x$ to the origin (so that~$\alpha+\beta=|t|$) and we combine~\eqref{OJHSD.hy769ijhk76GD98ewygievn7b7598unv6pojhgM10.2p} with the continuity of~$W$ at the origin, to see that
\begin{eqnarray*}&& \sup_{x\in (-\delta,\delta)}|W(x+t)-W(x)|\le
\sup_{x\in (-\delta,\delta)}\big( |W(x+t)-W(0)|+|W(0)-W(x)|\big)\\&&\qquad
\le \sup_{x\in (-\delta,\delta)}\left(\sup_{y\in(-\delta,0)}|W'(y)|\,\alpha+\sup_{y\in(0,\delta)}|W'(y)|\,\beta
\right)\le C\alpha+C\beta= C|t|,\end{eqnarray*}
thus establishing~\eqref{OJHSD.hy769ijhk76GD98ewygievn7b7598unv6pojhgM10.2}.

Hence, by means of~\eqref{OJHSD.hy769ijhk76GD98ewygievn7b7598unv6pojhgM10.2},
we obtain the desired result in~\eqref{OJHSD.hy769ijhk76GD98ewygievn7b7598unv6pojhgM10}
by choosing~$\delta_\epsilon:=\min\left\{\frac\epsilon{2C},\,\mu'-\delta\right\}$.

Now, thanks to~\eqref{OJHSD.hy769ijhk76GD98ewygievn7b7598unv6pojhgM10}, we can utilize
Theorem~\ref{C1uni} and deduce that the Fourier Series of~$W$ converges to~$W$ uniformly in~$(-\delta,\delta)$, namely
\begin{equation}\label{ojls.d5fgokmerft56hZX1.owfnIAN2} \lim_{N\to+\infty}\sup_{x\in(-\delta,\delta)}|W(x)-S_{N,W}(x)|=0\end{equation}
and therefore
\begin{eqnarray*}&&
\lim_{N\to+\infty}\sup_{x\in\left(-\delta,\delta\right)}
\left||h(x)-S_{N,h}(x)|-\frac{\ell_+-\ell_-}2|w(x)-S_{N,w}(x)|\right|
\\&&\qquad\le\lim_{N\to+\infty}\sup_{x\in\left(-\delta,\delta\right)}
\left|h(x)-S_{N,h}(x)-\frac{\ell_+-\ell_-}2\big(w(x)-S_{N,w}(x)\big)\right|\\&&\qquad\le\lim_{N\to+\infty}\sup_{x\in\left(-\delta,\delta\right)}|W(x)-S_{N,W}(x)|\\&&\qquad=0,
\end{eqnarray*}
which combined with~\eqref{ojls.d5fgokmerft56hZX1} entails that
$$ \lim_{N\to+\infty}\sup_{x\in\left(-\delta,\delta\right)}
\left|\frac{\ell_+-\ell_-}2|w(x)-S_{N,w}(x)|-|g(x)-S_{N,g}(x)|\right|=0.$$

Hence, by~\eqref{AMs-RLMAP-0} and~\eqref{AMs-RLMAP-1},
\begin{eqnarray*}
&& \lim_{N\to+\infty}\sup_{x\in\left(-\delta,\delta\right)}
\left|\frac{\ell_+-\ell_-}2
\big(S_{N,w}(x)-w(x)\big)
-\big(S_{N,f}(x)-f(x)\big)\right|\\
&&\qquad=\lim_{N\to+\infty}\sup_{x\in\left(-\delta,\delta\right)}
\left|\frac{\ell_+-\ell_-}2\big(S_{N,w}(x)-w(x)\big)-\big(S_{N,g}(x)-g(x)\big)\right|\\&&\qquad=0.
\end{eqnarray*}
This and~\eqref{GIB:FO:CO} yield the desired result in~\eqref{GIB-a1}.

Furthermore, since
$$g=h+\varphi=W+\frac{\ell_++\ell_-}2+ \frac{\ell_+-\ell_-}2 w+\varphi,$$
we have that
\begin{eqnarray*}
S_{N,g}=S_{N,W}+\frac{\ell_++\ell_-}2+ \frac{\ell_+-\ell_-}2 S_{N,w}+S_{N,\varphi}
\end{eqnarray*}
and accordingly
\begin{equation}\label{AMs-RLMAP-1cv}\lim_{N\to+\infty}
S_{N,g}\left(\frac1{2N}\right)=\lim_{N\to+\infty}
\left(S_{N,W}\left(\frac1{2N}\right)+\frac{\ell_++\ell_-}2+ \frac{\ell_+-\ell_-}2 S_{N,w}\left(\frac1{2N}\right)+S_{N,\varphi}\left(\frac1{2N}\right)\right).
\end{equation}

Besides, by means of~\eqref{ojls.d5fgokmerft56hZX1.owfnIAN} and~\eqref{ojls.d5fgokmerft56hZX1.owfnIAN2},
\begin{equation} \label{AMs-RLMAP-1cz}\limsup_{N\to+\infty}\left|S_{N,W}\left(\frac1{2N}\right)\right|\le
\limsup_{N\to+\infty}\left( \sup_{x\in{(-\delta,\delta)}}\left|S_{N,W}(x)-W(x)\right|+
\left|W\left(\frac1{2N}\right)\right|\right)=0.\end{equation}

On a similar note, in virtue of~\eqref{AMs-RLMAP-1.192eihdfn2}
(and recalling that~$\varphi(0)=0$, thanks to~\eqref{AMs-RLMAP-1.192eihdfn}),
\begin{equation*} \limsup_{N\to+\infty}\left|S_{N,\varphi}\left(\frac1{2N}\right)\right|\le
\limsup_{N\to+\infty}\left( \sup_{x\in\R}\left|S_{N,\varphi}(x)-\varphi(x)\right|+
\left|\varphi\left(\frac1{2N}\right)\right|\right)=0.\end{equation*}

Using this information, \eqref{AMs-RLMAP-1cv}, and~\eqref{AMs-RLMAP-1cz}, we gather that
\begin{equation*}\lim_{N\to+\infty}
S_{N,g}\left(\frac1{2N}\right)=\lim_{N\to+\infty}
\left(\frac{\ell_++\ell_-}2+ \frac{\ell_+-\ell_-}2 S_{N,w}\left(\frac1{2N}\right)\right).
\end{equation*}
For this reason and~\eqref{GIB:FO:CO.002} we conclude that
\begin{equation*}\lim_{N\to+\infty}
S_{N,g}\left(\frac1{2N}\right)=\frac{\ell_++\ell_-}2+ \frac{\ell_+-\ell_-}\pi\int_0^{\pi } \frac{\sin \tau }{\tau}\,d\tau.
\end{equation*}
{F}rom this and~\eqref{AMs-RLMAP-1} we obtain~\eqref{GIB-a1.las}, as desired.
\end{proof}

See e.g.~\cite[Section~1.6]{MR442564}, \cite[Chapter~II, Section~9 and Chapter~III, Section~11]{MR1963498}, \cite[pages~87--88]{MR2742531},
\cite[Sections~16 and~17]{MR4404761}, \cite[pages~181--182]{MR545506}, 
as well as the references therein,
for more details on the convergence properties of Fourier Series in the presence of jump discontinuities.
See also~\cite{MR1491051, MR1650415} for a general framework for the Gibbs phenomenon.

Several practical consequences of the Gibbs phenomenon will be presented in Section~\ref{GI:CONSECT}.

\begin{exercise}\label{9k:INTEHigii.0}
Let~$\varphi\in C^\infty_0((-1,1))$ with~$\varphi(0)=1$.
Let also~$ \varphi_N(x):=\varphi(Nx-1)$. 

Prove that, for all~$x\in\R$,
\begin{equation}\label{NOAKsm0.1mef1} \lim_{N\to+\infty}\varphi_N(x)=0,\end{equation}
but not uniformly.
\end{exercise}

\begin{exercise}\label{ALLSFAMNLS} Let~$\delta\in\left(0,\frac12\right)$.
Prove that
\begin{equation*}
\lim_{N\to+\infty}\,\sup_{0<m<2\delta(N+1)}
\left|\int_0^{\pi m} \frac{\sin \tau }{\pi(N+1)\sin\left(\frac{\tau}{2(N+1)}\right)}\,d\tau
-\frac2\pi\int_0^{\pi m} \frac{\sin \tau }{\tau}\,d\tau\right|=0.
\end{equation*}
\end{exercise}

\begin{exercise}\label{9k:INTEHigii} Let~$f\in L^1((0,1))$ be periodic of period~$1$. Let also~$p\in\R$.
The Gibbs Set \index{Gibbs Set} of~$f$ at~$p$ is the collection of values that can be attained by the Fourier Series of~$f$
through appropriate sub-sequences approaching the point~$p$.

More explicitly, $c\in\R$ belongs to the Gibbs Set of~$f$ at~$p$ if and only if there exists a sequence
of points~$p_N$ with~$p_N\to p$ as~$N\to+\infty$ such that
$$ \lim_{N\to+\infty} S_{N}(p_N)=c.$$

Prove that the Gibbs Set is an interval.
\end{exercise}

\begin{exercise} \label{KMSX:0oikjnhG-1sq0wuohgi}
Assume the setting of Theorem~\ref{IGIDE}.
Prove that the Gibbs Set (as defined in Exercise~\ref{9k:INTEHigii}) of~$f$ at~$x_0$ is the interval
$$ \big[\ell_- - \lambda(\ell_+-\ell_-),\,\ell_+ +\lambda(\ell_+-\ell_-)\big],$$
where
$$\lambda:=\frac1\pi\int_0^{\pi } \frac{\sin \tau }{\tau}\,d\tau-\frac{1}{2}\simeq 0.08948987.$$
\end{exercise}

\begin{exercise} \label{0-o2epkdjfm.019k4-9k:INTEHigii}
Convince yourself that the Gibbs phenomenon described in Theorem~\ref{IGIDE} and Lemma~\ref{MODELEM} is specific for Fourier Series.

For example, consider the sequence of functions
$$ \phi_N(x):=(-1)^{\lfloor 2x\rfloor} \big( \sin^2(2\pi x)\big)^{\frac{1}{N}},$$
let~$w$ be the square wave in Exercise~\ref{SQ:W},
and prove that, for all~$x\in\R\setminus(\Z/2)$,
\begin{equation}\label{smdkdKSjd02owdkHAScvkfdgdahJ-1} \lim_{N\to+\infty}\phi_N(x)=w(x)\end{equation}
and that, for each~$\delta\in\left(0,\frac12\right)$,
\begin{equation}\label{smdkdKSjd02owdkHAScvkfdgdahJ-2} \lim_{N\to+\infty}\sup_{x\in(0,\delta)}\big(\phi_N(x)-w(x)\big)=0.\end{equation}
\end{exercise}

\section{\faBomb Convergence issues}\label{EXCE}

Given the positive results presented so far about the convergence of the Fourier Series to the periodic function used to compute its Fourier coefficients, an optimist could believe that this always happens. But life is sometimes more complicated than expected and, as already mentioned in footnote~\ref{EXCE-foodf} on page~\pageref{EXCE-foodf} and also in the text of page~\pageref{0o2jrf98-l.1}, there are cases of divergent Fourier Series, so some care is needed when dealing with these delicate topics.

Also, using computers to visualise divergent Fourier Series of continuous functions is sometimes a bit tricky. Indeed, the function under consideration is typically continuous, but ``not too smooth'', otherwise the Fourier Series would converge (e.g., by Theorems~\ref{ONEGO} and~\ref{C1uni}). In this regard, a numerical simulation would need to plot a large number of points to provide a good picture. Moreover, the divergence of the Fourier Series implies that computing the values of an infinite sum just by adding more and more terms can fall beyond the numerical accuracy of the device.

Similarly, it may be difficult to explain ``geometrically'' or ``qualitatively'' the divergence of Fourier Series for continuous functions, and sometimes one needs to get their hands on some technical analytical details.

To familiarise with the lack of convergence of Fourier Series, we start with an explicit example, showing that the Fourier Series of a continuous function needs not to converge pointwise (though it does \label{ACarleson}
converges almost everywhere, due to the result of Carleson mentioned on page~\pageref{Carleson}):

\begin{theorem}\label{COMADIVESKMDWD}
There exist continuous functions of period~$1$ whose Fourier Series do not converge at a point.
\end{theorem}

\begin{proof} Here is an explicit example\footnote{For a slightly different methodology to construct
explicit examples that prove Theorem~\ref{COMADIVESKMDWD}, see Exercise~\ref{MSCOCOS}.
Another example will be given in Appendix~\ref{ANEXDIHALA}.} due to Fej\'er (see~\cite[page~264]{GOURDON}).

For all~$x\in\left[0,\frac12\right]$, let\footnote{Recalling the notational remark in footnote~\ref{PERGBDSE0232345021} on page~\pageref{PERGBDSE0232345021}, $2^{k^3}$ means, here and always, $2^{ ( k^3 ) }$.}
$$ f(x):=\sum_{k=1}^{+\infty}\frac1{k^2} \sin\left( \big(2^{ k^3 }+1\big)\pi x\right).$$
We stress that this is a valid definition, since the series above converges uniformly.
Since~$f(0)=0$, we can extend~$f$ by an even reflection across~$x=0$, and then take its periodic extension of period~$1$.
In this way, $f$ is continuous, even, and periodic of period~$1$. See Figure~\ref{LUNIikdnsEU} for a
(somewhat unavoidably rough) sketch of this situation.

\begin{figure}[h]
\includegraphics[height=3cm]{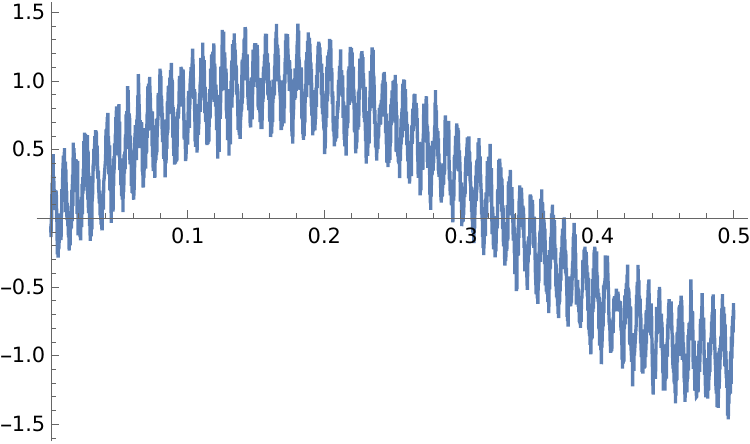}$\,\;$\includegraphics[height=3cm]{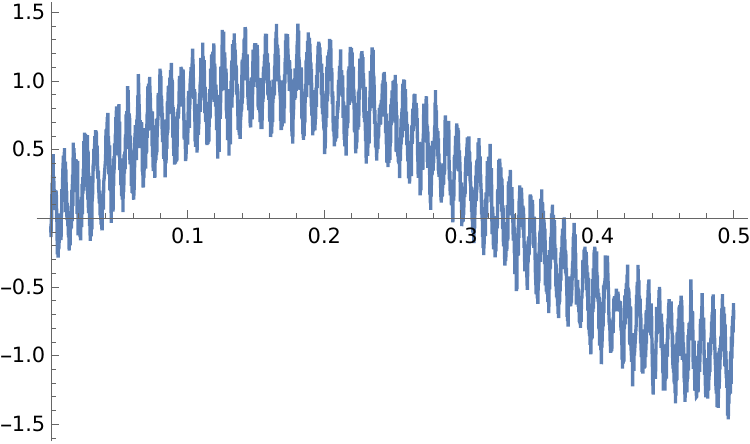}$\,\;$
\includegraphics[height=3cm]{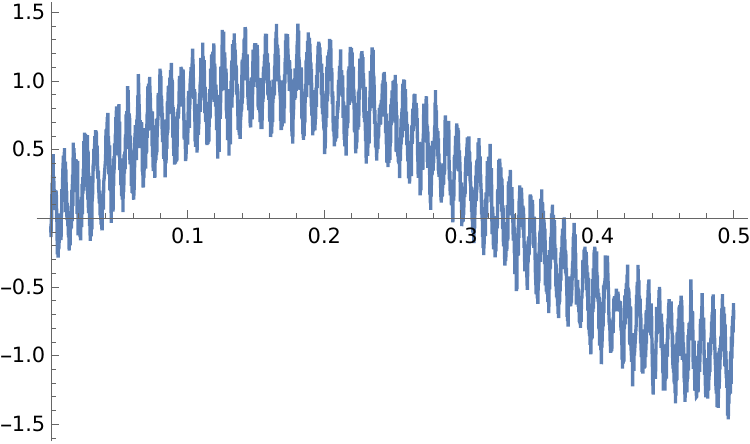}
\centering
\caption{Plot of~$\displaystyle\sum_{k=1}^{N}\frac1{k^2} \sin\left( \big(2^{ k^3 }+1\big)\pi x\right)$ with~$N\in\{10,\,50,\,100\}$.}\label{LUNIikdnsEU}
\end{figure}

We now write the Fourier Series of~$f$ in trigonometric form~\eqref{TRISI}.
To this end, since~$f$ is an even function, 
the Fourier Series of~$f$ (see Exercise~\ref{920-334PKSXu9o2fgfbsmos}) takes the form
\begin{equation}\label{OJSLM-2oekwf.13v523f23rfv02o8.2} \frac{a_0}2+\sum_{k=1}^{+\infty} a_k\cos(2\pi kx)\end{equation}
with (see Exercise~\ref{smc203e})
\begin{eqnarray*}a_k&=&4\int_0^{1/2} f(x)\,\cos(2\pi kx)\,dx\\&=&
\sum_{\ell=1}^{+\infty}\frac{4}{\ell^2}\int_0^{1/2} \sin\left( \big(2^{ \ell^3 }+1\big)\pi x\right)\,\cos(2\pi k x)\,dx,
\end{eqnarray*}
where the integral and series signs have been swapped thanks to the uniform convergence of the series.

From this (see Exercise~\ref{SPDCD-0.01anc}, used here with~$\alpha:=\big(2^{ \ell^3 }+1\big)\pi x$ and~$\beta:=2\pi kx$) we obtain that
$$ a_k=\sum_{\ell=1}^{+\infty}\frac{2}{\ell^2}\int_0^{1/2} 
\left[
\sin\left( \left(\big(2^{ \ell^3 }+1\big)-2 k\right)\pi x\right)+\sin\left(\left( \big(2^{ \ell^3 }+1\big)+2 k\right)\pi x\right)
\right]\,dx$$
and then, calculating the integrals,
\begin{eqnarray*} a_k&=&\sum_{\ell=1}^{+\infty}\frac{2}{\pi\,\ell^2}
\left(\frac{1}{ 1 + 2^{ \ell^3 } -2 k} + \frac{1}{1 + 2^{ \ell^3 }+ 2 k}\right)
.\end{eqnarray*}

As a result, recalling~\eqref{OJSLM-2oekwf.13v523f23rfv02o8.2},
\begin{equation}\label{OJSLM-2oekwf.13v523f23rfv02o8.3}\begin{split}
S_{N,f}(0)-\frac{a_0}2=\sum_{k=0}^{N} a_k
=\sum_{k=0}^{N}\sum_{\ell=1}^{+\infty}\frac{2}{\pi\,\ell^2}
\left(\frac{1}{ 1 + 2^{ \ell^3 } -2 k} + \frac{1}{1 + 2^{ \ell^3 }+ 2 k}\right).\end{split}
\end{equation}
We stress that the denominators above never vanish, because this would correspond to
$$ \N\ni k=\frac12+2^{ \ell^3 - 1 },$$
which is not an integer, since~$\ell\in\N\cap[1,+\infty)$.

Also, while the multiple summation index in~\eqref{OJSLM-2oekwf.13v523f23rfv02o8.3} may look scary at a first glance,
it just represents a finite sum (in~$k$) of finite quantities (the convergent series in~$\ell$, whose terms are also made by a sum of two objects).

As a consequence of~\eqref{OJSLM-2oekwf.13v523f23rfv02o8.3} (see Exercise~\ref{SENFYUJNSBKSxD}) we have that
\begin{equation}\label{OJSLM-2oekwf.13v523f23rfv02o8.4}\begin{split}
S_{N,f}(0)- \frac{a_0}2&=\sum_{\ell=1}^{+\infty}\sum_{k=0}^{N}\frac{2}{\pi\,\ell^2}
\left(\frac{1}{ 1 + 2^{ \ell^3 } -2 k} + \frac{1}{1 + 2^{ \ell^3 }+ 2 k}\right)\\
&=\frac2\pi \sum_{\ell=1}^{+\infty}\frac{\;\sigma_{N,2^{ \ell^3 }}\;}{\ell^2},\end{split}
\end{equation}
where, for all~$N$, $m\in\N\cap[1,+\infty)$, we set
$$ \sigma_{N,m}:=\sum_{k=0}^N
\left(\frac{1}{ 1 + 2m -2 k} + \frac{1}{1 + 2m+ 2 k}\right).$$

We observe that the summand above changes sign. Nevertheless, as we will see here below, we have that, for all~$N$, $m\in\N\cap[1,+\infty)$,
\begin{equation}\label{sigmaposikem-1}\sigma_{N,m}\ge0.
\end{equation}
Moreover, at least in one case the above term gets sufficiently large, namely, as we will now check,
\begin{equation}\label{sigmaposikem-2} \sigma_{m,m}\ge c\ln m,
\end{equation}
for some~$c>0$.

To prove~\eqref{sigmaposikem-1},
we first observe that~$\sigma_{N,m}\ge0$ when~$N\le m$, because in this case
$$ \frac{1}{ 1 + 2m -2 k} + \frac{1}{1 + 2m+ 2 k}=\frac{2(1+2m)}{ (1 + 2m )^2-4 k^2} $$ 
and, since~$(1 + 2m)^2 -4 k^2\ge
(1 + 2m)^2 -4 N^2\ge0$, the relevant summand is positive.

Hence, one has to prove~\eqref{sigmaposikem-1} when~$N>m$, as well as~\eqref{sigmaposikem-2}.
For this, we perform some algebraic manipulations 
on finite sums. Namely, we use the index substitutions~$j_1:=m-k$, $j_2:=k-m-1$ and~$j_3:=m+k$ to
observe that
\begin{eqnarray*}
\sigma_{N,m}&=&\sum_{k=0}^N\frac{1}{ 1 + 2m -2 k} +\sum_{k=0}^N \frac{1}{1 + 2m+ 2 k}\\
&=&\sum_{k=0}^m\frac{1}{ 1 + 2m -2 k}+\sum_{k=m+1}^{ N}\frac{1}{ 1 + 2m -2 k}
+\sum_{k=0}^N \frac{1}{1 + 2m+ 2 k}\\
&=&\sum_{j_1=0}^{m}\frac{1}{ 1 + 2j_1}-\sum_{j_2=0}^{ N-m-1}\frac{1}{1+2j_2}
+\sum_{j_3=m}^{N+m} \frac{1}{1 + 2j_3}\\&\ge&\sum_{j=0}^{N+m}\frac{1}{ 1 + 2j}-\sum_{j_2=0}^{ N-m-1}\frac{1}{1+2j_2}
\\&=&\sum_{j=N-m}^{N+m}\frac{1}{ 1 + 2j},
\end{eqnarray*}
that implies~\eqref{sigmaposikem-1}.

Moreover, the above observation for~$N:=m$ also gives that
\begin{eqnarray*}
\sigma_{m,m}\ge\sum_{j=0}^{2m}\frac{1}{ 1 + 2j},
\end{eqnarray*}
which leads to~\eqref{sigmaposikem-2}.

Now we combine~\eqref{sigmaposikem-1} and~\eqref{sigmaposikem-2} with~\eqref{OJSLM-2oekwf.13v523f23rfv02o8.4}
(picking~$N:=2^{ M^3 }$) and we find that
\begin{equation*}\begin{split}&
S_{2^{ M^3 },f}(0)- \frac{a_0}2=\frac2\pi \sum_{\ell=1}^{+\infty}\frac{\;\sigma_{2^{ M^3 },2^{ \ell^3 }}\;}{\ell^2}\ge \frac{\;2\sigma_{2^{ M^3 },2^{ M^3 }}\;}{\pi\,M^2}\ge\frac{2c\ln (2^{ M^3 })}{\pi\,M^2}=\frac{2cM\ln 2}\pi
,\end{split}
\end{equation*}
which diverges as~$M\to+\infty$, thus showing that this Fourier Series does not converge at~$x=0$.
\end{proof}

For further insight on periodic continuous functions whose Fourier Series diverges at a point,
see~\cite[Chapter~II, Section~2]{MR2039503},
\cite[Exercises 3.3.6 and 3.5.9]{MR3243734}, and~\cite[Chapter~18]{MR4404761}
(in particular, in~\cite[Remark after proof of Theorem 2.1]{MR2039503} one also finds a variation
of the arguments presented here leading to a continuous function whose Fourier Series diverges on a dense set).

Interestingly, the divergence of Fourier Series is also related to
the existence of continuous functions whose Fourier Series converges pointwise but not uniformly (see Exercise~\ref{MSCOCOS}, as well as~\cite[Chapter VIII]{MR1963498}
and the literature referred to there).

For other examples of divergent Fourier Series, see e.g.~\cite[Chapters~4--6]{MR44660} and~\cite[Part~1, Sections~2.1 and~2.3]{MR1408905}.

The reader already familiar with functional analysis (specifically, with the Baire Category Theorem) may like\footnote{The reader not
already familiar with the Baire Category Theorem may wish to skip Theorem~\ref{BaireCategoryTh} for the moment.}
to see that the example in Theorem~\ref{COMADIVESKMDWD} is not just an isolated pathology. Instead, given a point~$x_0\in\R$,
the set of continuous functions whose Fourier Series converges at~$x_0$ is somewhat ``meager'' and divergence is instead a ``generic phenomenon''. More precisely, we have that:

\begin{theorem}\label{BaireCategoryTh}
Let~$x_0\in\R$. Let~$C_{\text{per}}$ be 
the space of continuous functions, periodic of period~$1$, 
endowed with the $L^\infty(\R)$-norm.

Let
$$ {\mathcal{Z}}:=\left\{ f\in C_{\text{per}} {\mbox{ s.t. }} \sup_{N\in\N}|S_{N,f}(x_0)|<+\infty\right\}.$$

Then, the set~$ {\mathcal{Z}}$ has empty interior in~$C_{\text{per}}$.
\end{theorem}

\begin{proof} For every~$j\in\N$, we let
$$ {\mathcal{Z}}_j:=\left\{ f\in C_{\text{per}} {\mbox{ s.t. }} \sup_{N\in\N}|S_{N,f}(x_0)|\le j\right\}.$$
Since
$$ {\mathcal{Z}}=\bigcup_{j\in\N}{\mathcal{Z}}_j,$$
in light of the Baire Category Theorem (see e.g.~\cite[Theorem~5.6]{MR924157}) it suffices to show that 
\begin{equation}\label{CSDVAF}
{\mbox{each set~${\mathcal{Z}}_j$ is closed and has empty interior in~$C_{\text{per}}$.}}\end{equation}

For this, let
$$ {\mathcal{Z}}_{j,N}:=\left\{ f\in C_{\text{per}} {\mbox{ s.t. }}|S_{N,f}(x_0)|\le j\right\}.$$
Let~$f\in C_{\text{per}}\setminus {\mathcal{Z}}_{j,N}$ and~$\epsilon>0$. Then, $ |S_{N,f}(x_0)|>j$
and if~$g\in C_{\text{per}}$ with~$0<\|f-g\|_{L^\infty(\R)}\le\epsilon$, setting~$h:=\frac{f-g}{\|f-g\|_{L^\infty(\R)}}$ we have 
that
\begin{eqnarray*}
&&|S_{N,g}(x_0)|\ge|S_{N,f}(x_0)|-|S_{N,f-g}(x_0)|
=|S_{N,f}(x_0)|-\|f-g\|_{L^\infty(\R)}\,\left|S_{N,h}(x_0)\right|\\
&&\qquad\ge |S_{N,f}(x_0)|-\epsilon\sigma_N,
\end{eqnarray*}
where
$$\sigma_N:=\sup_{{\phi\in C_{\text{per}}}\atop{\|\phi\|_{L^\infty(\R)}=1}} |S_{N,\phi}(x_0)|.$$
We stress that~$\sigma_N<+\infty$ (see~\eqref{TACLCSIMOMNDADEKPD.1} in Exercise~\ref{SRTAMAYBHJAOPLKAZXAQNBKNEMINSNC}) and therefore if~$\epsilon\in\left(0,\frac{|S_{N,f}(x_0)|-j}{\sigma_N}\right)$ we have that~$|S_{N,g}(x_0)|>j$, whence~$g\in C_{\text{per}}\setminus{\mathcal{Z}}_{j,N}$.

This yields that~$C_{\text{per}}\setminus{\mathcal{Z}}_{j,N}$ is open and accordingly~${\mathcal{Z}}_{j,N}$ is closed in~$C_{\text{per}}$.

On this account, the set~${\mathcal{Z}}_j$ (which is the countable intersection of the closed sets~${\mathcal{Z}}_{j,N}$
over~$N\in\N$) is also closed.

Thus, to complete the proof of~\eqref{CSDVAF}, we need to show that~${\mathcal{Z}}_j$ has empty interior in~$C_{\text{per}}$.
To this end, we argue by contradiction and assume that there exist~$j\in\N$,
$f_j\in {\mathcal{Z}}_j$ and~$\epsilon_j>0$ such that if~$\psi\in C_{\text{per}}$ with~$\|f_j-\psi\|_{L^\infty(\R)}\le\epsilon_j$, then~$\psi\in{\mathcal{Z}}_j$. 

Hence, for every~$f\in C_{\text{per}}$ with~$\|f\|_{L^\infty(\R)}=1$, we take~$\psi_j:=f_j-\epsilon_j f$. In this way, we have that
$$\|f_j-\psi_j\|_{L^\infty(\R)}=\epsilon_j\|f\|_{L^\infty(\R)}=\epsilon_j$$
and therefore~$\psi_j\in{\mathcal{Z}}_j$.

In virtue of this fact,
$$ j\ge\sup_{N\in\N}|S_{N,\psi_j}(x_0)|\ge \sup_{N\in\N}|S_{N,\epsilon_j f}(x_0)|- \sup_{N\in\N}|S_{N,f_j}(x_0)|
\ge\epsilon_j \sup_{N\in\N}|S_{N,f}(x_0)|-j$$
and thus, for all~$N\in\N$,
$$ |S_{N,f}(x_0)|\le\frac{2j}{\epsilon_j}.$$
Since this is valid for every~$f\in C_{\text{per}}$ with~$\|f\|_{L^\infty(\R)}=1$, we conclude that
$$ \sup_{{f\in C_{\text{per}}}\atop{\|f\|_{L^\infty(\R)}=1}} |S_{N,f}(x_0)|\le\frac{2j}{\epsilon_j},$$
but this is absurd
(see~\eqref{TACLCSIMOMNDADEKPD.2} in Exercise~\ref{SRTAMAYBHJAOPLKAZXAQNBKNEMINSNC}).\end{proof}

The divergence phenomena for Fourier Series are a very difficult, but fascinating topic.
Historically, it dates back to~1876, when du Bois-Reymond constructed the first, and rather complicated,
example of a continuous function whose Fourier Series diverges at one point; since then, several other explicit examples have been provided, also to exhibit even wilder behaviours.

Among the many classical examples, one of the most famous was found by
Kolmogoroff (at the time, a third-year student)
who first constructed an explicit example of Fourier Series that does not converge in a set of full measure
(see~\cite{UL10115} and~\cite[Chapter~VIII, Theorem~(3$\cdot$1)]{MR1963498} for full details) and then refined this result
by constructing another explicit example of Fourier Series that converges nowhere
(see~\cite[Chapter~VIII, Theorem~(4$\cdot$1)]{MR1963498} for full details;
see also~\cite{MR1009439} for a scientific biography of Kolmogoroff).

The result by Kolmogoroff was also retaken by Marcinkiewicz~\cite{zbMATH03022296}, who constructed an example of
Fourier Sum which remains bounded in a set of full measure, without converging in this set
(see also~\cite{MR96068}, in which it is shown that one can prescribe any divergent sequence~$N_k$ for which~$|S_{N_k,f}(x)|\to+\infty$
for almost every~$x\in\R$ as~$k\to+\infty$, i.e. one can 
generalise Kolmogoroff's example by
even arbitrarily picking a divergent sequence of Fourier Sums; for related results, see~\cite{zbMATH03903268}).

For continuous periodic functions, the convergence properties\footnote{A somewhat related question regarding the Fourier Series of continuous periodic functions is: how fast do the Fourier coefficients of these functions decay? On the one hand, we know from Bessel's Inequality~\eqref{AJSa} that if~$f$ is continuous and periodic of period~$1$ (whence in~$L^2((0,1))$), then necessarily
$$ \sum_{k\in\Z}|\widehat f_k|^2<+\infty.$$
On the other hand, one can show that this decay cannot in general be improved: namely, for every sequence~$c_k$ such that
$$ \sum_{k\in\Z}|c_k|^2<+\infty,$$
one can find a continuous function~$f$, periodic of period~$1$, such that~$|\widehat f_k|\ge |c_k|$, for all~$k\in\Z$. See~\cite[Theorem~2.1 on page~278]{MR2039503} for full details on this.}
of Fourier Series are also rich and very complicated:
for example, it is proved in~\cite{MR199633}
that for any set~$E\subset(0,1)$ of zero measure, there is a continuous function~$f$, periodic of period~$1$,
whose Fourier Series diverges on~$E$ (and this result is optimal, since this type of Fourier Series must converge
almost everywhere, due to the result by
Carleson mentioned on page~\pageref{Carleson}.

\begin{exercise}\label{SENFYUJNSBKSxD}
Consider a sequence of real numbers~$\{\mu_{k,\ell}\}_{k,\ell\in\N}$ such that, for all~$k\in\N$,
\begin{equation}\label{SENFYUJNSBKSxD0}
{\mbox{the series }}\sum_{\ell=0}^{+\infty}\mu_{k,\ell}{\mbox{ is convergent.}}\end{equation}
Prove that, for every~$N\in\N$,
$$ \sum_{k=0}^N\sum_{\ell=0}^{+\infty}\mu_{k,\ell}=\sum_{\ell=0}^{+\infty}\sum_{k=0}^N\mu_{k,\ell}.$$
\end{exercise}

\begin{exercise}\label{SRTAMAYBHJAOPLKAZXAQNBKNEMINSNC}
Let~$x_0\in\R$. Let~$C_{\text{per}}$ be as in Theorem~\ref{BaireCategoryTh}.

Prove that, for all~$N\in\N$,
\begin{equation}\label{TACLCSIMOMNDADEKPD.1}
\sup_{{f\in C_{\text{per}}}\atop{\|f\|_{L^\infty(\R)}=1}} |S_{N,f}(x_0)|<+\infty
\end{equation}
and that
\begin{equation}\label{TACLCSIMOMNDADEKPD.2}\lim_{N\to+\infty}
\sup_{{f\in C_{\text{per}}}\atop{\|f\|_{L^\infty(\R)}=1}} |S_{N,f}(x_0)|=+\infty.
\end{equation}\end{exercise}

\begin{exercise}\label{SRTAMAYBHJAOPLKAZXAQNBKNEMINSNC.L1pe}
Let~$L^1_{\text{per}}$ be the class of periodic functions~$f$, of period~$1$, such that~$f\in L^1((0,1))$.

Prove that, for all~$N\in\N$,
\begin{equation}\label{TACLCSIMOMNDADEKPD.1pe}
\sup_{{f\in L^1_{\text{per}}}\atop{\|f\|_{L^1((0,1))}=1}} \|S_{N,f}\|_{L^1((0,1))}<+\infty
\end{equation}
and that
\begin{equation}\label{TACLCSIMOMNDADEKPD.2pe}\lim_{N\to+\infty}
\sup_{{f\in L^1_{\text{per}}}\atop{\|f\|_{L^1((0,1))}=1}}\|S_{N,f}\|_{L^1((0,1))}=+\infty.
\end{equation}\end{exercise}

\begin{exercise}\label{KDKAMSDW:SDCMILCSKMc}
Endow~$L^1_{\text{per}}$ with the norm~$\|\cdot\|_{L^1((0,1))}$ and prove that the set of functions~$f\in L^1_{\text{per}}$
for which~$S_{N,f}$ converges to~$f$ in~$L^1((0,1)$ has empty interior in~$L^1_{\text{per}}$.
\end{exercise}

\begin{exercise}\label{CLURVB}
Let~$\beta>0$. Prove that
$$ \sup_{\alpha\ge\beta}\left|\int_{\alpha}^1\frac{\cos(\pi t)}{t}\,dt\right|<+\infty.$$
\end{exercise}

\section{Functions of arbitrary periods}\label{ANY}

For concreteness, we have focused here on periodic functions of period~$1$.
The case of periodic functions with different periodicity is completely similar and offers no conceptual difficulties: in fact, for a periodic function~$f\in L^1((0,P))$ with periodicity of~$P\in(0,+\infty)$, one can consider the function~$ g(x):=f(Px)$, which is periodic of period~$1$, develop all the theory for~$g$ and reconstruct from this the analogous results for the original function~$f$.

In particular (see Exercise~\ref{fc:PkT-0.3}), for a periodic function of period~$P$
the Fourier coefficients in~\eqref{FOUCO} become
\begin{equation}\label{fc:PkT-0.1}
\widehat f_k={\frac{1}{P}}\int _{0}^P f(x)\,e^{-{\frac{2\pi ikx}{P}}}\,dx\end{equation}
and the corresponding Fourier coefficients in trigonometric form in~\eqref{jasmx23er} are
\begin{equation}\label{fc:PkT-0.2}\begin{split}
&a_k={\frac{2}{P}}\int _{0}^P f(x)\,\cos \left({\frac{2\pi kx}{P}}\right)\,dx\\
{\mbox{and }}\qquad&
b_k={\frac{2}{P}}\int _{0}^P f(x)\,\sin \left({\frac{2\pi kx}{P}}\right)\,dx.
\end{split}\end{equation}

From~\eqref{fc:PkT-0.1} one can appreciate the notable property that~$\widehat f_0$ is actually the average value of the function~$f$ over its periodicity domain.

In the same vein, the formal Fourier Series in~\eqref{FOSE} for a function of period~$P$ takes the form
\begin{equation}\label{fc:PkT-0.45}
\sum _{k\in\Z}\widehat f_{k}\, e^{\frac{2\pi ikx}{P}} 
\end{equation}
and the corresponding trigonometric representation~\eqref{TRISI} is
\begin{equation}\label{fc:PkT-0.55}
\frac{a_0}2+\sum _{k=1}^{\infty }\left(a_k\cos \left({\frac{2\pi kx}{P}}\right)+b_k\sin \left({\frac{2\pi kx}{P}}\right)\right).
\end{equation} 

\begin{exercise}\label{fc:PkT-0.3} Check that the Fourier coefficients
of a periodic function of period~$P$
take the form presented in~\eqref{fc:PkT-0.1} and~\eqref{fc:PkT-0.2},
and that the corresponding Fourier Series are the ones in~\eqref{fc:PkT-0.45} and~\eqref{fc:PkT-0.55}.
\end{exercise}

\begin{exercise}
Let~$L>0$, $\lambda\in\R$ and~$f:\R\to\R$ be a continuous and odd function, periodic of periodic~$2L$, with~$f(0)=f(L)=0$ and~$f\left(\frac{L}4\right)=\lambda$. 

Suppose that the graph of~$f$ coincides with a straight line in~$\left(0,\frac{L}4\right)$ and in~$\left(\frac{L}4,L\right)$.

Calculate the Fourier coefficients of~$f$.\end{exercise}

\begin{exercise}\label{013p24rtg.01324r5t4y56.0}
Let
$$f(x):=\begin{dcases}1&{\mbox{ if }}x\in\Q,\\
0&{\mbox{ if }}x\in\R\setminus\Q.
\end{dcases}$$
Prove that every positive rational number is a period
for~$f$.
\end{exercise}

\begin{exercise}\label{013p24rtg.01324r5t4y56}
Let~$f:\R\to\R$ be a continuous, nonconstant, periodic function.
Prove that the set
\begin{equation}\label{013p24rtg.01324r5t4y56.e} \big\{ {\mbox{$T>0$ s.t. $T$ is a period for $f$}}\big\}\end{equation}
admits a minimum~$T_0$ (which is called in jargon \index{fundamental period}
the ``fundamental period'' of~$f$).
\end{exercise}

\begin{exercise}\label{013p24rtg.01324r5t4y56.b} In the setting of Exercise~\ref{013p24rtg.01324r5t4y56}, prove that if~$T>0$ is any period of~$f$, then~$T=kT_0$, for some~$k\in\N\cap[1,+\infty)$.
\end{exercise}

\begin{exercisesk}\label{013p24rtg.01324r5t4y56.bb}
Prove an analogue of Exercise~\ref{013p24rtg.01324r5t4y56}
in Lebesgue spaces. Namely, let~$\tau>0$ and~$f:\R\to\R$ be a periodic function of period~$\tau$, with~$f\in L^1([0,\tau])$.

Prove that either the set in~\eqref{013p24rtg.01324r5t4y56.e} 
admits a minimum (hence~$f$ possesses a fundamental period),
or~$f$ is constant (in the sense of Lebesgue,
namely there exist~$c\in\R$ and~$Z\subset\R$ such that~$Z$
has null Lebesgue measure and~$f(x)=c$ for every~$x\in\R\setminus Z$).
\end{exercisesk}

\begin{exercise}
Show that the function $$\R\ni x\longmapsto
\cos(8\pi x)+\cos(12\pi x)+\sin^2(2\pi x)+7$$
is periodic and calculate its fundamental period.\end{exercise}

\begin{exercise}\label{QUASDFG-qaudwe4y6.PREC}
Let~$\phi:\R\to\R$ be continuous and periodic with periods~$T_1$, $T_2>0$.
Assume that~$\frac{T_1}{T_2}$ is an irrational number.

Prove that~$\phi$ is constant.
\end{exercise}

\begin{exercise}\label{QUASDFG-qaudwe4y6}
Let~$f_1$, $f_2:\R\to\R$ be continuous, nonconstant functions. Suppose that~$f_1$ is periodic of period~$T_1>0$ and that~$f_2$ is periodic of period~$T_2>0$.
Assume also that~$\frac{T_1}{T_2}$ is an irrational number and let~$g:=f_1+f_2$.

Prove that~$g$ is {\em not} a periodic function.
\end{exercise}

\section{Fourier Series in any dimension}\label{SEC:ANUYDIMEAGa}

Now we briefly discuss what happens to Fourier Series in higher dimension. The quick statement in this case is that everything
carries over, except what doesn't!

Let's see why.

A function~$f:\R^n\to\R$ is said to be $\Z^n$-periodic \index{$\Z^n$-periodic} if,
for every~$k\in\Z^n$, \label{PKSMDODKUIOIYFDYPIJBN8.17y2ert0012ihfbbgrce4gtdk}
$$ f(x+k)=f(x).$$

In analogy with~\eqref{FOUCO}, given~$k\in\Z^n$, the $k$th Fourier of a $\Z^n$-periodic function~$f\in L^1((0,1)^n)$ is defined by
$$\widehat f_k:=\int_{(0,1)^n} f(x)\,e^{-2\pi i k\cdot x}\,dx$$
and, similarly to~\eqref{FOSE}, the corresponding Fourier Series is defined as
$$ \sum_{k\in\Z^n} \widehat f_k \,e^{2\pi i k\cdot x}.$$

One may be interested in proving the convergence of this series to the original function, but, as we will see, there is some caveat
since in higher dimension even the definition of Fourier Sum requires a
careful specification of the precise order of summation to reconstruct the Fourier Series from a Fourier Sum. For example,
a natural multidimensional extension of~\eqref{FourierSum} consists in
\begin{equation*}
S_{N,f,{\scriptscriptstyle\mbox{cir}}}(x):=\sum_{{k=(k_1,\dots,k_n)\in\Z^n}\atop{\sqrt{k_1^2+\dots k_n^2}\le N}} \widehat f_k \,e^{2\pi i k\cdot x},
\end{equation*}
which is known as \emph{circular Fourier Sum} \index{circular Fourier Sum}
(where ``circular'' reflects the idea that the norm used on the indices~$k$ is the Euclidean one,
whose level sets are circles).

However, another possible multidimensional extension of~\eqref{FourierSum} consists in
\begin{equation}\label{SKMELSQFSS}
S_{N,f,{\scriptscriptstyle\mbox{squ}}}(x):=\sum_{{k=(k_1,\dots,k_n)\in\Z^n}\atop{|k_1|\le N,\,\dots,\,|k_n|\le N}} \widehat f_k \,e^{2\pi i k\cdot x},
\end{equation}
which is instead known as \emph{square Fourier Sum} \index{square Fourier Sum}
(where ``square'' describes the level sets on the norm used on the indices~$k$; the nomenclature used varies anyway from
author to author).

And, in general, one could consider possible multidimensional extensions of~\eqref{FourierSum} by taking into account
different norms on~$\R^n$, see e.g.~\cite{MR2943169} for a survey on this topic.

The difference between these possible definitions of Fourier Sums is very deep, since conceptually different phenomena may arise
and several one-dimensional results are unknown in multiple dimensions or require a careful specification on summation assumptions. For example, the equivalent of Carleson's result mentioned on pages~\pageref{Carleson}
and~\pageref{ACarleson} is still unknown for circular Fourier Sums (and it is false in general for ``rectangular'' Fourier Sums~\cite{MR279529})
but it has been established
by Fefferman in~\cite{MR435724, MR340926, MR1469321} for square Fourier Sums and for other types of ``polygonal'' Fourier Sums.

So, the bottom line is that choosing a notion of Fourier Sum in higher dimension may affect the corresponding results;
likely, the square Fourier Sums having some convenience over the others, since they usually entail nice factorisations into products of one-dimensional objects (see e.g. Exercise~\ref{SQDITJMS0x12S}).

On a positive note, most of the standard results about decay of Fourier coefficients and convergence of Fourier Series for functions in the appropriate sense remain valid, and one can easily deduce them by extending the proofs presented so far in the one-dimensional case to the higher-dimensional setting, e.g. using a convenient multi-index notation (or by deducing the high-dimensional results from the one-dimensional ones).

In any case, we will not dive here into the full development of the general theory of Fourier Series in
multiple dimensions.
Rather, we will provide just some very basic results which will be sufficient for the goals of this book (and the methods of proof can
be also extended to obtain sharper results). We start with 
a global convergence result in a very smooth setting.

For this we say that a sequence~$\{{\mathcal{U}}_\ell\}_{\ell\in\N}$ of finite subsets of~$\Z^n$ is an invasion of~$\Z^n$ if, for all~$\ell\in\N$,
\begin{equation}\label{INVA-01}
{\mathcal{U}}_\ell\subseteq{\mathcal{U}}_{\ell+1}\end{equation}
and
\begin{equation}\label{INVA-02} \bigcup_{\ell\in\N}{\mathcal{U}}_\ell=\Z^n.\end{equation}

With this notation, we have:

\begin{theorem}\label{SELINOPOR} Let~$\{{\mathcal{U}}_\ell\}_{\ell\in\N}$ be an invasion of~$\Z^n$
and~$f\in C^{n+1}(\R^n)$ be~$\Z^n$-periodic.

Then,
\begin{equation}\label{PRECON56789-BIS.lo8} \lim_{N\to+\infty}\sup_{x\in\R^n}\left|f(x)-\sum_{k\in{\mathcal{U}}_N} \widehat f_k \,e^{2\pi i k\cdot x}\right|=0.\end{equation}
\end{theorem}

\begin{proof} One could go through the arguments in Section~\ref{UNIFORMCO:SECTI} and check that they can be suitably extended to higher dimension, but, to provide a more self-contained exposition, we argue directly as follows.

We observe that
\begin{equation}\label{EUGUNOCO}\begin{split}&
{\mbox{if we know the result in~\eqref{PRECON56789-BIS.lo8} for a given invasion~$\{{\mathcal{U}}_\ell\}_{\ell\in\N}$ of~$\Z^n$,}}\\&{\mbox{then we also know it for any other invasion~$\{{\mathcal{V}}_\ell\}_{\ell\in\N}$ of~$\Z^n$.}}\end{split}\end{equation}
For this, one observes that, for all~$k\in\Z^n$ and all~$\mu=(\mu_1,\dots,\mu_n)\in\N^n$ with~$\mu_1+\dots+\mu_n\le n+1$, we have that
$$ \widehat{D^\mu f}_k=(2\pi i k)^\mu\widehat f_k,$$
where, as customary, we use the multi-index notation~$k^\mu=k_1^{\mu_1}\dots k_n^{\mu_n}$.
This can be obtained by proceeding as in the proof of Theorem~\ref{SMXC22b}, just using multi-indices.

As a consequence (see Exercise~\ref{PRECON56789-BIS.lo}) we obtain that
\begin{equation}\label{oslm0-AmdfiNBNAcH} |\widehat f_k|\le\frac{C}{|k|^{n+1}},\end{equation}
for some~$C\ge0$, where~$|\cdot|$ is the Euclidean norm.

Now we recall~\eqref{INVA-01} and~\eqref{INVA-02}, and,
given~$L>0$, to be taken as large as we wish in what follows, we take~$N_L\in\N$ sufficiently large such that
for all~$N\ge N_L$
$$ \big\{ k\in \Z^n{\mbox{ s.t. }}|k|\le L\big\}\subseteq\bigcup_{\ell=0}^N{\mathcal{U}}_\ell={\mathcal{U}}_N$$
and similarly
$$ \big\{ k\in \Z^n{\mbox{ s.t. }}|k|\le L\big\}\subseteq{\mathcal{V}}_N.$$

In particular, taking~$N\ge N_L$, we have that if~$k\in\Z^n\setminus{\mathcal{U}}_N$, or if~$k\in\Z^n\setminus{\mathcal{V}}_N$, then necessarily~$|k|>L$.

On this account, and using~\eqref{PRECON56789-BIS.lo}, we see that, if~$N\ge N_L$,
\begin{eqnarray*}&&
\left|\sum_{k\in{\mathcal{U}}_N} \widehat f_k \,e^{2\pi i k\cdot x}-\sum_{k\in{\mathcal{V}}_N} \widehat f_k \,e^{2\pi i k\cdot x}\right|\le
\sum_{k\in{\mathcal{U}}_N\setminus{\mathcal{V}}_N} |\widehat f_k |+\sum_{k\in{\mathcal{V}}_N\setminus{\mathcal{U}}_N} |\widehat f_k |\\&&\qquad\le2
\sum_{{k\in\Z^n}\atop{|k|>L}} |\widehat f_k |\le2C
\sum_{{k\in\Z^n}\atop{|k|>L}}\frac{1}{|k|^{n+1}},
\end{eqnarray*}
which is the tail of a convergent series.

Therefore, if we know~\eqref{PRECON56789-BIS.lo8} for the invasion~$\{{\mathcal{U}}_\ell\}_{\ell\in\N}$ of~$\Z^n$, we also know that
\begin{eqnarray*}
&&\lim_{N\to+\infty}\sup_{x\in\R^n}\left|f(x)-\sum_{k\in{\mathcal{V}}_N} \widehat f_k \,e^{2\pi i k\cdot x}\right|\\
&&\qquad\le
\lim_{N\to+\infty}\left(\sup_{x\in\R^n}\left|f(x)-\sum_{k\in{\mathcal{U}}_N} \widehat f_k \,e^{2\pi i k\cdot x}\right|
+\sup_{x\in\R^n}\left|\sum_{k\in{\mathcal{U}}_N}-\sum_{k\in{\mathcal{V}}_N} \widehat f_k \,e^{2\pi i k\cdot x}\right|\right)
\\&&\qquad\le0+2C
\sum_{{k\in\Z^n}\atop{|k|>L}}\frac{1}{|k|^{n+1}}.
\end{eqnarray*}
Therefore, sending~$L\to+\infty$,
$$\lim_{N\to+\infty}\sup_{x\in\R^n}\left|f(x)-\sum_{k\in{\mathcal{V}}_N} \widehat f_k \,e^{2\pi i k\cdot x}\right|=0,$$
which is~\eqref{PRECON56789-BIS.lo8} for the invasion~$\{{\mathcal{V}}_\ell\}_{\ell\in\N}$ of~$\Z^n$, thus establishing~\eqref{EUGUNOCO}.

Thus, in view of~\eqref{EUGUNOCO}, it suffices to check~\eqref{PRECON56789-BIS.lo8} using the notion of square Fourier Sum
in~\eqref{SKMELSQFSS} (which indeed provides an invasion of~$\Z^n$ with sets~${\mathcal{U}}_j:=\{k\in\Z^n$ s.t. $|k_1|\le j,\dots,
|k_n|\le j\}$).

What is more, one can focus on pointwise convergence, namely we claim that it suffices to check that, for every~$x\in\R^n$,
\begin{equation}\label{UNIFORMCO:SECTIc}
\lim_{N\to+\infty}\sum_{{{k=(k_1,\dots,k_n)\in\Z^n}\atop{|k_1|\le N,\,\dots,\,|k_n|\le N}}} \widehat f_k \,e^{2\pi i k\cdot x}=f(x).
\end{equation}
Indeed, suppose that~\eqref{UNIFORMCO:SECTIc} holds true. Then, owing to~\eqref{oslm0-AmdfiNBNAcH},
we observe that, for all~$M> N$,
\begin{eqnarray*}&&
\left|\sum_{{{k=(k_1,\dots,k_n)\in\Z^n}\atop{|k_1|\le M,\,\dots,\,|k_n|\le M}}} \widehat f_k \,e^{2\pi i k\cdot x}-\sum_{{{k=(k_1,\dots,k_n)\in\Z^n}\atop{|k_1|\le N,\,\dots,\,|k_n|\le N}}} \widehat f_k \,e^{2\pi i k\cdot x}\right|\le
\sum_{j=1}^n\sum_{{{{k=(k_1,\dots,k_n)\in\Z^n}\atop{|k_1|\le M,\,\dots,\,|k_n|\le M}}}\atop{|k_j|\ge N+1}}| \widehat f_k|\\&&\qquad\le 
\sum_{j=1}^n\sum_{{{{k=(k_1,\dots,k_n)\in\Z^n}\atop{|k_1|\le M,\,\dots,\,|k_n|\le M}}}\atop{|k_j|\ge N+1}}\frac{C}{|k|^{n+1}}
\le\sum_{{{{k=(k_1,\dots,k_n)\in\Z^n}\atop{|k|\ge N+1}}}}\frac{Cn}{|k|^{n+1}},
\end{eqnarray*}
which is the tail of a convergent series.

Therefore, one obtains the uniform convergence of the associated series, i.e. there exists a function~$g$ such that
$$ \lim_{N\to+\infty}\sup_{x\in\R^n}\left|g(x)-\sum_{{k\in\Z^n}\atop{\|k\|\le N}} \widehat f_k \,e^{2\pi i k\cdot x}\right|=0.$$
On grounds of this and~\eqref{UNIFORMCO:SECTIc}, it follows that~$g$ coincides with~$f$, and this gives the desired result in~\eqref{PRECON56789-BIS.lo8} for the square Fourier Sum.

These considerations show that we can focus on the proof of~\eqref{UNIFORMCO:SECTIc}.
For this, we write that
\begin{eqnarray*}
S_{N,f,{\footnotesize\mbox{squ}}}(x)&
=&\sum_{{k=(k_1,\dots,k_n)\in\Z^n}\atop{|k_1|\le N,\,\dots,\,|k_n|\le N}} \widehat f_k \,e^{2\pi i k\cdot x}
\\&=&
\sum_{{k=(k_1,\dots,k_n)\in\Z^n}\atop{|k_1|\le N,\,\dots,\,|k_n|\le N}} \int_{(0,1)^n}f(y) \,e^{2\pi i k\cdot (x-y)}\,dy
\\&=&\int_{(0,1)^n}f(y) \,D_N(x_1-y_1)\dots D_N(x_n-y_n)\,dy,
\end{eqnarray*}
where~$D_N$ is the one-dimensional Dirichlet Kernel in~\eqref{PAKSw-L4} (compare with Exercise~\ref{SQDITJMS0x12S}).

Hence (see Exercise~\ref{fr12}) we have that
$$ S_{N,f,{\footnotesize\mbox{squ}}}(x)=\int_{(0,1)^n}f(x-y) \,D_N(y_1)\dots D_N(y_n)\,dy.$$

We also remark that, by Exercise~\ref{K-1PIO},
$$ \int_{(0,1)^n} D_N(y_1)\dots D_N(y_n)\,dy=\Bigg(\int_0^1 D_N(y_1)\,dy_1\Bigg)\dots\Bigg(\int_0^1 D_N(y_n)\,dy_n\Bigg)=1$$
and therefore
\begin{eqnarray*}
S_{N,f,{\footnotesize\mbox{squ}}}(x)-f(x)=\int_{(0,1)^n}\big(f(x-y)-f(x)\big) \,D_N(y_1)\dots D_N(y_n)\,dy.
\end{eqnarray*}

Hence, by Lemma~\ref{KASMqwdfed123erDKLI},
\begin{eqnarray*}&&
S_{N,f,{\footnotesize\mbox{squ}}}(x)-f(x)\\&&\qquad=\int_{(0,1)^n}\big(f(x-y)-f(x)\big) \,
\frac{\sin\big((2N+1)\pi y_1\big)\dots\sin\big((2N+1)\pi y_n\big)}{\sin(\pi y_1)\dots\sin(\pi y_n)}
\,dy.
\end{eqnarray*}

Given~$x\in\R^n$, we now consider the function
{\footnotesize{$$ \R\ni y_1\longmapsto \phi^{(x)}(y_1):=
\int_{(0,1)^{n-1}}\big(f(x-y)-f(x)\big) \,
\frac{\sin\big((2N+1)\pi y_2\big)\dots\sin\big((2N+1)\pi y_n\big)}{\sin(\pi y_1)\dots\sin(\pi y_n)}
\,dy_2\,\dots,dy_n$$}}
and we see that
\begin{eqnarray*}
S_{N,f,{\footnotesize\mbox{squ}}}(x)-f(x)=\int_{0}^1
\phi^{(x)}(y_1)\,\sin\big((2N+1)\pi y_1\big)\,dy_1,
\end{eqnarray*}
which is the imaginary part of a Fourier coefficient of the bounded\footnote{The reader may notice that
less regularity is needed for this argument, thus highlighting how square Fourier Sums are often easier to deal with.} function~$\phi^{(x)}$.

Accordingly, the Riemann-Lebesgue Lemma (see Theorem~\ref{RLjoqwskcdc}) returns that
$$ \lim_{N\to+\infty} S_{N,f,{\footnotesize\mbox{squ}}}(x)-f(x)=0,$$
that is~\eqref{UNIFORMCO:SECTIc}.
\end{proof}

Now we briefly address the theory of Fourier Series in~$L^2((0,1)^n)$:

\begin{theorem}\label{L2INDIN0igwi92jgm.2}
Let~$\{{\mathcal{U}}_\ell\}_{\ell\in\N}$ be an invasion of~$\Z^n$ and~$f\in L^2((0,1)^n)$ be~$\Z^n$-periodic.

Then,
\begin{equation}\label{PARSENEN9c1903x} \lim_{N\to+\infty}
\left\|f-\sum_{k\in{\mathcal{U}}_N} \widehat f_k \,e^{2\pi i k\cdot x}\right\|_{L^2((0,1)^n)}=0.\end{equation}

Moreover, the following \index{Parseval's Identity} Parseval's Identity holds true:
\begin{equation}\label{PARSENEN} \sum_{{k\in\Z^n}} |\widehat f_k|^2=\|f\|_{L^2((0,1)^n)}^2.\end{equation}
\end{theorem}

\begin{proof} One could extend the methods put forth in Sections~\ref{L2pe}, \ref{SFOL2}, and~\ref{COL2},
but, to provide a more self-contained exposition, we argue directly as follows.

We consider a finite set~$Z\subseteq\Z^n$ and we observe that
\begin{equation}\label{FCONJHALKDD-1}\begin{split}&
\left\| \sum_{k\in Z} \widehat f_k \,e^{2\pi i k\cdot x}\right\|_{L^2((0,1)^n)}^2=
\sum_{{k,h\in Z}}
\widehat f_k\,\overline{\widehat f_h}\int_{(0,1)^n} e^{2\pi i (k-h)\cdot x}\,dx
=\sum_{{k\in Z}}|\widehat f_k|^2.
\end{split}
\end{equation}
As a result,
\begin{equation}\label{PARSENEN9c1903x2}\begin{split}&
\left\| f-\sum_{k\in Z} \widehat f_k \,e^{2\pi i k\cdot x}\right\|_{L^2((0,1)^n)}^2\\&\qquad=
\|f\|_{L^2((0,1)^n)}^2+\sum_{{k\in Z}}
|\widehat f_k|^2\\&\qquad\qquad
-\sum_{{k\in Z}}
\overline{\widehat f_k}\int_{(0,1)^n} f(x)\,e^{-2\pi ik\cdot x}\,dx
-\sum_{{k\in Z}}
\widehat f_k\,\int_{(0,1)^n} f(x)\,e^{2\pi i k\cdot x}\,dx
\\&\qquad=\|f\|_{L^2((0,1)^n)}^2+\sum_{{k\in Z}}
|\widehat f_k|^2-2\sum_{{k\in Z}}
|\widehat f_k|^2\\&\qquad=\|f\|_{L^2((0,1)^n)}^2-\sum_{{k\in Z}}
|\widehat f_k|^2
\end{split}
\end{equation}
and in particular
\begin{equation}\label{BEGH0-019340.1} \sum_{{k\in \Z^n}}|\widehat f_k|^2\le\|f\|_{L^2((0,1)^n)}^2,\end{equation}
which can be seen as a multi-dimensional analogue of Bessel's Inequality~\eqref{AJSa}.

Now, let~$\epsilon>0$ and (see e.g.~\cite[Theorem~9.6]{MR3381284}) pick a~$\Z^n$-periodic function~$f_\epsilon\in C^\infty(\R^n)$ such that~$\|f-f_\epsilon\|_{L^2((0,1)^n)}\le\epsilon$. 

We define~$g_\epsilon:=f-f_\epsilon$ and deduce from~\eqref{BEGH0-019340.1} that
$$ \sum_{{k\in \Z^n}}|\widehat g_{\epsilon,k}|^2\le\|g_\epsilon\|_{L^2((0,1)^n)}^2=\|f-f_\epsilon\|_{L^2((0,1)^n)}^2\le\epsilon^2.$$

For this reason and~\eqref{FCONJHALKDD-1},
\begin{eqnarray*}
&&\left\| \sum_{k\in Z} \widehat f_{\epsilon,k} \,e^{2\pi i k\cdot x}-\sum_{k\in Z} \widehat f_k \,e^{2\pi i k\cdot x}\right\|_{L^2((0,1)^n)}=\left\| \sum_{k\in Z} \widehat g_{\epsilon,k} \,e^{2\pi i k\cdot x}\right\|_{L^2((0,1)^n)}^2\le\epsilon^2
\end{eqnarray*}
and therefore
\begin{eqnarray*}
\left\| f-\sum_{k\in Z} \widehat f_k \,e^{2\pi i k\cdot x}\right\|_{L^2((0,1)^n)}&\le&
\|f-f_\epsilon\|_{L^2((0,1)^n)}+
\left\| f_\epsilon-\sum_{k\in Z} \widehat f_{\epsilon,k} \,e^{2\pi i k\cdot x}\right\|_{L^2((0,1)^n)}\\&&\qquad
+\left\| \sum_{k\in Z} \widehat f_{\epsilon,k} \,e^{2\pi i k\cdot x}-\sum_{k\in Z} \widehat f_k \,e^{2\pi i k\cdot x}\right\|_{L^2((0,1)^n)}\\
&\le&2\epsilon+\left\|f_\epsilon- \sum_{k\in Z} \widehat f_{\epsilon,k} \,e^{2\pi i k\cdot x}\right\|_{L^2((0,1)^n)}.
\end{eqnarray*}
Hence, we can take~$Z:={\mathcal{U}}_N$, send~$N\to+\infty$, recall~\eqref{PRECON56789-BIS.lo8}
(used here on the smooth function~$f_\epsilon$),
then send~$\epsilon\searrow0$, thus obtaining~\eqref{PARSENEN9c1903x}, as desired.

From~\eqref{PARSENEN9c1903x} and~\eqref{PARSENEN9c1903x2}, we also obtain~\eqref{PARSENEN}.
\end{proof}

In spite of the care required in general to deal with Fourier Series in higher dimension, let us mention that counterexamples typically carry over nicely from dimension~$1$ to~$n$. Namely, the one-dimensional examples provided in Section~\ref{EXCE} to highlight the lack of convergence for Fourier Series can be easily transferred to the $n$-dimensional case, just by constantly extending the given function in the other variables (i.e., by considering the function~$F(x_1,\dots,x_n):=f(x_1)$).

See also~\cite[Chapter~3]{MR3243734} for an approach to Fourier Series starting directly from dimension~$n$,
\cite[Chapter~79]{MR4404761} for a generalisation of the basic theory of Fourier Series from dimension~$1$ to~$n$,
and~\cite[Part~2]{MR1408905} for additional information.

Classical Fourier analysis can also be seen as a special case of the representation theory on homogeneous spaces, see e.g.~\cite{MR1143783}.

\begin{exercise}\label{SQDITJMS0x12S} The \emph{square Dirichlet Kernel} \index{square Dirichlet Kernel} is the function $$\R^n\ni x\longmapsto \sum_{{k=(k_1,\dots,k_n)\in\Z^n}\atop{|k_1|\le N,\,\dots,\,|k_n|\le N}}e^{-2\pi i k\cdot x}.$$
Prove that it can be written as
$$\prod_{j=1}^n D_N(x_j),$$ where~$D_N$ is the one-dimensional Dirichlet Kernel in~\eqref{PAKSw-L4}.\end{exercise}

\begin{exercise}\label{PRECON56789}
Let~$E$ be a bounded, measurable subset of~$\R^n$.

Consider the characteristic function of~$E$,
\begin{equation}\label{CHARAFA} \chi_E(x):=\begin{dcases}
1&{\mbox{ if }}x\in E,\\0&{\mbox{ if }}x\in\R^n\setminus E.
\end{dcases}\end{equation}
Let
$$ f(x):=\sum_{m\in\Z^n}\chi_E(2(x+m)).$$
Prove that the series above consists of a finite sum, hence~$f$ is well-defined.

Prove also that~$f$ is $\Z^n$-periodic and that, for every~$k\in\Z^n$,
$$\widehat f_k=\frac1{2^n}\int_{E}e^{-\pi ik\cdot x}\,dx.$$
\end{exercise}

\begin{exercise}\label{PRECON56789-BIS}
Let~$f$ be as in Exercise~\ref{PRECON56789}. Calculate~$\|f\|_{L^2((0,1)^n)}$.\end{exercise}

\begin{exercise}\label{PRECON56789-BIS.lo}Let~$\ell\in\N$.
Prove that there exists~$C>0$, depending only on~$n$ and~$\ell$, such that for all~$k=(k_1,\dots,k_n)\in\Z^n$
there exists~$\mu=(\mu_1,\dots,\mu_n)\in\N^n$ with~$\mu_1+\dots+\mu_n=\ell$ such that
$$ |k_1|^{\mu_1}\dots|k_n|^{\mu_n}\ge \big(k_1^2+\dots+k_n^2\big)^{\frac\ell2}.$$\end{exercise}


\chapter{Fourier Series: applications}\label{APPL-SEZ1}

The applications of Fourier Series are limitless.
Without aiming at being exhaustive, we list some of them here below.

\section{The old-fashioned epicycle theory}\label{ASTRP}

The Greeks, upon observing the night sky, formulated the hypothesis that all celestial motions could be explained as superpositions of uniform circular motions. This concept traces its origins back to figures like Apollonius and Hipparchus and was consolidated in Ptolemy's Almagest~\cite{MR1656361}. 

\begin{figure}[h]
\includegraphics[width=8cm]{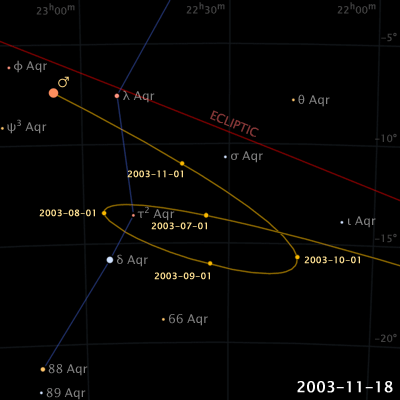}
\centering
\caption{Computer simulation of the night sky showing the apparent retrograde motion of the planet Mars in August and September of 2003 in the constellation Aquarius (work by Eugene Alvin Villar,
image from Wikipedia, CC BY-SA 4.0).}\label{RETP}
\end{figure}

This hypothesis was not only motivated by the elegance and simplicity of circular motions but it also turned out to be remarkably effective in making accurate predictions. The core of the theory is that the movements of heavenly bodies, at a first approximation, appear
as uniform circular orbits around the Earth's centre. That is, one could estimate the position of a celestial body by envisioning it moving uniformly along a circle with centre at the Earth. 

This average motion was \index{deferent}
called \emph{deferent}. However in reality, celestial bodies rarely occupy precisely this average position; rather, they slightly deviate from it, either moving ahead or lagging behind (see e.g. Figure~\ref{RETP}).

\begin{figure}[h]
\includegraphics[width=8cm]{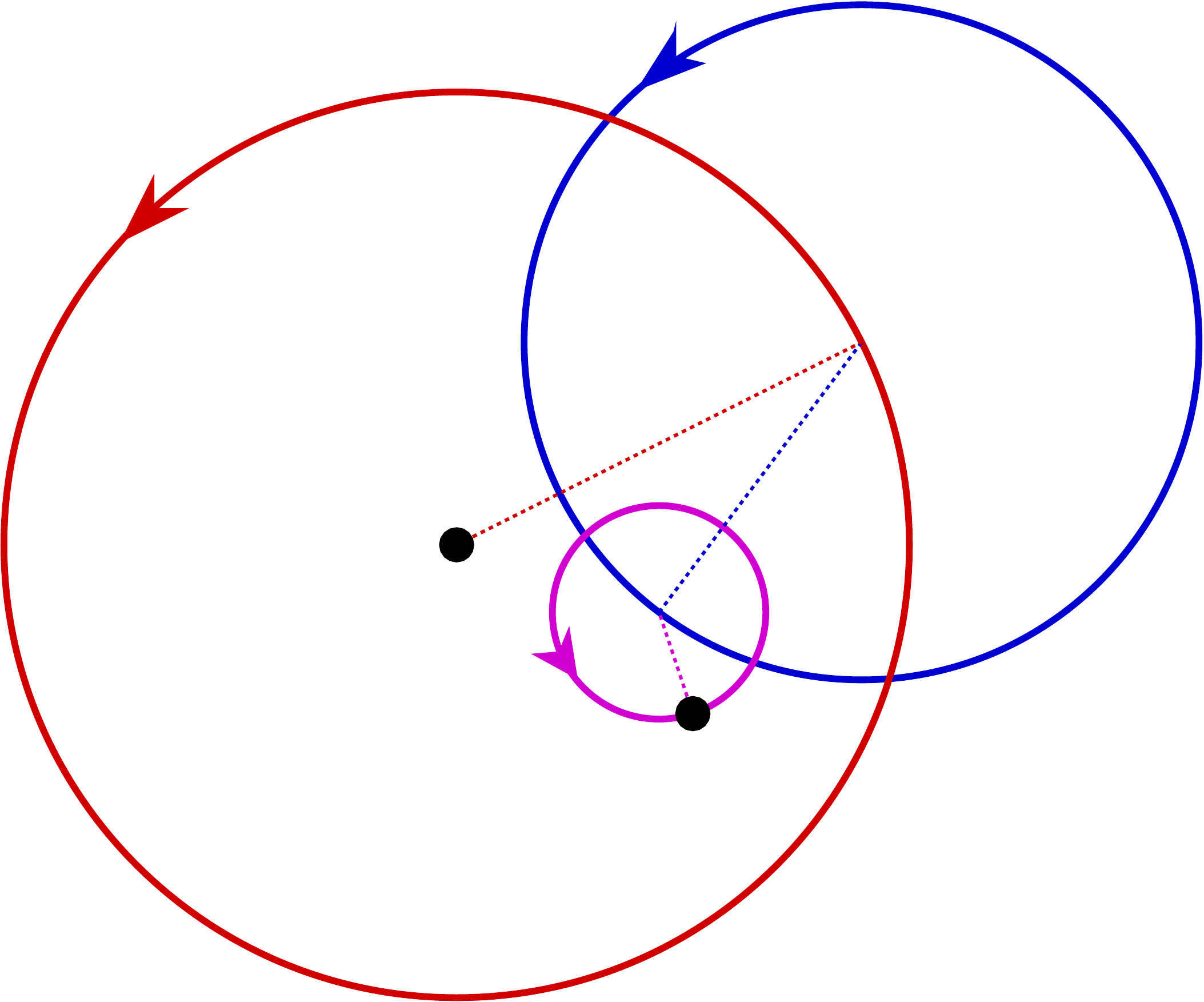}
\centering
\caption{The deferent motion (in red), a first epicycle (in blue), a second epicycle (in magenta).}\label{RETP2}
\end{figure}

To account for this discrepancy, the Greeks introduced a second circular motion, known as \index{epicycle}
\emph{epicycle}, centred at the point on the deferent where the celestial body was located. 
The idea of correcting circular motions by further circular motions in order to better fit the astronomical data can be repeated over and over, by adding a second epicycle centred at the first one, a third epicycle centred at the second one, and so on, see Figure~\ref{RETP2}.

\begin{figure}[h]
\includegraphics[height=0.35\textwidth]{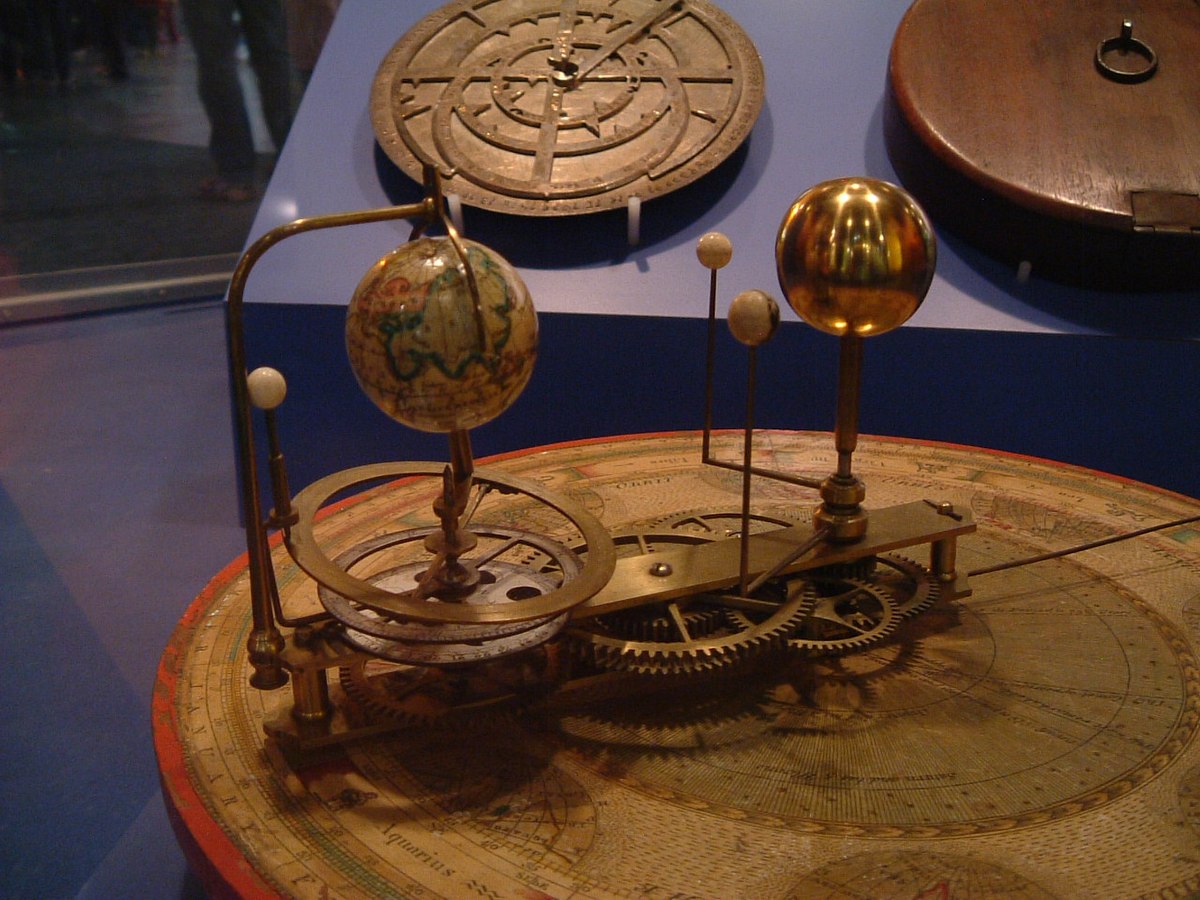}$\,$
\includegraphics[height=0.35\textwidth]{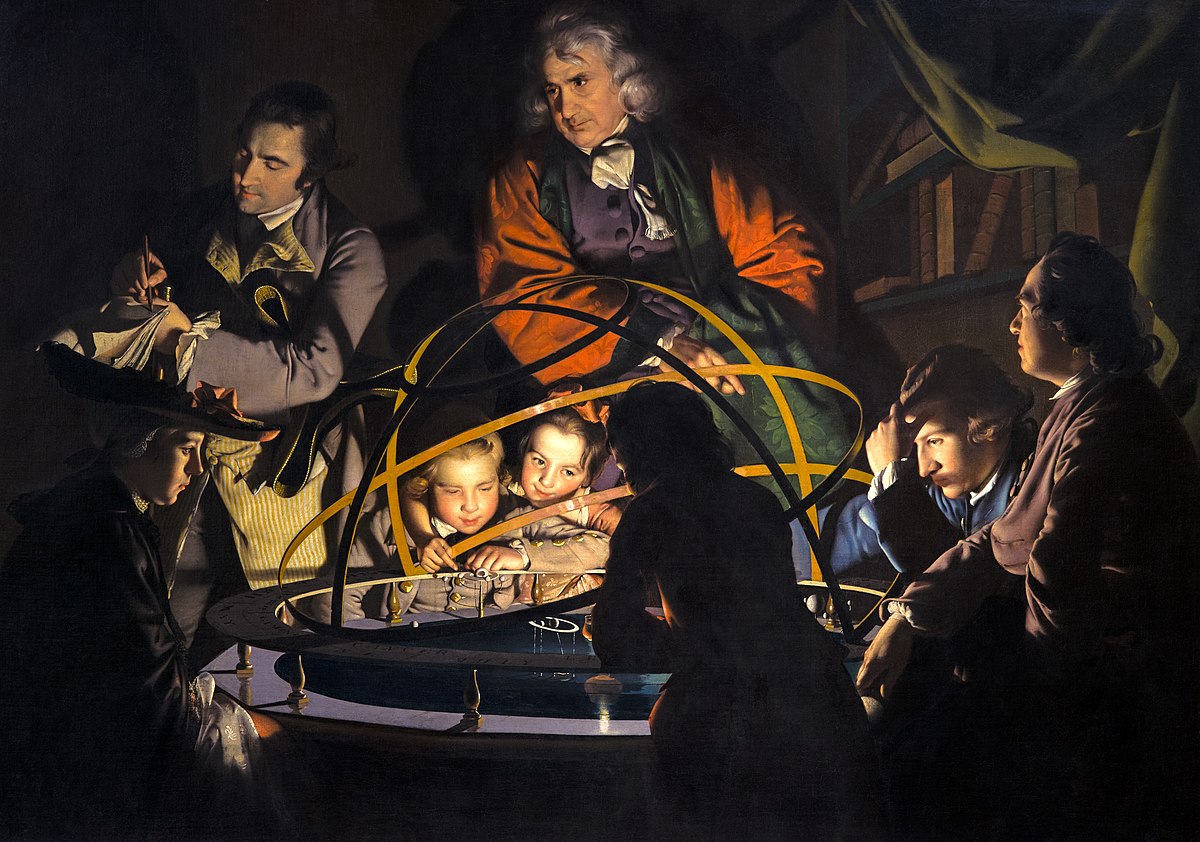}
\centering
\caption{Left, an orrery \index{orrery}
showing Earth and the inner planets (image from Wikipedia by Kaptain Kobold, CC BY 2.0). Right, ``A Philosopher Giving that Lecture on the Orrery, in which a Lamp is put in place of the Sun or The Orrery'', by artist Joseph Wright of Derby (public domain image from Wikipedia).}\label{RETPx}
\end{figure}

Of course, computing all these ``elementary'' circular motions and superimposing them at a time in which computers had not been invented yet was a tour-de-force, but the method turned out to be very fruitful. In fact, the reader may rightly be sceptical about the theory exposed so far, in which the Greeks placed the Earth, and not the Sun, at the centre of their reasoning, but:
\begin{itemize}
\item The system of deferents and epicycles obtained a great descriptive and predictive precision at that time, resulting much more accurate from the experimental point of view than other theories, such as the one by Aristarchus, that used to place the Sun, instead of the Earth, at the centre of the argument (indeed, Greeks also developed ``heliocentrism'', but the ``geocentric'' calculations at that time were far more in agreement with experimental measurements!),
\item The issue was likely not so much deciding what body has to be in the centre (the Sun is not at any centre of the universe after all, and maybe, 
according to contemporary observations, no matter how we try to define and identify it, there is simply nothing like the centre of the universe). The issue is more on describing the relative motion of the celestial bodies as seen from the Earth, since, in practice, this is the information that relates more directly to observations,
\item The system of epicycles could be constructed in practice with rather spectacular mechanical systems,
such as sets of levers and motors rotating at constant speed forming clockwork carrying the planets around
(these devices are called in the jargon ``orreries'', see Figure~\ref{RETPx}).
\end{itemize}

The reason for which we are discussing the epicycle theory in this set of notes is that there is a clear correspondence between it and the theory of Fourier Series. To see this, suppose that the deferent circle has radius~$R_0$ and that
the motion on it runs with an angular velocity~$\omega_0$, that is, if we identify the orbital plane with the complex plane~$\C$, the deferent motion takes the form~$R_0 e^{i\omega_0t}$. Similarly, suppose that the circular motion along the first epicycle takes the form~$R_1 e^{i\omega_1t}$, for some radius~$R_1$ and angular velocity~$\omega_1$, that the circular motion along the second epicycle takes the form~$R_2 e^{i\omega_2t}$, and so on.

Then, the motion of the celestial body obtained as a superposition of these circular motions can be written as
$$ R_0 e^{i\omega_0t}+R_1 e^{i\omega_1t}+R_2 e^{i\omega_2 t}+\dots$$
and we suppose, for simplicity, that a rational relation between the frequencies takes place, namely that
the frequency ratios~$\frac{\omega_i}{\omega_j}$ are rational\footnote{On the one hand, this is not an unreasonable assumption, since the rationals are a dense subset of the real numbers (i.e., the fact that every real number has rational numbers arbitrarily close to it suggests that one can assume the rationality of the frequency ratio up to an arbitrary degree of approximation). On the other hand, the rationals have zero Lebesgue measure, so they should somewhat appear with zero probability. This is indeed a tricky issue, and, in practical circumstances, it may be affected by complicated systems of resonances, tidal phenomena, etc. In fact, a good chunk of chaos theory stems out of this.}
numbers.

Thus, we write
$$\frac{\omega_1}{\omega_0}=\frac{k_1}{k_0},\qquad\frac{\omega_2}{\omega_0}=\frac{k_2}{k_0},\dots$$
for suitable integers~$k_0$, $k_1$, $k_2\dots$.

In this way, changing conveniently the units of measurements of time (i.e., looking at the rescaled time~$\tau:=\frac{k_0}{m_0} t$), the motion of a celestial body is described by
$$ R_0 e^{i k_0\tau}+R_1 e^{ik_1\tau}+R_2 e^{ik_2 \tau}+\dots$$
and the analogy with the Fourier Series is now apparent!

From this perspective, it is no wonder that the old-fashioned epicycle theory was so successful and accurate: after all,
the motion of most of the celestial bodies seen from Earth is (almost) periodic (e.g., the Moon takes about
27 days to complete a revolution around the Earth), and we know (e.g., from Theorem~\ref{DINITS})
that every ``reasonable'' periodic motion can be well-approximated by a Fourier expansion (whence, by a sufficiently large number
of epicycles). So, in some sense, the whole theory of epicycles can be considered as Fourier analysis ``ante litteram''
(yet, Fourier analysis is much more than that).

See~\cite{MR1898455} for a beautiful discussion about ancient and modern theories of celestial bodies' motions, their links to Fourier Series and chaos theory, and much more (see also the forthcoming Section~\ref{ODCABRBAGN} for another application of Fourier methods to celestial mechanics).

\begin{figure}[h]
\includegraphics[height=3.2cm]{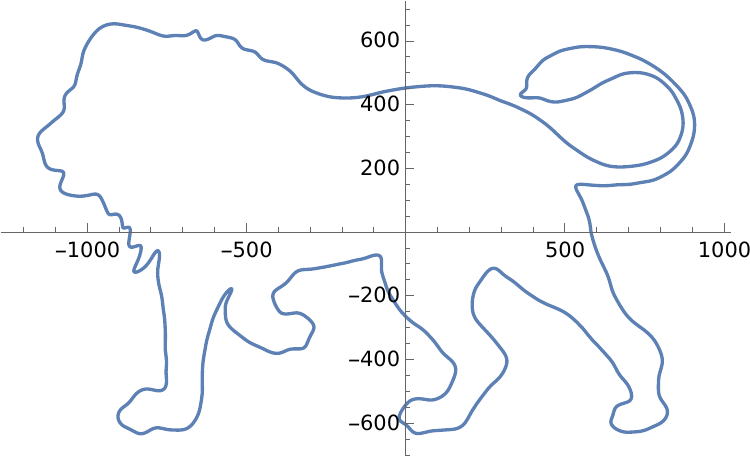}$\,$
\includegraphics[height=3.2cm]{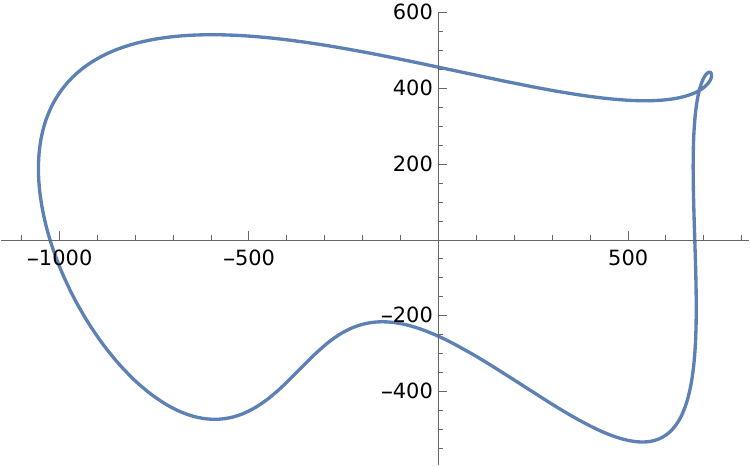}$\,$
\includegraphics[height=3.2cm]{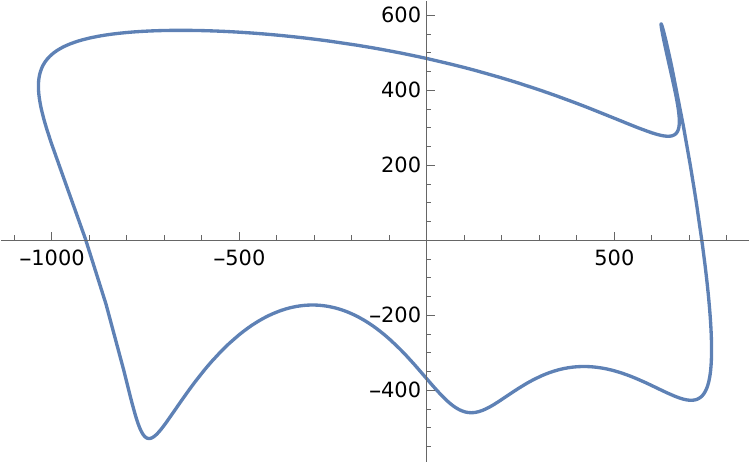}\\
\includegraphics[height=3.2cm]{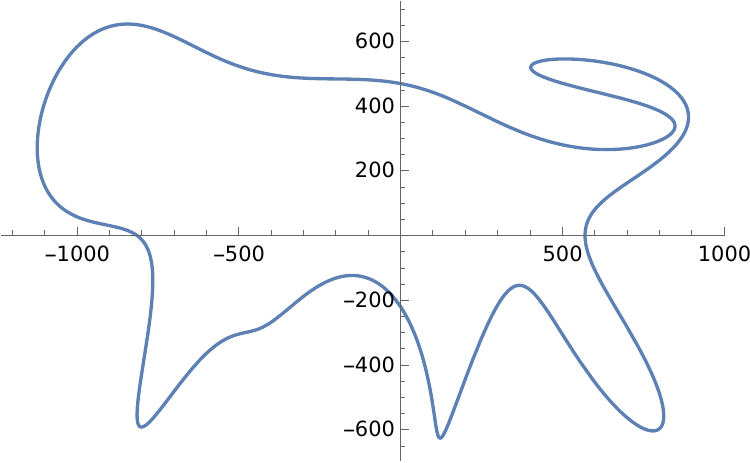}$\,$
\includegraphics[height=3.2cm]{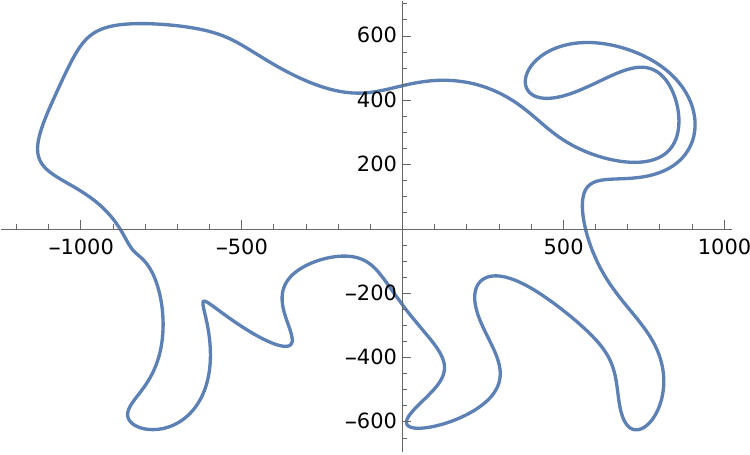}$\,$
\includegraphics[height=3.2cm]{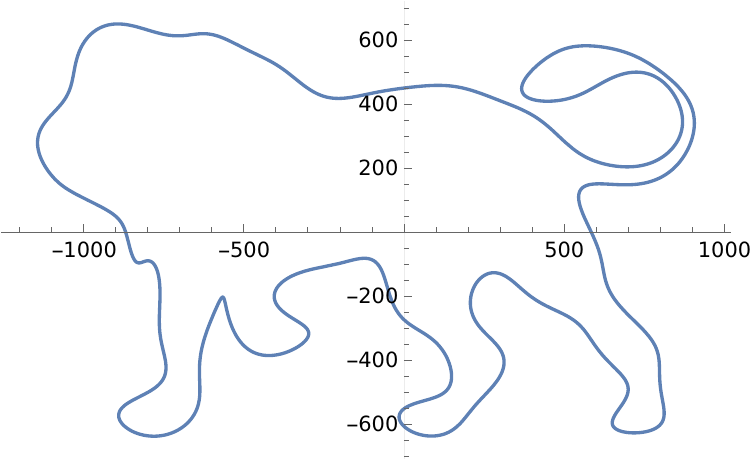}
\centering
\caption{The ``lion curve'' and its approximation by Fourier Series with frequencies up to~$3$, $5$, $10$, $20$, and~$40$.}\label{qTAF.fbuLLDOc3rcbD2.mET24Oikt.9dnsEU0}
\end{figure}

To confirm the capability of epicycles to virtually capture any periodic curve, one may play around with some Wolfram Alpha popular curves (see \protect\url{https://www.wolframalpha.com/examples/mathematics/geometry/curves-and-surfaces/popular-curves/} and check how well they can be approximated by epicycles, i.e. Fourier Sums.
To start with, we consider the ``lion curve'' (available on mathematica via the command {\texttt{Entity["PopularCurve", "LionCurve"][EntityProperty["PopularCurve", "ParametricEquations"]]}}). This curve is actually described by a sum of sines (with high frequencies and suitable phases) therefore it is a (rather long) Fourier Sum. Figure~\ref{qTAF.fbuLLDOc3rcbD2.mET24Oikt.9dnsEU0} shows the approximation of this curve by Fourier Sums with frequencies up to~$3$, $5$, $10$, $20$, and~$40$.

\begin{figure}[h]
\includegraphics[height=4.5cm]{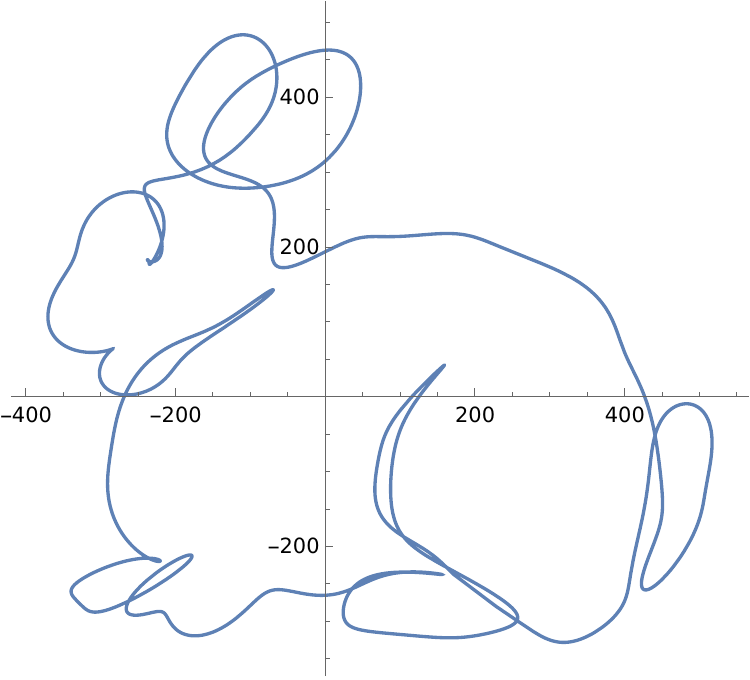}$\,$
\includegraphics[height=4.5cm]{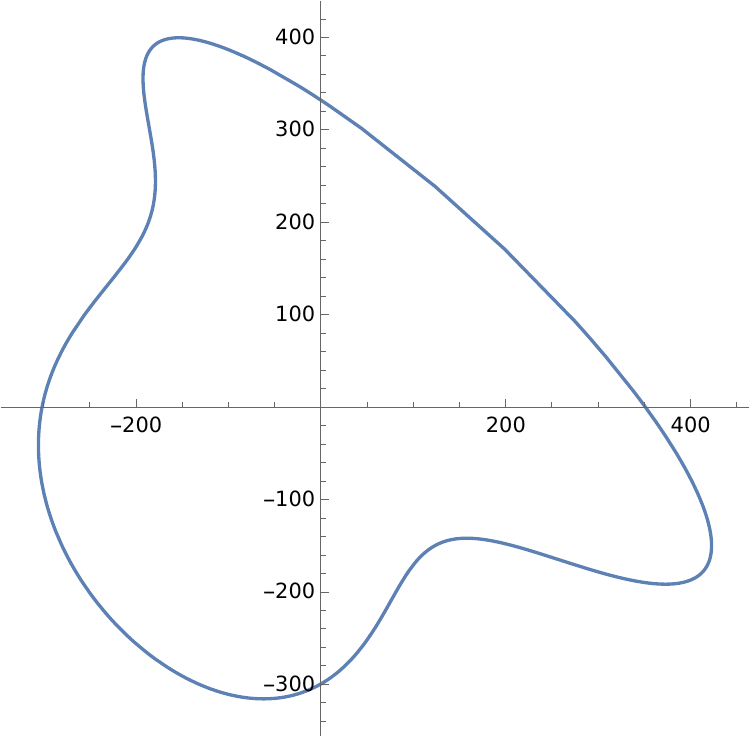}$\,$
\includegraphics[height=4.5cm]{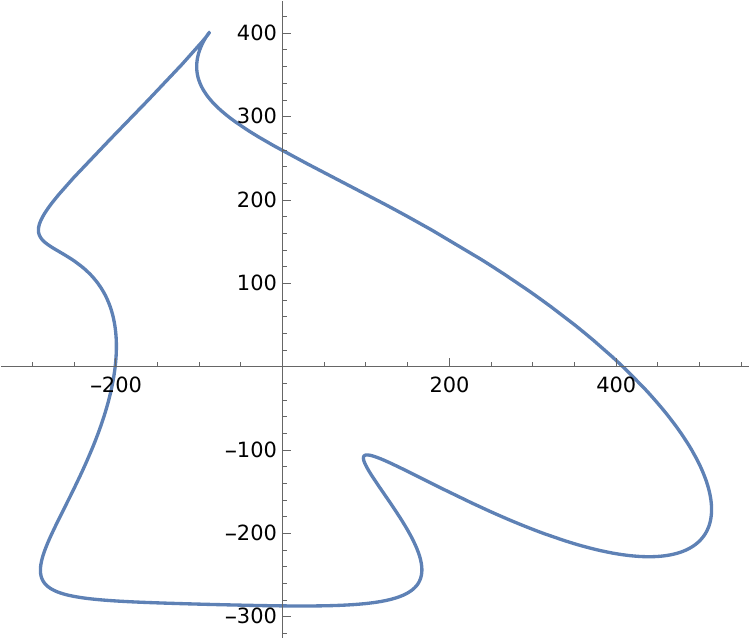}\\
\includegraphics[height=4.5cm]{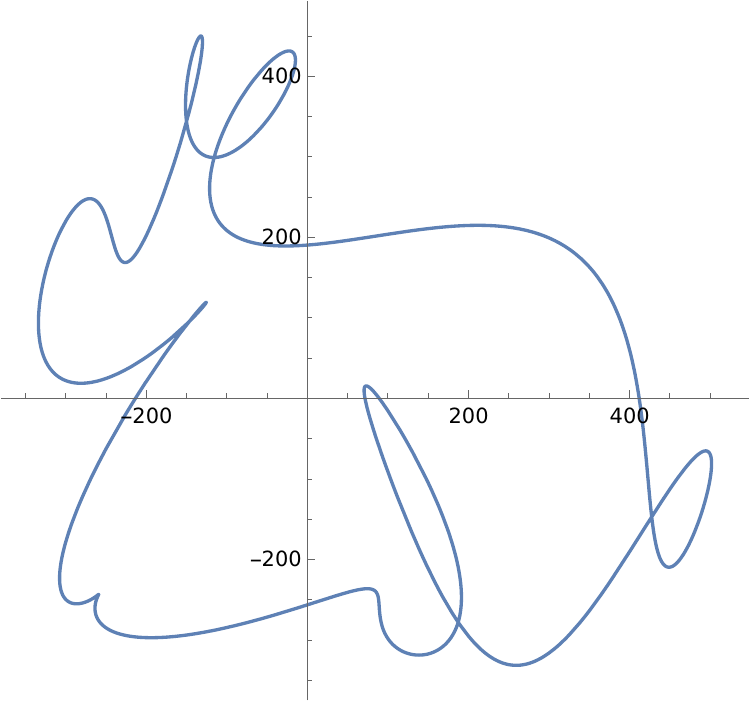}$\,$
\includegraphics[height=4.5cm]{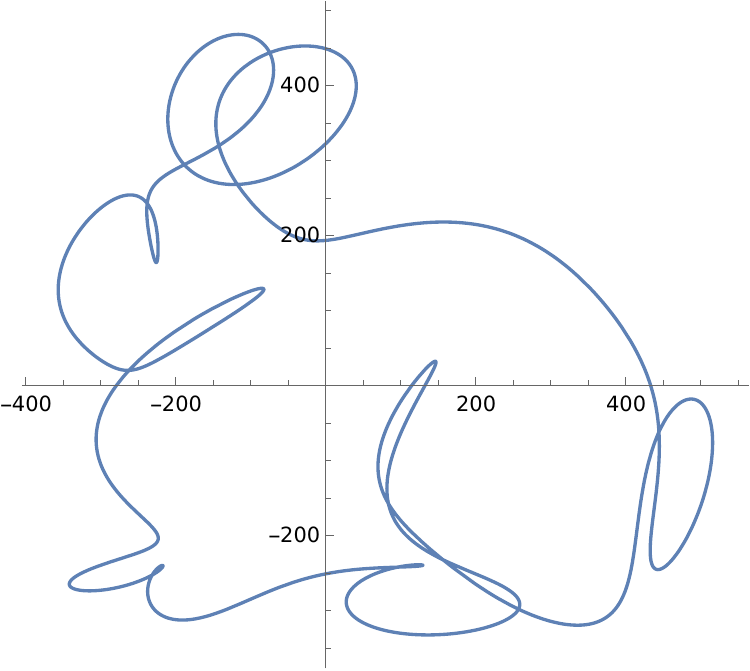}$\,$
\includegraphics[height=4.5cm]{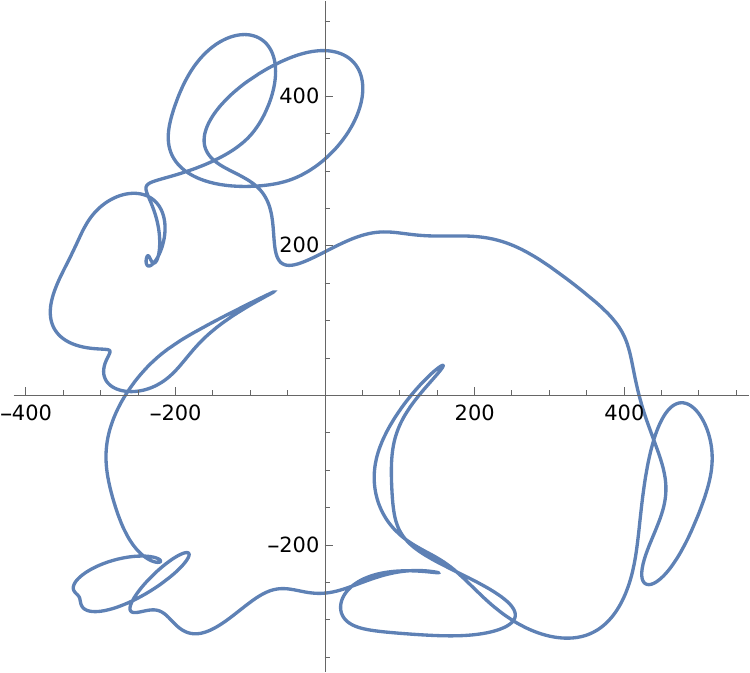}
\centering
\caption{The ``bunny curve'' and its approximation by Fourier Series with frequencies up to~$3$, $5$, $10$, $20$, and~$40$.}\label{qTAF.fbuLLDOc3rcbD2.mET24Oikt.9dnsEU02}
\end{figure}

The ``bunny curve'' is also a nice example of popular curve obtained by a finite sum of sines with high frequencies and phases (available on mathematica via the command {\texttt{Entity["PopularCurve", "BunnyCurve"]}}). Figure~\ref{qTAF.fbuLLDOc3rcbD2.mET24Oikt.9dnsEU02} shows the approximation of this curve by Fourier Sums with frequencies up to~$3$, $5$, $10$, $20$, and~$40$.

\section{Does anyone need a computer?}\label{ODCABRBAGN}

It would be desirable to frame Fourier methods in a way which makes it possible to implement them in computers. This is an important topic, better suited into the general framework of Fourier Transform, which requires suitable discretisation techniques (leading to the concept of \index{Discrete Fourier Transform} \emph{Discrete Fourier Transform}, or DFT) and smart ways of making the required calculations quicker to save computing time and optimise practical efficiency (leading to the concept of \index{Fast Fourier Transform} \emph{Fast Fourier Transform}, or FFT).

Here, we will not fully explore this territory and for a deeper insight we refer the reader to~\cite{MR735963}, \cite{VANHAL},
\cite[Chapter~7]{MR1970295}, \cite[Sections~3.8--3.10]{MR3616140},
\cite[Chapter~94--99]{MR4404761}, \cite[Chapter~IV]{PICARDELLO}, and the references therein.

\begin{figure}[h]
\includegraphics[height=6.2cm]{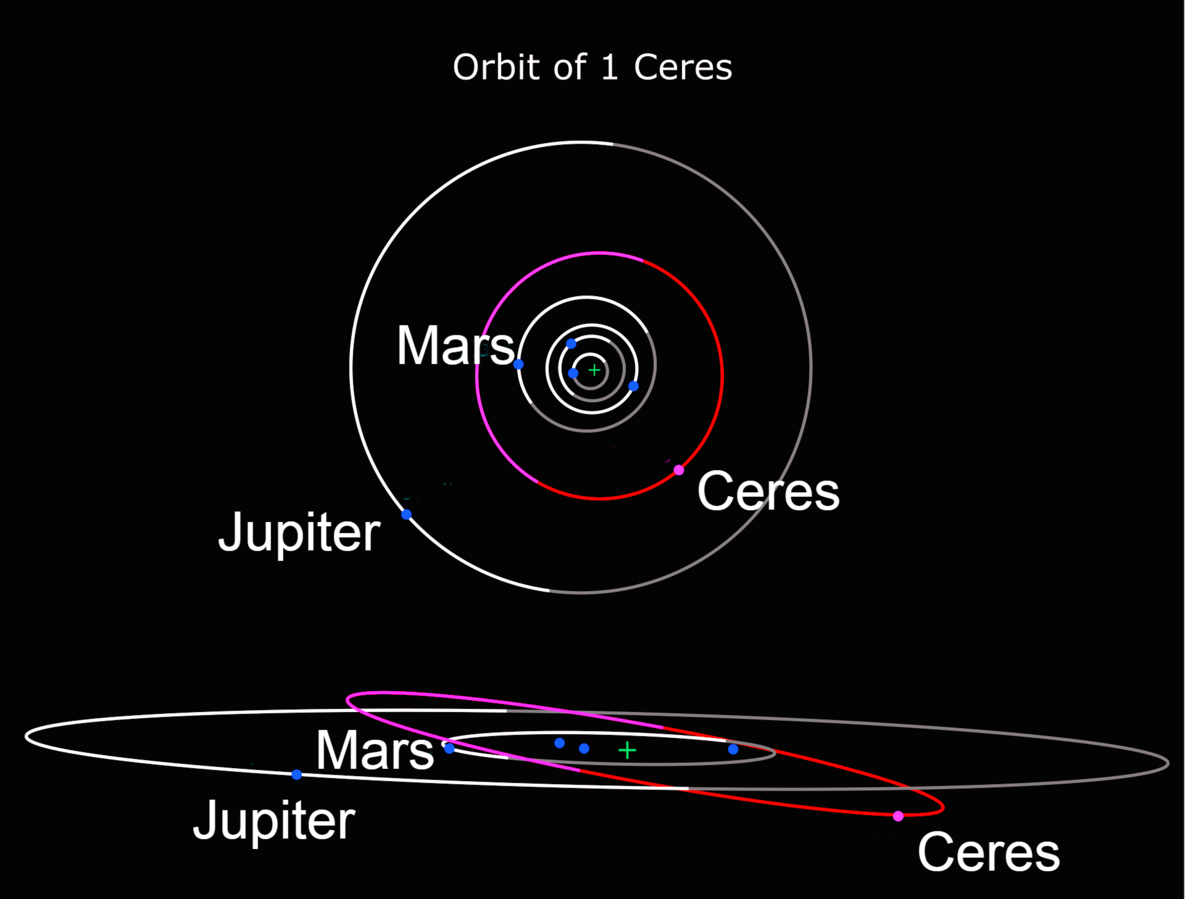}$\quad$\includegraphics[height=6.2cm]{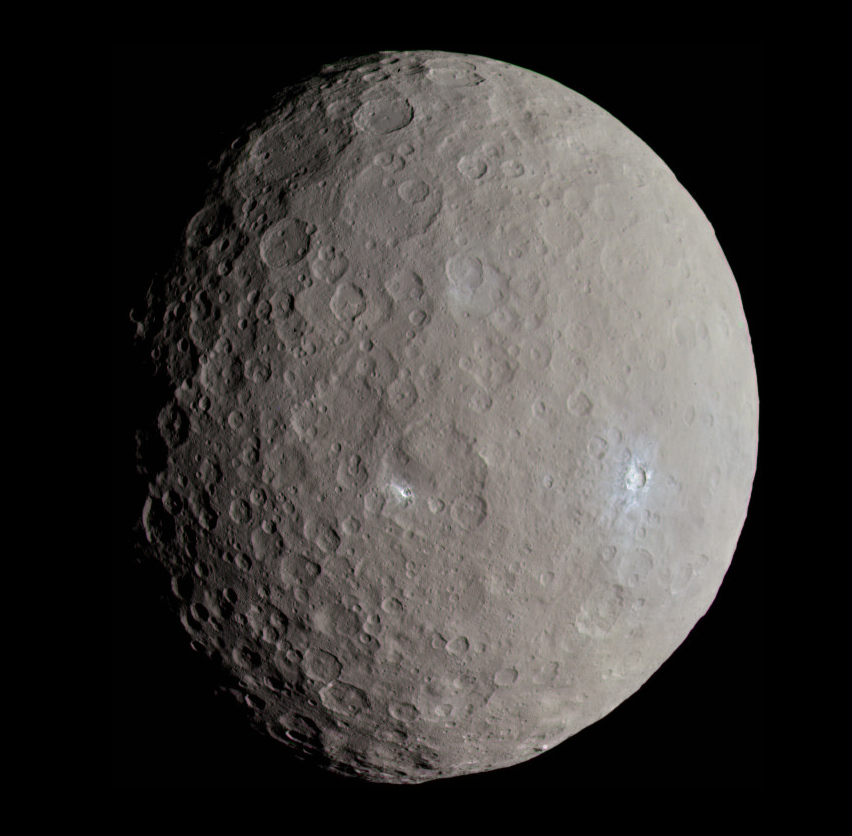}
\centering
\caption{Left: orbit of Ceres
(image from Wikipedia, based on work by Orionist and Amitchell, licensed under the Creative Commons Attribution-ShareAlike 4.0 License).
Right: image of Ceres in~2015 (by Justin Cowart, public domain image from
Wikipedia).}\label{c3.02poelf-1XrcbD2.mET24Oikt.9dnsEUvg50}
\end{figure}

\begin{figure}[h]
\includegraphics[height=6.4cm]{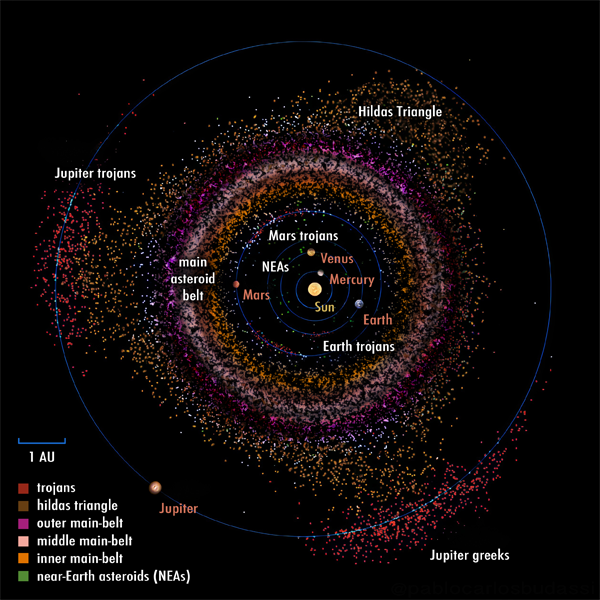}
\centering
\caption{A view of the asteroid belt (image from Wikipedia, based on work by Pablo Carlos Budassi, licensed under the Creative Commons Attribution-Share Alike 4.0 International License).}\label{AYS6yc3.02poelf-1XrcbD2.mET24Oikt.9dnsEUvg50}
\end{figure}

For the humble purposes of this book however it can be instructive to recall the first occurrence in which scientists have invented algorithms to suitably approximate Fourier coefficients by a method involving only a finite number of elementary operations such as additions and multiplications. Quite surprisingly, this predates Fourier (as well as modern computers) and dates back to the beginning of the nineteenth century, when some ``minor planets'' (or ``asteroids'', as we would say today) had been\footnote{The reason for which astronomers at that time started investigating consistently the region of space between Mars and Jupiter lies in the so-called \index{Titius-Bode Law} Titius-Bode Law,
according to which the semi-major axis of a planet orbit
should be~$0.4+0.3\cdot 2^{n}$ astronomical units with~$n\in\{-\infty,0,1,2,\dots\}$.

At that time, this heuristic formula returned reasonable results for Mercury (corresponding to~$n=-\infty$,
predicted distance by the Titius-Bode Law~$0.4$, 
actual distance about~$0.39$ astronomical units), Venus
(corresponding to~$n=0$,
predicted distance by the Titius-Bode Law~$0.7$, 
actual distance about~$0.72$ astronomical units), 
Earth
(corresponding to~$n=1$,
predicted and 
actual distance of~$1$ astronomical units, which is actually the definition of astronomical unit), Mars (corresponding to~$n=2$,
predicted distance by the Titius-Bode Law~$1.6$, 
actual distance about~$1.52$ astronomical units), Jupiter (corresponding to~$n=4$,
predicted and
actual distance about~$5.2$ astronomical units), and
Saturn (corresponding to~$n=5$,
predicted distance by the Titius-Bode Law~$10$, 
actual distance about~$9.54$ astronomical units).

The popularity of Titius-Bode Law further increased after the discovery of 
Uranus  in 1781 (which provided a rather good agreement with the law, with~$n=6$, corresponding
to a predicted distance~$19.6$, actual distance~$19.18$ astronomical units).

It is intriguing to conjecture about a ``missing planet'' that should have corresponded to the choice~$n=3$, or to find valid explanations for the absence
of such a planet: for instance, Immanuel Kant had conjectured that the gap between Mars and Jupiter
could have been created by 
Jupiter's gravity.

Ceres was discovered on 1 January 1801, by astronomer and mathematician (as well as priest of the Theatine order)
Giuseppe Piazzi. Announced at first as a new planet, Ceres was later reclassified as an asteroid and nowadays as a dwarf planet. The agreement with Titius-Bode Law was quite intriguing
(in spite of the fact that the orbit is slightly elongated, its
semi-major axis is about~$2.77$ astronomical units, to be compared with the
predicted value of~$2.8$).

To be perfectly honest, maybe the Titius-Bode Law was not really the main drive for Piazzi. It seems that Piazzi,
with the intent of
compiling a stellar catalogue,
was simply searching for an already known star, but found that this start ``was preceded by another", so science sometimes does rely on a good dose of serendipity. 

The name of Ceres was proposed by
Piazzi (in honour of the Roman goddess of agriculture, 
whose most ancient sanctuary was located in Sicily). Actually, the name proposed by Piazzi was Ceres Ferdinandea (the latter, in honour of King Ferdinand III of Sicily), but the name did not enjoy much success in other nations and therefore ended up to be dropped.

By the way, the brines of Ceres provide a potential habitat for microbial life, which is a classic theme for a science fiction novel that the reader is invited to invent.

Nowadays, we know that the region between Mars and Jupiter is occupied
by an ``asteroid belt'', lying between~$2.2$ and~$3.2$ astronomical units from the Sun, see Figure~\ref{AYS6yc3.02poelf-1XrcbD2.mET24Oikt.9dnsEUvg50}.
It is believed that the asteroid belt originated from a primordial group of planetesimals,
formed out of cosmic dust grains, which, due to the strong gravitational effect of Jupiter, did not manage to give rise to a full planet (so, maybe Kant was not completely wrong after all).

Till now, there is no solid theoretical explanation supporting the Titius-Bode Law, but it is often believed that laws of this type could arise as a consequence of orbital resonances: see~\cite[pages~88--91]{CELL1} and the references therein for more information.} observed.
The story begins with the observation
of the celestial body called Ceres, see Figure~\ref{c3.02poelf-1XrcbD2.mET24Oikt.9dnsEUvg50},
which was detected twenty-four times
between 1 January and 11 February 1801.
By this time, the apparent position of Ceres had changed,
due to Earth's motion around the Sun, becoming too close to the Sun to be observed.
Later on, Ceres was expected to become visible again, yet it was difficult to predict its exact position
after a long time with no direct observations.
While many astronomers tried in vain to relocate Ceres' position relying on telescopes, Gauss decided to use some mathematics for this goal. By December 1801, Gauss was able to predict the exact position of Ceres, allowing astronomers for their coveted observation.
Gauss repeated his success in 1802, with another asteroid, called Pallas.

Likely, the method of Gauss relied on previous work  by Euler, Clairaut, and Lagrange, with the idea of describing a periodic orbit of a celestial body as a trigonometric polynomial, focusing on astute techniques to recover the coefficients of this trigonometric polynomial (the Fourier coefficients, we would say nowadays) from experimental observations.

Concretely, one can think that a given periodic function~$f$, normalised to have period~$1$,
is approximated by (or is equal to, for simplicity) a trigonometric polynomial, thus writing
\begin{equation}\label{GACOC0N-1} f(t)=\sum_{k=0}^{N-1}C_k\,e^{2\pi ikt},\end{equation}
for some unknown complex coefficients~$C_0,\dots,C_{N-1}$.

One supposes that the values of~$f$ have been observed at time~$t_j:=\frac{j}{N}$, with~$j\in\{0,\dots,N-1\}$ and would like to recover the coefficients~$C_0,\dots,C_{N-1}$ from this information
(with this, all the values of~$f$ would be determined by~\eqref{GACOC0N-1}, the analogy with Ceres
being to be able to recover its position at any time~$t$ from previous observations at time~$t_0,\dots,t_{N-1}$).

To recover the desired coefficients, one could infer from the know experimental data in~\eqref{GACOC0N-1} that, for each~$j\in\{0,\dots,N-1\}$,
$$f(t_j)=\sum_{k=0}^{N-1}C_k\,e^{2\pi ikt_j}=\sum_{k=0}^{N-1}C_k\,e^{\frac{2\pi ikj}{N}}.
$$
Hence, for all~$h\in\{0,\dots,N-1\}$,
$$ f(t_j) \,e^{\frac{-2\pi ihj}{N}}=\sum_{k=0}^{N-1}C_k\,e^{\frac{2\pi i(k-h)j}{N}},$$
leading (see Exercise~\ref{GABNMSOIPEDCFQMw0-1}) to
$$ \sum_{j=0}^{N-1}f(t_j) \,e^{\frac{-2\pi ihj}{N}}=\sum_{k,j=0}^{N-1}C_k\,e^{\frac{2\pi i(k-h)j}{N}}=
N\sum_{k=0}^{N-1}C_k\,\delta_{k,h}=NC_h.$$

On that account, the coefficients in~\eqref{GACOC0N-1} are reconstructed by the formula,
for each~$h\in\{0,\dots,N-1\}$,
\begin{equation}\label{CGCHGA} C_h=\frac1N\sum_{j=0}^{N-1}f(t_j) \,e^{\frac{-2\pi ihj}{N}}.\end{equation}

But this is not the end of the story. This formula works fine, in principle, and would reconstruct all the values
of the function~$f$ via~\eqref{GACOC0N-1}. Fine. But the computational time required for the task would be rather long.
Indeed, to calculate each of the coefficients in~\eqref{CGCHGA}, one needs to perform~$N$ multiplications, then sum~$N$ objects, and then make a division. 

Namely, to compute all coefficients under consideration in~\eqref{CGCHGA}, \begin{equation}\label{CGCHGA.3}
{\mbox{$O(N^2)$ elementary operations are needed.}}
\end{equation}

Gauss spotted a remarkable simplification in~\eqref{CGCHGA} which dramatically reduced the computational effort. The gist is to consider the case in which~$N$ is a composite number (and most of the numbers are composite anyway, so having a simplification for them would be already good enough!), say when~$N=N_1\,N_2$, with~$N_1$, $N_2\in\N\cap[2,+\infty)$.
As we will see, in this case the simplification is the outcome of noticing that many multiplications return the number~$1$
and can thereby be disregarded. 

To see these cancellations, let us consider the index~$h\in\{0,\dots,N-1\}=\{0,\dots,N_1 N_2-1\}$ and let us perform a division by~$N_1$. In this way, we can write that
$$ h=N_1 h_1+h_2,$$
for some~$h_1=h_1(h)\in\N$ and~$h_2=h_2(h)\in\{0,\dots,N_1-1\}$.

We also observe that
$$ N_1N_2>N_1N_2 -1\ge h\ge N_1 h_1$$
and therefore~$h_1<N_2$, yielding that~$h_1\in\{0,\dots,N_2-1\}$.

Similarly, we can divide each index~$j\in\{0,\dots,N-1\}=\{0,\dots,N_1N_2-1\}$ in~\eqref{CGCHGA}
by~$N_2$ and obtain that
$$ j=N_2 j_1+j_2,$$
for some~$j_1=j_1(j)\in\N$ and~$j_2=j_2(j)\in\{0,\dots,N_2-1\}$.

In this situation,
$$ N_1N_2>N_1N_2 -1\ge j\ge N_2 j_1,$$
yielding that~$j_1<N_1$ and thus that~$j_1\in\{0,\dots,N_1-1\}$.

On that account, we can consider the~$N$th root of unity~$\omega:=e^{\frac{-2\pi i}N}$
(notice indeed that~$\omega^N= e^{-2\pi i}=1$) and,
when~$N=N_1 N_2$, rewrite~\eqref{CGCHGA} in the form
\begin{equation} \label{CGCHGA.2}\begin{split}
C_h&=\frac1N\sum_{j=0}^{N-1}f(t_j) \,e^{\frac{-2\pi ihj}{N}}\\
&=\frac1N\sum_{j=0}^{N-1}f(t_j) \,\omega^{hj}
\\&=\frac1N\sum_{j_1=0}^{N_1-1}\sum_{j_2=0}^{N_2-1} f(t_{N_2 j_1+j_2}) \,\omega^{(N_1 h_1+h_2)(N_2 j_1+j_2)}
\\&=\frac1N\sum_{j_1=0}^{N_1-1}\sum_{j_2=0}^{N_2-1} f(t_{N_2 j_1+j_2}) \,\omega^{N_1N_2 h_1 j_1}
\,\omega^{N_1 h_1j_2 + N_2 h_2 j_1+h_2j_2}\\&=\frac1N\sum_{j_1=0}^{N_1-1}\sum_{j_2=0}^{N_2-1} f(t_{N_2 j_1+j_2}) 
\,\omega^{N_1 h_1j_2 + N_2 h_2 j_1+h_2j_2}\\&=\frac1N\sum_{j_2=0}^{N_2-1}
\underbrace{\left(\sum_{j_1=0}^{N_1-1} f(t_{N_2 j_1+j_2}) 
\,\omega^{N_2 h_2 j_1}\right)}_{B_{j_2,h_2}}\,\omega^{h j_2}.
\end{split}\end{equation}
While not obvious at first glance, a significant speed-up is achieved from this two-step (double sum) calculation. The array of $B_{j_2,h_2}$ coefficients has $N_1N_2=N$ entries but each now requires only $N_1$ operations (inner sum), resulting in a complexity of $O(NN_1)$ to construct $B_{j_2,h_2}$. Then, the $N$ coefficients $C_h$ are evaluated from the $B_{j_2,h_2}$ coefficients but only requiring $N_2$ operations (outer sum). Overall, the algorithm thus requires
\begin{equation*}
{\mbox{only $O(N(N_1+N_2))$ elementary operations.}}
\end{equation*}
This is a significant improvement from~\eqref{CGCHGA.3} when~$N_1$ and~$N_2$ are large. The procedure can be repeated recursively to break $N_1$ and $N_2$ into smaller composite numbers. In cases where $N$ is "extremely composite", for example $N=2^M$ for $M\in \N$, such \emph{divide-and-conquer} strategy approaches a complexity of $O(N\log_2 N)$, which is a tremendous reduction compared to $O(N^2)$. 

Algorithms similar to the one originally introduced by Gauss were later independently rediscovered, see e.g.~\cite{MR178586}.
See also~\cite{MR484905} for a close comparison between the method of Gauss and similar ones introduced later on. The FFT algorithm is unanimously considered to be one of the most important in scientific computing.

\begin{exercise}\label{GABNMSOIPEDCFQMw0-1}
Let~$k$, $h\in\N\cap[0,N-1]$. Prove that
$$\frac1N\sum_{j=0}^{N-1} e^{\frac{2\pi ij(k-h)}N}=\delta_{k,h}.$$
\end{exercise}

\section{Predicting tides}\label{TIDES}

William Thomson, 1st Baron Kelvin, was a rather precocious kid. At the age of 15, he wrote
\emph{An Essay on the Figure of the Earth}, containing so advanced and original mathematical analysis to
be awarded a gold medal from the University of Glasgow. This early work was a fruitful companion
for Lord Kelvin throughout his life and it is commonly reported that, as a coping strategy during times of personal stress,
he used to come back to the problems raised in the essay to develop further studies on those topics.

\begin{figure}[h]
\includegraphics[height=0.4\textwidth]{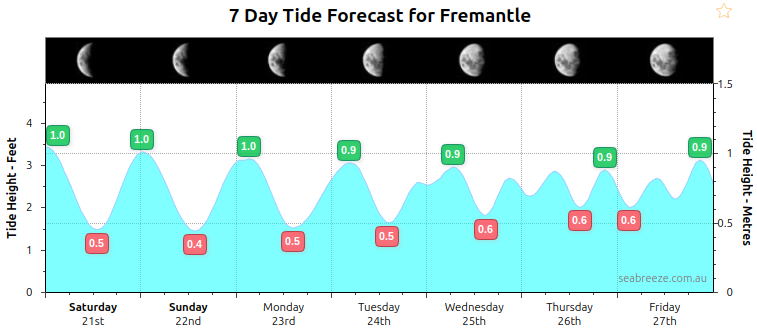}
\centering
\caption{Tide forecasts from \protect\url{https://www.seabreeze.com.au}.}\label{SEAo90pjHAhanaokdfRETPx}
\end{figure}

Interestingly, the essay opened with a quotation of Alexander Pope:

``Go, wondrous creature! mount where Science guides;

Go measure earth, weigh air, and state the tides''.

Predicting tides has indeed been a major scientific problem, and it is a typical example of very complex system,
being influenced by the simultaneous effects of the gravitational forces exerted by the Moon and the Sun, the rotation of the Earth, the seasonal changes, the shape of the coastlines, etc., and being also affected by transient meteorological effects such as changes of barometric pressure, wind, and storms. 

Besides Pope and tides, another important ingredient in Kelvin's training was played by Fourier analysis.
While the work of Fourier at the time was being heavily criticised by several scientists
(so, we can well remain serene when our own work gets criticised too), Kelvin felt enthusiastic
about Fourier's mathematical approach towards the physical world. 

Hence, when Kelvin got involved with problems related to navigation, not only he made good use of Fourier methods
in the prediction theory of tides, but also he designed and built several mechanical machines for this scope.

The theoretical idea is, after all, relatively simple: as a first approximation, one can describe the height~$h$
of the tide at a certain point at time~$t$ and postulate that it is influenced
by forces which change periodically (such as the gravitational effects of the Moon, Sun and Earth, and seasonal effects),
and also assume that the effect of  these forces is additive.
For instance, if~$h_1$ takes into account the effect of the Moon, $h_2$ of the Sun, and so on, we write
\begin{equation}\label{0kmws-MAjnedfimsT.1} h(t)=\sum_{j=1}^M h_j(t),\end{equation}
for some~$M\in\N$.

\begin{figure}[h]
\includegraphics[height=0.5\textwidth]{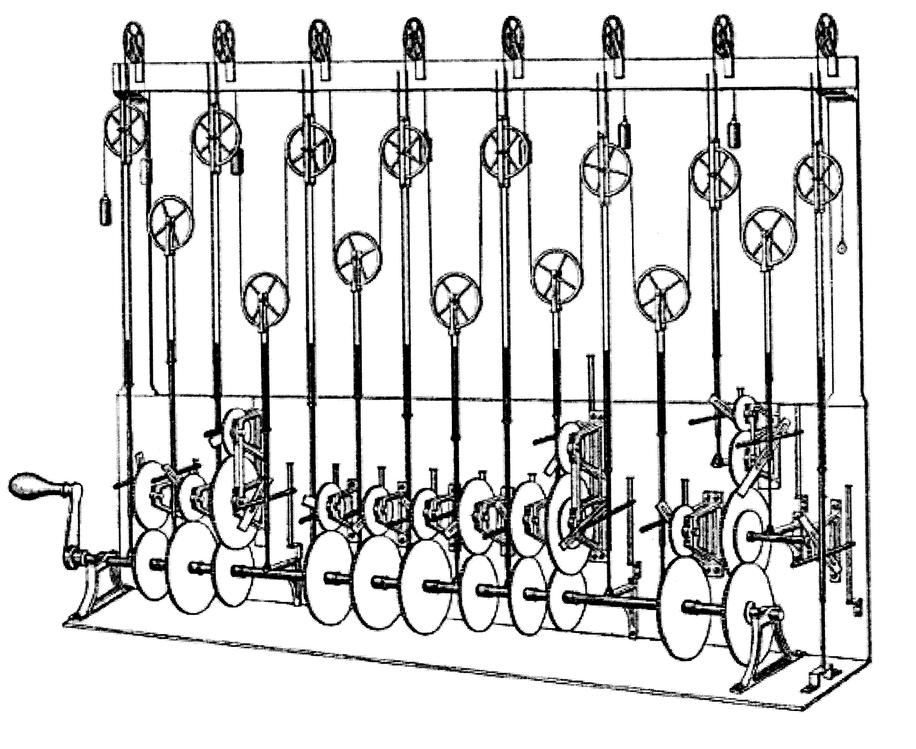}
\centering
\caption{Kelvin's
design sketch for a tide-predicting machine
(scanned by Terry0051 and released into the public domain).}\label{o90pjHAhanaokdfR134ETPx.22}
\end{figure}

If all these forces follow a periodic pattern, we can hope to approximate them with their Fourier Sum. Hence, for all~$j\in\N$, we assume that
\begin{equation}\label{0kmws-MAjnedfimsT.2} h_j(t)=\sum_{{k\in\Z}\atop{|k|\le N}} \widehat h_{j,k}\,e^{i\omega_k t},\end{equation}
for some~$N\in\N$ and
for a suitable frequency~$\omega_k$.

We stress that these frequencies are in general different one from another
(since, for instance, the period of rotation of the Earth with respect to the Moon is different from the period of rotation of the Earth with respect to the Sun, and the mismatch of frequency is certainly expected to cause a number of mathematical difficulties). The superpositions of effects with different periods is clearly visible in tide forecasts, see Figure~\ref{SEAo90pjHAhanaokdfRETPx}.

By formally inserting~\eqref{0kmws-MAjnedfimsT.2} into~\eqref{0kmws-MAjnedfimsT.1} we obtain
\begin{equation}\label{21er.TIDH} h(t)=\sum_{{1\le j\le M}\atop{{k\in\Z}\atop{|k|\le N}}} \widehat h_{j,k}\,e^{i\omega_k t}.\end{equation}

It is therefore highly desirable to:
\begin{itemize} \item determine the Fourier coefficients of~$h$,
\item in practice, given these coefficients, trace out efficiently the predicted height of the tide, \item in practice,
compute these coefficients from the experimental record of the past tidal height.\end{itemize}

\begin{figure}[h]
\includegraphics[height=0.45\textwidth]{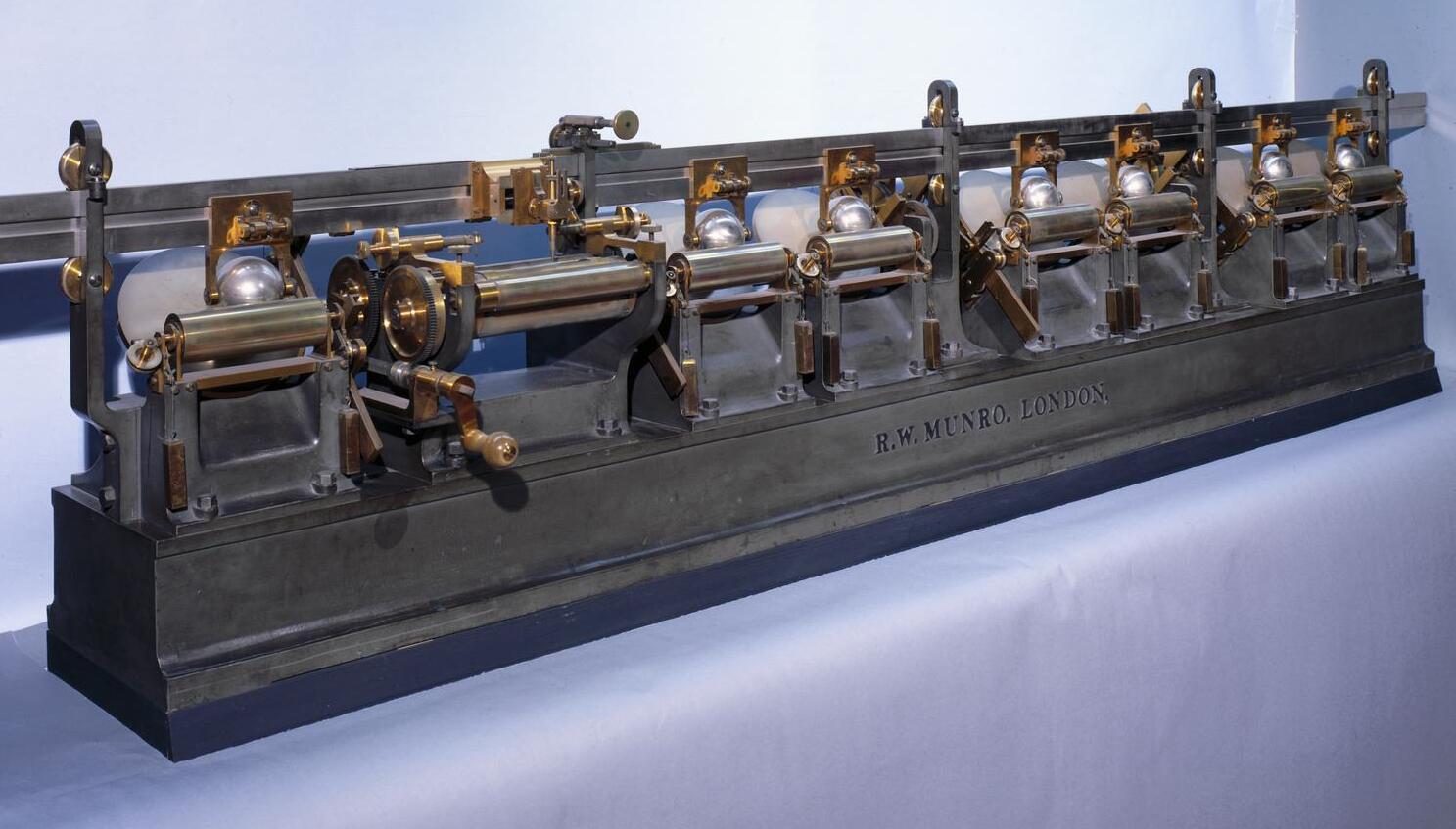}
\centering
\caption{Kelvin's harmonic analyser, assembled in 1878 in London.
Science Museum Group Collection \protect\url{https://www.sciencemuseumgroup.org.uk} (Creative Commons Attribution-NonCommercial-ShareAlike 4.0).}\label{o90pjHAhanaokdfRETPx}
\end{figure}

The first task, i.e. the determination of the Fourier coefficients, can be accomplished thanks to
an elegant mathematical formula (see Exercise~\ref{TIDEX}). Namely, if
$$ \widehat h_k:=\sum_{j=1}^M \widehat h_{j,k},$$
one can rewrite~\eqref{21er.TIDH} in the form
\begin{equation}\label{0kmwcvBs-MAjnedfimsT.2ghtaska0} h(t)=\sum_{{{k\in\Z}\atop{|k|\le N}}} \widehat h_{k}\,e^{i\omega_k t},\end{equation}
consider any initial time~$t_0$,
and, having a record of the tidal height~$h(t)$ over a long period of time~$[t_0, t_0 + T]$, the desired Fourier coefficients
can be obtained as a limit of a time average, as in the expression
$$\widehat h_{k}=\lim_{T\to+\infty}\frac1T\int_{t_0}^{t_0+T}h(t)\,e^{-i\omega_k t}\,dt.$$

The practical issue of a formula of this type is that it involves a long time limit, and in fact, to reach high accuracy,
the time length considered has to be quite large, given the substantial gap between the frequencies
involved (compare the fortnight lunar tidal period with half a solar year).
Hence, to accomplish the other tasks, namely to use the Fourier coefficients to predict the height of the tide
and effectively compute these coefficients from the experimental data, Kelvin invented
a variety of technological instruments specifically designed for this scope,
known under the various names of ball-and-disk integrators, \index{ball-and-disk integrator}
tide-predicting machines, \index{tide-predicting machine}
harmonic analysers, \index{harmonic analyser}
etc.

Specifically, plotting the vertical position on a curve on a moving band of paper to record
oscillatory motions required a new device, in which
a cord, anchored at one end, goes up over and down under a series of pulleys, which
move up and down,
either taking in or letting out the cord depending on how it moves,
and connecting via a system of gears to a drive shaft.
According to Kelvin, this rather sophisticated apparatus was
inspired to him by engineer Beauchamp Tower
(see Figure~\ref{o90pjHAhanaokdfR134ETPx.22} for its sketch).

\begin{figure}[h]
\includegraphics[height=0.45\textwidth]{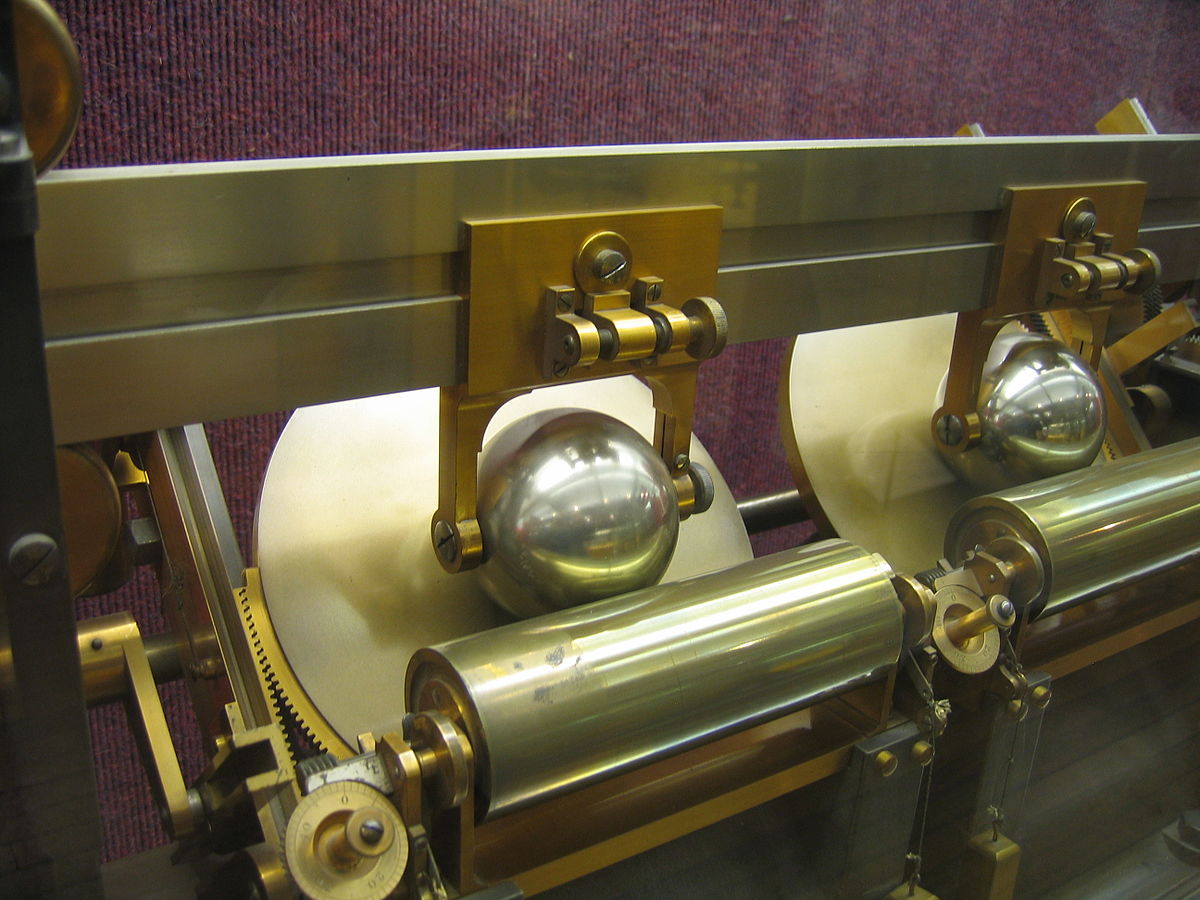}
\centering
\caption{Ball-and-disk integrators in Kelvin's tidal calculators,
showing the input spinning disk, the ball, and the cylindrical output shaft
(photo by Andy Dingley,
Creative Commons Attribution 3.0 Unported License).}\label{o90pjHAhanaokdfR134ETPx.2}
\end{figure}

As for the reconstruction of the Fourier modes,
the mechanism of the harmonic analyser (see Figure~\ref{o90pjHAhanaokdfRETPx})
traces a point along a curve, taking graphical records of daily changes in quantities of interest, whose oscillation causes many discs to rotate, connected to rolling spheres which communicate this motion to the recording rollers.
The selection of different frequencies was performed by moving the balls to the left or right along the rack
(see Figure~\ref{o90pjHAhanaokdfR134ETPx.2}): namely,
when the bearing is at the centre of the spinning disk, the ball does not impart any
motion to the output shaft, which remains stationary, but as the carriage moves the bearing away from the centre 
and towards the edge of the spinning disk, the output shaft rotates faster and faster.

As it often happens, inventions with a civil impact were turned into military devices: ball-and-disk integrators have been repeatedly used, especially before the advent of digitisation, in naval ballistics and missile systems, to efficiently compute the relative motion of the target.

\begin{exercise}\label{TIDEX} Let~$h$ be as in~\eqref{0kmwcvBs-MAjnedfimsT.2ghtaska0}, with~$\omega_k\ne\omega_m$ whenever~$k\ne m$.

Prove that, for all~$t_0\in\R$,
$$ \widehat h_{k}=\lim_{T\to+\infty}\frac1T\int_{t_0}^{t_0+T}h(t)\,e^{-i\omega_k t}\,dt.$$
\end{exercise}

\section{Partial differential equations}\label{PDESE}

A classical application of Fourier Series consists in finding explicit solutions for some important partial differential equations. The gist of the method is often to use\footnote{The method described here is often referred to with the name of \index{separation of variables} {\em separation of variables}. In spite broadly used, this technique is not exactly a mathematical theory, since its applicability in concrete cases is limited by a number of features, such as the structure of the equation under consideration, the domain, and the given data. See however~\cite{MR1567202} for a precise approach to this matter and its link with the notion of symmetry.}
separation of variables, and expansions in Fourier Series to be differentiated. Some care is needed to go through such a program, since some technicalities may arise in connection with the convergence properties of the infinite series involved, with the notion of boundary or initial conditions, and with the uniqueness of the solution for the equation under consideration.

Here, we do not aim at exhausting this very deep topic and we rather limit ourselves to some explicit examples (see e.g.~\cite{MR3380662}, \cite{MR3497072}, and~\cite{ELEM1} for an introduction to the partial differential equations considered here and a discussion about their importance in concrete applications).

We start with the \index{heat equation}
heat equation along a one-dimensional bar, for a given initial temperature, with endpoints held at a constant temperature.

\begin{theorem}\label{HEAT} Let~$\ell>0$, $T>0$, and~$f\in C([0,\ell])\cap C^1((0,\ell))$ be such that
\begin{equation}\label{OJSNILCESNpfa.EXPGF7.1}
f(0)=f(\ell)=0\end{equation}
and
\begin{equation}\label{OJSNILCESNpfa.EXPGF7.2}
\sup_{(0,\ell)}|f'|<+\infty.\end{equation}

Then, there exists one and only one continuous function~$u:[0,\ell]\times[0,T)\to\R$ such that
\begin{eqnarray}\label{OJSNILCESNpfa.EXPGF7}
&&{\mbox{for each~$t\in(0,T)$, we have that $u(\cdot,t)\in C^2((0,\ell))$, }}\\
\label{OJSNILCESNpfa.EXPGF8}&&{\mbox{and for each~$x\in(0,\ell)$, we have that $u(x,\cdot)\in C^1((0,T))$,}}\end{eqnarray} solving
$$ \begin{dcases}
\partial_t u(x,t)=\partial^2_xu(x,t) &\;{\mbox{ for all }}(x,t)\in(0,\ell)\times(0,T),\\
u(0,t)=u(\ell,t)=0&\;{\mbox{ for all }}t\in[0,T),\\ u(x,0)=f(x)
&\;{\mbox{ for all }}x\in[0,\ell].
\end{dcases}$$

Furthermore, the above solution~$u$ can be explicitly written as
\begin{equation}\label{OJSNILCESNpfa.EXPGF2} u(x,t)=\sum_{k=1}^{+\infty} b_k\,e^{-\frac{\pi^2 k^2 t}{\ell^2}} \,\sin\left(\frac{\pi kx}\ell\right),\end{equation}
where the coefficients~$b_k$ are obtained from the expansion
\begin{equation}\label{OJSNILCESNpfa.EXPGF}
f(x)=\sum_{k=1}^{+\infty} b_k\, \sin\left(\frac{\pi kx}\ell\right).
\end{equation}
\end{theorem}

Before proving Theorem~\ref{HEAT} let us point out that it is easy to guess the form that the solution~$u$ must take:
indeed, if, for a given~$t\in[0,+\infty)$, we could develop the function~$u(\cdot,t)$ in Fourier Series as
$$ u(x,t)=\sum_{k=1}^{+\infty} \zeta_k(t)\,\sin\left(\frac{\pi kx}\ell\right),$$
and if we were allowed to freely take derivatives inside the summation sign, we would find that
$$\sum_{k=1}^{+\infty} \partial_t\zeta_k(t)\,\sin\left(\frac{\pi kx}\ell\right)=\partial_tu(x,t)=\partial^2_xu(x,t)=
-\sum_{k=1}^{+\infty} \frac{\pi^2 k^2\,\zeta_k(t)}{\ell^2}\,\sin\left(\frac{\pi kx}\ell\right),$$
which could be solved by taking, for all~$k\in\N\cap[1,+\infty)$,
$$\frac{d}{dt}\zeta_k(t)=- \frac{\pi^2 k^2\,\zeta_k(t)}{\ell^2},$$
leading to
$$\zeta_k(t)=\mu_k\,e^{-\frac{\pi^2 k^2 t}{\ell^2}},$$
for some~$\mu_k\in\R$ to be determined.

Now, formally writing that
$$ \sum_{k=1}^{+\infty} b_k\, \sin\left(\frac{\pi kx}\ell\right)=f(x)=u(x,0)=\sum_{k=1}^{+\infty} \zeta_k(0)\,\sin\left(\frac{\pi kx}\ell\right)=\sum_{k=1}^{+\infty} \mu_k\,\sin\left(\frac{\pi kx}\ell\right),$$
it is tempting to choose~$\mu_k:=b_k$.

These heuristics lead to the expression in~\eqref{OJSNILCESNpfa.EXPGF2}, showing that a function~$u$ in that form
is a reasonable candidate for a solution.

But are these steps rigorously justifiable? And if so, is the solution constructed in this way the only possible one?
Theorem~\ref{HEAT} indeed gives a positive answer to this question,
and now we provide a proof for it.

\begin{proof}[Proof of Theorem~\ref{HEAT}] First of all, we consider the odd reflection of~$f$ and we then extend this function periodically:
in this way, owing to~\eqref{OJSNILCESNpfa.EXPGF7.1} and~\eqref{OJSNILCESNpfa.EXPGF7.2}
we have obtained a function, still denoted by~$f$ for simplicity, which is globally Lipschitz continuous and periodic of period~$2\ell$. 

Therefore (see 
Theorem~\ref{C1uni}, Exercises~\ref{920-334PKSXu9o2fgfbsmos}, \ref{DINI1UNI}, and~\ref{DINI2UNI}, as well as~\eqref{fc:PkT-0.55}) we can expand~$f$
in Fourier Series, with uniform convergence, as in~\eqref{OJSNILCESNpfa.EXPGF}.

Now we check that defining~$u$ as in~\eqref{OJSNILCESNpfa.EXPGF2} one indeed finds a solution of the problem under consideration. To this end, by the Riemann-Lebesgue Lemma (see Theorem~\ref{RLjoqwskcdc}) we know that~$b_k\to0$ as~$k\to+\infty$, and therefore~$|b_k|\le C$ for some~$C>0$.

On this account, for all~$t\in(0,T)$, we obtain that
\begin{equation}\label{OJSNILCESNpfa.EXPGF9} |b_k|\,e^{-\frac{\pi^2 k^2 t}{\ell^2}} \le C\,e^{-\frac{\pi^2 k^2 t}{\ell^2}}. \end{equation}
From this, by the uniform convergence of the series of the derivatives involved
(having just an extra factor of the order~$k^2$ thrown
in, which would be hit anyway by the exponential decay), we obtain~\eqref{OJSNILCESNpfa.EXPGF7}
and~\eqref{OJSNILCESNpfa.EXPGF8}, and also that
termwise differentiation of the series is possible, leading to
\begin{eqnarray*}
&&\partial_t u(x,t)-\partial^2_xu(x,t)=\sum_{k=1}^{+\infty}
\partial_t\left(
b_k\,e^{-\frac{\pi^2 k^2 t}{\ell^2}} \,\sin\left(\frac{\pi kx}\ell\right)\right)-\sum_{k=1}^{+\infty}
\partial_x^2\left(
b_k\,e^{-\frac{\pi^2 k^2 t}{\ell^2}} \,\sin\left(\frac{\pi kx}\ell\right)\right)=0.
\end{eqnarray*}

We also infer from~\eqref{OJSNILCESNpfa.EXPGF9} that, for all~$\tau\in(0,T)$,
$$ \lim_{[0,\ell]\ni [0,T)\ni (x,t)\to(0,\tau)}u(x,t)=
\sum_{k=1}^{+\infty}
\lim_{[0,\ell]\ni [0,T)\ni (x,t)\to(0,\tau)}
b_k\,e^{-\frac{\pi^2 k^2 t}{\ell^2}} \,\sin\left(\frac{\pi kx}\ell\right)=0$$
and similarly
$$ \lim_{[0,\ell]\ni [0,T)\ni (x,t)\to(\ell,\tau)}u(x,t)=0.$$

To finish checking that~$u$ satisfies all the desired properties, we need to show that, for all~$\eta\in[0,\ell]$,
\begin{equation}\label{OJSNILCESNpfa.EXPGFC}
\lim_{[0,\ell]\ni [0,T)\ni (x,t)\to(\eta,0)}u(x,t)=f(x).
\end{equation}
To this end, we point out that, 
in light of Exercise~\ref{G0ilMajsx912e-14},
$$ \sum_{k=1}^{+\infty}|b_k|<+\infty.$$
As a consequence, given~$\epsilon>0$, we pick~$N_\epsilon\in\N$ such that
$$ \sum_{k=N_\epsilon+1}^{+\infty}|b_k|<\epsilon$$
and we see that
\begin{eqnarray*}&&
\lim_{[0,\ell]\ni [0,T)\ni (x,t)\to(\eta,0)}\big|u(x,t)-f(x)\big|\\&&\qquad=
\lim_{[0,\ell]\ni [0,T)\ni (x,t)\to(\eta,0)}\left|
\sum_{k=1}^{+\infty} b_k\,e^{-\frac{\pi^2 k^2 t}{\ell^2}} \,\sin\left(\frac{\pi kx}\ell\right)-
\sum_{k=1}^{+\infty} b_k \,\sin\left(\frac{\pi kx}\ell\right)\right|\\
&&\qquad\le\lim_{[0,\ell]\ni [0,T)\ni (x,t)\to(\eta,0)}
\left[
\left|
\sum_{k=1}^{N_\epsilon} b_k\,e^{-\frac{\pi^2 k^2 t}{\ell^2}} \,\sin\left(\frac{\pi kx}\ell\right)-
\sum_{k=1}^{N_\epsilon} b_k \,\sin\left(\frac{\pi kx}\ell\right)\right|+2\sum_{N_\epsilon+1}^{+\infty}| b_k|\right]\\
&&\qquad=2\sum_{N_\epsilon+1}^{+\infty}| b_k|\\&&\qquad\le2\epsilon.
\end{eqnarray*}
By taking the limit as~$\epsilon\searrow0$, we thereby complete the proof of~\eqref{OJSNILCESNpfa.EXPGFC}
and we have thus verified that~$u$
satisfies all the desired properties to be a solution of the problem under consideration.

We also observe that the solution that we have just considered is the only possible one:
indeed if two solutions~$u$ and~$v$ are given, one could define~$w:=u-v$ and check that~$w$ must necessarily vanish identically (see Exercise~\ref{UNIZER}).
\end{proof}

See Sections~\ref{ILVCEMDlamcu} and~\ref{POISSONKERN-ex4}, as well as
Exercise~\ref{BESSEL-FC-EX4-DIFFB-HYPB}, for other links between Fourier methods and partial differential equations.
See also~\cite[Chapter~4]{MR1145236}, \cite[Chapters~2, 3, and~5]{MR3497072},
and the references therein for further instructive readings.

\begin{exercise}\label{UNIZER} Let~$\ell$, $T>0$.
Prove that the only continuous function~$w:[0,\ell]\times[0,T)\to\R$ such that
\begin{eqnarray*}
&&{\mbox{for each~$t\in(0,T)$, we have that $w(\cdot,t)\in C^2((0,\ell))$, }}\\
&&{\mbox{and for each~$x\in(0,\ell)$, we have that $w(x,\cdot)\in C^1((0,T))$,}}\end{eqnarray*} solving
$$ \begin{dcases}
\partial_t w(x,t)=\partial^2_x w(x,t) &\;{\mbox{ for all }}(x,t)\in(0,\ell)\times(0,T),\\
w(0,t)=w(\ell,t)=0&\;{\mbox{ for all }}t\in[0,+\infty),\\ w(x,0)=0
&\;{\mbox{ for all }}x\in[0,\ell],
\end{dcases}$$
is the function that vanishes identically.
\end{exercise}

\begin{exercise} Let~$\ell$, $T>0$. Can one apply, under ``reasonable assumptions'' on the data, Fourier methods to solve the \index{wave equation}
wave equation
$$ \begin{dcases}
\partial_t^2 u(x,t)=\partial^2_x u(x,t) &\;{\mbox{ for all }}(x,t)\in(0,\ell)\times(0,T),\\
u(0,t)=u(\ell,t)=0&\;{\mbox{ for all }}t\in[0,+\infty),\\ u(x,0)=f(x)
&\;{\mbox{ for all }}x\in[0,\ell],\\ \partial_tu(x,0)=g(x)
&\;{\mbox{ for all }}x\in[0,\ell]?
\end{dcases}$$
\end{exercise}

\begin{exercise} Let~$a_1>b_1$, $\dots$, $a_n>b_n$.
Can one apply, under ``reasonable assumptions'' on the data, Fourier methods to solve the 
\index{elastic membrane equation}
elastic membrane equation
$$ \begin{dcases}
\partial^2_{x_1} u(x)+\dots+\partial^2_{x_n} u(x)=f(x) &\;{\mbox{ for all }}x=(x_1,\dots,x_n)\in{\mathcal{R}},\\
u(x)=0&\;{\mbox{ for all }}x\in\partial{\mathcal{R}},
\end{dcases}$$
where~${\mathcal{R}}:=(a_1,b_1)\times\dots\times(a_n,b_n)$?
\end{exercise}

\begin{exercise} Let~$\ell$, $T>0$. Can one apply, under ``reasonable assumptions'' on the data, Fourier methods to solve the \index{telegraph equation}
telegraph equation
$$ \begin{dcases}
\partial_t^2 u(x,t)+2\partial_tu(x,t)+u(x,t)=\partial^2_x u(x,t) &\;{\mbox{ for all }}(x,t)\in(0,\ell)\times(0,T),\\
u(0,t)=u(\ell,t)=0&\;{\mbox{ for all }}t\in[0,+\infty),\\ u(x,0)=f(x)
&\;{\mbox{ for all }}x\in[0,\ell],\\ \partial_tu(x,0)=g(x)
&\;{\mbox{ for all }}x\in[0,\ell]?
\end{dcases}$$
\end{exercise}

\begin{exercise}\label{NSSEEV}
The velocity field~$v:\R^3\times\R\to\R^3$ of an incompressible homogeneous fluid is often described by the \index{Navier-Stokes equation}
Navier-Stokes equation\footnote{When~$\nu=0$, \eqref{NS6THTgBplaISw} is called the \index{Euler equation}
Euler equation.}
\begin{equation}\label{NS6THTgBplaISw} 
\begin{dcases} \frac{\partial v }{\partial t}+(v \cdot \nabla )v =\nu \,\Delta v -\nabla p,\\ \nabla\cdot v=0
,\end{dcases}\end{equation}
where~$p:\R^3\to\R$ is the pressure and~$\nu\in[0,+\infty)$ is the viscosity of the fluid
(the first equation above accounting for momentum conservation, the second for fluid incompressibility).

Assume that~$v=v(x,t)$ and~$p=p(x,t)$ are~$\Z^3$-periodic in the space variable~$x$ (according to the notation introduced in Section~\ref{SEC:ANUYDIMEAGa}) and obtain formal\footnote{The convergence issues related to the Navier-Stokes equation are deep and highly problematic, therefore we do not aim at making Exercise~\ref{NSSEEV} anything more than a heuristic calculation to train with Fourier methods and appreciate their power even in situations where rigorous justifications are lacking.}
relations among their Fourier coefficients (which depend on the time variable~$t$).
\end{exercise}

\begin{exercise}\label{2NSSEEV}
In the notation of Exercise~\ref{NSSEEV}, the quantity
\begin{equation}\label{WTTGfdgtchofrTLA:0-2} {\mathcal{V}}(t):=\int_{(0,1)^3} |\nabla\times v(x,t)|^2\,dx\end{equation}
is called the \emph{enstrophy} \index{enstrophy}
and it accounts for the \emph{vorticity} \index{vorticity} $\nabla\times v$ of the fluid.

One can also consider the rate of fluid deformation
\begin{equation}\label{WTTGfdgtchofrTLA:0-2.9rjHNX} {\mathcal{D}}(t):=\sum_{j,m=1}^3\int_{(0,1)^3} \left|\frac{\partial v_m}{\partial x_j}(x,t)\right|^2\,dx.\end{equation}

Calculate, formally and without attempting to take care of any convergence issue, the Fourier Series of these quantities.
\end{exercise}

\begin{exercise}\label{2NSSEEV.bis}
In the notation of Exercise~\ref{NSSEEV}, the quantity
\begin{equation}\label{WTTGfdgtchofrTLA:0-1} {\mathcal{E}}(t):=\frac12\int_{(0,1)^3} |v(x,t)|^2\,dx\end{equation}
accounts for the total kinetic energy of the fluid.

Calculate, formally and without attempting to take care of any convergence issue, the derivative of the energy in Fourier modes and relate it to the enstrophy.

Is it true that if the viscosity~$\nu$ is zero then the energy is conserved?\end{exercise}

\section{Fourier's cellar}\label{ILVCEMDlamcu}

An interesting classical topic in science deals with the determination of the temperature in the interior of the Earth. Besides the intriguing nature of such a question, the problem has also some practical applications: for instance, for obvious reasons, in mining, but also in environmental sciences, since in regions subject to severe winter or snowfalls, many animals (e.g., earthworms) must find a way to survive the season and one possibility is to dig tunnels underground to find warmer spots
(and of course similar arguments hold for regions with great temperature excursion between day and night,
in which case the animals need to find a suitable refuge on a circadian basis).

Wine connoisseurs may also have an interest in this problem. Indeed, quick temperature fluctuations  can negatively impact wines, while medium and long term storage, and more generally ageing precious wines over a long period of time, strictly demand a thermically stable environment to retain, and even improve, the quality and flavours of wine. A traditional solution for this is to construct a wine cellar keeping wine at a stable, ideal temperature (as side benefits, it would also provide appropriate humidity for longevity and ensure that the wine is stored out of the way of bright light).

\begin{figure}[h]
\includegraphics[height=6.94cm]{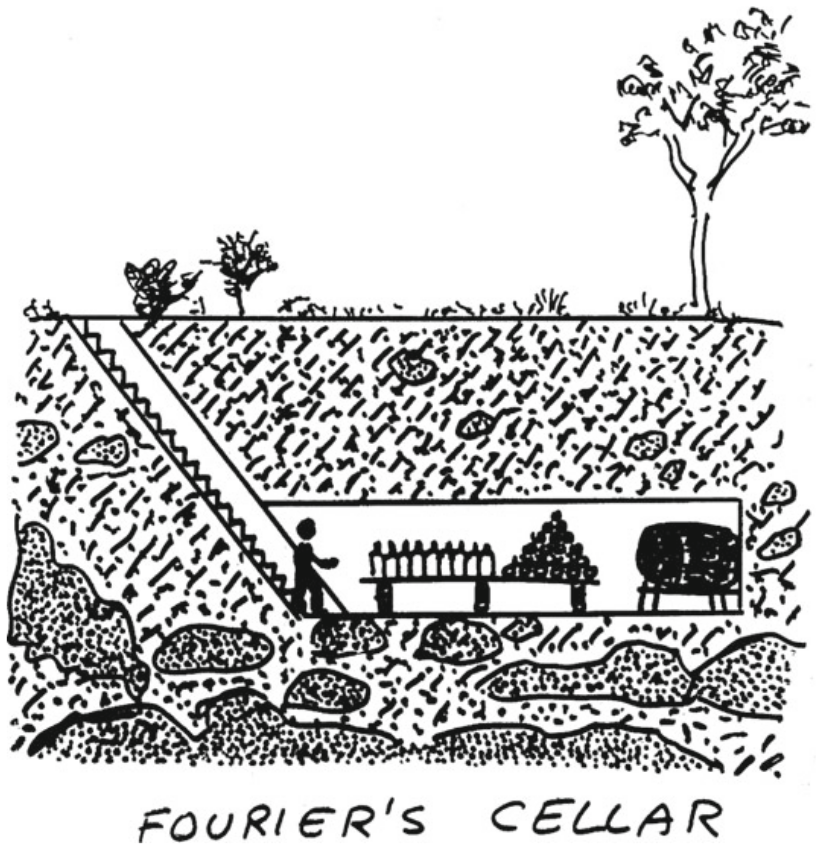}
\centering
\caption{Fourier descending into a cellar: drawing by Enrico Bombieri
(as reproduced in~\cite{MR3616140}).}\label{DECELT24OikdnsEU0}
\end{figure}

The question is therefore at what depth a cellar should be built for it to be at a cool temperature, with only minor seasonal variations, see Figure~\ref{DECELT24OikdnsEU0}. To address this problem, we denote by~$y$ the depth coordinate under the soil (i.e., $y=0$ corresponds to ground level, $y=1$ to 1 meter underground, and so on). Also, since realistically a cellar cannot be placed too far underground, we disregard any potential source of geothermal energy (such as magma, geysers, hot springs, steam vents, etc.).

The temperature at depth~$y$ and time~$t$ will be denoted by~$u(y, t)$ and we assume that the surface temperature is prescribed by environmental conditions as~$u(0, t) = u_0 (t)$.

We assume that heat diffuses according to the heat equation (see Section~\ref{PDESE})
\begin{equation}\label{HEDICOMbrtEAUNVDMDBS} \partial_t u=c \partial^2_y u,\end{equation}
for some~$c>0$ accounting for the thermal conductivity of the soil. For ordinary soil, it is customary to take
\begin{equation}\label{HEDICOMbrtEAUNVDMDBS.c}c=
2 \cdot 10^{-3},\end{equation}
measured in~cm$^2/$sec, which makes~\eqref{HEDICOMbrtEAUNVDMDBS} dimensionally consistent.

Ideally, we could suppose that the temperature~$u$ is periodic in time, with period~$P$ of one year. This allows us (under the Ansatz that the temperature is a nice function) to expand~$u$ in its Fourier Series in the variable~$t$, for a given~$y\in[0,+\infty)$. Thus, in view of~\eqref{fc:PkT-0.45} we have that
\begin{equation*} u(y,t)=\sum _{k\in\Z}\widehat u_{k}(y)\, e^{\frac{2\pi ikt}{P}} \end{equation*}
and similarly, for the surface temperature,
\begin{equation*} u_0(t)=\sum _{k\in\Z}\widehat u_{0,k}\, e^{\frac{2\pi ikt}{P}}. \end{equation*}

Hence, equation~\eqref{HEDICOMbrtEAUNVDMDBS} prescribes that
\begin{eqnarray*}
\sum _{k\in\Z} {\frac{2\pi ik\,\widehat u_{k}(y)}{P}}\, e^{\frac{2\pi ikt}{P}}=
\partial_t u(y,t)=c \partial^2_y u(y,t)=c\sum _{k\in\Z}\widehat u_{k}''(y)\, e^{\frac{2\pi ikt}{P}}.
\end{eqnarray*}
A solution for this equation (and, in fact, essentially the only possible solution, in view of the uniqueness results presented in Section~\ref{SEC:UNIQ:3}) accounting for the surface condition leads to the system
$$ \begin{dcases}
\widehat u_{k}''(y)={\frac{2\pi ik\,\widehat u_{k}(y)}{c P}},\\
\widehat u_{k}(0)=\widehat u_{0,k}.
\end{dcases}$$
This and the observation that
$$ (1+i)^2=2i\qquad{\mbox{and}}\qquad (1-i)^2=-2i$$
yield that
$$ \widehat u_k(y)= \alpha_{1,k}\exp\left( \omega_{1,k}\,y\right)+\alpha_{2,k}
\exp\left(\omega_{2,k}\,y\right),$$
with~$\alpha_{1,k}$, $\alpha_{2,k}\in\C$ such that~$\alpha_{1,k}+\alpha_{2,k}=\widehat u_{0,k}$, and also
$$ \omega_{1,k}:=\begin{dcases} (1+i) \,\sqrt{{\frac{\pi k}{c P}}}&{\mbox{ if }}k\ge0,\\
(1-i) \,\sqrt{{\frac{\pi |k|}{c P}}}&{\mbox{ if }}k<0
\end{dcases}\qquad{\mbox{and}}\qquad\omega_{2,k}:=\begin{dcases} -(1+i) \,\sqrt{{\frac{\pi k}{c P}}}&{\mbox{ if }}k\ge0,\\
(i-1) \,\sqrt{{\frac{\pi |k|}{c P}}}&{\mbox{ if }}k<0.
\end{dcases}$$

Now, to be physically coherent, we take~$\alpha_{1,k}:=0$, otherwise~$\Re\widehat u_k(y)\to+\infty$ as~$y\to+\infty$, and we obtain that
$$ \widehat u_k(y)= \widehat u_{0,k}\,
\exp\left(\omega_{2,k}\,y\right)=\widehat u_{0,k}\,
\exp\left(-\sqrt{{\frac{\pi| k|}{c P}}}\,y\right)\,
\exp\left(-i{\operatorname{sign}}(k)\,\sqrt{{\frac{\pi| k|}{c P}}}\,y\right),
$$
where the sign function notation in~\eqref{lmq-21iDgbnmpeTArhmewxibis} has been adopted.

As a result,
\begin{equation*} u(y,t)=\sum _{k\in\Z}
\widehat u_{0,k}\,
\exp\left(-\sqrt{{\frac{\pi| k|}{c P}}}\,y\right)\,
\exp\left(-i{\operatorname{sign}}(k)\,\sqrt{{\frac{\pi| k|}{c P}}}\,y+\frac{2\pi ikt}{P}\right) .\end{equation*}

Given the exponential decay in the depth~$y$, we make the Ansatz that a good approximation of our solution~$u$ is provided by the lower harmonics, i.e. we restrict to the indices such that~$|k|\le1$.
With this assumption, and recalling the reality condition in~\eqref{fasv}, we gather that
\begin{equation} \label{lmq-21iDgbnmpeTArhmewxiBGBDM-2}
\begin{split}u(y,t)&\simeq 
\widehat u_{0,0}+
\widehat u_{0,1}\,\exp\left(-\sqrt{{\frac{\pi}{c P}}}\,y\right)\,
\exp\left(-i\,\sqrt{{\frac{\pi}{c P}}}\,y+\frac{2\pi it}{P}\right)\\&\qquad\qquad
+\overline{\widehat u_{0,1}}\,
\exp\left(-\sqrt{{\frac{\pi}{c P}}}\,y\right)\,
\exp\left(i\,\sqrt{{\frac{\pi}{c P}}}\,y-\frac{2\pi it}{P}\right)\\&=\widehat u_{0,0}+2\exp\left(-\sqrt{{\frac{\pi}{c P}}}\,y\right)\,\Re\left[
\widehat u_{0,1}\,
\exp\left(-i\,\sqrt{{\frac{\pi}{c P}}}\,y+\frac{2\pi it}{P}\right)\right]\\&=\widehat u_{0,0}+
A_0\,\exp\left(-\sqrt{{\frac{\pi}{c P}}}\,y\right)\,\cos\left(\frac{2\pi t}{P}-\sqrt{{\frac{\pi}{c P}}}\,y+\phi_0\right),
\end{split}\end{equation}
for suitable~$A_0$, $\phi_0$ depending only on the surface condition.

We stress that the term~$\sqrt{{\frac{\pi}{c P}}}\,y$ in the last line of~\eqref{lmq-21iDgbnmpeTArhmewxiBGBDM-2} plays the role of a phase shift in the cosine function.

Thus, to construct a cellar with a stable temperature, we may decide to pick its depth~$y$ in such a way that the heat and cold waves get there with a delay of six months
(say, the cooling arriving during the summer and the warming during the winter). For this, we want to find~$y_\star>0$ such that the wave at depth~$y=y_\star$ corresponds to a shift of size~$\pi$ of the cosine function, namely
$$\sqrt{{\frac{\pi}{c P}}}\,y_\star=\pi,$$
and thus
$$ y_\star=\sqrt{\pi c P}.$$
Using the value in~\eqref{HEDICOMbrtEAUNVDMDBS.c}
(and writing the period of one year as~$P=3600\cdot24\cdot365$ sec),
we obtain
$$ y_\star=\sqrt{\pi \cdot 2\cdot 3.6\cdot24\cdot365}\simeq 445\quad\text{cm}.$$
That is, a cellar located at about four meters and a half underground would probably provide a sufficiently stable temperature for properly ageing some quality wine.

The damping factor in the temperature, corresponding in~\eqref{lmq-21iDgbnmpeTArhmewxiBGBDM-2} to the term
$$ \exp\left(-\sqrt{{\frac{\pi}{c P}}}\,y_\star\right)=e^{-\pi}$$
would also contribute to maintain the thermal environment conveniently stable. Indeed, from the last line in~\eqref{lmq-21iDgbnmpeTArhmewxiBGBDM-2}, we have that the oscillation of temperature at depth~$y$ is given by
$$ A_0\exp\left(-\sqrt{{\frac{\pi}{c P}}}\,y\right).$$
Hence, $A_0$ represents the temperature oscillation at the surface, while~$A_0e^{-\pi}$ corresponds to the temperature variation in the cellar over the year. To have an idea of the effectiveness of this cellar, in Perth, Australia, the winter
minimum temperatures usually sit around 3 degrees Celsius and the maximum in the summer maybe 40 degrees:
this would correspond to an oscillation on the ground of about~$40-3=37$ degrees, but in the cellar of only~$37\,e^{-\pi}$ degrees, that is slightly more than one degree and a half.

In more extreme whether conditions, such as Verkhoyannsk, in Siberia, allegedly the yearly temperature can pass from~$-70$
to~$34$ Celsius degrees, with a remarkable excursion of $34+70=104$: in this case, the temperature oscillation in the cellar is predicted to be about~$104\,e^{-\pi}$ degrees, namely less than four and a half degrees, which is perhaps not ideal for the finest wines (the acceptable level of temperature fluctuation being considered to be about 2 to 3 Celsius degrees in a year)
but likely still within the tolerance range of more commercial wines.

See~\cite[Section~2.6]{MR3616140} and the references therein for more information about this fascinating topic.

\section{The Dirichlet problem on the two-dimensional disk}\label{POISSONKERN-ex4}

The classical Dirichlet problem considers an open domain~$\Omega\subset\R^n$ and some given continuous function~$f$, assigned on the boundary of~$\Omega$, and requires to find a continuous function~$u$ which is twice continuously differentiable in~$\Omega$ and continuous up to the boundary of~$\Omega$ and such that~$u$ is harmonic in~$\Omega$ and~$ u = f$ on the boundary of~$\Omega$.

The Dirichlet problem is one of the major milestones in mathematics
(and physics, and engineering, and biology, and virtually any human discipline) since it studies a question arising,
under different forms, in complex analysis,
harmonic analysis, partial differential equations,
calculus of variations, gravitation, electrostatics, elasticity,
population dynamics, social networks, and so on. See e.g.~\cite{ELEM2} and the references therein for
more information on the Dirichlet problem.

Here we will confine ourselves to the simplest possible case (besides that of an interval in~$\R$), namely the one in which~$\Omega$ is a disk in the plane.
For this, we will employ the \index{periodic Poisson Kernel} periodic Poisson Kernel~${\mathcal{P}}$ introduced in
Exercise~\ref{POISSONKERN-ex2} and use polar coordinates~$(r,\theta)\in[0,1)\times\R$ to describe the unit disk in the plane. In this way, a continuous function~$f:\partial B_1\to\R$ can be thought as~$f=f(\theta)$ and we consider the
function, defined in~$B_1$, 
\begin{equation*} [0,1)\times\R\ni(r,\theta)\longmapsto u_f(r,\theta):=\int_0^1 f(\theta-\tau)\,{\mathcal{P}}(r,\tau)\,d\tau.\end{equation*}

We have that~$u_f$ solves the Dirichlet problem, as remarked by the following result:

\begin{theorem}\label{POIPBDIRHDMOPSMAS}
The function~$u_f$ belongs to~$C^2(B_1)\cap C(\overline{B_1})$ and solves the equation
$$ \begin{dcases}
\Delta u_f=0&{\mbox{ in }}B_1,\\u_f=f&{\mbox{ on }}\partial B_1.
\end{dcases}$$
\end{theorem}

\begin{proof} By virtue of~\eqref{POISSONKERN-ex2-form3}, $u_f$ coincides with the Abel mean of~$f$, as introduced in
Exercise~\ref{POISSONKERN-ex1}. 

Hence, on account of Exercise~\ref{POISSONKERN-ex3}, for all~$\theta\in\R$,
$$ \lim_{r\nearrow1}u_f(r,\theta)=f(\theta),$$
giving that~$u_f$ converges to~$f$ along~$\partial B_1$.

Moreover, by Exercise~\ref{POISSONKERN-ex2}, using the notation~$z:=re^{2\pi ik\theta}$,
\begin{equation}\label{MjAPLSvbiyKamdfRbr} u_f(r,\theta)=\sum_{{k\in\Z}} r^{|k|}\,\widehat f_k\,e^{2\pi ik\theta}=\Re\left(\widehat f_0+2\sum_{k=1}^{+\infty} r^k\,\widehat f_k\,e^{2\pi ik\theta}\right)=\Re\left(\widehat f_0+2\sum_{k=1}^{+\infty} \widehat f_k\,z^k\right).\end{equation}
Being the real part of an analytic function on the complex unit disk, we have (see e.g.~\cite[Theorem~11.4]{MR924157})
that~$u_f$ is harmonic, i.e.~$\Delta u_f=0$
in~$B_1$.
\end{proof}

See e.g.~\cite[Chapter~11]{MR924157} for further reading about harmonic functions in the plane.

\begin{exercise}\label{SPKEMLEMFCf01}
Prove that the function~$u_f$ in Theorem~\ref{POIPBDIRHDMOPSMAS} satisfies the following maximum principle. If~$f\ge0$, then~$u_f\ge0$. \index{maximum principle}
Similarly, if~$f\le0$, then~$u_f\le0$.
\end{exercise}

\begin{exercise}\label{U3fNIDIRIY.1}
Prove that there exists one and only one solution of the Dirichlet Problem in Theorem~\ref{POIPBDIRHDMOPSMAS}.
\end{exercise}

\begin{exercise}\label{U3fNIDIRIY.1bisa}
Prove that the function~$u_f$ in Theorem~\ref{POIPBDIRHDMOPSMAS} satisfies the following mean value property: for every~$p_0\in B_1$ and every~$r\in(0,1-|p_0|)$, we have that \index{mean value property}
\begin{equation}\label{U3fNIDIRIY.1bisa-EQ1} u_f(p_0)=\frac1{\pi r^2}\int_{B_r(p_0)} u_f\end{equation}
and
\begin{equation}\label{U3fNIDIRIY.1bisa-EQ2} u_f(p_0)=\frac1{2\pi r}\int_{\partial B_r(p_0)} u_f.\end{equation}
\end{exercise}

\section{The Dirichlet problem on a square and the behaviour near corners}\label{TEojemy04i7u:Asdweg54uk:028tuyjh}

\begin{figure}[h]
\includegraphics[height=6.5cm]{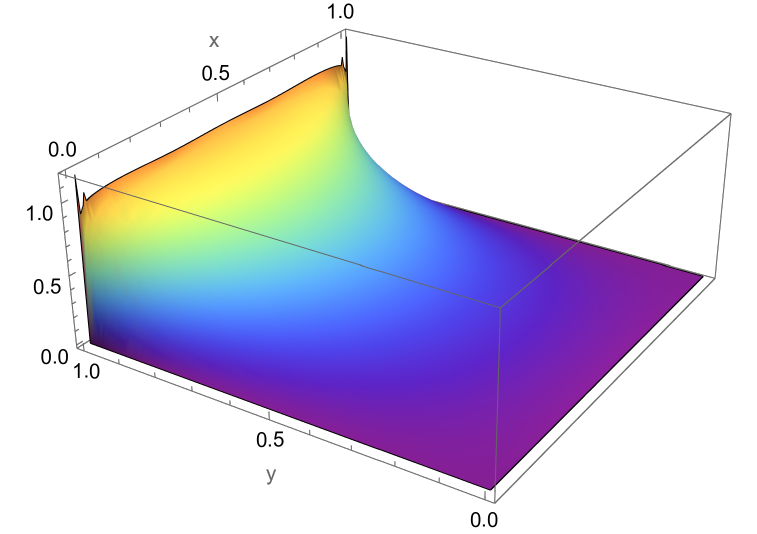}
$\,$\includegraphics[height=6.5cm]{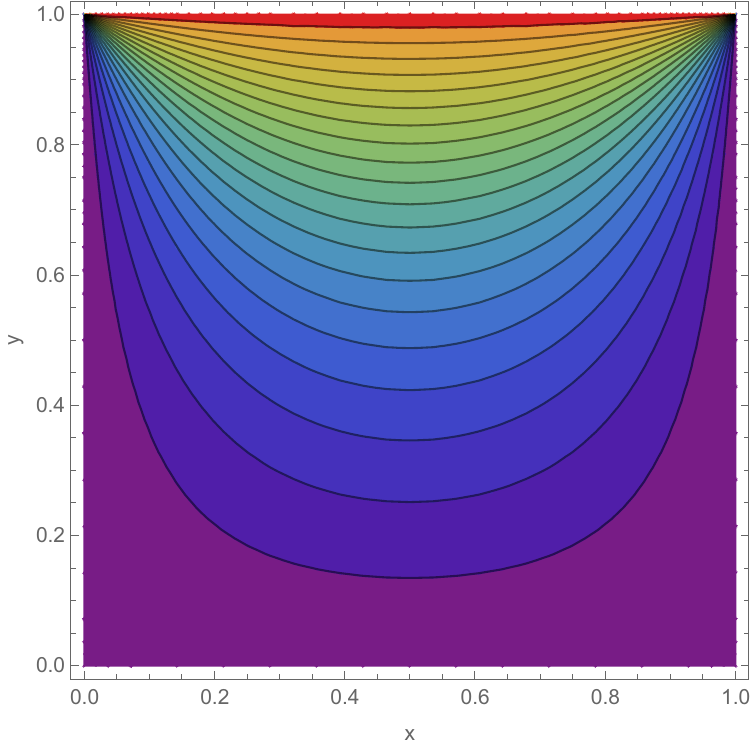}
\centering
\caption{Plot of the function~$\displaystyle
\frac{4}{\pi} \sum_{k=0}^{100} \frac{1}{2k+1}\frac{\sinh((2k+1)\pi y)}{\sinh((2k+1)\pi)}\sin((2k+1)\pi x)$ and sketch of its
level sets.}\label{qTAF.fc3rcbD2.mEcornT24Oikt.9dnsEU0}
\end{figure}

Sometimes, the naive use of mathematics yields surprisingly useful and physically sound results. We consider the problem of finding the steady-state of the heat equation in a unit square in two-dimensions, say~$(0,1)\times(0,1)$, where three edges are maintained at a reference temperature of $0$ and the fourth edge at a temperature of $1$. This is expressed via the boundary value problem
\begin{equation}\label{u:funcorners:2}
\begin{dcases}
  \Delta u(x,y) = 0& {\mbox{for all }}(x,y)\in (0,1)\times(0,1),\\
  u(0,y)=u(1,y)=0& {\mbox{for all }} y\in(0,1),\\
  u(x,0)=0 & {\mbox{for all }} x\in(0,1),\\
  u(x,1)=1 & {\mbox{for all }} x\in(0,1).
\end{dcases}
\end{equation}
Fourier's method of separation of variables works beautifully to produce the solution
\begin{equation}\label{u:funcorners}
u(x,y) = \frac{4}{\pi} \sum_{k=0}^{+\infty} \frac{1}{2k+1}\frac{\sinh((2k+1)\pi y)}{\sinh((2k+1)\pi)}\sin((2k+1)\pi x),
\end{equation}
see Exercise~\ref{jlwmeotrjho590iukjiqdhfpirnHSNdlm0ujt1}
(see also Figure~\ref{qTAF.fc3rcbD2.mEcornT24Oikt.9dnsEU0}
for a sketch of an approximation of this function).

As expected, the convergence properties of the above series are 
somewhat problematic at the corners~$(0,1)$ and~$(1,1)$,
where the boundary datum is discontinuous. Any truncation will poorly represent the behaviour of the temperature function, and will suffer from a two-dimensional variant of Gibbs phenomenon in the corner regions. We need to verify that there are no issues with the model description, and that the mathematics does not breakdown. It may be physically impossible to arrange to
have~$u(x,0)=0$ on one edge near the corner and~$u(x,1)=1$ on the adjacent edge. Alternatively, it may be that the series solution is ill-suited to describe the sharp transition near the corner.

Fortunately, mathematics is strong enough to clarify the situation.
To this end, to appreciate the structure
of~$u$ near the corners at which the boundary datum is discontinuous
(say, focusing on the corner~$(0,1)$) we use polar coordinates
\begin{equation}\label{POCAUSDM:D} x= r\sin\theta\qquad{\mbox{and}}\qquad y= 1 - r\cos\theta\end{equation}
with~$r\in(0,1)$ (to be taken small
when we want to extract information at the corner) and~$\theta\in\left(0,\frac\pi2\right)$.

Also, given~$\rho\in\left[0,\frac12\right)$, we set
\begin{equation}\label{lkajd9358uyjh:pkdfo8mnb7m0KSdjfgwrnb7m80}
Q:=[0,1]\times[0,1], \quad
Q_{1,\rho}:=[0,\rho]\times[1-\rho,1]\quad{\mbox{and}}\quad
Q_{2,\rho}:=[1-\rho,1]\times[1-\rho,1].\end{equation}
Then, we have:

\begin{theorem}
The series in~\eqref{u:funcorners} converges locally uniformly in~$(0,1)\times(0,1)$
and, for every~$\rho\in\left[0,\frac12\right)$, uniformly in~$Q\setminus(Q_{1,\rho}\cup Q_{2,\rho})$.

Moreover, for small~$r>0$,
\begin{equation}\label{LOGREPDO343562.23}
u(x,y)=\frac{2\theta}{\pi}+O(r^2)+v(x,y),\end{equation}
where~$v\in C^m([0,1]\times[0,1])$ for every~$m\in\N$.
\end{theorem}

\begin{proof}
To simplify the analysis, we first observe that the series in~\eqref{u:funcorners} can be rewritten in the form
\begin{equation}\label{u:funcorners:W}
u(x,y) = \frac{4}{\pi} \sum_{k=0}^{+\infty} \frac{e^{(2k+1)\pi (y-1)}\,\sin((2k+1)\pi x)}{2k+1}+v(x,y),
\end{equation}
where~$v\in C^m([0,1]\times[0,1])$ for every~$m\in\N$, see Exercise~\ref{jlwmeotrjho590iukjiqdhfpirnHSNdlm0ujt1BIS}.

In complex notation, we can also rephrase~\eqref{u:funcorners:W} in the form
\begin{equation*}
u(x,y) = \frac{4}{\pi} \sum_{k=0}^{+\infty} \Im\left(\frac{e^{(2k+1)\pi \zeta}}{2k+1}\right)+v(x,y),\qquad{\mbox{where }}\;
\zeta:=(y-1)+ ix\in\C
\end{equation*}
or simply
\begin{equation}\label{u:funcorners:Z}
u(x,y) = \frac{4}{\pi} \sum_{k=0}^{+\infty} \Im\left(\frac{Z^{2k+1}}{2k+1}\right)+v(x,y),\qquad{\mbox{where }}\;
Z:=e^{\pi \zeta}\in\C.
\end{equation}

We stress that, if~$(x,y)\in(0,1)\times(0,1)$,
\begin{equation} \label{9hf934tkbm79VYASjow83thnb1k04-4yuj9}|Z|=e^{\pi(y-1)}<1,\end{equation}
ensuring the locally uniform convergence of the above series in~$(0,1)\times(0,1)$,
as well as the uniform convergence in the set~$Q\setminus(Q_{1,\rho}\cup Q_{2,\rho})$, see Exercise~\ref{jlwmeotrjho590iukjiqdhfpirnHSNdlm0ujt2}.

Another consequence of~\eqref{u:funcorners:Z}
and~\eqref{9hf934tkbm79VYASjow83thnb1k04-4yuj9}
consists in the complex logarithm representation\begin{equation}\label{LOGREPDO343562}\begin{split}
u(x,y)&=\frac{4}{\pi} \Im\left(\sum_{k=0}^{+\infty} \frac{Z^{2k+1}}{2k+1}\right)+v(x,y)\\&=
\frac{2}{\pi} \Im\left(\sum_{j=1}^{+\infty} \frac{Z^{j}}{j}-
\sum_{j=1}^{+\infty} \frac{(-1)^jZ^{j}}{j}
\right)+v(x,y)\\&=\frac{2}{\pi} \Im\left(
\operatorname{Log}(1+Z)-\operatorname{Log}(1-Z)
\right)+v(x,y).\end{split}\end{equation}

Now, using the polar coordinate notation
in~\eqref{POCAUSDM:D}, for small~$r>0$ it follows that
\begin{equation}\label{SDMLSPMCAUSDM:D}
Z=e^{-\pi r(\cos\theta- i\sin\theta)}
.\end{equation}

We want to compute the argument~$\alpha_\pm$
of the complex number~$1\pm Z$. This will be useful
in comparison with~\eqref{LOGREPDO343562},
because if~$1\pm Z=\rho_\pm\,e^{i\alpha_\pm}$ with~$\rho_\pm>0$, we have that~$\operatorname{Log}(1\pm Z)=\ln\rho_\pm +i\alpha_\pm$ and therefore
\begin{equation}\label{REWSDV:qwdfgvbweEQscFsdRWEGRBN-1}
\Im\left(
\operatorname{Log}(1+Z)-\operatorname{Log}(1-Z)
\right)=\Im\left(
\ln\rho_+ +i\alpha_+
-\ln\rho_- -i\alpha_-
\right)=\alpha_+-\alpha_-.
\end{equation}

Now we point out that (see Exercise~\ref{jlwmeotrjho590iukjiqdhfpirnHSNdlm0ujt2hgd})
\begin{equation}\label{oadhfweiroeghmtinvikm9iryhtgb02oegk90toh4657yh9203tmtr7u5}
\alpha_+= \frac{\pi r \sin\theta}2+O(r^2)\qquad{\mbox{and}}\qquad
\alpha_-=\theta +\frac{\pi r \sin\theta}2  + O(r^2). 
\end{equation}
This and \eqref{REWSDV:qwdfgvbweEQscFsdRWEGRBN-1}
yield that
$$\Im\left(
\operatorname{Log}(1+Z)-\operatorname{Log}(1-Z)
\right)=\theta+O(r^2).$$
{F}rom this and~\eqref{LOGREPDO343562} we deduce~\eqref{LOGREPDO343562.23}, as desired.\end{proof}

\begin{figure}[h]
\includegraphics[height=6.5cm]{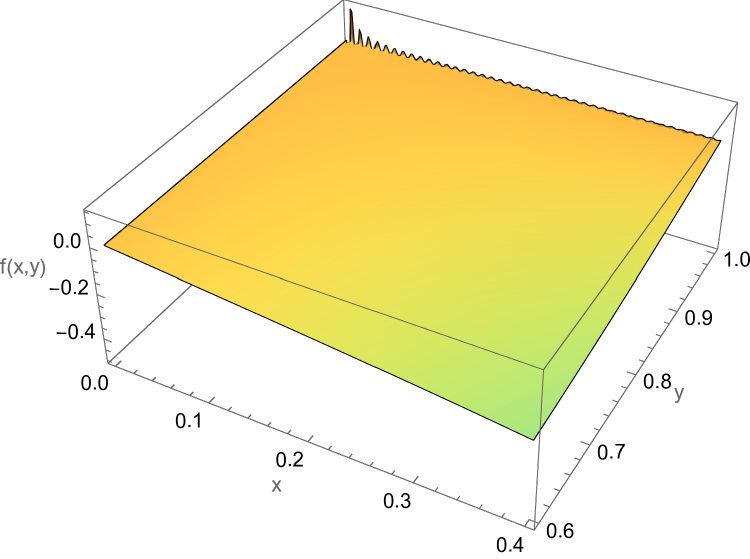}
$\,$\includegraphics[height=6.5cm]{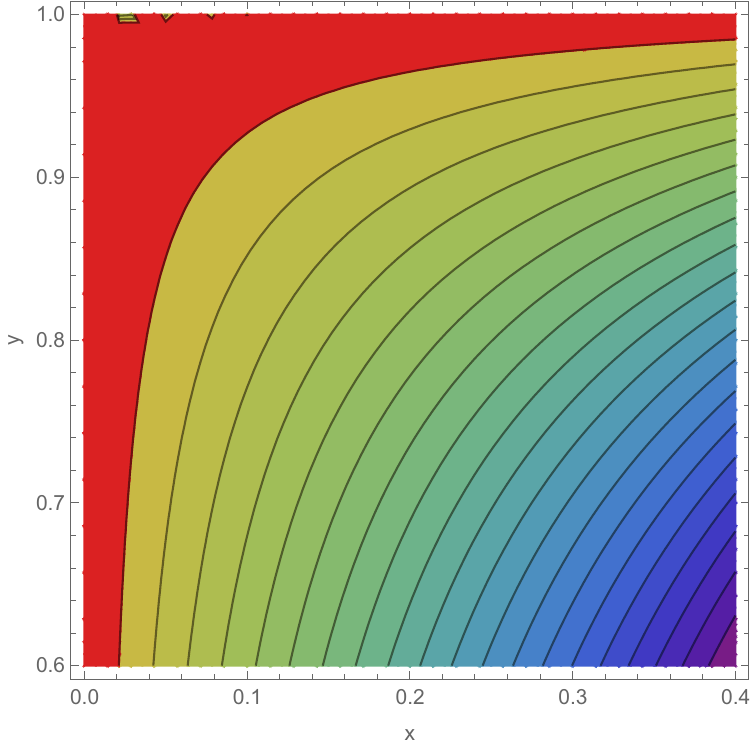}
\centering
\caption{Plot of the function~$\displaystyle
\frac{4}{\pi} \sum_{k=0}^{100} \frac{1}{2k+1}\frac{\sinh((2k+1)\pi y)}{\sinh((2k+1)\pi)}\sin((2k+1)\pi x)-\frac{2}\pi \arctan\left(\frac{x}{1-y}\right)$ in the vicinity of~$(0,1)$
and sketch of its
level sets.}\label{qTAF.fc3rcbD2.mEcornT24Oikt.9dnsEU0.0oeirjfg2}
\end{figure}

We remark that~\eqref{LOGREPDO343562.23}
is interesting, since it identifies
the behaviour of the solution~$u$ near the corner
as a (rescaled) angular variable, up to a small correction.
See Figure~\ref{qTAF.fc3rcbD2.mEcornT24Oikt.9dnsEU0.0oeirjfg2} for a sketch of the quadratic
separation between the solution~$u$
and the angular fuction.
See~\cite{MR1297119} for further
information about this model
and related problems.

\begin{exercise}\label{jlwmeotrjho590iukjiqdhfpirnHSNdlm0ujt1}
Check that the function in~\eqref{u:funcorners}
is a solution of~\eqref{u:funcorners:2}. 
\end{exercise}

\begin{exercise}\label{jlwmeotrjho590iukjiqdhfpirnHSNdlm0ujt1BIS}
Check that~\eqref{u:funcorners}
can be rewritten as~\eqref{u:funcorners:W}.
\end{exercise}

\begin{exercise}\label{jlwmeotrjho590iukjiqdhfpirnHSNdlm0ujt1BIS.04}
Let~$\delta\in\left(0,\frac12\right)$. Suppose that~$Z=r e^{\pi i\phi}$, with~$r\in(\delta,1]$ and~$\phi\in[\delta,1-\delta]$.

Prove that there exists~$C>0$, depending only on~$\delta$, such that, for all~$t\in(0,+\infty)$,
$$ \left|\frac{1-Z^t}{1-Z}\right|\le C.$$
\end{exercise}

\begin{exercise}\label{jlwmeotrjho590iukjiqdhfpirnHSNdlm0ujt2}
Assume the notation in~\eqref{lkajd9358uyjh:pkdfo8mnb7m0KSdjfgwrnb7m80}.
Check that the series in~\eqref{u:funcorners} and~\eqref{u:funcorners:Z} converge uniformly in the set~$Q\setminus(Q_{1,\rho}\cup Q_{2,\rho})$.
\end{exercise}

\begin{exercise}\label{jlwmeotrjho590iukjiqdhfpirnHSNdlm0ujt2hgd}
Prove~\eqref{oadhfweiroeghmtinvikm9iryhtgb02oegk90toh4657yh9203tmtr7u5}.
\end{exercise}

\section{Calculating the exact value of a series}\label{Pkjqn90-902}

A beautiful, important, and rather surprising application of Fourier methods consists in finding explicitly the exact value of an infinite series, a task which is very rarely possible to achieve. The paradigmatic example is the following one:

\begin{proposition}\label{BASEL}
We have that
\begin{equation}\label{BAFO} \sum_{k=1}^{+\infty }\frac{1}{k^2}=\frac{\pi^2}6.\end{equation}
\end{proposition}

This result is not a mere curiosity, it actually provides interesting information about a number of several problems. As a matter of fact,
the question of finding the value of the above series was posed already by Pietro Mengoli in 1650.
Though many prominent mathematicians attacked the problem from different angles
(for instance, in 1665 John Wallis computed the series to three decimal places finding the approximated result of~$1.645$), the question remained as a major unsolved problem until the rise of a young (28 years old at that time) genius, named\footnote{The question of precisely evaluating the infinite series on the left-hand side of~\eqref{BAFO} was called the \emph{Basel problem}.
\index{Basel problem} Euler was born in Basel. However, the name of the question is probably due to the fact that Basel was also the hometown of the Bernoulli family, who extensively worked on it, and especially the publishing location of Jacob Bernoulli's ``Ars conjectandi, opus posthumum. Accedit Tractatus de seriebus infinitis'' (published in 1689 and available at \url{http://www.kubkou.se/pdf/mh/jacobB.pdf}).
In particular, Jacob (sometimes spelled as Jakob, and also known as James or Jacques) Bernoulli included a reference to the problem and
proved that the series in~\eqref{BAFO} is less than~$2$ by an elegant telescopic trick (see Exercise~\ref{TELES}). Euler probably heard of the question from his mentor, Johann (also known as Jean or John) Bernoulli, brother of Jacob
(curiously, Jacob had previously taught Euler's father at the University of Basel).

After an initial period of collaboration, a bitter rivalry arose between Jacob and Johann.
Jacob died about thirty years before Euler found the solution to the Basel problem.
When Johann heard the news of Euler's solution to the problem, he vowed 
``My brother's most ardent wish is satisfied! If only my brother were still alive!'' (see~\cite[page~22]{MR256841}).

See also~\cite{MR2916494} and the references therein for more information about the history of the Basel problem.} Leonhard Euler.

To deal with the question, Euler, endowed with a remarkable computational capacity, first sharpened the
approximated result to six decimal places 
and later to the still approximated, but
very accurate, value\footnote{For the approximation to six decimal places, Euler
transformed the given series into two new series
which converge extremely rapidly. For the subsequent approximation, he introduced a summation formula which is called nowadays
Euler-Maclaurin formula, see~\cite[Section~9.5]{MR1397498}.}
of~$1.64493406684822643647$, till, in 1734, he discovered\footnote{In Euler's own words:
``quite unexpectedly I have found an
elegant formula involving the quadrature of the circle'',
the quadrature of the circle, on this occasion, referring to the number~$\pi$ (the area, or ``quadrature'',
of the circle).} the elegant result in Proposition~\ref{BASEL}.
The solution found by Euler in~1735 was incredibly pioneering at that time and required about a century of solid mathematical work to be completely justified (see~\cite{MR2321397} and~\cite[pages~37-40]{MR2561962} for a detailed account of this story and a discussion about Euler's original line of reasoning).

The result in Proposition~\ref{BASEL} is actually a cornerstone linking special functions, harmonic analysis,
number theory, and probability, as we will see in the exercises below
and in Section~\ref{NUTESE}. A unified framework to this type of results is given by the \index{Riemann zeta function}
Riemann zeta function, which is defined, for all~$s\in\C$ with~$\Re s>1$, as
\begin{equation}\label{Zzetas} \zeta (s):=\sum_{k=1}^{+\infty }{\frac{1}{k^{s}}}.\end{equation}
For other complex values,
the Riemann zeta function can also be defined via analytic continuation.

Not many explicit values of the Riemann zeta function are known, but remarkably Proposition~\ref{BASEL} states that~$\zeta(2)=\frac{\pi^2}6$.

Deeply understanding the Riemann zeta function is probably one the hardest challenges in mathematics: in particular, the so-called
Riemann Hypothesis is the conjecture that the zeros of the Riemann zeta function can only occur
at the negative even integers and at complex numbers with real part equal to~$\frac12$.
Likely, giving an answer to this question will make anyone become the most famous mathematician in the world
(plus receiving, as a small token of appreciation, a million US dollars from the Clay Mathematics Institute).
See e.g.~\cite{MR2060134, MR2063737} and the references therein for a fascinating and accessible introduction\footnote{Here are some lyrics by Tom Apostol (see~\cite[pages 394--395]{MR2063737} for the full song):

{\em Where are the zeros of zeta of s? 

G. F. B. Riemann has made a good guess:

``They're all on the critical line,'' stated he,

``And their density's one over two pi log T.''

This statement of Riemann's has been like a trigger,  

And many good men, with vim and with vigor,

Have attempted to find, with mathematical rigor,

What happens to zeta as mod t gets bigger.}} to
the \index{Riemann Hypothesis} Riemann Hypothesis.

\begin{proof}[Proof of Proposition~\ref{BASEL}] For every~$x\in[0,1)$ we let~$f(x):=x(1-x)$ and consider its periodic extension~$f_{\text{per}}$ of period~$1$, as defined in~\eqref{FPER}.

We know (see Exercise~\ref{ka-2}) that the Fourier Series of~$f_{\text{per}}$ equals
\begin{equation}\label{FAJHMNS02r} \frac16-\sum_{k\in\Z\setminus\{0\}}\frac{1}{2\pi^2 k^2}e^{2\pi ikx}.\end{equation}
We also note that~$f_{\text{per}}$ is Lipschitz and therefore it fulfils Dini's Condition~\eqref{DINI} at every point~$x$
(see Exercises~\ref{DINI1} and~\ref{DINI2}).

Consequently, by the pointwise convergence of Fourier Series (see Theorem~\ref{DINITS}) we deduce that~$f_{\text{per}}(x)$ equals~\eqref{FAJHMNS02r}
for every~$x\in\R$. In particular, taking~$x:=0$,
$$ 0=f_{\text{per}}(0)=\frac16-\sum_{k=1}^{+\infty}\frac{1}{\pi^2 k^2},$$
from which we obtain~\eqref{BAFO}.
\end{proof} 

Many different proofs of Proposition~\ref{BASEL} using different methods have been put forth.
Some of these proofs based on analytic techniques will be recalled in the forthcoming exercises (see also~\cite{MR2826455} and~\cite[Problem~2.3.25]{MR2961901} for a proof based on probability theory).

See also~\cite[Chapter~12]{zbMATH05014111} for additional information about this topic and related ones.

\begin{exercise}\label{TELES}
Prove that
$$ \sum _{k=1}^{+\infty }\frac{1}{k^2}\le2$$
by using a telescopic trick.
\end{exercise}

\begin{exercise}\label{TELES-INTR}
Prove that
$$ \sum _{k=1}^{+\infty }\frac{1}{k^2}\le2$$
by using an integral estimate.
\end{exercise}

\begin{exercise}\label{AOJsnkwe9rgeotg2-34rt-p0}
Give a proof of Proposition~\ref{BASEL} by using
the convergence theory of Fourier Series in~$L^2((0,1))$, as presented in Section~\ref{COL2}.\end{exercise}

\begin{exercise}\label{AOJsnkwe9rgeotg2-34rt}
Give a proof of Proposition~\ref{BASEL} by using Euler's sine product formula (see Exercise~\ref{PROSI}).\end{exercise}

\begin{exercise}\label{AOJsnkwe9rgeotg2-34rt-DOU}
Give a proof of Proposition~\ref{BASEL} by computing in two different ways the double integral
$$ \int_0^1\int_0^1 \frac{dx\,dy}{1-xy}.$$\end{exercise}

\begin{exercise}\label{SOMMS1}
The $n$th \index{harmonic number}
harmonic number, denoted by~$H_n$, is defined as the sum of the reciprocals of the first~$n$ natural numbers, namely
$$ H_n:=\sum_{j=1}^n \frac1j.$$
Prove that
$$ \sum_{n=1}^{+\infty}\frac{H_n}{n(n+1)}=\frac{\pi^2}{6}.$$\end{exercise}

\begin{exercise}\label{SOMMS2} Prove that
$$ \sum_{n=1}^{+\infty} \frac{\displaystyle 1+\frac12+\dots+\frac1n}{1+2+\dots+n}
=\frac{\pi^2}{3}.$$
\end{exercise}

\begin{exercise}\label{EUFO}
Euler's product formula states that, if~$s\in\C$ and~$\Re s>1$,
\begin{equation}\label{87uye80217yrg2353fhv5h6h434125} \zeta(s)=\prod_{{p{\text{ prime}}}\atop{p>1}}{\frac{1}{1-p^{-s}}}.\end{equation}
Prove this formula.
\end{exercise}

\begin{exercise}\label{EUFOanc}
Prove that
\begin{equation}\label{mnhgreg-190r2f439uut9032yf936hgv6b} \frac{\pi^2}6=\prod_{{p{\text{ prime}}}\atop{p>1}}{\frac{p^2}{p^2-1}}.\end{equation}
\end{exercise}

\begin{exercise}\label{EUCL}
Give five different proofs of the fact that there are infinitely many primes.\end{exercise}

\begin{exercise}\label{HAPRO-0} Prove that
$$\prod_{{p{\text{ prime}}}\atop{p>1}}\frac1{1-p^{-1}}=+\infty.$$
\end{exercise}

\begin{exercise}\label{HAPRO} Euler established that the sum of the reciprocals of all prime numbers diverges, i.e.
\begin{equation}\label{02jf90ewnv4578-213660} \sum_{{p{\text{ prime}}}}\frac1p=+\infty.\end{equation}
Prove this statement.
\end{exercise}

\begin{exercise}\label{x1kdo021ps6lungajkd}
Prove that
\begin{equation}\label{VACPAI} \sum_{k=1}^{+\infty}\frac1{k^4}=\frac{\pi^4}{90}\end{equation}
by using Exercise~\ref{0ei02r3290234pfp42}.
\end{exercise}

\begin{exercise}\label{KALMS0owk3klxslkfw-e}
Give another proof of~\eqref{VACPAI} by using Exercise~\ref{ka-2}.
\end{exercise}

\begin{exercise} Calculate the Fourier Series of the function
$$ f(x):=\big|\cos(\pi x)\big|,$$
discuss its pointwise convergence, and use this information to prove that
$$\sum_{k=1}^{+\infty} \frac{(-1)^{k-1}}{4k^2-1}=\frac{\pi}4-\frac{1}{2}.$$\end{exercise}

\begin{exercise} \label{aujsoJs893405}
Let
$$ \left[-\frac12,\frac12\right)\ni x\longmapsto f(x):=\cosh(2x)$$
and extend~$f$ to a periodic function of period~$1$.

Using the Fourier Series of~$f$, calculate
\begin{equation}\label{90OJSldfw0eojrf} \sum_{k=1}^{+\infty} \frac{(-1)^k}{k^2(\pi^2 k^2+1)}.\end{equation}
\end{exercise}

\section{Number theory}\label{NUTESE}

We have already seen a glimpse of number theory on several occasions
(such as Exercise~\ref{PROSIW} and in many instances in Section~\ref{Pkjqn90-902}). Here we recall this very alluring result, which is also related to probability:

\begin{theorem}
The probability that two random numbers are co-prime is~$\frac{6}{\pi^2}$.
\end{theorem}

\begin{proof} We have to calculate the probability that two natural numbers~$a$, $b$,
randomly and independently chosen, 
are co-prime.

We observe that the integer~$a$ has probability~$\frac1p$ to be divisible by a given prime~$p$
(because~$a=kp+j$, with~$k\in\N$ and~$j\in\{0,\dots,p-1\}$, that is one integer out of~$p$ is divisible by~$p$).

The same is true for the integer~$b$, hence, if~$a$ and~$b$ are independently chosen, the probability that they are both simultaneously divisible by~$p$ is~$\frac1{p^2}$.

In other words, the probability that~$a$ and~$b$ are not simultaneously divisible by a given prime~$p$ equals~$1-\frac1{p^2}=1-p^{-2}$.

Accordingly, the probability that~$a$ and~$b$ are co-prime (i.e., not
simultaneously divisible by any prime~$p$) is
$$ \prod_{{{p{\text{ prime}}}}}(1-p^{-2})=:{\mathcal{P}}.$$

Therefore, by Euler's product formula (see Exercise~\ref{EUFO}, used here with~$s:=2$),
$$ \frac1{\mathcal{P}}=\prod_{{{p{\text{ prime}}}}}\frac{1}{1-p^{-2}}=\zeta(2).$$
This and Proposition~\ref{BASEL} give that~$\frac1{\mathcal{P}}=\frac{\pi^2}6$ and then~${\mathcal{P}}=\frac6{\pi^2}$.
\end{proof}

See also Sections~\ref{PelapeammilonS} and~\ref{SE:MINK} for further applications to number theory.

\section{Inequalities of analytic flavour}\label{CSLE:A:SEZZ0-1}

Fourier methods are often useful to establish convenient functional inequalities.
As an example, we recall the following \index{Poincar\'e-Wirtinger Inequality} Poincar\'e-Wirtinger Inequality:

\begin{theorem}\label{CSLE:A:SEZZ0-1-pwin}
Let~$f\in C^1(\R)$, periodic of period~$1$, and such that
\begin{equation}\label{lk2340-012.cnh.ll-d} \int_0^1 f(x)\,dx=0.\end{equation}

Then,
\begin{equation}\label{lk2340-012.cnh} \int_0^1 |f(x)|^2\,dx\le\frac1{4\pi^2}\int_0^1 |f'(x)|^2\,dx.\end{equation}
\end{theorem}

\begin{proof} Among the many available in the literature, a possible proof proceeds as follows. We develop~$f$ in Fourier Series, owing to Theorem~\ref{C1uni}
(recall also Exercises~\ref{DINI1UNI} and~\ref{DINI2UNI}) and we observe that, by reason of~\eqref{lk2340-012.cnh.ll-d}, we have that~$\widehat f_0=0$.

Furthermore, on account of Theorem~\ref{THCOL2FB}, we know that
Parseval's Identity in~\eqref{PARS} holds true (both for~$f$ and for~$f'$).

Besides, by periodicity,
$$ \widehat{(f')}_0=\int_0^1f'(x)\,dx=f(1)-f(0)=0.$$

As a result, recalling Theorem~\ref{SMXC22b} we see that
\begin{equation}\label{BNCDVCuPXasjmdc.1}\begin{split}& \int_0^1 |f(x)|^2\,dx=\sum_{k\in\Z} |\widehat f_k|^2 
=\sum_{k\in\Z\setminus\{0\}}  |\widehat f_k|^2\le \sum_{k\in\Z\setminus\{0\}}|k \widehat f_k|^2\\&\qquad
=\frac{1}{4\pi^2}
\sum_{k\in\Z\setminus\{0\}}  |\widehat{Df}_k|^2=\frac{1}{4\pi^2}
\sum_{k\in\Z}  |\widehat{(f')}_k|^2=
\frac{1}{4\pi^2}\int_0^1 |f'(x)|^2\,dx,\end{split}\end{equation}
as desired.
\end{proof}

See~\cite[Section~4.4]{MR1412143} and~\cite{MR1814364} for several versions of the Poincar\'e-Wirtinger Inequality and related results.
See also~\cite{MR322688} for discrete counterparts of the Poincar\'e-Wirtinger Inequality.

\begin{exercise}\label{BNCDVCuPXasjmdc.1e} For which functions (if any) does the inequality in~\eqref{lk2340-012.cnh} become an equality?
\end{exercise}

\begin{exercise}\label{BNCDVCuPXasjmdc.1e.102i} A variant (among the many) of the 
Poincar\'e-Wirtinger Inequality in~\eqref{lk2340-012.cnh} goes as follows.
Let~$f\in C^1([0,1])$ with~$f(0)=f(1)=0$.

Then,
\begin{equation}\label{lk2340-012.cnh.102} \int_0^1 |f(x)|^2\,dx\le\frac1{\pi^2}\int_0^1 |f'(x)|^2\,dx.\end{equation}

Prove this statement.\end{exercise}

\begin{exercise}\label{BNCDVCuPXasjmdc.1e.102} For which functions (if any) does the inequality in~\eqref{lk2340-012.cnh.102} become an equality?
\end{exercise}

\begin{exercisesk}\label{CSLE:A}
The following result is known as Ste\v{c}kin's Lemma.

Suppose that~$f:\R\to\R$ is a trigonometric polynomial of degree~$N\ge1$, with~$x_0\in\R$ such that
$$ f(x_0)=\max_{x\in\R}|f(x)|.$$
Then, for every~$y\in\left[-\frac{1}{2N},\frac{1}{2N}\right]$,
\begin{equation}\label{lk2340-012} f(x_0+y)\ge f(x_0)\,\cos(2\pi Ny).\end{equation}

Prove this result.\end{exercisesk}

\begin{exercisesk}\label{CSLE:B}
The following result\footnote{See~\cite[page~12]{MR545506}, \cite[pages~97--109]{MR1261635}
 and the references therein for further information on this topic.} is known as Bernstein's Inequality.

Suppose that~$f:\R\to\C$ is a complex trigonometric polynomial of degree~$N$.

Then,
\begin{equation}\label{BVELE} \max_{x\in\R}|f'(x)|\le 2\pi N\max_{x\in\R}|f(x)|.\end{equation}
Prove this result.
\end{exercisesk}

\begin{exercise}\label{BVELEb} Prove that the estimate in~\eqref{BVELE} is sharp by providing a trigonometric polynomial of degree~$N$
attaining the equality.
\end{exercise}

\section{Inequalities of geometric flavor}\label{CSLE:A:GEO}

Along the lines of Section~\ref{CSLE:A:SEZZ0-1}, Fourier methods can be profitably used to establish inequalities with a geometric significance. We start with a version of the \index{Isoperimetric Inequality} Isoperimetric Inequality in the plane, with a proof relying
on Fourier Series which\footnote{Actually, the regularity requirements on the curve
can be relaxed and are assumed here just for the sake of simplicity, mainly to avoid
the introduction of generalised notions of length, area, and normal, as well as to use the Fourier methods in their simplest
possible formulation, see also~\cite[Section~4.1]{MR1412143} to appreciate how to set this approach into a general framework. See also Exercise~\ref{kmsdcJAx.iton4.143.e}
for an approach to the Isoperimetric Inequality in Theorem~\ref{kmsdcJAx.iton4.143}
relying on the Poincar\'e-Wirtinger Inequality.} was proposed in~\cite{zbMATH02662735}.

\begin{theorem}\label{kmsdcJAx.iton4.143}
Let $I$ be a closed interval.
Let $\gamma:I\to\R^2$ be the parametrisation of a regular, closed curve of class~$C^\infty$
and consider its length
\begin{equation} \label{LEN:TH}{\mathcal{L}}:=\int_I |\gamma'(t)|\,dt\end{equation}
and the area~${\mathcal{A}}$ enclosed by the curve.

Then,
\begin{equation}\label{POKSJND:lasdmcv} {\mathcal{L}}^2\ge 4\pi{\mathcal{A}}.\end{equation}

Also, equality in~\eqref{POKSJND:lasdmcv} holds true if and only if the curve is a circle.
\end{theorem}

\begin{proof} Since~\eqref{POKSJND:lasdmcv} is invariant under dilation, we can suppose that~${\mathcal{L}}=1$.
Moreover, to simplify computations, it is convenient to suppose, without loss of generality, that the curve is parameterised by its arc-length, i.e.~$I=[0,{\mathcal{L}}]=[0,1]$ and~$|\gamma'(t)|=1$ for all~$t\in[0,1]$.

We write the Cartesian coordinates in the plane as~$X=(x,y)\in\R^2$ and the components of the curve under consideration as~$\gamma(t)=(x(t),y(t))$. We also denote by~$\Omega$ the planar region encircled by the curve.
We assume, without loss of generality, that the curve is travelled anticlockwise  and denote by~$\nu$ its outward unit normal.
With this notation, we have that
$$ \nu(\gamma(t))=\frac{\big(y'(t),\,-x'(t)\big)}{|\gamma'(t)|}=\big(y'(t),\,-x'(t)\big).$$

Then, by the Divergence Theorem,
\begin{equation}\label{SMXC22.12wedmMKAx}\begin{split}& {\mathcal{A}}=\int_\Omega dX=\frac12\int_\Omega {\operatorname{div}}_X\,X\,dX=
\frac12\int_{\partial \Omega} X\cdot\nu=\frac12\int_0^1 
\big(x(t),\,y(t)\big)\cdot\big(y'(t),\,-x'(t)\big)\,dt\\&\qquad\qquad=
\frac12\int_0^1 
\big(x(t)y'(t)-y(t)x'(t)\big)\,dt.
\end{split}\end{equation}

Since the curve under consideration is closed, we can expand its coordinates in Fourier Series
(recall e.g. Theorem~\ref{C1uni}). In this way,
\begin{equation}\label{LASMN02rujgrmtg0o0u}
x(t)=\sum _{k\in\Z}\widehat x_{k}\, e^{2\pi ikt}\qquad{\mbox{and}}\qquad y(t)=\sum _{k\in\Z}\widehat y_{k}\, e^{2\pi ikt}
\end{equation}
and, by virtue of Exercise~\ref{OJSNILCESNpfa.sdwpqoed-23wedf},
\begin{equation}\label{kkidxehatyk}
x'(t)=\sum _{k\in\Z}
{{2\pi ik\,\widehat x_{k}}}\, e^{2\pi ikt}\qquad{\mbox{and}}\qquad y'(t)=\sum _{k\in\Z}{{2\pi ik\,\widehat y_{k}}}\, e^{2\pi ikt}.
\end{equation}
Also, we have uniform convergence of the series above (due e.g. to Theorems~\ref{SMXC22} and~\ref{BASw}) and accordingly~\eqref{SMXC22.12wedmMKAx} leads to
\begin{eqnarray*}
{\mathcal{A}}&=&\frac12\sum _{k,h\in\Z}
\int_0^1 \left(
{{2\pi ik\,\widehat x_h \,\widehat y_{k}}}\, e^{{2\pi i(k+h)t}}
-{{2\pi ik\,\widehat y_h \,\widehat x_{k}}}\, e^{{2\pi i(k+h)t}}
\right)\,dt
\\&=&\pi i\sum _{k\in\Z}
\left({ {k\,\widehat x_{-k} \;\widehat y_{k}}}
-{ {k\,\widehat y_{-k} \;\widehat x_{k}}}
\right).
\end{eqnarray*}
Thus, in light of the reality condition in~\eqref{fasv},
\begin{equation}\label{kkidxehatyk.2}
{\mathcal{A}}=\pi i\sum _{k\in\Z}
k\big(\overline{\widehat x_{k}} \;\widehat y_{k}-\widehat x_{k} \;\overline{\widehat y_{k}}\big).\end{equation}

Furthermore, by~\eqref{LEN:TH} and the arc-length parametrisation for the curve being discussed, we have that
$$ 1={\mathcal{L}}=\int_0^1 |\gamma'(t)|^2\,dt=\int_0^1 \big(|x'(t)|^2+|y'(t)|^2\big)\,dt.$$
For this reason and the Parseval's Identity (recall~\eqref{L2THM.0-02} and~\eqref{kkidxehatyk}), we see that
$$ 1=4\pi^2\sum _{k\in\Z} k^2\,\big(|\widehat x_{k}|^2+|\widehat y_{k}|^2\big).$$

As a result, recalling~\eqref{kkidxehatyk.2},
\begin{equation}\label{JMSiskdcm98PJDIMma.2}
 {\mathcal{L}}^2- 4\pi{\mathcal{A}}=1- 4\pi{\mathcal{A}}=
4\pi^2  \sum _{k\in\Z}\Big(
k^2\,\big(|\widehat x_{k}|^2+|\widehat y_{k}|^2\big)-
ki\big(\overline{\widehat x_{k}} \;\widehat y_{k}-\widehat x_{k} \;\overline{\widehat y_{k}}\big)
\Big).\end{equation}
We also remark that
\begin{eqnarray*}
\big|k\widehat x_{k}-i\widehat y_{k}\big|^2=k^2
|\widehat x_{k}|^2+|\widehat y_{k}|^2
-ki\big(\overline{\widehat x_{k}} \;\widehat y_{k}-\widehat x_{k} \;\overline{\widehat y_{k}}\big).\end{eqnarray*}

We plug this identity into~\eqref{JMSiskdcm98PJDIMma.2} and we obtain that
\begin{equation*}
 {\mathcal{L}}^2- 4\pi{\mathcal{A}}=
4\pi^2  \sum _{k\in\Z}\Big(
(k^2-1) |\widehat y_{k}|^2+\big|k\widehat x_{k}-i\widehat y_{k}\big|^2
\Big)\ge0.\end{equation*}

Also, the above quantity equals zero if and only if~$\widehat y_{k}=0$ for all~$k\in\Z\setminus\{-1,1\}$
and~$k\widehat x_{k}=i\widehat y_{k}$ for all~$k\in\Z$, i.e. if and only if~$\widehat x_{k}=\widehat y_{k}=0$ for all~$k\in\Z\setminus\{-1,1\}$, $\widehat x_{1}=i\widehat y_{1}$, and~$\widehat x_{-1}=-i\widehat y_{-1}$.

This condition is equivalent to rewriting~\eqref{LASMN02rujgrmtg0o0u} in the form
\begin{equation*}\begin{split}&
x(t)=2\big( \Re(\widehat x_{1})\cos(2\pi t)-\Im(\widehat x_{1})\sin(2\pi t)\big)\\
{\mbox{and}}\qquad &y(t)=2\big( \Re(\widehat x_{1})\sin(2\pi t)+\Im(\widehat x_{1})\cos(2\pi t)\big),
\end{split}\end{equation*}
which is the circle~$(x(t))^2+(y(t))^2=4|\widehat x_1|^2$.
\end{proof}

See~\cite[Chapter~4]{MR1412143} and~\cite[Chapters~3 and~4]{MR1970295} for further information on the topics
discussed in these pages.

\begin{exercise}\label{kmsdcJAx.iton4.143.e}
Give a proof of the Isoperimetric Inequality in Theorem~\ref{kmsdcJAx.iton4.143}
by using the Poincar\'e-Wirtinger Inequality of Theorem~\ref{CSLE:A:SEZZ0-1-pwin}.
\end{exercise}

\begin{exercisesk}\label{GEOEX:per} Consider a regular, closed curve of class~$C^\infty$. Assume that this curve is strictly convex, i.e. it encircles a strictly convex set. Prove that there exists a smooth parametrisation~$\gamma:\R\to\R^2$ of this curve such that
\begin{equation}\label{RTSm}
\gamma'(t)=i\mu(t)\,e^{it},
\end{equation}
for some~$\mu:\R\to(0,+\infty)$ periodic of period~$2\pi$.
\end{exercisesk}

\begin{exercisesk}\label{GEOEX} The diameter of a convex curve controls its length. More precisely,
the following result holds true.

Let $I$ be a closed interval and~$\gamma:I\to\R^2$ be the parametrisation of a regular, closed curve of class~$C^\infty$.

The diameter of this curve is
$$ {\mathcal{D}}:= \sup_{t_1,t_2\in I}|\gamma(t_1)-\gamma(t_2)|.$$
The length of this curve is as in~\eqref{LEN:TH}.
Assume that this curve is convex, i.e. it encircles a convex set. 

Then, \begin{equation}\label{GEOEX:INEQ}
{\mathcal{L}}\le \pi {\mathcal{D}}.\end{equation}

Prove this statement.\end{exercisesk}

\begin{exercise}\label{NOAKMSDceikfMKsCfDtgbdcAtVkmdcYI-1} Can one drop the convexity requirement in Exercise~\ref{GEOEX}?
That is: do non-convex curves also fulfil the inequality in~\eqref{GEOEX:INEQ}?.\end{exercise}

\begin{exercise}\label{NOAKMSDceikfMKsCfDtgbdcAtVkmdcYI-2} Can one replace the convexity assumption
in Exercise~\ref{GEOEX} with a star-shapedness hypothesis?\end{exercise}

\begin{exercise}\label{NOAKMSDceikfMKsCfDtgbdcAtVkmdcYI-2b} Is there a convex curve satisfying
the equality in~\eqref{GEOEX:INEQ}?\end{exercise}

\begin{exercise}\label{NOAKMSDceikfMKsCfDtgbdcAtVkmdcYI-3} Do all convex curve satisfy
the equality in~\eqref{GEOEX:INEQ}?\end{exercise}

\section{The Weierstrass Approximation Theorem}\label{CSLE:A:WEIE}

An important aspect in all branch of mathematics is whether or not a given function,
or a class of functions, can be approximated to any degree of accuracy by ``better'', or ``simpler'' functions.

It would be also very convenient if this approximation occurred in a uniform fashion (so that errors are controlled uniformly, and integrations can be freely performed). And it would be highly desirable that the class of approximating functions happened to be so small to constitute a finite dimensional space.

This ideal situation takes place for continuous functions on closed and bounded intervals: such a result,
known as \index{Weierstrass Approximation Theorem} \emph{Weierstrass Approximation Theorem}, goes as follows:

\begin{theorem}\label{WEO:T}
Any continuous function on a closed and bounded interval can be uniformly approximated on that interval by polynomials to any degree of accuracy.

That is, let~$a<b$ and~$f:[a,b]\to\R$ be a continuous function. Let~$\epsilon>0$.

Then, there exists a polynomial~$P$ such that
$$ \sup_{x\in[a,b]} |f(x)-P(x)|<\epsilon.$$\end{theorem}

\begin{proof} Up to replacing~$f(x)$ with~$f\big(2{(b-a)x}+a\big)$
we can reduce to the case in which~$a=0$ and~$b=\frac12$.

Also, we can continuously extend~$f$ in~$\left(\frac12,1\right]$ in such a way that~$f(1)=f(0)$.

Then, we can continuously extend~$f$ from~$[0,1]$ to the whole of~$\R$ making it a continuous function, periodic of period~$1$.

From this, the desired result follows from Exercise~\ref{ojdkfnvioewyr098765rewe67890iuhgvgyuuygfr43wdfgh8bft5xs}.
\end{proof}

An alternative proof of Theorem~\ref{WEO:T} elegantly use the Fej\'er Kernel and goes as follows:

\begin{proof}[Another proof of Theorem~\ref{WEO:T}] Up to replacing~$f(x)$ with~$f\big({(b-a)x}+a\big)$
we can reduce to the case in which~$a=0$ and~$b=1$.

For all~$x\in\R$, set
$$ g(x):=f\big(|\cos(2\pi x)|\big).$$
We notice that~$g$ is continuous and periodic of period~$1$ and that~$g(-x)=g(x)$.

As a result (see Exercises~\ref{smc203e2rf436b.2900rjm4on} and~\ref{920-334PKSXu9o2fg})
we have that the Fourier Sum of~$g$ contains only cosines, that is
$$ S_{N,g}(x)= \frac{a_0}2+\sum_{j=1}^{N} a_j\cos(2\pi jx),$$
for suitable coefficients~$\{a_0,\dots,a_N\}$.

Therefore,
\begin{eqnarray*}&& \frac{1}{N}\sum_{k=0}^{N-1}S_{k,g}(x)=
\frac{1}{N}\sum_{k=0}^{N-1}\left(
\frac{a_0}2+\sum_{j=1}^{k} a_j\cos(2\pi jx)\right)\\&&\qquad=\frac{a_0}2+\sum_{j=1}^{N-1}\sum_{k=j}^{N-1}\frac{a_j}N\,\cos(2\pi jx)=\frac{a_0}2+\sum_{j=1}^{N-1}\frac{a_j\,(N-j)}N\,\cos(2\pi jx).
\end{eqnarray*}

Then, in light of Theorem~\ref{FeJ-th.2}, given~$\epsilon>0$ we can find~$N_\epsilon\in\N$ such that, for all~$N\ge N_\epsilon$,
$$ \epsilon>\sup_{x\in\R}\left|
\frac1N\sum_{k=0}^{N-1}S_{k,g}(x)-g(x)
\right|=\sup_{x\in\R}\left|\frac{a_0}2+\sum_{j=1}^{N-1}\frac{a_j\,(N-j)}N\,\cos(2\pi jx)-g(x)\right|.$$

Hence, writing~$\cos(N\theta)=T_N(\cos\theta)$ for a suitable polynomial~$T_N$ of degree~$N$
(as provided by Exercise~\ref{LMSTNKNOLZOJOqdU2}) we find that
\begin{eqnarray*}
\epsilon&>&\sup_{x\in\R}\left|\frac{a_0}2+\sum_{j=1}^{N-1}\frac{a_j\,(N-j)}N\,T_j\big(\cos(2\pi x)\big)-g(x)\right|\\
&=&\sup_{x\in\R}\left|\frac{a_0}2+\sum_{j=1}^{N-1}\frac{a_j\,(N-j)}N\,T_j\big(\cos(2\pi x)\big)-f\big(|\cos(2\pi x)|\big)
\right|.
\end{eqnarray*}
Now, for all~$y\in[0,1]$ we consider the unique~$x_y\in\left[0,\frac14\right]$ for which~$\cos(2\pi x_y)=y$
and we conclude that
\begin{eqnarray*}
\epsilon&>&\sup_{y\in[0,1]}\left|\frac{a_0}2+\sum_{j=1}^{N-1}\frac{a_j\,(N-j)}N\,T_j(y)-f(y)
\right|,
\end{eqnarray*}
giving the desired result.
\end{proof}

For more information about the topics treated here, see~\cite[Sections~4 and~5]{MR4404761}.

\begin{exercise}\label{LMSTNKNOLZOJOqdU}
Let~$m\in\R$. Prove that, for every~$\theta\in\R$,
$$\cos(m\theta)+\cos\big((m-2)\theta\big)-2\cos\theta\,\cos\big((m-1)\theta\big)=0.$$
\end{exercise}

\begin{exercise}\label{LMSTNKNOLZOJOqdU2} Prove that, for every~$N\in\N$, there exists a polynomial of the form
$$ T_N(t)=\sum_{j=0}^N \kappa_j t^j,$$
with
\begin{equation}\label{eTGBSq01ijm43g-01i24} \kappa_N:=\begin{dcases} 1&{\mbox{ if }}N=0,\\
2^{N-1}&{\mbox{ if }}N\ge1,
\end{dcases}\end{equation}
such that, for every~$\theta\in\R$,
$$ T_N(\cos\theta)=\cos(N\theta).$$
\end{exercise}

\section{\faBomb The Radon Transform}\label{RADOSE}

\begin{figure}[h]
\includegraphics[height=5.6cm]{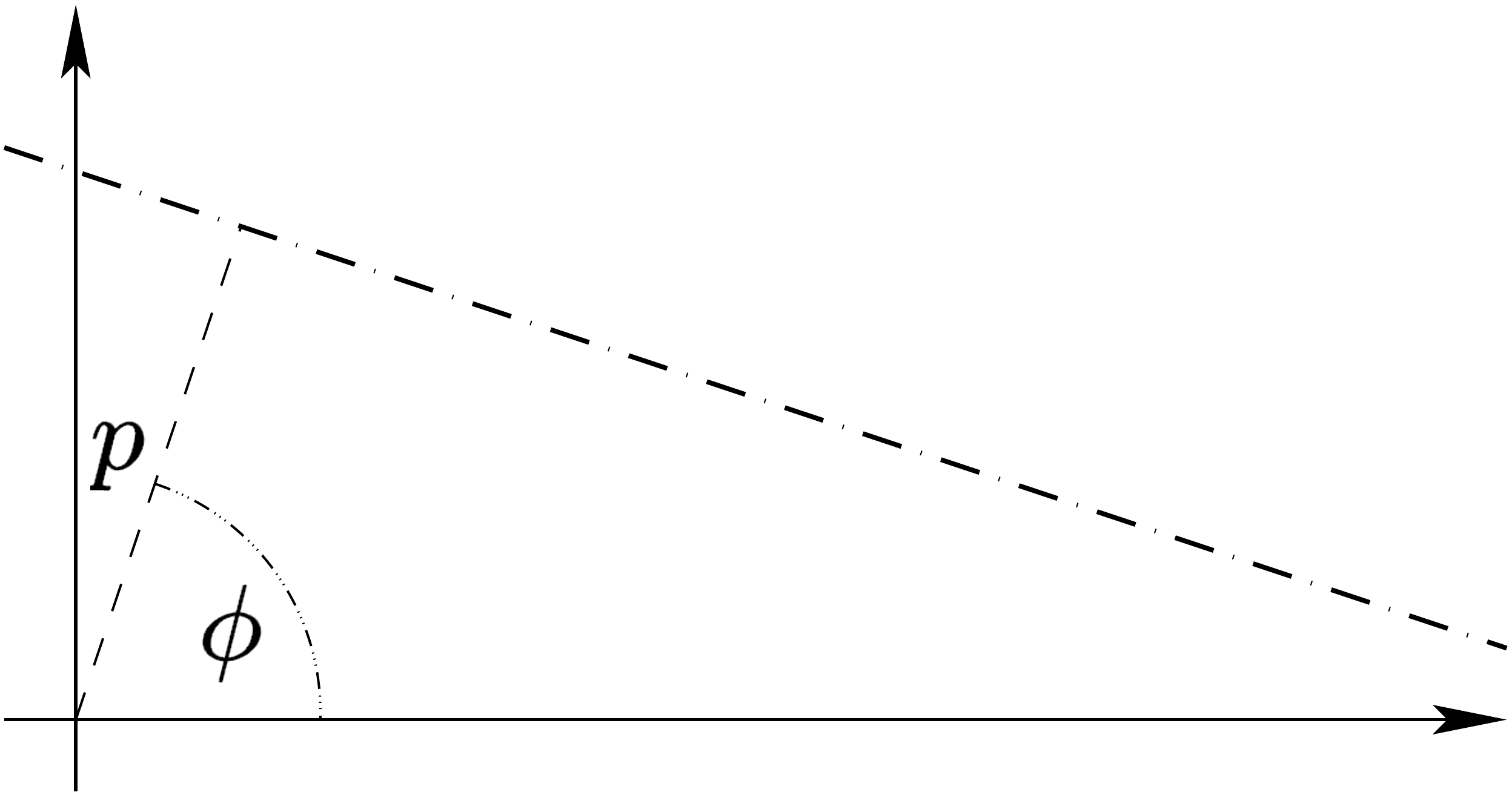}
\centering
\caption{Parametrising straight lines as in~\eqref{RAD:TRA:DE1}.}\label{qTAF.fc3rcbD2.mET24Oikt.9dnsEU0}
\end{figure}

The \index{Radon Transform} Radon Transform\footnote{For simplicity, we focus here on the Radon Transform on the plane,
but higher-dimensional versions are also possible, and useful for practical purposes.
For example, while in the plane the Radon Transform considers a line integral,
in three dimensions it considers an integral over two-dimensional planes. Higher-dimensional generalisations
are possible in the same vein and one can consider integrals
along more complicated objects in the context of integral geometry. There is also a complex analogue of the Radon Transform,
called Penrose Transform, which is used in theoretical physics.

See e.g.~\cite{MR1274701, MR1723736, MR3616276} and the references therein for more information about the Radon Transform and its
various applications.

Let us also mention that one can also apply the projection and reconstruction methods based on the Radon Transform for ``virtual objects'', such as images, e.g. for pattern recognition when searching for similar images in vast archives. In this case, the information stored in the image (for instance, colour or luminosity of pixels) is projected along various directions and gathered in short binary vectors (the collection of which is the virtual analogue of the Radon Transform), see e.g.~\cite{HOANG2012271} and the literature cited therein.} of a function~$f:\R^2\to\R$
is the integral of~$f$ along a given straight line.
Several approaches are possible to formalise this concept. Here we follow a strategy which is nicely
compatible with polar coordinates. Namely, a straight line in the plane is determined
by its distance to the origin, that will be denoted by~$p\in[0,+\infty)$, and the angle formed 
between the normal vector to the straight line
and the horizontal axis of the Cartesian coordinate system, that will be denoted by~$\phi\in\R$, see Figure~\ref{qTAF.fc3rcbD2.mET24Oikt.9dnsEU0}.

More precisely, from the analytic point of view, this straight line has the arc-length parametrisation
\begin{equation}\label{RAD:TRA:DE1} \R\ni t\longmapsto r_{p,\phi}(t):
=(-\sin\phi,\cos\phi)t+(\cos\phi,\sin\phi)p.\end{equation}

With this notation, we introduce the Radon Transform of a function~$f$ as the integral along the straight line above, i.e.
\begin{equation}\label{RAD:TRA:DE} {\mathcal{R}}_f(p,\phi):=\int_{-\infty}^{+\infty}f\big(r_{p,\phi}(t)\big)\,dt,\end{equation}
whenever~$f$ is regular enough for this integral to make sense.

While this may sound a bit abstract so far, the motivation for such a procedure is quite practical,
since it allows to retrieve information on a certain object without having easy access to its full volume.
Suppose for instance that one wants to know what is in our brain (besides mathematics). For this, one may
send an X-ray and detect how much energy passes through our brain. In this sense, the X-ray is mathematically
described by the straight line in~\eqref{RAD:TRA:DE1}, the function~$f$ represents the density of our brain
(or, say, its capacity to block some fraction of the energy of the ray) and the quantity in~\eqref{RAD:TRA:DE}
accounts for the energy absorbed by our brain in this procedure. Notice that this quantity can be experimentally
measured (knowing, for instance, the initial energy of the ray and measuring the residual energy on a screen after the ray passes through our brain).

In practice, one can repeat this measurement by acquiring heaps of X-ray images from different directions
(e.g., from the left side of our skull, from the front, from the right side, etc.): this is precisely\footnote{Actually, not precisely, more or less. For instance the attenuation of X-rays is not linear, but exponential,
in the quantity in~\eqref{RAD:TRA:DE}, see e.g.~\cite[equation~(2.2)]{MR1274701}, but this would not change qualitatively the rough sketch that we have given here.

Also, X-rays are just a cool example, but the same method is used in other contexts as well,
e.g. using ultrasonic pulses, beam of electrons (as in the electronic microscope), or microwave radiation (with applications in astronomy). See again~\cite{MR1274701}
for full details on these topics.}
how a tomography (or CT scan) works. The challenge is now to reconstruct the original function~$f$, which describes the physical properties of our brain, out of all these images (i.e., from all the values~${\mathcal{R}}_f(p,\phi)$ obtained
by sending X-rays of the form~$r_{p,\phi}$). In a nutshell, while
the Radon Transform collects the output of a series of tomographic scans,
inverting the Radon Transform allows one to reconstruct the original density from the projection data.

The inversion of the Radon Transform is a very delicate problem and we surely do not aim at exhausting its complexity
in these notes. This question is deeply related with inverse problems and Fourier methods,
and, to emphasise its link with the topics exposed here, we recall that a sufficiently nice
function in the plane can be reconstructed by its line integrals: for this, we follow an argument introduced in~\cite{zbMATH03191492}, based on writing the function as a Fourier Series in polar coordinates.

To this end, we use polar coordinates~$(\rho,\theta)\in[0,+\infty)\times\R$ to describe points in the plane
and denote by~$F$ the representation of~$f$ in this polar system of coordinates, namely
\begin{equation}\label{DERA6tghk:LA-1} F(\rho,\theta):=f(\rho\cos\theta,\rho\sin\theta).\end{equation}
Since, for each~$\rho\ge0$, the function~$\R\ni\theta\mapsto F(\rho,\theta)$ is periodic of period~$2\pi$,
we can consider its Fourier Series and write
\begin{equation} \label{DERA6tghk:LA-2}
F(\rho,\theta)=\sum _{k\in\Z}\widehat F_{k}(\rho)\, e^{ ik\theta} ,\end{equation}
where
$$\widehat F_k(\rho)={\frac{1}{2\pi}}\int _{0}^{2\pi} F(\rho,\lambda)\,e^{- ik\lambda}\,d\lambda,$$
see Theorem~\ref{ONEGO}, \eqref{fc:PkT-0.1}, and~\eqref{fc:PkT-0.45}.

We also consider~$(p,\phi)\in[0,+\infty)\times\R$ a set of polar coordinates for the function~${\mathcal{R}}_f$ in~\eqref{RAD:TRA:DE}: more specifically, since, for each~$p\ge0$, the map~$\R\ni\phi\mapsto {\mathcal{R}}_f(p,\phi)$ is periodic of period~$2\pi$, we can consider its Fourier Series and write
\begin{equation} \label{wcmRA}
{\mathcal{R}}_f(p,\phi)=\sum _{k\in\Z}\widehat{ {\mathcal{R}}}_{f,k}(p)\, e^{ ik\phi} ,\end{equation}
where
$$\widehat{{\mathcal{R}}}_{f,k}(p)={\frac{1}{2\pi}}\int _{0}^{2\pi} {\mathcal{R}}_f(p,\lambda)\,e^{- ik\lambda}\,d\lambda.$$

Moreover, we recall that, for all~$m\in\N$, the \index{Chebyshev polynomial} $m$th Chebyshev polynomial of the first kind~$T_m$ is defined as
\begin{equation}\label{RADODE1} \R\ni x\longmapsto {\mathcal{T}}_m(x):=\begin{dcases}
(-1)^m\cosh\big( m\arccosh(-x)\big)&{\mbox{ if }}x\in(-\infty,-1],\\
\cos\big( m\arccos x\big)&{\mbox{ if }}x\in(-1,1),\\
\cosh\big( m\arccosh x\big)&{\mbox{ if }}x\in[1,+\infty).
\end{dcases}\end{equation}

With this notation, we have the following reconstruction method
(see also~\cite{zbMATH03191492} for more general results):

\begin{theorem} Let~$f\in C^\infty(\R^2)$ with~$f=0$ outside~$B_1$.

Then, $f$ can be reconstructed by its Radon Transform via~\eqref{DERA6tghk:LA-1}, \eqref{DERA6tghk:LA-2}, and the identity
$$ \widehat F_{k}(\tau)=\begin{dcases}
\displaystyle
-\frac{1}{\pi}\frac{d}{d\tau}
\int_\tau^1 \frac{\tau\,{\mathcal{T}}_k\left(\frac{p}\tau\right)\;\widehat{ {\mathcal{R}}}_{f,k}(p)}{p\,\sqrt{p^2-\tau^2}}\,dp
& {\mbox{ if }}\tau\in[0,1),\\
0&{\mbox{ otherwise.}}\end{dcases}$$
\end{theorem}

\begin{proof} If~$(\rho_{p,\phi,t},\theta_{p,\phi,t})\in[0,+\infty)\times\R$ are the polar coordinates in the plane corresponding
to a point~$r_{p,\phi}(t)$ as in~\eqref{RAD:TRA:DE1}, 
since the vectors~$(-\sin\phi,\cos\phi)$ and~$(\cos\phi,\sin\phi)$ are orthogonal and of unit length
we have that
\begin{eqnarray*}
\rho_{p,\phi,t}^2=t^2+p^2.
\end{eqnarray*}
This also gives that
\begin{eqnarray*}
\rho_{p,\phi,t}=\sqrt{t^2+p^2}.
\end{eqnarray*}

Besides, considering the first coordinate of~$r_{p,\phi}(t)$, we have that
\begin{eqnarray*}
\sqrt{t^2+p^2}\,\cos\theta_{p,\phi,t}=\rho_{p,\phi,t}\,\cos\theta_{p,\phi,t}=-
t\sin\phi +p\cos\phi.
\end{eqnarray*}
Hence, setting
$$\alpha:=\arccos\frac{t}{\sqrt{t^2+p^2}},$$
we conclude that
$$ \sin\alpha=\sqrt{1-\cos^2\alpha}=\sqrt{1-\frac{t^2}{{t^2+p^2}}}=\frac{p}{\sqrt{t^2+p^2}}$$
and that
\begin{eqnarray*}&& \cos\theta_{p,\phi,t}=-\frac{t}{\sqrt{t^2+p^2}}\sin\phi+
\frac{p}{\sqrt{t^2+p^2}}\cos\phi=-\cos\alpha\,\sin\phi+\sin\alpha\,\cos\phi\\&&\qquad=\sin(\alpha-\phi)=
\cos\left(\frac\pi2-\alpha+\phi\right).
\end{eqnarray*}

Similarly, looking at the second coordinate of~$r_{p,\phi}(t)$, 
\begin{eqnarray*}
\sqrt{t^2+p^2}\,\sin\theta_{p,\phi,t}=\rho_{p,\phi,t}\,\sin\theta_{p,\phi,t}=
t\cos\phi +p\sin\phi
\end{eqnarray*}
and thus
\begin{eqnarray*}&&
\sin\theta_{p,\phi,t}=
\frac{t}{\sqrt{t^2+p^2}}\cos\phi +\frac{p}{\sqrt{t^2+p^2}}\sin\phi=
\cos\alpha\,\cos\phi +\sin\alpha\,\sin\phi\\&&\qquad=\cos(\alpha-\phi)=\sin\left(\frac\pi2-\alpha+\phi\right).
\end{eqnarray*}

From these observations, we arrive at
\begin{equation}\label{MSDVBILMSJ9O1NdASHDYS}
\theta_{p,\phi,t}=\frac\pi2-\alpha+\phi=\frac\pi2-\arccos\frac{t}{\sqrt{t^2+p^2}}+\phi.
\end{equation}
Since, for all~$\lambda\in[-1,1]$,
$$ \arccos\lambda +\arccos(-\lambda)=\pi,$$
we also find that\footnote{One could have also obtained~\eqref{MMLSdxcsd823TRfMA} (up to multiples of~$2\pi$) by trigonometric
considerations, looking at the triangle with vertexes the origin, $(\cos\phi,\sin\phi)p$, and~$r_{p,\phi}(t)$,
and comparing it with the (congruent)
triangle with vertexes the origin, $(\cos\phi,\sin\phi)p$, and~$r_{p,\phi}(-t)$.}
\begin{equation}\label{MMLSdxcsd823TRfMA}\begin{split}
\theta_{p,\phi,t}+
\theta_{p,\phi,-t}&=\pi-\arccos\frac{t}{\sqrt{t^2+p^2}}-\arccos\left(-\frac{t}{\sqrt{t^2+p^2}}\right)+2\phi\\&=2\phi.\end{split}\end{equation}

Hence (using the decay of the Fourier coefficients discussed in Section~\ref{DECAY:e:sFOL}
to swap the summation and integral signs) we can rewrite~\eqref{RAD:TRA:DE} in the form
\begin{eqnarray*}
{\mathcal{R}}_f(p,\phi)&=&\int_{-\infty}^{+\infty}f\big(\rho_{p,\phi,t}\,\cos\theta_{p,\phi,t},\,
\rho_{p,\phi,t}\,\sin\theta_{p,\phi,t}\big)\,dt\\&=&\int_{-\infty}^{+\infty}F\big(\rho_{p,\phi,t},\theta_{p,\phi,t}\big)\,dt\\
&=&\sum _{k\in\Z}\int_{-\infty}^{+\infty} \widehat F_{k}(\rho_{p,\phi,t})\, e^{ ik\theta_{p,\phi,t}}\,dt\\
&=&\sum _{k\in\Z}\int_{0}^{+\infty} \Big(\widehat F_{k}(\rho_{p,\phi,t})\, e^{ ik\theta_{p,\phi,t}}+
\widehat F_{k}(\rho_{p,\phi,-t})\, e^{ ik\theta_{p,\phi,-t}}\Big)\,dt\\&=&\sum _{k\in\Z}\int_{0}^{+\infty}
\widehat F_{k}(\rho_{p,\phi,t})\, \Big(e^{ ik\theta_{p,\phi,t}}+ e^{ -ik(2\phi-\theta_{p,\phi,t})}\Big)\,dt.
\end{eqnarray*}

It is now useful to notice that, for all~$\theta\in\R$,
\begin{eqnarray*}
&& e^{ ik\theta}+ e^{ik(2\phi-\theta)}=
e^{ik\phi}\big( e^{ ik(\theta-\phi)}+e^{ ik(\phi-\theta)}\big)=2e^{ik\phi}\,\cos\big(k(\theta-\phi)\big).
\end{eqnarray*}

Therefore, we conclude that
\begin{equation}\label{MSDVBILMSJ9O1NVPRXam7}
{\mathcal{R}}_f(p,\phi)=\sum _{k\in\Z}\int_{0}^{+\infty}
2\,\widehat F_{k}(\rho_{p,\phi,t})\, \cos\big(k(\theta_{p,\phi,t}-\phi)\big)\,e^{ik\phi}\,dt.
\end{equation}

We also deduce from~\eqref{MSDVBILMSJ9O1NdASHDYS} that
$$ \theta_{p,\phi,t}-\phi\in\left[0,\frac\pi2\right],$$
whence
$$ \cos(\theta_{p,\phi,t}+\phi)>0,$$
and that
\begin{eqnarray*}&&
\cos^2\left(\theta_{p,\phi,t}-\phi\right)=
\sin^2\left(\frac\pi2-\theta_{p,\phi,t}+\phi\right)=
1-\cos^2\left(\frac\pi2-\theta_{p,\phi,t}+\phi\right)\\&&\qquad=1-\frac{t^2}{{t^2+p^2}}=\frac{p^2}{{t^2+p^2}}.
\end{eqnarray*}
This gives that
$$ \cos(\theta_{p,\phi,t}-\phi)=\frac{p}{\sqrt{t^2+p^2}}.$$

Thus, recalling~\eqref{RADODE1},
$$ \cos\big(k(\theta_{p,\phi,t}-\phi)\big)=\cos\left(k\arccos\left(\frac{p}{\sqrt{t^2+p^2}}\right)\right)={\mathcal{T}}_k
\left(\frac{p}{\sqrt{t^2+p^2}}\right)$$
and therefore, in light of~\eqref{MSDVBILMSJ9O1NVPRXam7},
$$ {\mathcal{R}}_f(p,\phi)=\sum _{k\in\Z}\int_{0}^{+\infty}
2\,\widehat F_{k}(\rho_{p,\phi,t})\;{\mathcal{T}}_k
\left(\frac{p}{\sqrt{t^2+p^2}}\right)\,e^{ik\phi}\,dt.$$

Thus, comparing with~\eqref{wcmRA},
\begin{eqnarray*}
\widehat{ {\mathcal{R}}}_{f,k}(p)=
2\int_0^{+\infty}\widehat F_{k}(\rho_{p,\phi,t})\; {\mathcal{T}}_k
\left(\frac{p}{\sqrt{t^2+p^2}}\right)\,dt.
\end{eqnarray*}

Now we use the change of variable~$s:=\rho_{p,\phi,t}=\sqrt{t^2+p^2}$, corresponding to~$t=\sqrt{s^2-p^2}$ and
$$ dt=\frac{s}{\sqrt{s^2-p^2}}\,ds.
$$
In this way, we find that
\begin{equation}\label{SFANQPDCCED2}
\widehat{ {\mathcal{R}}}_{f,k}(p)=
2\int_p^{+\infty}\frac{s\,\widehat F_{k}(s)\; {\mathcal{T}}_k
\left(\frac{p}{s}\right)}{\sqrt{s^2-p^2}}\,ds.
\end{equation}

Also, since~$f$ vanishes outside the unit ball, we have that if~$s\ge1$ then
$$\widehat F_k(s)={\frac{1}{2\pi}}\int _{0}^{2\pi} F(s,\lambda)\,e^{- ik\lambda}\,d\lambda=0$$
and therefore~\eqref{SFANQPDCCED2} boils down to
\begin{equation}\label{SFANQPDCCED}
\widehat{ {\mathcal{R}}}_{f,k}(p)=
2\int_p^{1}\frac{s\,\widehat F_{k}(s)\; {\mathcal{T}}_k
\left(\frac{p}{s}\right)}{\sqrt{s^2-p^2}}\,ds.
\end{equation}

For all~$m\in\N$ and~$a\ge1$, we now define
\begin{equation}\label{RADODE1.gfbv} I_m(a):=\int_1^a \frac{a\,{\mathcal{T}}_m(x)\,{\mathcal{T}}_m\left(\frac{x}a\right)}{x\,\sqrt{(a^2-x^2)(x^2-1)}}\,dx.\end{equation}

With this notation, given~$\tau\in(0,1)$,
multiplying both sides of~\eqref{SFANQPDCCED}
by~$\frac{\tau\,{\mathcal{T}}_k\left(\frac{p}\tau\right)}{p\,\sqrt{p^2-\tau^2}}$, after an integration over~$p\in(\tau,1)$
and the substitution~$x:=\frac{p}\tau$,
we see that
\begin{eqnarray*}
\int_\tau^1 \frac{\tau\,{\mathcal{T}}_k\left(\frac{p}\tau\right)\;\widehat{ {\mathcal{R}}}_{f,k}(p)}{p\,\sqrt{p^2-\tau^2}}\,dp&=&
2\int_\tau^1\left(\int_{p}^{1}\frac{\tau s\,\widehat F_{k}(s)\; {\mathcal{T}}_k\left(\frac{p}\tau\right)\;{\mathcal{T}}_k
\left(\frac{p}{s}\right)}{p\,\sqrt{(s^2-p^2)(p^2-\tau^2)}}\,ds\right)\,dp\\&=&
2\int_\tau^1\left(\int_{1}^{s/\tau}\frac{s\,\widehat F_{k}(s)\; {\mathcal{T}}_k(x)\;{\mathcal{T}}_k
\left(\frac{\tau x}{s}\right)}{x\tau\,\sqrt{\left(\frac{s^2}{\tau^2}-x^2\right)(x^2-1)}}\,dx\right)\,ds\\&=&
2\int_\tau^1 \widehat F_{k}(s)\,I_k\left(\frac{s}\tau\right)\,ds.
\end{eqnarray*}
For this reason (see Exercise~\ref{RADOEXX}) we conclude that
$$ \int_\tau^1 \frac{\tau\,{\mathcal{T}}_k\left(\frac{p}\tau\right)\;\widehat{ {\mathcal{R}}}_{f,k}(p)}{p\,\sqrt{p^2-\tau^2}}\,dp=
\pi\,\int_\tau^1\widehat F_{k}(s)\,ds.$$
Taking a derivative in~$\tau$, we obtain the desired result.
\end{proof}

\begin{exercise}\label{RADOEXILE} Compute the Radon Transform of the Gaussian function~$\R^2\ni X\mapsto e^{-|X|^2}$.\end{exercise}

\begin{exercise}\label{RADOEXILE.2} Can one compute the Radon Transform of the Dirac Delta Function?\end{exercise}

\begin{exercise}\label{ANSBEHINGFVA6tYTRAhJAS}
Let~$u\in C^\infty_0(\R^2)$. Prove that the Radon Transform of the Laplacian of~$u$ is equal to the second derivative of the Radon Transform of~$u$ with respect to~$p$, namely $${\mathcal{R}}_{\Delta u}(p,\phi)=\partial^2_p{\mathcal{R}}_u(p,\phi).$$
\end{exercise}

\begin{exercise}\label{RADOEX1} In the notation of~\eqref{RADODE1.gfbv},
prove that, for all~$a\ge1$,
$$I_0(a)=I_1(a)=\frac\pi2.$$
\end{exercise}

\begin{exercisesk}\label{RADOEX2} In the notation of~\eqref{RADODE1.gfbv},
prove that, for all~$m\in\N$ and~$a\ge1$,
$$I_{m+2}(a)=I_m(a).$$
\end{exercisesk}

\begin{exercise}\label{RADOEXX} In the notation of~\eqref{RADODE1.gfbv},
prove that, for all~$m\in\N$ and~$a\ge1$,
$$I_m(a)=\frac\pi2.$$
\end{exercise}

\begin{exercisesk}\label{RADOEXXCR} Given a smooth curve~$\gamma$ in the plane, $p\in[0,+\infty)$, and~$\phi\in\R$, let~$N_\gamma(p,\phi)\in\N\cup\{+\infty\}$ be the number of points at which the curve~$\gamma$ and the straight line~$r_{p,\phi}$ in~\eqref{RAD:TRA:DE1} intersect.

A beautiful formula, which paves the way of integral geometry, reconstructs the length of the curve out of the average of this intersection number. More precisely, the Crofton Formula \index{Crofton Formula} states\footnote{The Crofton Formula
is not just a mathematical curiosity. For example, with that one can
approximately calculate the length of a given curve without having to know its equation, since
a proxy for the integral in~\eqref{LA97529SFOISMOsmSjA} can be obtained by 
considering a family of parallel straight lines at a small distance, then rotating this family by small angles, thus obtaining four families of finely distributed straight lines: for instance, in the spirit of Riemann Sums,
for small~$\epsilon>0$ and large~$m\in\N$,
the quantity $$\sum_{{\phi\in \left\{ 0,\,\frac{\pi}m,\,\frac{2\pi}m,\,\dots,\,\frac{(m-1)\pi}m\right\}}\atop{p\in \epsilon\N}}
\frac{\epsilon\,\pi\,N_\gamma(p,\phi)}{2m}$$
is often used as a proxy for the integral in~\eqref{LA97529SFOISMOsmSjA}.

Alternatively, one can
obtain a proxy for the integral in~\eqref{LA97529SFOISMOsmSjA} by considering a large number of random straight lines,
which provides a very intriguing link between analysis, geometry, and probability (see also Exercise~\ref{GECOB}
for further connections among these fields).

An explicit example implementing this approximation strategy is worked out on pages~45--46 of~\cite{MR394451}
to estimate the length of a DNA molecule pictured by an electron micrograph.

See also~\cite[Corollary~2.5]{MR2167254} 
and~\cite[Chapter~3]{MR3497381}
for higher-dimensional counterparts of the Crofton Formula.

Let us also mention that the regularity assumptions needed for the Crofton Formula can be relaxed,
see e.g.~\cite{MR76373, MR3069683, MR3981295}.}
that the length of~$\gamma$ equals
\begin{equation}\label{LA97529SFOISMOsmSjA}
\frac12\iint_{[0,+\infty)\times[0,2\pi)} N_\gamma(p,\phi)\,dp\,d\phi.\end{equation}

Prove this formula.
\end{exercisesk}

\begin{exercise}\label{RADOEXXISP.02.234}
This is a sanity check for the Crofton Formula in Exercise~\ref{RADOEXXCR}: prove that
if~$\gamma$ is the unit circle, the quantity in~\eqref{LA97529SFOISMOsmSjA} indeed equals~$2\pi$.
\end{exercise}

\begin{exercise}\label{RADOEXXISP.02.234.RSG}
This is another sanity check for the Crofton Formula in Exercise~\ref{RADOEXXCR}: prove that
if~$\gamma$ is a segment, the quantity in~\eqref{LA97529SFOISMOsmSjA} indeed is equal to the length of the segment.
\end{exercise}

\begin{exercise}\label{RADOEXXISP.dp} Let~$\gamma_1$ and~$\gamma_2$ be two smooth, closed curves in the plane.
Assume that~$\gamma_1$ is convex (i.e., it encircles a convex region)
and that it lies in the region encircled by~$\gamma_2$.

Let~${\mathcal{L}}_j$ be the length of~$\gamma_j$, for~$j\in\{1,2\}$, and prove that~${\mathcal{L}}_1\le{\mathcal{L}}_2$.
\end{exercise}

\begin{exercise}\label{RADOEXXISP.dpVAR} 
Let~$\gamma_1$ and~$\gamma_2$ be two smooth curves in the plane.
Assume that~$\gamma_2$ is closed and
that~$\gamma_1$ lies in the region encircled by~$\gamma_2$.

Let~${\mathcal{L}}_j$ be the length of~$\gamma_j$, for~$j\in\{1,2\}$.

Prove that there exists a straight line intersecting~$\gamma_1$ at least~$\left\lfloor\frac{2{\mathcal{L}}_1}{{\mathcal{L}}_2}\right\rfloor$ times.
\end{exercise}

\begin{exercisesk}\label{GECOB}
Georges-Louis Leclerc, Comte de Buffon, once posed this question, which is nowadays known as the \index{Buffon's needle problem} Buffon's needle problem.

Suppose that we have a floor made of parallel strips of wood, each of unit width, and we drop a needle of length~$\ell\in(0,1]$ onto the floor. What is the probability that this needle, randomly tossed, intersects a line between two strips?

Prove that the answer is~$\frac{2\ell}\pi$.
\end{exercisesk}

\section{Filters}\label{SE:FIL:023er-1}

The mathematical theory of filters is more accurately explained in terms of Fourier Transform and can be more effectively understood in terms of the theory of distributions, nevertheless one can somewhat catch a glimpse of this theory also after the modest material that has been put forth in this book.

The filter is a device that stops certain frequencies in a given signal, ideally leaving all the other frequencies unaltered.
For simplicity, we focus here on an ideal \emph{low-pass} filter, i.e. \index{filter, low-pass}
a filter\footnote{A system which removes frequencies lower than a given threshold \index{filter, high-pass}
is called instead \emph{high-pass} filter. A filter only allowing frequencies in a given band is \index{filter, band-pass}
called \emph{band-pass} filter.} which removes completely frequencies higher than a given threshold.
This can be also a convenient procedure to remove unnecessary information from a signal: for instance,
given the fact that the human ear can perceive sounds with frequencies up to 18000 Hz, 
a signal used in telephony does not need to carry frequencies higher than that (actually, the signal used in telephony
is often cut to about 4000 Hz and the results seem adequate to transmit intelligible
speech, compare also with footnote~\label{VOCEHZ} on page~\pageref{VOCEHZ}).

Though the theory of filtering is rather sophisticated, the simplest analogy that we can consider in the setting of Fourier Series  is to apply a low-pass filter~${\mathcal{L}}_{k_0}$ which, given the Fourier Series of some function~$f\in L^1((0,1))$ periodic of period~$1$, maintains only the Fourier Sum up to a given frequency value~$k_0$, that is
$$ {\mathcal{L}}_{k_0}(f):=\sum_{{k\in\Z}\atop{|k|\le k_0}} \widehat f_k \,e^{2\pi i k x}.$$

The natural question is whether there is a ``simple'' operation on the function~$f$ which produces precisely this ideal effect.
For this, the notion of periodic convolution introduced in Exercise~\ref{pecoOSJHN34jfgj} comes in handy.
In particular, we recall that if $$
f\star g(x):=\int_0^1 f(x-y)\,g(y)\,dy,$$
then
\begin{equation}\label{KLfcsckdqmenpfcs9mA}
\widehat{(f\star g)}_k=\widehat f_k\,\widehat g_k.\end{equation}

As a result, if we consider the trigonometric polynomial
$$ \ell_{k_0}(x):=\sum_{{k\in\Z}\atop{|k|\le k_0}} e^{2\pi i k x}=1+2\sum_{k=0}^{k_0} \cos(2\pi k x),$$
we have (see Exercise~\ref{PKS0-3-21}) that the corresponding Fourier coefficients are
$$ \widehat\ell_{k_0,k}:=\begin{dcases}1&{\mbox{ if }}|k|\le k_0,\\0&{\mbox{ if }}|k|>k_0,
\end{dcases}$$
yielding, in light of~\eqref{KLfcsckdqmenpfcs9mA}, that
$$ \widehat{(f\star\ell_{k_0})}_k=\begin{dcases}\widehat f_k&{\mbox{ if }}|k|\le k_0,\\0&{\mbox{ if }}|k|>k_0
\end{dcases}$$
and consequently that the Fourier Series of~$f\star\ell_{k_0}$ is actually a finite sum, equal to
$$ \sum_{{k\in\Z}} \widehat{ (f\star\ell_{k_0})}_k \,e^{2\pi i k x}=\sum_{{k\in\Z}\atop{|k|\le k_0}} \widehat f_k \,e^{2\pi i k x}=
{\mathcal{L}}_{k_0}(f).$$
This confirms the relation between filters and convolution operations.

\begin{figure}[h]
\includegraphics[height=2.26cm]{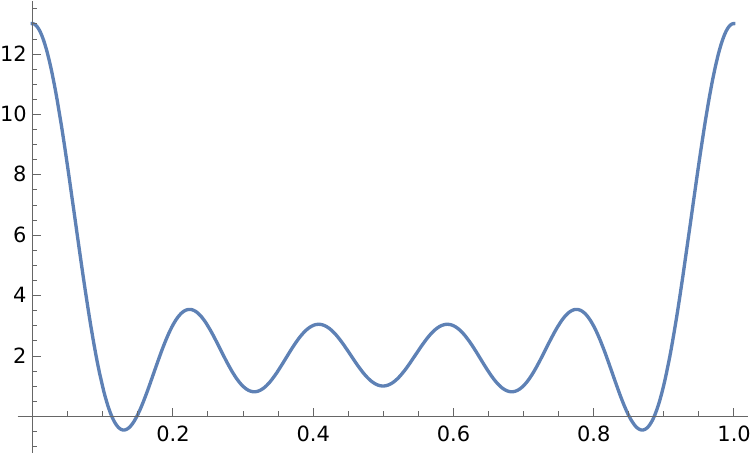}$\,\;$\includegraphics[height=2.26cm]{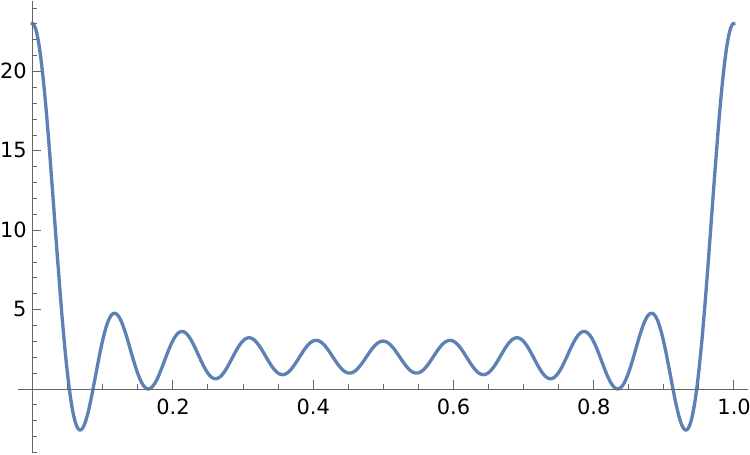}$\,\;$
\includegraphics[height=2.26cm]{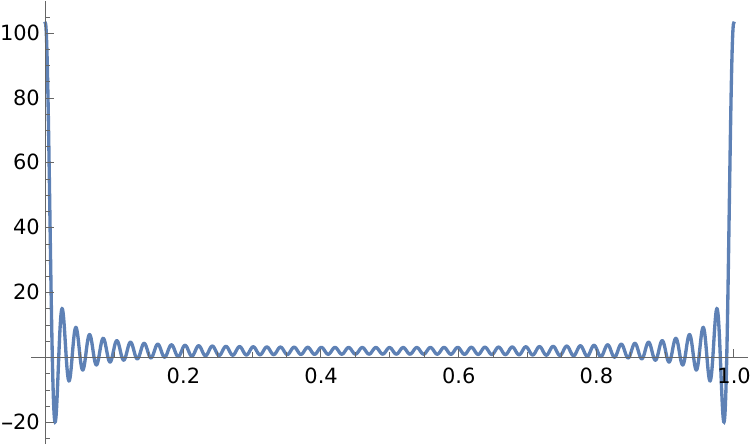}$\,\;$\includegraphics[height=2.26cm]{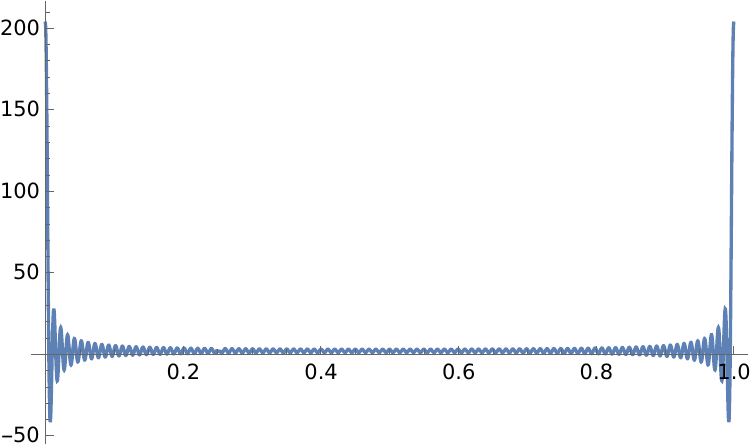}
\centering
\caption{Sketch of the function~$\displaystyle1+2\sum_{k=0}^{k_0} \cos(2\pi k x)$
for~$k_0\in\{5,10,50,100\}$.}\label{c3rcbD2.mET24Oikt.9dnsEU0}
\end{figure}

See Figure~\ref{c3rcbD2.mET24Oikt.9dnsEU0} for a diagram of the filtering function~$\ell_{k_0}$.

We will see in Section~\ref{MCNV:SNK-1md.a} some more application of filters.
Besides, for additional, and more precise, information about filters and related topics see e.g.~\cite{MR1850949}, \cite[Section~3.6]{MR3616140},
\cite[Sections~13, 17, 20, and 23]{PICARDELLO}, and the literature referred to there.

\section{Edge detection}\label{DEDESECTDET}

\begin{figure}[h]
\includegraphics[height=0.32\textwidth]{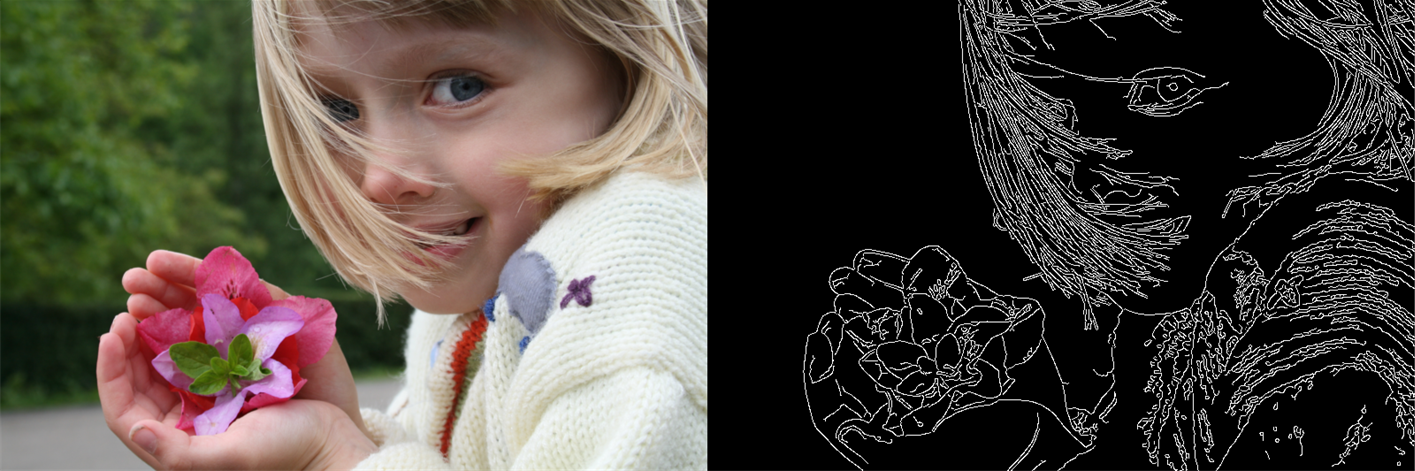}
\centering
\caption{Edge detection applied to a photograph
(work by JonMcLoone, image from Wikipedia,
licensed under the Creative Commons Attribution-Share Alike 3.0 Unported License).}\label{N2BCFOTEDDE1023oBROAM}\end{figure}

Edges and contours are essential ingredients for the processing and recognition of an image.
From a quantitative point of view, the edges appear as curves at which the brightness or intensity of an image changes sharply, see Figure~\ref{N2BCFOTEDDE1023oBROAM}.

In principle, detecting an edge should be a cinch: one could just look at a single row (or column) of the image
and consider, for any point~$x$ of this row (or column), the corresponding intensity, say~$f(x)$, then compute~$f'(x)$
and mark the points where~$f'$ reaches the greatest size.

The issue is that, in practice, the image could be blurred, or affected by some noise. These defects may alter the intensity~$f$ only very mildly, but they could highly impact the derivative of~$f$. Indeed,
a noise could be a continuous function, with small but extremely fast oscillations (e.g., a small multiple of a ``monster'' function
such as in Exercises~\ref{JSLD.023oler-NOMPO} and~\ref{DDEF8AMPMvS1}), and the addition of such noise could make the intensity function nowhere differentiable, thus making it impossible to detect edges.

To get around this sort of difficulties, several methods have been introduced.
The actual procedure of edge detection is rather sophisticated and in practice uses methods of Fourier Transform and wavelets
(see footnote~\ref{DEDESECTDET:F} on page~\pageref{DEDESECTDET:F}) but, to give a clue of it, we highlight here a simple method based on Fourier Series which can be seen as a baby version of more accurate techniques (such as
the Marr-Hildreth algorithm, \index{Marr-Hildreth algorithm} see~\cite{EDGEDETE}).

To this end, for all~$x\in\R$ and~$\epsilon>0$, we consider the \index{Gauss-Weierstrass Kernel} Gauss-Weierstrass Kernel $${\mathcal{G}}(x,\epsilon):=\frac1{\sqrt{\pi\epsilon}}e^{-\frac{x^2}{\epsilon}}.$$

We recall that the Gauss-Weierstrass Kernel can be conveniently used to approximate a given function by a smooth one.
In particular, if~$g\in L^1(\R)$ and continuous at the point~$x$, then
\begin{equation}\label{hGB.klS.mwerfXykLnamsdTYnsy71} \lim_{\epsilon\searrow0}\int_{-\infty}^{+\infty} g(x-\eta)\,{\mathcal{G}}(\eta,\epsilon)\,d\eta=g(x),\end{equation}
see~\cite[page~223]{MR3381284}.

Also, it is worth observing that the Gauss-Weierstrass Kernel satisfies the heat equation, in the sense that \begin{equation}\label{GWHEQmpec.1}
4\partial_\epsilon {\mathcal{G}}(x,\epsilon)=\partial^2_x {\mathcal{G}}(x,\epsilon),\end{equation} see Exercise~\ref{GWHEQmpec.2}.

In the edge detection applications, the Gauss-Weierstrass Kernel can be used to smooth out the singularities produced by the noise. To fit the problem into the framework of Fourier Series, up to scaling we can suppose that the image width is the interval~$[0,1)$ (hence~$f:[0,1)\to\R$) and we identify~$f$ with its periodic extension
(in practice, this can produce some little error on the sides of the image, since we are relating points at two boundary components far away, but we neglect this nuisance here).

We also consider a periodic modification of the Gauss-Weierstrass Kernel via the setting
\begin{equation}\label{PS.mwerfXykLnamsdTYnsy71} \widetilde{\mathcal{G}}(x,\epsilon):=\sum_{k\in\Z} {\mathcal{G}}(x+k,\epsilon)\end{equation}
and we stress that this series converges uniformly, thanks to the spatial decay of the Gauss-Weierstrass Kernel.

Then, in the notation of Exercise~\ref{pecoOSJHN34jfgj}, we consider a smoothed version of the intensity function~$f$ by looking at the function~$f_\epsilon:=f\star\widetilde{\mathcal{G}}(\cdot,\epsilon)$, i.e.
$$ f_\epsilon(x)=\int_0^1 f(x-y)\,\widetilde{\mathcal{G}}(y,\epsilon)\,dy.$$

We notice that, by the uniform convergence in~\eqref{PS.mwerfXykLnamsdTYnsy71},
\begin{eqnarray*}&& f_\epsilon(x)=\sum_{k\in\Z}\int_0^1 f(x-y)\,{\mathcal{G}}(y+k,\epsilon)\,dy=
\sum_{k\in\Z}\int_k^{k+1} f(x-\eta+k)\,{\mathcal{G}}(\eta,\epsilon)\,d\eta\\&&\qquad\qquad\qquad=\sum_{k\in\Z}\int_k^{k+1} f(x-\eta)\,{\mathcal{G}}(\eta,\epsilon)\,d\eta
=\int_{-\infty}^{+\infty} f(x-\eta)\,{\mathcal{G}}(\eta,\epsilon)\,d\eta\end{eqnarray*}
and therefore~$f_\epsilon$ can be thought as a good proxy for~$f$, due to~\eqref{hGB.klS.mwerfXykLnamsdTYnsy71}.

So, we could hope to have some indication about the edges of the image by looking at the location of the maxima of~$f_\epsilon'$. Now, if a point~$x_\epsilon$ maximises~$f'_\epsilon$, then~$f''_\epsilon(x_\epsilon)=0$ and consequently
the edges of the image could be sought among the points~$x_\epsilon$ that satisfy
\begin{equation}\label{pecoOSJHN34jfgj.01woiejfhiruyghtm3hUYng} 0=f''_\epsilon(x_\epsilon)=\int_{-\infty}^{+\infty} f(\eta)\;\partial_x^2{\mathcal{G}}(x-\eta,\epsilon)\,d\eta,\end{equation}
namely the problem of edge detection is reduced to the zero-crossing of a graph obtained by the convolution
with the second derivative of the Gauss-Weierstrass Kernel.

There are also a couple of stratagems which can help simplify the calculations involved in this method. First of all,
in virtue of Exercise~\ref{pecoOSJHN34jfgj}, the Fourier coefficients of~$ f_\epsilon$ (and its derivatives)
can be computed by multiplications of the Fourier coefficients of~$f$ and those of the periodic Gauss-Weierstrass Kernel
(and its derivatives).

Moreover, the second derivative of the Gauss-Weierstrass Kernel in~\eqref{pecoOSJHN34jfgj.01woiejfhiruyghtm3hUYng}
can be replaced by~$4\partial_\epsilon {\mathcal{G}}$, thanks to~\eqref{GWHEQmpec.1}.
This highlights that algorithms of this type basically involve  the subtraction of Gauss-Weierstrass Kernels with
slightly differing widths and they are thereby referred with the name of \index{Difference-of-Gaussians}
Difference-of-Gaussians operators (or DoG for short).

\begin{figure}[h] 
\includegraphics[height=0.29\textwidth]{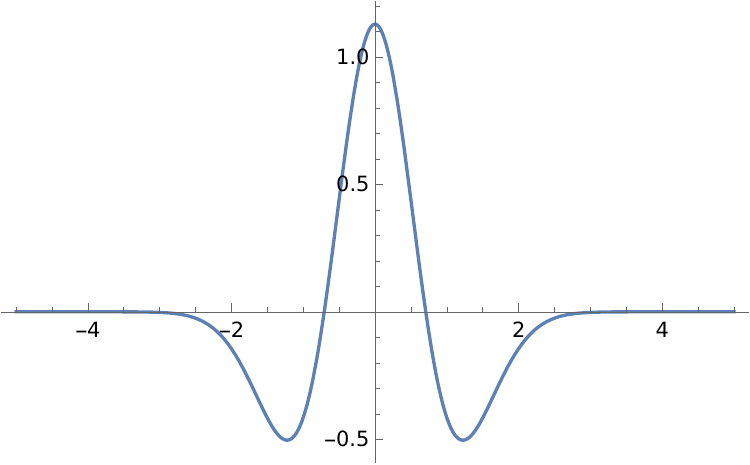}
\centering
\caption{Plot of~$-\partial_x^2{\mathcal{G}}(x,1)$.}\label{1MeX23NBCBROAM}\end{figure}

Interestingly, the second derivative of the Gauss-Weierstrass Kernel coincides (possibly up to a sign and a normalisation)
with a wavelet which is usually called, due to its shape, \emph{Mexican hat},
\index{Mexican hat wavelet} see Figure~\ref{1MeX23NBCBROAM}.

Notwithstanding its classical flavour, this method of edge detection presents some drawbacks and nowadays more accurate algorithms are used in practice. For instance, one limitation of the zero-crossing methods is that they can
locate the edge but by themselves they do not provide any quantitative interpretation of
the ``strength of the edge'' (say, how large the gradient of the intensity is).
Also, since the stationarity of a function is only a necessary, but not a sufficient, condition for a maximum to take place,
zero-crossing methods can also produce ``false edges'', and this error may be severe at curved edges, when additional critical points of the derivative of the intensity function can be created. For more information about edge detection and related topics, see e.g.~\cite{zbMATH05530592} and the literature referred to there.

Let us also highlight an algorithm for edge detection creatively relying on Fourier Series. To catch sight of this technique, which is fully explained in~\cite{MR2416251},
let us suppose that our intensity function~$f$ is periodic of period~$1$ and piecewise smooth,
presenting jump discontinuities at some points~$x_1,\dots,x_\ell\in(0,1)$, with jumps~$\sigma_1,\dots,\sigma_\ell\in\R\setminus\{0\}$, respectively.
More specifically, let us assume that~$x_1<\dots<x_\ell$, that~$f$ is differentiable in all the intervals~$(0,x_1)$, $(x_1,x_2)$, $\dots$, $(x_{\ell-1},x_\ell)$, $(x_\ell,1)$, that~$|f'(x)|\le M$ for all~$x$ belonging to any of these intervals, and that, for each~$j\in\{1,\dots,\ell\}$,
\begin{equation}\label{EDG:GWHEQmpec.2.DE9} \lim_{x\searrow x_j}f(x)-\lim_{x\nearrow x_j}f(x)=\sigma_j.\end{equation}

In this situation, a useful edge detector is provided by the function
\begin{equation}\label{EnaENSIADCac} E_N(x):=\frac{i\pi}{N}\sum_{{k\in\Z}\atop{|k|\le N}} k\widehat f_k\,e^{2\pi i kx},\end{equation}
because, on the one hand,
\begin{equation}\label{EDG:GWHEQmpec.2.DE1}
{\mbox{if~$x\in[0,1)\setminus\{x_1,\dots,x_\ell\}$, then }} \lim_{N\to+\infty}E_N(x)=0,
\end{equation}
and, on the other hand,
\begin{equation}\label{EDG:GWHEQmpec.2.DE3}
{\mbox{for all~$j\in\{1,\dots,\ell\}$, }} \lim_{N\to+\infty}E_N(x_{j})=\sigma_{j},
\end{equation}
see Exercise~\ref{EDG:GWHEQmpec.2.DE4}.

That is, according to~\eqref{EDG:GWHEQmpec.2.DE1} and~\eqref{EDG:GWHEQmpec.2.DE3}, the edge detector deviates from zero precisely at the jump points of the intensity functions, thus revealing the presence of an edge.

\begin{exercise}\label{GWHEQmpec.2}
Check the identity claimed in~\eqref{GWHEQmpec.1}.
\end{exercise}

\begin{exercise}\label{EDG:GWHEQmpec.2.DE4} Prove the claims in~\eqref{EDG:GWHEQmpec.2.DE1} and~\eqref{EDG:GWHEQmpec.2.DE3}.\end{exercise}

\begin{exercise}\label{SASCaCmnmsPcplstqE}
Write the function~$E_N$ in~\eqref{EnaENSIADCac} in
a convenient way to deal with Fourier Series in trigonometric form.
\end{exercise}

\begin{exercise}\label{SASCaCmnmsPcplstqE-vis} Try to visualise the edge detection theory summarised in~\eqref{EnaENSIADCac}, \eqref{EDG:GWHEQmpec.2.DE1}, and~\eqref{EDG:GWHEQmpec.2.DE3} for the example presented in Exercise~\ref{NU90i3orjf:12oeihfnvZMAS}.
\end{exercise}

\section{Denoising}\label{MCNV:SNK-1md.a}

The classical problem of denoising consists of recovering or estimating an ideal image or sound, given only one sample of it. Of course, posed in this generality, the problem turns out to be impossible, but a glimpse of a viable solution comes if one assumes that the ideal signal should not possess too high frequencies, and that the possible presence of these spurious frequencies in the sample is due to noise (see Section~\ref{DEDESECTDET} for similar considerations). 
With this regard,one could use a filter (see Section~\ref{SE:FIL:023er-1}) with the objective of masking parts of the Fourier spectrum of a signal to obtain a denoised one (this method is called in jargon \index{frequency cut-off} \emph{frequency cut-off}). For instance (while more conveniently applied with Fourier Transforms and wavelets), in the context of Fourier Series, this would correspond to replacing a signal modelled by a periodic function~$f$ with a Fourier Sum~$S_{N_0,f}$, with~$N_0\in\N$ suitably chosen to play the role of a frequency threshold telling apart the relevant information and the disturbances.

In a similar vein, one could assume that the dominant components of the signal present a sufficiently high intensity, while noisy data are of small amplitude. In this sense, one can choose an intensity threshold below which all noise components should remain. In the context of Fourier Series, this would produce modified series of the form
$$ \sum_{k\in\Z} \widetilde f_k\,e^{2\pi ikx},$$
where
$$ \widetilde f_k:=\begin{dcases}  \widehat f_k &{\mbox{ if }}| \widehat f_k|\ge \alpha_0,\\ 0&{\mbox{ otherwise,}}\end{dcases}$$
where~$\alpha_0>0$ is appropriately chosen as an amplitude threshold to separate valuable information from noise (this method is called in jargon \index{amplitude thresholding} \emph{amplitude thresholding}).

See e.g.~\cite{zbMATH02043275} and the references therein for more information on these topics.

\section{Who listens to the radio?}\label{TWHOLIKSMDRA902}

The title of this section is taken from a 1978 song by the Australian rock band The Sports. This song claims that logarithms can make someone go berserk (which we totally agree with) and answers the question in the title with ``AM or FM, I listen to both of them''.

\begin{figure}[h]
\includegraphics[height=0.3\textwidth]{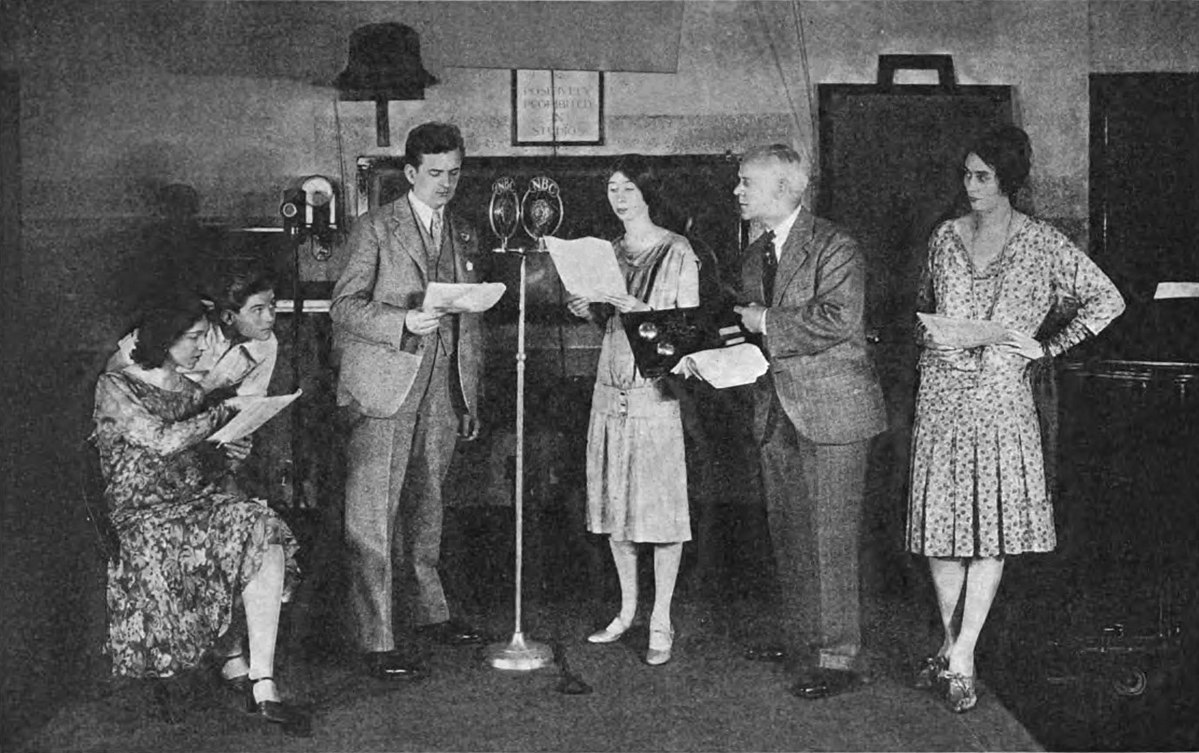}$\;\;$\includegraphics[height=0.3\textwidth]{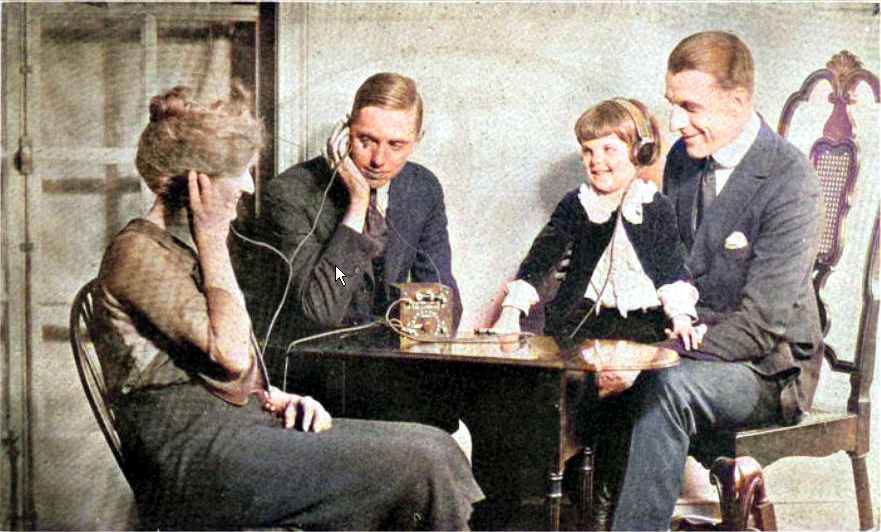}
\centering
\caption{Left: AM radio broadcasting, NBC Studios, New York; right: 1922 advertisement appeared in the ``Radio World'' magazine, depicting an American family listening to a radio (Public Domain Images from Wikipedia).}\label{NBCBROAM}\end{figure}

\begin{figure}[h] 
\includegraphics[height=0.14\textwidth]{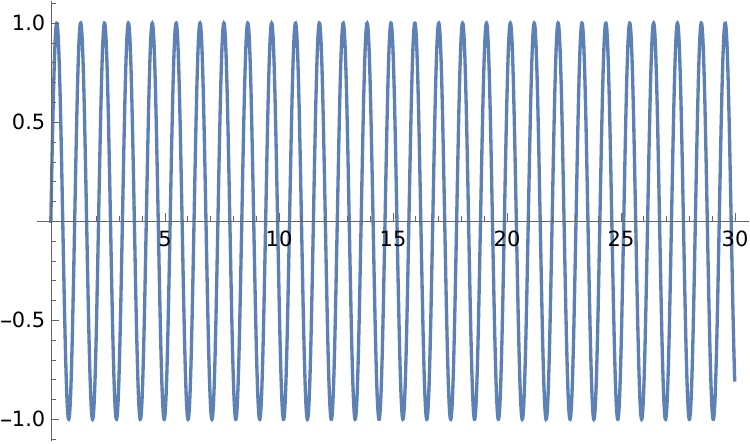}$\;\;\;\;$\includegraphics[height=0.14\textwidth]{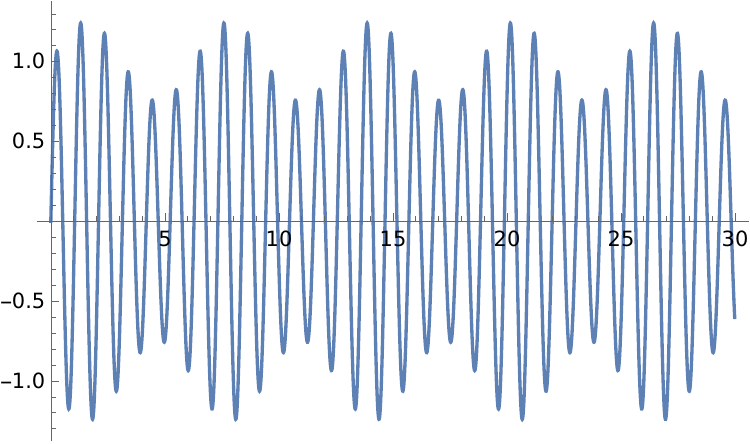}$\;\;\;\;$
\includegraphics[height=0.14\textwidth]{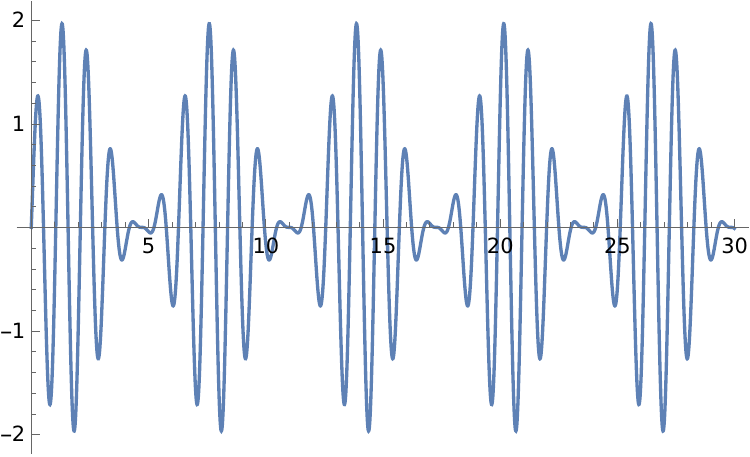}$\;\;\;\;$\includegraphics[height=0.14\textwidth]{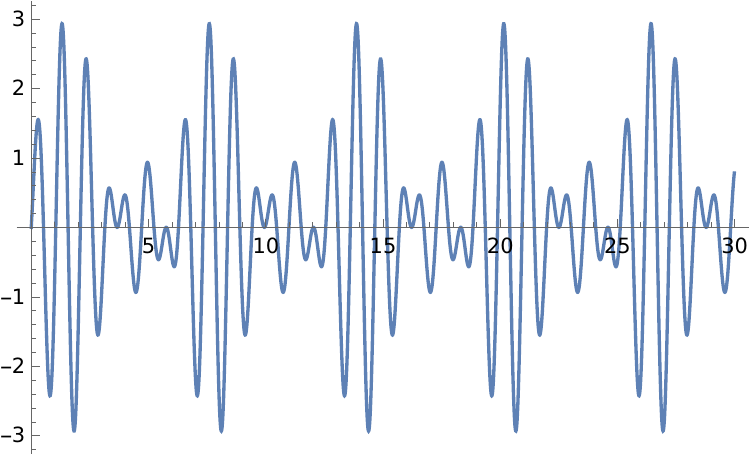}
\centering
\caption{Plots of~$\sin(6 t)$, of~$\left(1+\frac14 \sin t\right)\,\sin(6 t)$, of~$\left(1+\sin t\right)\,\sin(6 t)$,
and of~$\left(1+2\sin t\right)\,\sin(6 t)$.}\label{123NBCBROAM}\end{figure}

So, what's this story of AM and FM? Well, the first transmissions of music and entertainment by radio date back to the beginning of the twentieth century (see Figure~\ref{NBCBROAM}). The technology behind these transmissions relied on the \emph{amplitude modulation} \index{amplitude modulation} (or AM, for short). Roughly speaking, one can imagine that the easiest way to transmit the voice of a speaker consists in ``transforming'' the voice into an electromagnetic wave by tracing its ``highs'' and ``lows''.  
In this sense, the ``amplitude'' is an indication of ``how strong'' a signal is: for instance, computer diagrams of sounds typically show how the amplitude varies over time.
In this way, the easiest way to send out a radio transmission
is to ``modulate'' the amplitude of the signal: the gist is to use a single tone which can be increased or decreased in amplitude, for instance by a sinusoidal ``carrier wave'' of carrier frequency~$\omega$ and variable amplitude~$A(t)$, i.e.$$ A(t)\,\sin(2\pi\omega t).$$
For example, if the amplitude is also modulated by a (small) sinusoidal wave with modulation frequency~$\omega'$, we can consider the signal
\begin{equation}\label{PR0mTRsy1}\big(1+\epsilon \sin(2\pi\omega' t)\big)\,\sin(2\pi\omega t),\end{equation}
see Figure~\ref{123NBCBROAM}.

In practice, the antenna of our AM radio resonates at the appropriate frequency of the selected radio station, essentially ignoring (at least ideally) every other signal, thus producing the appropriate variations of voltage to activate the speakers and (again, ideally)
reproduce the original sound.

Actually, by Prosthaphaeresis, one can rewrite the expression in~\eqref{PR0mTRsy1} in the form
\begin{equation}\label{PR0mTRsy1ZVaC} \sin(2\pi\omega t)+\frac\epsilon2\Big( \cos\big( 2\pi(\omega-\omega')t\big)-\cos\big( 2\pi(\omega+\omega')t\big)\Big),\end{equation}
which separates the role of the carrying wave by two additional ``sidebands'' with frequencies~$\omega-\omega'$ and~$\omega+\omega'$ (i.e., one below and one above the carrier frequency~$\omega$): this specific pattern of amplitude modulation is sometimes called \emph{double-sideband amplitude modulation} \index{amplitude modulation, double-sideband}(or DSBAM, for short).

All in all, the advantages of the amplitude modulation are a simple conceptual idea for transmission and a relatively cheap cost to implement this strategy.

A drawback, however, is that the receiver detects noise in equal proportion to the signal: for instance,
adverse atmospheric conditions, such as thunderstorms, or electric currents and external sources of electric fields, also produce alterations in the amplitude of the radio wave. These perturbations can affect the quality of transmission and reception in amplitude modulation and the variations in the amplitude of the signal caused by environmental effects may produce
significant disturbances.

To address these inconveniences, American electrical engineer Edwin Howard Armstrong considered the idea of sending out
signals with constant amplitude, but with a varying frequency. Indeed, since atmospheric perturbations typically act as a forcing or damping term on the amplitudes, these signals are expected to be less affected by disturbances.

\begin{figure}[h] 
\includegraphics[height=0.24\textwidth]{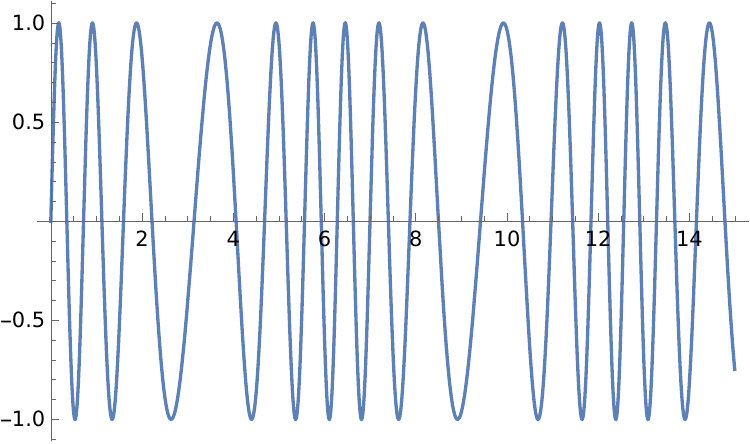}
\centering
\caption{Plot of~$\sin\big(6 t+3 \sin t\big)$.}\label{FM1AMFM23NBCBROAM}\end{figure}

The technique of encoding information in a carrier wave by varying the instantaneous frequency of the wave is called
\emph{frequency modulation} \index{frequency modulation} (FM, for short). In this case, the gist is to leave unaltered the amplitude of, say, a sinusoidal carrier wave and to alter\footnote{We stress that both the carrier and the modulating waves do not need to be necessarily sinusoidal waves: here we presented the case of the sines just for the sake of simplicity, but other types of waves can be taken into account as well. In fact, it is interesting to ``hear'' the sound of frequency modulation on standard waves, such as sawtooth, triangle, square waveforms, see e.g. \url{https://www.youtube.com/watch?v=hQYtRhQPcyY} and other resources easily available online.} the carrier frequency~$\omega$ via a modulating frequency~$\omega'$. This mathematical model, in the simplest situation in which also the frequency modulation occurs through a sinusoidal function, reads
\begin{equation}\label{PR0mTRsy1.b} \sin\big(2\pi\omega t+\epsilon \sin(2\pi\omega' t)\big),\end{equation}
see Figure~\ref{FM1AMFM23NBCBROAM} (the reader can actually compare the frequency modulation formula in~\eqref{PR0mTRsy1.b} with the amplitude modulation formula in~\eqref{PR0mTRsy1}, as well as the sketch of a frequency modulation in Figure~\ref{FM1AMFM23NBCBROAM} with the amplitude modulation in Figure~\ref{123NBCBROAM}).

We recall that the expression in~\eqref{PR0mTRsy1.b} can be rewritten in the form
\begin{equation}\label{PR0mTRsy1ZVaC2}
\sum_{{k\in\Z}} J_k(\epsilon)\,\sin\big(2\pi(\omega+k\omega')t\big),
\end{equation}
where~$J_k$ is the \index{Bessel Function} Bessel Function of order~$k$, see Exercises~\ref{BESSEL-FC-EX1} and~\ref{BESSEL-FC-EX4}
(see also~\cite{MR10746}, \cite{MR193286}, \cite[Appendix~A]{MR3380662},
and~\cite[Section~2.7]{MR3497072} for more information about Bessel Functions). The advantage of the expression in~\eqref{PR0mTRsy1ZVaC2} over the one in~\eqref{PR0mTRsy1.b} is that it makes it more explicit that a frequency modulation can be written as a superposition of standard waves (whose amplitude is determined by the Bessel Functions of the amplitude of the modulating signal).

It is interesting to observe that FM signals carry significantly more bandwidth than AM signals (indeed, only three frequencies appear in~\eqref{PR0mTRsy1ZVaC}, while, formally, in~\eqref{PR0mTRsy1ZVaC2} one has infinitely many frequencies): therefore the channels usually adopted in FM transmissions need to be more separated than the AM ones, and
the frequency used for FM waves is usually considerably higher\footnote{Actually, AM radio \label{RADIOFRE} frequencies range from 530 to 1710 KHz while FM ones range in a much higher spectrum, namely from about 65 to 108 MHz, the precise assigned frequency range varying from country to country.

To set a comparison with frequencies close to our intuition, let us recall
the frequency range of human voice. \label{VOCEHZ}
A female voice typically covers the frequency range from 350 Hz to 17KHz and a male voice from 100Hz to 8KHz.
Hence, when thinking about radio transmissions, the frequency of the carrier wave is typically much higher than
the frequency of the original signal to be sent out.} than the one used for AM waves.

A useful byproduct of the wide bandwidth carried by FM signals is the possibility of storing additional
information in the signal, for example paving the way to stereophony. To give a glimpse of how to send a stereo signal
using frequency modulation, one could pick a certain frequency, called \emph{pilot tone} \index{pilot tone}, which we denote by~$\omega_0$. This frequency (e.g., 19 KHz) acts as a marker, or ``reference signal'', indicating that there is some stereophonic information at 
the second harmonic~$2\omega_0$ of the pilot tone (e.g., 38 KHz): in this way, if a pilot tone is present,
the receiver is triggered to demodulate the stereo information (if no information is stored at the pilot tone frequency~$\omega_0$, then any information in frequencies comparable\footnote{The ``numerology'' of the frequencies used here is merely exemplificative, but it comes from FM channels in the USA. Historically,
channels were set to be~200 KHz apart to prevent overlap and interference. Hence, a deviation of~75 KHz was considered sufficiently inside the frequency gap. As we all know, dividing~75 in half results in~38 KHz, which is used as the frequency where to store the stereo information. The pilot tone would be then half of this frequency, namely at~19 KHz. \label{PILOTFPPOO}

A reasonable frequency span for the signals~$L$ and~$R$ (and so for the signals~$\sigma$ and~$\delta$) was considered to be~15 KHz, therefore the usual information was considered to lie in the range~38$\pm$15 KHz, corresponding precisely to~23--53 KHz. Nothing set in stone anyway.}
to~$2\omega_0$, e.g. 23--53 KHz, is ignored).

Concretely, suppose to have two signals, a ``left'' signal~$L=L(t)$, and a ``right'' signal~$R=R(t)$.
The information stored in~$L$ and~$R$ separately is certainly the same as that stored in the signals~$\sigma:=R+L$ and~$\delta:=R-L$, however it is more convenient to send out these signals rather than the original ones, so that
any ``mono'' receiver, such as a radio not endowed with stereo facilities, will be able to reproduce both the left and the right channels through a single loudspeaker just by reproducing the signal~$\sigma$ (this is called ``mono compatibility'',
which allows mono receivers to reproduce stereo sounds too, and it also permits stereo receivers to switch to mono listening when reception conditions\footnote{For instance, an appropriate reproduction of stereo channels can be jeopardised by
interferences coming from signals reflected by obstacles (such as hills, buildings, and trees) and it is often impossible for
a receiver to separate the primary signal from the secondary reflections of the same signal, which end up arriving out of phase,
causing distortions.  
In these cases, to rectify the signal, a possible solution is to employ a directional antenna designed to accept signals from only one direction (namely, the one of the primary station).}
are less than ideal).

The idea is now to transmit the signals~$\sigma$ and~$\delta$ over the same channel, to shift each signal to a different carrier frequency and then sum them together: this produces a composite signal, often called \emph{multiplex} \index{multiplex}
(or~MPX for short). To catch a sight of this method of transmission, one could describe the multiplexed signal as
\begin{equation} \label{VOCEHZ.e1}\mu(t):= \sigma(t)+\delta(t)\sin(4\pi\omega_0 t)+\kappa\sin(2\pi\omega_0 t),\end{equation}
where~$\kappa$ denotes the amplitude of the pilot tone.

Notice indeed that in~\eqref{VOCEHZ.e1} the sum signal~$\sigma$ remained unprocessed in order to be available for monophonic reception as well and the difference signal~$\delta$ travels on a ``subcarrier'' with frequency~$2\omega_0$.

\begin{figure}[h]
\includegraphics[height=2.8cm]{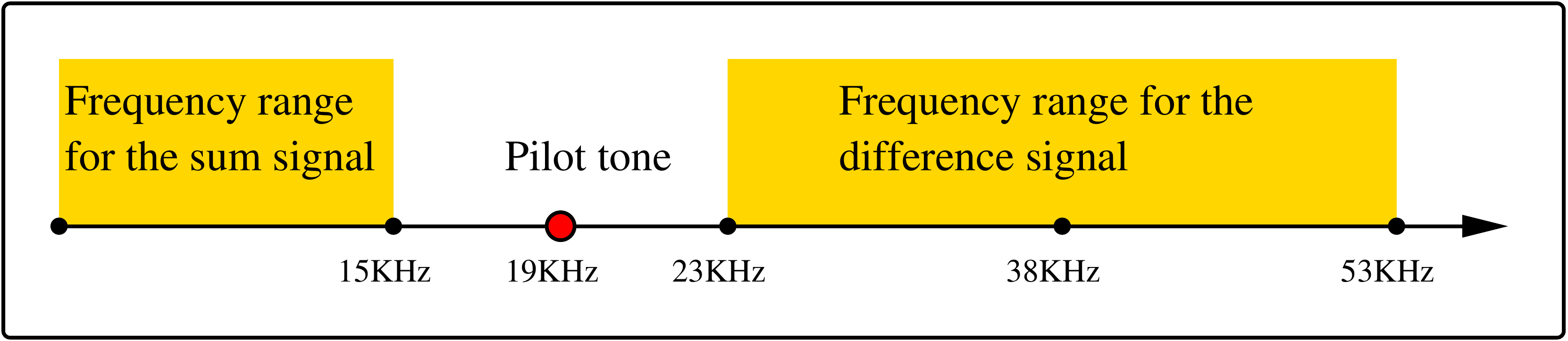}
\centering
\caption{The frequency of the pilot tone remains protected from the ones of the stereo sub-carriers.}\label{PILODET24OikdnsEU0}
\end{figure}

The expression in~\eqref{VOCEHZ.e1} is also consistent with a frequency separation between the different signals:
for instance, if the left and right signals~$L$ and~$R$ (and therefore the signals~$\sigma$ and~$\delta$), have frequency up to
a certain threshold~$\nu$ (e.g., 15 KHz), the information regarding~$\sigma$ is supported in the frequency band~$[0,\nu]$
(e.g., up to~15 KHz) and, by Prosthaphaeresis or algebraic manipulations with complex exponentials, the one related to~$\delta$
lies in the frequency band~$[2\omega_0-\nu, 2\omega_0+\nu]$ (e.g., when the frequency~$\omega_0$
of the pilot tone is 19 KHz, corresponding to the range between~$2\cdot 19-15=23$ and~$2\cdot 19+15=53$ KHz)
and these frequency bands do not interfere with the pilot note, see Figure~\ref{PILODET24OikdnsEU0} and compare with footnote~\ref{PILOTFPPOO}.

The stereo composite signal in~\eqref{VOCEHZ.e1} can be conveniently encoded in the frequency modulation transmission,
as described in~\eqref{PR0mTRsy1.b}: for this, for instance, one could consider the FM signal $$\sin\big(2\pi\omega t+\mu( t)\big).$$

\begin{figure}[h] 
\includegraphics[height=0.54\textwidth]{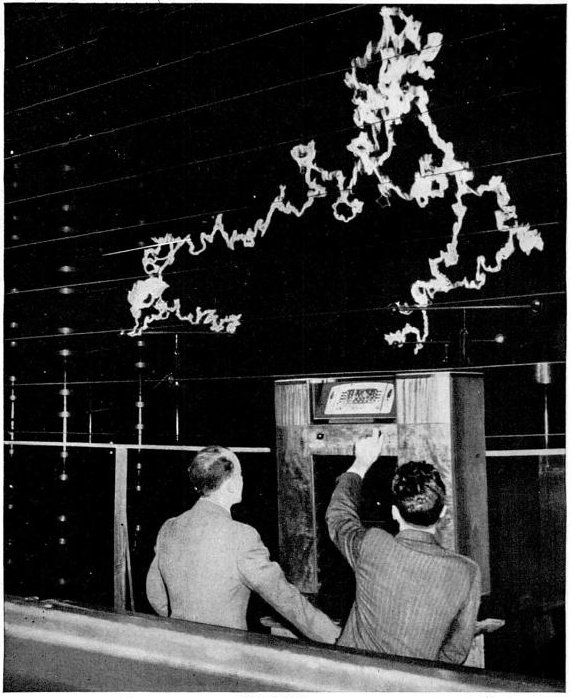}
\centering
\caption{Amplitude modulation versus frequency modulation.}\label{1AMFM23NBCBROAM}\end{figure}

As a historical note, we recall the rather dramatic public relations demonstration, portrayed in Figure~\ref{1AMFM23NBCBROAM}, that took place in 1940: on that occasion,
a lightning-like arc from a million volt three phase-transformer was mounted in a high voltage lab in New York
to act as a source of interference behind a radio having both an AM and FM receiver: while from the AM receiver
only an indistinct static roar could be heard, a music program came through rather clearly from the FM receiver, showcasing the superior audio quality of the frequency modulation and its capability of avoiding noise distortions.

In the same vein, let us remark that when two AM signals are received at the same frequency 
(such as two hypothetical AM radio stations at the same frequency), both will be heard in the receiver, since the amplitude of one will sum up with the amplitude of the other, thus causing an ``interference''.
Instead, if the receiver captures two FM signals, it would naturally amplify the stronger and basically ignore the weaker: in this case,
rather than having an interference caused by the superposition of two signals, the receiver would just select one of the two (this is the reason for which sometimes, while travelling by car, the radio may automatically jump from one station to another
as one drives around). The bottom line
is that AM transmission is prone to interference from other AM stations on the same frequency, while FM transmission provides a natural selection process among competing signals.

Of course, at this point, a natural question is: if FM is so much superior, why do we still have AM too?
Well, about this, one aspect to consider is business:
especially at the beginning, the frequency modulation appeared as a more complex technology, with higher costs.
Moreover, interestingly, an advantage of AM over FM is that typically AM can travel farther than FM,
thanks to the lower frequency of the system.

Let us try to have a glimpse of insight on this phenomenon. First of all, air molecules are prone to ``scatter'' electromagnetic
signals, forcing waves to deviate from their path. Roughly speaking, this deviation is due to the interactions of the wave
with air molecules. In this, not all the frequencies scatter equally, and usually higher frequencies scatter in greater amounts than the lower frequencies: for instance, the sunset looks reddish because the high frequencies (corresponding to the other colours) get scattered more than the low frequencies (corresponding to the red colour).

Though the scattering phenomenon is rather complex, we can imagine that the impact of the scattering on the highest frequencies is the consequence of the scale of the transmission medium with respect to the crests of the wave: as a rule of thumb,
when the frequency of the signal is higher, there are more oscillations between one crest and another,
which thus interact with more gas molecules, being deflected in all directions: this is also the reason for which
the sky appears blue during daytime, given that high frequencies (corresponding to blueish colours)
are deflected so many times that we see them arriving from virtually every direction.

Thus, coming back to AM and FM transmissions, since FM frequencies are considerably higher (recall footnote~\ref{RADIOFRE}),
AM signals typically travel farther than FM signals. This played an important role in the history of the radio communications
since while AM radios could easily cover a vast territory, the FM signals had shorter ranges and needed to be re-tuned while travelling.

\begin{figure}[h]
\includegraphics[height=5.9cm]{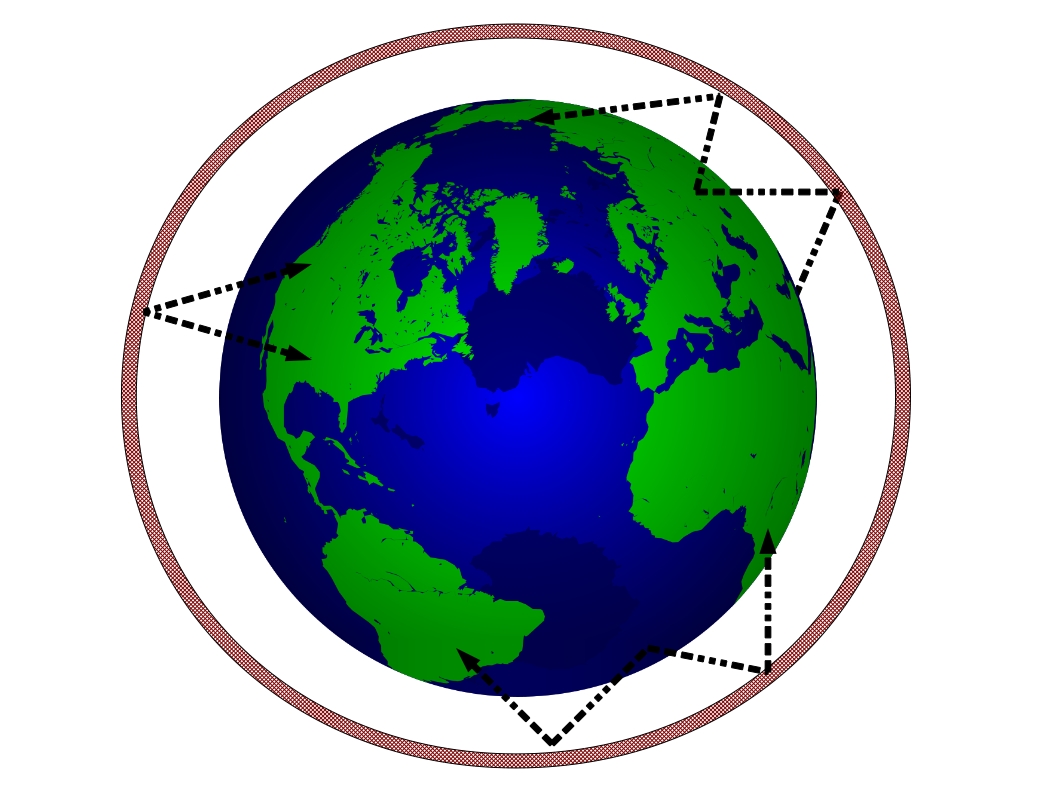}
\centering
\caption{Skywave propagation, with 
radio waves reflecting off the ionosphere (image from Wikipedia
by Kf4yfd, Noldoaran, and Augiasstallputzer$\sim$commonswiki,
licensed under the Creative Commons Attribution-Share Alike 3.0 Unported License).}\label{SKYW01FS.234.45fSIS.K9iuygfy876trfdfghjjnn2meI202}
\end{figure}

Interestingly, the further reach of AM signals is also enhanced especially at night, on the basis of the ionospheric reflection,
which makes the AM signal (but not the FM one) bounce back after arriving at a high layer in the atmosphere, reaching beyond the horizon line, see Figure~\ref{SKYW01FS.234.45fSIS.K9iuygfy876trfdfghjjnn2meI202}.

This \emph{skywave} \index{skywave} effect takes place at the uppermost region of the atmosphere, namely at the level of the ionosphere, where the air is ionised by photons and cosmic rays.
When electromagnetic signals enter the ionosphere, if the peak ionisation is strong enough, the ions and free electrons vibrate and radiate the energy back, and a reflected wave exits the bottom of the ionospheric layer earthwards. In this way, when the radio antenna points the signal close to the horizon, the signal reaches the ionosphere at a very low angle and returns back to the Earth after a considerable distance.

The reflection occurs mostly at night, since the radiation from the sun causes an excessive ionisation in some atmospheric layers and these free electrons and ions end up dampening or absorbing the signal. Also, the reflection of the ionosphere only takes place for lower frequencies: indeed, for frequencies which are too large in comparison with the electron density of the ionosphere, the signal just passes through, is not returned back to the earth, and escapes into space, hence skywave transmission occurs for AM signals, allowing them to bounce off the ionosphere and and get past the horizon line, but not for FM signals.
As a result, FM receiver should ``see'' the transmitter on a straight line of sight, which is restricted by the Earth curvature,
while AM signals can propagate along the surface of the Earth.

In fact, ionospheric propagation was the one initially used by Guglielmo Marconi in his transoceanic radio propagation experiments from Europe to America and led to the discovery of the ionosphere itself (allegedly, Marconi mistakenly believed that
the radio waves were following the curvature of the Earth, but, short after Marconi's experiments, Oliver Heaviside proposed the existence of a suitable high-quote atmospheric layer to allow radio signals to be transmitted around the Earth's curvature).

For more information about radio (and TV) signal transmission, see e.g.~\cite{zbMATH03858920, Wozencraft-Jacobs, LATHIojdc, SINCLAIR201147, MR3616140, NAHIN}.

\begin{exercise}\label{BESSEL-FC-EX1}
For every~$k\in\N$ and~$x\in\R$, let
\begin{equation}\label{BEDIEQG.010} J_k(x):=\sum_{j=0}^{+\infty}\frac{(-1)^j x^{2j+k}}{2^{2j+k}\,j!\,(j+k)!}.\end{equation}
Prove that this series is uniformly convergent in every bounded interval. 

Also, for every~$k\in\Z\cap(-\infty,-1]$, let
\begin{equation}\label{BEDIEQGnega}J_{k}(x):=(-1)^{k}J_{-k}(x).\end{equation}
Prove that, for all~$k\in\Z$,
\begin{equation}\label{BEDIEQG} x^2 J''_k(x)+x\,J'_k(x)+(x^2-k^2)\,J_k(x)=0.\end{equation}

The function~$J_k$ is called\footnote{For completeness, we mention that it
is possible to define Bessel Functions of every order~$\alpha\in\C$. For the humble purposes of this book,
it will suffice to consider integer orders. Also, sometimes the functions~$J_k$ are called Bessel Functions
\emph{of the first kind}, to distinguish them from other Bessel Functions \emph{of the second kind} (anyway, here we will not use other types
of Bessel Functions).} Bessel Function of order~$k$
and~\eqref{BEDIEQG} is called \index{Bessel's Differential Equation} Bessel's Differential Equation.
\end{exercise}

\begin{exercise}\label{BESSEL-FC-EX2fac}
In the notation of Exercise~\ref{BESSEL-FC-EX1}, prove that, for all~$k\in\Z$, the Bessel Function~$J_k$
is an even function when~$k$ is even and is an
odd function when~$k$ is odd.\end{exercise}

\begin{exercise}\label{BESSEL-FC-EX2}
In the notation of Exercise~\ref{BESSEL-FC-EX1}, prove that, for all~$x\in\R$ and~$z\in\C$,
\begin{equation}\label{GENEBESS}
\exp\left(\frac{x}2\left(z-\frac1z\right)\right)=\sum_{{k\in\Z}}J_k(x)\,z^k.
\end{equation}

The function on the left-hand side of~\eqref{GENEBESS} is called the \emph{generating function} of the Bessel Functions.
\end{exercise}

\begin{exercise}\label{BESSEL-FC-EX3}
In the notation of Exercise~\ref{BESSEL-FC-EX1}, prove that, for all~$\beta$, $\theta\in\R$,
\begin{eqnarray*}
&&\cos(\beta\sin \theta)  = J_0(\beta) + 2\sum_{m=1}^{+\infty}J_{2m}(\beta)\,\cos(2m\theta) \\ {\mbox{and }}&&\sin(\beta\sin \theta)  = 2\sum_{m=0}^{+\infty}J_{2m+1}(\beta)\,\sin((2m+1)\theta). \end{eqnarray*}\end{exercise}

\begin{exercise}\label{BESSEL-FC-EX4}
In the notation of Exercise~\ref{BESSEL-FC-EX1},
prove that the quantities in~\eqref{PR0mTRsy1.b} and~\eqref{PR0mTRsy1ZVaC2} are equal.
\end{exercise}

\begin{exercise}\label{BESSEL-FC-EX4-DIFFB.imp}
In the notation of Exercise~\ref{BESSEL-FC-EX1}, prove that
$$ J_0(x)=\frac1\pi\int_0^\pi \cos(x\sin\theta)\,d\theta=\frac2\pi\int_0^1\frac{\cos(xt)}{\sqrt{1-t^2}}\,dt.$$\end{exercise}

\begin{exercisesk}\label{BESSEL-FC-EX4-DIFFBmahUAHSNd01}
Prove that~$J_0$ changes its sign at integer multiples of~$\pi$.\end{exercisesk}

\begin{exercise}\label{BESSEL-FC-EX4-DIFFB}
In the notation of Exercise~\ref{BESSEL-FC-EX1},
prove that~$J_0$ has infinitely many positive zeros and that these zeros are isolated.
\end{exercise}

\begin{exercise}\label{BESSEL-FC-EX4-DIFFB-HYPB}
The displacement~$u=u(x,y,t)$ of a vibrating membrane is described by the equation \begin{equation}\label{BESSEL-FC-EX4-DIFFB-HYPB.e}\partial^2_t u=c\,\big(\partial^2_x u+\partial^2_y u\big),\end{equation}
where~$c>0$ is a parameter taking into account the elasticity of the membrane, $t\ge0$ is time, and~$(x,y)\in\R^2$.

We consider the case of a rotationally symmetric membrane clamped along the unit circle.
Namely, we consider polar coordinates~$(r,\theta)\in[0,1]\times\R$ and write the displacement as~$u(x,y,t)=U(r,t)$. 

The membrane's clamping condition then reads~$U(1,t)=0$. We also assume that the initial velocity of the membrane is null, namely~$\partial_t U(r,0)=0$.

Prove that a solution of this problem is given by
\begin{equation}\label{BESSEL-FC-EX4-DIFFB-HYPB.ui} U(r,t)=J_0(\lambda r)\,\cos(\sqrt{c}\,\lambda t),\end{equation}
where~$\lambda$ is a positive zero of~$J_0$.
\end{exercise}

\section{Overdrives, clipping, distortion, and ringing artifacts}\label{GI:CONSECT}

The Gibbs phenomenon showcased in Section~\ref{SEC:GIBBS-PH} is a truly fascinating mathematical concept, but it does present some concrete pitfalls. Indeed, a general feature of this phenomenon is to present some overshooting
and we now briefly discuss what sort of hindrances (but sometimes also benefits) can arise from that.

\begin{figure}[h]
\includegraphics[height=0.28\textwidth]{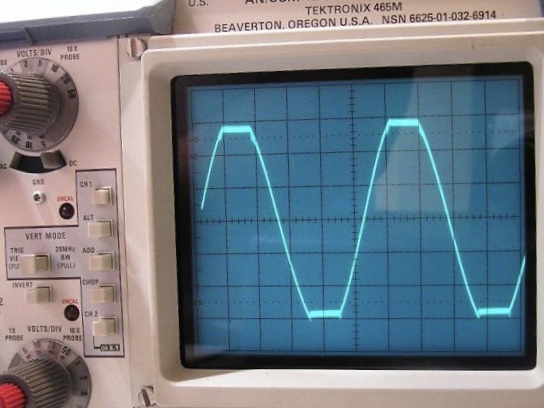}$\qquad$
\includegraphics[height=0.28\textwidth]{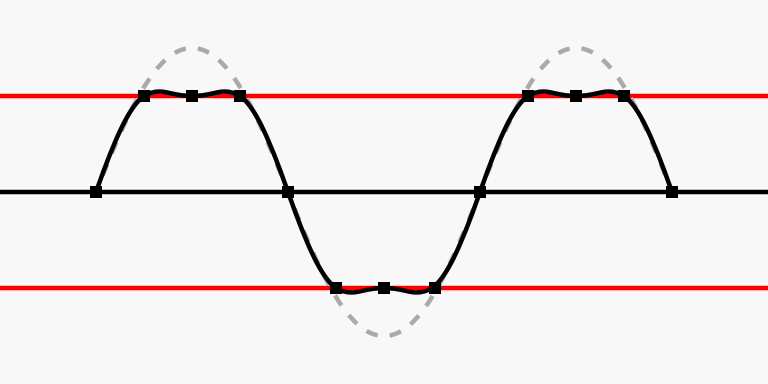}
\centering
\caption{Occurrences of clipping (images from Wikipedia; left: assumed work by Lgreen$\sim$commonswiki,
Creative Commons Attribution-Share Alike 3.0 Unported License; right:
work by Gutten p{\aa} Hemsen, Creative Commons Attribution-Share Alike 4.0 International License).}\label{RSTCNdqU.1ZKRETPx}
\end{figure}

One typical instance of overdrive occurs in audio signal and it is called \index{clipping} \emph{clipping}. In a nutshell,
when an audio signal is pushed beyond the maximum limit of the reproducing device,
it forces the speaker to a higher output voltage than what its power supply can produce.
The sound gets thereby altered or distorted, may produce an unpleasant listening experience,
and can potentially damage the device. In these situations, the signal usually ``cuts'' at the maximum capacity of the amplifier:
the result is a distorted waveform, often similar to a square wave (recall Exercise~\ref{SQ:W} and see Figure~\ref{RSTCNdqU.1ZKRETPx}).

\begin{figure}[h]
\includegraphics[height=0.18\textwidth]{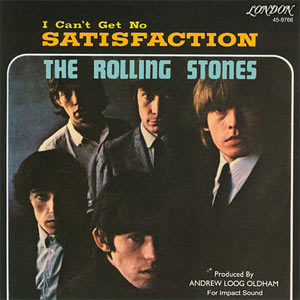}$\,$
\includegraphics[height=0.18\textwidth]{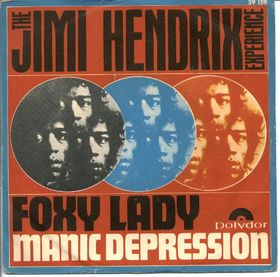}$\,$
\includegraphics[height=0.18\textwidth]{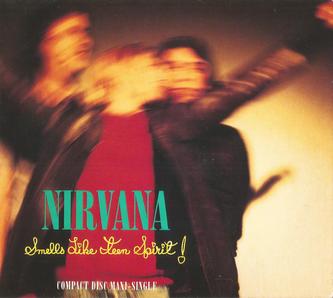}$\,$
\includegraphics[height=0.18\textwidth]{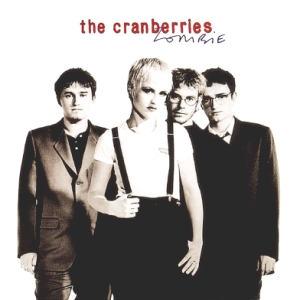}$\,$
\includegraphics[height=0.18\textwidth]{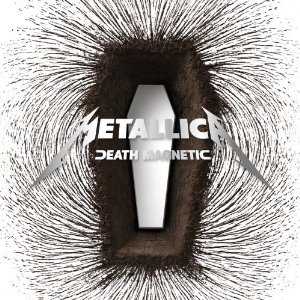}
\centering
\caption{Covers of some famous musical hits (images from Wikipedia;
it is believed that the use of the images of such covers
to illustrate the audio recording in question qualifies as fair use in terms of copyright).}\label{RSTCNZKRETPx}
\end{figure}

But let us remain positive: every cloud has a silver lining and clipping can turn out sometimes to be a precious resource.
Indeed, several famous songs  reached their popularity also thanks to an original and bewitching use of clipping, overdrive, ``fuzz'', and all sort of distortion effects, especially on electric guitars, which, in various forms, emphasise upper frequencies, possibly cutting away the middle frequencies, often transforming the audio signal into an almost square or triangular or sawtooth wave (or all sort of funny shape waves): just to mention few examples of popular songs, {\em (I Can't Get No) Satisfaction} by the Rolling Stones (1965), {\em Foxy Lady} by Jimi Hendrix (1967), {\em Smells Like Teen Spirit} by Nirvana (1991), {\em Zombie} by the Cranberries (1994), the album {\em Death Magnetic} by Metallica (2008),
and so on, see Figure~\ref{RSTCNZKRETPx}

Perhaps the introduction of distortion in music was a ``happy accident'', dating back to 1961: it seems that Glenn Snoddy, a studio engineer of a recording session, accidentally produced some overtones after having inadvertently turned up the volume way too high, also damaging some of the amplifiers, but, allegedly, rather than fixing the ``problem'', the musicians liked the results. After this, Snoddy invented a box, the Maestro FZ-1 Fuzz-Tone, released in 1962, to intentionally distort the guitar sound (without damaging the amplifiers).

We stress that the Gibbs phenomenon is one of the possible causes of clipping, but clipping may be the outcome of many other situations (including serendipity in the anecdote about Snoddy). In the same vein, the occurrence of clipping is not limited to audio signals. For instance, a visual clipping takes place when an image, or part of it, present ``desaturated''
pure white areas due to the clip of all colour components, or to regions with inaccurate colour reproduction,
caused by the clip of an individual colour channel.

\begin{figure}[h]
\includegraphics[height=0.48\textwidth]{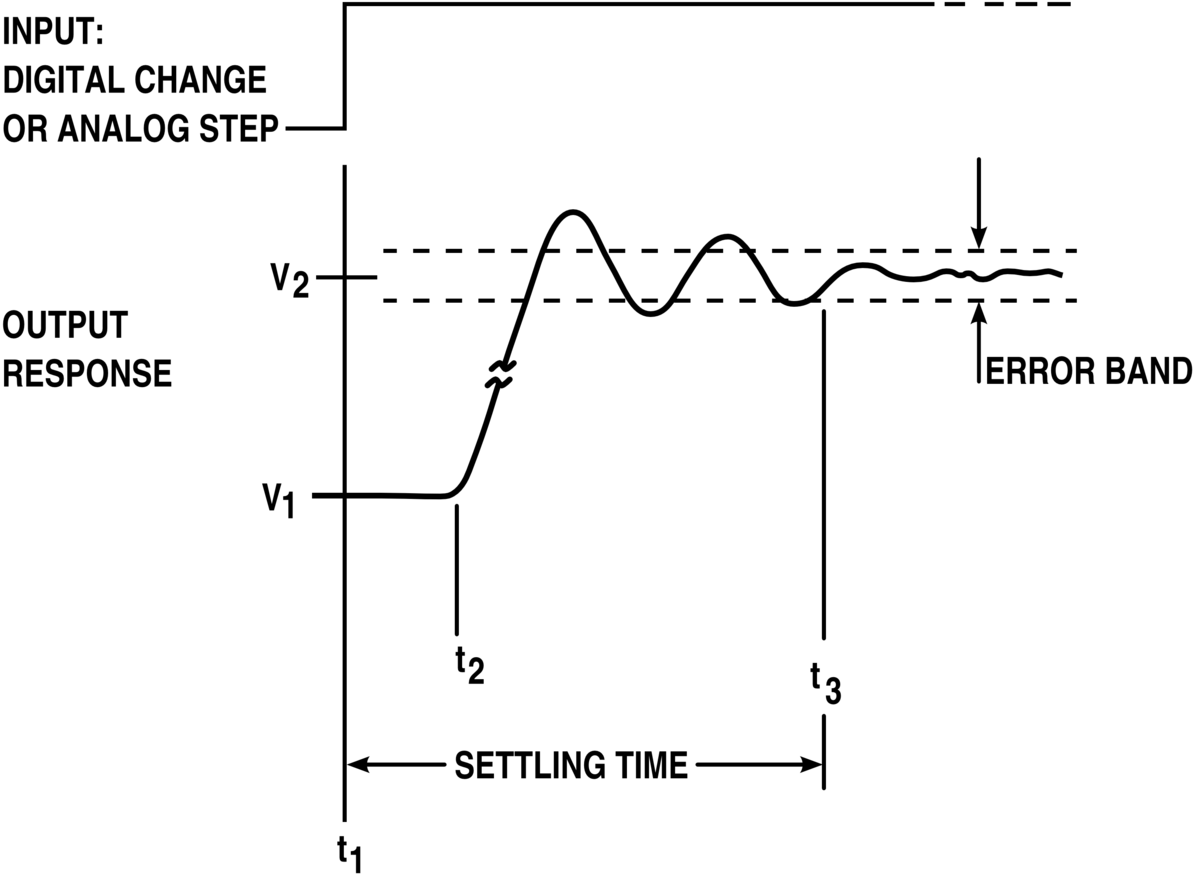}
\centering
\caption{Diagram of a step response of a filter to a discontinuous input, producing oscillatory patterns giving rise to ringing artifacts (image from Wikipedia, original work by Howard K. Schoenwetter, released into the Public Domain).}\label{Rin8STCNZKRETPx}
\end{figure}

Another spurious effect that, in the realm of sound and image processing, can be caused by overshooting patterns such as the ones produced by the Gibbs phenomenon is related to \index{ringing artifacts} the \emph{ringing artifacts}.
These artifacts can be caused, for instance, by the response of a filter to
a sudden, ideally discontinuous, input: in view of the Gibbs phenomenon, the filter may 
replace the sharp transition of the input with an output signal oscillating at a fading rate, see Figure~\ref{Rin8STCNZKRETPx}.

\begin{figure}[h]
\includegraphics[height=0.28\textwidth]{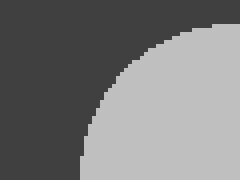}$\qquad$\includegraphics[height=0.28\textwidth]{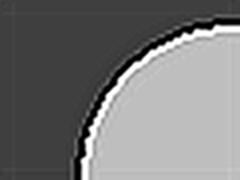}
\centering
\caption{Left: an illustration of a quarter-circle of light grey on dark grey.
Right: same image with ringing artifacts produced by upscaling (work by Nils R. Barth,
image from Wikipedia, made available under the Creative Commons CC0 1.0 Universal Public Domain Dedication).}\label{Rin8STCNZKRETPx.1}
\end{figure}

The name of these artifacts is due to to the ``bands'' or ``rings'' appearing in an image near sharp edges and caused
by the spurious oscillations that we have just described, see Figure~\ref{Rin8STCNZKRETPx.1}.

\begin{figure}[h]
\includegraphics[height=0.34\textwidth]{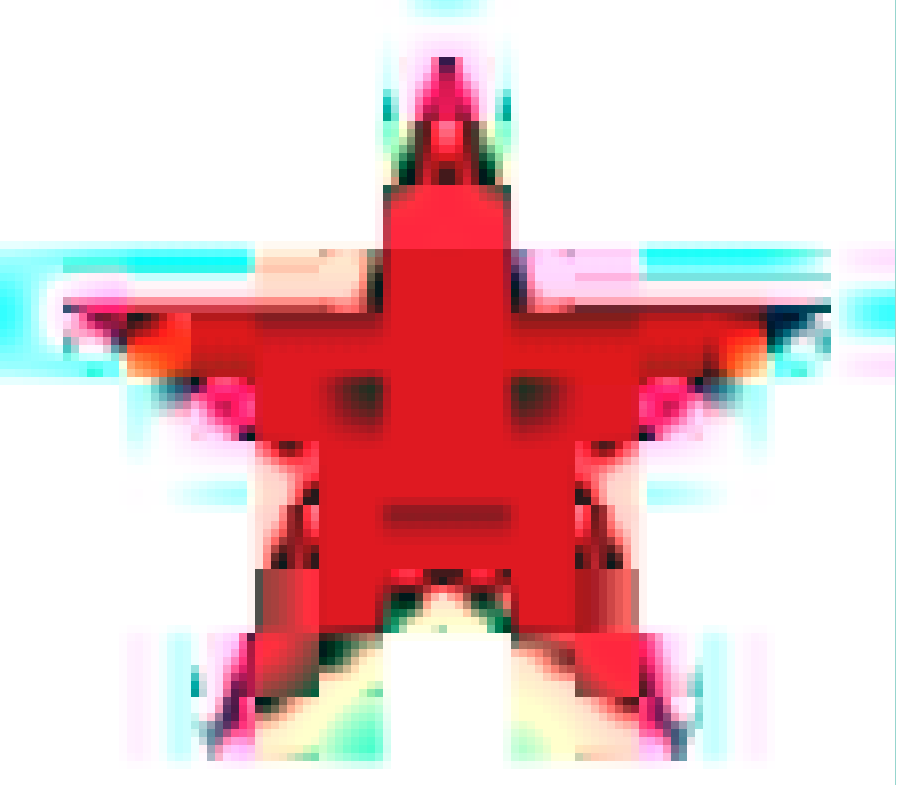}
\centering
\caption{A red star that is distorted by JPEG artifacts (work by Roger McLassus,
image from Wikipedia, made available under the  Creative Commons Attribution-Share Alike 3.0 Unported License).}\label{Rin8STCNZKRETPx.2}
\end{figure}

Ringing artifacts at sharp transitions can also occur as a result of a JPEG compression, since these algorithm use
discrete approximations of Fourier-type quantities (compare with Section~\ref{ODCABRBAGN})
and a block segmentation of the image. This discretisation can lose accuracy in high frequency components
and the poor rendering can be enhanced when the edges of the image do not match well with the subdivision into blocks,
see Figure~\ref{Rin8STCNZKRETPx.2}.

In audio signals, the ringing artifacts can occur especially in the presence of percussion instruments.
In this case, the sharp struck of the percussion is complemented by spurious echos and pre-echos.

\section{Linear motions on tori}\label{LIMOTORSE}

Linear motions on tori are a fascinating topic combining arithmetic, geometry, number theory, mathematical analysis, and mathematical physics.
Let's see how.

Given~$p\in\R^n$ and~$\omega\in\R^n$, we can consider the linear motion
$$ \phi_0(t):=p+\omega t.$$
This is just the description of a trajectory moving on a straight line.
We can also consider the projection of this motion on the torus, namely
$$\phi(t):=\{\phi_0(t)\}=\{p+\omega t\},$$
where~$\{\cdot\}$ denotes the fractional part, as introduced on page~\pageref{FRaPARTADE},
and the notation has to be read ``componentwise'', i.e., if~$p=(p_1,\dots,p_n)$ and~$\omega=(\omega_1,\dots,\omega_n)$,
$$ \phi(t)=\Big( \{p_1+\omega_1 t\},\dots,\{p_n+\omega_n t\}\Big).$$

Motions of this type are especially relevant in the study of many dynamical systems (including celestial mechanics and plasma confinement).
For this, it is often desirable to understand ``how much'' of the torus gets filled by the linear motion~$\phi$.
A prototypical result in this setting goes as follows:

\begin{theorem}\label{NORESAT}
Suppose that
\begin{equation}\label{NORESA}
{\mbox{for every~$k\in\Z^n\setminus\{0\}$ we have that~$\omega\cdot k\ne0$.}}
\end{equation}

Then, the orbit of~$\phi$ is dense on the whole torus.

More precisely, 
\begin{equation}\label{NORESA2}
\overline{\big\{ \phi(t),\;t\ge 0\big\}}=[0,1)^n.\end{equation}
\end{theorem}

In the jargon, when condition~\eqref{NORESA} is fulfilled, one says that the frequency~$\omega$ is
\index{nonresonant}
\emph{nonresonant}, or \emph{rationally independent}. \index{rationally independent}
Theorem~\ref{NORESAT} is beautiful, since it links the ``arithmetic'' condition~\eqref{NORESA}
to a dynamical and topological property, namely~\eqref{NORESA2}.

\begin{figure}[h]
\includegraphics[height=2.94cm]{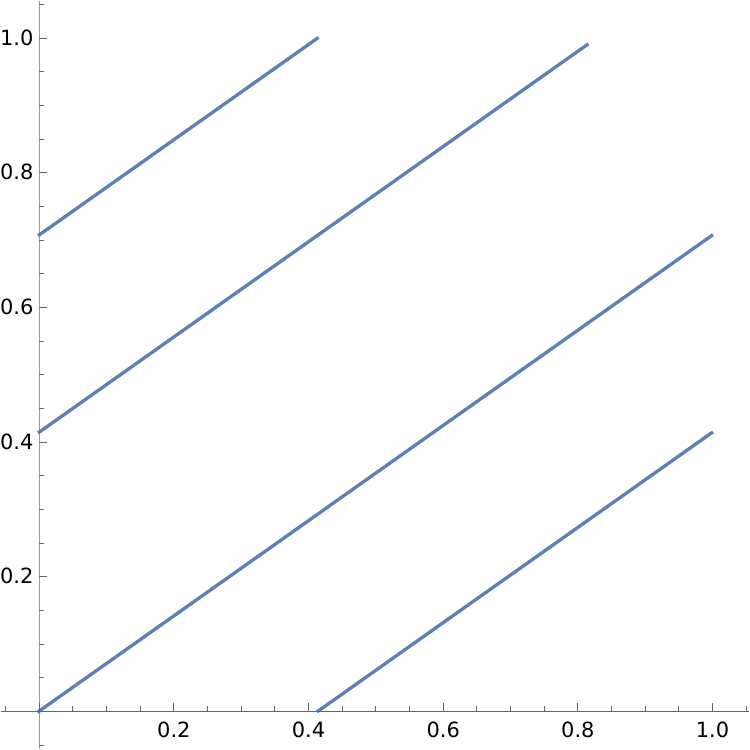}$\,\;$\includegraphics[height=2.94cm]{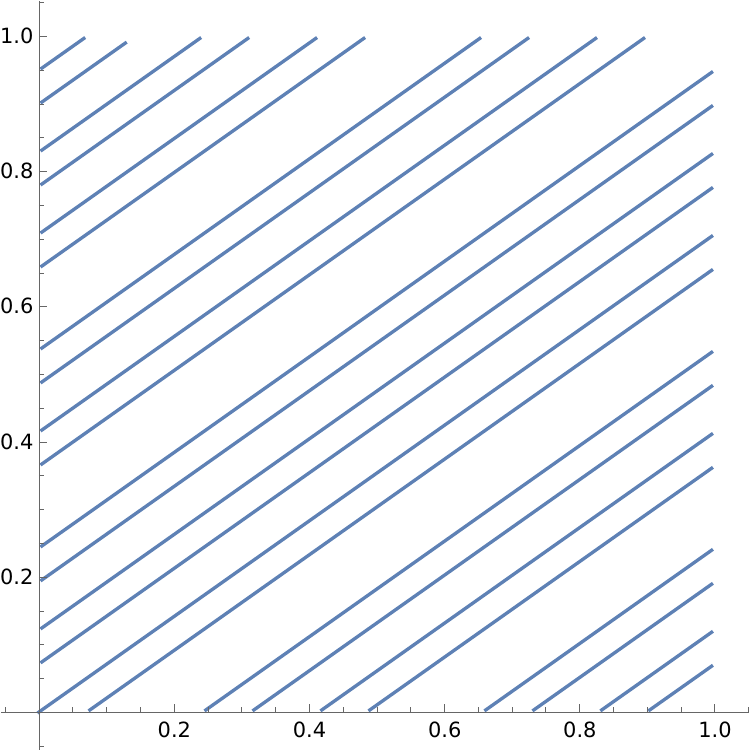}$\,\;$
\includegraphics[height=2.94cm]{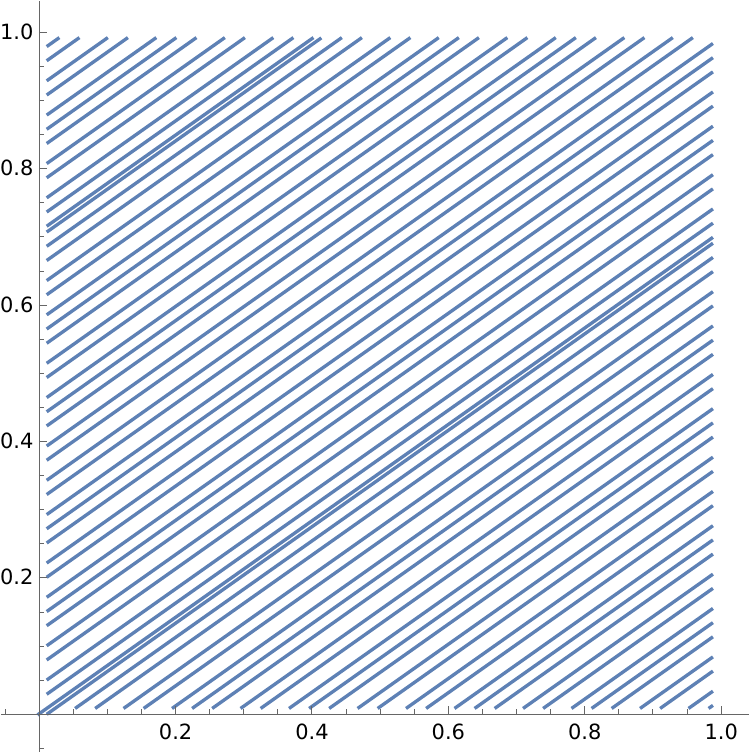}$\,\;$\includegraphics[height=2.94cm]{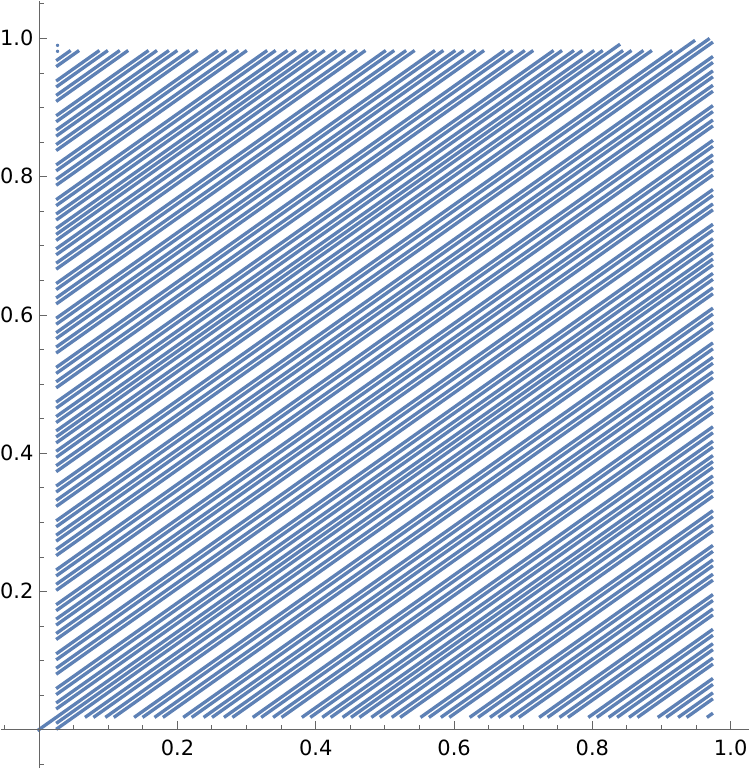}
\centering
\caption{The evolution of a rationally independent linear motion on a torus.}\label{DET24OikdnsEU0}
\end{figure}

\begin{figure}[h]
\includegraphics[height=2.94cm]{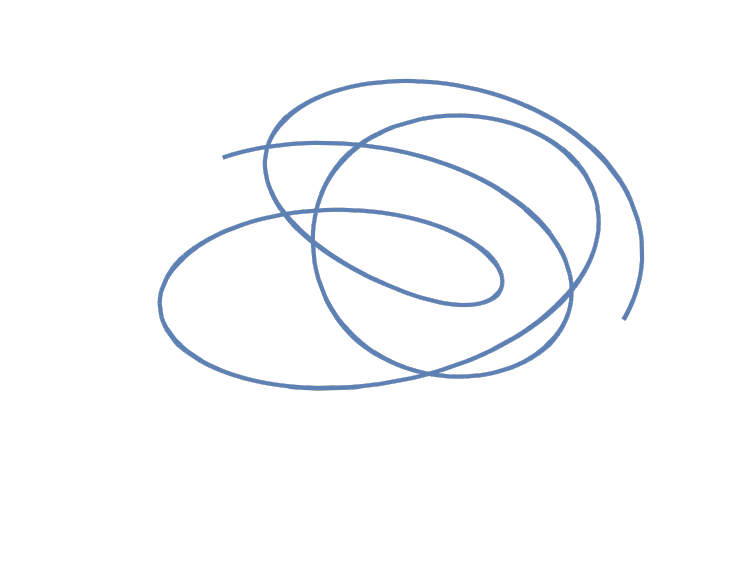}$\,\;$\includegraphics[height=2.94cm]{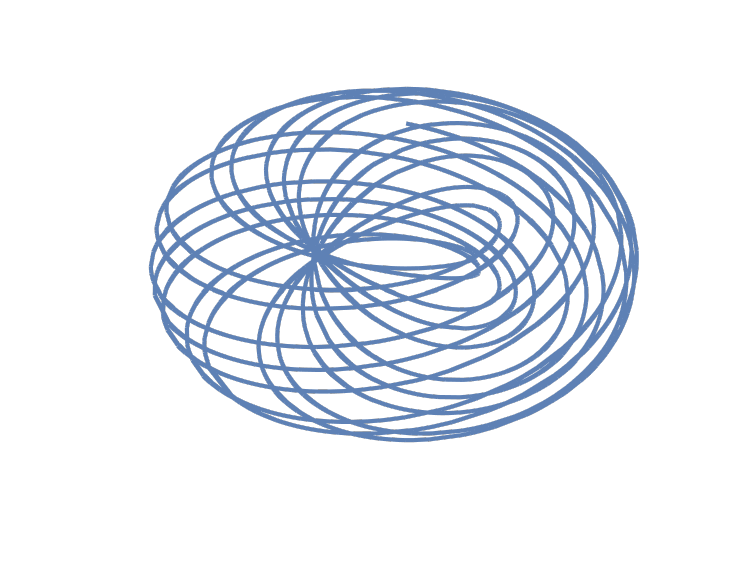}$\,\;$
\includegraphics[height=2.94cm]{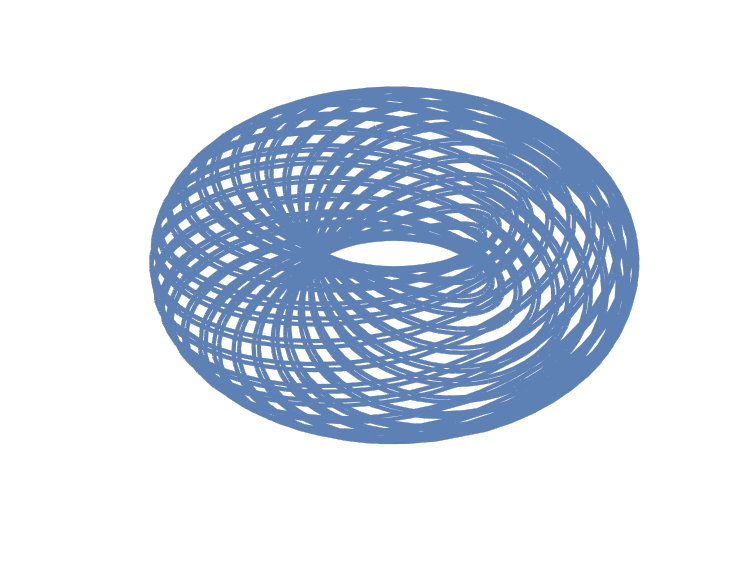}$\,\;$\includegraphics[height=2.94cm]{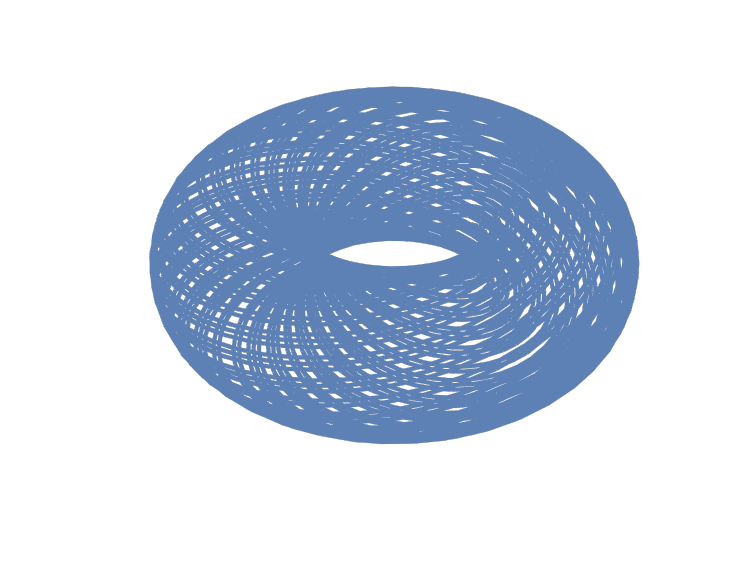}
\centering
\caption{The evolution of a rationally independent linear motion on a torus.}\label{DET24OikdnsEU}
\end{figure}

See Figures~\ref{DET24OikdnsEU0} and~\ref{DET24OikdnsEU} to visualise how rationally independent linear motions densely fill the torus.

To prove Theorem~\ref{NORESAT}, we introduce some elements of ergodic theory, with the aim of comparing
\emph{time averages} \index{time averages}
and \emph{phase averages}. \index{phase averages}
Specifically, in the notation of page~\pageref{PKSMDODKUIOIYFDYPIJBN8.17y2ert0012ihfbbgrce4gtdk},
given a $\Z^n$-periodic function~$f:\R^n\to\R$, we define the time average of~$f$ with respect to the linear motion~$\phi$ by
\begin{equation}\label{CqTr64RQL} \langle f \rangle_{\mbox{time}}:=\lim_{T\to+\infty}\frac1T\int_0^T f(\phi(t))\,dt,\end{equation}
provided that this limit exists.

Roughly speaking the time average of~$f$ measures the long-time asymptotic mean value of~$f$ recorded on the travelled trajectory.

Also, the phase average is simply the average of~$f$ on its fundamental domain, that is
\begin{equation}\label{CqTr64RQL.lamerpwoekfm} \langle f \rangle_{\mbox{phase}}:=\int_{(0,1)^n} f(x)\,dx.\end{equation}
Notice that this quantity is completely independent of the dynamics.

For instance, in the particular case in which we consider a measurable set~$A\subseteq[0,1)^n$ and the corresponding characteristic function
$$ \chi_A(x)=\begin{dcases} 1&{\mbox{ if }}x\in A,\\0&{\mbox{ otherwise}},
\end{dcases}$$
we see that~$\langle \chi_A \rangle_{\mbox{phase}}$ is simply the Lebesgue measure of~$A$
(or, if one prefers, since~$[0,1)^n$ has unit measure, the proportion of the measure of~$A$ in the fundamental domain under consideration).

Instead, $\langle \chi_A \rangle_{\mbox{time}}$ measures ``how often the trajectory stays in the set~$A$''.
Remarkably, for rationally independent linear flows on tori, these two conceptually very different quantities coincide:

\begin{theorem} \label{ERGOB}
If~$f$ is a $\Z^n$-periodic and continuous function and~$\omega$ is a rationally independent frequency, then the limit in~\eqref{CqTr64RQL} exists for every initial datum~$p$ and
$$ \langle f\rangle_{\mbox{time}}=\langle f \rangle_{\mbox{phase}}\,.$$
\end{theorem}

\begin{proof} Let~$\epsilon>0$. Let~$f_\epsilon\in C^\infty(\R^n)$ be $\Z^n$-periodic and such that~$\|f-f_\epsilon\|_{L^\infty(\R^n)}\le\epsilon$ (see e.g.~\cite[Theorem~9.8]{MR3381284} for this type of approximation). Let also~$N_\epsilon\in\N$ so large that~$\|f_\epsilon-S_{N_\epsilon,f_\epsilon}\|_{L^\infty(\R^n)}\le\epsilon$ (see Theorem~\ref{SELINOPOR}).

In this way,
\begin{equation}\label{ERGOB.eqRoscjh} \|f-S_{N_\epsilon,f_\epsilon}\|_{L^\infty(\R^n)}\le2\epsilon.\end{equation}
Then, writing
$$ S_{N_\epsilon,f_\epsilon}(x)=\sum_{{k\in\Z^n}\atop{|k|\le N_\epsilon}} \widehat f_{\epsilon,k}\,e^{2\pi i k\cdot x},$$ 
and observing that, by the periodicity of the complex exponential,
$$ e^{2\pi i k\cdot \phi(t)}=e^{2\pi i k\cdot (p+\omega t)},$$
we have that, for all~$T>0$,
\begin{eqnarray*}&&
\left|\frac1T\int_0^T f(\phi(t))\,dt-\int_{(0,1)^n} f(x)\,dx\right|\\&&\qquad\le
\left|\frac1T\int_0^T S_{N_\epsilon,f_\epsilon}(\phi(t))\,dt-\int_{(0,1)^n} S_{N_\epsilon,f_\epsilon}(x)\,dx\right|
+4\epsilon\\&&\qquad=
\left|\frac1T
\sum_{{k\in\Z^n}\atop{|k|\le N_\epsilon}}\widehat f_{\epsilon,k}
\int_0^T e^{2\pi i k\cdot \phi(t)}\,dt-\sum_{{k\in\Z^n}\atop{|k|\le N_\epsilon}}\widehat f_{\epsilon,k}\int_{(0,1)^n}e^{2\pi i k\cdot x}\,dx\right|
+4\epsilon\\&&\qquad=
\left|\frac1T
\sum_{{k\in\Z^n}\atop{|k|\le N_\epsilon}}\widehat f_{\epsilon,k}
\int_0^T e^{2\pi i k\cdot (p+\omega t)}\,dt-\widehat f_{\epsilon,0}\right|
+4\epsilon\\&&\qquad=
\left|\frac1T
\sum_{{k\in\Z^n}\atop{0<|k|\le N_\epsilon}}\frac{\widehat f_{\epsilon,k}\,e^{2\pi i k\cdot p}\big(
e^{2\pi i k\cdot \omega T}-1\big)}{2\pi k\cdot\omega}\right|
+4\epsilon\\&&\qquad\le
\frac{2N_\epsilon}{\pi T}\,\sup_{{k\in\Z^n}\atop{0<|k|\le N_\epsilon}}\frac{|\widehat f_{\epsilon,k}|}{| k\cdot\omega|}
+4\epsilon,
\end{eqnarray*}
and we stress that the above supremum is finite thanks to the assumption that~$\omega$ is rationally independent.

We can now send~$T\to+\infty$ and deduce that
$$ \lim_{T\to+\infty}\left|\frac1T\int_0^T f(\phi(t))\,dt-\int_{(0,1)^n} f(x)\,dx\right|\le4\epsilon.$$
Taking~$\epsilon$ as small as we wish, we obtain the desired result.
\end{proof}

With this, we can now complete the proof of Theorem~\ref{NORESAT}.

\begin{proof}[Proof of Theorem~\ref{NORESAT}] Suppose, for a contradiction, that~\eqref{NORESA2}
is not satisfied. Then, there exists a ball~$B$ contained in~$[0,1)^n\setminus
\overline{\big\{ \phi(t),\;t\ge 0\big\}}$. We pick~$f\in C^\infty_0(B,[0,1])$, with~$f=1$ in a set of positive measure.

Hence, by Theorem~\ref{ERGOB},
\begin{eqnarray*}&&
0<\int_B f(x)\,dx=\langle f \rangle_{\mbox{phase}}=
\langle f\rangle_{\mbox{time}}
\\&&\qquad=\lim_{T\to+\infty}\frac1T\int_0^T f(\phi(t))\,dt\le
\lim_{T\to+\infty}\frac1T\int_0^T \chi_B(\phi(t))\,dt=0.
\end{eqnarray*}This is a contradiction and the proof of~\eqref{NORESA2} is thereby complete.
\end{proof}

For more information about linear motions on tori see e.g.~\cite[Section~11.7]{MR3185347}. See also~\cite[Chapter~II, Section 5]{MR493419} and the references therein for further reading about ergodic theory. See also Section~\ref{PelapeammilonS}
for an application of the topics highlighted here to number theory.

\begin{exercise}\label{PelapeammilonDF} Let~$\gamma\in\R\setminus\Q$ and~$f:\R\to\R$ be continuous and periodic of period~$1$.

Prove, using\footnote{A slightly different approach to Exercise~\ref{PelapeammilonDF}
will be outlined in Exercise~\ref{Pelapeammilon5}.} Fourier methods, that
$$\lim_{K\to+\infty}\frac1K\sum_{k=0}^{K-1} f(\gamma k)=\int_0^1f(x)\,dx.$$
\end{exercise}

\begin{exercise}\label{PelapeammilonDF2} Can one drop the assumption that~$f$ is continuous in Exercise~\ref{PelapeammilonDF}?
\end{exercise}

\begin{exercisesk}\label{DIpFA:0}
Given~$\tau\in(0,+\infty)$, let\footnote{Exercise~\ref{DIpFA:0} is about the so-called \index{Diophantine vectors} Diophantine vectors, which are of great importance in dynamical systems see e.g.~\cite[Section~5.9]{MR698947} and~\cite[Section~12.3]{MR3185347}.

{F}rom the arithmetic perspective, not only these Diophantine vectors are
rationally independent (or nonresonant, according to the notation in~\eqref{NORESA}), but also they satisfy
a ``quantitative'' bound from the resonance.}
\begin{equation}\label{DEFIDIOF}\begin{split} {\mathcal{D}}_{n,\tau}&:=\Bigg\{ \omega\in\R^n{\mbox{ s.t. there exists $\gamma>0$ s.t. for every $k\in\Z^n\setminus\{0\}$}}\\&\qquad\qquad\qquad{\mbox{ we have that }}
|\omega\cdot k|\ge\frac{\gamma}{|k|^\tau}\Bigg\}.\end{split}\end{equation}

Prove that
\begin{itemize}
\item[{(i).}] When~$\tau\in(0,n-1)$ the set~${\mathcal{D}}_{n,\tau}$ is empty.
\item[{(ii).}] When~$\tau\in(n-1,+\infty)$ the set~$\R^n\setminus {\mathcal{D}}_{n,\tau}$ is of null Lebesgue measure.
\end{itemize}

Also, if~$n\ge2$, prove there exists~$C>0$ such that for all~$\omega\in\R^n$ there exist infinitely many~$k\in\Z^n\setminus\{0\}$ such that
\begin{equation}\label{DIpFA:0ladj01-211}
|\omega\cdot k|\le\frac{C\,|\omega|}{|k|^{n-1}}.
\end{equation}
\end{exercisesk}

\begin{exercisesk}\label{DIpFA:1}
Let~$n\in\N\setminus\{0\}$, $r\in\R$ and consider
\begin{equation}\label{ER:jDi01} \sum_{k\in\Z^n} \big(\cos(k\cdot x)\big)^{|k|^r}.\end{equation}

Prove that\footnote{The case~$n=1$ in Exercise~\ref{DIpFA:1} is already
quite nontrivial.}
\begin{itemize}
\item[{(i).}] When~$r\in(-\infty,2n]$ the series in~\eqref{ER:jDi01} diverges for every~$x\in\R^n$.
\item[{(ii).}] When~$r\in(2n,+\infty)$ the series in~\eqref{ER:jDi01}
converges absolutely for almost every~$x\in\R^n$ but diverges on a dense set in~$\R^n$.
\end{itemize}
\end{exercisesk}

\section{Weyl's Equidistribution Theorem}\label{PelapeammilonS}

A beautiful result, known under the name of \index{Weyl's Equidistribution Theorem}
Weyl's Equidistribution Theorem, states that irrational rotations are more or less uniformly scattered over the circle:

\begin{theorem}\label{PelapeammilonT} Let~$\gamma\in\R\setminus\Q$.

For all~$K\in\N$ and~$a$, $b\in\R$ with~$a\le b$, let~${\mathcal{N}}_{[a,b]}(K)$ be the number of elements of the form~$\gamma k$ with~$k\in\N\cap[0,K-1]$ such that
the fractional part of~$\gamma k$ lies in~$[a,b]$.

Then, for all~$a$, $b\in[0,1]$,
$$\lim_{K\to+\infty}\frac{{\mathcal{N}}_{[a,b]}(K)}{K}=b-a.$$
\end{theorem}

To prove this result, we start with an estimate from above:

\begin{lemma}\label{PelapeammilonTwLEm} Let~$a$, $b\in\R$ with~$a\le b$.
In the notation of Theorem~\ref{PelapeammilonT}, we have that
\begin{equation}\label{PelapeammilonTw2.2} \limsup_{K\to+\infty}\frac{{\mathcal{N}}_{[a,b]}(K)}{K}\le b-a.\end{equation}
\end{lemma}

\begin{proof} Without loss of generality, we can suppose that
\begin{equation}\label{PelapeammilonTw}
a,b\in[0,1].\end{equation}
Indeed, setting~$a':=\max\{a,0\}$ and~$b':=\min\{b,1\}$ we see that~$[a,b]\cap[0,1]=[a',b']$,
therefore~${\mathcal{N}}_{[a,b]}(K)={\mathcal{N}}_{[a',b']}(K)$. Hence, if~\eqref{PelapeammilonTw2.2} holds true for~$a'$ and~$b'$, we have that
$$  \limsup_{K\to+\infty}\frac{{\mathcal{N}}_{[a,b]}(K)}{K}=\limsup_{K\to+\infty}\frac{{\mathcal{N}}_{[a',b']}(K)}{K}\le b'-a'\le b-a.$$
Thus, we assume now that~\eqref{PelapeammilonTw} holds true and we let~$\delta\in(0,1)\cap\left(0,\frac1{4\gamma}\right)$,
$$ X:=\bigcup_{\ell\in\Z}[\ell+a-2\gamma\delta,\ell+b+2\gamma\delta]\qquad{\mbox{and}}\qquad
Y:=\bigcup_{\ell\in\Z}[\ell,\ell+\delta].$$
We define
$$ \R^2\ni(x,y)\longmapsto f(x,y):=\begin{dcases} 1 &{\mbox{ if $x\in X$ and~$y\in Y$,}}\\
0&{\mbox{ otherwise.}}
\end{dcases}$$  

We have thus reduced ourselves to the setting of Section~\ref{LIMOTORSE}: indeed, $f$ is~$\Z^2$-periodic
and
\begin{equation}\label{INTSERGOB}{\mbox{$\omega:=(\gamma,1)$ is nonresonant (in the sense of~\eqref{NORESA}),}} \end{equation}
since if~$(k_1,k_2)\in\Z^2$ is such that~$\gamma k_1+k_2=0$, then necessarily~$k_1=k_2=0$, due to the fact that~$\gamma$ is irrational).
We can thereby utilize Theorem~\ref{ERGOB} and, in view of~\eqref{CqTr64RQL} and~\eqref{CqTr64RQL.lamerpwoekfm}, conclude that
\begin{equation}\label{TCCqTr64RQL.sLO}
\begin{split}&
\lim_{K\to+\infty}\frac1{K}\int_0^K f(\gamma t,t)\,dt=
\langle f\rangle_{\mbox{time}}=\langle f \rangle_{\mbox{phase}}=
\iint_{(0,1)\times(0,1)} f(x,y)\,dx\,dy.\end{split}
\end{equation}

We also observe that 
\begin{equation}\label{CqTr64RQLmnH}\begin{split}&
{\mbox{if~$(x,y)\in(0,1)\times(0,1)$ and~$f(x,y)\ne0$, then necessarily}}\\&{\mbox{$x\in 
[0,2\gamma\delta]\cup[1-2\gamma\delta,1]\cup[a-2\gamma\delta,b+2\gamma\delta]
$ and~$y\in[0,\delta]$.}}\end{split}\end{equation}
To check this, let~$(x,y)$ as above. Then, there exist~$\ell_x$, $\ell_y\in\Z$ such that
$$\ell_x+a-2\gamma\delta\le x<1,\qquad 0<x\le\ell_x+b+2\gamma\delta,\qquad
\ell_y\le y<1,\qquad{\mbox{and}}\qquad 0<y\le\ell_y+\delta .$$
In particular, we have that~$1>\ell_y>-\delta>-1$, which says that~$\ell_y=0$ and thus~$y\in[0,\delta]$.
Moreover, by~\eqref{PelapeammilonTw},
\begin{eqnarray*}&& -2<-1-2\gamma\delta\le-b-2\gamma\delta\le\ell_x\\{\mbox{and }}&&\ell_x\le1-a+2\gamma\delta\le1+2\gamma\delta<2,\end{eqnarray*}
yielding that~$\ell_x\in\{-1,0,1\}$.

Correspondingly, when~$\ell_x=0$, we have that~$x\in[a-2\gamma\delta,b+2\gamma\delta]$.
When instead~$\ell_x=-1$ it follows that~$x<-1+b+2\gamma\delta\le2\gamma\delta$,
while when~$\ell_x=1$ we find that~$x\ge1+a-2\gamma\delta\ge1-2\gamma\delta$.
The proof of~\eqref{CqTr64RQLmnH} is thereby complete.

As a consequence of~\eqref{CqTr64RQLmnH}, we obtain that
\begin{equation}\label{TCCqTr64RQL.sLO2}\iint_{(0,1)\times(0,1)} f(x,y)\,dx\,dy
\le\delta\,(b-a+8\gamma\delta).\end{equation}

Moreover, for each~$k\in\N\cap[0,K-1]$, substituting for~$\tau:=t-k$ we deduce from the~$\Z^2$-periodicity of~$f$ that
\begin{equation*}\begin{split}&
\int_k^{k+1} f(\gamma t,t)\,dt=\int_k^{k+1} f(\gamma t,t-k)\,dt=\int_0^{1} f(\gamma \tau+\gamma k,\tau)\,d\tau=\int_0^{\delta} f(\gamma \tau+\gamma k,\tau)\,d\tau\end{split}
\end{equation*}
and consequently
\begin{equation}\label{BAL:PsBMoanr-daq-2.i2jd.qkwndc}\int_0^K f(\gamma t,t)\,dt=\sum_{k=0}^{K-1}
\int_k^{k+1} f(\gamma t,t)\,dt=\sum_{k=0}^{K-1}
\int_0^{\delta} f(\gamma \tau+\gamma k,\tau)\,d\tau.\end{equation}

Now we claim that, for all~$\tau\in[0,\delta]$,
\begin{equation}\label{BAL:PsBMoanr-daq-2.i2jd}
f(\gamma \tau+\gamma k,\tau)\ge\begin{dcases}
1&{\mbox{ if $\{\gamma k\}\in[a,b]$,}}\\
0&{\mbox{ otherwise,}}\end{dcases}
\end{equation}where we have used the fractional part notation~$\{\cdot\}$, as discussed on page~\pageref{FRaPARTADE}.

To check this, suppose that~$\{\gamma k\}\in[a,b]$. Then, there exists~$\ell\in\Z$ such that~$\gamma k\in[\ell+a,\ell+b]$
and therefore, for all~$\tau\in[0,\delta]$, we see that~$\gamma \tau+\gamma k\in[\ell+a-2\gamma\delta,\ell+b+2\gamma\delta]\subseteq X$, from which we gather that~$f(\gamma \tau+\gamma k,\tau)=1$.

Having completed the proof of~\eqref{BAL:PsBMoanr-daq-2.i2jd}, we deduce from this information that
\begin{eqnarray*}
\int_0^{\delta} f(\gamma \tau+\gamma k,\tau)\,d\tau\ge\begin{dcases}
\delta&{\mbox{ if $\{\gamma k\}\in[a,b]$,}}\\
0&{\mbox{ otherwise}}\end{dcases}
\end{eqnarray*}
and thus, by virtue of~\eqref{BAL:PsBMoanr-daq-2.i2jd.qkwndc}, that
$$ \int_0^K f(\gamma t,t)\,dt\ge\delta{\mathcal{N}}_{[a,b]}(K).$$

From this, \eqref{TCCqTr64RQL.sLO}, and~\eqref{TCCqTr64RQL.sLO2} we arrive at
$$ \limsup_{K\to+\infty}\frac{\delta{\mathcal{N}}_{[a,b]}(K)}{K}\le
\limsup_{K\to+\infty}\frac1{K}\int_0^K f(\gamma t,t)\,dt=
\iint_{(0,1)\times(0,1)} f(x,y)\,dx\,dy\le \delta\,(b-a+8\gamma\delta).$$
We now divide by~$\delta$ and send~$\delta\searrow0$ to obtain the desired result.
\end{proof}

With this preliminary work, we can now complete the proof of Theorem~\ref{PelapeammilonT}.

\begin{proof}[Proof of Theorem~\ref{PelapeammilonT}] 
We will exploit Lemma~\ref{PelapeammilonTwLEm} and slightly modify its proof.
Let~$\epsilon>0$, to be taken as small as we wish in what follows, and
$$ \R^2\ni(x,y)\longmapsto f(x,y):=\begin{dcases} 1 &{\mbox{ if $\{x\}\in[a,b]$ and~$\{y\}\in[0,\delta]$,}}\\
0&{\mbox{ otherwise,}}
\end{dcases}$$ where we have used the fractional part notation~$\{\cdot\}$, as discussed on page~\pageref{FRaPARTADE}.

As in~\eqref{INTSERGOB}, we have thus reduced ourselves to the setting of Section~\ref{LIMOTORSE}: indeed, $f$ is~$\Z^2$-periodic and we can thereby utilise Theorem~\ref{ERGOB}. Thus, by~\eqref{CqTr64RQL} and~\eqref{CqTr64RQL.lamerpwoekfm}, we conclude that
\begin{equation} \label{ma12253RGOB}\begin{split}&
\lim_{K\to+\infty}\frac1{K}\int_0^K f(\gamma t,t)\,dt=
\langle f\rangle_{\mbox{time}}=\langle f \rangle_{\mbox{phase}}\\&\qquad=
\iint_{(0,1)\times(0,1)} f(x,y)\,dx\,dy=\delta\,(b-a).\end{split}
\end{equation}

Moreover (recall~\eqref{BAL:PsBMoanr-daq-2.i2jd.qkwndc}) we have that
\begin{equation}\label{BAL:PsBMoanr-daq-0}\int_0^K f(\gamma t,t)\,dt=\sum_{k=0}^{K-1}
\int_0^{\delta} f(\gamma \tau+\gamma k,\tau)\,d\tau.\end{equation}

Now we observe that, for all~$\tau\in[0,\delta]$,
\begin{equation}\label{BAL:PsBMoanr-daq-1}
f(\gamma \tau+\gamma k,\tau)\le\begin{dcases}
1&{\mbox{ if $\{\gamma k\}\in[0,2\gamma\delta]\cup[1- 2\gamma\delta,1]\cup[a-2\gamma\delta,b+2\gamma\delta]$,}}\\
0&{\mbox{ otherwise.}}
\end{dcases}
\end{equation}
To check this, suppose that~$\{\gamma k\}\not\in[0,2\gamma\delta]\cup[1- 2\gamma\delta,1]\cup[a-2\gamma\delta,b+2\gamma\delta]$. Hence, to prove~\eqref{BAL:PsBMoanr-daq-1}, it suffices to check that~$\{\gamma \tau+\gamma k\}\not\in[a,b]$. 

Suppose, for a contradiction, that~$\{\gamma \tau+\gamma k\}\in[a,b]$. Then, if~$\{\gamma k\}\in(\gamma\delta,1-\gamma\delta)$, it follows that~$\gamma k=\nu+\{\gamma k\}$ for some~$\nu\in\Z$ and~$\gamma k
+\gamma\tau=\nu+\{\gamma k\}+\gamma\tau$, with~$\{\gamma k\}+\gamma\tau\in(0,1)$.
This yields that~$[a,b]\ni \{\gamma \tau+\gamma k\}=
\{\gamma k\}+\gamma\tau$ and therefore~$\{\gamma k\}\in [a-\gamma\delta,b+\gamma\delta]$, which is absurd. 

For this reason, we have that necessarily either~$\{\gamma k\}\in[0,\gamma\delta]$ or~$\{\gamma k\}\in[1-\gamma\delta,1]$ (or both), but this is also in contradiction with our assumptions and the proof of~\eqref{BAL:PsBMoanr-daq-1} is complete.

Thus, in light of~\eqref{BAL:PsBMoanr-daq-0} and~\eqref{BAL:PsBMoanr-daq-1},
\begin{eqnarray*}\\&&\frac1\delta
\int_0^K f(\gamma t,t)\,dt\\&&\qquad\le
{\mathcal{N}}_{[0,2\gamma\delta]}(K)+{\mathcal{N}}_{[1-2\gamma\delta,1]}(K)+
{\mathcal{N}}_{[a-2\gamma\delta,b+2\gamma\delta]}(K)
\\&&\qquad\le{\mathcal{N}}_{[0,2\gamma\delta]}(K)+{\mathcal{N}}_{[1-2\gamma\delta,1]}(K)+
{\mathcal{N}}_{[a-2\gamma\delta,a]}(K)+
{\mathcal{N}}_{[b,b+2\gamma\delta]}(K)+{\mathcal{N}}_{[a,b]}(K)
\\&&\qquad={\mathcal{N}}_{[a,b]}(K)
+\sum_{p\in{\mathcal{S}}}{\mathcal{N}}_{[p,p+2\gamma\delta]}(K)
,\end{eqnarray*}
where~${\mathcal{S}}:=\{0,\,a-2\gamma\delta,\,b-2\gamma\delta,\,1-2\gamma\delta\}$.

Therefore, we deduce from~\eqref{ma12253RGOB} that
\begin{equation*}\begin{split}&b-a=\frac1\delta
\liminf_{K\to+\infty}\frac1{K}\int_0^K f(\gamma t,t)\,dt\le\liminf_{K\to+\infty}\frac1{K}\left(
{\mathcal{N}}_{[a,b]}(K)
+\sum_{p\in{\mathcal{S}}}{\mathcal{N}}_{[p,p+2\gamma\delta]}(K)\right)
.\end{split}\end{equation*}
Since, by Lemma~\ref{PelapeammilonTwLEm},
$$ \limsup_{K\to+\infty}\frac{{\mathcal{N}}_{[p,p+2\gamma\delta]}(K)}{K}\le2\gamma\delta,$$
we gather that
\begin{equation*}\begin{split}&b-a\le\liminf_{K\to+\infty}\frac{{\mathcal{N}}_{[a,b]}(K)}{K}+8\gamma\delta
.\end{split}\end{equation*}
We now send~$\delta\searrow0$ and we obtain that
$$ b-a\le\liminf_{K\to+\infty}\frac{{\mathcal{N}}_{[a,b]}(K)}{K}.$$
The desired result now follows by combining this bound with that obtained in Lemma~\ref{PelapeammilonTwLEm}.
\end{proof}

\begin{exercise}\label{Pelapeammilon-14.1} As a small variant of Theorem~\ref{PelapeammilonT}, let~$I$ be
any interval of the form~$[a,b]$, $[a,b)$, $(a,b]$ or~$(a,b)$ and
let~${\mathcal{N}}_{I}(K)$ be the number of elements of the form~$\gamma k$ with~$k\in\N\cap[0,K-1]$ such that
the fractional part of~$\gamma k$ lies in~$I$, and prove that, for all~$a$, $b\in[0,1]$,
$$\lim_{K\to+\infty}\frac{{\mathcal{N}}_{I}(K)}{K}=b-a.$$
\end{exercise}

\begin{exercise}\label{Pelapeammilon4} Give another proof of Theorem~\ref{PelapeammilonT} by using Exercise~\ref{PelapeammilonDF}.\end{exercise}

\begin{exercise}\label{Pelapeammilon5} With regard to Exercise~\ref{Pelapeammilon4},
would it be possible to give another proof of the claim in Exercise~\ref{PelapeammilonDF}
by assuming Theorem~\ref{PelapeammilonT} (thus showing that
Exercise~\ref{PelapeammilonDF} and Theorem~\ref{PelapeammilonT} are equivalent)?\end{exercise}

\section{Minkowski's Theorem on convex sets}\label{SE:MINK}

Convexity is a beautiful, and very useful, mathematical property, finding applications in geometry,
optimisation, functional analysis, inequalities, machine learning, etc., see e.g.~\cite{MR1451876, MR3497790, MR4211776, MR4501642}.

The following result was proved by Hermann Minkowski in 1889 \index{Minkowski's Theorem on convex sets} and had a great impact in the application of geometric methods to study number theory.

To state it, we recall that a set~$E$ is said to be symmetric about the origin when a point~$x$ belongs to~$E$ if and only if~$-x$ belongs to~$E$.

If~$E\ne\varnothing$ is convex and symmetric about the origin, then~$E$ necessarily contains the origin
(because it contains a point~$x$, whence~$-x$, and therefore the segment joining these two points).

Minkowski's Theorem on convex sets states that, if the volume of the set is large enough, than the set must contain also another point of the lattice~$\Z^n$. More precisely:

\begin{theorem} \label{MINK}
Let~$E\subseteq\R^n$ be a convex set. Assume that~$E$ is symmetric about the origin and that the Lebesgue measure of~$E$ is strictly larger than~$2^n$.

Then, $E$ must contain a point of~$\Z^n$ different from the origin.
\end{theorem}

\begin{proof} Without loss of generality, we can assume that~$E$ is also bounded (see Exercise~\ref{MINKE2}).

Hence, we define
$$ f(x):=\sum_{m\in\Z^n}\chi_E(2(x+m))$$
and we have (see Exercise~\ref{PRECON56789}) that the series above consists of a finite sum and, for every~$k\in\Z^n$,
$$\widehat f_k=\frac1{2^n}\int_{E}e^{-\pi ik\cdot x}\,dx.$$

Moreover (see Exercise~\ref{PRECON56789-BIS}) we have that
$$ \|f\|_{L^2((0,1)^n)}^2={\frac1{2^n}\sum_{j\in\Z^n}\int_{\R^n}\chi_E(x)\,\chi_E(x+2j)\,dx}.$$

From these observations and Parseval's Identity~\eqref{PARSENEN} it follows that
\begin{equation}\label{082ei0218r903v652weder6fye5r6fbv82JEFj}\begin{split}&
\sum_{k\in\Z^n}\left|\frac1{2^n}\int_{E}e^{-\pi ik\cdot x}\,dx\right|^2
=\sum_{k\in\Z^n}|\widehat f_k|^2\\&\qquad=
\|f\|_{L^2((0,1)^n)}^2={\frac1{2^n}\sum_{j\in\Z^n}\int_{\R^n}\chi_E(x)\,\chi_E(x+2j)\,dx}.\end{split}\end{equation}
Accordingly, choosing~$k=0$ in the first sum of~\eqref{082ei0218r903v652weder6fye5r6fbv82JEFj},
\begin{equation}\label{082ei0218r903v652weder6fye5r6fbv82JEFjbi}\frac{|E|^2}{2^{2n}}=\left|
\frac1{2^n}\int_{E}\,dx\right|^2\le{\frac1{2^n}\sum_{j\in\Z^n}\int_{\R^n}\chi_E(x)\,\chi_E(x+2j)\,dx},\end{equation}
where~$|E|$ denotes the Lebesgue measure of~$E$.

Now suppose, for a contradiction, that
\begin{equation}\label{vfe3dfbu9gyv652weder6fye5r6fbv82JEFjbi}
{\mbox{among the points of~$\Z^n$, the set~$E$ contains only the origin.}}\end{equation}
In this case, we claim that, for all~$x\in\R^n$ and~$j\in\Z^n\setminus\{0\}$,
\begin{equation}\label{9ij8u7-0-r903v652weder6fye5r6fbv82JEFjbi}
\chi_E(x)\,\chi_E(x+2j)=0.
\end{equation}
Indeed, if~$x\not\in E$, it follows that~$\chi_E(x)=0$ and we are done.

Similarly, if~$x+2j\not\in E$, then~$\chi_E(x+2j)=0$ and we are done.

We can thereby focus on the case in which~$x\in E$ and~$x+2j\in E$. Since~$E$ is symmetric about the origin, we have that~$-x\in E$. Since~$E$ is convex, the midpoint between~$-x$ and~$x+2j$ must also belong to~$E$, thus
$$ E\ni\frac{-x+(x+2j)}2=j,$$
in contradiction with~\eqref{vfe3dfbu9gyv652weder6fye5r6fbv82JEFjbi}. This completes the proof of~\eqref{9ij8u7-0-r903v652weder6fye5r6fbv82JEFjbi}.

In view of~\eqref{082ei0218r903v652weder6fye5r6fbv82JEFjbi} and~\eqref{9ij8u7-0-r903v652weder6fye5r6fbv82JEFjbi}, we deduce that
$$ \frac{|E|^2}{2^{2n}}\le\frac1{2^n}\int_{\R^n}\chi_E(x)\,\chi_E(x+0)\,dx=\frac1{2^n}\int_{E}\,dx=\frac{|E|}{2^{n}}$$
and consequently~$|E|\le2^n$, against our assumptions.
\end{proof}

\begin{exercise}\label{MINKE}
Prove that the volume assumption in Theorem~\ref{MINK} is optimal by providing an example of convex set which is
symmetric about the origin, has Lebesgue measure equal to~$2^n$ but does not contain any point of~$\Z^n$ other than the origin.
\end{exercise}

\begin{exercise}\label{MINKE2}
Let~$E\subseteq\R^n$ be convex, symmetric about the origin and with Lebesgue measure larger than~$2^n$.
Prove that there exists~$R>0$ large enough such that~$E\cap B_R$ is also
convex, symmetric about the origin and with Lebesgue measure larger than~$2^n$.
\end{exercise}

\section{The Central Limit Theorem}\label{CLRTNSSECMAS}

The Central Limit Theorem states, roughly speaking, that the distribution of the sample mean is, under suitable assumptions, approximately normal and, as such, is one of the pillars of probability and statistics.

In our case, without diving into probabilistic technicalities, we consider functions (actually, \index{random variable} \emph{random variables} in the probabilistic lingo) from some set~$\Omega$ to\footnote{Usually, random variables are taken with values in the reals. Here, considering random variables with values in the circle is just a little trick to cast the problem into the framework of Fourier Series (the case of random variables with values in the reals would require Fourier Transform instead).
\label{ESILVRAs}

See e.g.~\cite{MR1828667, MR1836122, MR3752655} for probabilistic and statistical theories on the circle.}
the circle, identified as the quotient~$\R/\Z$, with representatives in~$[0,1)$.

The set~$\Omega$ is known as the \emph{sample space} \index{sample space} in the probabilistic terminology
and there is a measure~${\mathbb{P}}$ defined on (a suitable $\sigma$-algebra of) subsets\footnote{This $\sigma$-algebra of measurable sets according to the measure~${\mathbb{P}}$ is called \index{event space} \emph{event space} by probabilists, and the corresponding subsets of~$\Omega$ are called \index{event} \emph{events}.} of~$\Omega$, with~${\mathbb{P}}(\Omega)=1$. This measure is called \index{probability function} probability function.

That is, roughly speaking, with this formalism, one is able to speak about the probability of some event.
For our purposes, what we need is that given a random variable~$X:\Omega\to[0,1)$
one can speak about the probability that~$X$ takes value in some subset~$A$ of the circle, that is
$$ {\mathbb{P}}(X\in A):={\mathbb{P}}\Big(\big\{\omega\in\Omega {\mbox{ s.t. }}X(\omega)\in A\big\}\Big).$$

In this notation, two random variables~$X_1$ and~$X_2$ are defined\footnote{Again, in case of
random variables taking values in the reals, the set~$A$ would belong to the reals, rather than the circle (recall footnote~\ref{ESILVRAs}).
A similar modification would be required if one considers random variables taking values in a different space.} to be \emph{identically distributed} \index{random variables, identically distributed} if, for each set~$A$ in the circle,
\begin{equation}\label{DSTAFIPMAkcHienGHaBe} {\mathbb{P}}(X_1\in A)= {\mathbb{P}}(X_2\in A).\end{equation}

Also, $X_1$ and~$X_2$ are defined to be \emph{independent} \index{random variables, independent} 
if, for all sets~$A_1$ and~$A_2$ in the circle,
\begin{equation}\label{DSTAFIPMAkcHi} {\mathbb{P}}(X_1\in A_1,\,X_2\in A_2)= {\mathbb{P}}(X_1\in A_1)\, {\mathbb{P}}(X_2\in A_2),\end{equation}
where
$$ {\mathbb{P}}(X_1\in A_1,\,X_2\in A_2):={\mathbb{P}}\Big(\big\{\omega\in\Omega {\mbox{ s.t. }}X_1(\omega)\in A_1
{\mbox{ and }}X_2(\omega)\in A_2\big\}\Big).$$
Relatedly, given a sequence of random variables~$\{X_j\}_{j\in\N}$ one says that they
are mutually independent if, for any~$N\in\N$ and any sets~$A_0,\dots,A_N$ in the circle,
$$ {\mathbb{P}}(X_0\in A_0,\dots,\,X_N\in A_N)= {\mathbb{P}}(X_0\in A_0)\dots {\mathbb{P}}(X_N\in A_N).$$

We also consider the \emph{wrapped Gaussian distribution} \index{wrapped Gaussian distribution}
\begin{equation}\label{WRAPGA} {\mathcal{W}}(x):=\frac1{\sqrt{2\pi}}\sum_{k\in\Z} e^{-\frac{(x+k)^2}2}.\end{equation}

In this framework, one possible version of the \index{Central Limit Theorem}
Central Limit Theorem goes\footnote{Concerning assumptions~\eqref{SCCSICLA1qGHondcv-1} and~\eqref{SCCSICLA1qGHondcv-2}, we remark that, being the sequence of
random variables identically distributed, we have that, for all~$j\in\N$,
all the~$X_j$'s have the same expected value, that is
$$ \int_\Omega X_j(\omega)\,d{\mathbb{P}}_\omega=\int_\Omega X_0(\omega)\,d{\mathbb{P}}_\omega,$$
see Exercise~\ref{MNVCEC-aIN:2}.

On this basis, assumption~\eqref{SCCSICLA1qGHondcv-1} is simply stating that this common expected value is actually zero.

As a matter of fact, if the common expected value of the random variables~$X_j$ is~$\mu$
and~$\int_\Omega\big( X_j(\omega)\big)^2\,d{\mathbb{P}}_\omega=\varsigma^2$, one can always reduce to the setting in~\eqref{SCCSICLA1qGHondcv-1} and~\eqref{SCCSICLA1qGHondcv-2} by replacing~$X_j$ with~$\frac{X_j-\mu}{\sqrt{\varsigma^2-\mu^2}}$.} as follows:

\begin{theorem}\label{SCCSICLA1q} Let~$\{X_j\}_{j\in\N}$ be a sequence of random variables.
Suppose that they are all identically distributed and mutually independent.

Assume also that, for all~$j\in\N$,
\begin{equation}\label{SCCSICLA1qGHondcv-1}
\int_\Omega X_j(\omega)\,d{\mathbb{P}}_\omega=0
\end{equation}
and
\begin{equation}\label{SCCSICLA1qGHondcv-2}
\int_\Omega\big( X_j(\omega)\big)^2\,d{\mathbb{P}}_\omega=1.
\end{equation}

Then, if
$$ Z_N:=\frac1{\sqrt{N}}\sum_{j=0}^{N-1}X_j,$$
we have that, for every interval~$J$ in the circle,
\begin{equation}\label{iksd67.lo0pIIJlimNtfty}\lim_{N\to+\infty} {\mathbb{P}}(Z_N\in J)=
\int_J {\mathcal{W}}(x)\,dx.
\end{equation}
\end{theorem}

A salient feature of this kind of results is that no assumption is taken on~${\mathbb{P}}(X_k\in J)$,
namely whatever the form of the distribution of the original random variables~$X_k$,
the limit distribution of their normalised partial
sums is the same, as stated in~\eqref{iksd67.lo0pIIJlimNtfty}, highlighting
a noteworthy ``universal'' characteristic of the wrapped Gaussian distribution:
significantly, in a nutshell, this gives that any a sequence of
independent ``trials'' approaches a wrapped Gaussian distribution distribution as the number of
observations increases.

See~\cite{MR2569085, MR2743162} for more information about the importance and
the fascinating history of the Central Limit Theorem.
See Exercise~\ref{02ktml-PMSdPR-2lXC135.30568emahbP} for the construction
of a sequence of independent identically distributed random variables.
See also~\cite[page~41]{MR1852999}
for the existence of probability spaces that admit a sequence of independent identically distributed random variables with specified distributions.

\begin{proof}[Proof of Theorem~\ref{SCCSICLA1q}] Given a random variable~$X$, for all~$\xi\in\R$ and~$\omega\in\Omega$ we have that
$$|e^{2\pi i\xi X(\omega)}|=1$$ and therefore~$e^{2\pi i\xi X(\cdot)}\in L^1(\Omega,{\mathbb{P}})$. 

Hence, for every~$\xi\in\R$, we can define\footnote{The function in~\eqref{CHARAFA:VA} is often called the
\emph{characteristic function} of the random variable~$X$. Not to be confused with the characteristic function of a set,
as in~\eqref{CHARAFA}.}
\begin{equation}\label{CHARAFA:VA} {\mathcal{C}}_X(\xi):=\int_\Omega e^{2\pi i\xi X(\omega)}\,d{\mathbb{P}}_\omega.\end{equation}

We see that
\begin{equation*}
\begin{split}
{\mathcal{C}}_{Z_N}(\xi)&=\int_\Omega e^{2\pi i\,\frac\xi{\sqrt{N}}\big( X_0(\omega)+\dots+X_{N-1}(\omega)\big)}\,d{\mathbb{P}}_\omega\\&=\int_\Omega e^{2\pi i\,\frac\xi{\sqrt{N}}\,X_0(\omega)}\,
\dots\, e^{2\pi i\,\frac\xi{\sqrt{N}}\,X_{N-1}(\omega)}\,d{\mathbb{P}}_\omega.
\end{split}
\end{equation*}

Hence, using Exercise~\ref{MNVCEC-aIN},
\begin{equation*}
{\mathcal{C}}_{Z_N}(\xi)=\int_\Omega e^{2\pi i\,\frac\xi{\sqrt{N}}\,X_0(\omega)}\,d{\mathbb{P}}_\omega\,
\dots\, \int_\Omega e^{2\pi i\,\frac\xi{\sqrt{N}}\,X_{N-1}(\omega)}\,d{\mathbb{P}}_\omega.
\end{equation*}
Therefore, on the wake of Exercise~\ref{MNVCEC-aIN:2},
\begin{equation}\label{SCCSICLA1qGHondcv-3}
{\mathcal{C}}_{Z_N}(\xi)=\left(\int_\Omega e^{2\pi i\,\frac\xi{\sqrt{N}}\,X_0(\omega)}\,d{\mathbb{P}}_\omega\right)^N.
\end{equation}

Now, for a given~$\xi\in\R$, we perform a Taylor expansion of the complex exponential and we see that, as~$N\to+\infty$,
\begin{equation} \label{aojsn028rRTutgh9irghiougertq3gedfgbn34r58tuhjhJ34.n} e^{2\pi i\,\frac\xi{\sqrt{N}}\,X_0(\omega)}=1+\frac{2\pi i\,\xi\,X_0(\omega)}{\sqrt{N}}-\frac{2\pi^2\xi^2\,\big(X_0(\omega)\big)^2}{N}+\frac{\Upsilon_N(\omega,\xi)}N,\end{equation}
where
$$ \Upsilon_N(\omega,\xi):=
2\pi^2\xi^2\,\big(X_0(\omega)\big)^2\,\big(1-e^{i\upsilon_N(\omega,\xi)}\big),$$
with$$|\upsilon_N(\omega,\xi)|\le\frac{2\pi|\xi|\,|X_0(\omega)|}{\sqrt{N}}.$$

In particular,
$$ \lim_{N\to+\infty}\sup_{\xi\in J}|\upsilon_N(\omega,\xi)|=0$$
and, as a result, if
$$ \Upsilon_N(\omega):=\sup_{\xi\in J}| \Upsilon_N(\omega,\xi)|,$$
we have that
\begin{equation} \label{aojsn028rRTutgh9irghiougertq3gedfgbn34r58tuhjhJ34}
\lim_{N\to+\infty}\Upsilon_N(\omega)\le
2\pi^2\,\big(X_0(\omega)\big)^2\,\lim_{N\to+\infty}\sup_{\xi\in J}\left(\xi^2
\left|1-\exp\left(i\upsilon_N(\omega,\xi)\right)\right|\right)=0.
\end{equation}
We also stress that, for some~$C>0$,
$$|\Upsilon_N(\omega)|\le\sup_{\xi\in J}
4\pi^2\xi^2\,\big(X_0(\omega)\big)^2\le C\,\big(X_0(\omega)\big)^2,$$ whose integral over~$\Omega$ with respect to the measure~${\mathbb{P}}$ is finite, thanks to~\eqref{SCCSICLA1qGHondcv-2}.

This, \eqref{aojsn028rRTutgh9irghiougertq3gedfgbn34r58tuhjhJ34} and the Dominated Convergence Theorem (see e.g.~\cite[Theorem~10.31]{MR3381284}) yield that
$$ \lim_{N\to+\infty}\int_\Omega\Upsilon_N(\omega)\,d{\mathbb{P}}_\omega=0.$$

Therefore, in light of~\eqref{SCCSICLA1qGHondcv-1}, \eqref{SCCSICLA1qGHondcv-2}, and~\eqref{aojsn028rRTutgh9irghiougertq3gedfgbn34r58tuhjhJ34.n},
$$ \int_\Omega e^{2\pi i\,\frac\xi{\sqrt{N}}\,X_0(\omega)}\,d{\mathbb{P}}_\omega=
1-\frac{2\pi^2\xi^2}{N}+o\left(\frac1N\right).$$

Plugging this information into~\eqref{SCCSICLA1qGHondcv-3} we conclude that
\begin{equation}\label{CHARAFA:VA-5}
\lim_{N\to+\infty}{\mathcal{C}}_{Z_N}(\xi)=\lim_{N\to+\infty}\left(1-\frac{2\pi^2\xi^2}{N}+o\left(\frac1N\right)\right)^N=e^{-2\pi^2\xi^2}.
\end{equation}

Now we consider a smooth function~$f$, periodic of period~$1$, and, owing to Theorems~\ref{SMXC22} and~\ref{BASw},
we have that
$$f(x)=\sum_{k\in\Z}\widehat f_k\,e^{2\pi ikx},$$
with uniform convergence of the series and
\begin{equation}\label{CHARAFA:VA-56}
\sum_{k\in\Z}|\widehat f_k|<+\infty.
\end{equation}

Hence, by~\eqref{CHARAFA:VA},
\begin{equation}\label{CHARAFA:VA-99}\begin{split}&
\sum_{k\in\Z}\widehat f_k\,{\mathcal{C}}_{Z_N}(k)=\sum_{k\in\Z}\int_\Omega \widehat f_k\,e^{2\pi ik Z_N(\omega)}\,d{\mathbb{P}}_\omega\\&\qquad=\int_\Omega \sum_{k\in\Z}\widehat f_k\,e^{2\pi ik Z_N(\omega)}\,d{\mathbb{P}}_\omega=\int_\Omega f(Z_N(\omega))\,d{\mathbb{P}}_\omega.\end{split}
\end{equation}

We also deduce from~\eqref{CHARAFA:VA} that, for every random variable~$X$ and~$\xi\in\R$,
$$ |{\mathcal{C}}_X(\xi)|\le \int_\Omega \big|e^{2\pi i\xi X(\omega)}\big|\,d{\mathbb{P}}_\omega={\mathbb{P}}(\Omega)=1$$
and consequently, in virtue of~\eqref{CHARAFA:VA-5} and~\eqref{CHARAFA:VA-56},
$$ \lim_{N\to+\infty}\sum_{k\in\Z}\widehat f_k\,{\mathcal{C}}_{Z_N}(k)
=\sum_{k\in\Z}\widehat f_k\,\lim_{N\to+\infty}{\mathcal{C}}_{Z_N}(k)=\sum_{k\in\Z}\widehat f_k\,e^{-2\pi^2 k^2}.
$$

This and~\eqref{CHARAFA:VA-99} yield that, for every smooth function~$f$, periodic of period~$1$,
\begin{equation}\label{iksd67.lo0pIIJlimNtfty.09} \lim_{N\to+\infty}
\int_\Omega f(Z_N(\omega))\,d{\mathbb{P}}_\omega=
\sum_{k\in\Z}\widehat f_k\,e^{-2\pi^2 k^2}.
\end{equation}

Now, to prove~\eqref{iksd67.lo0pIIJlimNtfty}, we consider an interval~$J$ in the circle.
If~$J$ coincides with the full circle, both sides in~\eqref{iksd67.lo0pIIJlimNtfty} equals~$1$ and we are done,
so we may assume that~$J$ is a proper sub-interval of the circle. Thus,
given~$\epsilon>0$ sufficiently small, we can find intervals~$G_\epsilon$ and~$H_\epsilon$,
obtained respectively by shrinking or enlarging~$J$ by a length~$\epsilon$ on both sides,
and smooth periodic functions~$g_\epsilon$ and~$h_\epsilon$ on the circle,
non-negative and bounded by~$1$, and such that~$g_\epsilon=1$ in~$G_\epsilon$, $g_\epsilon=0$ outside~$J$,
$h_\epsilon=1$ in~$J$, $h_\epsilon=0$ outside~$H_\epsilon$.

In this way,
\begin{equation}\label{iksd67.lo0pIIJlimNtfty.013}
\begin{split}&
|\widehat g_{\epsilon,k}-\widehat\chi_{J,k}|=
\left| \int_J g_\epsilon(x) e^{-2\pi ikx}\,dx
-\int_J e^{-2\pi ikx}\,dx\right|\\&\qquad=
\left| \int_{J\setminus G_\epsilon}\big( g_\epsilon(x)-1\big)\, e^{-2\pi ikx}\,dx\right|\le
|J\setminus G_\epsilon|=2\epsilon
\end{split}\end{equation}
and similarly
\begin{equation}\label{iksd67.lo0pIIJlimNtfty.091}
|\widehat h_{\epsilon,k}-\widehat\chi_{J,k}|\le2\epsilon.
\end{equation}

Also, since~$g_\epsilon\le \chi_J\le h_\epsilon$,
\begin{equation*}
\int_\Omega g_\epsilon(Z_N(\omega))\,d{\mathbb{P}}_\omega\le
\int_\Omega \chi_J(Z_N(\omega))\,d{\mathbb{P}}_\omega\le
\int_\Omega h_\epsilon(Z_N(\omega))\,d{\mathbb{P}}_\omega.
\end{equation*}
From this, \eqref{iksd67.lo0pIIJlimNtfty.09}, and~\eqref{iksd67.lo0pIIJlimNtfty.013} we deduce that
\begin{equation}\label{YOTHCLCHINDMAKSTABS:a}\begin{split}&
\sum_{k\in\Z}\big(\widehat \chi_{J,k}-2\epsilon\big)\,e^{-2\pi^2 k^2}\le
\sum_{k\in\Z}\widehat g_{\epsilon,k}\,e^{-2\pi^2 k^2}=
\lim_{N\to+\infty}\int_\Omega g_\epsilon(Z_N(\omega))\,d{\mathbb{P}}_\omega\\&\qquad\le\lim_{N\to+\infty}
\int_\Omega \chi_J(Z_N(\omega))\,d{\mathbb{P}}_\omega=\lim_{N\to+\infty} {\mathbb{P}}(Z_N\in J)\end{split}
\end{equation}
and likewise, using~\eqref{iksd67.lo0pIIJlimNtfty.091},
\begin{equation*}
\sum_{k\in\Z}\big(\widehat \chi_{J,k}+2\epsilon\big)\,e^{-2\pi^2 k^2}\ge\lim_{N\to+\infty} {\mathbb{P}}(Z_N\in J).
\end{equation*}

Comparing this and~\eqref{YOTHCLCHINDMAKSTABS:a}, taking the limit as~$\epsilon\searrow0$
and using the Dominated Convergence Theorem, we obtain that
\begin{equation*}\begin{split}&
\lim_{N\to+\infty} {\mathbb{P}}(Z_N\in J)=
\sum_{k\in\Z}\widehat \chi_{J,k}\,e^{-2\pi^2 k^2}
=\sum_{k\in\Z}\int_0^1 \chi_{J}(x)\,e^{-2\pi ikx-2\pi^2 k^2}\,dx\\&\qquad=\int_0^1 \chi_{J}(x)\,\left(\sum_{k\in\Z}e^{-2\pi ikx-2\pi^2 k^2}\right)\,dx
.\end{split}\end{equation*}
The proof of~\eqref{iksd67.lo0pIIJlimNtfty} is thereby complete, thanks to Exercise~\ref{MNVCEC-aIN:4}.
\end{proof}

\begin{exercise} \label{MNVCEC-aIN}
Let~$X_1$ and~$X_2$ be independent random variables and~$f$ be a continuous function, periodic of period~$1$.

Prove that
$$ \int_\Omega f\big( X_1(\omega)\big)\,f\big( X_2(\omega)\big)\,d{\mathbb{P}}_\omega=
\int_\Omega f\big( X_1(\omega)\big)\,d{\mathbb{P}}_\omega\,
\int_\Omega f\big( X_2(\omega)\big)\,d{\mathbb{P}}_\omega.$$
\end{exercise}

\begin{exercise} \label{MNVCEC-aIN:2}
Let~$X_1$ and~$X_2$ be identically distributed random variables and~$f$ be a continuous function, periodic of period~$1$.

Prove that
$$ \int_\Omega f\big( X_1(\omega)\big)\,d{\mathbb{P}}_\omega=
\int_\Omega f\big( X_2(\omega)\big)\,d{\mathbb{P}}_\omega.$$
\end{exercise}

\begin{exercise} \label{MNVCEC-aIN:3} Given~$a\in(0,+\infty)$ and~$\zeta\in\R$, prove that
$$ \int_{-\infty}^{+\infty} e^{-ax^2-2\pi i\zeta x}\,dx=\sqrt{\frac{\pi}a}\,e^{-\frac{\pi^2\zeta^2}{a}}
.$$\end{exercise}

\begin{exercise}\label{MNVCEC-aIN:4}
In the notation of~\eqref{WRAPGA}, prove that
$$ {\mathcal{W}}(x)= \sum_{k\in\Z}e^{-2\pi ikx-2\pi^2 k^2}.$$
\end{exercise}

\begin{exercise}\label{02ktml-PMSdPR-2lXC135.30568emahbP}
Let~$\Omega:=[0,1)$. For each~$\omega\in\Omega$, we consider its expansion in base~$2$, say~$\omega=0.\omega_0\omega_1\omega_2\dots$, with~$\omega_j\in\{0,1\}$.
For each~$j\in\N$, let
$$X_j(\omega):=\frac{\omega_j}2.$$
Prove that this is a sequence of identically distributed and mutually independent random variables with respect to the Lebesgue measure on~$[0,1)$, satisfying
$$ \int_\Omega X_j(\omega)\,d\omega=\frac14\qquad{\mbox{ and }}\qquad
\int_\Omega X_j^2(\omega)\,d\omega=\frac18\qquad{\mbox{for all }}j\in\N.$$\end{exercise}

\chapter{What Comes Next?}

The reader who has patiently walked with us thus far may wonder what new territory lies ahead.

Of course, there is no single answer to this question, as the natural progression depends on the reader's preferences, interests, and goals. 

Nevertheless, a standard path to follow after Fourier Series is to study the Fourier Transform (to which we intend to devote a book~\cite{NEXT}) and related methods such as the Laplace Transform (see~\cite{MR842186, MR1716143, MR2001192, MR3185433}).

Fourier methods also find many natural applications in Partial Differential Equations~\cite{MR2597943} and play a key role in functional analysis~\cite{MR1216137, MR1892228}, providing a deep understanding of infinite-dimensional spaces and operators.

These areas naturally lead to the development of suitable functional spaces, such as Sobolev spaces~\cite{MR2597943, MR2895178, MR2944369, MR3726909, MR4567945}
and to the introduction of concepts that go beyond classical functions and differential operators (such as distributions~\cite{MR80707, MR115085, MR209834, MR3469458,
MR3469849, MR3468845},
and pseudo-differential operators~\cite{MR1852334, MR597144, MR597145, MR2884718}).

Significantly, Fourier theory is also extremely influential in complex analysis~\cite{MR924157}, geometry~\cite{MR2132704, MR2244106}, combinatorics~\cite{MR2289012},
number theory~\cite{MR434929, MR2061214, MR3307692}, probability and statistics~\cite{MR2839251}, etc.

The applications of Fourier analysis are vast and include digital signal processing~\cite{MR4789586}, wavelets and time-frequency analysis~\cite{MR1162107, MR4201879},
numerical analysis~\cite{MR4696613}, quantum mechanics~\cite{MR2722363, MR2797644}, and more:
as Hamlet reminds us, there are more things in heaven and earth than are dreamt of in our philosophy.

\chapter{Solutions to selected exercises}\label{SOLUTIONS}

\section{Solutions to selected exercises of Section~\ref{SEC:PER:FUNC:913uo}}

\paragraph{Solution to Exercise~\ref{PRD1}.}
We integrate by parts twice and we see that, for every~$n\ge2$,
\begin{eqnarray*}&&
I_n(\omega)\\&=&\frac1{2\pi\omega}\int_0^{1/4} \cos^n(2\pi t)\,\frac{d}{dt}\big(\sin(2\pi\omega t)\big)\,dt\\
&=&\frac1{2\pi\omega}\Bigg(
\cos^n\left(\frac\pi2\right)\,\sin\left(\frac{\pi\omega}2\right)-\cos^n(0)\,\sin(0)\\&&\qquad+
2\pi n\int_0^{1/4} \cos^{n-1}(2\pi t)\,\sin(2\pi t)\,\sin(2\pi\omega t)\,dt
\Bigg)\\&=&\frac{n}{\omega}\int_0^{1/4} \cos^{n-1}(2\pi t)\,\sin(2\pi t)\,\sin(2\pi\omega t)\,dt
\\&=& -\frac{n}{2\pi\omega^2}\int_0^{1/4} \cos^{n-1}(2\pi t)\,\sin(2\pi t)\,\frac{d}{dt}\left(\cos(2\pi\omega t)\right)\,dt\\&=& -\frac{n}{2\pi\omega^2}\Bigg(
\cos^{n-1}\left(\frac\pi2\right)\,\sin\left(\frac\pi2\right)\,\cos\left(\frac{\pi\omega}2\right)
- \cos^{n-1}(0)\,\sin(0)\,\cos(0)\\&&-
\int_0^{1/4} \Big(2\pi \cos^{n-1}(2\pi t)\,\cos(2\pi t)
-2\pi(n-1)\cos^{n-2}(2\pi t)\,\sin^2(2\pi t)\Big)\,
\cos(2\pi\omega t)\,dt\Bigg)\\&=& \frac{n}{\omega^2}
\int_0^{1/4} \Big( \cos^{n}(2\pi t)
-(n-1)\cos^{n-2}(2\pi t)\,\big(1-\cos^2(2\pi t)\big)\Big)\,
\cos(2\pi\omega t)\,dt\\&=& \frac{n}{\omega^2}
\int_0^{1/4} \Big( n\cos^{n}(2\pi t)
-(n-1)\cos^{n-2}(2\pi t)\Big)\,
\cos(2\pi\omega t)\,dt\\&=&\frac{n^2}{\omega^2}\,I_{n}(\omega)-\frac{n(n-1)}{\omega^2}\,I_{n-2}(\omega),
\end{eqnarray*}
from which one obtains the desired result by moving one term to the other side.

\paragraph{Solution to Exercise~\ref{PRD2}.}
We argue by induction. When~$m=0$, the result is true, since the product in the right-hand side is void (hence equal to~$1$).

We suppose now the result to be true for some~$m\in\N$ and we prove it for~$m+1$. To this end, we employ Exercise~\ref{PRD1} with~$n:=2(m+1)$ and we see that
$$ I_{2(m+1)}(\omega)=\frac{2(m+1)(2m+1)}{(2(m+1))^2-\omega^2}\,I_{2m}(\omega).$$
Then, we use the inductive assumption to compute that
\begin{eqnarray*} &&\frac{I_{2(m+1)}(\omega)}{I_{2(m+1)}(0)}\,\prod_{j=1}^{m+1}\frac{(2j)^2-\omega^2}{(2j)^2}
\\&=&\frac{2(m+1)(2m+1)}{(2(m+1))^2-\omega^2}\cdot\frac{(2(m+1))^2}{2(m+1)(2m+1)}\cdot
\frac{I_{2m}(\omega)}{I_{2m}(0)}\,\prod_{j=1}^{m+1}\frac{(2j)^2-\omega^2}{(2j)^2}\\
&=&\frac{(2(m+1))^2}{(2(m+1))^2-\omega^2}\cdot
\frac{I_0(\omega)}{I_0(0)}\cdot\frac{(2(m+1))^2-\omega^2}{(2(m+1))^2}\\&=&\frac{I_0(\omega)}{I_0(0)},
\end{eqnarray*}as desired.

\paragraph{Solution to Exercise~\ref{PRD3}.}
Since~$I_n(\omega)$ is even in~$\omega$, we can focus on the case~$\omega\in\left[0,\frac1{2\pi}\right]$.
For all~$t\in\left[0,\frac14\right]$ we have that~$\cos(2\pi t)\ge0$.
Additionally, for~$\omega\in\left[0,\frac1{2\pi}\right]$,
we have that~$\cos(2\pi\omega t)\ge0$ and therefore, for all~$n\ge1$,
$$ \cos^n(2\pi t)\,\cos(2\pi\omega t) \le \cos^n(2\pi t).$$
Integrating this relation, we find that
\begin{equation}\label{oqwdnf94}I_n(\omega)\le I_n(0).\end{equation}

Moreover, by the monotonicity of the cosine function in the interval~$\left[0,\frac\pi2\right]$,
$$ \cos(2\pi\omega t)\ge\cos(2\pi t) $$
and therefore $$\cos^n(2\pi t)\,\cos(2\pi\omega t) \ge \cos^{n+1}(2\pi t)\ge\cos^{n+2}(2\pi t).$$
Integrating this relation, we find that
\begin{equation}\label{IENNE}I_n(\omega)\ge I_{n+2}(0).\end{equation}
Moreover, in light of Exercise~\ref{PRD1},
$$ I_{n+2}(0)=\frac{n+1}{n+2}\,I_{n}(0).$$
This and~\eqref{IENNE} yield that
\begin{equation*}I_n(\omega)\ge \frac{n+1}{n+2}\,I_{n}(0).\end{equation*}
The desired result now follows from this inequality and~\eqref{oqwdnf94}.

\paragraph{Solution to Exercise~\ref{PROSI}.}
We stress that both sides of~\eqref{EUJNC-P} are holomorphic functions (by uniform convergence,
see e.g.~\cite[Theorem~10.28]{MR924157}), hence, by the uniqueness of power series (see e.g.~\cite[Theorem~10.18]{MR924157})
it suffices to
check~\eqref{EUJNC-P} in the set~$x\in\left(\frac1{100},\frac1{100}\right)\times\{0\}$.

Hence, we pick~$x\in\left(\frac1{100},\frac1{100}\right)$.
Let~$I_n$ be as in Exercise~\ref{PRD1}. Then, for all~$\omega\in\left[-\frac1{2\pi},\frac1{2\pi}\right]$,
$$ I_0(\omega)=\int_0^{1/4} \cos(2\pi\omega t)\,dt=\frac{\sin\left(\frac{\pi \omega}2\right)}{2 \pi\omega}$$
and
$$ I_0(0)=\frac14.$$
Therefore, the result of Exercise~\ref{PRD2} can be restated as
$$ \frac{2\sin\left(\frac{\pi \omega}2\right)}{\pi\omega}=\frac{I_{2m}(\omega)}{I_{2m}(0)}\,\prod_{j=1}^m\frac{(2j)^2-\omega^2}{(2j)^2}$$
or equivalently, choosing~$x:=\frac{ \omega}2$,
\begin{equation*} \frac{\sin(\pi x)}{\pi x}=\frac{I_{2m}(2 x)}{I_{2m}(0)}\,\prod_{j=1}^m\left(1-\frac{ x^2}{ j^2}\right).\end{equation*}
We now send~$m\to+\infty$ and deduce the desired result from Exercise~\ref{PRD3}.

\begin{figure}[h]
\includegraphics[height=3cm]{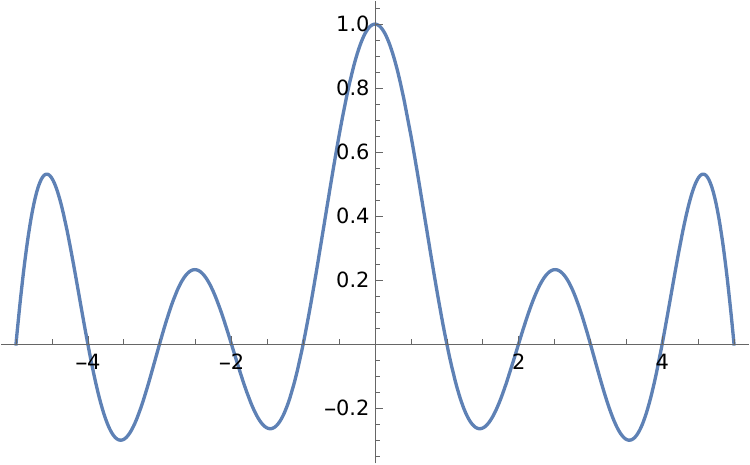}$\,\;$\includegraphics[height=3cm]{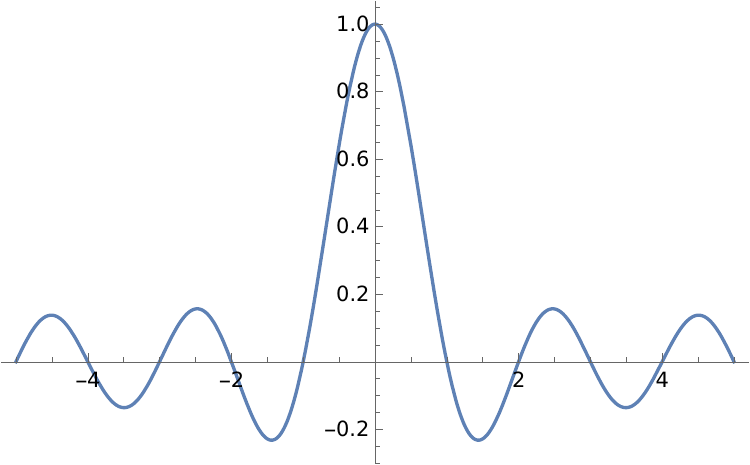}$\,\;$
\includegraphics[height=3cm]{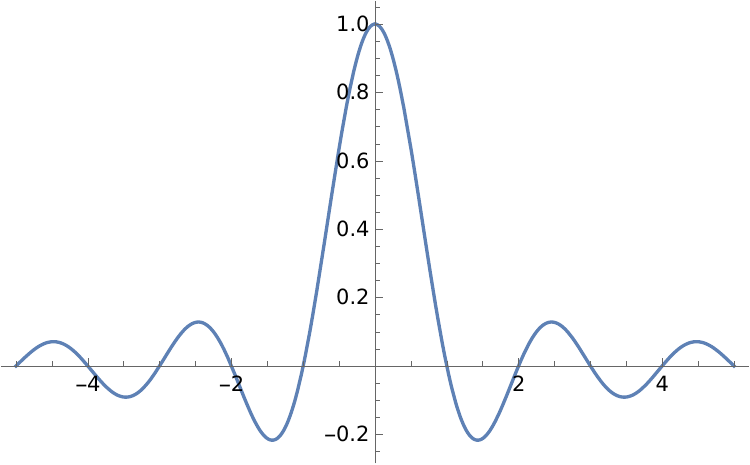}
\centering
\caption{Left: plot of~$\displaystyle\prod_{j=1}^{10}\left(1-\frac{x^2}{ j^2}\right)$.
Centre: plot of~$\displaystyle\prod_{j=1}^{30}\left(1-\frac{x^2}{ j^2}\right)$.
Right: plot of~$\displaystyle\frac{\sin (\pi x)}{\pi x}$.}\label{ikdnsEU}
\end{figure}

See Figure~\ref{ikdnsEU} for a sketch of the approximation procedure given by
Euler's sine product formula. See also~\cite[pages~12--18]{MR1483074}, \cite{MR3383900}
and the references therein for more details on Euler's sine product formula.

\paragraph{Solution to Exercise~\ref{PROSIW}.} We employ Euler's sine product formula
of Exercise~\ref{PROSI} with~$z:=\frac12$ and we see that
$$ \frac2\pi=\frac{2\sin \left(\frac\pi2\right)}{\pi}=\prod_{j=1}^{+\infty}\left(1-\frac{1}{(2 j)^2}\right)=
\prod_{j=1}^{+\infty}\frac{4j^2-1}{4j^2}
.$$ The desired formula now follows by taking the reciprocal of this identity.

\begin{figure}[h]
\includegraphics[height=5.3cm]{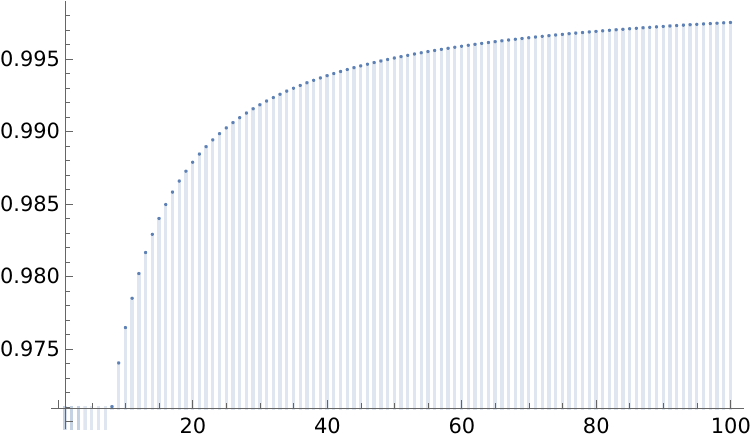}
\centering
\caption{Plot of~$\displaystyle\frac2\pi \prod_{j=1}^n \frac{4 j^2}{4 j^2 - 1}$
for~$n$ up to~$100$ (to be compared also with Figure~\ref{STISamdi}).}\label{ikdns-0oiuygfd98uytf0iu-09EU}
\end{figure}

For additional readings, see e.g.~\cite{MR2519493} to see a link between Wallis's formula
and the integral of the Gaussian.

\paragraph{Solution to Exercise~\ref{STIRL-va}.}
We point out that
\begin{eqnarray*} &&\prod_{j=1}^{n}{\frac{4j^{2}}{4j^{2}-1}}=4^n \prod_{j=1}^{n}{\frac{j^{2}}{(2j+1)(2j-1)}}=
4^{2n} \prod_{j=1}^{n}{\frac{j^{4}}{(2j+1)(2j)(2j-1)(2j)}}\\&&\qquad=
\frac{4^{2n} (n!)^4}{\displaystyle
\prod_{k=1}^{n}(2k+1)(2k)\,\prod_{h=1}^{n}(2h-1)(2h)}=\frac{4^{2n} (n!)^4}{\displaystyle
(2n+1)!\,(2n)!}=\frac{4^{2n} (n!)^4}{\displaystyle
(2n+1)((2n)!)^2}.
\end{eqnarray*}
Hence, by Wallis's formula in Exercise~\ref{PROSIW},
$$\sqrt{\frac{\pi }{2}}=\lim_{n\to+\infty}\sqrt{\prod_{j=1}^{n }{\frac{4j^{2}}{4j^{2}-1}}}
=\lim_{n\to+\infty}
\frac{2^{2n} \,(n!)^2}{\sqrt{2n+1}\,(2n)!}=\frac1{\sqrt2}\lim_{n\to+\infty}
\frac{2^{2n} \,(n!)^2}{\sqrt{n}\;(2n)!},$$
as desired.

\paragraph{Solution to Exercise~\ref{STIRL}.}
Let
$$ \alpha_n:=\frac{n!}{\sqrt{n}}{\left({\frac{e}{n}}\right)}^{n}$$
and observe that
\begin{eqnarray*}
\frac{\alpha_n^4}{\alpha_{2n}^2}=\frac{2^{4n+1}\,(n!)^4}{n\,((2n)!)^2}.
\end{eqnarray*}
Hence, by Exercise~\ref{STIRL-va},
\begin{equation}\label{0qowekd-34rt} \lim_{n\to+\infty}\frac{\alpha_n^4}{\alpha_{2n}^2}=2\pi.\end{equation}

We now aim at showing that~$\alpha_n$ approaches a well-defined limit as~$n\to+\infty$.
To this end, we remark that
\begin{eqnarray*}\beta_n:=
\ln \alpha_n=\ln(n!)-\frac{\ln n}2+n-n\ln n=
\sum_{k=2}^n\ln k-\frac{\ln n}2+n-n\ln n.
\end{eqnarray*}
That is, since
$$ \int_1^n \ln x \,dx = n \ln n - n+1,$$
we see that
\begin{eqnarray*}
\beta_n=
\sum_{k=2}^{n}\ln k-\frac{\ln n}2-\int_1^n \ln x \,dx-1.
\end{eqnarray*}
As a result, for large~$n$,
\begin{eqnarray*}
|\beta_{n+1}-\beta_n|&=&\left|
\frac{\ln (n+1)+\ln n }2-\int_n^{n+1} \ln x \,dx
\right|\\&=&\left|\frac{\ln (n+1)+\ln n }2+n \ln n - (n + 1) \ln(n + 1) + 1\right|\\&=&\left|
1 - \left(n + \frac12\right) \ln\left(1+\frac{1}{n}\right) \right|\\&=&\left|
1 - \left(n + \frac12\right)\left(\frac{1}{n}-\frac{1}{2n^2}+O\left(\frac1{n^3}\right)
\right) \right|\\&=&O\left(\frac1{n^2}\right).
\end{eqnarray*}
Consequently, there exists~$C>0$ such that, for all~$m\in\N$, $$
|\beta_{n+m}-\beta_n|\le\sum_{j=n}^{n+m-1}
|\beta_{j+1}-\beta_{j}|\le C\sum_{j=n}^{n+m-1}\frac1{j^2}.$$
Since the latter is the tail of a convergent series, it follows that~$\beta_n$ is a Cauchy sequence,
whence~$\beta_n$ converges  to some~$\beta\in\R$ as~$n\to+\infty$.

This gives that
$$ \lim_{n\to+\infty}\alpha_n=\lim_{n\to+\infty} e^{\beta_n}=e^\beta=:L\in(0,+\infty).$$

As a result, in light of~\eqref{0qowekd-34rt}, 
$$ L^2=\frac{L^4}{L^2}=\lim_{n\to+\infty}\frac{\alpha_n^4}{\alpha_{2n}^2}=2\pi,$$
from which the desired result follows.

\begin{figure}[h]
\includegraphics[height=5.3cm]{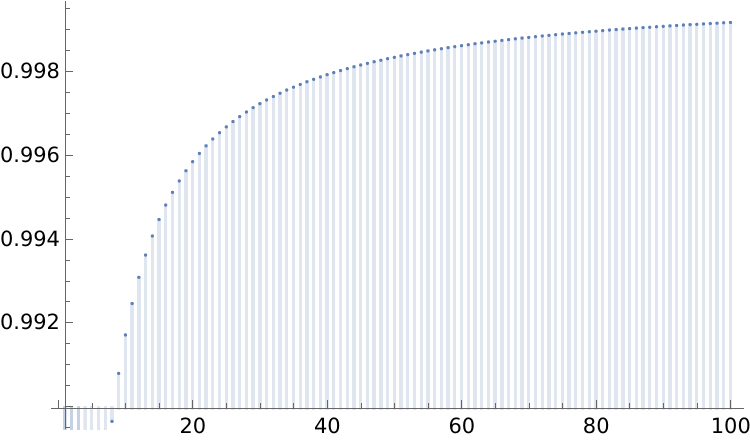}
\centering
\caption{Plot of~$\displaystyle\frac{\sqrt {2\pi n}}{n!}{\left({\frac{n}{e}}\right)}^{n}$ for~$n$ up to~$100$
(to be compared also with Figure~\ref{ikdns-0oiuygfd98uytf0iu-09EU}).}\label{STISamdi}
\end{figure}

See Figure~\ref{STISamdi} for a diagram sketching Stirling's approximation for the factorial.
For further readings on this topic, see e.g.~\cite[Theorem~15.19]{MR248290}.

\paragraph{Solution to Exercise~\ref{SnsdFrhndiYUJS-2}.} Let~$P$ be a non-constant polynomial with complex coefficients. Suppose by contradiction that~$P(z)\ne0$ for each~$z\in\C$. Then, for all~$r\in[0,+\infty)$ and~$x\in\R$, we can consider the complex-valued function
$$ [0,+\infty)\times\R\ni (r,x)\longmapsto \frac{1}{P(r e^{2\pi i x})}.$$
We remark that
\begin{eqnarray*}
&& \frac{d}{dr}\left(\frac{1}{P(r e^{2\pi i x})}\right)=
\frac{e^{2\pi i x}\,P'(r e^{2\pi i x})}{\big(P(r e^{2\pi i x})\big)^2}
\\ {\mbox{and }}&&\frac{d}{dx}\left(\frac{1}{P(r e^{2\pi i x})}\right)=
\frac{2\pi ir e^{2\pi i x}\,P'(r e^{2\pi i x})}{\big(P(r e^{2\pi i x})\big)^2}
\end{eqnarray*}
and therefore
\begin{equation}\label{7.0} \frac{d}{dx}\left(\frac{1}{P(r e^{2\pi i x})}\right)=
2\pi ir\,\frac{d}{dr}\left(\frac{1}{P(r e^{2\pi i x})}\right).\end{equation}

We write
$$P(z)=\sum_{k=0}^N c_k \,z^k,$$ for some~$N\ge1$, $c_0,\dots,c_N\in\C$, and~$c_N\ne0$.

Accordingly, if~$r$ is sufficiently large,
$$ \sup_{x\in\R} \left|\frac{1}{r\,P(r e^{2\pi i x})}\right|\le\frac{1}{r\left(\displaystyle|c_N|\,r^N-
\sum_{k=1}^{N-1} |c_k |\,r^k\right)}\le\frac{1}{2|c_N|\,r^{N+1}}\le\frac2{|c_N|\,r^{2}}.
$$
Hence, given~$\delta>0$, to be taken as small as we wish here below, we can define
$$ f(x):=\int_\delta^{+\infty} \frac{1}{r\,P(r e^{2\pi i x})}\,dr$$
and we observe that~$f$ is periodic of period~$1$.

We claim that
\begin{equation}\label{7.1}
f'(x)=\int_\delta^{+\infty}\frac1r \,\frac{d}{dx}\left(\frac{1}{P(r e^{2\pi i x})}\right)\,dr.
\end{equation}
To check this, we observe that
\begin{eqnarray*}&&
\lim_{\epsilon\to0} \left| \frac{f(x+\epsilon)-f(x)}\epsilon-\int_\delta^{+\infty}\frac1r \,\frac{d}{dx}\left(\frac{1}{P(r e^{2\pi i x})}\right)\,dr\right|\\&\le&
\lim_{\epsilon\to0} \int_\delta^{+\infty}\frac1r\left| 
\frac1\epsilon\left(\frac{1}{P(r e^{2\pi i (x+\epsilon)})}-\frac{1}{P(r e^{2\pi i x})}\right)
-\frac{d}{dx}\left(\frac{1}{P(r e^{2\pi i x})}\right)\right|\,dr\\&=&
\lim_{\epsilon\to0} \int_\delta^{+\infty}\frac1r\left| 
\frac1\epsilon\int_0^\epsilon \frac{d}{d\tau}
\left(\frac{1}{P(r e^{2\pi i (x+\tau)})}\right)\,d\tau
-\frac{d}{dx}\left(\frac{1}{P(r e^{2\pi i x})}\right)\right|\,dr\\&=&
\lim_{\epsilon\to0} \int_\delta^{+\infty}\frac1r\left| 
\frac1\epsilon\int_0^\epsilon \frac{d}{dx}
\left(\frac{1}{P(r e^{2\pi i (x+\tau)})}\right)\,d\tau
-\frac{d}{dx}\left(\frac{1}{P(r e^{2\pi i x})}\right)\right|\,dr\\&\le&\lim_{\epsilon\to0} \int_\delta^{+\infty}\frac1{\epsilon r}\left[ 
\int_0^\epsilon\left| \frac{d}{dx}
\left(\frac{1}{P(r e^{2\pi i (x+\tau)})}\right)
-\frac{d}{dx}\left(\frac{1}{P(r e^{2\pi i x})}\right)\right| \,d\tau\right]\,dr\\&\le&\lim_{\epsilon\to0} \int_\delta^{+\infty}\frac1{\epsilon r}\left[ 
\int_0^\epsilon\left( \int_0^\tau\left| \frac{d^2}{dx^2}
\left(\frac{1}{P(r e^{2\pi i (x+\theta)})}\right)\right| \,d\theta\right)\,d\tau\right]\,dr.
\end{eqnarray*}
Since, for all~$r\ge0$,
\begin{eqnarray*}&& \sup_{x\in\R}\left| \frac{d^2}{dx^2}
\left(\frac{1}{P(r e^{2\pi i x})}\right)\right|\le(2\pi )^2 r\sup_{x\in\R}\left(
\frac{r\,|P''(r e^{2\pi i x})|}{|P(r e^{2\pi i x})|^2}+\frac{2r\,|P'(r e^{2\pi i x})|^2}{|P(r e^{2\pi i x})|^3}+\frac{|P'(r e^{2\pi i x})|}{|P(r e^{2\pi i x})|^2}\right)\\&&\qquad\le \frac{C}{1+r^N}\le \frac{C}{r},\end{eqnarray*}
for some~$C\ge0$, we conclude that
\begin{eqnarray*}&&
\lim_{\epsilon\to0} \left| \frac{f(x+\epsilon)-f(x)}\epsilon-\int_\delta^{+\infty}\frac1r \,\frac{d}{dx}\left(\frac{1}{P(r e^{2\pi i x})}\right)\,dr\right|\\&\le&
C\lim_{\epsilon\to0} \int_\delta^{+\infty}\frac1{\epsilon r^{2}}\left[ 
\int_0^\epsilon\left( \int_0^\tau\,d\theta\right)\,d\tau\right]\,dr
\\&=&C\lim_{\epsilon\to0}\frac{\epsilon}{2\delta}
\\&=&0,\end{eqnarray*}
proving~\eqref{7.1}.

Now we observe that, on the one hand, by~\eqref{7.0} and~\eqref{7.1},
\begin{eqnarray*} \int_0^1 f'(x)\,dx&=&\int_0^1\left[\int_\delta^{+\infty}\frac1r \,\frac{d}{dx}\left(\frac{1}{P(r e^{2\pi i x})}\right)\,dr\right]\,dx
\\&=&2\pi i\int_0^1\left[\int_\delta^{+\infty}\frac{d}{dr}\left(\frac{1}{P(r e^{2\pi i x})}\right)\,dr\right]\,dx
\\&=&-2\pi i\int_0^1\frac{1}{P(\delta e^{2\pi ix})}\,dx.
\end{eqnarray*}
On the other hand, since~$f$ is periodic of period~$1$,
$$ \int_0^1 f'(x)\,dx= f(1)-f(0)=0.$$

We have thereby shown that
$$ \int_0^1\frac{1}{P(\delta e^{2\pi ix})}\,dx=0.$$
We can now use the Dominated Convergence Theorem (recall that~$P$ is supposed to never vanish) and conclude that
$$ \frac{1}{P(0)}=
\int_0^1\frac{1}{P(0)}\,dx=
\lim_{\delta\searrow0}\int_0^1\frac{1}{P(\delta e^{2\pi ix})}\,dx=0.$$
We have thus obtained the desired contradiction.

\paragraph{Solution to Exercise~\ref{JSLD.023oler-NOMPO}.}
Since~$|\cos(100^k\pi x)|\le1$, the series in~\eqref{JSLD.023oler-NOMPOn} converges uniformly, thus defining a continuous function.

Now, seeking a contradiction, suppose that~$f$ is monotonic in some interval~$I$. Let~$\alpha<\beta$ be the endpoints of this interval and pick~$p\in\N$, $p\ge1$, sufficiently large such that
\begin{equation}\label{91wieyhd0ritojgno5nbfdw4r2tgh62yh1j28-00}
100^{2-p}<\beta-\alpha.\end{equation}

Let also~$m$ be the smallest integer for which
\begin{equation}\label{91wieyhd0ritojgno5nbfdw4r2tgh62yh1j28-0} \frac{m}{100^p}\ge\frac{\alpha+\beta}2.\end{equation}

We consider the points
\begin{equation*} a:=\frac{m}{100^p},\qquad \underline a:=a-\frac{1}{100^{p+1}}\qquad{\mbox{and}}\qquad\overline a:=a+\frac{1}{100^{p+1}}.\end{equation*}
Notice that
\begin{equation}\label{91wieyhd0ritojgno5nbfdw4r2tgh62yh1j28-1}
\underline a< a<\overline a.
\end{equation}

We claim that
\begin{equation}\label{91wieyhd0ritojgno5nbfdw4r2tgh62yh1j28-2}
\underline a>\alpha.
\end{equation}
Indeed, by~\eqref{91wieyhd0ritojgno5nbfdw4r2tgh62yh1j28-00} and~\eqref{91wieyhd0ritojgno5nbfdw4r2tgh62yh1j28-0},
$$\underline a-\alpha=\frac{m}{100^p}-\frac{1}{100^{p+1}}-\alpha\ge\frac{\alpha+\beta}2-\frac{1}{100^{p+1}}-\alpha=
\frac{\beta-\alpha}2-\frac{1}{100^{p+1}}\ge\frac{\beta-\alpha}2-\frac{1}{100^{p}}>0,$$
showing the validity of~\eqref{91wieyhd0ritojgno5nbfdw4r2tgh62yh1j28-2}.

We also claim that
\begin{equation}\label{91wieyhd0ritojgno5nbfdw4r2tgh62yh1j28-3}
\overline a<\beta.
\end{equation}
For this, we use the fact that~$m$ is the smallest possible integer fulfilling~\eqref{91wieyhd0ritojgno5nbfdw4r2tgh62yh1j28-0} to see that
$$\frac{m-1}{100^p}\le\frac{\alpha+\beta}2$$
and thus, by~\eqref{91wieyhd0ritojgno5nbfdw4r2tgh62yh1j28-00},
$$\overline a-\beta =\frac{m}{100^p}+\frac{1}{100^{p+1}}-\beta\le
\frac{1}{100^p}+\frac{\alpha+\beta}2+\frac{1}{100^{p+1}}-\beta
\le\frac{1}{100^{p-1}}-\frac{\beta-\alpha}2<0,
$$
which is~\eqref{91wieyhd0ritojgno5nbfdw4r2tgh62yh1j28-3}.

As a consequence of~\eqref{91wieyhd0ritojgno5nbfdw4r2tgh62yh1j28-1}, \eqref{91wieyhd0ritojgno5nbfdw4r2tgh62yh1j28-2}, and~\eqref{91wieyhd0ritojgno5nbfdw4r2tgh62yh1j28-3}, we obtain that the three points~$\underline a<a<\overline a$ belong to the supposed monotonicity interval~$I$.

Hence, a contradiction will be reached once we show that
\begin{equation}\label{91wieyhd0ritojgno5nbfdw4r2tgh62yh1j28-4}
{\mbox{$f(a)-f(\underline a)$ and $f(\overline a)-f(a)$ have opposite signs.}}
\end{equation}
To this end, it is useful to notice that
$$ \sum_{k=0}^{+\infty}\frac1{2^k}=\frac1{1-\frac12}=2$$
and so we replace~$f$ in~\eqref{JSLD.023oler-NOMPOn} with~$
\widetilde f:=2-f$,
which would maintain the same monotonicity properties
of the original~$f$, that is we look at a new function
in the form
\begin{equation}\label{JSLD.023oler-NOMPOn1} \widetilde f(x)=\sum_{k=0}^{+\infty}\frac{1-\cos(100^k\pi x)}{2^k}.\end{equation}
The technical advantage of~\eqref{JSLD.023oler-NOMPOn1}
over~\eqref{JSLD.023oler-NOMPOn} is that the summand is now non-negative, thus making it easier to spot algebraic simplifications.

We observe that, when~$k\ge p+1$,
$$ 100^k\pi a=100^{k-p} m\pi\in 2\pi\Z$$
and thus~$\cos(100^k\pi a)=1$.

This and~\eqref{JSLD.023oler-NOMPOn1} yield that
\begin{equation}\label{JSLD.023oler-NOMPOn1z3z} \widetilde f(a)=\sum_{k=0}^{p}\frac{1-\cos(100^k\pi a)}{2^k}.\end{equation}

Moreover, when~$k\ge p+2$,
$$ 100^k\pi\left(a\pm \frac{1}{100^{p+1}}\right)=
100^{k-p} m\pi\pm 100^{k-p-1}\pi\in2\pi\Z$$
and thus, by~\eqref{JSLD.023oler-NOMPOn1},
\begin{equation*} \widetilde f(\underline a)=\sum_{k=0}^{p+1}\frac{1-\cos(100^k\pi \underline a)}{2^k}\end{equation*}
and
\begin{equation*} \widetilde f(\overline a)=\sum_{k=0}^{p+1}\frac{1-\cos(100^k\pi \overline a)}{2^k}.\end{equation*}

From these observations and~\eqref{JSLD.023oler-NOMPOn1z3z} it follows that
\begin{equation}\label{ABDFMCSCKRT-01}\begin{split}&
\widetilde f(\overline a)-\widetilde f(a)=\sum_{k=0}^{p}\frac{\cos(100^k\pi a)-\cos(100^k\pi \overline a)}{2^k}
+\frac{1-\cos(100^{p+1}\pi \overline a)}{2^{p+1}}
\end{split}\end{equation}
and similarly
\begin{equation}\label{ABDFMCSCKRT-02}\begin{split}&
\widetilde f(\underline a)-\widetilde f(a)= 
\sum_{k=0}^{p}\frac{\cos(100^k\pi a)-\cos(100^k\pi \underline a)}{2^k}
+\frac{1-\cos(100^{p+1}\pi \underline a)}{2^{p+1}}.
\end{split}\end{equation}

Now, since
$$ |\overline a-a|=|\underline a-a|=\frac{1}{100^{p+1}}$$
and, for all~$\alpha$, $\beta\in\R$,
$$ |\cos\alpha-\cos\beta|=\left| \int^{\alpha}_\beta\sin t\,dt\right|\le|\alpha-\beta|,$$
we see that
\begin{equation}\label{DFAVCKTMZACMD-01} \left|\sum_{k=0}^{p}\frac{\cos(100^k\pi a)-\cos(100^k\pi \overline a)}{2^k}\right|\le
\sum_{k=0}^{p}50^k \pi | \overline a-a|\le\frac{50^{p+1}\pi}{49\cdot100^{p+1}}=\frac{\pi}{49\cdot2^{p+1}}\end{equation}
and, in a similar fashion,
\begin{equation} \label{DFAVCKTMZACMD-02}\left|\sum_{k=0}^{p}\frac{\cos(2^k\pi a)-\cos(2^k\pi \overline a)}{2^k}\right|\le\frac{\pi}{49\cdot2^{p+1}}.\end{equation}

In addition, 
$$ 100^{p+1}\pi \overline a=
100^{p+1}\pi\left(\frac{m}{100^p}+\frac{1}{100^{p+1}}\right)=
100 m\pi +\pi\equiv\pi\qquad({\mbox{mod }}2\pi)
$$
and therefore
\begin{equation}\label{DFAVCKTMZACMDx-01}
\frac{1-\cos(100^{p+1}\pi \overline a)}{2^{p+1}}=
\frac{1-\cos\pi}{2^{p+1}}=\frac{1}{2^{p}}.
\end{equation}
In the same vein,
\begin{equation}\label{DFAVCKTMZACMDx-02}
\frac{1-\cos(100^{p+1}\pi \underline a)}{2^{p+1}}= \frac{1}{2^{p}}.\end{equation}

We deduce from~\eqref{ABDFMCSCKRT-01}, \eqref{DFAVCKTMZACMD-01}, and~\eqref{DFAVCKTMZACMDx-01} that
\begin{equation}\label{DFAVCKTMZACMDbb-01}
\widetilde f(\overline a)-\widetilde f(a)\ge -\frac{\pi}{49\cdot2^{p+1}}+\frac{1}{2^{p}}>0.
\end{equation}
Likewise, we infer from~\eqref{ABDFMCSCKRT-02}, \eqref{DFAVCKTMZACMD-02}, and~\eqref{DFAVCKTMZACMDx-02} that
\begin{equation}\label{DFAVCKTMZACMDbb-02}
\widetilde f(\underline a)-\widetilde f(a) >0.
\end{equation}

By virtue of~\eqref{DFAVCKTMZACMDbb-01} and~\eqref{DFAVCKTMZACMDbb-02},
the proof of~\eqref{91wieyhd0ritojgno5nbfdw4r2tgh62yh1j28-4} is thereby complete.

\begin{figure}[h]
\includegraphics[height=2.8cm]{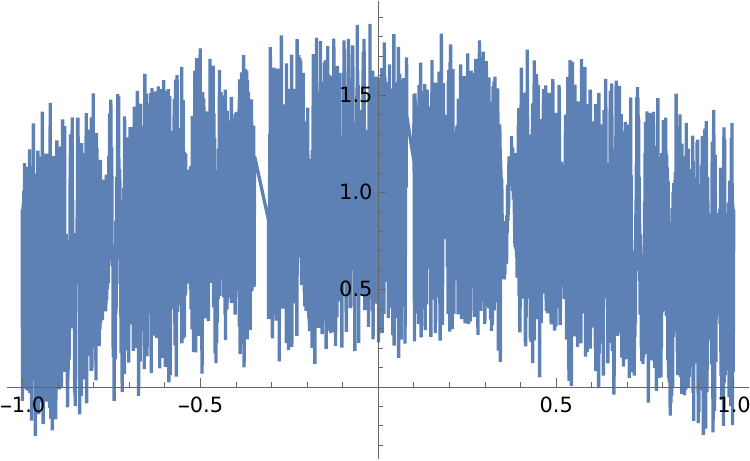}$\,\;\quad$\includegraphics[height=2.8cm]{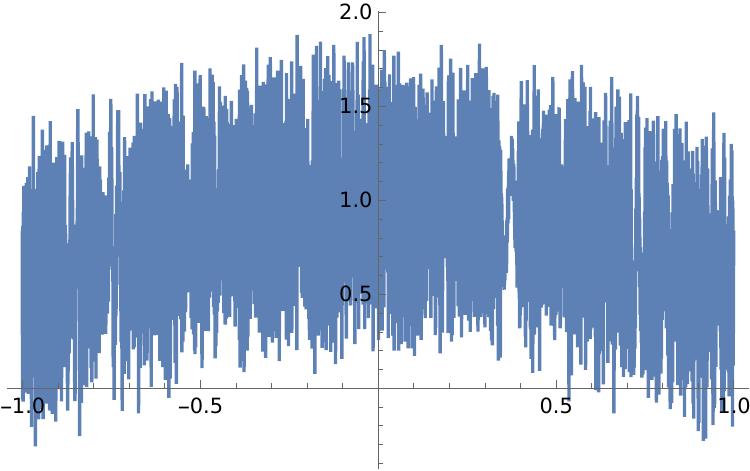}$\,\;\quad$
\includegraphics[height=2.8cm]{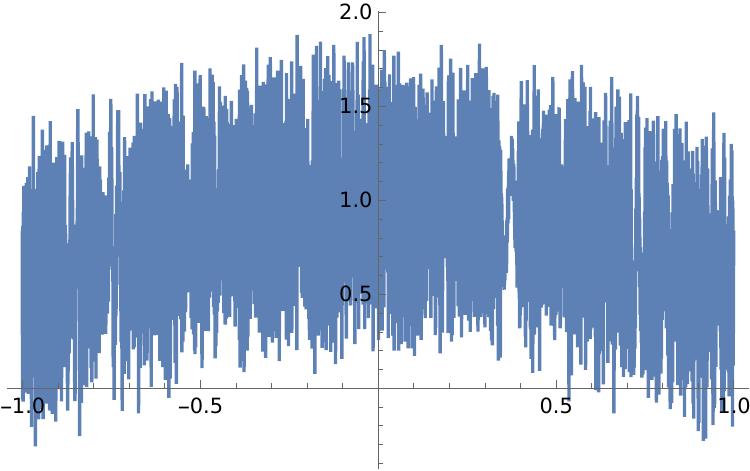}
\centering
\caption{Plot of the function~$\displaystyle\sum_{k=0}^{N}\frac{ \cos(100^k\pi x)}{2^k}$
with~$N\in\{3,20,100\}$.}\label{MO9ikm9400mOntrtm0g.1vcl8b9mb8n.ikdnsEU}
\end{figure}

Of course, visualising ``monsters'' as the functions considered here is tricky, even with computers, because wild oscillations at small scales may just ``thicken the graph'' without being clearly visible with the naked eye.
In any case, a rough sketch of the situation considered in this exercise is given in Figure~\ref{MO9ikm9400mOntrtm0g.1vcl8b9mb8n.ikdnsEU}.

These functions are special cases of a class of functions which was first introduced by Weierstrass, see~\cite[pages~29, 38--39]{MR1996162} and~\cite[pages~48--49]{MR655599} (see also Exercise~\ref{DDEF8AMPMvS1}).

\paragraph{Solution to Exercise~\ref{DDEF8AMPMvS1}.} This is a slight sharpening of the argument in Exercise~\ref{JSLD.023oler-NOMPO}. Full, self-contained details are provided for convenience.
To find a contradiction, suppose that~$f$ is differentiable at some point~$x_0\in\R$.
Then,
$$ \lim_{h\to0}\frac{f(x_0+h)-f(x_0)}{h}=f'(x_0)$$
and thus for every~$\epsilon>0$ there exists~$\delta_\epsilon>0$ such that if~$h\in[-\delta_\epsilon,\delta_\epsilon]$ then
$$\left|\frac{f(x_0+h)-f(x_0)}{h}-f'(x_0)\right|\le\epsilon.$$

As a result, if~$\underline x\in[x_0-\delta_\epsilon,x_0]$ and~$\overline x\in[x_0,x_0+\delta_\epsilon]$, with~$\underline x<\overline x$,
we find that
\begin{equation}\label{PSJL209iuh-9iuhbv3d.M-eps2}\begin{split}&
\left|\frac{f(\overline x)-f(\underline x)}{\overline x-\underline x}-f'(x_0)\right|=
\left|\frac{f(\overline x)-f(x_0)}{\overline x-\underline x}-\frac{f(\underline x)-f(x_0)}{\overline x-\underline x}-f'(x_0)\right|\\
&\qquad=\left|\frac{\overline x-x_0}{\overline x-\underline x}\left(
\frac{f(\overline x)-f(x_0)}{\overline x-x_0}-f'(x_0)\right)-\frac{\underline x-x_0}{\overline x-\underline x}\left(
\frac{f(\underline x)-f(x_0)}{\underline x-x_0}-f'(x_0)\right)\right|\\&\qquad\le
\epsilon\,\left( 
\frac{|\overline x-x_0|}{|\overline x-\underline x|}+\frac{|\underline x-x_0|}{|\overline x-\underline x|}
\right)\\&\qquad\le2\epsilon.
\end{split}\end{equation}

We will take~$\epsilon:=1$ and call~$\delta:=\delta_1$.
We also take~$p\in\N$ large enough such that~$\frac2{100^p}<\delta$
and we look at the function~$\R\ni x\mapsto\cos\big(100^p\pi(x_0+x)\big)$, which is periodic of period~$\frac2{100^p}$.
Hence, we find points~$\underline x\in \left[x_0-\frac2{100^p},x_0\right)$ and~$\overline x\in \left(x_0,x_0+\frac2{100^p}\right]$ such that~$\cos(100^p\pi\underline x)=-1$ and~$\cos(100^p\pi\overline x)=1$.

Namely, $100^p\pi\underline x\equiv\pi$ and~$100^p\pi\overline x\equiv0$
$({\mbox{mod }}2\pi)$.

Therefore, if~$k\ge p+1$ then~$100^k\pi\underline x=
100^{k-p}100^p\pi\underline x\equiv100^{k-p}\pi\equiv0$ and~$100^k\pi\overline x=100^{k-p}100^p\pi\overline x\equiv0$
$({\mbox{mod }}2\pi)$.

As a result,
\begin{eqnarray*}
f(\overline x)-f(\underline x)&=&\sum_{k=0}^{p-1}\frac{ \cos(100^k\pi \overline x)-\cos(100^k\pi\underline x)}{2^k}
+\frac2{2^p}\\
&\ge&-\sum_{k=0}^{p-1}\frac{ 100^k\pi |\overline x-\underline x|}{2^k}
+\frac2{2^p}\\&=&-\frac{ (50^p-1)\pi |\overline x-\underline x|}{49}
+\frac2{2^p}\\&\ge&-\frac{ 4\pi}{49\cdot2^p}
+\frac2{2^p}\\&=&\frac{c}{2^p},
\end{eqnarray*}
for some~$c>0$.

Consequently,
$$ \frac{f(\overline x)-f(\underline x)}{\overline x-\underline x}\ge\frac{c}{2^p(\overline x-\underline x)}\ge\frac{c\;50^p}4$$
and thus, comparing with~\eqref{PSJL209iuh-9iuhbv3d.M-eps2},
\begin{eqnarray*}2\ge
\left|\frac{f(\overline x)-f(\underline x)}{\overline x-\underline x}-f'(x_0)\right|
\ge \frac{f(\overline x)-f(\underline x)}{\overline x-\underline x}-f'(x_0)\ge\frac{c\;50^p}4-f'(x_0).\end{eqnarray*}
We can now take~$p$ as large as we want and obtain the desired contradiction.

\section{Solutions to selected exercises of Section~\ref{L2pe}}

\paragraph{Solution to Exercise~\ref{PKS0-3-21}.}
We deduce from~\eqref{GL4} that the~$j$th Fourier coefficient of the function~$f_k(x):=
e^{2\pi ikx}$ is~$ \widehat f_{k,j}=\delta_{k,j}$.

Hence (see Exercise~\ref{LINC-bisco}) we have that
the~$j$th Fourier coefficient of the function~$f$ is
$$ \widehat f_j=\sum_{{k\in\Z}\atop{|k|\le n}} c_k\,\widehat f_{k,j}=\sum_{{k\in\Z}\atop{|k|\le n}} c_k\,\delta_{k,j}=
\begin{dcases}
c_j &{\mbox{ if }}|j|\le n,\\
0&{\mbox{ if }}|j|> n.
\end{dcases}$$

For this reason,
$$ \sum_{{j\in\Z}\atop{|j|\le N}}\widehat f_j\,e^{2\pi ijx}=\sum_{{j\in\Z}\atop{|j|\le \min\{n,N\}}}c_j\,e^{2\pi ijx}
$$
and, considering the case when~$N\ge n$, the desired result follows.

\paragraph{Solution to Exercise~\ref{NICELI}.}
We write~$F$ in its power series representation in the open disk of radius~$R$ centred at~$z_0$
(see e.g.~\cite[Corollary~4.3 and Theorem~4.4]{MR1976398}), namely for all~$z\in\C$ with~$|z-z_0|<R$ we have that
\begin{equation}\label{DnKOCKSTF.0-1} F(z)=\sum_{j=0}^{+\infty}\frac{D^jF(z_0)}{j!}(z-z_0)^j,\end{equation}
with, for all~$j\in\N$ and~$\rho\in(0,R)$,
\begin{equation}\label{DnKOCKSTF.0-2} \frac{D^jF(z_0)}{j!}\le\frac{M}{\rho^j},\end{equation}
for a suitable~$M\ge0$.

From now on, we will suppose that~$r\in(0,R)$ is given and we take~$\rho\in(r,R)$. Given~$\theta\in\R$, we define
$$ f(\theta):=F\big(z_0+re^{2\pi i\theta}\big)$$
and we deduce from~\eqref{DnKOCKSTF.0-1} that
$$ f(\theta)=\sum_{j=0}^{+\infty}\frac{D^jF(z_0)\;r^j}{j!} e^{2\pi ij\theta}.$$
We stress that, by virtue of~\eqref{DnKOCKSTF.0-2}, the above series is uniformly convergent for~$\theta$ in~$\R$
and accordingly~$f$ is a continuous function, periodic of period~$1$.

We can thereby compute its Fourier coefficients by means of~\eqref{FOUCO}, finding that, for all~$k\in\N$,
$$ \widehat f_k=\int_0^1 f(\theta)\,e^{-2\pi i k\theta}\,d\theta
=\sum_{j=0}^{+\infty}\int_0^1\frac{D^jF(z_0)\;r^j}{j!} e^{2\pi i(j-k)\theta}\,d\theta,$$
where the swapping of summation and integral signs is made possible by the uniform convergence of the series.

Hence, calculating the integral or recalling~\eqref{GL4}, we arrive at
\begin{equation}\label{AJSa.0olXZX2} \widehat f_k=\frac{D^k F(z_0)\;r^k}{k!}.\end{equation}

Moreover, we can parameterise~${\mathcal{C}}_r$ by~$\{z_0+re^{2\pi i\theta}$, $\theta\in[0,1)\}$ and thus obtain that, for all~$k\in\N$,
$$ \oint_{{\mathcal{C}}_r}\frac{F(z)}{(z-z_0)^{k+1}}\,dz=2\pi ir
\int_0^1 \frac{F\big(z_0+re^{2\pi i\theta}\big)}{\big(re^{2\pi i\theta}\big)^{k+1}}\,e^{2\pi i\theta}d\theta=
\frac{2\pi i}{r^k}
\int_0^1 f(\theta)\,e^{-2\pi ik\theta}d\theta=
\frac{2\pi i}{r^k}\,\widehat f_k.
$$
Comparing this with~\eqref{AJSa.0olXZX2} we see that
$$ \frac{D^k F(z_0)\;r^k}{k!}=\widehat f_k=\frac{r^k}{2\pi i}\oint_{{\mathcal{C}}_r}\frac{F(z)}{(z-z_0)^{k+1}}\,dz,$$
leading to the desired result after simplifying~$r^k$.

\paragraph{Solution to Exercise~\ref{NONSepCFGDV2}.} 
By the Cauchy Integral Formula (see e.g.~\cite[equation~(4.44)]{MR3024399}),
for all~$j\in\N$ and~$r\in(0,1)$,
\begin{equation}\label{PKSM0-i55kvbmS-103}
D^j F(0)=\frac{j!}{2\pi i}\oint_{{\mathcal{C}}_r}\frac{F(z)}{z^{j+1}}\,dz,\end{equation}
where~${{\mathcal{C}}_r}$
is the circle of radius~$r$ in the complex plane, travelled anticlockwise.

Also, we let~$u:=\Re F$ and~$v:=\Im F$.
Thus, in view of~\eqref{PKSM0-i55kvbmS-103}, we see that, for all~$r\in(0,1)$,
\begin{equation}\label{PKSM0-i55kvbmS-104}\begin{split}
\frac{r^k}{k!}\,\big| D^k F(0) \big|&=\left| 
\int_0^1 F(re^{2\pi ix})\,e^{-2\pi ikx}\,dx\right|\\&=\left| 
\int_0^1 \big(u(re^{2\pi ix})+iv(re^{2\pi ix})\big)\,e^{-2\pi ikx}\,dx\right|.
\end{split}\end{equation}

Moreover, for all~$k\in\N\cap[1,+\infty)$, the function~$z^{k-1} F(z)$ is holomorphic in the unit disk and consequently,
by the Cauchy Integral Theorem (see e.g.~\cite[Theorem 4.2.11]{MR3024399}), 
$$ 0=\frac1{2\pi ir^{k}}\oint_{{\mathcal{C}}_r}z^{k-1} F(z)\,dz=
\int_0^1 \big(u(re^{2\pi ix})+iv(re^{2\pi ix})\big)\,e^{2\pi ikx}\,dx.
$$
Taking the complex conjugate of this identity, we conclude that, for all~$k\in\N\cap[1,+\infty)$,
$$ 0=\int_0^1 \overline{\big(u(re^{2\pi ix})+iv(re^{2\pi ix})\big)}\,e^{-2\pi ikx}\,dx
=\int_0^1 \big(u(re^{2\pi ix})-iv(re^{2\pi ix})\big)\,e^{-2\pi ikx}\,dx.
$$

From this and~\eqref{PKSM0-i55kvbmS-104} we gather that
\begin{equation}\label{PKSM0-i55kvbmS-105}\begin{split}
\frac{r^k}{k!}\,\big| D^k F(0) \big|=2\left| 
\int_0^1 u(re^{2\pi ix})\,e^{-2\pi ikx}\,dx\right|.
\end{split}\end{equation}

Now we define
$$ f_r(x):= u(r e^{2\pi ix})$$
and we denote by~$\widehat f_{r,k}$ the $k$th Fourier coefficient of~$f_r$.

Hence, we infer from~\eqref{PKSM0-i55kvbmS-105} that
$$ \frac{r^k}{k!}\,\big| D^k F(0) \big|=2\left|\widehat f_{r,k}\right|.$$
This and Exercise~\ref{NONSepCFGDV} yield that
\begin{equation}\label{PKSM0-i55kvbmS-106} \frac{r^k}{k!}\,\big| D^k F(0) \big|=2\widehat f_{r,0}.\end{equation}

Furthermore, on account of~\eqref{PKSM0-i55kvbmS-103},
\begin{eqnarray*}&& \Re F(0)=\Re\left(\frac{1}{2\pi i}\oint_{{\mathcal{C}}_r}\frac{F(z)}{z}\,dz\right)=\Re\left(\int_0^1F(re^{2\pi ix})\,dx\right)\\&&\qquad=
\int_0^1u(re^{2\pi ix})\,dx=\int_0^1f_r(x)\,dx=\widehat f_{r,0}.
\end{eqnarray*}
We insert this information in~\eqref{PKSM0-i55kvbmS-106} and we conclude that
$$ \frac{r^k}{k!}\,\big| D^k F(0) \big|\le2\Re F(0).$$
The desired result now follows by sending~$r\nearrow1$.

We point out that this exercise shows a nice interplay between analytic functions (which are typically related to
Taylor Series) and Fourier coefficient (which are typically related to Fourier Series).
See also footnote~\ref{LAFOOTTF} on page~\pageref{LAFOOTTF}.

\paragraph{Solution to Exercise~\ref{NONSepCFGDV2-cont}.} No. Take for instance~$F(z):=i$.

\paragraph{Solution to Exercise~\ref{ILTRI-FAA}.}
The desired result can be proved by induction over~$\ell$.

Alternatively, one can use Fa\`a di Bruno's Formula for the derivative of the composition (see e.g.~\cite[Theorem~1.3.2]{MR1916029}).

\paragraph{Solution to Exercise~\ref{ILTRI}.}
Suppose by contradiction that there exist a trigonometric polynomial~$f$ of degree~$N$, non-negative integers~$k_1,\dots,k_\ell$,
and distinct points~$\{x_j\}_{j\in\{1,\dots,L\}}\in I$ with~$D^\ell f(x_j)=0$ for all~$\ell\in\{0,\dots,k_j\}$ and
\begin{equation}\label{TROP} \sum_{j=1}^L (k_j+1)> 2N.\end{equation}
We write
$$ f(x)=\sum_{{h\in\Z}\atop{|h|\le N}}c_h\,e^{2\pi ihx}.$$
We let~$C_m:=c_{m-N}$ and consider the complex polynomial
$$ P(z):=z^N\sum_{{h\in\Z}\atop{|h|\le N}}c_h\,z^h=
\sum_{m=0}^{2N} C_m z^m.$$
Let also~$z_j:=e^{2\pi i x_j}$ and note that they are different complex numbers for~$j\in\{1,\dots,L\}$.

Moreover,
$$ P(e^{2\pi ix})=e^{2\pi iNx} f(x)$$
and accordingly, for all~$j\in\{1,\dots,L\}$ and~$\ell\in\{0,\dots,k_j\}$, the Product Rule for derivatives give that
\begin{eqnarray*}\Psi_\ell(x):=
D^\ell_x\big(P(e^{2\pi ix})\big)=D^\ell_x\big(e^{2\pi iNx} f(x)\big)=
\sum_{q=0}^\ell \left( {\ell}\atop{q}\right) D^{\ell-q}_x(e^{2\pi iNx})\,D^qf(x),
\end{eqnarray*}
yielding that
\begin{equation*}
\Psi_\ell(x_j)=0.
\end{equation*}

In addition, owing to Exercise~\ref{ILTRI-FAA},
$$ \Psi_\ell(x)=e^{2\pi i\ell x}\,D^\ell P(e^{2\pi ix})+\sum_{j=1}^{\ell-1} c_{\ell,j}\,e^{2\pi ijx}\,D^j P(e^{2\pi ix}),$$
for suitable coefficients~$c_{\ell,1},\dots,c_{\ell,\ell-1}$, yielding that, for all~$j\in\{1,\dots,L\}$ and~$\ell\in\{0,\dots,k_j\}$,
\begin{equation}\label{ILTRI-FAA-eq01}\begin{split}&
D^\ell P(z_j)+\sum_{j=1}^{\ell-1} c_{\ell,j}\,z_j^{j-\ell}\,D^j P(z_j)\\&\qquad=
D^\ell P(e^{2\pi ix_j})+\sum_{j=1}^{\ell-1} c_{\ell,j}\,e^{2\pi i(j-\ell)x_j}\,D^j P(e^{2\pi ix_j})=0.\end{split}
\end{equation}

We claim that, for all~$j\in\{1,\dots,L\}$ and~$\ell\in\{0,\dots,k_j\}$,
\begin{equation}\label{ILTRI-FAA-eq02}
D^\ell P(z_j)=0.
\end{equation}
The proof of this is by induction over~$\ell$ (for any given~$j\in\{1,\dots,L\}$). When~$\ell=0$, we deduce from~\eqref{ILTRI-FAA-eq01} that~$P(z_j)=0$, as desired.
Hence, to perform the inductive step, we assume that~\eqref{ILTRI-FAA-eq02} holds true for all the indices~$0,\dots,\ell-1$
and we aim at establishing it for the index~$\ell$. To this end, we deduce from the inductive hypothesis that
$$ \sum_{j=1}^{\ell-1} c_{\ell,j}\,e^{2\pi i(j-\ell)x_j}\,D^j P(e^{2\pi ix_j})=0.$$
Accordingly, in virtue of~\eqref{ILTRI-FAA-eq01}, we find that~$D^\ell P(z_j)=0$, thus completing the proof of~\eqref{ILTRI-FAA-eq02}.

In light of~\eqref{TROP} and~\eqref{ILTRI-FAA-eq02}, we conclude that the complex polynomial~$P$ of degree~$2N$ has more than~$2N$ roots, counted with multiplicity,
in contradiction with the Fundamental Theorem of Algebra (see Exercise~\ref{SnsdFrhndiYUJS-2}).

\paragraph{Solution to Exercise~\ref{ILTRI-22}.} No. For example, the trigonometric polynomial~$2+\cos(2\pi Nx)$
has degree~$N$ and does not possess any zero.

\paragraph{Solution to Exercise~\ref{920-334PKSXu9o2fgfbsmos}.}
We prove the first claim, the second being similar. By way of~\eqref{jasmx23er}
and Exercise~\ref{fr12},
\begin{eqnarray*}
&& b_k=
2\int_0^1 f(x)\,\sin(2\pi kx)\,dx=2\int_{-1/2}^{1/2} f(x)\,\sin(2\pi kx)\,dx,\end{eqnarray*}
which vanishes since the integrand is an odd function integrated over a symmetric domain.

\paragraph{Solution to Exercise~\ref{FO:DE:MAGIB}.}
By Exercise~\ref{ojld03-12d},
$$ \sum_{{k\in\Z}{|k|\le N}}\widehat w_k\,e^{2\pi ikx}=\sum_{j=0}^{\widetilde N} \frac4{\pi(2j+1)}\sin(2\pi(2j+1)x),$$where
$$\widetilde N:=\begin{dcases}\displaystyle\frac{N-1}2&{\mbox{ if $N$ is odd,}}\\
\displaystyle\frac{N-2}2&{\mbox{ if $N$ is even.}}
\end{dcases}$$
Hence, for all~$x\in\left(0,\frac12\right)$,
$$ F_N(x):=\sum_{{k\in\Z}{|k|\le N}}\widehat w_k\,e^{2\pi ikx}-w(x)=\sum_{j=0}^{\widetilde N} \frac4{\pi(2j+1)}\sin(2\pi(2j+1)x)-1.$$
On this account, for all~$x\in\left(0,\frac12\right)$,
\begin{equation}\label{EFDEN} \begin{split}F_N'(x)&=8\sum_{j=0}^{\widetilde N} \cos(2\pi(2j+1)x)\\&=
8\,\Re \left(\sum_{j=0}^{\widetilde N} e^{2\pi(2j+1)ix}\right)\\&=
8\,\Re \left( e^{2\pi ix} \sum_{j=0}^{\widetilde N} (e^{4\pi ix})^j\right)\\&=
8\,\Re \left( \frac{e^{2\pi ix}( e^{4\pi i(\widetilde N+1)x}-1)}{ e^{4\pi ix}-1}\right)\\&=
8\,\Re \left( \frac{ e^{4\pi i(\widetilde N+1)x}-1}{ e^{2\pi ix}-e^{-2\pi ix}}\right)\\&=
8\,\Re \left( \frac{ \cos(4\pi (\widetilde N+1)x)+i\sin(4\pi (\widetilde N+1)x)-1}{ 2i\sin(2\pi x)}\right)\\&=
\frac{4\sin(4\pi (\widetilde N+1)x)}{\sin(2\pi x)}.
\end{split}\end{equation}

\begin{figure}[h]
\includegraphics[height=3cm]{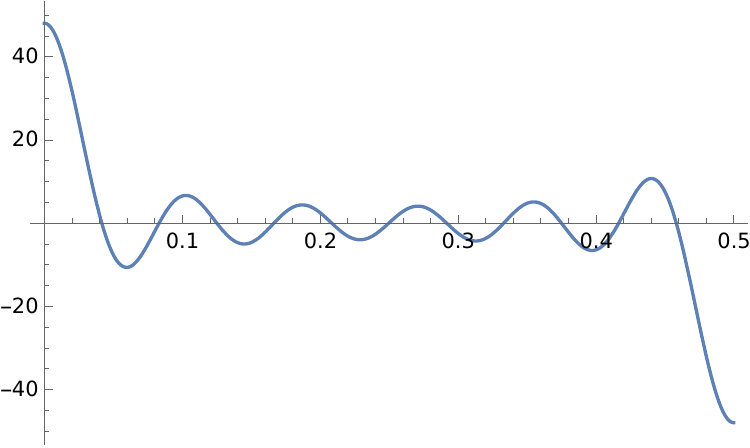}$\quad$\includegraphics[height=3cm]{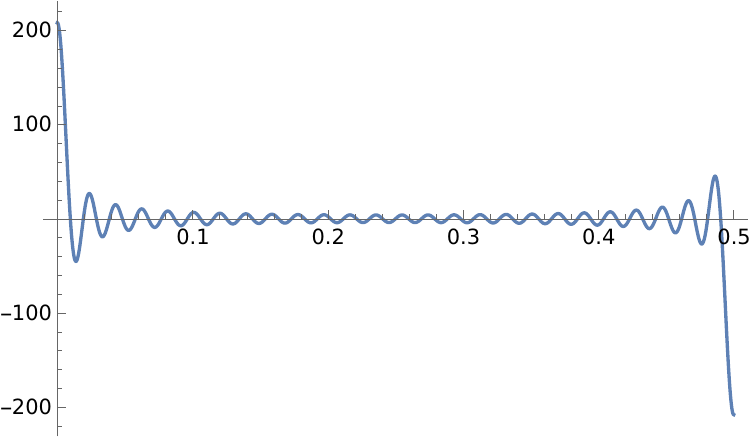}$\quad$\includegraphics[height=3cm]{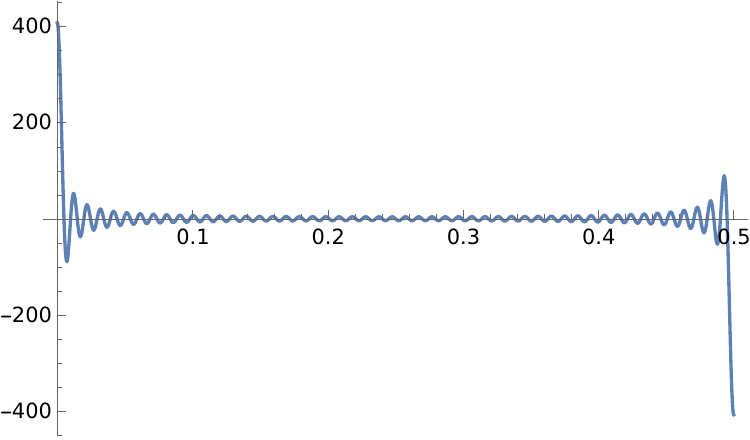}
\centering
\caption{Plot of the function in~\eqref{EFDEN} when~$\widetilde N\in\{5,25,50\}$.}\label{a13e23rt43.Anmdi}
\end{figure}

We observe that, for large~$N$, the function in~\eqref{EFDEN} develops higher and higher oscillations near the origin
(and in general at all points of the form~$\frac{m}2$ with~$m\in\Z$), see Figure~\ref{a13e23rt43.Anmdi}.
This situation will play an important role for the Gibbs phenomenon in Section~\ref{SEC:GIBBS-PH}.

\paragraph{Solution to Exercise~\ref{NU90i3orjf:12oeihfnvZMAS}.}
We have that
$$ \int_{-1/3}^{1/3} x^2\,\sin(2\pi kx)\,dx=0$$
and
$$ \int_{-1/3}^{1/3} x^2\,\cos(2\pi kx)\,dx=\begin{dcases}
\displaystyle\frac2{81}&{\mbox{ if }}k=0,\\
\displaystyle\frac{(2\pi^2k^2-9)S_k+6\pi kC_k}{18\pi^3 k^3} &{\mbox{ if }}k\ne0,\\
\end{dcases}$$ where $$ S_k:=\begin{dcases}
\displaystyle0 &{\mbox{ if }}k\in3\Z,\\
\displaystyle\frac{\sqrt3}{2} &{\mbox{ if }}k\in3\Z+1,\\
\displaystyle-\frac{\sqrt3}{2} &{\mbox{ if }}k\in3\Z+2
\end{dcases}$$ and $$ C_k:=\begin{dcases}
\displaystyle1&{\mbox{ if }}k\in3\Z,\\
\displaystyle-\frac{1}{2} &{\mbox{ if }}k\in(3\Z+1)\cup(3\Z+2).
\end{dcases}$$ 

\begin{figure}[h]
\includegraphics[height=2.3cm]{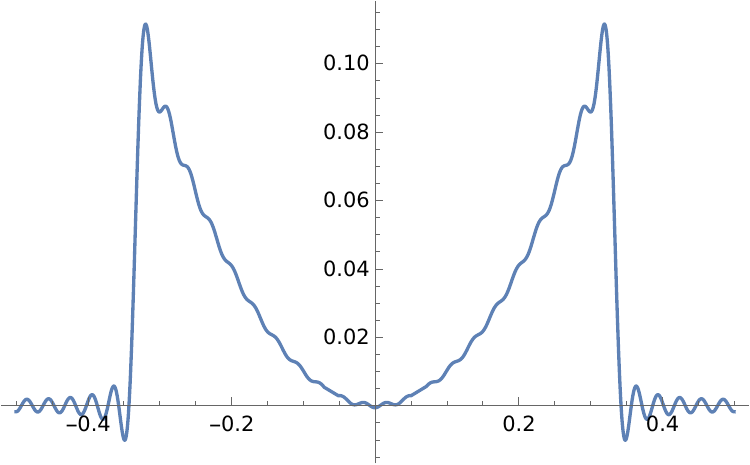}$\,$\includegraphics[height=2.3cm]{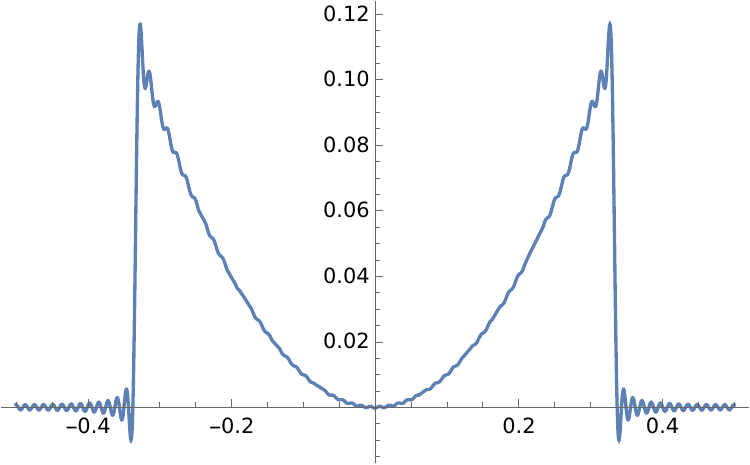}$\,$
\includegraphics[height=2.3cm]{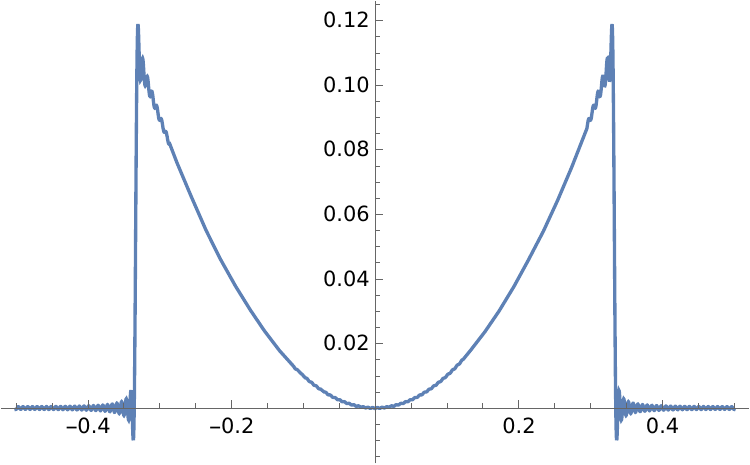}$\,$\includegraphics[height=2.3cm]{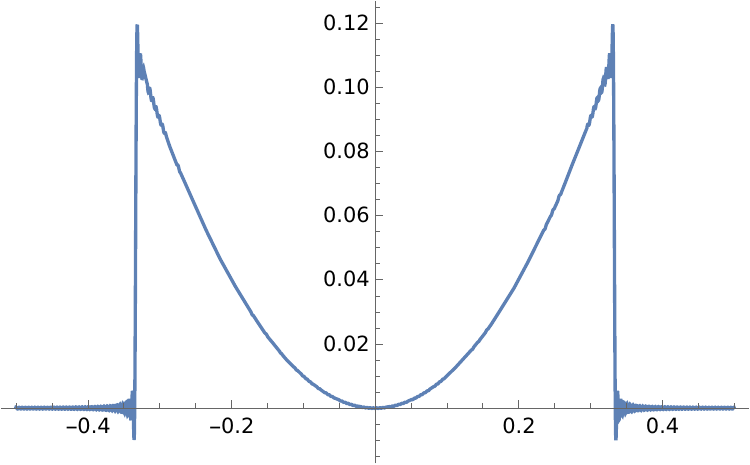}
\centering
\caption{Plot of the sum up to~$N$ approximating the series in~\eqref{02938urfjiikmsiiklo0PSnm}
with~$N\in\{10,\,25,\,50,\,75\}$.}\label{C7n7okjntghbrdGHJKFRO903.lk90fv0-12920045-23-5283012P}
\end{figure}

On this account, the Fourier Series of~$f$ in trigonometric form is \begin{equation}\label{02938urfjiikmsiiklo0PSnm}\begin{split}&
\frac{2}{81}+\sum_{j=1}^{+\infty}
\frac{2}{27\pi^2 j^2}\,\cos(6\pi jx)\\&\qquad+\sum_{j=0}^{+\infty}
\frac{(2\pi^2(3j+1)^2-9)\sqrt{3}-6\pi (3j+1)}{18\pi^3 (3j+1)^3}\,\cos(2\pi(3 j+1)x)\\&\qquad-\sum_{j=0}^{+\infty}
\frac{(2\pi^2(3j+2)^2-9)\sqrt{3}+6\pi (3j+2)}{18\pi^3 (3j+2)^3}\,\cos(2\pi(3 j+2)x),\end{split}\end{equation}
to be compared with Figure~\ref{C7n7okjntghbrdGHJKFRO903.lk90fv0-12920045-23-5283012P}.

\paragraph{Solution to Exercise~\ref{OJSNILCESNpfa.sdwpqoed-23wedf}.}
The $k$th Fourier coefficient of~$f'$ is
\begin{eqnarray*}&& \int_0^1 f'(x)\,e^{-2\pi ikx}\,dx=\int_0^1 
\big( f(x)\,e^{-2\pi ikx}\big)'\,dx+2\pi ik\int_0^1 f(x)\,e^{-2\pi ikx}\,dx\\&&\qquad=
f(1)e^{-2\pi ik}-f(0)+2\pi ik\widehat f_k=
f(1)-f(0)+2\pi ik\widehat f_k=2\pi ik\widehat f_k,\end{eqnarray*}
where the continuity and periodicity of~$f$ has been used in the last step.

\paragraph{Solution to Exercise~\ref{OJSNILCESNpfa.sdwpqoed-23wedf.1}.}
The assumption that~$f$ is continuous is necessary. For instance, if~$f$ is the sawtooth waveform
in Exercise~\ref{SA:W}, we know that, for all~$x\in(0,1)$, $f(x)=x-\frac12$ and therefore~$f'(x)=1$,
giving that the Fourier Series of~$f'$ (periodically extended outside~$(0,1)$) is~$1$.

But in this case the quantity in~\eqref{pqjdwlfmaonsveaN} takes the form
$$  -\sum_{{k\in\Z}\setminus\{0\}}\,e^{2\pi ikx}.$$

\paragraph{Solution to Exercise~\ref{G0ilMajsx912e}.} 
We let~$g:=f'\in L^2((0,1))$ and we consider the Fourier coefficients of~$g$.
Indeed, by Exercise~\ref{OJSNILCESNpfa.sdwpqoed-23wedf}, we know that~$\widehat g_k=2\pi ik\widehat f_k$.

Moreover, we can apply Bessel's Inequality~\eqref{AJSa} and infer that $$
\sum_{k\in\Z} |\widehat g_k|^2 \leq \|g\|^2_{L^2((0,1))}<+\infty.$$

All in all,
$$ \sum_{k\in\Z\setminus\{0\}} |\widehat f_k|=
\frac1{2\pi}\sum_{k\in\Z\setminus\{0\}} \frac{|\widehat g_k|}{|k|}
\le\frac1{2\pi}\sqrt{\sum_{k\in\Z\setminus\{0\}}|\widehat g_k|^2}\;\sqrt{
\sum_{k\in\Z\setminus\{0\}} \frac{1}{|k|^2}}<+\infty,$$ from which the desired result follows.

\paragraph{Solution to Exercise~\ref{G0ilMajsx912e-14}.}
We revisit Exercise~\ref{G0ilMajsx912e} (actually, Exercise~\ref{G0ilMajsx912e} corresponds to the case~$M=1$).

For each~$j\in\{1,\dots,M\}$, we define
$$ g_j(x):=\begin{dcases}
f'(x)&{\mbox{ if }}x\in(a_{j-1},a_j),\\
0&{\mbox{ otherwise}}.
\end{dcases}$$

We claim that, for all~$x\in[0,1)$,
\begin{equation}\label{ANIVAPCSMBNA}
f(x)=f(0)+\sum_{j=1}^M \int_0^x g_j(t)\,dt.
\end{equation}
To check this, we let~$j_x\in\{1,\dots,M\}$ be such that~$x\in[a_{j_x-1},a_{j_x})$ and, on the account of a telescopic sum
and of the Fundamental Theorem of Calculus, we observe that
\begin{eqnarray*}
f(x)-f(0)&=&f(x)-f(a_{j_x-1})+\sum_{j=1}^{j_x-1}\big(f(a_{j})-f(a_{j-1})\big)\\&=&
\int^x_{a_{j_x-1}} f'(t)\,dt+\sum_{j=1}^{j_x-1}\int^{a_j}_{a_{j-1}} f'(t)\,dt\\&=&
\int^x_{a_{j_x-1}} g_{j_x}(t)\,dt+\sum_{j=1}^{j_x-1}\int^{a_j}_{a_{j-1}} g_j(t)\,dt.
\end{eqnarray*}
Hence, recalling that~$g_j$ vanishes outside~$(a_{j-1},a_j)$,
\begin{eqnarray*}
f(x)-f(0)&=&
\int^x_0 g_{j_x}(t)\,dt+\sum_{j=1}^{j_x-1}\int^{x}_{0} g_j(t)\,dt\\&=&\sum_{j=1}^M \int_0^x g_j(t)\,dt,
\end{eqnarray*}
proving~\eqref{ANIVAPCSMBNA}.

Taking the limit as~$x\nearrow1$ in~\eqref{ANIVAPCSMBNA}
and using the periodicity of~$f$, we also see that
$$0=f(1)-f(0)=\sum_{j=1}^M \int_0^1 g_j(t)\,dt.$$

In the wake of~\eqref{ANIVAPCSMBNA} and the latter identity, we find that, for all~$k\in\Z\setminus\{0\}$,
\begin{eqnarray*}
\widehat f_k&=&\int_0^1\left(f(0)+\sum_{j=1}^M \int_0^x g_j(t)\,dt\right)\,e^{-2\pi ikx}\,dx\\
&=&\sum_{j=1}^M\int_0^1\left( \int_t^1 g_j(t)\,e^{-2\pi ikx}\,dx\right)\,dt\\&=&
\frac1{2\pi ik}\sum_{j=1}^M\int_0^1 g_j(t)\,\big(e^{-2\pi ikt}-1\big)\,dt\\&=&
\frac1{2\pi ik}\sum_{j=1}^M\left(\widehat g_{j,k}-\int_0^1 g_j(t)\,dt\right)
\\&=&\frac1{2\pi ik}\sum_{j=1}^M\widehat g_{j,k}.
\end{eqnarray*}

Furthermore, we can apply Bessel's Inequality~\eqref{AJSa} and infer that $$
\sum_{k\in\Z} |\widehat g_{j,k}|^2 \leq \|g_j\|^2_{L^2((0,1))}<+\infty$$
and consequently
$$ \sum_{k\in\Z\setminus\{0\}} |\widehat f_k|=
\frac1{2\pi}\sum_{{k\in\Z\setminus\{0\}}\atop{1\le j\le M} }\frac{|\widehat g_{j,k}|}{|k|}
\le\frac1{2\pi}\sqrt{\sum_{{k\in\Z\setminus\{0\}}\atop{1\le j\le M}}|\widehat g_{j,k}|^2}\;\sqrt{
\sum_{{k\in\Z\setminus\{0\}}\atop{1\le j\le M}} \frac{1}{|k|^2}}<+\infty,$$ from which the desired result follows.

\paragraph{Solution to Exercise~\ref{LADEBD}.}
We stress that the definition of~$f$ is well-posed, because the series in~\eqref{AMs-ccf} is uniformly convergent, thanks to~\eqref{AMs-ccg}. This uniform convergence also allows us to exchange the integral and summation signs and conclude that, for every~$j\in\Z$,
\begin{eqnarray*}
\widehat f_j=\int_0^1 f(x)\,e^{-2\pi ijx}\,dx=\int_0^1 \sum_{k\in\Z} c_k\,e^{2\pi i(k-j)x}\,dx=\sum_{k\in\Z}\int_0^1  c_k\,e^{2\pi i(k-j)x}\,dx=c_j,
\end{eqnarray*}
which is what one is asked to prove.

\paragraph{Solution to Exercise~\ref{FO:DE:MA}.} We consider the periodic extension of~$\phi_a$
(still denoted by~$\phi_a$ for simplicity) and we observe that it is an even function.

As a result, for all~$k\in\N$,
$$\int_{-1/2}^{1/2}\phi_a(x)\,\sin(2\pi kx)\,dx=0$$
and
\begin{eqnarray*}
&& \int_{-1/2}^{1/2}\phi_a(x)\,\cos(2\pi kx)\,dx=2\int_{0}^{1/2}\phi_a(x)\,\cos(2\pi kx)\,dx\\
&&\qquad
=\frac2a\int_{0}^{1/2}\phi\left(\frac{x}a\right)\,\cos(2\pi kx)\,dx=2\int_{0}^{1/(2a)}\phi(y)\,\cos(2\pi ak y)\,dy\\
&&\qquad=2\int_{0}^{1/(2a)}\max\{ 0,\,1-y\}\,\cos(2\pi ak y)\,dy
=2\int_{0}^{1}(1-y)\,\cos(2\pi ak y)\,dy\\&&\qquad=\begin{dcases}\displaystyle\frac{\sin^2(ka\pi)}{ (k \pi a)^2}
&{\mbox{ if }}k\ne0,\\
1&{\mbox{ if }}k=0.\end{dcases}
\end{eqnarray*}

From these considerations, \eqref{jasmx23er} and Exercise~\ref{fr12}, we conclude that the Fourier Series of~$\phi_a$ has the trigonometric form
$$ 1+\sum_{k=1}^{+\infty}\frac{2\sin^2(ka\pi)}{ (k \pi a)^2}\cos(2\pi kx).$$

\begin{figure}[h]
\includegraphics[height=3.8cm]{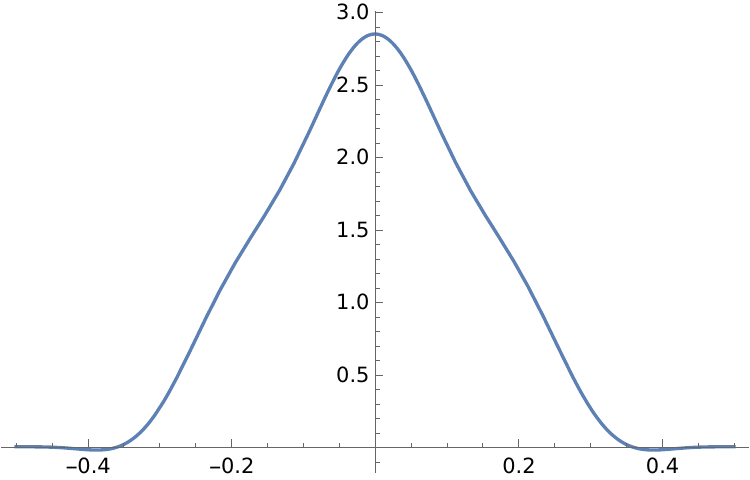}$\,\;\qquad\quad$\includegraphics[height=3.8cm]{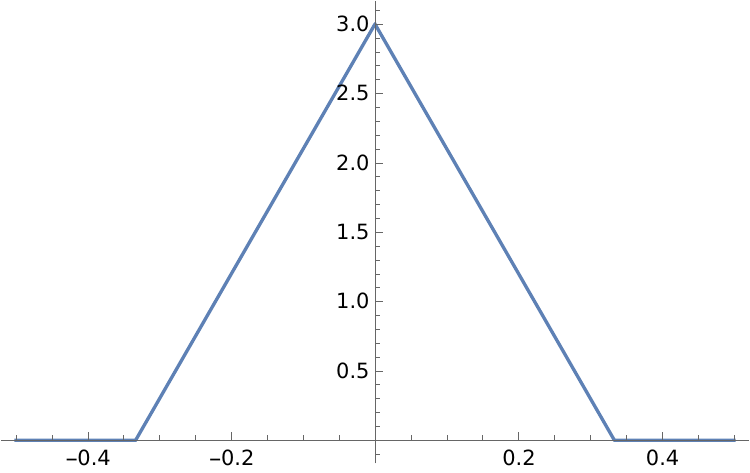}\centering
\caption{Left: plot of the function~$\displaystyle
1+\sum_{k=1}^{5}\frac{2\sin^2\left(\frac{k\pi}3\right)}{ (k \pi /3)^2}\cos(2\pi kx)$.
Right: plot of the function~$\displaystyle3\max\{0,1-3|x|\}$.}\label{C7okjntghbrdfv0-12920045-23-5283012P}
\end{figure}

See Figure~\ref{C7okjntghbrdfv0-12920045-23-5283012P} for a sketch of the situation when~$a=\frac13$.

\paragraph{Solution to Exercise~\ref{RESID}.} If~$f$ denotes the function in this exercise,
we have that, for each~$k\in\Z$,
\begin{eqnarray*}
\widehat f_k&=&\int_0^1 \frac{\cos(2\pi x)+1}{\cos(2\pi x)+5\pi}\,e^{-2\pi ikx}\,dx\\
&=& \int_0^1 \frac{e^{2\pi ix}+e^{-2\pi ix}+2}{e^{2\pi ix}+e^{-2\pi ix}+10\pi}\,e^{-2\pi i(k+1)x}\,e^{2\pi ix}\,dx\\
&=&\frac1{2\pi i}\oint_{\mathcal{C}} F(z)\,dz,
\end{eqnarray*}
where~${\mathcal{C}}$ is the unit circle in the complex plane, travelled anticlockwise, and we have introduced the meromorphic function
$$ \C\ni z\mapsto F(z):=
\frac{z+z^{-1}+2}{z+z^{-1}+10\pi}\,z^{-(k+1)}=
\frac{(z+1)^2}{z^{k+1}(z^2+10\pi z+1)}=\frac{(z+1)^2}{z^{k+1}(z-\alpha)(z-\beta)},$$
with
$$\alpha:= -5 \pi+\sqrt{25 \pi^2 - 1}\qquad{\mbox{and}}\qquad
\beta:= -5 \pi - \sqrt{25 \pi^2 - 1}.$$
We note that~$|\alpha|<1<|\beta|$. 

It is now tempting to use the Residue Theorem (see e.g.~\cite[Theorem~19 in Section~5.1]{MR4506522}) to compute~$\widehat f_k$.

For instance, if~$k=0$, the singularities of~$F$ in the unit disk
are~$z=\alpha$ and~$z=0$, which are simple poles.

As a result,
the sum of the residues of~$F$ inside the unit disk equals
$$ \frac{(\alpha+1)^2}{\alpha(\alpha-\beta)}+\frac{1}{\alpha\beta}=1+\frac{\big(5\pi-1\big)\Big(5\pi-\sqrt{25\pi^2-1}\Big)}{25\pi^2-1-5\pi\sqrt{25\pi^2-1}}$$
and consequently, by the Residue Theorem,
$$ \widehat f_0=1+\frac{\big(5\pi-1\big)\Big(5\pi-\sqrt{25\pi^2-1}\Big)}{25\pi^2-1-5\pi\sqrt{25\pi^2-1}}.$$

When~$k\ge1$, the calculation of the residues of~$F$ inside the unit disk is possibly more complicated, due to the presence of higher order singularities. It is therefore convenient to compute the $k$th Fourier
coefficient when~$k\le-1$ and then obtain the~$k$th Fourier
coefficient when~$k\ge1$ by using~\eqref{fasv}.

More specifically, when~$k\le-1$, the only singularity of~$F$ inside the unit disk is~$z=\alpha$, which is a simple pole, hence the corresponding residue equals
$$ \frac{(\alpha+1)^2}{\alpha^{k+1}(\alpha-\beta)}=\frac{\Big(1 - 5 \pi + \sqrt{25 \pi^2 - 1}\Big)^2 
}{2 \sqrt{25 \pi^2 - 1}\,(\sqrt{25 \pi^2 - 1} - 5 \pi)^{k+1}}.$$
As a consequence, for all~$k\le-1$,
$$ \widehat f_k=\frac{\Big(1 - 5 \pi + \sqrt{25 \pi^2 - 1}\Big)^2 
}{2 \sqrt{25 \pi^2 - 1}\,(\sqrt{25 \pi^2 - 1} - 5 \pi)^{k+1}}.$$

From this and~\eqref{fasv}, we also find that, for all~$k\ge1$,
$$ \widehat f_k=\overline{\widehat f_{-k}}=\frac{\Big(1 - 5 \pi + \sqrt{25 \pi^2 - 1}\Big)^2 
}{2 \sqrt{25 \pi^2 - 1}\,(\sqrt{25 \pi^2 - 1} - 5 \pi)^{1-k}}.$$

\begin{figure}[h]
\includegraphics[height=3.6cm]{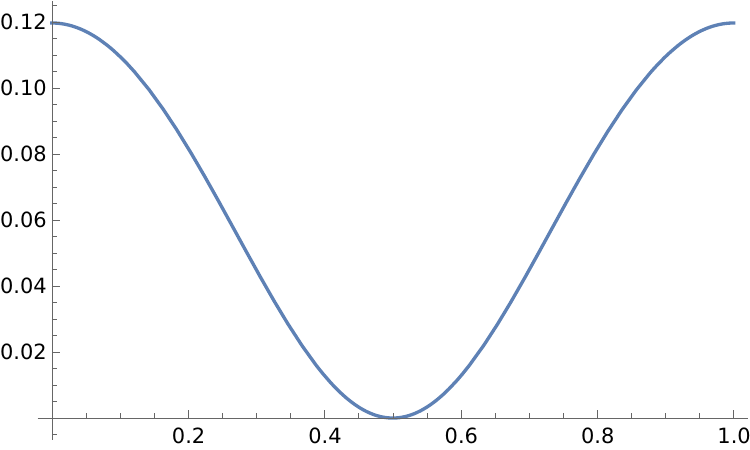}$\,\;\qquad\qquad$\includegraphics[height=3.6cm]{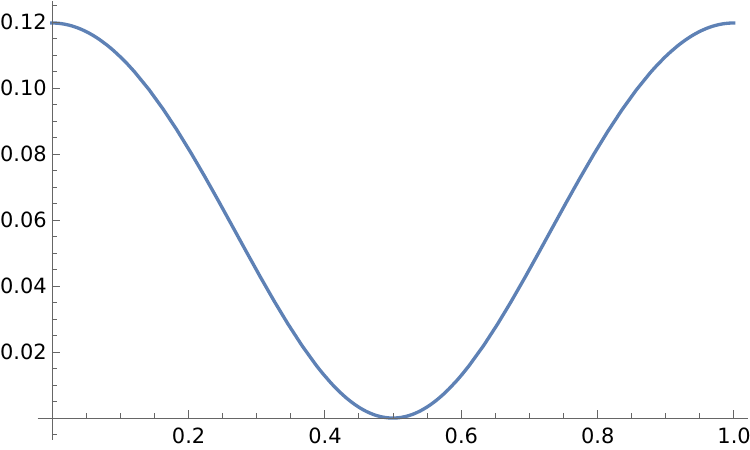}\centering
\caption{Left: sketch of the graph of the function~$\displaystyle
1+\frac{\big(5\pi-1\big)\Big(5\pi-\sqrt{25\pi^2-1}\Big)}{25\pi^2-1-5\pi\sqrt{25\pi^2-1}}
+\frac{\Big(1 - 5 \pi + \sqrt{25 \pi^2 - 1}\Big)^2 
}{\sqrt{25 \pi^2 - 1}}\,\sum_{k=1}^{3} \Big(\sqrt{25 \pi^2 - 1} - 5 \pi\Big)^{k-1}\,\cos(2\pi kx)$. Right:
sketch of the graph of the function~$\displaystyle \frac{\cos(2\pi x)+1}{\cos(2\pi x)+5\pi}$.}\label{CADna23792004drt5-23-5283012P}
\end{figure}

Thus, recalling~\eqref{jasmx23er}, the corresponding Fourier Series in trigonometric form is
$$ 1+\frac{\big(5\pi-1\big)\Big(5\pi-\sqrt{25\pi^2-1}\Big)}{25\pi^2-1-5\pi\sqrt{25\pi^2-1}}
+\frac{\Big(1 - 5 \pi + \sqrt{25 \pi^2 - 1}\Big)^2 
}{\sqrt{25 \pi^2 - 1}}\,\sum_{k=1}^{+\infty} \Big(\sqrt{25 \pi^2 - 1} - 5 \pi\Big)^{k-1}\,\cos(2\pi kx),$$
see Figure~\ref{CADna23792004drt5-23-5283012P}.

\paragraph{Solution to Exercise~\ref{AKxz3eMS-90394585865}.} The requested Fourier Series is
$$\frac13 \sum_{k\in\Z} \frac1{2^{|k|}}\, e^{2\pi i kx}.$$
Indeed, we can see~$f$ as a complex function, by defining
$$ g(z):=\frac{1}{5-2(z+z^{-1})}$$
and noticing that
\begin{eqnarray*}
f(x)=\frac{1}{5-2(e^{2\pi ix}+e^{-2\pi i x})}=g(e^{2\pi ix}).
\end{eqnarray*}
In this way, by a Taylor expansion, for all~$z\in\C$ with~$1/2<|z|<2$, we see that
\begin{eqnarray*}
g(z)&=&\frac13\left(\frac{1}{1-(z/2)}+\frac{z^{-1}}{2(1-(z^{-1}/2))}\right)\\
&=&\frac13\left(
\sum_{k=0}^{+\infty}\left(\frac{z}2\right)^k+\frac{z^{-1}}{2}
\sum_{k=0}^{+\infty}\left(\frac{z^{-1}}2\right)^k
\right)\\&=&\frac13\left(
\sum_{k=0}^{+\infty}\left(\frac{z}2\right)^k+\sum_{k=0}^{+\infty}\left(\frac{z^{-1}}2\right)^{k+1}
\right)\\&=&\frac13\left(
\sum_{k=0}^{+\infty}2^{-k}z^k+\sum_{k=0}^{+\infty} 2^{-(k+1)} z^{-(k+1)}
\right)\\&=&\frac13\left(
\sum_{j\ge0}2^{-j}z^j+\sum_{j\le-1} 2^{j} z^{j}
\right)\\&=&\frac13\sum_{j\in\Z}2^{-|j|}z^j
\end{eqnarray*}
and consequently
\begin{equation} \label{9902rj32-24ru9bvmbmm} f(x)=g(e^{2\pi ix})=\frac13\sum_{j\in\Z}2^{-|j|} e^{2\pi ijx}.\end{equation}
Now, 
$$ \sum_{j\in\Z}2^{-|j|}<+\infty,$$
therefore the series on the right-hand side of~\eqref{9902rj32-24ru9bvmbmm}
is precisely the Fourier Series of~$f$ (see Exercise~\ref{LADEBD}).

\begin{figure}[h]
\includegraphics[height=3.8cm]{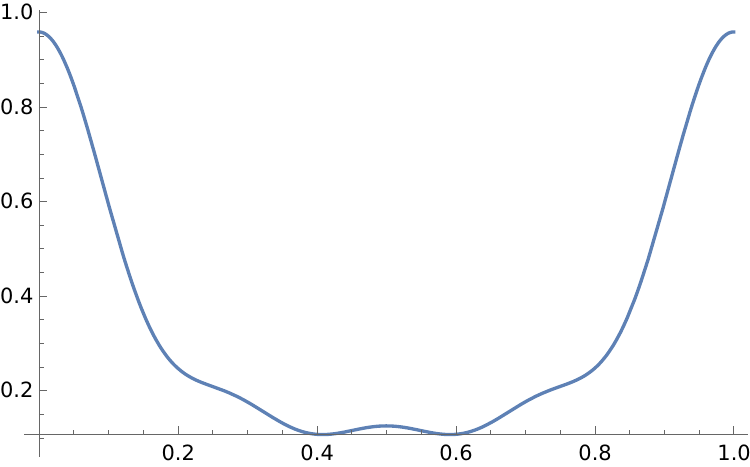}$\,\;\qquad\quad$\includegraphics[height=3.8cm]{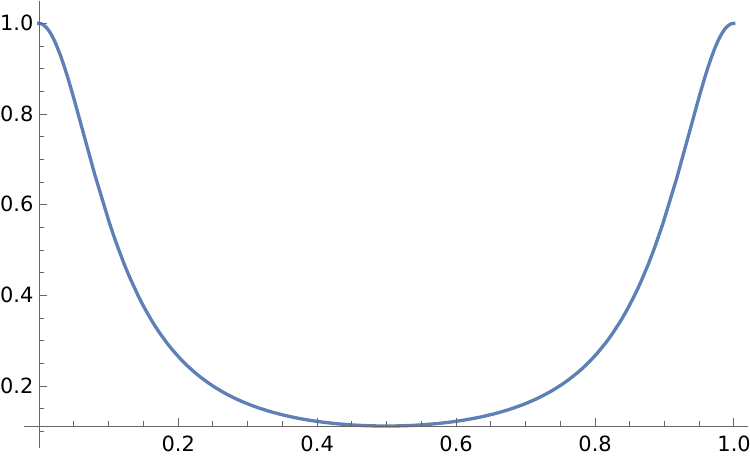}\centering
\caption{Left: plot of~$\displaystyle\frac13 \sum_{{k\in\Z}\atop{|k|\le4}} \frac1{2^{|k|}}\, e^{2\pi i kx}$.
Right: plot of~$\displaystyle\frac{1}{5-4\cos(2\pi x)}$.}\label{C7920045-23-5283012P}
\end{figure}

See Figure~\ref{C7920045-23-5283012P} for a diagram comparing the original function with a suitable trigonometric polynomial.

\paragraph{Solution to Exercise~\ref{ASKSMc-1}.} By elementary trigonometry,
\begin{equation}\label{Aspdwfe44} \cos^2 x=\frac12 +\frac12\,\cos(2 x)\end{equation}
and the latter is a Fourier series (in fact, a trigonometric polynomial) in trigonometric form. Then, the Fourier coefficients of the two functions on both sides of~\eqref{Aspdwfe44} must agree (see Exercise~\ref{LADEBD}) and therefore the right-hand side of~\eqref{Aspdwfe44} represents the Fourier Series in trigonometric form of the function on the left-hand side of~\eqref{Aspdwfe44}.

\paragraph{Solution to Exercise~\ref{FERR-2}.} First assume that
\begin{equation}\label{ILPO}
{\mbox{$f(x)>0$ for all $x\in\R$.}}
\end{equation}
We write
$$ f(x)=\sum_{{k\in\Z}\atop{|k|\le N}}c_k \,e^{2\pi ikx}$$
and we assume, without loss of generality, that 
\begin{equation}\label{2134SXD-2}
{\mbox{either~$c_N\ne0$ or~$c_{-N}\ne0$.}}\end{equation}
Since~$f$ takes values in the reals, by~\eqref{fasv} and Exercise~\ref{PKS0-3-21}
we know that
\begin{equation}\label{cnsaodp} c_{-k}=\widehat f_{-k}=\overline{\widehat f_k}=\overline{c_k}.\end{equation}
This and~\eqref{2134SXD-2} give that
\begin{equation}\label{NORO}
c_N\ne0\qquad{\mbox{and}}\qquad c_{-N}\ne0.\end{equation}

We also define
$$ P(z):=\sum_{{k\in\Z}\atop{|k|\le N}}c_k \,z^{N+k}.$$
We notice that~$P$ is a complex polynomial of degree~$2N$ (and not of a lower degree, thanks to~\eqref{NORO})
and that, by~\eqref{NORO},
\begin{equation}\label{cnsaodp2.0} P(0)=c_{-N}\ne0.\end{equation}
We also remark that 
\begin{equation}\label{cnsaodp2} P(z)\ne0 {\mbox{ for every $z\in\C$ with~$|z|=1$,}}\end{equation}
because if~$|z|=1$ we can find~$x_z\in\R$ for which~$z=e^{2\pi ix_z}$ and then
$$ P(z)=z^N\,\sum_{{k\in\Z}\atop{|k|\le N}}c_k \,z^{k}=z^N\,\sum_{{k\in\Z}\atop{|k|\le N}}c_k \,e^{2\pi ikx_z}=z^N\,f(x_z)\ne0,$$
owing to~\eqref{ILPO}.

Also, using~\eqref{cnsaodp},
\begin{equation}\label{cnsaodp3}\begin{split}&
z^{2N}\; \overline{P\left(\frac1{\overline z}\right)} -P(z)=z^{2N}\;\overline{\;\;\sum_{{k\in\Z}\atop{|k|\le N}}c_k \,\frac1{\overline{z}^{N+k}}\;\;}-
\sum_{{k\in\Z}\atop{|k|\le N}}c_k \,z^{N+k}\\&\qquad=z^{2N}\;\sum_{{k\in\Z}\atop{|k|\le N}}\overline{c_k} \,\frac{\overline{z}^{N+k}}{|z|^{2(N+k)}} -
\sum_{{k\in\Z}\atop{|k|\le N}}c_k \,z^{N+k}\\&\qquad=z^{2N}\;\sum_{{k\in\Z}\atop{|k|\le N}}{c_{-k}} \,\frac{\overline{z}^{N+k}}{|z|^{2(N+k)}} -
\sum_{{k\in\Z}\atop{|k|\le N}}c_k \,z^{N+k}\\&\qquad=z^{2N}\;\sum_{{k\in\Z}\atop{|k|\le N}}{c_k} \,\frac{\overline{z}^{N-k}}{|z|^{2(N-k)}} -
\sum_{{k\in\Z}\atop{|k|\le N}}c_k \,z^{N+k}\\&\qquad=0.
\end{split}\end{equation}
As a consequence, if~$z_\star\in\C$ is a root of~$P$, then~$z_\star\ne0$, due to~\eqref{cnsaodp2.0}, $|z_\star|\ne1$, due to~\eqref{cnsaodp2},
and also~$\frac1{\overline{z_\star}}$ is a root of~$P$, due to~\eqref{cnsaodp3}
(and we stress that~$z_\star\ne\frac1{\overline{z_\star}}$ because~$|z_\star|\ne1$).

Accordingly, we can factorise~$P$ through the Fundamental Theorem of Algebra (see Exercise~\ref{SnsdFrhndiYUJS-2}) and conclude that
$$ P(z)= c\prod_{j=1}^N \left(z-z_j\right)\left(z-\frac{1}{\overline{z_j}}\right),$$
for some~$c\in\C$, $z_1,\dots,z_N\in\C\setminus\{0\}$ (not necessarily distinct) with~$|z_j|<1$ for all~$j\in\{1,\dots,N\}$.

Therefore, setting
$$ C:=c(-1)^N\prod_{j=1}^N\frac1{\overline z_j},$$
we gather that, for all~$z\in\C\setminus\{0\}$,
\begin{eqnarray*}&&
\sum_{{k\in\Z}\atop{|k|\le N}}c_k \,z^{k}=cz^{-N}\prod_{j=1}^N \left(z-z_j\right)\left(z-\frac{1}{\overline{z_j}}\right)=C(-1)^Nz^{-N}\prod_{j=1}^N \left(z-z_j\right)\left(\overline{z_j}z-1\right)\\
&&\qquad\qquad\qquad\qquad=C\prod_{j=1}^N \left(z-z_j\right)\left(\frac1z-\overline{z_j}\right).
\end{eqnarray*}

Hence, we write~$z=e^{2\pi ix}$ and~$z_j=r_j e^{2\pi ix_j}$ for some~$x_j\in\R$ and~$r_j\in(0,1)$, and we see that
\begin{equation*}
\begin{split}
f(x)&=C\prod_{j=1}^N \left(e^{2\pi ix}-r_j e^{2\pi ix_j}\right)\left(e^{-2\pi ix}-r_j e^{-2\pi ix_j}\right)\\&
=C\prod_{j=1}^N \left(e^{2\pi ix}-r_j e^{2\pi ix_j}\right)\overline{\left(e^{2\pi ix}-r_j e^{\pi ix_j}\right)}\\&
=C\prod_{j=1}^N \left|e^{2\pi ix}-r_j e^{2\pi ix_j}\right|^2.\end{split}
\end{equation*}
By~\eqref{ILPO}, we have that~$C>0$, and hence we can define
$$ g(x):=\sqrt{C}\,\prod_{j=1}^N \left(e^{2\pi ix}-r_j e^{2\pi ix_j}\right),$$
we have that~$g$ is a trigonometric polynomial of degree~$N$ and that~$f=|g|^2$.
This proves the desired result under the additional assumption~\eqref{ILPO}.

Let us now deal with the general case of a non-negative~$f$. For this, we pick~$\epsilon\in(0,1)$, and use the result for~$f_\epsilon:=f+\epsilon$, which now satisfies~\eqref{ILPO}. This gives that there exists a trigonometric polynomial~$g_\epsilon$ such that~$f_\epsilon=|g_\epsilon|^2$. Notice that the degree of~$g_\epsilon$ does not exceed the degree of~$f$.

In particular,
$$ \sup_{x\in\R}|g_\epsilon(x)|\le\sqrt{\sup_{x\in\R}|f_\epsilon(x)|}\le\sqrt{\sup_{x\in\R}|f(x)|+1}.$$
Since the space of trigonometric polynomials of a certain degree is finite-dimensional (a trigonometric
polynomial of degree~$N$ being identified by~$2N+1$ complex coefficients), this bound gives
that, by compactness, up to a sub-sequence, $g_\epsilon$ converge to some trigonometric polynomial~$g$ of degree~$N$.

Consequently,
$$ f=\lim_{\epsilon\searrow0} f_\epsilon=\lim_{\epsilon\searrow0} |g_\epsilon|^2=|g|^2,$$
as desired.

\section{Solutions to selected exercises of Section~\ref{SFOL2}}

\paragraph{Solution to Exercise~\ref{SBPF2}.}
For further applications of the method proposed in this exercise, see e.g. Abel's Theorem in~\cite[Theorem~3 in Section~2.5]{MR4506522}.
  
Let~$x\in(0,1)$, $\rho\in[0,1]$ and \begin{equation}\label{DEBERHO}\beta_k(\rho):=
\sum_{1\le j\le k-1} (\rho e^{2\pi ix})^j=
\frac{\rho e^{2 i \pi x}-\rho^k e^{2 i \pi k x} }{1-\rho e^{2 i \pi x}  }.
\end{equation}
We recall the ``summation by parts'' formula (see Exercise~\ref{SBPF}, used here with~$\alpha_k:=\frac1k$ and~$\beta_k:=\beta_k(\rho)$) to find that
\begin{equation}\label{BESU}\begin{split}&\sum_{k=1}^{N}\frac{\rho^k e^{2\pi ikx}}k=
\sum_{k=1}^{N}\frac1k\big(\beta_{k+1}(\rho)-\beta_{k}(\rho)\big)\\&\qquad=\left(\frac1{N}\beta_{N+1}(\rho)-\beta_{1}(\rho)\right)-\sum_{k=2}^{N}\beta_{k}(\rho)\left(\frac1k-\frac1{k-1}\right)\\&\qquad=
\frac{\beta_{N+1}(\rho)}{N}+\sum_{k=2}^{N}\frac{\beta_{k}(\rho)}{k(k-1)}.
\end{split}\end{equation}

Moreover,
\begin{equation} \label{MU1SOSD}\mu:=\inf_{\rho\in[0,1]}|1-\rho e^{2 i \pi x}|>0.\end{equation}
Indeed, otherwise there would exist a sequence~$\rho_k$
for which~$|1-\rho_k e^{2 i \pi x}|\to0$ as~$k\to+\infty$. Up to a sub-sequence, $\rho_k$
converges to some~$\bar\rho\in[0,1]$ as~$k\to+\infty$, therefore~$|1-\bar\rho e^{2 i \pi x}|=0$.
This gives necessarily that~$\bar\rho=1$ and~$ x\in\Z$, against our assumption,
thus proving~\eqref{MU1SOSD}.

From~\eqref{MU1SOSD}, we find that
\begin{equation}\label{ojskds9LO-3} \beta_k(\rho)\le
\frac{2 }{|1-\rho e^{2 i \pi x}  |}\le\frac2\mu,\end{equation}
which is bounded uniformly in~$k$.

As a result, we can pass to the limit as~$N\to+\infty$ in~\eqref{BESU} and find that, for every~$x\in(0,1)$ and~$\rho\in[0,1]$,
\begin{equation}\label{ojskds9LO1}
\sum_{k=1}^{+\infty}\frac{\rho^k e^{2\pi ikx}}k=\sum_{k=2}^{+\infty}\frac{\beta_{k}(\rho)}{k(k-1)}.
\end{equation}

Now we recall the logarithmic power series, valid for all~$z\in\C$ with~$|z|<1$,
\begin{equation}\label{ojskds9LO}
-\ln(1-z)=\sum_{k=1}^{+\infty}\frac{z^k}k,\end{equation}
see e.g.~\cite[formula~(29) on page~101]{MR4577812}. Here, the notation~$\ln$ stands for the
principal branch of the complex logarithm (i.e. the branch of~$\ln(1+z)$ which is equal to~$0$ when~$z=0$).

\begin{figure}[h]
\includegraphics[height=3.8cm]{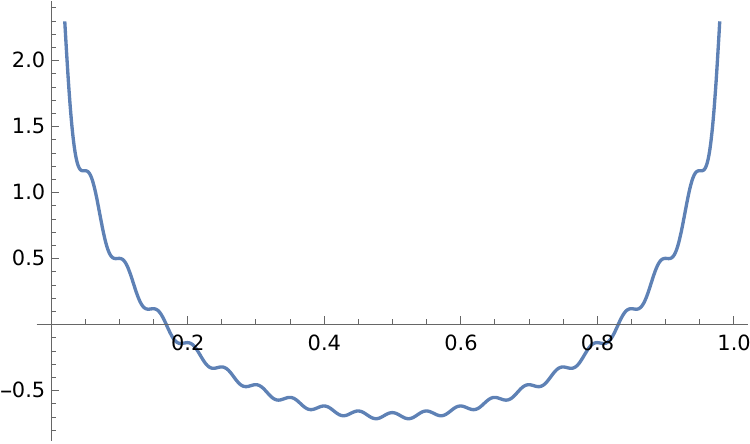}$\,\;\qquad\quad$\includegraphics[height=3.8cm]{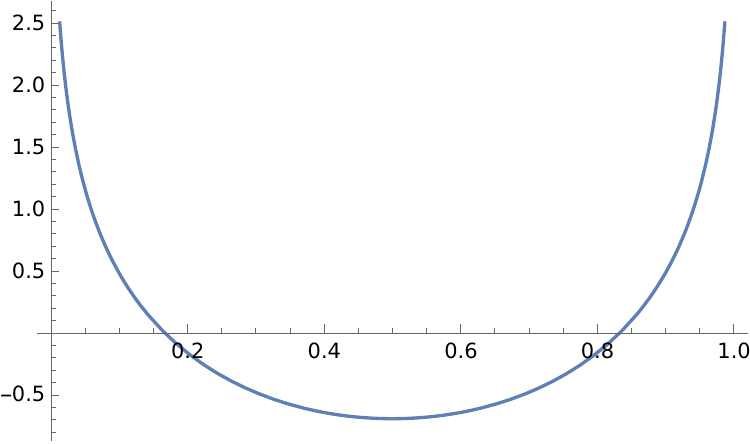}\centering
\caption{Left: plot of~$\displaystyle\sum_{k=1}^{20}\frac{\cos(2\pi kx)}k$.
Right: plot of~$\displaystyle\frac12\Big(\ln 2+\ln(1-\cos(2\pi x))\Big)$.}\label{C7283012P}
\end{figure}

In particular, we can use~\eqref{ojskds9LO} for~$z:=\rho e^{2\pi ix}$ when~$\rho\in[0,1)$
and deduce from~\eqref{ojskds9LO-3} and~\eqref{ojskds9LO1} that, for every~$K\in\N$, to be taken as large as we wish in what follows,
\begin{eqnarray*}&&
\left|- \ln (1-e^{2\pi ix})-\sum_{k=1}^{+\infty}\frac{e^{2\pi ikx}}k\right|=
\lim_{\rho\nearrow1}\left|-\ln (1-\rho e^{2\pi ix})-\sum_{k=1}^{+\infty}\frac{e^{2\pi ikx}}k\right|\\
&&\qquad=
\lim_{\rho\nearrow1}\left| \sum_{k=1}^{+\infty}\frac{\rho^k e^{2\pi ikx}}k-\sum_{k=1}^{+\infty}\frac{e^{2\pi ikx}}k\right|=\lim_{\rho\nearrow1}\left|
\sum_{k=2}^{+\infty}\frac{\beta_{k}(\rho)}{k(k-1)}-\sum_{k=2}^{+\infty}\frac{\beta_{k}(1)}{k(k-1)}
\right|\\&&\qquad\le
\lim_{\rho\nearrow1}\left[ 
\sum_{k=2}^{K}\frac{\big|\beta_{k}(\rho)-\beta_k(1)\big|}{k(k-1)}
+
\sum_{k=K+1}^{+\infty}\frac{4}{\mu\,k(k-1)}\right]\\&&\qquad=\frac4\mu\sum_{k=K+1}^{+\infty}\frac{1}{k(k-1)}.
\end{eqnarray*}
The latter is the tail of a convergent sequence, hence we can send~$K\to+\infty$ and obtain that
\begin{equation} \label{pksldmc-4}
\sum_{k=1}^{+\infty}\frac{e^{2\pi ikx}}k=-\ln (1-e^{2\pi ix}).\end{equation}
Interestingly, this says that we can somewhat ``pass the logarithmic expansion~\eqref{ojskds9LO}'' to the boundary of the disk.

We now write~$1-e^{2\pi ix}$ in polar form.
Namely, we have that~$1-e^{2\pi ix}=r(x) e^{i\theta(x)}$, with
\begin{eqnarray*}&& r(x):=\sqrt{2(1-\cos(2\pi x))}\\
{\mbox{and}}\qquad&&\tan\theta(x)=-\frac{\sin(2\pi x)}{1-\cos(2\pi x)}=-\frac1{\tan(\pi x)}=\tan\left(\pi\left(x-\frac12\right)\right).\end{eqnarray*}
As a consequence, for all~$x\in(0,1)$,
$$ \ln (1-e^{2\pi ix})=\ln r(x)+i\pi \left(x-\frac12\right),$$
where the addendum~$\frac12$ comes from the fact that we took the principal branch of the logarithm
(hence the imaginary part of~$\ln (1-e^{2\pi ix})$ must vanish when~$1-e^{2\pi ix}=2$).

\begin{figure}[h]
\includegraphics[height=3.8cm]{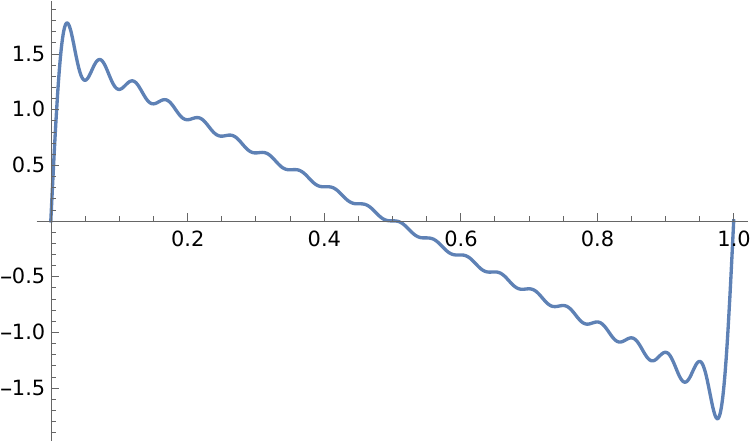}$\,\;\qquad\quad$\includegraphics[height=3.8cm]{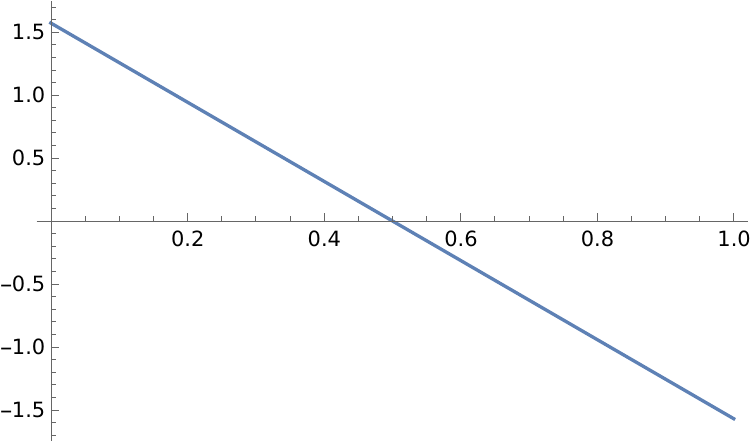}\centering
\caption{Left: plot of~$\displaystyle\sum_{k=1}^{20}\frac{\sin(2\pi kx)}k$.
Right: plot of~$\displaystyle\pi\left(\frac12- x\right)$.}\label{S7283012P}
\end{figure}

Therefore, taking the real and imaginary parts of~\eqref{pksldmc-4},
\begin{eqnarray*}&&\sum_{k=1}^{+\infty}\frac{\cos(2\pi kx)}k=
\Re\left(\sum_{k=1}^{+\infty}\frac{e^{2\pi ikx}}k\right)=-\Re\big(\ln (1-e^{2\pi ix})\big)=-\ln r(x)\\&&\qquad=-\frac12\Big(\ln 2+\ln(1-\cos(2\pi x))\Big)
\end{eqnarray*}
and
\begin{eqnarray*}&&\sum_{k=1}^{+\infty}\frac{\sin(2\pi kx)}k=\Im\left(
\sum_{k=1}^{+\infty}\frac{e^{2\pi ikx}}k\right)=-\Im\big(\ln (1-e^{2\pi ix})\big)=\pi\left(\frac12- x\right),
\end{eqnarray*}
yielding the desired result.

See Figures~\ref{C7283012P} and~\ref{S7283012P} for a visualisation of these identities.

\paragraph{Solution to Exercise~\ref{NOCAP}.}
By~\eqref{GL4}, if~$k\ne j\in\Z$,
\begin{eqnarray*}&&
\| e^{2\pi ikx}-e^{2\pi ijx}\|_{L^2((0,1))}^2=2,
\end{eqnarray*}
hence the sequence~$\{e^{2\pi ikx}\}_{k\in\Z}$ belongs to the closed unit ball of~$L^2((0,1))$ but none of its sub-sequence converges in~$L^2((0,1))$.

\paragraph{Solution to Exercise~\ref{LEGEPO}.} First off, we observe that~$\left(\frac{d}{dx}\right)^k(x^2-1)^k$ is a polynomial of degree~$k$ (since one takes~$k$ derivatives of a polynomial of degree~$2k$).

Now, integrating by parts~$k$ times,
\begin{eqnarray*}
\int_{-1}^1 |P_k(x)|^2\,dx&=&\frac{2k+1}{2^{2k+1}\,(k!)^2}\int_{-1}^1\left|\left(\frac{d}{dx}\right)^k(x^2-1)^k\right|^2\,dx\\
&=&\frac{(2k+1)\,(-1)^k}{2^{2k+1}\,(k!)^2}\int_{-1}^1(x^2-1)^k\,\left(\frac{d}{dx}\right)^{2k}(x^2-1)^k\,dx\\
&=&\frac{(2k+1)\,(-1)^k\,(2k)!}{2^{2k+1}\,(k!)^2}\int_{-1}^1(x^2-1)^k\,dx\\&=&\frac{(2k+1)!}{2^{2k}\,(k!)^2}\int_{0}^1(1-x^2)^k\,dx.
\end{eqnarray*}
Therefore, to check that
$$ \int_{-1}^1 |P_k(x)|^2\,dx=1,$$
it suffices to establish that, for all~$k\in\N$,
\begin{equation}\label{0oqlkewf-0.012l3xeX}
\int_{0}^1(1-x^2)^k\,dx=\frac{2^{2k}\,(k!)^2}{(2k+1)!}.
\end{equation}
We prove this by induction. The claim in~\eqref{0oqlkewf-0.012l3xeX} is obviously true for~$k=0$, hence we assume it for the index~$k-1$ in order to prove it for the index~$k$. To this end, we integrate by parts to see that
\begin{eqnarray*}
&& I_k:=\int_{0}^1(1-x^2)^k\,dx=\int_{0}^1(1-x^2)(1-x^2)^{k-1}\,dx\\&&\qquad=
\int_{0}^1(1-x^2)^{k-1}\,dx-\int_{0}^1x^2(1-x^2)^{k-1}\,dx
\\&&\qquad=I_{k-1}+\frac1{2k}\int_{0}^1x\frac{d}{dx}(1-x^2)^{k}\,dx
\\&&\qquad=I_{k-1}-\frac1{2k}\int_{0}^1 (1-x^2)^{k}\,dx\\&&\qquad=I_{k-1}-\frac{I_k}{2k},
\end{eqnarray*}
yielding that
$$I_k=\frac{2k}{2k+1}I_{k-1}.$$
Using the inductive assumption, we thus find that
\begin{eqnarray*}&&
I_k=\frac{2k}{2k+1}\cdot\frac{2^{2k-2}\,((k-1)!)^2}{(2k-1)!}=\frac{2^{2k}\,(k!)^2}{(2k+1)!},
\end{eqnarray*}completing the proof of~\eqref{0oqlkewf-0.012l3xeX}.

Now we take~$k$, $m\in\N$ with~$k>m$.
To prove that
$$ \int_{-1}^1 P_k(x)\,P_m(x)\,dx=0$$
we integrate by parts~$m+1$ times in the following calculation:
\begin{eqnarray*}&&
\int_{-1}^1 \left(\frac{d}{dx}\right)^k(x^2-1)^k\,\left(\frac{d}{dx}\right)^m(x^2-1)^m\,dx\\&&\qquad=
(-1)^{m+1}\int_{-1}^1 \left(\frac{d}{dx}\right)^{k-m-1}(x^2-1)^k\,\left(\frac{d}{dx}\right)^{2m+1}(x^2-1)^m\,dx=0
,\end{eqnarray*}
which gives the desired result.

\paragraph{Solution to Exercise~\ref{HAARBA}.}
Let us first deal with the case~$n_1=n_2=:n$ and~$k_1=k_2=:k$. In this case, the quantity in the left-hand side of~\eqref{po09iuyhglsqd-2.1} boils down to
$$ \int_{\frac{k-1}{2^n}}^{\frac{k}{2^n}} 2^n\,dx=1,$$
as desired.

Let us now suppose that~$n_1\ne n_2$. Up to swapping these two indices, we may suppose that~$n_2> n_1$.
Then, the the quantity in the left-hand side of~\eqref{po09iuyhglsqd-2.1} becomes
\begin{equation}\label{TCJHVLCDERSSRICMDJ-1}
2^{n_1/2}\left(\int_{\frac{k_1-1}{2^{n_1}}}^{\frac{k_1-\frac12}{2^{n_1}}} h_{k_2,n_2}(x)\,dx-
\int_{\frac{k_1-\frac12}{2^{n_1}}}^{\frac{k_1}{2^{n_1}}}\,h_{k_2,n_2}(x)\,dx\right).
\end{equation}
We observe that both these terms are zero, namely
\begin{equation}\label{TCJHVLCDERSSRICMDJ}
\int_{\frac{k_1-1}{2^{n_1}}}^{\frac{k_1-\frac12}{2^{n_1}}} h_{k_2,n_2}(x)\,dx=0\qquad{\mbox{and}}\qquad
\int_{\frac{k_1-\frac12}{2^{n_1}}}^{\frac{k_1}{2^{n_1}}}\,h_{k_2,n_2}(x)\,dx.
\end{equation}
Indeed, suppose not: say for instance that
\begin{equation}\label{lXa-203.112} \int_{\frac{k_1-1}{2^{n_1}}}^{\frac{k_1-\frac12}{2^{n_1}}} h_{k_2,n_2}(x)\,dx\ne0.\end{equation}
Then necessarily~$\left(\frac{k_2-1}{2^{n_2}},\frac{k_2}{2^{n_2}}\right)\cap\left({\frac{k_1-1}{2^{n_1}}},{\frac{k_1-\frac12}{2^{n_1}}}\right)\ne\varnothing$. 

Therefore,
$$ \frac{k_2}{2^{n_2}}>{\frac{k_1-1}{2^{n_1}}}\qquad{\mbox{and}}\qquad\frac{k_2-1}{2^{n_2}}<{\frac{k_1-\frac12}{2^{n_1}}}.$$
Hence,
$$ k_2>2^{n_2-n_1}(k_1-1)\in\Z\qquad{\mbox{and}}\qquad k_2-1<2^{n_2-n_1}\left(k_1-\frac12\right)\in\Z,
$$
yielding that
$$k_2\ge 1+2^{n_2-n_1}(k_1-1) \qquad{\mbox{and}}\qquad k_2\le2^{n_2-n_1}\left(k_1-\frac12\right).$$

As a result,
\begin{eqnarray*}
\frac{k_2-1}{2^{n_2}}\ge\frac{1+2^{n_2-n_1}(k_1-1)-1}{2^{n_2}}={\frac{k_1-1}{2^{n_1}}}
\end{eqnarray*}
and
\begin{eqnarray*}
\frac{k_2}{2^{n_2}}\le\frac{2^{n_2-n_1}\left(k_1-\frac12\right)}{2^{n_2}}={\frac{k_1-\frac12}{2^{n_1}}}.
\end{eqnarray*}

We thereby conclude that
$$ \left(\frac{k_2-1}{2^{n_2}},\frac{k_2}{2^{n_2}}\right)\subseteq\left({\frac{k_1-1}{2^{n_1}}},{\frac{k_1-\frac12}{2^{n_1}}}\right)$$
and accordingly~\eqref{lXa-203.112} becomes
$$ 0\ne\int_{\frac{k_2-1}{2^{n_2}}}^{\frac{k_2}{2^{n_2}}} h_{k_2,n_2}(x)\,dx=
2^{n_2/2}\left(-
\int_{\frac{k_2-1}{2^{n_2}}}^{\frac{k_2-\frac12}{2^{n_2}}} \,dx
+
\int_{\frac{k_2-\frac12}{2^{n_2}}}^{\frac{k_2}{2^{n_2}}} \,dx\right)=0
.$$
This is a contradiction and we have thus completed the proof of~\eqref{TCJHVLCDERSSRICMDJ}.

From~\eqref{TCJHVLCDERSSRICMDJ-1} and~\eqref{TCJHVLCDERSSRICMDJ} we conclude that the quantity in the left-hand side of~\eqref{po09iuyhglsqd-2.1} vanishes whenever~$n_1\ne n_2$. 

It remains now to consider the case in which~$n_1=n_2=:n$ and~$k_1\ne k_2$.
But then
$$ \left({\frac{k_1-1}{2^{n}}},{\frac{k_1}{2^{n}}}\right)\cap\left(\frac{k_2-1}{2^{n}},\frac{k_2}{2^{n}}\right)
=\varnothing$$
and, as a consequence, the quantity in the left-hand side of~\eqref{po09iuyhglsqd-2.1} must vanish.

\paragraph{Solution to Exercise~\ref{CLA:AUNIFO}.}
We use the notation in~\eqref{DEBERHO}, with~$\beta_k:=\beta_k(1)$.
We know from~\eqref{BESU} (used here with~$\rho:=1$) that
\begin{eqnarray*}
\sum_{k=1}^{N}\frac{\cos(2\pi kx)}k+i\sum_{k=1}^{N}\frac{\sin( 2\pi kx)}k=
\sum_{k=1}^{N}\frac{\rho^k e^{2\pi ikx}}k=
\frac{\beta_{N+1}}{N}+\sum_{k=2}^{N}\frac{\beta_{k}}{k(k-1)}.
\end{eqnarray*}
Thus, since
$$ \sup_{x\in[a,1-a]}|\beta_k|\le\sup_{x\in[a,1-a]}\frac2{|1-e^{2\pi ix}|},$$
the desired uniform converge would follow once we show that
\begin{equation}\label{09-02wd} \inf_{x\in[a,1-a]} |1-e^{2\pi ix}|>0.\end{equation}
So, we focus on the proof of~\eqref{09-02wd}. We argue by contradiction: if not,
by continuity, there would exist~$x_\star\in[a,1-a]$ such that~$|1-e^{2\pi ix_\star}|=0$.
But this would give that~$x_\star\in \Z$, against our assumption, and the proof of~\eqref{09-02wd}
is thereby complete.

\paragraph{Solution to Exercise~\ref{1-CLA:AUNIFO}.}
We recall the ``summation by parts'' formula (see Exercise~\ref{SBPF}, used here with~$\alpha_k:=\frac1{k\,\ln k}$ and~$\beta_k:=\beta_k(1)$ as in~\eqref{DEBERHO}) and we obtain
\begin{eqnarray*}
\sum_{k=m}^{n}\frac{e^{2\pi ikx}}{k\,\ln k}=\left(\frac{e^{2n\pi ix}}{n\,\ln n}-\frac{e^{2(m-1)\pi ix}}{m\,\ln m}\right)-\sum_{k=m+1}^{n} e^{2(k-1)\pi ix}\left(\frac1{k\,\ln k}-\frac1{(k-1)\,\ln( k-1)}\right).
\end{eqnarray*}
Recognising a telescopic sum, we thus obtain that
\begin{eqnarray*}&&
\sup_{x\in\R}\left|\sum_{k=m}^{n}\frac{e^{2\pi ikx}}{k\,\ln k}\right|\le
\frac{1}{n\,\ln n}+\frac{1}{m\,\ln m}+
\sum_{k=m+1}^{n} \left(\frac1{(k-1)\,\ln( k-1)}-\frac1{k\,\ln k}\right)\\
&&\qquad\le\frac{1}{n\,\ln n}+\frac{2}{m\,\ln m},
\end{eqnarray*}
which is as small as we wish provided that~$m$ and~$n$ are sufficiently large.

Since\begin{eqnarray*}
\left|\sum_{k=m}^{n}\frac{e^{2\pi ikx}}{k\,\ln k}\right|&=&
\left|\sum_{k=m}^{n}\frac{\cos(2\pi kx)}{k\,\ln k}
+i\sum_{k=m}^{n}\frac{\sin(2\pi kx)}{k\,\ln k}\right|\\&\ge&\max\left\{
\left|\sum_{k=m}^{n}\frac{\cos(2\pi kx)}{k\,\ln k}\right|,\,\left|\sum_{k=m}^{n}\frac{\sin(2\pi kx)}{k\,\ln k}\right|
\right\},\end{eqnarray*}
the desired result follows.

\paragraph{Solution to Exercise~\ref{MASCHE}.} We use the ``summation by parts'' formula (see Exercise~\ref{SBPF}). Specifically, we employ~\eqref{BYPA} with~$\alpha_k:=\frac1k$ and~$\beta_k:=k+1$ and we find that\begin{eqnarray*}&&\sum_{k=1}^{{\lfloor x\rfloor}}\frac1k=
\sum_{k=1}^{{\lfloor x\rfloor}}\alpha_{k}(\beta_{k+1}-\beta_{k})=\left(\alpha_{{\lfloor x\rfloor}}\beta_{{\lfloor x\rfloor}+1}-\alpha_{1}\beta_{1}\right)-\sum_{k=2}^{{\lfloor x\rfloor}}\beta_{k}(\alpha_{k}-\alpha_{k-1})\\&&\qquad=\left(\frac{{\lfloor x\rfloor}+2}{{\lfloor x\rfloor}}-2\right)-\sum_{k=2}^{{\lfloor x\rfloor}}(k+1)\left(\frac1{k}-\frac1{k-1}\right)\\&&\qquad=-\frac{{\lfloor x\rfloor} - 2}{\lfloor x\rfloor}
-\sum_{k=2}^{{\lfloor x\rfloor}}\left(\frac1k -\frac2{k - 1}\right).
\end{eqnarray*}

Moreover, for all~$x>1$,
\begin{eqnarray*}&&
\ln x+\gamma =\ln x+\int_{1}^{+\infty }\left({\frac{1}{\lfloor y\rfloor }}-{\frac{1}{y}}\right)\,dy\\&&\qquad=
\ln\frac{x}{\lfloor x\rfloor}+\ln{\lfloor x\rfloor}
+\int_{1}^{{\lfloor x\rfloor}}\left({\frac{1}{\lfloor y\rfloor }}-{\frac{1}{y}}\right)\,dy
+\int_{{\lfloor x\rfloor}}^{+\infty }\left({\frac{1}{\lfloor y\rfloor }}-{\frac{1}{y}}\right)\,dy\\&&\qquad=
\ln\frac{x}{\lfloor x\rfloor}+\int_{1}^{{\lfloor x\rfloor}}{\frac{dy}{\lfloor y\rfloor }}
+\int_{{\lfloor x\rfloor}}^{+\infty }\left({\frac{1}{\lfloor y\rfloor }}-{\frac{1}{y}}\right)\,dy\\
\\&&\qquad=
\ln\frac{x}{\lfloor x\rfloor}+\sum_{k=2}^{{\lfloor x\rfloor}}
\int_{k-1}^{k}{\frac{dy}{\lfloor y\rfloor }}
+\int_{{\lfloor x\rfloor}}^{+\infty }\left({\frac{1}{\lfloor y\rfloor }}-{\frac{1}{y}}\right)\,dy\\&&\qquad=
\ln\frac{x}{\lfloor x\rfloor}+\sum_{k=2}^{{\lfloor x\rfloor}}
\int_{k-1}^{k}\frac{dy}{k-1}
+\int_{{\lfloor x\rfloor}}^{+\infty }\left({\frac{1}{\lfloor y\rfloor }}-{\frac{1}{y}}\right)\,dy\\&&\qquad=
\ln\frac{x}{\lfloor x\rfloor}+\sum_{k=2}^{{\lfloor x\rfloor}}
\frac{1}{k-1}
+\int_{{\lfloor x\rfloor}}^{+\infty }\left({\frac{1}{\lfloor y\rfloor }}-{\frac{1}{y}}\right)\,dy.
\end{eqnarray*}

From the observations above, one deduces that
\begin{eqnarray*}
&&\lim_{x\to+\infty} \ln x+\gamma-\sum_{k=1}^{\lfloor x\rfloor}\frac1k\\
&=&\lim_{x\to+\infty} \ln\frac{x}{\lfloor x\rfloor}+\sum_{k=2}^{{\lfloor x\rfloor}}\frac1{k-1}
+\int_{{\lfloor x\rfloor}}^{+\infty }\left({\frac{1}{\lfloor y\rfloor }}-{\frac{1}{y}}\right)\,dy\\&&\qquad\qquad\qquad+\frac{{\lfloor x\rfloor} - 2}{\lfloor x\rfloor}
+\sum_{k=2}^{{\lfloor x\rfloor}}\left(\frac1k -\frac2{k - 1}\right)\\
&=&\lim_{x\to+\infty} \int_{{\lfloor x\rfloor}}^{+\infty }\left({\frac{1}{\lfloor y\rfloor }}-{\frac{1}{y}}\right)\,dy+1
+\sum_{k=2}^{{\lfloor x\rfloor}}\left(\frac1k -\frac1{k - 1}\right)\\
&=&\lim_{x\to+\infty} \int_{{\lfloor x\rfloor}}^{+\infty }\left({\frac{1}{\lfloor y\rfloor }}-{\frac{1}{y}}\right)\,dy+1+\frac1{{\lfloor x\rfloor}} -1\\&=&\lim_{x\to+\infty} \int_{{\lfloor x\rfloor}}^{+\infty }\left({\frac{1}{\lfloor y\rfloor }}-{\frac{1}{y}}\right)\,dy.
\end{eqnarray*}
Hence, to complete the proof of the desired result, it remains to show that
\begin{equation}\label{EUMACLA}
\lim_{x\to+\infty} \int_{{\lfloor x\rfloor}}^{+\infty }\left({\frac{1}{\lfloor y\rfloor }}-{\frac{1}{y}}\right)\,dy=0.\end{equation}
For this, we observe that, for all~$y\ge{\lfloor x\rfloor}$,
$$ \left|{\frac{1}{\lfloor y\rfloor }}-{\frac{1}{y}}\right|=
\left|\frac{y-\lfloor y\rfloor}{\lfloor y\rfloor \,y}\right|\le\frac{1}{\lfloor y\rfloor \,y}\le\frac1{(y-1)^2}$$
and consequently
$$\left|\int_{{\lfloor x\rfloor}}^{+\infty }\left({\frac{1}{\lfloor y\rfloor }}-{\frac{1}{y}}\right)\,dy\right|\le
\int_{{\lfloor x\rfloor}}^{+\infty }\frac1{(y-1)^2}\,dy,$$
which is infinitesimal as~$x\to+\infty$, thus establishing~\eqref{EUMACLA}.

\paragraph{Solution to Exercise~\ref{ABEDI1}.}
We use the summation by parts formula in Exercise~\ref{SBPF}, used here with
$$ \alpha_k:=\tau_k\qquad{\mbox{and}}\qquad\beta_k:=\sum_{N\le j\le k-1} \sigma_j,$$
where~$N$ will be chosen appropriately large.

In this way, we have that
$$ \beta_{k+1}-\beta_{k}=\sigma_k$$
and therefore, in view of~\eqref{BYPA}, for all~$m$, $n\in\N$ with~$n>m$,
\begin{eqnarray*} \sum_{k=m}^{n}\sigma_k\tau_{k}=\left(\tau_{n}
\sum_{N\le j\le n} \sigma_j
-\tau_{m}\sum_{N\le j\le m-1} \sigma_j\right)-\sum_{k=m+1}^{n}
\left(\sum_{N\le j\le k-1} \sigma_j\right)
(\tau_{k}-\tau_{k-1}).\end{eqnarray*}

Hence, given~$\epsilon>0$, we pick~$N$ sufficiently large such that, for all~$M\ge N$,
$$ \left| \sum_{N\le j\le M} \sigma_j\right|\le\epsilon$$
and we obtain that, if~$n$ and~$m$ are sufficiently large,
\begin{eqnarray*} \left|\sum_{k=m}^{n}\sigma_k\tau_{k}\right|\le2\epsilon\sup_{j\in\N}|\tau_{j}|+\epsilon\sum_{k=m+1}^{n}
|\tau_{k}-\tau_{k-1}|.\end{eqnarray*}
The monotonicity of~$\{\tau_k\}_{k\in\N}$ gives that the latter is a telescopic sum, hence
\begin{eqnarray*} \left|\sum_{k=m}^{n}\sigma_k\tau_{k}\right|\le4\epsilon\sup_{j\in\N}|\tau_{j}|.\end{eqnarray*}
From this, one obtains the desired result.

\paragraph{Solution to Exercise~\ref{ABEDI2}.}
We use the summation by parts formula in Exercise~\ref{SBPF}, used here with
$$ \alpha_k:=\tau_k\qquad{\mbox{and}}\qquad\beta_k:=\sum_{0\le j\le k-1} \sigma_j.$$
In this way, we have that
$$ \beta_{k+1}-\beta_{k}=\sigma_k$$
and therefore, in view of~\eqref{BYPA}, for all~$m$, $n\in\N$ with~$n>m$,
\begin{equation*} \sum_{k=m}^{n}\tau_{k}\sigma_k=\tau_{n}\left(\sum_{0\le j\le n} \sigma_j\right)
-\tau_{m}\left(\sum_{0\le j\le m-1} \sigma_j\right)-\sum_{k=m+1}^{n}
\left(\sum_{0\le j\le k-1} \sigma_j\right)
(\tau_{k}-\tau_{k-1}).\end{equation*}
Consequently,
\begin{eqnarray*}\left|
\sum_{k=m}^{n}\tau_{k}\sigma_k\right|\le M\big(|\tau_{n}|+|\tau_{m}|\big)+M\sum_{k=m+1}^{n}
|\tau_{k}-\tau_{k-1}|
\end{eqnarray*}
The monotonicity of~$\{\tau_k\}_{k\in\N}$ gives that the latter is a telescopic sum, hence
\begin{eqnarray*}\left|
\sum_{k=m}^{n}\tau_{k}\sigma_k\right|\le M\big(|\tau_{n}|+|\tau_{m}|\big)+M
|\tau_{n}-\tau_{m}|\le2 M\big(|\tau_{n}|+|\tau_{m}|\big).
\end{eqnarray*}
The fact that~$\tau_k$ is infinitesimal gives the desired result.

\section{Solutions to selected exercises of Section~\ref{DECAY:e:sFOL}}

\paragraph{Solution to Exercise~\ref{LASESQ}.}
Given~$f\in L^1((0,1))$, we define
$$ g_M(x):=\begin{cases}
f(x) & {\mbox{ if }}f(x)\in(- M,M),\\
0&{\mbox{ otherwise}}\end{cases}$$
and
$$ h_M(x):=f(x)-g_M(x)=\begin{cases}
0 & {\mbox{ if }}f(x)\in(- M,M),\\
f(x)&{\mbox{ otherwise.}}\end{cases}$$

Notice that~$|g_M(x)|\le M$, hence~$g_M\in L^\infty((0,1))$.
 
We also stress that since~$f\in L^1((0,1))$, we have that~$f$ is finite a.e. (see e.g.~\cite[Theorem~5.22]{MR3381284})
 and accordingly~$h_M\to0$ a.e. in~$(0,1)$ as~$M\to+\infty$. Moreover, $|h_M(x)|\le|f(x)|$ for all~$x\in(0,1)$. As a result, in light of the Dominated Convergence Theorem (see e.g.~\cite[Theorem~5.36]{MR3381284}),
$$ \lim_{M\to+\infty}\int_0^1 |h_M(x)|\,dx=0,$$
as desired.

\paragraph{Solution to Exercise~\ref{BEDD-pre}.}
By~\eqref{FOUCO} and~\eqref{CARGi},
\begin{eqnarray*}
|\widehat f_k|=\frac{|\widehat{D^m f}_k|}{|2\pi k|^m}=\frac{1}{|2\pi k|^m}\,
\left|\int_0^1 D^mf(x)\,e^{-2\pi i kx}\,dx\right|\le
\frac{1}{|2\pi k|^m}\int_0^1 |D^mf(x)|\,dx,
\end{eqnarray*}
giving the desired result.

\paragraph{Solution to Exercise~\ref{BEDD}.}
Let~$g:=D^mf$. Notice that~$g$ is bounded, and in particular~$g\in L^2((0,1))$. Thus, we can apply Bessel's Inequality~\eqref{AJSa} and infer that
$$\sum_{k\in\Z} |\widehat g_k|^2 \leq \|g\|^2_{L^2((0,1))}. $$
This and~\eqref{CARGi} imply the desired result.

\paragraph{Solution to Exercise~\ref{0-1-303kuHNAsq1-1}.} 
By Lemma~\ref{Le-ojqdwn23E},
we have that~$f\in C^m(\R)$ and
$$D^mf(x)=(2\pi)^m \sum_{\ell=0}^{+\infty} \frac{X(2^{\ell+1}\pi x)}{2^{\ell\alpha}}.$$
where~$X$ stands for either~$\sin$, $\cos$, $-\sin$, $-\cos$.
 
As a consequence, using Theorem~\ref{SMXC22b} and the uniform convergence of the above series,
\begin{eqnarray*}
&&\limsup_{k\to\pm\infty} |k|^{m+\alpha}\,|\widehat f_k|=\limsup_{k\to\pm\infty}
\frac{|k|^\alpha}{(2\pi)^m}\,|\widehat{D^m f}_k|\ge
\limsup_{j\to+\infty} \frac{2^{j\alpha}}{(2\pi)^m}\,\big|\widehat{D^m f}_{2^j}\big|\\&&\qquad=
\limsup_{j\to+\infty} 2^{j\alpha}\,\left|
\sum_{\ell=0}^{+\infty} \int_0^1 \frac{X(2^{\ell+1}\pi x)\,e^{-2^{j+1}\pi ix}}{2^{\ell\alpha}}\,dx
\right|=\limsup_{j\to+\infty} 2^{j\alpha}\,
\sum_{\ell=0}^{+\infty}\frac{\delta_{\ell,j}}{2^{\ell\alpha}}=1.
\end{eqnarray*}

\paragraph{Solution to Exercise~\ref{S91uejfnv10moOrEo023.4fSIAL1}}
First we prove the desired result for~$f\in C^1(\R)$ periodic of period~$1$.
This follows by integrating by parts, namely
\begin{eqnarray*}
|\widehat f_k|=\frac{1}{2\pi k}\left|\int_0^1 f(x)\,\frac{d}{dx}\big(e^{-2\pi i kx}\big)\,dx\right|=\frac{1}{2\pi k}
\left|\int_0^1 f'(x)\,e^{-2\pi i kx}\,dx\right|\le
\frac{\|f'\|_{L^\infty((0,1))}}{2\pi k},
\end{eqnarray*}which is infinitesimal as~$k\to+\infty$, as desired.

Let us now prove the desired result for the general case.
For this, let~$\epsilon>0$ and~$g=g_\epsilon\in C^1(\R)$
be a continuous function, periodic of period~$1$, such that~$\|f-g\|_{L^1((0,1))}\le\epsilon$ (see e.g.~\cite[Theorem~9.6]{MR3381284}).
Then, by Exercise~\ref{FBA},
$$ |\widehat f_k-\widehat g_k|\le\epsilon.$$
{F}rom the previous discussion, we also know that~$\widehat g_k\to0$ as~$k\to+\infty$ and therefore
$$ \lim_{k\to+\infty}|\widehat f_k|\le
\lim_{k\to+\infty}\big(|\widehat f_k-\widehat g_k|+
|\widehat g_k|\big)\le\epsilon.$$
We can now take~$\epsilon$ as small as we like and obtain the result in
Theorem~\ref{RLjoqwskcdc}.

\paragraph{Solution to Exercise~\ref{S91uejfnv10moOrEo023.4fSIAL2}}
Suppose first that, for all~$x\in[0,1]$, we have that~$f(x)=\chi_I$ for some interval~$
I$ with endpoints~$a$ and~$b$. Then,
$$ \widehat f_k=\int_a^b e^{-2\pi i kx}\,dx
=\frac{e^{-2\pi i ka}-e^{-2\pi i kb}}{2\pi ik},$$
which is infinitesimal as~$k\to+\infty$.

This is Theorem~\ref{RLjoqwskcdc} for characteristic functions of intervals. By linearity (see Exercise~\ref{LINC-bisco}) the same holds for step functions of the form
\begin{equation}\label{DFCPTNnAA}\sum_{j=1}^N c_j\chi_{I_j},\end{equation}
with~$c_j\in\R$ and~$I_j$ intervals in~$[0,1]$.

Now, for the general case, we can use the density of the step functions in~$L^1((0,1))$,
see e.g.~\cite[Theorems~2.26 and~2.41]{MR1681462}: namely, given~$\epsilon>0$,
one can find a step function~$g=g_\epsilon$ as in~\eqref{DFCPTNnAA}
such that~$\|f-g\|_{L^1((0,1))}\le\epsilon$. As a byproduct (see Exercise~\ref{FBA}),
$$|\widehat f_k-\widehat g_k|\le\epsilon.$$

Also, by the previous observation, we know that~$\widehat g_k\to0$ as~$k\to+\infty$. All in all,
$$ \lim_{k\to+\infty}|\widehat f_k|\le
\lim_{k\to+\infty}\big(|\widehat f_k-\widehat g_k|+
|\widehat g_k|\big)\le\epsilon.$$
We can now take~$\epsilon$ as small as we like and obtain the result in
Theorem~\ref{RLjoqwskcdc}.

\paragraph{Solution to Exercise~\ref{S91uejfnv10moOrEo023.4fSIAL3}}
Let~$\epsilon>0$. We observe that
$$\int_0^1 f(x+\epsilon)\,e^{-2\pi i kx}\,dx=
\int_\epsilon^{1+\epsilon} f(y)\,e^{-2\pi i k(y-\epsilon)}\,dy
=e^{2\pi i k\epsilon}
\int_0^{1} f(y)\,e^{-2\pi i ky}\,dy=e^{2\pi i k\epsilon}\,\widehat f_k,$$ thanks to Exercise~\ref{fr12},
and therefore
$$ \big(1-e^{2\pi i k\epsilon}\big)\widehat f_k=
\int_0^1 \big(f(x)-f(x+\epsilon)\big)\,e^{-2\pi i kx}\,dx.$$

We now choose~$\epsilon:=\epsilon_k:=\frac{1}{2k}$ and we find
$$ 2\widehat f_k=\big(1-e^{\pi i }\big)\widehat f_k=
\int_0^1 \big(f(x)-f(x+\epsilon_k)\big)\,e^{-2\pi i kx}\,dx.$$

Note also that~$\epsilon_k$ is infinitesimal as~$k\to+\infty$.
Consequently, by the continuity of the translations in Lebesgue spaces (see e.g.~\cite[Theorem~8.19]{MR3381284}),
$$ 2\lim_{k\to+\infty}|\widehat f_k|=\lim_{k\to+\infty}\left|
\int_0^1 \big(f(x)-f(x+\epsilon_k)\big)\,e^{-2\pi i kx}\,dx\right|
\le
\lim_{k\to+\infty}
\int_0^1 \big|f(x)-f(x+\epsilon_k)\big|\,dx=0.$$

\section{Solutions to selected exercises of Section~\ref{CONVESob}}

\paragraph{Solution to Exercise~\ref{SPDCD-0.02}.} By Exercise~\ref{SPDCD-0.01}, used here with~$\alpha:=2\pi Nx$ and~$\beta:=2\pi jx$,
$$ \frac{\cos(2\pi(N-j)x)}{j}-\frac{\cos(2\pi(N+j)x)}{j}=
\frac{2\sin(2\pi Nx)\,\sin(2\pi jx)}j $$
and accordingly
\begin{eqnarray*}
&& \sum_{{j\in\Z\setminus\{0\}}\atop{|j|\le\ell}}\frac{\cos(2\pi(N-j)x)}{j}=
\sum_{j=1}^{\ell}\frac{\cos(2\pi(N-j)x)}{j}-\sum_{j=1}^{\ell}\frac{\cos(2\pi(N+j)x)}{j}
\\&&\qquad=\sum_{j=1}^{\ell}\left(\frac{\cos(2\pi(N-j)x)}{j}-\frac{\cos(2\pi(N+j)x)}{j}\right)=\sum_{j=1}^{\ell}\frac{2\sin(2\pi Nx)\,\sin(2\pi jx)}j .\end{eqnarray*}

\section{Solutions to selected exercises of Section~\ref{CONVES}}

\paragraph{Solution to Exercise~\ref{DINI-addC-1}.} No, Dini's Condition~\eqref{DINI} rules out jump discontinuities, because, if, say, $f(x+t)\ge f(x)+a$ for all~$t\in(0,\delta_0)$ for some~$a$, $\delta_0>0$, then, setting~$\delta_1:=\min\{\delta,\delta_0\}$, we see that
$$ \int_{-\delta}^{\delta}\left|\frac{f (x + t) - f(x)}{t}\right|\,dt\ge
\int_{0}^{\delta_1}\frac{f (x + t) - f(x)}{t}\,dt\ge
\int_{0}^{\delta_1}\frac{a}{t}\,dt=+\infty.$$

\paragraph{Solution to Exercise~\ref{DINI-addC-2}.} No, $f$ is not necessarily continuous\footnote{And, similarly, $f$ does not necessarily coincide with a continuous function almost everywhere.
This also highlights that the specific value of~$f$ at the point~$x$ in Dini's Condition~\eqref{DINI} has to be taken in the pointwise sense (not messing around with sets of null Lebesgue measure): for instance, the function vanishing identically
obviously satisfies Dini's Condition~\eqref{DINI}, and
the function vanishing identically except at the origin in which takes value~$1$
obviously does not satisfy Dini's Condition~\eqref{DINI}.}
at the point~$x$ fulfilling Dini's Condition~\eqref{DINI}.

To construct an instructive example, we can proceed as follows.
By an induction argument, one observes\footnote{For completeness, we show the inductive step:
if we know that~$2^{k+1}(k+2)<(k+1)!$, then
$$ 2^{k+2}(k+3)=\frac{ 2(k+3)\,2^{k+1}(k+2)}{k+2}<
\frac{ 2(k+3)\,(k+1)!}{k+2}\le(k+2)!$$
as desired.}
that, for all~$k\in\N\cap[5,+\infty)$,
$$ 2^{k+1}(k+2)<(k+1)!$$
and therefore, the intervals~$\left(\frac1{2^k}-\frac{1}{k!},\frac1{2^k}+\frac{1}{k!} \right)$
are disjoint.

Thus, given a sequence of real numbers~$a_k\in [1,k]$, we can define
$$ \left(-\frac12,\frac12\right]\ni x\longmapsto f(x):=\begin{dcases}
a_k & {\mbox{ if $x\in \displaystyle\left(\frac1{2^k}-\frac{1}{k!},\frac1{2^k}+\frac{1}{k!} \right)$ for some~$k\in\N\cap[5,+\infty)$,}}\\
0 &{\mbox{ otherwise.}}
\end{dcases}$$
We can also extend~$f$ periodically and identify it with a function of period~$1$ and we remark that
$$ \int_{-1/2}^{1/2} |f(x)|\,dx\le\sum_{k=5}^{+\infty} \frac{2a_k}{k!}\le
\sum_{k=5}^{+\infty} \frac{2}{(k-1)!}<+\infty,$$
showing that~$f\in L^1((0,1))$.

Notice also that~$f(0)=0$ and
$$ \liminf_{k\to+\infty} f\left(\frac1{2^k}\right)\ge1$$
and consequently~$f$ is not continuous at~$0$ (not even bounded if~$a_k:=k$).

However, $f$ satisfies Dini's Condition~\eqref{DINI} at~$0$ with~$\delta:=\frac{1}4$, because
\begin{eqnarray*}&&
\int_{-1/4}^{1/4}\left|\frac{f (t) - f(0)}{t}\right|\,dt=
\int_{0}^{1/4} \frac{f (t) }{t}\,dt=
\sum_{k=5}^{+\infty} \int_{\frac1{2^k}-\frac{1}{k!}}^{\frac1{2^k}+\frac{1}{k!}}
\frac{a_k}{t}\,dt\\&&\qquad\le
\sum_{k=5}^{+\infty}
\frac{2a_k}{k!\displaystyle\left(\frac1{2^k}-\frac{1}{k!}\right)}
\le\sum_{k=5}^{+\infty}
\frac{2^{k+1}k}{k!-2^k}<+\infty.
\end{eqnarray*}

\section{Solutions to selected exercises of Section~\ref{SEC:UNIQ:3}}

\paragraph{Solution to Exercise~\ref{ojqdwn23EEE}.} Notice that the series defining~$f$ converges uniformly, hence the definition of~$f$ is well-posed. Let
$$ c_k:=\begin{dcases}\displaystyle
\frac1{2|k|^{m+1}} &{\mbox{ if }}k\ne0,\\0&{\mbox{ if }}k=0.
\end{dcases}$$
We note that
$$ \sum_{k\in\Z} c_k\,e^{2\pi ikx}=\sum_{k=1}^{+\infty} c_k\,\big(e^{2\pi ikx}+e^{-2\pi ikx}\big)=
\sum_{k=1}^{+\infty} 2c_k\,\cos(2\pi kx)=f(x).$$
Moreover,
$$ \sum_{k\in\Z}|k|^{m-1}|c_k|=\sum_{k\in\Z\setminus\{0\}}\frac1{2|k|^2}<+\infty.$$
Thus, by Lemma~\ref{Le-ojqdwn23E} (used here with~$m-1$ in the place of~$m$) we know that~$f\in C^{m-1}(\R)$ and periodic of period~$1$ (and also that~$c_k=\widehat f_k$).

But~$f\not\in C^m(\R)$. To check this, we argue by contradiction.
If~$f\in C^m(\R)$, Theorem~\ref{SMXC22b} returns that~$ \widehat{D^m f}_k=(2\pi i k)^m\widehat f_k$
and therefore
\begin{eqnarray*}&& \widehat{D^m f}_k=(2\pi i k)^mc_k=\begin{dcases}\displaystyle
\frac{(2\pi i k)^m}{2|k|^{m+1}} &{\mbox{ if }}k\ne0,\\0&{\mbox{ if }}k=0.
\end{dcases}
\end{eqnarray*}

Let us now distinguish two cases. If~$m$ is odd, we have that~$|k|^{m+1}=k^{m+1}$ and that~$i^m\in\{i,-i\}$, thus we obtain that
\begin{eqnarray*}&& \widehat{D^m f}_k=\begin{dcases}\displaystyle
\pm\frac{(2\pi)^m i}{2k} &{\mbox{ if }}k\ne0,\\0&{\mbox{ if }}k=0.
\end{dcases}
\end{eqnarray*}
That is, the Fourier coefficients of~$D^m f$ coincide (up to a multiplicative factor) with those of the sawtooth waveform in Exercise~\ref{SA:W} and consequently, by Theorem~\ref{UNIQ}, we have that~$D^m f$ is a sawtooth waveform, which
is discontinuous, against our assumptions.

Instead, if~$m$ is even, we have that~$k^m=|k|^{m}$ and that~$i^m\in\{1,-1\}$, yielding that
\begin{eqnarray*}&& \widehat{D^m f}_k=\begin{dcases}\displaystyle
\pm\frac{(2\pi)^m}{2|k|} &{\mbox{ if }}k\ne0,\\0&{\mbox{ if }}k=0.
\end{dcases}
\end{eqnarray*}
Hence, up to a multiplicative constant, the function~$D^m f$ has the same Fourier coefficients
as the function~$\phi$ in Exercise~\ref{ojqdwn23EEE0921ef-LN}.

Consequently, by Theorem~\ref{UNIQ}, we have that~$D^m f$ coincides with~$\phi$, up to a 
multiplicative constant, and therefore, by~\eqref{GIJMPPrvU},
$$ \lim_{x\to0}D^mf(x)=+\infty.$$
But this shows that~$D^mf$ is not a continuous function, against our assumptions.

\paragraph{Solution to Exercise~\ref{anjqnxn}.}
We use the Taylor expansion of the exponential, the uniform convergence of a power series and the Fubini-Tonelli Theorem to see that, for each~$k\in\Z$,
\begin{eqnarray*}&& \widehat f_k=\int_0^1 f(x)\,e^{-2\pi i kx}\,dx=
\int_0^1 f(x)\,\sum_{n=0}^{+\infty} \frac{(-2\pi i kx)^n}{n!}\,dx\\&&\qquad=\sum_{n=0}^{+\infty} \frac{(-2\pi i k)^n}{n!}
\int_0^1 f(x)\,x^n\,dx=0.
\end{eqnarray*}
Since also the Fourier coefficients of the zero function vanish, we deduce from Theorem~\ref{UNIQ}
that~$f$ and the zero function coincide a.e. in~$\R$.

\paragraph{Solution to Exercise~\ref{PKS0-3-21bismo}.} Let~$g:=S_{N,f}$.
By~\eqref{GL4}, we know that, for all~$k\in\Z$,
$$ \widehat g_k= \sum_{{h\in\Z}\atop{|h|\le N}}\int_0^1 \widehat f_k\, e^{2\pi i(h-k)x}=
\sum_{{h\in\Z}\atop{|h|\le N}} \widehat f_k\,\delta_{h,k}=\begin{dcases}
\widehat f_k&{\mbox{ if }}|k|\le N,\\ 0&{\mbox{ otherwise,}}
\end{dcases}$$
and therefore, by our assumptions, $ \widehat g_k= \widehat f_k$ for all~$k\in\Z$.

Hence, by way of Theorem~\ref{UNIQ}, up to a set of null Lebesgue measure,
$f=g=S_{N,f}$.

\paragraph{Solution to Exercise~\ref{PKS0-3-21bismo2}.}
This follows from Exercises~\ref{PKS0-3-21} and~\ref{PKS0-3-21bismo}.

\paragraph{Solution to Exercise~\ref{PKS0-3-21bismo2.NECE}.}
Let~$j_0\in\N$ be such that, whenever~$j\ge j_0$,
$$ \sup_{x\in\R}|f_j(x)-f(x)|\le1.$$
Then,
$$ \int_0^1|f(x)|\,dx\le \int_0^1|f_{j_0}(x)|\,dx+\int_0^1|f_{j_0}(x)-f(x)|\,dx\le\int_0^1|f_{j_0}(x)|\,dx+1,$$
whence~$f\in L^1((0,1))$.

Also, for all~$x\in\R$,
$$ f(x+1)-f(x)=\lim_{j\to+\infty} \big(f_j(x+1)-f_j(x) \big)=0,$$
showing that~$f$ is periodic of period~$1$.

Moreover, for all~$k\in\Z$,
\begin{eqnarray*}&& |\widehat f_k|=\left|\int_0^1 f(x)\,e^{-2\pi ikx}\,dx\right|
\le\left|\int_0^1 f_j(x)\,e^{-2\pi ikx}\,dx\right|+\left|\int_0^1 \big(f_j(x)-f(x)\big)\,e^{-2\pi ikx}\,dx\right|\\&&\qquad\qquad\qquad
\le |\widehat f_{j,k}|+\sup_{x\in\R}|f_j(x)-f(x)|.\end{eqnarray*}

From this and~\eqref{ICDCMMANHSMCGIEVIMN} we infer that, for all~$k\in\Z$ with~$|k|\ge M$, 
$$ |\widehat f_k|\le\lim_{j\to+\infty}\left(|\widehat f_{j,k}|+\sup_{x\in\R}|f_j(x)-f(x)|\right)=0.$$

This gives that the Fourier Series of~$f$ contains only finitely many terms. Accordingly,
$f$ is a trigonometric polynomial (see Exercise~\ref{PKS0-3-21bismo}).

\paragraph{Solution to Exercise~\ref{KPSLM.01e2ourfhvb0odc.2ewrfevgbx3xdnfv0-1}.} We argue by induction over~$m$.
When~$m=0$, we have that
$$  \phi_0(\sigma_0)=0\cdot\sigma_0=0$$
and therefore
$$ \sum_{{\sigma_0\in\{-1,1\}}}\frac{e^{2\pi i\phi_0 (\sigma) x}}{2^0}=1=f_0(x),$$
which is the basis of the induction.

Suppose now recursively that the claim is true for the index~$m$. Then,
\begin{eqnarray*}f_{m+1}(x)&=&\prod_{j=0}^{m+1} \cos(2\pi jx)\\&=&
f_m(x)\,\cos(2\pi(m+1)x)\\&=&
\sum_{{\sigma_0,\dots,\sigma_m\in\{-1,1\}}}
\frac{e^{2\pi i\phi_m(\sigma) x}}{2^m}
\cdot\frac{e^{2\pi i(m+1)x}+e^{-2\pi i(m+1)x}}{2}\\
&=&
\sum_{{\sigma_0,\dots,\sigma_m\in\{-1,1\}}}
\frac{e^{2\pi i\phi_m(\sigma) x}}{2^{m+1}}
\cdot\sum_{{\sigma_{m+1}\in\{-1,1\}}} e^{2\pi i\sigma_{m+1}(m+1)x}
\\&=& \sum_{{\sigma_0,\dots,\sigma_m,\sigma_{m+1}\in\{-1,1\}}}\frac{e^{2\pi i(\phi_m(\sigma)+\sigma_{m+1}(m+1)) x}}{2^{m+1}}
\\&=&\sum_{{\sigma_0,\dots,\sigma_{m+1}\in\{-1,1\}}}\frac{e^{2\pi i\phi_{m+1}(\sigma) x}}{2^{m+1}},
\end{eqnarray*}
which is the desired claim for the index~$m+1$.

\paragraph{Solution to Exercise~\ref{phimwlxc2}.}
We use the modular arithmetic, taking notice that both~$1\equiv 1$ $({\mbox{mod }}2)$ and~$-1\equiv 1$ $({\mbox{mod }}2)$. 

As a result, for every~$k\in\{0,\dots,m\}$, we have that~$\sigma_k\equiv 1$ $({\mbox{mod }}2)$.

On this account and~\eqref{phimwlxc},
$$ 0= \phi_m(\sigma)=
\sum_{k=0}^m k\sigma_k\equiv\sum_{k=0}^m k
=\frac{m(m+1)}2\qquad
({\mbox{mod }}2).$$

That being so, it follows that~$m(m+1)\equiv0$ $({\mbox{mod }}4)$.

Hence, by looking at the cases~$m\in\{0,1,2,3\}$, necessarily either~$m\equiv0$ $({\mbox{mod }}4)$ or~$m\equiv3$ $({\mbox{mod }}4)$, which is what we were asked to prove.

\paragraph{Solution to Exercise~\ref{phimwlxc3}.} Assume first that~$\frac{m}4\in\N$, i.e.~$m=4\ell$ for some~$\ell\in\N$.

Then,
\begin{eqnarray*}&&\phi_m(\sigma)=\sum_{k=0}^{4\ell} k\sigma_k=\sum_{k=1}^{4\ell} k\sigma_k
=\sum_{j=0}^{\ell-1}\; \;\sum_{k\in\{4j+1,4j+2,4j+3,4j+4\}} k\sigma_k\\&&\quad=
\sum_{j=0}^{\ell-1} \Big( (4j+1)\sigma_{4j+1}
+(4j+2)\sigma_{4j+2}+(4j+3)\sigma_{4j+3}+(4j+4)\sigma_{4j+4}\Big).\end{eqnarray*}

Hence, by virtue of~\eqref{SIKMSDgmwiedj1},
$$\phi_m(\sigma)=
\sum_{j=0}^{\ell-1} \Big( (4j+1)-
(4j+2)-(4j+3)+(4j+4)\Big)=\sum_{j=0}^{\ell-1}0=0,$$
as desired.

On a similar note, assume now that~$\frac{m-3}4\in\N$, i.e.~$m=4\ell+3$ for some~$\ell\in\N$.

Then,
\begin{eqnarray*}&&\phi_m(\sigma)=\sum_{k=0}^{4\ell+3} k\sigma_k
=\sum_{j=0}^{\ell} \;\;\sum_{k\in\{4j,4j+1,4j+2,4j+3\}} k\sigma_k\\&&\quad=
\sum_{j=0}^{\ell} \Big(4j\sigma_{4j}+
(4j+1)\sigma_{4j+1}+(4j+2)\sigma_{4j+2}+(4j+3)\sigma_{4j+3}\Big).\end{eqnarray*}

As above, on the strength of~\eqref{SIKMSDgmwiedj12},
$$\phi_m(\sigma)=\sum_{j=0}^{\ell} \Big(4j-
(4j+1)-(4j+2)+(4j+3)\Big)=\sum_{j=0}^{\ell-1}0=0,$$
as desired.

\paragraph{Solution to Exercise~\ref{LAVppnUNOMASMAJONC0989okj744124IJN}.}
The functions involved in the finite product in the right-hand side of~\eqref{LAVppnUNOMASMAJONC0989okj744124IJN.0}
are all continuous and periodic of period~$1$, as well as even, thus so is the function~$f_m$.

Thus, \eqref{LAVppnUNOMASMAJONC0989okj744124IJN.1} follows from Exercise~\ref{920-334PKSXu9o2fgfbsmos}.

Furthermore, by~\eqref{jasmx23er} and Exercise~\ref{KPSLM.01e2ourfhvb0odc.2ewrfevgbx3xdnfv0-1}, for all~$k\in\N$,
\begin{eqnarray*}
a_k&=&2\int_0^1 f_m(x)\,\cos(2\pi kx)\,dx\\
&=&\int_0^1
\sum_{{\sigma_0,\dots,\sigma_m\in\{-1,1\}}}
\frac{e^{2\pi i\phi_m(\sigma) x}}{2^{m-1}}\,\cos(2\pi kx)\,dx\\&=&\frac1{2^m}
\int_0^1
\sum_{{\sigma_0,\dots,\sigma_m\in\{-1,1\}}}
e^{2\pi i\phi_m(\sigma) x}\big( e^{2\pi ikx}+e^{-2\pi ikx}\big)\,dx\\
\\&=&\frac1{2^m} \sum_{{\sigma_0,\dots,\sigma_m,\tau\in\{-1,1\}}}
\int_0^1
e^{2\pi i (\phi_m(\sigma) +\tau k )x}\,dx.
\end{eqnarray*}

Hence, since, for every~$r\in\R$,
$$ \int_0^1 e^{2\pi i r x}\,dx=\begin{dcases} 1 &{\mbox{ if }}r=0,\\
0&{\mbox{ otherwise,}}
\end{dcases}$$
we conclude that, for all~$k\in\N$,
\begin{equation}\label{RRSDFV0kSKNEMSHPRLMCLIYAHT}\begin{split}
a_k=\frac{\nu_{m,k}}{2^m} ,\end{split}\end{equation}
where we denoted by~$\nu_{m,k}$ the number of possible choices of~$\sigma_0,\dots,\sigma_m,\tau\in\{-1,1\}$
for which~$\phi_m(\sigma) +\tau k =0$.

Since, by~\eqref{phimwlxc}, $$|\phi_m(\sigma)|\le \sum_{k=0}^m k|\sigma_k|=\sum_{k=0}^m k=\frac{m(m+1)}{2},$$
whenever~$k>\frac{m(m+1)}{2}$ it follows that
$$ |\phi_m(\sigma) +\tau k|\ge|\tau| k-| \phi_m(\sigma)|=k-| \phi_m(\sigma)|\ge k-\frac{m(m+1)}{2}>0$$
and correspondingly~$a_k=0$, which establishes~\eqref{LAVppnUNOMASMAJONC0989okj744124IJN.2}.

As a byproduct of~\eqref{LAVppnUNOMASMAJONC0989okj744124IJN.2}, we also have that the series in~\eqref{LAVppnUNOMASMAJONC0989okj744124IJN.7} and~\eqref{LAVppnUNOMASMAJONC0989okj744124IJN.3} are finite sums (say, they only possess terms up to some index~$k_0\in\N$). As a result, by Exercise~\ref{PKS0-3-21bismo} and~\eqref{LAVppnUNOMASMAJONC0989okj744124IJN.1},
 we have that, for all~$x\in\R$,
\begin{eqnarray*}f_m(x)=
\frac{a_{0}}2+\sum_{k=1}^{k_0}\Big(a_{k}\cos(2\pi kx) + b_{k}\sin(2\pi kx)\Big)
=\frac{a_{0}}2+\sum_{k=1}^{k_0} a_{k}\cos(2\pi kx),\end{eqnarray*}
from which, choosing~$x:=0$, we obtain~\eqref{LAVppnUNOMASMAJONC0989okj744124IJN.3}, as desired.

We now discuss for which values of~$m$ it holds that~$a_0\ne0$. To this end, we retake~\eqref{RRSDFV0kSKNEMSHPRLMCLIYAHT} and we see that~$a_0\ne0$ if and only if
there exists at least a choice of~$\sigma_0,\dots,\sigma_m\in\{-1,1\}$
for which~$\phi_m(\sigma)=0$.
In light of Exercises~\ref{KPSLM.01e2ourfhvb0odc.2ewrfevgbx3xdnfv0-1}
and~\ref{phimwlxc2} this occurs if and only if either~$\frac{m}4\in\N$ or~$\frac{m-3}4\in\N$.

That is, $a_0\ne0$ if and only if either~$\frac{m}4\in\N$ or~$\frac{m-3}4\in\N$.

\section{Solutions to selected exercises of Section~\ref{DEC2}}

\paragraph{Solution to Exercise~\ref{ojqdwn23EEE0921ef-LN}.}
Given the summation formula in~\eqref{SBPF2-E1},
a natural candidate for the desired function is
$$ \phi(x):=-\frac12\Big(\ln 2+\ln(1-\cos(2\pi x))\Big).$$
 We show that indeed~$\phi\in L^1((0,1))$ and that~\eqref{OSK-1-013-1} holds true (the other properties can be easily checked).
 
To prove that~$\phi\in L^1((0,1))$, it suffices to show that
\begin{equation}\label{IF1IN} I:=\int_0^1 \big| \ln(1-\cos(2\pi x))\big|\,dx<+\infty.\end{equation}
To this end, we use the substitution~$y:=1-x$ and we find that
\begin{eqnarray*}
I&=&\int_0^{1/2} \big| \ln(1-\cos(2\pi x))\big|\,dx+\int_{1/2}^1 \big| \ln(1-\cos(2\pi x))\big|\,dx\\
&=&\int_0^{1/2} \big| \ln(1-\cos(2\pi x))\big|\,dx-\int_{1/2}^0 \big| \ln(1-\cos(2\pi (1-y)))\big|\,dy\\&=&2\int_0^{1/2} \big| \ln(1-\cos(2\pi x))\big|\,dx.
\end{eqnarray*}

Moreover,
$$ \lim_{x\to0}\frac{1-\cos(2\pi x)}{x^2}=2\pi^2$$
therefore there exists~$x_0\in\left(0,\frac14\right)$ such that, for all~$x\in(0,x_0)$,
$$\frac{1-\cos(2\pi x)}{x^2}\ge\pi^2$$
and accordingly, since the logarithm is monotone increasing,
$$ -|\ln(1-\cos(2\pi x))|=
\ln(1-\cos(2\pi x))\ge \ln(\pi^2 x^2)=2\ln(\pi x).$$

This observation, together with the piecewise monotonicity of the cosine function, yields that
\begin{eqnarray*}
I&\le&-4\int_0^{x_0} \ln(\pi x)\,dx-2\int_{x_0}^{1/2} \ln(1-\cos(2\pi x))\,dx\\
&\le&-4\int_0^{x_0} \ln(\pi x)\,dx-2\int_{x_0}^{1/2} \ln(1-\cos(2\pi x_0))\,dx\\
&=&4x_0 (1 -\ln(\pi x_0))-(1-2x_0)\ln(1-\cos(2\pi x_0)),
\end{eqnarray*}
which gives~\eqref{IF1IN}.

We now check the validity of~\eqref{OSK-1-013-1}.
To this end, we know from Exercises~\ref{SBPF2} and~\ref{CLA:AUNIFO} that, for all~$a\in\left(0,\frac12\right)$,
$$ \phi(x)=\sum_{k=1}^{+\infty}\frac{\cos(2\pi kx)}k,$$
with uniform convergence when~$x\in[a,1-a]$.

Consequently, we can swap the order of integration and find that, for every~$j\in\Z$,
\begin{eqnarray*}&&
\int_a^{1-a}\phi(x)\,e^{-2\pi ijx}\,dx=\sum_{k=1}^{+\infty}\int_a^{1-a}\frac{\cos(2\pi kx)}k\,e^{-2\pi ijx}\,dx.
\end{eqnarray*}
Also, since~$\phi\in L^1((0,1))$,
$$ \lim_{a\searrow0}\int_0^a|\phi(x)|\,dx+\int_{1-a}^1|\phi(x)|\,dx=0,$$
and, as a result,
\begin{eqnarray*}
\widehat \phi_j=\lim_{a\searrow0}\int_a^{1-a}\phi(x)\,e^{-2\pi ijx}\,dx.
\end{eqnarray*}

These observations entail that
\begin{equation}\label{AKSmxcs-serf}
\widehat \phi_j=\lim_{a\searrow0}\sum_{k=1}^{+\infty}\int_a^{1-a}\frac{\cos(2\pi kx)}k\,e^{-2\pi ijx}\,dx=\lim_{a\searrow0}\sum_{k=1}^{+\infty}
c_k(j,a),
\end{equation}
where
$$ c_k(j,a):=\int_a^{1-a}\frac{\cos(2\pi kx)}k\,e^{-2\pi ijx}\,dx.$$
Moreover, integrating by parts,
\begin{eqnarray*}&&
|c_k(j,a)|\\&=&\left|
\int_a^{1-a}\frac{d}{dx}\left(\frac{\sin(2\pi kx)}{2\pi ik^2}\right)\,e^{-2\pi ijx}\,dx
\right|\\&=&\left|
\frac{\sin(2k\pi (1-a))}{2\pi ik^2}\,e^{-2\pi ij(1-a)}-
\frac{\sin(2k\pi a)}{2\pi ik^2}\,e^{-2\pi ija}
+\int_a^{1-a}\frac{j\sin(2\pi kx)}{k^2}\,e^{-2\pi ijx}\,dx
\right|\\&\le&\frac{2+j}{2\pi k^2}
\end{eqnarray*}
and this ensures that the last series in~\eqref{AKSmxcs-serf} converges uniformly in~$a$.

We can therefore swap the limit and the summation sign in~\eqref{AKSmxcs-serf} and obtain that
\begin{eqnarray*}&&
\widehat \phi_j=\sum_{k=1}^{+\infty}\lim_{a\searrow0}
c_k(j,a)=\sum_{k=1}^{+\infty}\lim_{a\searrow0}\int_a^{1-a}\frac{\cos(2\pi kx)}k\,e^{-2\pi ijx}\,dx\\&&\qquad
=\sum_{k=1}^{+\infty}\int_0^{1}\frac{\cos(2\pi kx)}k\,e^{-2\pi ijx}\,dx=\sum_{k=1}^{+\infty}\int_0^{1}\frac{e^{2\pi ikx}+e^{-2\pi ikx}}{2k}\,e^{-2\pi ijx}\,dx\\&&\qquad
=\begin{dcases}
\displaystyle\frac1{2j}& {\mbox{ if }}j>0,\\
\displaystyle-\frac1{2j}& {\mbox{ if }}j<0,\\
0 & {\mbox{ if }}j=0.
\end{dcases}
\end{eqnarray*}
This completes the proof of~\eqref{OSK-1-013-1}.

\paragraph{Solution to Exercise~\ref{ui80-0}.} For all~$x>0$, let~$f(x):=(a x)^j e^{-\sigma x}$.
We observe that
$$ f'(x)=aj\,(a x)^{j-1} e^{-\sigma x}-\sigma (a x)^j e^{-\sigma x}$$
and therefore~$f$ is increasing in~$\left(0,\frac{j}\sigma\right)$
decreasing in~$\left(\frac{j}\sigma,+\infty\right)$.

As a consequence,
\begin{eqnarray*}&&
\int_{\lfloor j/\sigma\rfloor+1}^{+\infty}f(x)\,dx=\sum_{k=\lfloor j/\sigma\rfloor+1}^{+\infty}\int_{k}^{k+1}f(x)\,dx
\ge \sum_{k=\lfloor j/\sigma\rfloor+1}^{+\infty}\int_{k}^{k+1}f(k+1)\,dx\\&&\qquad
=\sum_{k=\lfloor j/\sigma\rfloor+1}^{+\infty} f(k+1)=\sum_{k=\lfloor j/\sigma\rfloor+2}^{+\infty} f(k)=
\sum_{k=\lfloor j/\sigma\rfloor+2}^{+\infty}(a k)^j e^{-\sigma k}.
\end{eqnarray*}

Furthermore, if~$\frac{j}\sigma\ge2$,
\begin{eqnarray*}&& \int_1^{\lfloor j/\sigma\rfloor-1}f(x)\,dx=
\sum_{k=1}^{\lfloor j/\sigma\rfloor-2}\int_{k}^{k+1}f(x)\,dx\ge\sum_{k=1}^{\lfloor j/\sigma\rfloor-2}\int_{k}^{k+1}f(k-1)\,dx\\&&\qquad=\sum_{k=1}^{\lfloor j/\sigma\rfloor-2}f(k-1)=\sum_{k=0}^{\lfloor j/\sigma\rfloor-3}f(k)=\sum_{k=0}^{\lfloor j/\sigma\rfloor-3}(ak)^je^{-\sigma k}.
\end{eqnarray*}

Hence,
\begin{equation}\label{AMSGA-1}
\sum_{k=0}^{+\infty}(ak)^je^{-\sigma k}\le 6\left(
\int_0^{+\infty} f(x)\,dx+f\left(\frac{j}\sigma\right)
\right)
\end{equation}

Besides, for all~$j\in\N$,
\begin{equation}\label{MA01} j!\ge\left(\frac{j}e\right)^j.
\end{equation}
This can be checked by induction over~$j$. Indeed, when~$j=0$ the claim is obvious. Suppose now the claim to be true for~$j$ and let us prove it for~$j+1$. To this end, we recall (see e.g.~\cite[formula~(14) on page~64]{MR385023})
that~$\left(1+\frac1j\right)^j<e$ and therefore
\begin{equation*} (j+1)!=(j+1)j!\ge(j+1)
\left(\frac{j}e\right)^j=\left(\frac{j+1}e\right)^{j+1}\,\frac{e}{\left(1+\frac1j\right)^j}
\ge\left(\frac{j+1}e\right)^{j+1}.
\end{equation*}
The inductive step is thereby complete and~\eqref{MA01} is proved.

From~\eqref{MA01}, one infers that
\begin{equation}\label{AMSGA-2}
f\left(\frac{j}\sigma\right)=\frac1{e^j}\left(\frac{aj}\sigma\right)^j\le\left(\frac{a}\sigma\right)^j j!\,.
\end{equation}

Moreover, we can integrate by parts~$j$ times to see that
$$\int_0^{+\infty} f(x)\,dx=a^j\int_0^{+\infty}x^j e^{-\sigma x}\,dx=
\frac{a^j\,j!}{\sigma^j}\int_0^{+\infty} e^{-\sigma x}\,dx=\frac{a^j\,j!}{\sigma^{j+1}}.$$

Gathering this information, \eqref{AMSGA-1} and~\eqref{AMSGA-2}, we obtain the desired result.

\paragraph{Solution to Exercise~\ref{EPPmdcer0}.}
We have that
$$  \cos^2(2\pi kx+\xi_k)=
\frac14\left(e^{i(2\pi kx+\xi_k)}+e^{-i(2\pi kx+\xi_k)}\right)^2
=\frac14\left(e^{2i(2\pi kx+\xi_k)}+e^{-2i(2\pi kx+\xi_k)}+2\right).$$

Additionally, by the Riemann-Lebesgue Lemma (see Theorem~\ref{RLjoqwskcdc}, applied here to the function~$\chi_E$), we find that
$$ \lim_{k\to+\infty}\int_E e^{2i(2\pi kx+\xi_k)}\,dx
=\lim_{k\to+\infty}  e^{2i\xi_k}\int_0^1\chi_E(x)\, e^{4\pi ikx}\,dx=0$$
and similarly
$$ \lim_{k\to+\infty}\int_E e^{-2i(2\pi kx+\xi_k)}\,dx=0.$$

On this account,
\begin{eqnarray*}
&& \lim_{k\to+\infty}\int_E \cos^2(2\pi kx+\xi_k)\,dx= 
\frac14\lim_{k\to+\infty}\int_E\left(e^{2i(2\pi kx+\xi_k)}+e^{-2i(2\pi kx+\xi_k)}+2\right)\,dx= \frac{|E|}2.
\end{eqnarray*}

\paragraph{Solution to Exercise~\ref{EPPmdcer0.1}.}
Arguing for a contradiction, we let~$\rho_k:=\sqrt{a_k^2+b_k^2}$ and we suppose that~\eqref{EPPmdcer0.4} does not hold, whence~$\rho_k$ is not infinitesimal. Accordingly, there exist a sub-sequence~$k_j$ and some~$\epsilon_0>0$ such that, for all~$j\in\N$,
\begin{equation}\label{EPPmdcer0.2} \rho_{k_j}\ge\epsilon_0.\end{equation}

We stress that the point~$\left(\frac{a_k}{\rho_k},-\frac{b_k}{\rho_k}\right)$ lies on the unit circle and therefore we can find~$\xi_k\in\R$ such that~$\cos\xi_k=\frac{a_k}{\rho_k}$ and~$\sin\xi_k=-\frac{b_k}{\rho_k}$.

In this way,\begin{eqnarray*} \rho_k\cos(2\pi kx+\xi_k)&=&
\rho_k\cos(2\pi kx)\cos\xi_k-\rho_k\sin(2\pi kx)\sin\xi_k\\&=&a_k\cos(2\pi kx)+b_k\sin(2\pi kx).\end{eqnarray*}
Hence, on grounds of~\eqref{EPPmdcer0.3},
$$ \lim_{k\to+\infty}\rho_k\cos(2\pi kx+\xi_k)=0.$$
This and~\eqref{EPPmdcer0.2} entail that
$$ \lim_{j\to+\infty}\cos(2\pi k_jx+\xi_{k_j})=0.$$

As a result, by the Dominated Convergence Theorem,
$$ \lim_{j\to+\infty}\int_E \cos^2(2\pi k_jx+\xi_{k_j})\,dx=0,$$
but this is in contradiction with~\eqref{EPPmdcer0.3ns}.

See~\cite[Volume~I, Chapter~IX]{MR1963498} for further reading on this topic.

\paragraph{Solution to Exercise~\ref{EPPmdcer0.1-029.1}.} This problem is interesting for several reasons.
First of all, it is worth noticing that the function~$F$ defined in~\eqref{EPPmdcer0.1-029.231.4.23.1} satisfies the 
Schr{\"o}dinger equation~$i\partial_tF=\partial^2_xF$.
In this regard, the asymptotics in~\eqref{EPPmdcer0.1-029.231.4.23} aims at detecting the asymptotic behaviour of this type of solutions. The interest of~\eqref{EPPmdcer0.1-029.231.4.23} is that three distinct components are clearly visible:
first, a decay term, encoded, up to multiplicative constants, by~$\frac1{\sqrt{t}}$,
second, an oscillatory term, given by~$e^{-\frac{ix^2(t)}{2t}}$, and third, the profile of the solution
along the prescribed slope~$f\left(\frac{\rho}{2}\right)$.

The gist to solve the given exercise is now, philosophically speaking, to focus
on ``stationary phases'' in~\eqref{EPPmdcer0.1-029.231.4.23.1}, that is, omitting the dependence on~$t$, to let~$\phi(\xi):=\xi^2t-\xi x(t)$ and stare at the oscillatory integrand~$ e^{i(\xi^2t-\xi x(t))}=e^{i\phi(\xi)}$.
After the Riemann-Lebesgue Lemma in Theorem~\ref{RLjoqwskcdc}, we have developed the intuition that
rapidly varying phases will average out all contributions to the integral, because they will favor a ``destructive
interference'' between different contributions which have no reason to add up in any coherent fashion.
In this spirit, the main contribution to the integral should come from stationary phases, corresponding to the condition~$\phi'(\xi)=0$, since the slow variation of the phase in this regime may produce some ``constructive interferece'' which is not subject to cancellations.

Of course, this intuition needs to be nailed down by some serious mathematics. The strategy that we follow is indeed a special case of the so-called \index{stationary phase (method of)} \emph{method of stationary phase}, which finds plenty of applications in harmonic analysis, differential equations, and theoretical physics. For the humble purposes of this exercise, the use of this technique goes as follows.

We substitute for~$\eta:=\xi-\frac{x(t)}{2t}$ and, since
$$ \xi^2t-\xi x(t)=\eta^2 t-\frac{x^2(t)}{2t},$$
we see, recalling Exercise~\ref{fr12}, that
\begin{equation}\label{EPPmdcer0.1-029.231.2} F(x(t),t)=e^{-\frac{ix^2(t)}{2t}}\int_0^1 f\left(\eta+\frac{x(t)}{2t}\right)\,e^{i\eta^2 t}\,d\eta.\end{equation}

For all~$\eta\in\R\setminus\{0\}$ and~$t\in\R$, we let
$$ g(\eta,t):=\frac1{2\eta}\left[ f\left(\eta+\frac{x(t)}{2t}\right)-f\left(\frac{x(t)}{2t}\right)\right]$$
and we remark that
$$ \sup_{(\eta,t)\in(\R\setminus\{0\})\times\R}|g(\eta,t)|\le\frac{\|f'\|_{L^\infty((0,1))}}2.$$
We can thereby integrate by parts and conclude that
\begin{equation}\label{EPPmdcer0.1-029.231}\begin{split}&
\int_0^1 f\left(\eta+\frac{x(t)}{2t}\right)\,e^{i\eta^2 t}\,d\eta-\int_0^1 f\left(\frac{x(t)}{2t}\right)\,e^{i\eta^2 t}\,d\eta=
\int_0^1 2\eta g(\eta,t)\,e^{i\eta^2 t}\,d\eta\\&\qquad=
-\frac{i}t\int_0^1 g(\eta,t)\,\frac{d}{d\eta}\big(e^{i\eta^2 t}-1\big)\,d\eta\\&\qquad=
-\frac{i g(1,t)\,(e^{it}-1)}t
+\frac{i}t\int_0^1 \partial_\eta g(\eta,t)\,(e^{i\eta^2 t}-1)\,d\eta.
\end{split}\end{equation}

We also observe that
\begin{eqnarray*}
|\partial_\eta g(\eta,t)\,(e^{i\eta^2 t}-1)|\le C\,|\partial_\eta g(\eta,t)|\,\min\{1,\,\eta^2 t\}\le
C\,\|f\|_{C^1((0,1))}\min\{\eta^{-1},\,\eta t\},
\end{eqnarray*} for some~$C>0$ that we freely rename from line to line, and therefore, for all~$t\ge2$,
\begin{eqnarray*}
\int_0^1 |\partial_\eta g(\eta,t)\,(e^{i\eta^2 t}-1)|\,d\eta\le
C\int_0^{1/\sqrt{t}} \eta t\,d\eta
+C\int_{1/\sqrt{t}}^1\frac{d\eta}\eta\le C\ln t.
\end{eqnarray*}

Plugging this information into~\eqref{EPPmdcer0.1-029.231} we gather that, as~$t\to+\infty$,
$$ \int_0^1 f\left(\eta+\frac{x(t)}{2t}\right)\,e^{i\eta^2 t}\,d\eta-\int_0^1 f\left(\frac{x(t)}{2t}\right)\,e^{i\eta^2 t}\,d\eta=o\left(\frac1{\sqrt{t}}\right).$$
This, \eqref{EPPmdcer0.1-029.231.74.2489}, and~\eqref{EPPmdcer0.1-029.231.2} yield that
\begin{equation}\label{EPPmdcer0.1-029.231.4}
F(x(t),t)=e^{-\frac{ix^2(t)}{2t}}\int_0^1 f\left(\frac{x(t)}{2t}\right)\,e^{i\eta^2 t}\,d\eta+o\left(\frac1{\sqrt{t}}\right).
\end{equation}

Now we recall that, by a contour integration on the complex plane, 
$$ \int_0^{+\infty}\sin(\theta^2)\,d\theta=\frac{1}{2}\sqrt{\frac\pi2}=\int_0^{+\infty}\cos(\theta^2)\,d\theta,$$
see e.g.~\cite[page~205]{MR4506522} (these integrals are often \index{Fresnel integrals} called \emph{Fresnel integrals}).

Consequently, substituting for~$\theta:=\sqrt{t}\,\eta$, we see that \begin{eqnarray*}&&
\int_0^1 e^{i\eta^2 t}\,d\eta=\frac{1}{\sqrt{t}} \int_0^{\sqrt{t}} e^{i\theta^2}\,d\eta
=\frac{1}{\sqrt{t}} \int_0^{\sqrt{t}}\big( \cos(\theta^2)+i\sin(\theta^2)\big)\,d\eta\\&&\qquad
=\frac{1}{\sqrt{t}} \int_0^{+\infty}\big( \cos(\theta^2)+i\sin(\theta^2)\big)\,d\eta+o\left(\frac1{\sqrt{t}}\right)\\&&\qquad=
\frac{1+i}{2}\sqrt{\frac\pi{2t}}+o\left(\frac1{\sqrt{t}}\right).
\end{eqnarray*}
From this observation and~\eqref{EPPmdcer0.1-029.231.4} we arrive at 
\begin{equation*}
F(x(t),t)=\frac{1+i}{2}\sqrt{\frac\pi{2t}}\,e^{-\frac{ix^2(t)}{2t}} f\left(\frac{x(t)}{2t}\right)+o\left(\frac1{\sqrt{t}}\right)
.\end{equation*}
The desired result in~\eqref{EPPmdcer0.1-029.231.4.23} now follows from~\eqref{EPPmdcer0.1-029.231.74.2489}.

\begin{figure}[h]
\includegraphics[height=3.1cm]{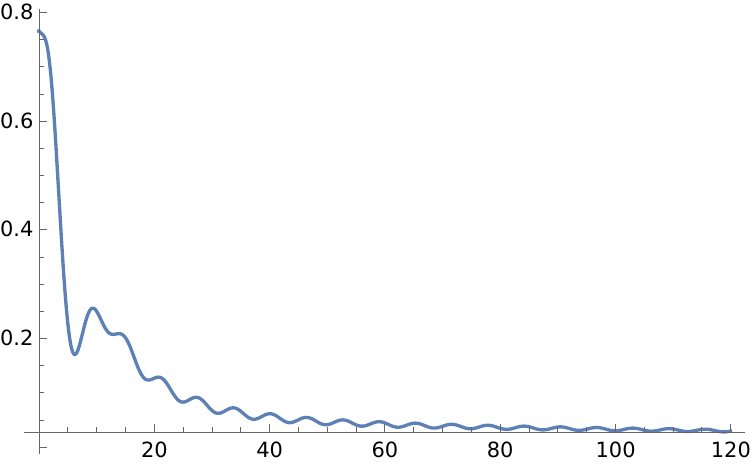}$\;\;$\includegraphics[height=3.1cm]{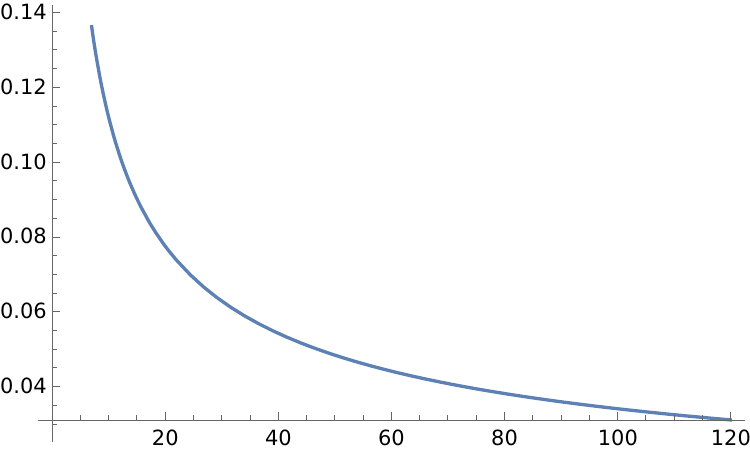}$\;\;$\includegraphics[height=3.1cm]{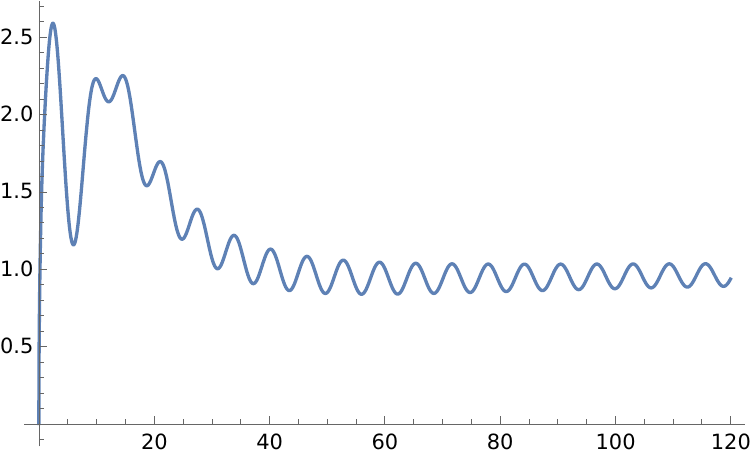}
\centering
\caption{Plot of the real part of the function in~\eqref{EPPmdcer0.1-029.231.4.23.1}, that of the function in~\eqref{EPPmdcer0.1-029.231.4.23}
(disregarding the ``little o''), and of the ratio of these,
when~$x(t):=\frac{t^2}{t^2+1}$ and~$f(\xi):=\cos(\cos(2\pi\xi))$.}\label{0.1.m2je9JJ90-1.ken70t0y3u-u1usqf70KaKcamdi}
\end{figure}

See Figure~\ref{0.1.m2je9JJ90-1.ken70t0y3u-u1usqf70KaKcamdi} for a sketch of the situation.
See also~\cite[Section~10.3]{MR2453734}, \cite[Section~5.10.6]{MR3497072},
\cite[Section~6.9]{MR4647607}, and the bibliography therein for further information about the method of stationary phase.

\paragraph{Solution to Exercise~\ref{WRO-eq00}.}
The statement in~\eqref{WRO-eq1} is certainly mistaken, since
(as in Exercise~\ref{ojqdwn23EEE}) 
one can construct examples of converging Fourier Series
(possibly with nontrivial Fourier coefficients only when~$k\ge0$)
of functions that are not analytic.

Yet, why is the argument provided incorrect? Well, first of all, the statement in~\eqref{WRO-eq3} is sloppy, because it does not specify the domains of the functions involved
(the statement in~\eqref{WRO-eq2} is similarly sloppy too;
also the notion of ``convergence''
in~\eqref{WRO-eq1} is not specified in detail). 

To be more accurate (see e.g.~\cite[Proposition~1.3.3]{MR1916029}), one should replace~\eqref{WRO-eq3} with:
\begin{equation}\label{WRO-eq200}
\begin{split}&
{\mbox{if~$U$ and~$V$ are {\em open} subsets of~$\C$,}}\\&{\mbox{$g:U\to\C$ and~$h:V\to\C$ are analytic, and~$h(V)\subseteq U$,}}\\ &{\mbox{then
the composition~$f:=g\circ h$ is analytic on V.}}
\end{split}
\end{equation}

In our framework, the argument provided tries to apply this statement with~$g$ as in~\eqref{WRO-eq20} and~$h$ the complex exponential map, and indeed~$h$ is entire, so one could, in principle, take any {\em open} set~$V$ in~\eqref{WRO-eq200}.
But since we want the original variable~$x$ to be in the reals,
\begin{equation}\label{WRO-eq2000}
{\mbox{$V$ must contain at least
an open interval~$J$ in the reals.}}\end{equation}

Also, in~\eqref{WRO-eq200} it is required that~$h(V)$ lies in an open set on which~$g$ is analytic. But this would require the power series defining~$g$ in~\eqref{WRO-eq20} to converge in an open neighborhood of~$h(V)$, that is, by~\eqref{WRO-eq2000}, in an open neighborhood of~$h(J)$, which is an arc of the unit circle.

The latter must contain some~$z_0\in\C$ with~$|z_0|>1$. Accordingly (by the Root Test or~\cite[Lemma 1.1.6]{MR1916029}),
\begin{equation}\label{h4okfh40yujLP8bg5t}
\limsup_{k\to+\infty}\big|\widehat f_k\big|^{\frac1k}<1.
\end{equation}

However, there is no reason for this to be true in general, so one cannot guarantee that one can apply~\eqref{WRO-eq200}.

In fact,
by Theorem~\ref{ojqdwn23E}, if~$f$ is not~$C^m(\R^n)$, then there exists a sequence~$k_j\to+\infty$ such that
$$\limsup_{j\to+\infty}|k_j|^{m+2}\big|\widehat f_{k_j}\big|=+\infty.$$
As a result, for large~$j$, we have that~$|k_j|^{m+2}\big|\widehat f_{k_j}\big|\ge1$, therefore~$\big|\widehat f_{k_j}\big|
\ge|k_j|^{-(m+2)}$ and accordingly
\begin{equation*}
\limsup_{j\to+\infty}|\widehat f_{k_j}|^{\frac1{k_j}}\ge
\limsup_{j\to+\infty}|k_j|^{-\frac{m+2}{k_j}}=1
.\end{equation*}This is in contradiction with~\eqref{h4okfh40yujLP8bg5t}
and it shows that, in this case, the hypotheses of~\eqref{WRO-eq200} are not fulfilled.

\section{Solutions to selected exercises of Section~\ref{UNIFORMCO:SECTI}}

\paragraph{Solution to Exercise~\ref{NEWLA.pre}.} This follows from Exercise~\ref{G0ilMajsx912e-14} and
Theorem~\ref{BASw}.

\paragraph{Solution to Exercise~\ref{NEWLA}.}
By~\eqref{PAKSw-L4},
\begin{eqnarray*}&&1+2\sum_{k=1}^N \cos(2\pi kx)=
\sum_{{k\in\Z}\atop{|k|\le N}} \cos(2\pi kx)=
\Re\left(\sum_{{k\in\Z}\atop{|k|\le N}} e^{-2\pi i kx}\right)\\&&\qquad=\Re\left(
\frac{\sin\big((2N+1)\pi x\big)}{\sin(\pi x)}\right)=\frac{\sin\big((2N+1)\pi x\big)}{\sin(\pi x)}.\end{eqnarray*}

\paragraph{Solution to Exercise~\ref{LGEGMCTLFDCD5D325E28A9N}.}
By~\eqref{KASMqwdfed123erDKLI}, 
$$ |D_N(x)|=\left|\sum_{{k\in\Z}\atop{|k|\le N}} e^{-2\pi i
kx}\right|\le\sum_{{k\in\Z}\atop{|k|\le N}} \big|e^{-2\pi i kx}\big|=2N+1=D_N(0).$$

\paragraph{Solution to Exercise~\ref{K-0PIO}.}
By~\eqref{PAKSw-L4},
\begin{eqnarray*}&& D_N\left(\frac12-x\right)-D_N\left(\frac12+x\right)
=\sum_{{k\in\Z}\atop{|k|\le N}} e^{-\pi ik+2\pi i kx}-\sum_{{k\in\Z}\atop{|k|\le N}} e^{-\pi ik-2\pi i kx}\\&&\qquad=
\sum_{{k\in\Z}\atop{|k|\le N}} e^{-\pi ik+2\pi i kx}-\sum_{{k\in\Z}\atop{|k|\le N}} e^{\pi ik+2\pi i kx}=
\sum_{{k\in\Z}\atop{|k|\le N}} (-1)^ke^{2\pi i kx}-\sum_{{k\in\Z}\atop{|k|\le N}} (-1)^ke^{2\pi i kx}
=0.\end{eqnarray*}

\paragraph{Solution to Exercise~\ref{K-1PIO}.}
By~\eqref{PAKSw-L4},
\begin{eqnarray*} \int_{0}^{1} D_N(x)\,dx= \sum_{{k\in\Z}\atop{|k|\le N}}\int_{0}^{1} e^{-2\pi i kx}\,dx=
1+\sum_{{k\in\Z\setminus\{0\}}\atop{|k|\le N}}\int_{0}^{1} e^{-2\pi i kx}\,dx=1.
\end{eqnarray*}

\paragraph{Solution to Exercise~\ref{K-1PIO-bis}.}
By~\eqref{PAKSw-L4},
\begin{eqnarray*}&& \int_{0}^{1/2} D_N(x)\,dx= \sum_{{k\in\Z}\atop{|k|\le N}}\int_{0}^{1/2} e^{-2\pi i kx}\,dx=
\frac12+\sum_{{k\in\Z\setminus\{0\}}\atop{|k|\le N}}\int_{0}^{1/2} e^{-2\pi i kx}\,dx
\\&&\qquad=\frac12+\sum_{k=1}^N\int_{0}^{1/2} \big(e^{2\pi ikx}+e^{-2\pi i kx}\big)\,dx=\frac12+2\sum_{k=1}^N\int_{0}^{1/2} \cos(2\pi kx)\,dx=\frac12.
\end{eqnarray*}

\paragraph{Solution to Exercise~\ref{K-2PIO}.}
By~\eqref{PAKSw-L4},
\begin{eqnarray*}&& \int_{0}^{1} \big(D_N(x)\big)^2\,dx=
\sum_{{k,h\in\Z}\atop{|k|,|h|\le N}}\int_{0}^{1} e^{-2\pi i (k+h)x}\,dx=
\sum_{{{k,h\in\Z}\atop{|k|,|h|\le N}}\atop{k+h=0}}1+
\sum_{{{k,h\in\Z}\atop{|k|,|h|\le N}}\atop{k+h\ne0}}\int_{0}^{1} e^{-2\pi i (k+h)x}\,dx
\\&&\qquad=\sum_{{{k\in\Z}\atop{|k|\le N}}}1+0=2N+1.
\end{eqnarray*}

\paragraph{Solution to Exercise~\ref{K-3PIO}.}
By~\eqref{PAKSw-L4}, \begin{eqnarray*}&&J:=
\int_{-1/2}^{1/2} |D_N(x)|\,dx=\int_{-1/2}^{1/2}
\frac{|\sin\big((2N+1)\pi x\big)|}{|\sin(\pi x)|}\,dx\ge\int_{-1/2}^{1/2}
\frac{|\sin\big((2N+1)\pi x\big)|}{\pi|x|}\,dx\\
&&\qquad=\frac2\pi\int_{0}^{1/2}
\frac{|\sin\big((2N+1)\pi x\big)|}{x}\,dx=\frac2{\pi}\int_{0}^{(2N+1)\pi/2}
\frac{|\sin y|}{y}\,dy\\&&\qquad\ge\frac2\pi\sum_{j=1}^{2N}
\int_{j\pi/2}^{(j+1)\pi/2}\frac{|\sin y|}{y}\,dy\ge
\frac4{\pi^2}\sum_{j=1}^{2N}\frac{1}{j+1}
\int_{j\pi/2}^{(j+1)\pi/2} |\sin y|\,dy.
\end{eqnarray*}

Also, using the substitution~$t:=y-\frac{j\pi}2$,
\begin{eqnarray*}&& \int_{j\pi/2}^{(j+1)\pi/2} |\sin y|\,dy=
\int_{0}^{\pi/2} \left|\sin \left(t+\frac{j\pi}2\right)\right|\,dt\\&&\qquad\ge
\int_{0}^{\pi/2} \min\big\{|\sin t|,\,|\cos t|\big\}\,dt=2-\sqrt{2}
\end{eqnarray*}

From these observations we arrive at
\begin{eqnarray*}&&J\ge
\frac{4(2-\sqrt{2})}{\pi^2}\sum_{j=1}^{2N}\frac{1}{j+1}=\frac{4(2-\sqrt{2})}{\pi^2}\sum_{j=1}^{2N}\frac{1}{j+1}\int_{j+1}^{j+2}dt
\ge\frac{4(2-\sqrt{2})}{\pi^2}\sum_{j=1}^{2N}\int_{j+1}^{j+2}\frac{dt}t\\&&\qquad\qquad=\frac{4(2-\sqrt{2})}{\pi^2}\int_{2}^{2N+2}\frac{dt}t=\frac{4(2-\sqrt{2})}{\pi^2}\,\ln(N+1),
\end{eqnarray*}
yielding the desired result.

\paragraph{Solution to Exercise~\ref{K-3PIO.eDhnZmasKvcTYhFA.1}.} Yes, it is, because, for all~$N\ge2$,
\begin{eqnarray*}&&
\int_{-1/2}^{1/2} |D_N(x)|\,dx=2\int_{0}^{1/2} |D_N(x)|\,dx=2\int_0^{1/2}
\left|\frac{\sin\big((2N+1)\pi x\big)}{\sin(\pi x)}\right|\,dx\\&&\quad
\le C\left( N\int_0^{1/N}\,dx
+\int_{1/N}^{1/2}
\left|\frac{1}{\sin(\pi x)}\right|\,dx\right)\le
 C\left( 1
+\int_{1/N}^{1/2}
\left|\frac{1}{x}\right|\,dx\right)\le C\ln N,
\end{eqnarray*}
for some~$C>0$ varying from line to line.

\paragraph{Solution to Exercise~\ref{IFEJ}.}
On grounds of Lemma~\ref{KASMqwdfed123erDKLI}
and the trigonometric identity in Exercise~\ref{SPDCD-0.01},
$$ D_N(x)=\frac{\sin\big((2N+1)\pi x\big)}{\sin(\pi x)}=
\frac{\sin\big((2N+1)\pi x\big)\,\sin(\pi x)}{\sin^2(\pi x)}=
\frac{\cos(2N\pi x)-\cos\big(2(N+1)\pi x\big)}{2\sin^2(\pi x)}.$$
Therefore, we recognise a telescopic sum and conclude that
$$F_N(x)=\frac{1}{2N\sin^2(\pi x)}\sum_{k=0}^{N-1}\big(
\cos(2k\pi x)-\cos\big(2(k+1)\pi x\big)\big)=
\frac{1}{2N\sin^2(\pi x)}\big(
1-\cos(2N\pi x)\big).$$
The Double-Angle Formula for the cosine now yields the desired resut.

\paragraph{Solution to Exercise~\ref{FEJPO}.} This is a straightforward consequence of
Exercise~\ref{IFEJ}.

\paragraph{Solution to Exercise~\ref{VjweWRFVSHEgKFB8ndVIAPMAQ}.}
By Exercise~\ref{FO:DE:MA}, we know that the Fourier Series of~$g_a$ is
\begin{equation}\label{22VjweWRFVSHEgKFB8ndVIAPMAQ} \sum_{k=1}^{+\infty}\frac{2\sin^2(ka\pi)}{ (k \pi a)^2}\cos(2\pi kx).\end{equation}
Since
$$\left|\frac{2\sin^2(ka\pi)}{ (k \pi a)^2}\cos(2\pi kx)\right|\le
\frac{2}{ (k \pi a)^2},$$
we can use Theorem~\ref{BASw} and conclude that the series in~\eqref{22VjweWRFVSHEgKFB8ndVIAPMAQ}
converges uniformly to~$g_a$.

We also observe that, by Exercise~\ref{fr12},
\begin{eqnarray*}&& \int_0^1 \phi_a(x)\,dx=\int_{-1/2}^{1/2} \phi_a(x)\,dx
=\frac1a\int_{-1/2}^{1/2}\phi\left(\frac{x}a\right)\,dx
\\&&\qquad=\int_{-1/(2a)}^{1/(2a)}\phi(y)\,dy
=2\int_0^1 (1-y)\,dy=1
\end{eqnarray*}
and therefore
\begin{equation}\label{OJALS02-L123} \int_0^1 g_a(x)\,dx=\int_0^1 \phi_a(x)\,dx-1=0.\end{equation}

Also, $\phi_a\ge0$ and consequently
\begin{equation}\label{OJALS02-L1}
\|g_a\|_{L^1((0,1))}\le\|\phi_a\|_{L^1((0,1))}+1=\int_0^1 \phi_a(x)\,dx+1=2.
\end{equation}

Now, for every~$N\in\N$, we define
\begin{equation}\label{KASM01oe2irjfniXqyhdgfouwie2}
G_N(x):=\sum_{3\le k\le N} p_k\,g_{1/k}(x).\end{equation}
We recall~\eqref{OJALS02-L1} and observe that, for all~$M\ge N\ge3$,
$$ \| G_N-G_M\|_{L^1((0,1))}
\le\sum_{N+1\le k\le M} p_k\,\|g_{1/k}\|_{L^1((0,1))}\le
2\sum_{N+1\le k\le M} p_k,$$
which is the tail of a convergent series.

This entails that
\begin{equation}\label{KASM01oe2irjfniXqyhdgfouwie}
{\mbox{$G_N$ converges to~$f$ in~$L^1((0,1))$}}\end{equation}
and, as a result, we have that~$f\in L^1((0,1))$
and it is even and periodic of period~$1$ (since so is~$g_a$).

Furthermore, in view of~\eqref{OJALS02-L123},
\begin{eqnarray*}&&
\left|\int_{0}^{1}f(x)\,dx\right|\le\lim_{N\to+\infty}\|f-G_N\|_{L^1((0,1))}+
\left|\int_{0}^{1}G_N(x)\,dx\right|\\&&\qquad=0+
\left|\sum_{3\le k\le N} p_k\int_{0}^{1}g_{1/k}(x)\,dx\right|=0,
\end{eqnarray*}
completing the proof of the desired result.

\paragraph{Solution to Exercise~\ref{VjweWRFVSHEgKFB8ndVIAPMAQ-bis}.}
By~\eqref{KASM01oe2irjfniXqyhdgfouwie2}, \eqref{KASM01oe2irjfniXqyhdgfouwie}, and Exercises~\ref{LINC-bisco} and~\ref{09iulappcoemmdba},
for all~$m\in\N$,
\begin{eqnarray*}
\widehat f_m=\lim_{N\to+\infty} \widehat G_{N,m}=\lim_{N\to+\infty} 
\sum_{3\le k\le N} p_k\,\widehat g_{1/k,m}.
\end{eqnarray*}

Hence, in light of~\eqref{22VjweWRFVSHEgKFB8ndVIAPMAQ}, the $m$th Fourier coefficient of~$f$
of sine type vanishes (as well as the average of~$f$) and the $m$th Fourier coefficient of~$f$ of cosine type is
$$ \sum_{k=3}^{+\infty} p_k\,\frac{2\sin^2(m \pi/k)}{ (m \pi /k)^2}
.$$
Therefore, the Fourier Series of~$f$ in trigonometric form is
\begin{equation}\label{oqdl93r-20404kkdd} \sum_{m=1}^{+\infty}\sum_{k=3}^{+\infty} p_k\,\frac{2\sin^2(m \pi/k)}{ (m \pi /k)^2}\,\cos(2m\pi x).\end{equation}

\begin{figure}[h]
\includegraphics[height=3.1cm]{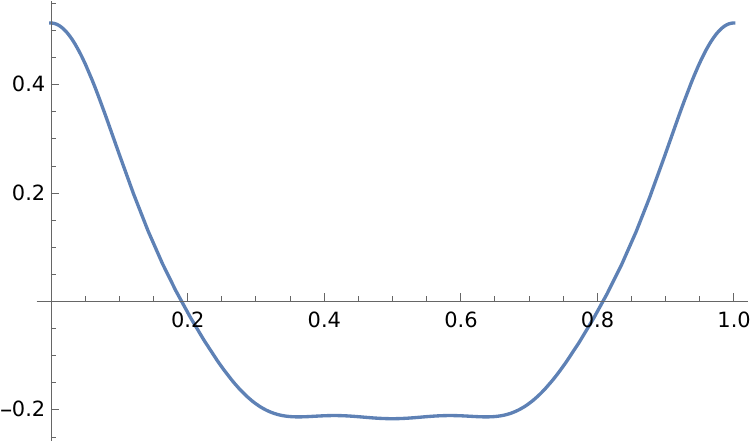}$\;$\includegraphics[height=3.1cm]{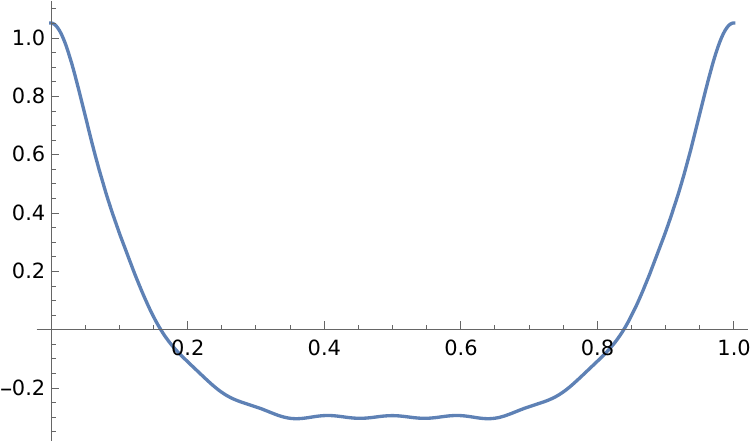}$\;$\includegraphics[height=3.1cm]{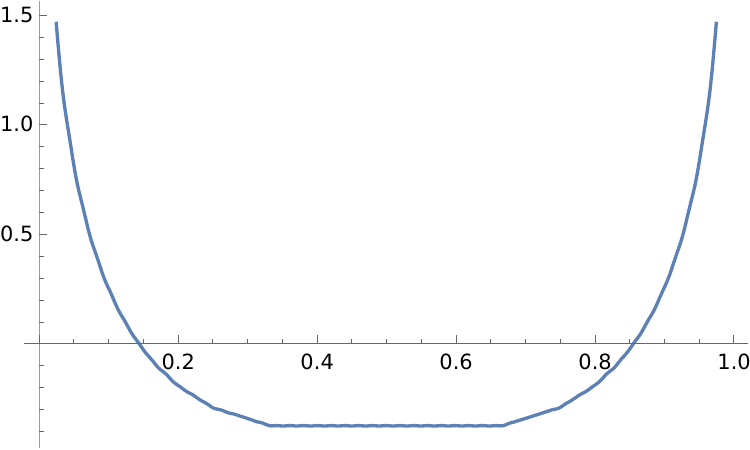}
\centering
\caption{Plot of~$\displaystyle\sum_{m=1}^{N}\sum_{k=3}^{N} p_k\,\frac{2\sin^2(m \pi/k)}{ (m \pi /k)^2}\,\cos(2m\pi x)$
with~$p_k:=1/k^2$ and~$N\in\{4,10,50\}$.}\label{0.1.m299ken70t0y3u-u1usqf70KaKcamdi}
\end{figure}

See Figure~\ref{0.1.m299ken70t0y3u-u1usqf70KaKcamdi} for a sketch of this situation.

\paragraph{Solution to Exercise~\ref{K-4PIO}.}
In light of~\eqref{PAKSw-L4} and Exercise~\ref{K-1PIO-bis},
\begin{equation}\label{oByu8mftb6m2345}\begin{split}&\frac12=\int_{0}^{1/2} D_N(x)\,dx= 
\int_0^{1/2} \frac{\sin\big((2N+1)\pi x\big)}{\sin(\pi x)}\,dx\\&\qquad
=\int_0^{1/2} \sin\big((2N+1)\pi x\big)\left[\frac{1}{\sin(\pi x)}-\frac1{\pi x}\right]\,dx
+\int_0^{1/2} \frac{\sin\big((2N+1)\pi x\big)}{\pi x}\,dx\\
&\qquad=\int_0^{1/2} f(x)\,\sin\big((2N+1)\pi x\big)\,dx
+\frac1\pi \int_0^{\frac{(2N+1)\pi}2} \frac{\sin y}{y}\,dy,\end{split}\end{equation}
where, for all~$x\in\left(0,\frac12\right]$, we set
$$ f(x):=\frac{1}{\sin(\pi x)}-\frac1{\pi x}.$$
Since
$$ \lim_{x\searrow0}f(x)=\frac1\pi\lim_{x\searrow0} \left(\frac{1}{x+O(x^3)}-\frac1{ x}\right)=\frac1\pi\lim_{x\searrow0} \frac1x\left(1+O(x^2)-1\right)=0,
$$
we see that~$f$ can be continuously extended to~$\left[0,\frac12\right]$ and then, by odd reflection,
to a uniformly continuous function in~$\left(-\frac12,\frac12\right]$.

Up to extending~$f$ periodically, we can thus assume that~$f$ is odd, periodic of period~$1$, with~$f\in L^1((0,1))$.

Hence, also~$g(x):=f(x)\cos(\pi x)$ and~$h(x):=f(x)\sin(\pi x)$ are periodic of period~$1$ and belong to~$L^1((0,1))$.

The Riemann-Lebesgue Lemma (see Theorem~\ref{RLjoqwskcdc}) thus gives that
\begin{eqnarray*}&& \lim_{N\to+\infty} \int_0^{1/2} f(x)\,\sin\big((2N+1)\pi x\big)\,dx=
\frac12\lim_{N\to+\infty} \int_{-1/2}^{1/2} f(x)\,\sin\big((2N+1)\pi x\big)\,dx\\&&\qquad=
\frac12\lim_{N\to+\infty} \int_{-1/2}^{1/2} f(x)\,\Big(
\sin(2N\pi x)\cos(\pi x)+\cos(2N\pi x)\sin(\pi x)
\Big)\,dx\\&&\qquad=\lim_{N\to+\infty}\left( \int_{-1/2}^{1/2} g(x)\,\sin(2N\pi x)\,dx+\int_{-1/2}^{1/2} h(x)\,\cos(2N\pi x)\,dx\right)=0.
\end{eqnarray*}

This and~\eqref{oByu8mftb6m2345} give that
\begin{equation}\label{ftinLsidyacht}\begin{split}
\frac\pi2&=\pi\lim_{N\to+\infty}\left(\int_0^{1/2} f(x)\,\sin\big((2N+1)\pi x\big)\,dx
+\frac1\pi \int_0^{\frac{(2N+1)\pi}2} \frac{\sin y}{y}\,dy\right)\\&=
\lim_{N\to+\infty} \int_0^{\frac{(2N+1)\pi}2} \frac{\sin y}{y}\,dy.\end{split}\end{equation}

This is almost the desired claim, since for the moment we have chosen a special sequence,
namely~${\frac{(2N+1)\pi}2}$.
To complete the desired claim, we need to take care of any divergent sequence. For this,
let~$\epsilon>0$. Then, if~$L>0$, we let~$N:=\left\lfloor \frac{L}{\pi}\right\rfloor-1$
and we use~\eqref{ftinLsidyacht} to see that, as long as~$L$ is large enough, 
\begin{eqnarray*}&&
\left| \int_0^L\frac{\sin y}{y}\,dy-\frac\pi2\right|\le
\left| \int_0^{\frac{(2N+1)\pi}2}\frac{\sin y}{y}\,dy-\frac\pi2\right|+\left| \int_{\frac{(2N+1)\pi}2}^L\frac{\sin y}{y}\,dy\right|\\&&\qquad
\le\epsilon+\int_{\frac{(2N+1)\pi}2}^L\frac{dy}{y}\le\epsilon+\int_{L-\frac\pi2}^L\frac{2\,dy}{(2N+1)\pi}=\epsilon+\frac{1}{2N+1}
\le2\epsilon,
\end{eqnarray*}
which establishes the desired result.

\paragraph{Solution to Exercise~\ref{CO213421FSBCA.02-ejNS12c2g2n4}.}
For~$m\in\N\setminus\{0\}$, let
$$ \sigma_m:=\int_{(m-1)\pi}^{m\pi}\frac{\sin x}{x}\,dx.$$
Looking at the positivity set of the sine function, we see that
\begin{equation}\label{CO213421FSBCA.02-ejNS12c2g2n4.E01}
{\mbox{$\sigma_m>0$ if~$m$ is odd and~$\sigma_m<0$ if~$m$ is even.}}\end{equation}

On this account, $|\sigma_m|=(-1)^{m+1}\sigma_m$ and consequently, using the substitution~$y:=x-m\pi$ and~$y:=x-(m-1)\pi$, we find that
\begin{equation}\label{CO213421FSBCA.02-ejNS12c2g2n4.E02}\begin{split}&
|\sigma_{m+1}|-|\sigma_m|=(-1)^{m+2}\sigma_{m+1}-(-1)^{m+1}\sigma_m=
(-1)^m(\sigma_{m+1}+\sigma_m)\\&\qquad=(-1)^m\left(
\int_{m\pi}^{(m+1)\pi}\frac{\sin x}{x}\,dx+
\int_{(m-1)\pi}^{m\pi}\frac{\sin x}{x}\,dx\right)\\&\qquad=(-1)^m\left(
\int_{0}^{\pi}\frac{\sin (y+m\pi)}{y+m\pi}\,dy+
\int_{0}^{\pi}\frac{\sin (y+(m-1)\pi)}{y+(m-1)\pi}\,dy\right)\\&\qquad=(-1)^m\left(
(-1)^m\int_{0}^{\pi}\frac{\sin y}{y+m\pi}\,dy+(-1)^{m-1}
\int_{0}^{\pi}\frac{\sin y}{y+(m-1)\pi}\,dy\right)\\&\qquad=
\int_{0}^{\pi}\frac{\sin y}{y+m\pi}\,dy-
\int_{0}^{\pi}\frac{\sin y}{y+(m-1)\pi}\,dy
\\&\qquad=-\pi
\int_{0}^{\pi}\frac{\sin y}{(y+(m-1)\pi)(y+m\pi)}\,dy\\&\qquad=-\pi
\int_{0}^{\pi}\frac{|\sin y|}{(y+(m-1)\pi)(y+m\pi)}\,dy\\&\qquad<0.
\end{split}\end{equation}

Let now
$$ \mu_m:=\int_0^{\pi m}\frac{\sin x}{x}\,dx=\sum_{j=1}^m \sigma_j.$$
If~$m$ is odd, we deduce from~\eqref{CO213421FSBCA.02-ejNS12c2g2n4.E01} that~$\mu_{m+1}=\mu_{m}+\sigma_{m+1}<\mu_m$
and therefore
\begin{equation}\label{CO213421FSBCA.02-ejNS12c2g2n4.E07}
\sup_{m\in\N\setminus\{0\}}\mu_m=\sup_{{m\in\N\setminus\{0\}}\atop{{\mbox{\tiny{$m$ odd}}}}}\mu_m=\sup_{j\in\N}\mu_{2j+1}.
\end{equation}

We also observe that, in view of~\eqref{CO213421FSBCA.02-ejNS12c2g2n4.E01} and~\eqref{CO213421FSBCA.02-ejNS12c2g2n4.E02},
$$ \mu_{2j+3}-\mu_{2j+1}=\sigma_{2j+3}+\sigma_{2j+2}=|\sigma_{2j+3}|-|\sigma_{2j+2}|<0.
$$
Combining this with~\eqref{CO213421FSBCA.02-ejNS12c2g2n4.E07}, we conclude that
$$ \sup_{m\in\N\setminus\{0\}}\mu_m=\mu_1,$$ which is the desired result in~\eqref{CO213421FSBCA.02-ejNS12c2g2n4.E010}.

Now we prove~\eqref{CO213421FSBCA.02-ejNS12c2g2n4.E011}. To this end, we define
$$ \nu_m:=\int_0^{2\pi m}\frac{\sin x}{x}\,dx$$
and we observe that, by~\eqref{CO213421FSBCA.02-ejNS12c2g2n4.E01} and~\eqref{CO213421FSBCA.02-ejNS12c2g2n4.E02},
\begin{eqnarray*}
\nu_{m+1}-\nu_m=\sigma_{2m+1}+
\sigma_{2m+2}=|\sigma_{2m+1}|-|\sigma_{2m+2}|>0
\end{eqnarray*}
and therefore~$\nu_m$ is increasing with respect to~$m$.

For this reason,
$$ \sup_{m\in\N\setminus\{0\}}\int_0^{2\pi m}\frac{\sin x}{x}\,dx=\sup_{m\in\N\setminus\{0\}}\nu_m=
\lim_{m\to+\infty}\nu_m=\int_0^{+\infty}\frac{\sin x}{x}\,dx.$$

Thanks to this and Exercise~\ref{K-4PIO},
the proof of~\eqref{CO213421FSBCA.02-ejNS12c2g2n4.E011} is complete.

\paragraph{Solution to Exercise~\ref{DINI-addC-2-maunif}.}
Yes, in this case~$f$ is necessarily continuous at all points of~$I$. To see this, one can argue as follows.
Suppose for the sake of contradiction that~$f$ is discontinuous at a certain point of the interval~$I$. Up to horizontal translations and reflections, we can suppose that the discontinuity point is the origin (in particular, $0\in I$) and that there exists a sequence of~$x_k$ of positive real numbers such that~$x_k\to0$ as~$k\to+\infty$, with~$|f(x_k)-f(0)|\ge \eta_0$, for some~$\eta_0>0$. In fact, up to replacing~$f$ with~$f-f(0)$, we can also assume without loss of generality that~$f(0)=0$.

By the uniform version of Dini's Condition~\eqref{DINI-UNO}, we can find~$\delta>0$ such that
$$\frac{\eta_0}{10}\ge
\sup_{x\in I}\int_{-\delta}^{\delta}\left|\frac{f (x + t) - f(x)}{t}\right|\,dt.$$
Since, given~$\delta'\in\left(0,\frac\delta4\right)$, for large~$k$ 
when~$t\in(-x_k+\delta',-x_k+2\delta')$
we have that
$$ |x_k+t|\ge |t|-|x_k|\ge \frac{|t|}2+\frac{\delta'}4-|x_k|\ge\frac{|t|}2,$$
we thereby find that
\begin{equation}\label{34rtghrbg.89sudhs9JMSynnIEmsdLjsdj.d}\begin{split}&
\frac{\eta_0}{10}\ge
\int_{-\delta}^{\delta}\left|\frac{f (x_k + t) - f(x_k)}{t}\right|\,dt\ge
\int_{-x_k+\delta'}^{-x_k+2\delta'}\frac{|f ( x_k+t)-f(x_k)|}{|t|}\,dt\\&\qquad\qquad\qquad\ge\frac12\int_{-x_k+\delta'}^{-x_k+2\delta'}\frac{|f ( x_k+t)-f(x_k)|}{|x_k+t|}\,dt
.\end{split}\end{equation}
In the same vein,
\begin{equation*}\begin{split}
&\frac{\eta_0}{10}\ge
\int_{-\delta}^{\delta}\left|\frac{f ( t) - f(0)}{t}\right|\,dt\ge
\int_{\delta'}^{2\delta'}\left|\frac{f ( t) - f(0)}{t}\right|\,dt
=\int_{\delta'}^{2\delta'}\left|\frac{f ( t) }{t}\right|\,dt\\&\qquad\qquad
=\int_{-x_k+\delta'}^{-x_k+2\delta'}\left|\frac{f ( x_k+\tau)}{x_k+\tau}\right|\,d\tau
\ge
\int_{-x_k+\delta'}^{-x_k+2\delta'}\frac{|f(x_k)|-|f ( x_k+\tau)-f(x_k)|}{|x_k+\tau|}\,d\tau\\&\qquad\qquad
\ge
\int_{-x_k+\delta'}^{-x_k+2\delta'}\frac{{\eta_0}}{|x_k+\tau|}\,d\tau
-
\int_{-x_k+\delta'}^{-x_k+2\delta'}\frac{|f ( x_k+\tau)-f(x_k)|}{|x_k+\tau|}\,d\tau.
\end{split}\end{equation*}
This and~\eqref{34rtghrbg.89sudhs9JMSynnIEmsdLjsdj.d} yield that, for large~$k$,
\begin{equation*}\begin{split}&
\int_{-x_k+\delta'}^{-x_k+2\delta'}\frac{{\eta_0}}{|x_k+\tau|}\,d\tau
\le\frac{\eta_0}{10}+
\int_{-x_k+\delta'}^{-x_k+2\delta'}\frac{|f ( x_k+\tau)-f(x_k)|}{|x_k+\tau|}\,d\tau
\le\frac{\eta_0}{10}+\frac{\eta_0}5=\frac{3\eta_0}{10}.
\end{split}\end{equation*}
As a result, dividing by~$\eta_0$ and taking the limit in~$k$,
$$ \ln2=
\int_{\delta'}^{2\delta'}\frac{{d\tau}}{\tau}=
\lim_{k\to+\infty}\int_{-x_k+\delta'}^{-x_k+2\delta'}\frac{{d\tau}}{|x_k+\tau|}\le\frac3{10},$$
which is a contradiction.

\paragraph{Solution to Exercise~\ref{C1uni-c-EXEOP}.} We recall the square wave in Exercise~\ref{SQ:W}, that is
$$ f(x):=\begin{dcases}
1&{\mbox{ if }}x\in\displaystyle\left[0,\frac12\right),\\
-1&{\mbox{ if }}x\in\displaystyle\left[\frac12,1\right).
\end{dcases}$$
We let~$I_1:=\left(0,\frac12\right)$, $I_2:=\left(\frac12,1\right)$, $g_1(x):=1$ and~$g_2(x):=-1$. We have that~$f=g_j$ in~$I_j$, for~$j\in\{1,2\}$.

We point out that~$S_{N,g_1}(x)=1$ and $S_{N,g_2}(x)=-1$ for all~$x\in\R$, therefore, up set of null measure,
$$ f(x)=\begin{dcases} S_{N,g_1} (x)&{\mbox{ if }}x\in I_1,\\S_{N,g_2}(x)&{\mbox{ if }}x\in I_2.\end{dcases}$$
Suppose now, for a contradiction, that
$$ \lim_{N\to+\infty}\sup_{x\in I_j}|S_{N,f}(x)-S_{N,g_j}(x)|=0,$$
for all~$j\in\{1,2\}$.

Then,
\begin{eqnarray*}&&
\lim_{N\to+\infty} \|S_{N,f}-f\|_{L^\infty(0,1)}=
\lim_{N\to+\infty} \max\big\{\|S_{N,f}-f\|_{L^\infty((0,1/2))},\,
\|S_{N,f}-f\|_{L^\infty((1/2,1))}\big\}\\&&\qquad
=\lim_{N\to+\infty} \max\big\{\|S_{N,f}-S_{N,g}\|_{L^\infty(I_1)},\,
\|S_{N,f}-S_{N,g}\|_{L^\infty(I_2)}\big\}=0.
\end{eqnarray*}
As a consequence, $S_{N,f}$ is a Cauchy sequence in~$L^\infty((0,1))$ and therefore it must converge uniformly to some function~$\phi\in L^\infty((0,1))$.
Since~$S_{N,f}$ is continuous, we obtain that~$\phi$ is continuous.
But by the uniqueness of the pointwise limit, up to set of measure zero, we have that~$\phi=f$, which has a jump discontinuity, contradiction.

See Section~\ref{SEC:GIBBS-PH} for more details on this phenomenon.

\paragraph{Solution to Exercise~\ref{7un.mSIjkaja21dcA.13}.} By Theorem~\ref{7un.mSIjkaja21dcA.13th}, we know that~$|\widehat f_k|\le Ce^{-\sigma|k|}$. Hence, by Theorem~\ref{BASw},
\begin{eqnarray*}
&&\sup_{x\in\R}|f(x)-S_N(x)|=\sup_{x\in\R}\left| \sum_{k\in\Z}\widehat f_k\,e^{2\pi ikx}-
\sum_{{k\in\Z}\atop{|k|\le N}}\widehat f_k\,e^{2\pi ikx}\right|
=\sup_{x\in\R}\left| \sum_{{k\in\Z}\atop{|k|> N}}\widehat f_k\,e^{2\pi ikx}\right|\\&&\qquad\le 
\sum_{{k\in\Z}\atop{|k|> N}}|\widehat f_k|\le
C\sum_{{k\in\Z}\atop{|k|> N}}e^{-\sigma|k|}
\le Ce^{-\sigma N},
\end{eqnarray*}
up to renaming~$C$.

\paragraph{Solution to Exercise~\ref{7un.mSIjkaja21dcA.13BIS}.} By Theorem~\ref{SMXC22},
we know that~$|\widehat f_k|\le \frac{C}{|k|^m}$. Hence, by Theorem~\ref{BASw},
\begin{eqnarray*}
&&\sup_{x\in\R}|f(x)-S_N(x)|=\sup_{x\in\R}\left| \sum_{k\in\Z}\widehat f_k\,e^{2\pi ikx}-
\sum_{{k\in\Z}\atop{|k|\le N}}\widehat f_k\,e^{2\pi ikx}\right|
=\sup_{x\in\R}\left| \sum_{{k\in\Z}\atop{|k|> N}}\widehat f_k\,e^{2\pi ikx}\right|\\&&\qquad\le 
\sum_{{k\in\Z}\atop{|k|> N}}|\widehat f_k|\le
C\sum_{{k\in\Z}\atop{|k|> N}}\frac1{|k|^m}
\le \frac{C}{N^{m-1}},
\end{eqnarray*}
up to renaming~$C$.

\paragraph{Solution to Exercise~\ref{PIPIO-010}.} We can assume that~$k-1\ge m$ (otherwise the summation runs over an empty family of indices and we are done). Then, if~$K:=k-1-m\ge0$,
\begin{eqnarray*}&&
\left|\sum_{{j\in\Z}\atop{m\le j\le k-1}}e^{2\pi ijx}\right|=
\left|e^{2\pi imx}\sum_{{j\in\Z}\atop{m\le j\le k-1}}e^{2\pi i(j-m)x}\right|=
\left|\sum_{{j\in\Z}\atop{m\le j\le k-1}}e^{2\pi i(j-m)x}\right|\\&&\qquad=
\left|\sum_{k=0}^Ke^{2\pi ikx}\right|=
\left|\frac{1-e^{2\pi i(K+1)x}}{1-e^{2\pi ix}}\right|=
\left|\frac{e^{-\pi ix}-e^{2\pi i\left(K+\frac12\right)x}}{e^{-\pi ix}-e^{\pi ix}}\right|\\
&&\qquad=\left|\frac{e^{-\pi ix}-e^{2\pi i\left(K+\frac12\right)x}}{2\sin(\pi x)}\right|\le
\frac{|e^{-\pi ix}|+\big|e^{2\pi i\left(K+\frac12\right)x}\big|}{2|\sin(\pi x)|}=\frac1{|\sin(\pi x)|}.
\end{eqnarray*}

\paragraph{Solution to Exercise~\ref{PIPIO-0}.}
We employ the ``summation by parts'' formula (see Exercise~\ref{SBPF}), used here with~$\alpha_k:=\gamma_k$ and
$$\beta_k:=\sum_{{j\in\Z}\atop{m\le j\le k-1}}e^{2\pi ijx}.$$
We observe that~$\beta_{k+1}-\beta_{k}=e^{2\pi ikx}$ and~$\beta_m=0$ (since no indices are present in the summation in this case).

Consequently, we gather from~\eqref{BYPA} that
\begin{equation}\begin{split}\label{IMNSITGNITNIHBU}&\sum_{{j\in\Z}\atop{m\le j\le n}}\gamma_{k}e^{2\pi ikx}=
\sum_{{j\in\Z}\atop{m\le j\le n}}\alpha_{k}(\beta_{k+1}-\beta_{k})=\left(\alpha_{n}\beta_{n+1}-\alpha_{m}\beta_{m}\right)
-\sum_{{j\in\Z}\atop{m+1\le j\le n}}\beta_{k}(\alpha_{k}-\alpha_{k-1})\\&\qquad=
\gamma_n\beta_{n+1}-\sum_{{j\in\Z}\atop{m+1\le j\le n}}\beta_{k}(\gamma_{k}-\gamma_{k-1})
\end{split}\end{equation}
and consequently, by our assumptions on~$\gamma_k$,
\begin{equation}\label{NOECX}
\left|\sum_{{j\in\Z}\atop{m\le j\le n}}\gamma_{k}e^{2\pi ikx}\right|\le
\gamma_n |\beta_{n+1}|+\sum_{{j\in\Z}\atop{m+1\le j\le n}}|\beta_{k}|(\gamma_{k-1}-\gamma_{k}).
\end{equation}

Besides, owing to Exercise~\ref{PIPIO-010},
\begin{equation}\label{HLPATOST} |\beta_k|=\left|\sum_{{j\in\Z}\atop{m\le j\le k-1}}e^{2\pi ijx}\right|
\le\frac{1}{|\sin(\pi x)|},
\end{equation}
hence we deduce from~\eqref{NOECX} that
\begin{eqnarray*}|\sin(\pi x)|\,
\left|\sum_{{j\in\Z}\atop{m\le j\le n}}\gamma_{k}e^{2\pi ikx}\right|&\le&
\gamma_n +\sum_{{j\in\Z}\atop{m+1\le j\le n}}(\gamma_{k-1}-\gamma_{k})\\&=&
\gamma_n+\gamma_m-\gamma_n\\&=&\gamma_m,
\end{eqnarray*}
which entails the desired result.

\paragraph{Solution to Exercise~\ref{PIPIOmaq}.} 
By Exercise~\ref{PIPIO-0},
$$\sup_{x\in(\epsilon,1-\epsilon)}\left|\sum_{{k\in\Z}\atop{m\le k\le n}}\gamma_k e^{2\pi ikx}\right|\le\sup_{x\in(\epsilon,1-\epsilon)} \frac{\gamma_m}{|\sin(\pi x)|}=
\frac{\gamma_m}{\sin(\pi \epsilon)}
.$$
For large~$m$, this becomes as small as we wish, giving the desired result.

\paragraph{Solution to Exercise~\ref{PIPIO}.} By the periodicity of the sine function, we can assume that~$x\in[0,1)$.
Also, since the desired result is obvious when~$x=0$, we can suppose without loss of generality that~$x\in(0,1)$.

Let now~$p\in\N$, with~$p\ge1$. Then, by Exercise~\ref{PIPIO-0},
$$ \left|\sum_{k=p}^N\gamma_k\sin(2\pi kx)\right|
=\left|\Im\left(\sum_{k=p}^N\gamma_k e^{2\pi ik x}\right)\right|\le\left|\sum_{k=p}^N\gamma_k e^{2\pi i kx}\right|\le
\frac{\gamma_p}{|\sin(\pi x)|}\le\frac{1}{p\,|\sin(\pi x)|}.$$
As a consequence,
\begin{equation}\label{me57uyjhn9uhgv}\begin{split}&
\left|\sum_{k=1}^N\gamma_k\sin(2\pi kx)\right|\le\left|\sum_{k=1}^{p-1}\gamma_k\sin(2\pi kx)\right|+\frac{1}{p\,|\sin(\pi x)|}\\
&\qquad\le \sum_{k=1}^{p-1}\gamma_k|\sin (2\pi kx)|+\frac{1}{p\,|\sin(\pi x)|}\le\sum_{k=1}^{p-1}\frac{|\sin(2\pi kx)|}k+\frac{1}{p\,|\sin(\pi x)|}.
\end{split}\end{equation}

We also remark that
$$ \left|\sin\left(2k\pi (1-x)\right)\right|=
|\sin(2k\pi -2\pi kx)|=|\sin(-2\pi kx)|=|\sin(2\pi kx)|
$$
and 
$$ |\sin(\pi (1-x))|=|\sin(-\pi x)|=|\sin(\pi x)|.$$
Therefore, if~$y:=\min\left\{x,\,1-x\right\}$, we have that~$y\in\left(0,\frac12\right]$ and we can write~\eqref{me57uyjhn9uhgv}
as
\begin{equation*}
\left|\sum_{k=1}^N\gamma_k\sin(2\pi kx)\right|\le\sum_{k=1}^{p-1}\frac{|\sin(2\pi ky)|}k+\frac{1}{p\,|\sin(\pi y)|}.
\end{equation*}

Since~$\big|\sin(2\pi ky)\big|\le 2\pi ky$, we thereby conclude that
\begin{equation}\label{ALZIK}\begin{split}&
\left|\sum_{k=1}^N\gamma_k\sin(2\pi kx)\right|
\le\sum_{k=1}^{p-1}\frac{2\pi ky}k+\frac{1}{p\,|\sin(\pi y)|}
= 2\pi(p-1)y+\frac{1}{p\sin(\pi y)}.\end{split}
\end{equation}

Choosing~$p:=1+ \left\lfloor\frac1{2y}\right\rfloor$, we have that
$$ 2\pi (p-1)y= 2\pi \left\lfloor\frac1{2y}\right\rfloor y\le\pi$$
and
$$ p\sin(\pi y)\ge\frac{\sin(\pi y)}{2y}.$$
These observations and~\eqref{ALZIK} give that
\begin{equation}\label{ALZIK-2}\begin{split}&
\left|\sum_{k=1}^N\gamma_k\sin(2\pi kx)\right|
\le\pi+\frac{2y}{\sin(\pi y)}.\end{split}
\end{equation}

Now we claim that, for all~$y\in\left(0,\frac12\right]$,
\begin{equation}\label{00efgPLJMNS}
\frac{\sin(\pi y)}{2y}\ge 1.
\end{equation}
To this end, we write
$$ \phi(y):=\frac{\sin(\pi y)}{\pi y}=\sum_{j=0}^{+\infty}\frac{(-1)^j (\pi y)^{2j}}{(2j+1)!},$$
whence
$$ \phi'(y)=\sum_{j=1}^{+\infty}\frac{(-1)^j 2j (\pi y)^{2j-1}}{(2j+1)!}=
\sum_{j=1}^{+\infty} (-1)^j\mu_j(y),$$
with, for~$j\ge1$,
$$ \mu_{j}(y):=\frac{2j (\pi y)^{2j-1}}{(2j+1)!}.$$
We observe that, for all~$y\in\left(0,\frac12\right]$,
$$ \frac{\mu_{j+1}(y)}{\mu_{j}(y)}=\frac{(\pi y)^2}{2j(2j+3)}\le\frac{(\pi y)^2}{10}\le\frac{\pi^2}{40}\le1.$$
As a result,
\begin{eqnarray*} &&\phi'(y)=\lim_{K\to+\infty}
\sum_{j=1}^{2K} (-1)^j\mu_j(y)=
-\big(\mu_1(y)- \mu_2(y)\big)-\ldots-\big( \mu_{2K-1}(y)-\mu_{2K}(y)\big)\le0.\end{eqnarray*}
This states that~$\phi$ is decreasing in~$\left(0,\frac12\right]$, and accordingly, for every~$y\in\left(0,\frac12\right]$,
$$ \frac{\sin(\pi y)}{\pi y}=\phi(y)\ge\phi\left(\frac12\right)=\frac{2}\pi,$$
yielding~\eqref{00efgPLJMNS}.

Thus, inserting~\eqref{00efgPLJMNS} into~\eqref{ALZIK-2} we find that
$$ \left|\sum_{k=1}^N\gamma_k\sin(2\pi kx)\right|
\le\pi+1,$$ as desired.

\paragraph{Solution to Exercise~\ref{RENOLLME}.}
No. For instance, taking~$\gamma_k:=\frac1{k+1}$ and~$x:=0$, we have that
$$ \sum_{k=1}^N\gamma_k\cos(2\pi kx)=\sum_{k=1}^N\frac1{k+1},$$
which diverges as~$N\to+\infty$.

\paragraph{Solution to Exercise~\ref{PIPIO-00-123-321}.} Let~$\epsilon>0$.
By assumption, there exists~$k_\epsilon\in\N$ such that for all~$k\ge k_\epsilon$ we have that~$k\gamma_k\le\epsilon$. We define
$$\widetilde\gamma_k:=\begin{dcases}\displaystyle\frac1k&{\mbox{ if }}k\in\{0,\dots,k_\epsilon-1\},\\
\\
\displaystyle\frac{\gamma_{k}}\epsilon&{\mbox{ if }}k\ge k_\epsilon.
\end{dcases}$$

We stress that
$$ \frac1{k_\epsilon-1}\ge\frac1{k_\epsilon}\ge\frac{\gamma_{k_\epsilon}}\epsilon$$
and therefore~$\widetilde\gamma_k\ge\widetilde\gamma_{k+1}\ge0$.

Moreover, $k\widetilde\gamma_k\le1$ and
we are thereby in the position of using Exercise~\ref{PIPIO}, concluding that, for all~$N\in\N$,
\begin{eqnarray*}\pi+1&\ge&
\left|\sum_{k=1}^N\widetilde\gamma_k\sin(2\pi k\delta)\right|\\
&=&\left|\sum_{0\le k<k_\epsilon}\frac{\sin(2\pi k\delta)}k
+\sum_{k_\epsilon\le k\le N}\frac{\gamma_k\,\sin(2\pi k\delta)}\epsilon\right|.
\end{eqnarray*}
Taking the limit as~$N\to+\infty$, we find that
\begin{eqnarray*}\pi+1&\ge&\left|\sum_{0\le k<k_\epsilon}\frac{\sin(2\pi k\delta)}k
+\sum_{k\ge k_\epsilon}\frac{\gamma_k\,\sin(2\pi k\delta)}\epsilon\right|\\
&\ge&\left|\sum_{k\ge k_\epsilon}\frac{\gamma_k\,\sin(2\pi k\delta)}\epsilon\right|-\left|\sum_{0\le k<k_\epsilon}\frac{\sin(2\pi k\delta)}k\right|.
\end{eqnarray*}
Sending now~$\delta\to0$, we obtain
\begin{eqnarray*} \pi+1&\ge&\lim_{\delta\to0}\left|\sum_{k\ge k_\epsilon}\frac{\gamma_k\,\sin(2\pi k\delta)}\epsilon\right|\\&\ge&\lim_{\delta\to0}\left(\left|\sum_{k=0}^{+\infty}\frac{\gamma_k\,\sin(2\pi k\delta)}\epsilon\right|-\left|
\sum_{0\le k<k_\epsilon}\frac{\gamma_k\,\sin(2\pi k\delta)}\epsilon
\right|\right)\\&=&\lim_{\delta\to0}\left|\sum_{k=0}^{+\infty}\frac{\gamma_k\,\sin(2\pi k\delta)}\epsilon\right|.
\end{eqnarray*}
This says that $$
\lim_{\delta\to0}\left|\sum_{k=0}^{+\infty} \gamma_k\,\sin(2\pi k\delta)\right|\le(\pi+1)\epsilon$$
and thus the desired result now follows by taking~$\epsilon$ as small as we wish.

\paragraph{Solution to Exercise~\ref{ESIPNo}.} The answer is negative. Suppose, for a contradiction, that such a function~$f$ exists. Then, by~\eqref{muFVPFbaGtpf}, for all~$k\ge2$,
\begin{equation}\label{muFVPFbaGtpf2} \int_0^1 f(x)\,\sin(2\pi kx)\,dx=\frac{1}{\ln k}.\end{equation}
Let now~$\epsilon\in(0,1)$. Then, there exists~$\delta_0\in\left(0,\frac12\right)$ such that, for every~$\delta\in(0,\delta_0)$,
$$ \int_{(0,\delta)\cup(1-\delta,1)}|f(x)|\,dx\le\epsilon.$$
For~$N\in\N$, let also$$ U_N(x):=\sum_{k=2}^N \frac{\sin(2\pi kx)}{k}.$$By Exercise~\ref{PIPIO}, used here\footnote{We observe that
this step would not remain valid if one had the cosine function instead of the sine function, see \label{REVCNCOFOmdfe8uj.O9023}
Exercise~\ref{RENOLLME}.} with~$\gamma_k:=1/k$, we know that, for all~$x\in\R$,
$$ |U_N(x)|\leq\pi+1.$$As a result,\begin{equation*}\begin{split}&\left| \int_0^1 f(x)\,U_N(x)\,dx-\int_\delta^{1-\delta} f(x)\,U_N(x)\,dx
\right|\\&\qquad\le \sup_{x\in\R}|U_N(x)|\int_{(0,\delta)\cup(1-\delta,1)}|f(x)|\,dx\\&\qquad\le(\pi+1)\epsilon.\end{split}\end{equation*}
Moreover,\begin{eqnarray*}\int_\delta^{1-\delta} f(x)\,U_N(x)\,dx=
\sum_{k=2}^N \frac{1}{k} \int_\delta^{1-\delta} f(x)\,\sin(2\pi kx)\,dx
\end{eqnarray*}and therefore\begin{eqnarray*}(\pi+1)\|f\|_{L^1((0,1))}&\ge&\int_0^1 f(x)\,U_N(x)\,dx\\&\ge&\int_\delta^{1-\delta} f(x)\,U_N(x)\,dx-(\pi+1)\epsilon\\&=&\sum_{k=2}^N \frac{1}{k} \int_\delta^{1-\delta} f(x)\,\sin(2\pi kx)\,dx-(\pi+1)\epsilon.\end{eqnarray*}Now we send~$\delta\searrow0$ and we use~\eqref{muFVPFbaGtpf2} to find that\begin{eqnarray*}(\pi+1)\|f\|_{L^1((0,1))}&\ge&\sum_{k=2}^N \frac{1}{k}\int_0^{1} f(x)\,\sin(2\pi kx)\,dx-(\pi+1)\epsilon
\\&=&\sum_{k=2}^N \frac{1}{k\,\ln k}-(\pi+1)\epsilon
.\end{eqnarray*}
And now we send~$N\to+\infty$, obtaining that, if~$\epsilon>0$ is sufficiently small,
\begin{eqnarray*}&& (\pi+1)\big(\|f\|_{L^1((0,1))}+1\big)\ge\sum_{k=2}^{+\infty} \frac{1}{k\,\ln k}=
\sum_{\ell=1}^{+\infty}\;\;\sum_{k=2^\ell}^{2^{\ell+1}-1}
\frac{1}{k\,\ln k}\ge\sum_{\ell=1}^{+\infty}\;\;\sum_{k=2^\ell}^{2^{\ell+1}-1}\frac{1}{2^{\ell+1}\,\ln (2^{\ell+1})}\\
&&\qquad=\sum_{\ell=1}^{+\infty}\frac{2^\ell}{2^{\ell+1}\,\ln (2^{\ell+1})}=\frac1{2\ln2}\sum_{\ell=1}^{+\infty}\frac{1}{\ell+1}=+\infty,
\end{eqnarray*}
which is a contradiction.

\begin{figure}[h]
\includegraphics[height=3.1cm]{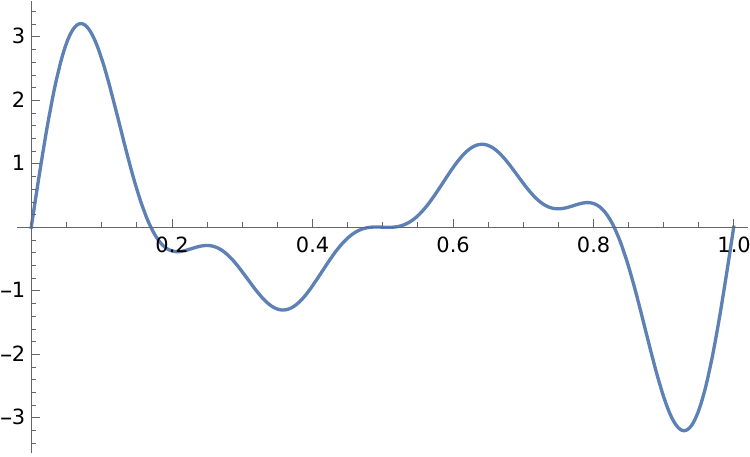}$\;$\includegraphics[height=3.1cm]{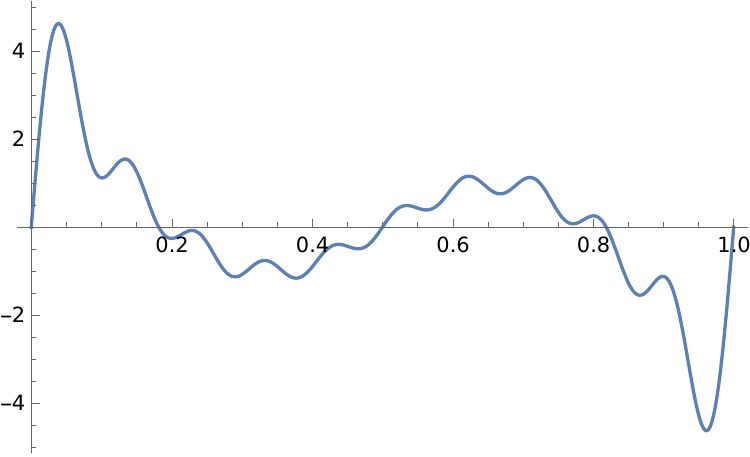}$\;$\includegraphics[height=3.1cm]{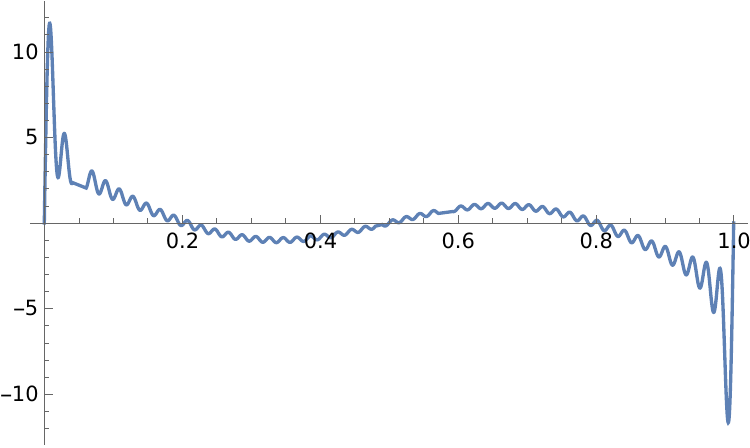}
\centering
\caption{Plot of~$\displaystyle\sum_{k=2}^{N} \frac{\sin(2\pi kx)}{\ln k}$
with~$N\in\{5,10,50\}$.}\label{01003-170KaKcamdi}
\end{figure}

See Figure~\ref{01003-170KaKcamdi} to visualise how things go wrong.

Examples of this type are interesting, since they show that there exist  trigonometric series with coefficients tending to zero which are not Fourier series of any integrable function: see~\cite[Volume~I, Section~30]{MR171116}, \cite[Chapter~7]{MR545506},
and~\cite[Chapter~I, Section~4.2]{MR2039503} for more information about this topic.

\paragraph{Solution to Exercise~\ref{ESIPNo-2}.}
The answer is yes (which is extremely\footnote{Maybe, it is slightly less surprising if one reconsiders the footnote
on page~\pageref{REVCNCOFOmdfe8uj.O9023}.} surprising, compared with Exercise~\ref{ESIPNo}).
The idea is that we cannot only ``integrate by parts'' once, but also twice (and here,
since we are dealing with series instead of integrals, we will replace ``integration by parts''
with ``summation by parts'', as put forth in Exercise~\ref{SBPF}, in any case, the gist is to perform this trick twice). A more thorough approach to this technique will be presented in Section~\ref{FEJERKESE}.

For all~$x\in(0,1)$, we define
$$ f(x):=\sum_{k=2}^{+\infty} \frac{\cos(2\pi kx)}{\ln k}.$$
We stress that this definition is well-posed, due to Exercise~\ref{PIPIOmaq}, which also tells us that the series above is uniformly convergent in~$(\epsilon,1-\epsilon)$, for any given~$\epsilon\in\left(0,\frac12\right)$.

Our objective to answer the question proposed is now to show that
\begin{equation}\label{TSI-09-1}
f\in L^1((0,1))
\end{equation}
and
\begin{equation}\label{TSI-09-2}
\int_0^1 f(x)\,\cos(2\pi kx)\,dx=
\begin{dcases}
0&{\mbox{ if }}k\in\{0,1\},\\
\displaystyle\frac{1}{\ln k}&{\mbox{ if }}k\ge2.
\end{dcases}\end{equation}

To this end, we observe that
\begin{eqnarray*}&& I_{\epsilon,k,\ell}:=\int_0^\epsilon \cos(2\pi kx)\,\cos(2\pi\ell x) \,dx
=\int_{1-\epsilon}^1\cos(2\pi kx)\,\cos(2\pi\ell x) \,dx\\&&\qquad\qquad\qquad
=\frac1{4\pi}\left(\frac{\sin(2 \pi (k - \ell) \epsilon)}{k - \ell} + \frac{\sin(2 \pi (k + \ell) \epsilon)}{k + \ell}\right).\end{eqnarray*}
Hence, by Exercise~\ref{PIPIO-00-123-321},
\begin{equation}\label{OAJSm-1-023oriuihgASI} \lim_{\epsilon\searrow0}\sum_{k=2}^{+\infty} \frac{I_{\epsilon,k,\ell}}{\ln k}=0.\end{equation}

In addition, by the Lagrange Trigonometric Identity (see Exercise~\ref{NEWLA}),
$$ 2\sum_{k=1}^N \cos(2\pi kx)=1+D_N(x).$$
Hence, by the summation by parts formula in Exercise~\ref{SBPF},
and omitting, for the sake of shortness, the dependence on~$x$ in the kernels,
\begin{equation}\label{102owieujrhf0NqaRndfKpoqrL0je01}\begin{split}\sum_{k=2}^{n}\frac{2\cos(2\pi kx)}{\ln k}&=
\sum_{k=2}^{n}\frac1{\ln k}(D_{k}-D_{k-1})\\&=\left(\frac{D_{n}}{\ln n}-\frac{D_{1}}{\ln 2}\right)-\sum_{k=3}^{n}D_{k-1}\left(\frac1{\ln k}-\frac1{\ln (k-1)}\right).
\end{split}\end{equation}

Now we use once more the summation by parts formula in Exercise~\ref{SBPF}, by defining
$$ F_{N}(x):=\frac{1}{N}\sum_{k=0}^{N-1}D_{k}(x).$$
This is the Fej\'er Kernel which was introduced in Exercise~\ref{IFEJ}
and will be studied in detail in Section~\ref{FEJERKESE}. For the moment, we use it in the summation by parts formula to obtain that, for every sequence~$\{\alpha_k\}_{k\in\N}$,
\begin{equation*} \sum_{k=3}^{n}D_{k-1}\alpha_{k}=
\sum_{k=3}^{n}(F_k-F_{k-1})\alpha_{k}
=\left(\alpha_{n}F_{n}-\alpha_{3}F_{2}\right)-\sum_{k=4}^{n}F_{k-1}(\alpha_{k}-\alpha_{k-1}).\end{equation*}

In particular, we find that
\begin{equation}\label{ILCOSEPOCm21plw1} \begin{split}&\sum_{k=3}^{n}D_{k-1}\left(\frac1{\ln k}-\frac1{\ln (k-1)}\right)=
\left(\frac1{\ln n}-\frac1{\ln (n-1)}\right)F_{n}-\left(\frac1{\ln 3}-\frac1{\ln 2}\right)F_{2}\\&\qquad
-\sum_{k=4}^{n}F_{k-1}\left(\frac1{\ln k}-\frac2{\ln (k-1)}+\frac1{\ln (k-2)}\right).\end{split}
\end{equation}
The latter is an interesting ``second order relation'' involving the inverse of the logarithm
that we now analyze further using a convexity argument. Specifically, by the convexity of the function
$$ [2,+\infty)\ni \tau\mapsto\phi(\tau):=\frac1{\ln \tau},$$
whose second derivative is
$$ [2,+\infty)\ni \tau\mapsto\phi''(\tau)=\frac{2 + \ln\tau}{\tau^2 \ln^3\tau}\ge0,$$
it follows that, for all~$k\ge4$,
\begin{eqnarray*}&&
\frac1{\ln k}-\frac2{\ln (k-1)}+\frac1{\ln (k-2)}=
\phi(k)-\phi(k-1)+\phi(k-2)-\phi(k-1)\\&&\qquad=\int_0^1 \phi'(k-t)\,dt-\int_0^1 \phi'(k-t-1)\,dt=
\iint_{(0,1)^2} \phi''(k-t-s)\,ds\,dt\ge0.
\end{eqnarray*}

Using this information and the fact that~$F_N(x)\ge0$ for all~$x\in\R$ (see Exercise~\ref{FEJPO}) we deduce from~\eqref{ILCOSEPOCm21plw1} that $$\sum_{k=3}^{n}D_{k-1}\left(\frac1{\ln k}-\frac1{\ln (k-1)}\right)\le
\left(\frac1{\ln n}-\frac1{\ln (n-1)}\right)F_{n}-\left(\frac1{\ln 3}-\frac1{\ln 2}\right)F_{2}\le-\left(\frac1{\ln 3}-\frac1{\ln 2}\right)F_{2}.$$
Therefore, in light of~\eqref{102owieujrhf0NqaRndfKpoqrL0je01},
\begin{equation*}\sum_{k=2}^{n}\frac{2\cos(2\pi kx)}{\ln k}\ge
\frac{D_{n}}{\ln n}-\frac{D_{1}}{\ln 2}+
\left(\frac1{\ln 3}-\frac1{\ln 2}\right)F_{2}.\end{equation*}
Hence, owing to~\eqref{PAKSw-L4}, for every~$x\in(0,1)$,
\begin{equation*}\sum_{k=2}^{n}\frac{2\cos(2\pi kx)}{\ln k}\ge-
\frac{1}{\ln n\,|\sin(\pi x)|}-\frac{\sin(3\pi x)}{\ln 2\,\sin(\pi x)}+\left(\frac1{\ln 3}-\frac1{\ln 2}\right)F_{2}\ge-\frac{1}{\ln n\,|\sin(\pi x)|}-C,
\end{equation*}
for some~$C>0$,
and thus, sending~$n\to+\infty$,
\begin{equation*}f(x)\ge-C.\end{equation*}

Consequently,
\begin{equation}\label{GMPDI}\begin{split}&\int_0^1|f(x)|\,dx\le
\int_0^1\left|f(x)+C\right|\,dx+C\\&\qquad=
\int_0^1\left(f(x)+C\right)\,dx+C=
\int_0^1 f(x)\,dx+2C.\end{split}\end{equation}

Now, to prove~\eqref{TSI-09-1} and~\eqref{TSI-09-2} we use the uniform convergence of the series in~$(\epsilon,1-\epsilon)$ to see that, for all~$\ell\in\N$,
\begin{eqnarray*}&& \int_\epsilon^{1-\epsilon} f(x)\,\cos(2\pi \ell x)\,dx=
\sum_{k=2}^{+\infty} \int_\epsilon^{1-\epsilon}\frac{\cos(2\pi kx)\,\cos(2\pi\ell x)}{\ln k}\,dx\\
&&\qquad=\sum_{k=2}^{+\infty}\left( \int_0^{1}\frac{\cos(2\pi kx)\,\cos(2\pi\ell x)}{\ln k}\,dx
-\frac{2I_{\epsilon,k,\ell}}{\ln k}\right)=
\frac{\chi_{[2,+\infty)}(\ell)}{\ln\ell}-\sum_{k=2}^{+\infty}\frac{2I_{\epsilon,k,\ell}}{\ln k}.
\end{eqnarray*}
Hence, by~\eqref{OAJSm-1-023oriuihgASI}, for all~$\ell\in\N$,
$$ \int_0^{1} f(x)\,\cos(2\pi \ell x)\,dx=
\frac{\chi_{[2,+\infty)}(\ell)}{\ln\ell}.$$
The proof of~\eqref{TSI-09-2} is thereby complete.

Also, taking~$\ell:=0$ here above and recalling~\eqref{GMPDI},
$$\int_0^1|f(x)|\,dx\le2C.$$
This proves~\eqref{TSI-09-1}.

\begin{figure}[h]
\includegraphics[height=3.1cm]{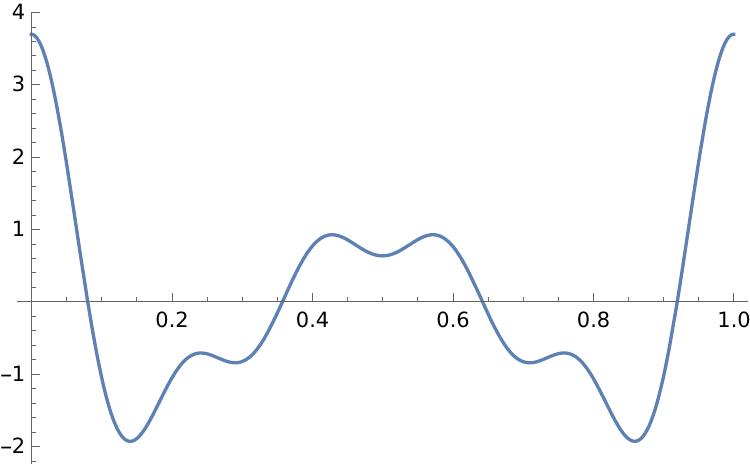}$\;$\includegraphics[height=3.1cm]{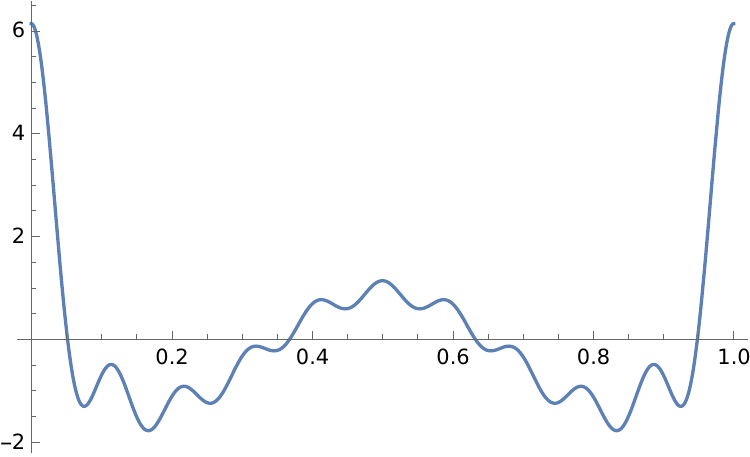}$\;$\includegraphics[height=3.1cm]{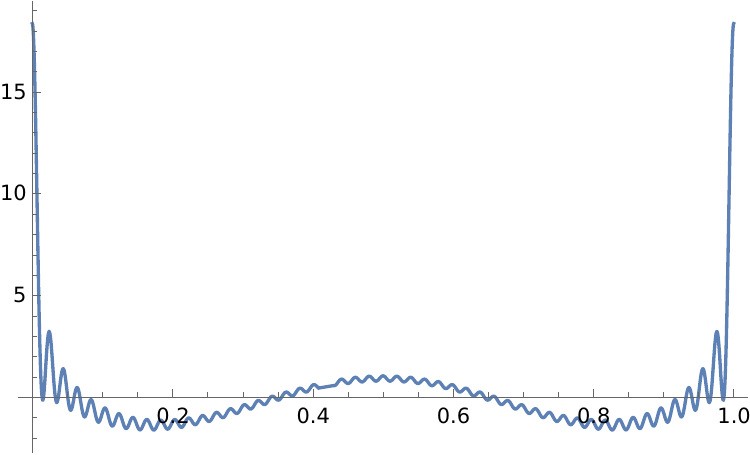}
\centering
\caption{Plot of~$\displaystyle\sum_{k=2}^{N} \frac{\cos(2\pi kx)}{\ln k}$
with~$N\in\{5,10,50\}$.}\label{01003-170KaKcamdi-2}
\end{figure}

See Figure~\ref{01003-170KaKcamdi-2} to visualise the situation in this case:
however, compared to Figure~\ref{01003-170KaKcamdi}, it is probably not immediate to detect just from
the pictures the striking differences between the situation here and that in Exercise~\ref{ESIPNo}.

See~\cite[Volume~I, Section~30]{MR171116} 
and~\cite[Chapter~7]{MR545506}
for additional information.

\paragraph{Solution to Exercise~\ref{LESENUEDGE}.} The proof goes along the same lines as that of Exercise~\ref{ESIPNo-2}.

For more information about this phenomenon, see~\cite[Section~7.3.1]{MR545506}.

\paragraph{Solution to Exercise~\ref{ESIPNo-3}.} The answer is yes (which is perhaps a bit surprising, compared with Exercise~\ref{ESIPNo}).
By Exercise~\ref{1-CLA:AUNIFO}, we can define
$$ f(x):=\sum_{k=2}^{+\infty}\frac{\sin(2\pi kx)}{k\,\ln k}$$
and deduce from the uniform convergence of the series that~$f$ is continuous and periodic of period~$1$.

Also, thanks to this uniform convergence, for every~$\ell\in\N\setminus\{0\}$,
\begin{eqnarray*}
\int_0^1f(x)\,\sin(2\ell\pi x)\,dx=\sum_{k=2}^{+\infty}\int_0^1\frac{\sin(2\pi kx)\,\sin(2\ell\pi x)}{k\,\ln k}\,dx
=\frac1{\ell\,\ln\ell}
\end{eqnarray*}
and similarly, for every~$\ell\in\N$,
\begin{eqnarray*}
\int_0^1f(x)\,\cos(2\ell\pi x)\,dx=\sum_{k=2}^{+\infty}\int_0^1\frac{\sin(2\pi kx)\,\cos(2\ell\pi x)}{k\,\ln k}\,dx
=0,
\end{eqnarray*}
from which we obtain that the Fourier Series of~$f$ in trigonometric form is precisely the one stated
in~\eqref{ITFS8901P13cde}.

\begin{figure}[h]
\includegraphics[height=3.1cm]{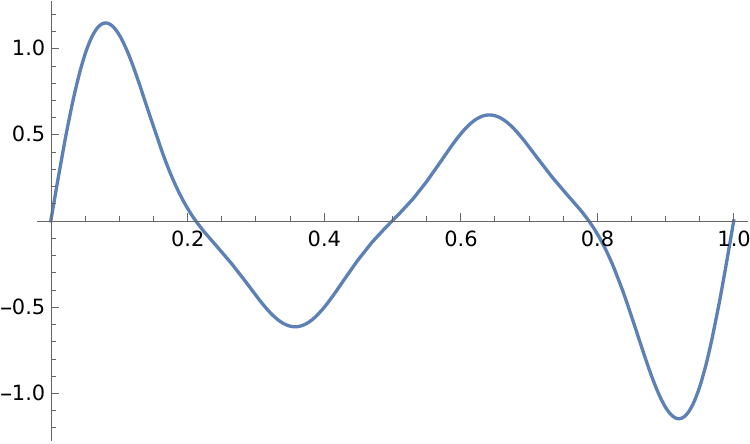}$\;$\includegraphics[height=3.1cm]{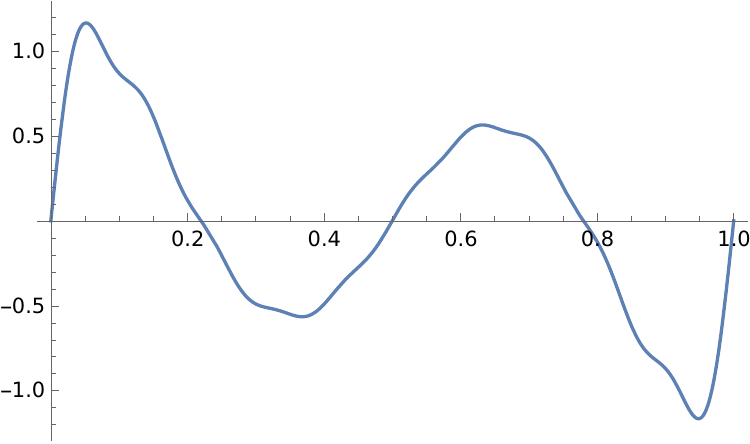}$\;$\includegraphics[height=3.1cm]{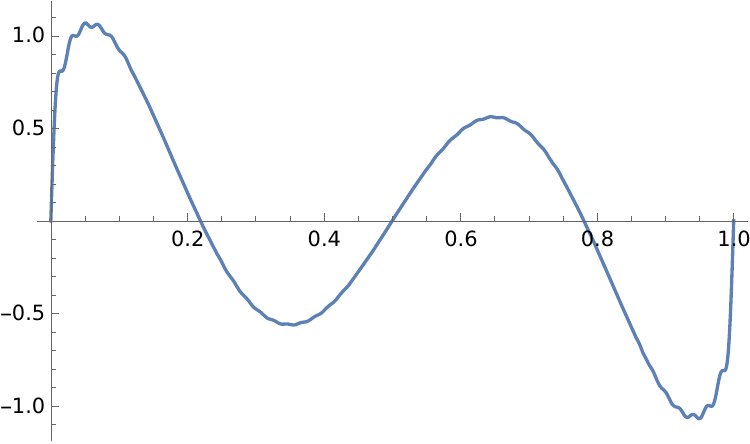}
\centering
\caption{Plot of~$\displaystyle\sum_{k=2}^{N} \frac{\sin(2\pi kx)}{k\,\ln k}$
with~$N\in\{5,10,50\}$.}\label{01003-17i8ngt56fd0KaKcamdi-2}
\end{figure}

See Figure~\ref{01003-17i8ngt56fd0KaKcamdi-2} for a sketch of this situation.

\paragraph{Solution to Exercise~\ref{ESIPNo-2-ap2}.}
The answer is in the negative (and this has to be compared with the positive answer in Exercise~\ref{ESIPNo-2}:
all in all, the function~$f$ happens to belong to~$L^1((0,1))$ but not to~$L^2((0,1))$.

Indeed, suppose by contradiction that there exists~$f\in L^2((0,1))$ whose Fourier Series is
$$ \sum_{k=2}^{+\infty} \frac{\cos(2\pi kx)}{\ln k}=
\sum_{{k\in\Z}\atop{|k|\ge2}}\frac{e^{2\pi ikx}}{2\ln |k|}.
$$
Then, by Bessel's Inequality~\eqref{AJSa},
$$\|f\|^2_{L^2((0,1))}\ge
\sum_{k\in\Z} |\widehat f_k|^2 = \sum_{{k\in\Z}\atop{|k|\ge2}}\frac{1}{4\ln^2 |k|}=+\infty,
$$which is a contradiction.

\paragraph{Solution to Exercise~\ref{FTCPASNo1}.} Let~$K_{-1}:=0$, $K_0:=1$, and~$K_1:=3$.
Then, for all~$\ell\in\N\cap[2,+\infty)$, we consider a sequence~$K_\ell\in\N$ such that 
\begin{equation}\label{VAGAES}
K_{\ell+1}\ge 1+2K_\ell
\end{equation}
and such that for all~$\ell\ge2$ and~$k\ge K_\ell$ we have that
\begin{equation}\label{SIKA} |\sigma_k|\le\frac{1}{\ell+1}.\end{equation}

Roughly speaking, the idea is now to define~$\gamma_{K_\ell}:=\frac{1}{\ell}$ for all~$\ell\in\N$ and then extend the definition
to all~$\gamma_k$ for~$k\in[K_{\ell-1}+1,K_\ell-1]$ by linear interpolation. This is a bit problematic for small values of~$\ell$, due to the appearance of zero divisors, but these are just finitely many instances which we can deal with directly.

Concretely, for all~$\ell\ge2$ and~$k\in\N\cap[K_{\ell-1},K_{\ell}]$, we define
$$ \gamma_k:=\frac{\ell K_\ell-(\ell-1)K_{\ell-1}-k}{\ell(\ell-1)(K_\ell-K_{\ell-1})}.$$

We observe that this defines~$\gamma_k$ only for~$k\ge K_1=3$, hence~$\gamma_0$, $\gamma_1$ and~$\gamma_2$ need to be defined directly. For this, we set
\begin{eqnarray*}&&\gamma_2:=\max\big\{\gamma_3,\,2\gamma_3-\gamma_4,\,|\sigma_2|\big\},\\&&
\gamma_1:=\max\big\{\gamma_2,\,2\gamma_2-\gamma_3,\,|\sigma_1|\big\},\\
{\mbox{and }}&&\gamma_0:=\max\big\{\gamma_1,\,2\gamma_1-\gamma_2,\,|\sigma_0|\big\}
\end{eqnarray*}

We observe that this sequence satisfies all the desired properties. Indeed, when~$k\in\{0,1,2\}$ one can proceed by direct inspection. Also, when~$k\in\N\cap[K_{\ell-1},K_{\ell}]\cap[3,+\infty)$, in view of~\eqref{SIKA}
we have that
\begin{eqnarray*}&&|\sigma_k|-\gamma_k=
|\sigma_k|-\frac{\ell K_\ell-(\ell-1)K_{\ell-1}-k}{\ell(\ell-1)(K_\ell-K_{\ell-1})}\\&&\qquad
\le |\sigma_k|-\frac{\ell K_\ell-(\ell-1)K_{\ell-1}-K_\ell}{\ell(\ell-1)(K_\ell-K_{\ell-1})}
=|\sigma_k|-\frac1\ell
\le0.\end{eqnarray*}

Moreover, if~$k$, $k+1\in\N\cap[K_{\ell-1},K_{\ell}]\cap[3,+\infty)$, then
\begin{eqnarray*}
\gamma_{k+1}-\gamma_k=-\frac{1}{\ell(\ell-1)(K_\ell-K_{\ell-1})}\le0
,\end{eqnarray*}
which establishes the desired monotonicity of the sequence under consideration.

Additionally, if~$k-1$, $k$, $k+1\in\N\cap[K_{\ell-1},K_{\ell}]\cap[3,+\infty)$, then
\begin{eqnarray*}
\gamma_{k+1}+\gamma_{k-1}-2\gamma_k
=0.\end{eqnarray*}

Also, by~\eqref{VAGAES},
\begin{eqnarray*}&&
\gamma_{K_{\ell}+1}+\gamma_{K_{\ell}-1}-2\gamma_{K_{\ell}}\\&&\qquad=
\frac{(\ell+1) K_{\ell+1}-\ell K_{\ell}-K_\ell-1}{(\ell+1) \ell (K_{\ell+1}-K_\ell)}
+\frac{\ell K_\ell-(\ell-1)K_{\ell-1}-(K_\ell-1)}{\ell(\ell-1)(K_\ell-K_{\ell-1})}\\&&\qquad\qquad\qquad\qquad
-2\frac{(\ell+1) K_{\ell+1}-\ell K_{\ell}-K_\ell}{(\ell+1)\ell(K_{\ell+1}-K_\ell)}
\\&&\qquad=\frac{1}{\ell}
-\frac{1}{(\ell+1) \ell (K_{\ell+1}-K_\ell)}+\frac{1}{\ell}+
\frac{1}{\ell(\ell-1)(K_\ell-K_{\ell-1})}
-\frac{2}\ell
\\&&\qquad=-\frac{1}{(\ell+1) \ell (K_{\ell+1}-K_\ell)}+
\frac{1}{\ell(\ell-1)(K_\ell-K_{\ell-1})}\\
&&\qquad\ge0.
\end{eqnarray*}
These arguments show the desired convexity of the sequence under consideration.

See also~\cite[Section~7.1.5]{MR545506} and~\cite[Lemma~3.3.2]{MR3243734}
for further details on convexity methods for sequences.

\paragraph{Solution to Exercise~\ref{FTCPASNo2}.}
This follows from Exercises~\ref{LESENUEDGE} and~\ref{FTCPASNo1}.

For further readings on this topic, see e.g.~\cite[Chapter~7]{MR545506}, \cite[Theorem 3.3.4]{MR3243734} and the references therein.

\paragraph{Solution to Exercise~\ref{FTCPASNo2CONYI}.}
We present here an example due to Iosif Pinelis.

We define
$$\widetilde\sigma_k:=\sup_{j\ge k}|\sigma_j|.$$

We let~$L\in\N$ such that if~$|x|\le\frac1L$ then~$\frac{\sin x}x\ge\frac12$.
We define
$$ p_k:=\begin{dcases}2(\widetilde\sigma_m-\widetilde\sigma_{m+1}) & {\mbox{ if there exists~$m\in\N$ such that~$k=4mL$,}}\\0&{\mbox{ otherwise.}}
\end{dcases}$$
We remark that~$\widetilde\sigma_{m+1}\le\widetilde\sigma_m$, whence~$p_k\ge0$.

Furthermore,
$$ \sum_{k=0}^{+\infty} p_k<+\infty,$$
since it is a telescopic series.

We let~$f$ be as in Exercises~\ref{VjweWRFVSHEgKFB8ndVIAPMAQ} and~\ref{VjweWRFVSHEgKFB8ndVIAPMAQ-bis}.
Then, by~\eqref{oqdl93r-20404kkdd}, the Fourier Series takes the trigonometric form in~\eqref{090-0987yhn-q0woieuxn83gfakowrynt0-w.MAQ-bis} with
\begin{eqnarray*}\gamma_m&:=&
\sum_{k=3}^{+\infty} p_k\,\frac{2\sin^2(m \pi/k)}{ (m \pi /k)^2}.
\end{eqnarray*}
Hence, by a telescopic series, we find that
\begin{eqnarray*}\gamma_m&\ge&\sum_{k=4mL}^{+\infty} p_k\,\frac{2\sin^2(m \pi/k)}{ (m \pi /k)^2} \\&\ge&\frac12\sum_{k=4mL}^{+\infty} p_k\\&=&\widetilde\sigma_m\\&\ge&|\sigma_m|,
\end{eqnarray*}as desired.

\paragraph{Solution to Exercise~\ref{ojdkfnvioewyr098765rewe67890iuhgvgyuuygfr43wdfgh8bft5xs}.}
Let~$f_\epsilon\in C^\infty(\R)$ be a periodic function such that
$$ \sup_{x\in\R}|f(x)-f_\epsilon(x)|<\frac\epsilon2,$$
see e.g.~\cite[Theorem~9.8]{MR3381284} for the construction of such a function.
Then, by Theorem~\ref{SMXC22} the $k$th Fourier coefficient~$\widehat{f_{\epsilon,k}}$ of~$f$ satisfies that
\begin{equation*} |\widehat {f_{\epsilon,k}}|\le \frac{\displaystyle\sup_{x\in[0,1)}|D^2 f_\epsilon(x)|}{|2k\pi |^2}\end{equation*}
and consequently
$$ \sum_{k\in\Z}|\widehat {f_{\epsilon,k}}|<+\infty.$$
We can thereby utilise Theorem~\ref{BASw} and find that, for~$N_\epsilon\in\N$ sufficiently large,
$$ \sup_{x\in\R}|f_\epsilon (x)-S_{N_\epsilon, f_\epsilon}(x)|<\frac\epsilon2.$$
Thus, $S_{N_\epsilon, f_\epsilon}$ is a trigonometric polynomial and satisfies
$$ \sup_{x\in\R}|f(x)-S_{N_\epsilon, f_\epsilon}(x)|<\epsilon,$$
as desired (for completeness, we mention that, in alternative to Theorems~\ref{SMXC22} and~\ref{BASw}, one can obtain the desired result also by using
Theorem~\ref{C1uni} and a covering argument).

Another solution of this exercise will be proposed in Exercise~\ref{ojdkfnvioewyr098765rewe67890iuhgvgyuuygfr43wdfgh8bft5xsBIS}.

\paragraph{Solution to Exercise~\ref{ojdkfnvioewyr098765rewe67890iuhgvgyuuygfr43wdfgh8bft5xs-EAN}.}
Given~$p\in[1,+\infty)$, $\epsilon>0$, and~$f\in L^p((0,1))$, we take a smooth function~$g$ such that~$\|f-g\|_{L^p((0,1))}<\epsilon$
(see e.g.~\cite[Theorem~9.6]{MR3381284}).

Then, by Exercise~\ref{ojdkfnvioewyr098765rewe67890iuhgvgyuuygfr43wdfgh8bft5xs}, we pick a trigonometric polynomial~$S$ such that~$\|g-S\|_{L^\infty((0,1))}<\epsilon$.

As a result, $\|f-S\|_{L^p((0,1))}\le\|f-g\|_{L^p((0,1))}+\|g-S\|_{L^p((0,1))}<\epsilon+\|g-S\|_{L^\infty((0,1))}<2\epsilon$.

\paragraph{Solution to Exercise~\ref{OAKSPiojlfS023a}.}
Let
$$ T_N(x):=\sum_{{k\in\Z}\atop{|k|\le N}} \mu_k\,e^{2\pi ikx}
=\sum_{\ell=0}^{\log_2N}\frac1{i\pi^{\ell}}\,\big(e^{2^{1+\ell}\pi ix}-e^{-2^{1+\ell}\pi ix}\big)
=\sum_{\ell=0}^{\log_2N}\frac{2}{\pi^{\ell}}\,\sin(2^{1+\ell}\pi x).$$
Since
$$ \sum_{\ell=0}^{+\infty}\frac{2}{\pi^{\ell}}<+\infty,$$
we have that~$T_N$ converges uniformly to the continuous function
$$ f(x):=\sum_{\ell=0}^{+\infty}\frac{2}{\pi^{\ell}}\,\sin(2^{1+\ell}\pi x)=
\sum_{\ell=0}^{+\infty}\frac1{i\pi^{\ell}}\,\big(e^{2^{1+\ell}\pi ix}-e^{-2^{1+\ell}\pi ix}\big)
.$$
Then, by Dominated Convergence Theorem,
$$\widehat f_k=\sum_{\ell=0}^{+\infty}\frac1{i\pi^{\ell}}\int_0^1\big(e^{2^{1+\ell}\pi ix}-e^{-2^{1+\ell}\pi ix}\big)e^{-2\pi ikx}\,dx
=\mu_k.$$
Hence, the answer to the first question is in the positive.

As for the second question, the idea is to use de l'H\^{o}pital's Rule to calculate the limit, take a derivative
inside the series and set~$x=0$. To justify this procedure, we argue as follows. We observe that, for all~$\tau\in\R$,
$$ 0\le1-\cos\tau\le\min\left\{2,\frac{\tau^2}2\right\},$$
hence, given~$\epsilon>0$, we have that, for all~$x\in(0,1)$,
\begin{eqnarray*}&&
\left|\frac{f(x)}{x}-\frac{4\pi^2}{\pi-2}\right|=
\left|\sum_{\ell=0}^{+\infty}\frac{2}{\pi^{\ell}}\cdot\frac{\sin(2^{1+\ell}\pi x)}{x}
-\frac{4\pi^2}{\pi-2}\right|\\&&\qquad=4\pi
\left|\sum_{\ell=0}^{+\infty}\frac{1}{2\pi^{1+\ell}}\cdot\frac{\sin(2^{1+\ell}\pi x)}{x}
-\frac{1}{1-(2/\pi)}\right|\\&&\qquad=4\pi
\left|\sum_{\ell=0}^{+\infty}\frac{1}{2\pi^{1+\ell}x}\int_0^{2^{1+\ell}\pi x}\cos\tau\,d\tau
-\sum_{\ell=0}^{+\infty}\left(\frac2\pi\right)^\ell\right|\\&&\qquad=4\pi
\left|\sum_{\ell=0}^{+\infty}\frac{1}{2\pi^{1+\ell}x}\int_0^{2^{1+\ell}\pi x}\big(\cos\tau-1\big)\,d\tau\right|\\&&\qquad\le4\pi
\sum_{\ell=0}^{+\infty}\frac{1}{2\pi^{1+\ell}x}\int_0^{2^{1+\ell}\pi x}\min\left\{2,\frac{\tau^2}2\right\}\,d\tau\\&&\qquad\le4\pi
\sum_{{\ell\in\N}\atop{2^\ell \pi x\le1}}\frac{1}{2\pi^{1+\ell}x}\int_0^{2^{1+\ell}\pi x}\frac{\tau^2}2\,d\tau+4\pi
\sum_{{\ell\in\N}\atop{2^\ell \pi x>1}}\frac{1}{\pi^{1+\ell}x}\int_0^{2^{1+\ell}\pi x}\,d\tau\\&&\qquad\le\frac{8\pi^3 x^2}3
\sum_{{\ell\in\N}\atop{\ell\le|\log_2(\pi x)|}} \left(\frac{8}\pi\right)^\ell 
+8\pi\sum_{{\ell\in\N}\atop{\ell>|\log_2(\pi x)|}}\left(\frac{2}{\pi}\right)^\ell\\&&\qquad\le24\pi^4 x^2\left(
\frac8\pi\right)^{|\log_2(\pi x)|}
+8\pi\sum_{{\ell\in\N}\atop{\ell>|\log_2(\pi x)|}}\left(\frac{2}{\pi}\right)^\ell\\&&\qquad=24\pi^{4-\log_2(8/\pi)} x^{2-\log_2(8/\pi)}
+8\pi\sum_{{\ell\in\N}\atop{\ell>|\log_2(\pi x)|}}\left(\frac{2}{\pi}\right)^\ell.
\end{eqnarray*}Therefore, since~$2-\log_2(8/\pi)>0$ and~$\frac{f(x)}x$ is an even function,
$$\lim_{x\to0}\frac{f(x)}{x}=\frac{4\pi^2}{\pi-2}.$$

\paragraph{Solution to Exercise~\ref{ZIMadmATnZ-1}.} Suppose, for a contradiction, that there exists~$x\in\R$
for which the series in~\eqref{LA:zkI1} converges. In particular, we have that
\begin{equation}\label{pjqldwm-1.0}\lim_{k\to+\infty}\cos(2\pi kx)=0.\end{equation}

We claim that also
\begin{equation}\label{pjqldwm-1.1}
\lim_{k\to+\infty}\sin(2\pi kx)=0.
\end{equation}
Indeed, suppose not. Then, by compactness, we can find a divergent sub-sequence~$h_j$ and~$\sigma\in[-1,1]\setminus\{0\}$ such that
$$ \lim_{j\to+\infty}\sin(2\pi h_jx)=\sigma.$$
For this reason, for all~$j\in\N$,
\begin{eqnarray*}0&=&\lim_{m\to+\infty} \cos(2\pi (h_m+h_j)x)\\&=&\lim_{m\to+\infty}\big(
\cos(2\pi h_mx)\cos(2\pi h_jx)-\sin(2\pi h_mx)\sin(2\pi h_jx)\big)\\&=&0-
\lim_{m\to+\infty}\sin(2\pi h_mx)\sin(2\pi h_jx)\\&=&-\sigma\sin(2\pi h_jx).
\end{eqnarray*}
We can now send~$j\to+\infty$ and conclude that~$0=-\sigma^2$, which is a contradiction,
whence the proof of~\eqref{pjqldwm-1.1} is complete.

It follows from~\eqref{pjqldwm-1.0} and~\eqref{pjqldwm-1.1} that there exists~$k_{x}\in\N$ large enough such that, if~$k\ge k_{x}$, then
$$ |\cos(2\pi kx)|<\frac{\sqrt2}2\qquad{\mbox{and}}\qquad|\sin(2\pi kx)|<\frac{\sqrt2}2.$$
This gives that~$2\pi k_x x$ must belong to two different regions of the trigonometric circle, which is impossible,
showing that the series in~\eqref{LA:zkI1} cannot converge.

\begin{figure}[h]
\includegraphics[height=3cm]{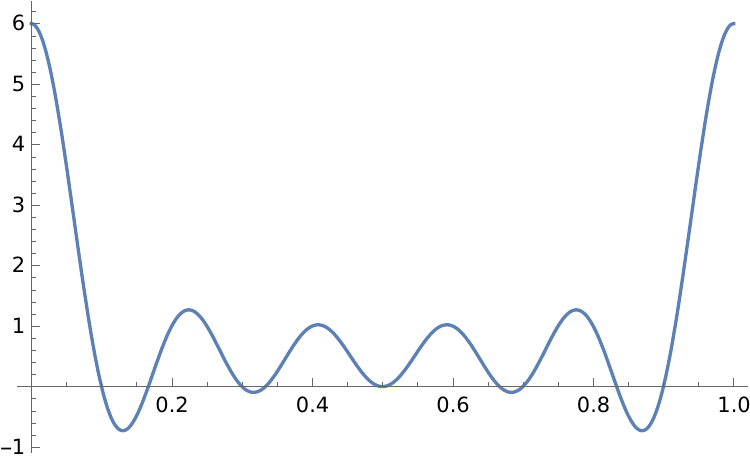}$\,\;$\includegraphics[height=3cm]{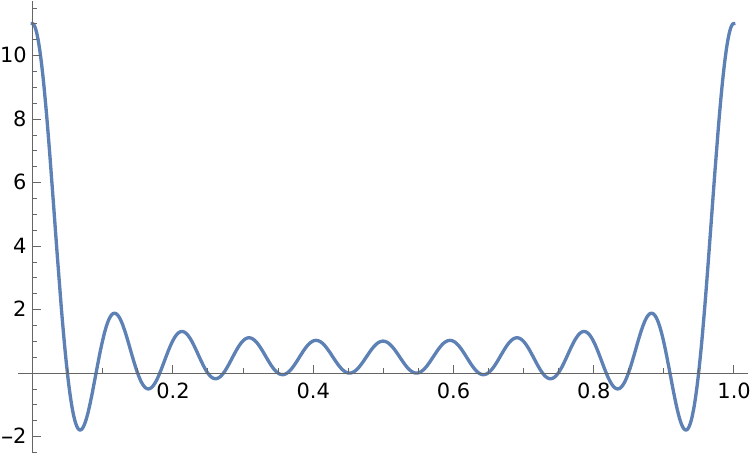}$\,\;$
\includegraphics[height=3cm]{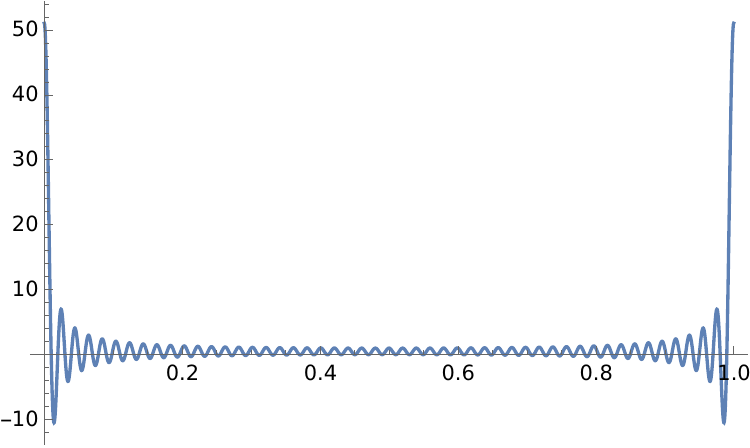}
\centering
\caption{Plot of~$\displaystyle\sum_{k=0}^N\cos(2\pi kx)$ with~$N\in\{5,10,50\}$.}\label{ikdnsE020e-U.01}
\end{figure}

Now we prove~\eqref{ZIMadmATnZ-1.01}. To this end, we use Exercise~\ref{PIPIO-010} to see that,
for all~$x\in(0,1)$,
$$\left|\sum_{j=0}^N\cos(2\pi jx)\right|=\left|\sum_{j=0}^N\Re\big(e^{2\pi ijx}\big)\right|=
\left|\Re\left(\sum_{{j\in\Z}\atop{0\le j\le N}}e^{2\pi ijx}\right)\right|\le\left|\sum_{{j\in\Z}\atop{0\le j\le N}}e^{2\pi ijx}\right|\le\frac1{|\sin(\pi x)|},$$
from which we obtain~\eqref{ZIMadmATnZ-1.01}.

See Figure~\ref{ikdnsE020e-U.01} to try to visualise the complications of this example.

See also~\cite[Section~7.6.5]{PICARDELLO} for further information on this kind of examples.

\paragraph{Solution to Exercise~\ref{ZIMadmATnZ-2}.}
This is a careful modification of Exercise~\ref{ZIMadmATnZ-1}.
Suppose, for a contradiction, that there exists~$x\in\left(0,\frac12\right)$
for which the series in~\eqref{LA:zkI2} converges. In particular, we have that
\begin{equation}\label{pjqldwm-1.0a}\lim_{k\to+\infty}\sin(2\pi kx)=0.\end{equation}

We claim that also
\begin{equation}\label{pjqldwm-1.1a}
\lim_{k\to+\infty}\cos(2\pi kx)=0.
\end{equation}
Indeed, suppose not. Then, by compactness, we can find a divergent sub-sequence~$h_j$ and~$\sigma\in[-1,1]\setminus\{0\}$ such that
$$ \lim_{j\to+\infty}\cos(2\pi h_jx)=\sigma.$$
For this reason, 
\begin{eqnarray*}0&=&\lim_{m\to+\infty} \sin(2\pi (h_m+1)x)\\&=&\lim_{m\to+\infty}\big(
\sin(2\pi h_mx)\cos(2\pi x)+\cos(2\pi h_mx)\sin(2\pi x)\big)\\&=&0+
\lim_{m\to+\infty}\cos(2\pi h_mx)\sin(2\pi x)\\&=&\sigma\sin(2\pi x).
\end{eqnarray*}
Since~$x\in\left(0,\frac12\right)$, we thus conclude that~$\sigma=0$, which is a contradiction,
whence the proof of~\eqref{pjqldwm-1.1a} is complete.

It follows from~\eqref{pjqldwm-1.0a} and~\eqref{pjqldwm-1.1a} that there exists~$k_{x}\in\N$ large enough such that, if~$k\ge k_{x}$, then
$$ |\cos(2\pi kx)|<\frac{\sqrt2}2\qquad{\mbox{and}}\qquad|\sin(2\pi kx)|<\frac{\sqrt2}2.$$
This gives that~$2\pi k_x x$ must belong to two different regions of the trigonometric circle, which is impossible,
showing that the series in~\eqref{LA:zkI2} cannot converge.

\begin{figure}[h]
\includegraphics[height=3cm]{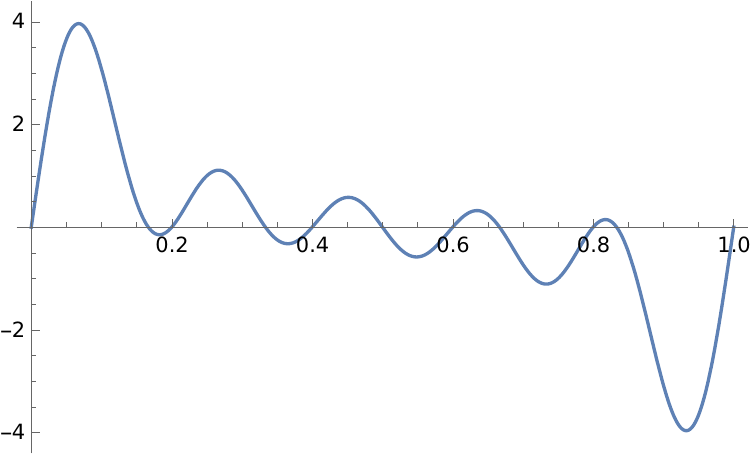}$\,\;$\includegraphics[height=3cm]{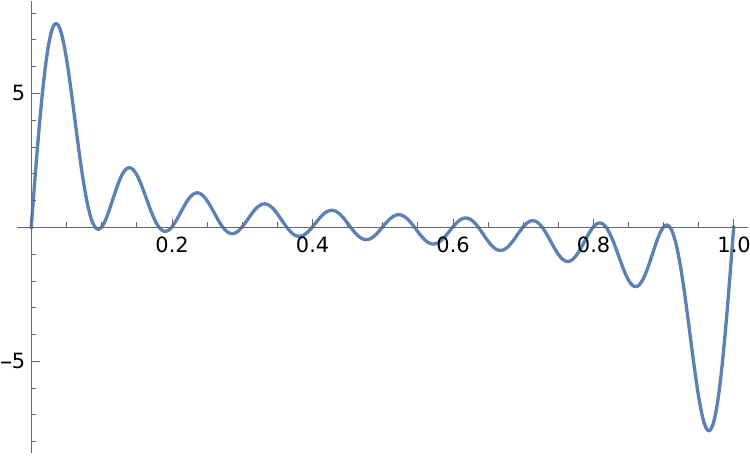}$\,\;$
\includegraphics[height=3cm]{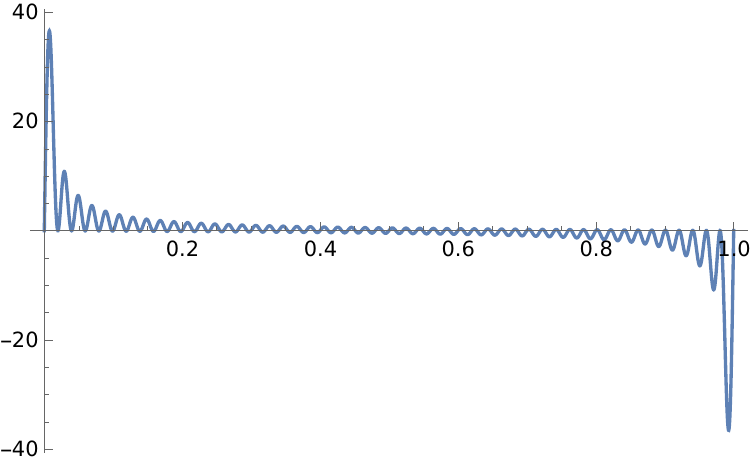}
\centering
\caption{Plot of~$\displaystyle\sum_{k=0}^N\sin(2\pi kx)$ with~$N\in\{5,10,50\}$.}\label{ikdnsE020e-U.02}
\end{figure}

Now we prove~\eqref{ZIMadmATnZ-1.02}. To this end, we observe that when~$x\in\Z$, the sum in~\eqref{ZIMadmATnZ-1.02} vanishes, so we can assume that~$x\not\in\Z$, and therefore~$\sin(\pi x)\ne0$.
Hence, we can use Exercise~\ref{PIPIO-010} to see that
$$\left|\sum_{j=0}^N\sin(2\pi jx)\right|=\left|\sum_{j=0}^N\Im\big(e^{2\pi ijx}\big)\right|=
\left|\Im\left(\sum_{{j\in\Z}\atop{0\le j\le N}}e^{2\pi ijx}\right)\right|\le\left|\sum_{{j\in\Z}\atop{0\le j\le N}}e^{2\pi ijx}\right|\le\frac1{|\sin(\pi x)|},$$
from which we obtain~\eqref{ZIMadmATnZ-1.02}.

See Figure~\ref{ikdnsE020e-U.02} to try to visualise the complications of this example.

\paragraph{Solution to Exercise~\ref{ZIMadmATnZ-3}.}
The answer is negative for both cases, because the corresponding Fourier coefficients are not infinitesimal,
against the Riemann-Lebesgue Lemma (see Theorem~\ref{RLjoqwskcdc}).

\section{Solutions to selected exercises of Section~\ref{COL2}}

\paragraph{Solution to Exercise~\ref{PARSCA}.} We observe that~$h:=f+ g\in L^2((0,1))$.
Hence, in light of~\eqref{L2THM.0-02} in Theorem~\ref{THCOL2FB}(i) (recall also Exercise~\ref{LINC-bisco}) we have that
\begin{eqnarray*}&&
\sum_{k\in\Z}\big(|\widehat f_k|^2+| \widehat g_k|^2+2 \widehat f_k\,\overline{\widehat g_k}\big)=
\sum_{k\in\Z}|\widehat f_k+ \widehat g_k|^2=
\sum_{k\in\Z} |\widehat h_k|^2 \\&&\qquad= \|h\|^2_{L^2((0,1))}=\|f+ g\|^2_{L^2((0,1))}\\&&\qquad=\|f\|^2_{L^2((0,1))}
+\| g\|^2_{L^2((0,1))}+2\int_0^1 f(x)\,g(x)\,dx\\&&\qquad=
\sum_{k\in\Z}\big(|\widehat f_k|^2+| \widehat g_k|^2\big)+2\int_0^1 f(x)\,g(x)\,dx,
\end{eqnarray*}
from which we obtain the desired result.

\paragraph{Solution to Exercise~\ref{ISOMel2}.} Given~$f\in L^2((0,1))$, we define~$T(f)$ to be the sequence~$\{\widehat f_k\}_{k\in\Z}$.
We notice that this is a linear map, thanks to Exercise~\ref{LINC-bisco}, and that~$T(f)\in \ell^2$,
due to~\eqref{fasv} and Bessel's Inequality~\eqref{AJSa}.

Moreover, the claim in~\eqref{L2THM.0-02-makdefg} holds true, as a consequence of~\eqref{L2THM.0-02}.

We also stress that~$T$ is injective, owing to Theorem~\ref{UNIQ}, and it is surjective, thanks to Theorem~\ref{THCOL2FB}(ii),
completing the proof of the desired result.

The observation contained in this exercise is structurally important, because it essentially states a general fact, according to which all separable infinite dimensional Hilbert spaces are linearly isometric with one another,
see e.g.~\cite[Theorems~8.33 and~8.35]{MR3381284} for more information.

\paragraph{Solution to Exercise~\ref{OJSNILCESNpfa.sdwpqoed-23wedf-87a}.}
We define
$$ T_N(x):=\int_0^x S_{N,f}(t)\,dt=
\widehat f_0\,x+\sum_{{k\in\Z}\atop{0<|k|\le N}} \frac{\widehat f_k}{2\pi ik}\,e^{2\pi ikx}$$
and we observe that
\begin{eqnarray*}&&
\sup_{x\in(0,1)}|F(x)-T_N(x)|=
\sup_{x\in(0,1)}\left|\int_0^x f(t) -S_{N,f}(t)\,dt\right|\\&&\qquad\le
\int_0^1| f(t) -S_{N,f}(t)|\,dt\le
\sqrt{\int_0^1| f(t) -S_{N,f}(t)|^2\,dt},
\end{eqnarray*}
which is infinitesimal, by virtue of~\eqref{L2THM.0-01}.

\section{Solutions to selected exercises of Section~\ref{COLP}}

\paragraph{Solution to Exercise~\ref{0pXe-244kf.1EX-PRE.0p23rl}.} The claim for~$S$ is a consequence of~\eqref{cjamsdperrea} and actually one can write
$$ S(x)=\frac{a_0}2+\sum_{k=1}^M a_k\, \cos(2\pi k x)
+b_k\,\sin(2\pi k x),$$
with~$a_k:=2\Re(c_k)$ and~$b_k:=-2\Im(c_k)$.

As for~${S}^\Diamond$, we recall~\eqref{cjamsdperrea} and~\eqref{0pXe-244kf.1EX-PRE.pe-djenfEF} to see that
\begin{eqnarray*}{S}^\Diamond(x)&=&-i\sum_{k=1}^M c_k \,e^{2\pi i k x}
+i\sum_{k=1}^M c_{-k} \,e^{-2\pi i k x}
\\&=&-i\sum_{k=1}^M c_k \,e^{2\pi i k x}
+i\sum_{k=1}^M \overline{c_{k}} \,e^{-2\pi i k x}\\
&=&-i\sum_{k=1}^M c_k \,\big(\cos(2\pi k x)+i\sin(2\pi k x)\big)
+i\sum_{k=1}^M \overline{c_{k}} \,\big(\cos(2\pi k x)-i\sin(2\pi k x)\big)\\
&=&i\sum_{k=1}^M \big(\overline{c_{k}}-c_k\big) \, \cos(2\pi k x)
+\sum_{k=1}^M \big(c_{k} +\overline{c_{k}}\big)\,\sin(2\pi k x)
\\&=&2\sum_{k=1}^M \Im(c_k) \, \cos(2\pi k x)
+2\sum_{k=1}^M \Re(c_{k})\,\sin(2\pi k x),
\end{eqnarray*}
which yields that~${S}^\Diamond$ takes values in the reals.

For completeness, let us remark that,
comparing with~\eqref{jasmx23er}, we have just found that
$$ {S}^\Diamond(x)=
\sum_{k=1}^M a_k\,\sin(2\pi k x)- b_k \, \cos(2\pi k x) ,$$
which can be considered as the trigonometric form of~${S}^\Diamond$.

\paragraph{Solution to Exercise~\ref{0pXe-244kf.1EX-PRE.0p23rlmanyOnt}.}
For instance, if~$S(x)=2\cos(2\pi x)=e^{2\pi ix}+e^{-2\pi ix}$, we have that~$P^+_{S}(x)=e^{2\pi ix}$ and~$P^-_{S}(x)=e^{-2\pi ix}$, which are complex-valued functions.

\paragraph{Solution to Exercise~\ref{0pXe-244kf.1EX-PRE.0p23rl.oqjsdlwnfvMSAS01maspe}.}
We will use the shift of indices~$h:=k+n$. In the notation of Exercise~\ref{0pXe-244kf.1EX-PRE.pe}, we have that
\begin{eqnarray*}
S^{(n)}(x)=\sum_{{k\in\Z}\atop{|k|\le M}} c_{k} \,e^{2\pi i( k+n) x}=\sum_{{k\in\Z}\atop{n-M\le h\le n+ M}} c_{h-n} \,e^{2\pi ih x}
=\sum_{{h\in\Z}\atop{|h|\le M_n}}\widetilde c_h\,e^{2\pi ihx},
\end{eqnarray*}
where
\begin{eqnarray*}&&\widetilde c_h:=\begin{dcases}c_{h-n}&{\mbox{ if }}h\in[n-M,n+M],\\ 0&{\mbox{ otherwise}}
\end{dcases}\\
{\mbox{and }}&&M_n:=\max\big\{ |n-M|,\,|n+M|\big\}.
\end{eqnarray*} 
Thus, by~\eqref{gbsdfnigbfa.S01maspe},
$$ P^+_{S^{(n)}}(x)=\sum_{h=1}^{M_n}\widetilde c_h\,e^{2\pi ihx}
=\sum_{{h\in\Z}\atop{1\le h\le n+ M}} c_{h-n} \,e^{2\pi ih x}=
\sum_{{k\in\Z}\atop{1-n\le k\le M}} c_k \,e^{2\pi i(k+n) x}.$$
From this, one obtains the desired result.

\paragraph{Solution to Exercise~\ref{0pXe-244kf.1EX-PRE.0p23rl.oqjsdlwnfvMSAS01maspe.bismdc}.} We exploit Exercise~\ref{0pXe-244kf.1EX-PRE.0p23rl.oqjsdlwnfvMSAS01maspe} with~$n:=N+1$ and~$n:=-N$, thus obtaining that
\begin{eqnarray*}&&
e^{-2\pi i(N+1)x}P^+_{S^{(N+1)}}(x)-e^{2\pi iNx}P^+_{S^{(-N)}}(x)=
\sum_{{k\in\Z}\atop{-N\le k\le M}} c_{k} \,e^{2\pi i k x}-\sum_{{k\in\Z}\atop{1+N\le k\le M}} c_{k} \,e^{2\pi i k x}\\&&\qquad=\sum_{{k\in\Z}\atop{-N\le k\le N}} c_{k} \,e^{2\pi i k x}
,\end{eqnarray*}
as desired.

\paragraph{Solution to Exercise~\ref{0pXe-244kf.1EX-PRE}.} By virtue of Exercise~\ref{0pXe-244kf.1EX-PRE.pe},
\begin{eqnarray*}&& \int_0^1 \big( S(x)+i{S}^\Diamond(x)\big)^{\ell}\,dx=
\int_0^1 \big( 2P^+_{S}(x)+c_0\big)^{\ell}\,dx\\&&\qquad=2^\ell\int_0^1 \big( P^+_{S}(x)\big)^{\ell}\,dx\\
&&\qquad=2^\ell
\sum_{{1\le k_1,\dots,k_\ell\le M}} c_{k_1}\dots c_{k_\ell} \,\int_0^1e^{2\pi i (k_1+\dots+k_\ell) x}\,dx\\&&\qquad=
2^\ell
\sum_{{1\le k_1,\dots,k_\ell\le M}} \frac{c_{k_1}\dots c_{k_\ell} \,\big(e^{2\pi i (k_1+\dots+k_\ell)}-1\big)}{2\pi i (k_1+\dots+k_\ell)}\\&&\qquad=0,
\end{eqnarray*}as desired.

\paragraph{Solution to Exercise~\ref{0pXe-244kf.1EX-PRE-012}.}
By Exercise~\ref{0pXe-244kf.1EX-PRE},
\begin{eqnarray*} &&0=\int_0^1 \big( S(x)+i{S}^\Diamond(x)\big)^{2\ell}\,dx=\int_0^1
\sum_{j=0}^{2\ell} i^{2\ell-j}\binom{2\ell}{j} \big(S(x)\big)^{j}\big({S}^\Diamond(x)\big)^{2\ell-j}\,dx.
\end{eqnarray*}
Taking the term with~$j=0$ to the left-hand side, we find that
$$ (-1)^{\ell+1}\int_0^1\big({S}^\Diamond(x)\big)^{2\ell}\,dx=
\sum_{j=1}^{2\ell} i^{2\ell-j}\binom{2\ell}{j} \int_0^1\big(S(x)\big)^{j}\big({S}^\Diamond(x)\big)^{2\ell-j}\,dx
$$
and consequently, using H\"older's Inequality with exponents~$\frac{2\ell}{j}$ and~$\frac{2\ell}{2\ell-j}$,
\begin{eqnarray*}
\int_0^1\big({S}^\Diamond(x)\big)^{2\ell}\,dx&\le&
\sum_{j=1}^{2\ell} \binom{2\ell}{j} \int_0^1\big|S(x)\big|^{j}\big|{S}^\Diamond(x)\big|^{2\ell-j}\,dx
\\&\le&\sum_{j=1}^{2\ell} \binom{2\ell}{j}\left( \int_0^1\big(S(x)\big)^{2\ell}\,dx\right)^{\frac{j}{2\ell}}\left(\int_0^1
\big({S}^\Diamond(x)\big)^{2\ell}\,dx\right)^{\frac{2\ell-j}{2\ell}}
,\end{eqnarray*}
yielding the desired result.

\paragraph{Solution to Exercise~\ref{0pXe-244kf.1EX-PRE-012.b}.}
If~$X<1$ we are done, hence we can suppose that~$X\ge1$. Thus, for all~$j\in\{1,\dots,2\ell\}$, we have that~$X^{2\ell-j}\le X^{2\ell-1}$ and therefore
$$ X^{2\ell}\le\sum_{j=1}^{2\ell} \binom{2\ell}{j} X^{2\ell-j}\le\sum_{j=1}^{2\ell} \binom{2\ell}{j}X^{2\ell-1}\le
2\ell\,M_\ell\, X^{2\ell-1},
$$
where
$$ M_\ell:= \max_{j\in\{1,\dots,2\ell\}}\binom{2\ell}{j}.$$
This gives that~$X\le2\ell\,M_\ell$, as desired.

\paragraph{Solution to Exercise~\ref{0pXe-244kf.1EX-PRE-012.c}.} We observe that
\begin{equation}\label{cjamsdperrea-24}\begin{split}&
{\mbox{it suffices to prove~\eqref{0pXe-244kf.1EX-PRE-012.d} under the additional assumption}}\\
&{\mbox{that the trigonometric polynomial~$S$ as in~\eqref{cjamsdperrea-23}
satisfies~$c_0=0$.}}\end{split}
\end{equation}
Indeed, one can define~$T:=S-c_0$ and observe that~${T}^\Diamond={S}^\Diamond$. If~\eqref{0pXe-244kf.1EX-PRE-012.d} is valid under the additional assumption above, then we know that
\begin{equation}\label{cjamsdperrea-25}
\|{S}^\Diamond\|_{L^{2\ell}((0,1))}=\|{T}^\Diamond\|_{L^{2\ell}((0,1))}\le C\|T\|_{L^{2\ell}((0,1))}
\le C\big(\|S\|_{L^{2\ell}((0,1))}+|c_0|\big).
\end{equation}

Moreover,
$$ c_0=\int_0^1 S(x)\,dx$$
and therefore, by H\"older's Inequality,
$$ |c_0|\le\int_0^1|S(x)|\,dx\le \|S\|_{L^{2\ell}((0,1))}.$$

Inserting this information into~\eqref{cjamsdperrea-25}, it follows that$$\|{S}^\Diamond\|_{L^{2\ell}((0,1))}
\le 2C\|S\|_{L^{2\ell}((0,1))},$$ hence~\eqref{0pXe-244kf.1EX-PRE-012.d} for~$S$ as well (up to renaming constants).
This proves~\eqref{cjamsdperrea-24}.

Now, on account of~\eqref{cjamsdperrea-24}, we can suppose that~$c_0=0$ and this will allow us to employ
Exercise~\ref{0pXe-244kf.1EX-PRE-012}.

Let~$\epsilon>0$ and
\begin{equation} \label{0pXe-244kf.1EX-PRE-012.e}
X:=\frac{\|{S}^\Diamond\|_{L^{2\ell}((0,1))}}{\|S\|_{L^{2\ell}((0,1))}+\epsilon} .\end{equation}

On grounds of Exercise~\ref{0pXe-244kf.1EX-PRE-012}, we know that
\begin{eqnarray*} X^{2\ell}&=&\frac{\|{S}^\Diamond\|_{L^{2\ell}((0,1))}^{2\ell}}{\big(\|S\|_{L^{2\ell}((0,1))}+\epsilon\big)^{2\ell}}
\\&\le&
\sum_{j=1}^{2\ell} \binom{2\ell}{j}\,\frac{\|S\|_{L^{2\ell}((0,1))}^{j}\,\|{S}^\Diamond\|_{L^{2\ell}((0,1))}^{2\ell-j}}{\big(\|S\|_{L^{2\ell}((0,1))}+\epsilon\big)^{2\ell}}\\&\le&
\sum_{j=1}^{2\ell} \binom{2\ell}{j}\,\frac{\|{S}^\Diamond\|_{L^{2\ell}((0,1))}^{2\ell-j}}{\big(\|S\|_{L^{2\ell}((0,1))}+\epsilon\big)^{2\ell-j}}\\&=&\sum_{j=1}^{2\ell}\binom{2\ell}{j} X^{2\ell-j}.
\end{eqnarray*}

Hence, in light of Exercise~\ref{0pXe-244kf.1EX-PRE-012.b},
$$ \frac{\|{S}^\Diamond\|_{L^{2\ell}((0,1))}}{\|S\|_{L^{2\ell}((0,1))}+\epsilon}=X\le C$$
and accordingly
$$ \|{S}^\Diamond\|_{L^{2\ell}((0,1))}\le C\Big(\|S\|_{L^{2\ell}((0,1))}+\epsilon\Big).$$
The desired result now follows by sending~$\epsilon\searrow0$.

As a side remark, let us mention that the role of~$\epsilon$ in the above computation was merely ancillary, to avoid possible divisions by zero in~\eqref{0pXe-244kf.1EX-PRE-012.e}; one could have avoided this small complication by treating separately the (trivial)
case in which~$\|S\|_{L^{2\ell}((0,1))}=0$.

\paragraph{Solution to Exercise~\ref{0pXe-244kf.1EX-SK}.}
This is a nice example of interpolation methods. First of all, we prove~\eqref{0pXe-244kf.1EX-PRE-012.d-96}
when~$p\in[2,+\infty)$. Thanks to
Exercise~\ref{0pXe-244kf.1EX-PRE-012.c}, we already know the desired result for~$p\in2\N\cap[2,+\infty)$, therefore
we can assume that~$p\in(2\ell,2\ell+2)$ for some~$\ell\in\N\cap[1,+\infty)$.

Now given two spaces~$X$ and~$Y$, we 
define~$X+Y:=\{x+y$, $x\in X$, $y\in Y\}$, and a linear operator
\begin{equation}\label{men01j45rf2}
T:L^{2\ell+2}((0,1))+L^{2\ell}((0,1))\longrightarrow L^{2\ell+2}((0,1))+L^{2\ell}((0,1))\end{equation} such that~$TS=
{S}^\Diamond$ for all trigonometric polynomials~$S$. To perform this operation, let~$f\in L^{2\ell+2}((0,1))+L^{2\ell}((0,1))$,
i.e., let~$f=f_1+f_2$ with~$f_1\in L^{2\ell+2}((0,1))$ and~$f_2\in L^{2\ell}((0,1))$. By the density of trigonometric polynomials
(see Exercise~\ref{ojdkfnvioewyr098765rewe67890iuhgvgyuuygfr43wdfgh8bft5xs-EAN}), we can find sequences of
trigonometric polynomials~$\{S_{1,j}\}_{j\in\N}$ and~$\{S_{2,j}\}_{j\in\N}$ such that
\begin{equation}\label{men01j45rf} \lim_{j\to+\infty}\| f_1-S_{1,j}\|_{L^{2\ell+2}((0,1))}=0\qquad{\mbox{and}}\qquad
\lim_{j\to+\infty}\| f_2-S_{2,j}\|_{L^{2\ell}((0,1))}=0.\end{equation}
We observe that, by Exercise~\ref{0pXe-244kf.1EX-PRE-012.c}, for all~$m$, $j\in\N$,
\begin{eqnarray*}&& \|{S}^\Diamond_{1,j}-{S}^\Diamond_{1,m}\|_{L^{2\ell+2}((0,1))}\le C\|S_{1,j}-S_{1,m}\|_{L^{2\ell+2}((0,1))}\\{\mbox{and}}\qquad&&
\|{S}^\Diamond_{2,j}-{S}^\Diamond_{2,m}\|_{L^{2\ell}((0,1))}\le C\|S_{2,j}-S_{2,m}\|_{L^{2\ell}((0,1))}.\end{eqnarray*}

This and~\eqref{men01j45rf} yield that~${S}^\Diamond_{1,j}$ is a Cauchy sequence in~$L^{2\ell+2}((0,1))$
and~${S}^\Diamond_{2,j}$ is a Cauchy sequence in~$L^{2\ell}((0,1))$ and accordingly we can define~$Tf_1$ as the limit
in~$L^{2\ell+2}((0,1))$ of~${S}^\Diamond_{1,j}$ and~$Tf_2$ as the limit
in~$L^{2\ell}((0,1))$ of~${S}^\Diamond_{2,j}$, as~$j\to+\infty$.
This provides the construction\footnote{Let us also stress that the definition of~$T$ is independent of
the approximating trigonometric polynomials chosen. For instance,
if one chose another sequence of
trigonometric polynomials~$\{\widetilde S_{1,j}\}_{j\in\N}$ such that
$$\lim_{j\to+\infty}\| f_1-\widetilde S_{1,j}\|_{L^{2\ell+2}((0,1))}=0 ,$$
then it would have followed from Exercise~\ref{0pXe-244kf.1EX-PRE-012.c} that
\begin{eqnarray*}&& \lim_{j\to+\infty}\| S_{1,j}^\Diamond-\widetilde S_{1,j}^\Diamond\|_{L^{2\ell+2}((0,1))}\le
C\lim_{j\to+\infty}\| S_{1,j}-\widetilde S_{1,j}\|_{L^{2\ell+2}((0,1))}\\&&\qquad\le
C\lim_{j\to+\infty}\big(\| S_{1,j}-f_1\|_{L^{2\ell+2}((0,1))}+\|f_1-\widetilde S_{1,j}\|_{L^{2\ell+2}((0,1))}\big)=0,\end{eqnarray*}
thus providing the same definition for~$Tf_1$.} of the operator in~\eqref{men01j45rf2}.

Moreover, again by Exercise~\ref{0pXe-244kf.1EX-PRE-012.c},
$$\| Tf_1\|_{L^{2\ell+2}((0,1))}=\lim_{j\to+\infty}
\| S_{1,j}^\Diamond\|_{L^{2\ell+2}((0,1))}\le C\lim_{j\to+\infty}
\| S_{1,j}\|_{L^{2\ell+2}((0,1))}=C\| f_1\|_{L^{2\ell+2}((0,1))}$$
and similarly
$$\| Tf_2\|_{L^{2\ell}((0,1))}\le C\| f_2\|_{L^{2\ell}((0,1))}.$$

Accordingly, we can apply the Riesz Interpolation Theorem
(see e.g.~\cite[Theorem~2.1 on pages~52--53]{MR2827930}, used here with~$q:=p$, $q_0:=p_0:=2\ell$, $q_1:=p_1:=2\ell+2$ and~$t:=\frac{(\ell+1)(p-2\ell)}{p}$) and deduce that, for all~$f\in L^p((0,1))$,
$$ \| Tf\|_{L^p((0,1))}\le C\| f\|_{L^p((0,1))},$$
up to renaming~$C$.

This completes the proof of~\eqref{0pXe-244kf.1EX-PRE-012.d-96}
when~$p\in[2,+\infty)$.

We now prove~\eqref{0pXe-244kf.1EX-PRE-012.d-96}
when~$p\in(1,2)$. In this case, we consider the dual exponent~$q$ of~$p$, namely we define~$q:=\frac{p}{p-1}$ and we remark that 
\begin{equation}\label{0pXe-244kf.1EX-PRE-012.d-962-i4k8gBN4}
q>2.
\end{equation}

Moreover, we can compute the norm of a function in~$L^p((0,1))$ by duality (see e.g.~\cite[Theorem~10.44]{MR3381284}), that is
$$ \|f\|_{L^p((0,1))}=\sup_{{g\in L^q((0,1))}\atop{\|g\|_{L^q((0,1))}=1}}\int_0^1 f(x)\,{g(x)}\,dx.$$

In particular, given~$\epsilon>0$, if~$S$ is a trigonometric polynomial as in~\eqref{cjamsdperrea-23},
we take~$g_{S,\epsilon}\in L^q((0,1))$, with~$\|g_{S,\epsilon}\|_{L^q((0,1))}=1$, such that
\begin{equation}\label{ojdkfnvioewyr098765rewe67890iuhgvgyuuygfr432134135tymTGbdHwdfgh8bft5xs-EAN}
\|{S}^\Diamond\|_{L^p((0,1))}\le\epsilon+\int_0^1 {S}^\Diamond(x)\,{g_{S,\epsilon}(x)}\,dx.
\end{equation}
We use again the density of trigonometric polynomials
(see Exercise~\ref{ojdkfnvioewyr098765rewe67890iuhgvgyuuygfr43wdfgh8bft5xs-EAN}) to find a trigonometric polynomial~$S_\epsilon$ such that~$\|g_{S,\epsilon}-S_\epsilon\|_{L^q((0,1))}\le\epsilon$.

In this way, H\"older's Inequality returns that
\begin{eqnarray*}
\int_0^1 {S}^\Diamond(x)\,\big( g_{S,\epsilon}(x)-S_\epsilon(x)\big)\,dx\le
\|{S}^\Diamond\|_{L^p((0,1))}\|g_{S,\epsilon}-S_\epsilon\|_{L^q((0,1))}\le\epsilon\|{S}^\Diamond\|_{L^p((0,1))},
\end{eqnarray*}
which, in combination with~\eqref{ojdkfnvioewyr098765rewe67890iuhgvgyuuygfr432134135tymTGbdHwdfgh8bft5xs-EAN}, produces that
\begin{equation}\label{ojdkfnvioewyr098765rewe67890iuhgvgyuuygfr432134135tymTGbdHwdfgh8bft5xs-EAN-mjk-1}
(1-\epsilon)\|{S}^\Diamond\|_{L^p((0,1))}\le\epsilon+\int_0^1 {S}^\Diamond(x)\, S_\epsilon(x)\,dx.
\end{equation}

Now we write~$S_\epsilon$ in the explicit form
$$ S_\epsilon(x)=\sum_{{k\in\Z}\atop{|k|\le N_\epsilon}} c_{\epsilon,k} \,e^{2\pi i k x},$$
for suitable coefficients~$c_{\epsilon,k}$ satisfying the reality condition~$c_{\epsilon,-k}=\overline{c_{\epsilon,k}}$,
and some~$N_\epsilon\in\N$ (without loss of generality, up to putting some zero coefficients into the sums,
we can also suppose that~$N_\epsilon\ge M$, where~$M$ is as in~\eqref{cjamsdperrea-23}).

Thus, using Exercise~\ref{PARSCA}, the definition of conjugacy in~\eqref{0pXe-244kf.1EX-PRE.pe-djenfEF}, and the fact that
trigonometric polynomials are in~$L^2((0,1))$, we see that
\begin{eqnarray*}&&
\int_0^1 {S}^\Diamond(x)\, S_\epsilon(x)\,dx=
-i\sum_{{k\in\Z}\atop{|k|\le M}} {\operatorname{sign}}(k)\, c_k\,c_{\epsilon,-k}
=-\int_0^1 {S}(x)\, S_\epsilon^\Diamond(x)\,dx.
\end{eqnarray*}
Consequently,
\begin{equation}\label{PARSCAp-kjedmf2c9yT}
\int_0^1 {S}^\Diamond(x)\, S_\epsilon(x)\,dx\le\|{S}\|_{L^p((0,1))}\|{S}_\epsilon^\Diamond\|_{L^q((0,1))}.
\end{equation}

Also, recalling~\eqref{0pXe-244kf.1EX-PRE-012.d-962-i4k8gBN4}, we can apply~\eqref{0pXe-244kf.1EX-PRE-012.d-96} for the exponent~$q$ and find that
$$ \|{S}_\epsilon^\Diamond\|_{L^q((0,1))}\le C\|{S}_\epsilon\|_{L^q((0,1))}\le
C\big(\|g_{S,\epsilon}\|_{L^q((0,1))}+\|g_{S,\epsilon}-S_\epsilon\|_{L^q((0,1))}\big)\le C(1+\epsilon).
$$

In the wake of this and~\eqref{PARSCAp-kjedmf2c9yT}, we see that
$$\int_0^1 {S}^\Diamond(x)\, S_\epsilon(x)\,dx\le C(1+\epsilon)\|{S}\|_{L^p((0,1))}.$$
Thus, recalling~\eqref{ojdkfnvioewyr098765rewe67890iuhgvgyuuygfr432134135tymTGbdHwdfgh8bft5xs-EAN-mjk-1},
$$ (1-\epsilon)\|{S}^\Diamond\|_{L^p((0,1))}\le\epsilon+C(1+\epsilon)\|{S}\|_{L^p((0,1))}.$$
We can now send~$\epsilon\searrow0$ and conclude that~$\|{S}^\Diamond\|_{L^p((0,1))}\le C\|{S}\|_{L^p((0,1))}$, as desired.

\paragraph{Solution to Exercise~\ref{0pXe-244kf.1EX-SK-ilcoM}.} 
One uses~\eqref{0pXe-244kf.1EX-PRE-012.d-96} to see that
\begin{eqnarray*}
\| {S}^\Diamond\|_{L^p((0,1),\C)}&=&\| {S}^\Diamond_1+i{S}^\Diamond_2\|_{L^p((0,1),\C)}\\&
\le& \| {S}^\Diamond_1\|_{L^p((0,1))}+\|{S}^\Diamond_2\|_{L^p((0,1))}\\&\le&
C\,\Big(\|S_1\|_{L^p((0,1))}+\|S_2\|_{L^p((0,1))}\Big)\\&\le&C\,\Big(\|S_1\|_{L^p((0,1))}^p+\|S_2\|_{L^p((0,1))}^p\Big)^{\frac1p}
\\&=&C\left(\int_0^1 \big( |S_1(x)|^p+|S_2(x)|^p\big)\,dx
\right)^{\frac1p}\\&\le&C\left(\int_0^1 \big( |S(x)|^p+|S(x)|^p\big)\,dx
\right)^{\frac1p}\\&\le&C\,\| {S}\|_{L^p((0,1))},
\end{eqnarray*}
up to renaming~$C$ line after line.

\paragraph{Solution to Exercise~\ref{0pXe-244kf.1EX-SK0-1o23lk-0923}.} First of all, we show that
\begin{equation}\label{0pXe-244kf.1EX-PRE-012.d-96-0923k2}\begin{split}&
{\mbox{it suffices to check~\eqref{0pXe-244kf.1EX-PRE-012.d-96-0923k} under the additional assumption}}\\ &{\mbox{that~$f$
is a trigonometric polynomial.}}\end{split}\end{equation}
To this end, let~$\epsilon>0$. By the density of trigonometric polynomials (see Exercise~\ref{ojdkfnvioewyr098765rewe67890iuhgvgyuuygfr43wdfgh8bft5xs-EAN}), we can find a trigonometric polynomial~$S_\epsilon$ such that~$\|f-S_\epsilon\|_{L^p((0,1))}\le\epsilon$.
Also, if~\eqref{0pXe-244kf.1EX-PRE-012.d-96-0923k} is valid for trigonometric polynomials, we have that
$$ \| S_{N,S_\epsilon}\|_{L^p((0,1))}\le C\,\|S_\epsilon\|_{L^p((0,1))}$$
and therefore
\begin{equation}\label{0pXe-244kf.1EX-PRE-012.d-96-0923k2-0o2eX} \| S_{N,S_\epsilon}\|_{L^p((0,1))}\le C\,\big(\|f\|_{L^p((0,1))}+\|f-S_\epsilon\|_{L^p((0,1))}\big)\le
C\,\big(\|f\|_{L^p((0,1))}+\epsilon\big).\end{equation}

Besides, if~$h\in L^p((0,1))$ is periodic of period~$1$,
$$ |\widehat h_k|\le\int_0^1|h(x)|\,dx\le\|h\|_{L^p((0,1))}$$
and therefore
$$ |S_{N,h}(x)|\le\sum_{{k\in\Z}\atop{|k|\le N}}|\widehat h_k|\le (2N+1)\|h\|_{L^p((0,1))}.$$
In particular,
$$ \|S_{N,h}\|_{L^p((0,1))}\le (2N+1)\|h\|_{L^p((0,1))}.$$
Applying this relation with~$h:=f-S_\epsilon$ we find that
$$ \|S_{N,f}-S_{N,S_\epsilon}\|_{L^p((0,1))}\le (2N+1)\epsilon.$$

The latter inequality and~\eqref{0pXe-244kf.1EX-PRE-012.d-96-0923k2-0o2eX} yield that
\begin{eqnarray*}
\| S_{N,f}\|_{L^p((0,1))}\le \|S_{N,S_\epsilon}\|_{L^p((0,1))}+\|S_{N,f}-S_{N,S_\epsilon}\|_{L^p((0,1))}\le
C\,\big(\|f\|_{L^p((0,1))}+\epsilon\big)+(2N+1)\epsilon.
\end{eqnarray*}
We can now send~$\epsilon\searrow0$ and conclude that~$
\| S_{N,f}\|_{L^p((0,1))}\le C\,\big(\|f\|_{L^p((0,1))}$, which is~\eqref{0pXe-244kf.1EX-PRE-012.d-96-0923k}.

The proof of~\eqref{0pXe-244kf.1EX-PRE-012.d-96-0923k2} is thereby complete and therefore we now
assume, without loss of generality, that~$f$ is a trigonometric polynomial.
Actually, in general, the difficult case to treat is when the degree of~$f$ is strictly greater than~$N$,
otherwise~$S_{N,f}$ would coincide with~$f$ and~\eqref{0pXe-244kf.1EX-PRE-012.d-96-0923k} would obviously hold true,
therefore we can suppose that
\begin{equation}\label{0pXe-244kf.1EX-PRE-012.d-96-0923k2-0o2eX-fg01} f(x)=\sum_{{k\in\Z}\atop{|k|\le M}}c_k\,e^{2\pi ikx},\end{equation}
for suitable coefficients~$c_k$ and~$M\in \N\cap[N+1,+\infty)$.

For all~$n\in\Z$, we let
$$ S^{(n)}(x):=e^{2\pi inx} f(x)$$
and, owing to Exercise~\ref{0pXe-244kf.1EX-PRE.0p23rl.oqjsdlwnfvMSAS01maspe.bismdc}, we see that, for all~$N\le M-1$,
$$ S_{N,f}(x)=\sum_{{k\in\Z}\atop{-N\le k\le N}} c_{k} \,e^{2\pi i k x}=e^{-2\pi i(N+1)x}P^+_{S^{(N+1)}}(x)-
e^{2\pi iNx}P^+_{S^{(-N)}}(x).$$
On this account,
\begin{equation}\label{0pXe-244kf.SCDOVERNGMXCDFKMsfdgbdifukg02o3kerpjhxe}\begin{split}
\| S_{N,f}\|_{L^p((0,1))}&\le\|e^{-2\pi i(N+1)x}P^+_{S^{(N+1)}}\|_{L^p((0,1),\C)}+\| e^{2\pi iNx}P^+_{S^{(-N)}}\|_{L^p((0,1),\C)}\\&=
\|P^+_{S^{(N+1)}}\|_{L^p((0,1),\C)}+\| P^+_{S^{(-N)}}\|_{L^p((0,1),\C)}.
\end{split}\end{equation}

Moreover, in virtue of Exercise~\ref{0pXe-244kf.1EX-PRE.pe}, we know that, for all~$n\in\N$,
$$ P^+_{S^{(n)}}=\frac{S^{(n)}+i(S^{(n)})^\Diamond-c_{n,0}}2,$$
with
$$ c_{n,0}:=\int_0^1 S^{(n)}(x)\,dx.$$
We stress that
$$ |c_{n,0}|\le\int_0^1 |S^{(n)}(x)|\,dx=
\int_0^1 |f(x)|\,dx\le\|f\|_{L^p((0,1))}.$$
Hence, off the back of Exercise~\ref{0pXe-244kf.1EX-SK-ilcoM}, and freely renaming~$C$ line after line,
\begin{eqnarray*}
\|P^+_{S^{(n)}}\|_{L^p((0,1),\C)}&\le&\frac{\|S^{(n)}\|_{L^p((0,1),\C)}+\|(S^{(n)})^\Diamond\|_{L^p((0,1),\C)}+|c_{n,0}|}2
\\&\le& C\,\Big(
\|S^{(n)}\|_{L^p((0,1),\C)}+\|f\|_{L^p((0,1))}
\Big)\\&=& C\,\Big(
\|e^{2\pi inx} f\|_{L^p((0,1),\C)}+\|f\|_{L^p((0,1))}
\Big)\\
&\le& C\,\|f\|_{L^p((0,1))}.
\end{eqnarray*}
Combining this information and~\eqref{0pXe-244kf.SCDOVERNGMXCDFKMsfdgbdifukg02o3kerpjhxe} we gather the desired result in~\eqref{0pXe-244kf.1EX-PRE-012.d-96-0923k}.

\paragraph{Solution to Exercise~\ref{0pXe-244kf.1EX}.}
Let~$\epsilon>0$. By reason of Exercise~\ref{ojdkfnvioewyr098765rewe67890iuhgvgyuuygfr43wdfgh8bft5xs-EAN},
we can find a trigonometric polynomial~$S_\epsilon$ such that~$\|S_\epsilon-f\|_{L^p((0,1))}\le\epsilon$.

Moreover, by Exercise~\ref{PKS0-3-21}, there exists~$N_\epsilon\in\N$ such that if~$N\ge N_\epsilon$ then~$S_{N,S_\epsilon}=S_\epsilon$.

Therefore, if~$N\ge N_\epsilon$,
\begin{eqnarray*} &&\|S_{N,f}-f\|_{L^p((0,1))}\le \|S_{N,f}-S_\epsilon\|_{L^p((0,1))}+\|S_\epsilon-f\|_{L^p((0,1))}\\&&\qquad
\le \|S_{N,f}-S_{N,S_\epsilon}\|_{L^p((0,1))}+\|S_\epsilon-f\|_{L^p((0,1))}\\&&\qquad=
\|S_{N,f-S_\epsilon}\|_{L^p((0,1))}+\|S_\epsilon-f\|_{L^p((0,1))}.\end{eqnarray*}

Hence, after Exercise~\ref{0pXe-244kf.1EX-SK0-1o23lk-0923},
$$ \|S_{N,f}-f\|_{L^p((0,1))}\le C\,\|S_\epsilon-f\|_{L^p((0,1))}+\|S_\epsilon-f\|_{L^p((0,1))}\le (C+1)\,\epsilon.$$
We can now send first~$N\to+\infty$ and then~$\epsilon\searrow0$ to conclude that
$$ \lim_{N\to+\infty} \|S_{N,f}-f\|_{L^p((0,1))}=0,$$ as desired.

\section{Solutions to selected exercises of Section~\ref{FEJERKESE}}

\paragraph{Solution to Exercise~\ref{CESA}.} Assume first that~$\gamma\in\R$. Let~$\epsilon>0$ and~$K_\epsilon\in\N$ be sufficiently large such that, for every~$k\ge K_\epsilon$,
$$ \gamma_k\in [\gamma-\epsilon,\gamma+\epsilon].$$
Let also
$$ \mu_\epsilon:=\max_{k\in\{0,\dots,K_\epsilon-1\}}|\gamma_k|.$$
Then, if~$N\ge K_\epsilon$,
\begin{eqnarray*}
\sum_{k=0}^{N-1}\gamma_k&\in&\left[ \sum_{k=0}^{K_\epsilon-1}\gamma_k +(N-K_\epsilon)(\gamma-\epsilon),\,
\sum_{k=0}^{K_\epsilon-1}\gamma_k +(N-K_\epsilon)(\gamma+\epsilon)
\right]\\&
\subseteq&\left[ -K_\epsilon\mu_\epsilon+(N-K_\epsilon)(\gamma-\epsilon),\,
K_\epsilon\mu_\epsilon+(N-K_\epsilon)(\gamma+\epsilon)
\right].\end{eqnarray*}
As a consequence,
\begin{eqnarray*}&&\gamma-\epsilon=
\liminf_{N\to+\infty}\left(-\frac{K_\epsilon\mu_\epsilon}N+\frac{(N-K_\epsilon)(\gamma-\epsilon)}N\right)\le
\liminf_{N\to+\infty}\frac1N\sum_{k=0}^{N-1}\gamma_k\\&&\qquad
\le\limsup_{N\to+\infty}\frac1N\sum_{k=0}^{N-1}\gamma_k
\le\limsup_{N\to+\infty}\left(\frac{K_\epsilon\mu_\epsilon}N+\frac{(N-K_\epsilon)(\gamma+\epsilon)}N\right)=\gamma+\epsilon.
\end{eqnarray*}
Since~$\epsilon>0$ can be taken arbitrarily small, this entails that
$$ \liminf_{N\to+\infty}\frac1N\sum_{k=0}^{N-1}\gamma_k=\limsup_{N\to+\infty}\frac1N\sum_{k=0}^{N-1}\gamma_k=\gamma,$$
as desired.

Let us now consider the case in which~$\gamma=-\infty$ (the case in which~$\gamma=+\infty$ being similar).
We take~$K_\epsilon\in\N$ sufficiently large such that, for every~$k\ge K_\epsilon$,
we have that~$\gamma_k\le-\frac1\epsilon$.

In this way, if~$N\ge K_\epsilon$,
$$ \sum_{k=0}^{N-1}\gamma_k\le \sum_{k=0}^{K_\epsilon-1}\gamma_k -\frac{N-K_\epsilon}\epsilon
\le K_\epsilon\mu_\epsilon-\frac{N-K_\epsilon}\epsilon$$
and, as a result,
$$ \limsup_{N\to+\infty}\frac1N\sum_{k=0}^{N-1}\gamma_k\le
\limsup_{N\to+\infty}\left(\frac{K_\epsilon\mu_\epsilon}N-\frac{N-K_\epsilon}{N\epsilon}\right)
=-\frac1\epsilon.$$
Thus, taking~$\epsilon>0$ as small as we wish,
$$\limsup_{N\to+\infty}\frac1N\sum_{k=0}^{N-1}\gamma_k\le-\infty,$$
as desired.

As for the requested example, consider for instance the sequence
$$ \gamma_k:=\begin{dcases}
k&{\mbox{ if there exists~$j\in\N$ such that }}k=2^j,\\
\displaystyle\frac1{e^k}&{\mbox{ otherwise}}
\end{dcases}$$ and observe that, if~$N-1\in[2^M,2^{M+1})$ for some~$M\in\N$,
$$ \sum_{k=0}^{N-1}\gamma_k\le \sum_{k=0}^{2^{M+1}}\gamma_k
\le\sum_{j=0}^{M+1} j+\sum_{k=0}^{2^{M+1}} \frac1{e^k}
\le (M+2)^2+\frac{e}{e-1}$$
and for this reason
$$ \lim_{N\to+\infty}\frac1N\sum_{k=0}^{N-1}\gamma_k\le
\lim_{M\to+\infty}\frac1{2^M}\left((M+2)^2+\frac{e}{e-1}\right)=0.$$

\paragraph{Solution to Exercise~\ref{NOINF-FeJ-th.3}.} For example, let~$f$ be the square wave in Exercise~\ref{SQ:W}. Then, $f$ is bounded and periodic of period~$1$, and we show that~\eqref{FeJ-th.3.NOIJMS.021324v4c25} cannot hold true with~$p=+\infty$.

Indeed, suppose by contradiction that
$$\lim_{N\to+\infty}\left\|
\frac1N\sum_{k=0}^{N-1}S_k-f
\right\|_{L^\infty((0,1))}=0.$$
Then, by this uniform convergence, $f$ is necessarily a continuous function.
This cannot be, since~$f$ is discontinuous at~$x=\frac12$.

\paragraph{Solution to Exercise~\ref{IFEJ-FORM2.bn-EX}.}
By Exercise~\ref{K-1PIO} we know that
$$ \int_{0}^{1} D_k(x)\,dx=1,$$
which, together with~\eqref{01foqld03rkfg.GnwedRA.a}, gives the desired result.

\paragraph{Solution to Exercise~\ref{GKFND}.} On the one hand, by Exercises~\ref{fr12} and~\ref{K-3PIO},
there exists~$c>0$ such that
\begin{eqnarray*}&& \sup_{N\in\N}\int_0^1 |D_N(x)|\,dx= \sup_{N\in\N}\int_{-1/2}^{1/2} |D_N(x)|\,dx\ge c \sup_{N\in\N}\ln N
=+\infty.\end{eqnarray*}
On the other hand, by~\eqref{IFEJ-FORM2} and Exercise~\ref{IFEJ-FORM2.bn-EX},
$$ \sup_{N\in\N}\int_0^1 |F_N(x)|\,dx=\sup_{N\in\N}\int_0^1 F_N(x)\,dx=1.$$

\paragraph{Solution to Exercise~\ref{UTSJIFPOAX2e921rhnr}.}
By~\eqref{IFEJ-FORM2} and~\eqref{023wefv4567uygfdfgyhuizo0iwhfg0eb5627e01.p1rebi},
\begin{eqnarray*}&& \left|\frac1N\sum_{k=0}^{N-1}S_k(x)\right|=\left|\int_0^1 f(y)\,F_N(x-y)\,dy\right|\\&&\qquad\le
\int_0^1 |f(y)|\,F_N(x-y)\,dy\le\|f\|_{L^\infty(\R)}\int_0^1F_N(x-y)\,dy.\end{eqnarray*}
Thus, the desired result follows by Exercises~\ref{fr12} and~\ref{IFEJ-FORM2.bn-EX}.

\paragraph{Solution to Exercise~\ref{UTSJIFPOA}.}
Using~\eqref{IFEJ-FORM.bi} and~\eqref{023wefv4567uygfdfgyhuizo0iwhfg0eb5627e01.p1rebi}, we have that
\begin{eqnarray*}&&\frac1N\sum_{k=0}^{N-1}S_{k}(x)=\int_0^1 f(y)\,F_N(x-y)\,dy=
\sum_{{k\in\Z}\atop{|k|\le N-1}}\int_0^1\left(1-\frac{|k|}{N}\right)\,f(y)\,e^{2\pi ik(x-y)}\,dy\\
&&\qquad=\sum_{{k\in\Z}\atop{|k|\le N-1}}\left(1-\frac{|k|}{N}\right)\,\widehat f_k\,e^{2\pi ik x}.
\end{eqnarray*}

\paragraph{Solution to Exercise~\ref{UNIQ:ALAmowifj30othgore0-1}.} Under the assumptions of Theorem~\ref{UNIQ},
in view of Exercise~\ref{UTSJIFPOA} we have that
\begin{equation}\label{0-kwsmd.10oiekd-12} \begin{split}&\frac1N\sum_{k=0}^{N-1}S_{k,f}(x)=\sum_{{k\in\Z}\atop{|k|\le N-1}}\left(1-\frac{|k|}{N}\right)\,\widehat f_k\,e^{2\pi ik x}\\&\qquad
=\sum_{{k\in\Z}\atop{|k|\le N-1}}\left(1-\frac{|k|}{N}\right)\,\widehat g_k\,e^{2\pi ik x}=\frac1N\sum_{k=0}^{N-1}S_{k,g}(x).\end{split}\end{equation}
In light of Theorem~\ref{FeJ-th.3}, we know that the first term in~\eqref{0-kwsmd.10oiekd-12} converges to~$f$ in~$L^1((0,1))$ and that the last term in~\eqref{0-kwsmd.10oiekd-12} converges to~$g$ in~$L^1((0,1))$, therefore necessarily~$f(x)=g(x)$ as functions of~$L^1((0,1))$ and almost everywhere.

\paragraph{Solution to Exercise~\ref{UNIQ:ALAmowifj30othgore0-1ILBVIS}.}Given~$\epsilon>0$ we take~$L_\epsilon\in\N$ such that, for all~$L\ge L_\epsilon$,
$$\sum_{{k\in\Z}\atop{|k|\ge L}}|\widehat f_k|\le\epsilon.$$
Hence, by virtue of Exercise~\ref{UTSJIFPOA}, if~$N> M> L_\epsilon$,
\begin{eqnarray*}&&
\left|\frac1N\sum_{k=0}^{N-1}S_{k}(x)-\frac1M\sum_{k=0}^{M-1}S_{k}(x)\right|=\left|
\sum_{{k\in\Z}\atop{|k|\le N-1}}\left(1-\frac{|k|}{N}\right)\,\widehat f_k\,e^{2\pi ik x}
-\sum_{{k\in\Z}\atop{|k|\le M-1}}\left(1-\frac{|k|}{M}\right)\,\widehat f_k\,e^{2\pi ik x}
\right|\\&&\qquad\le\sum_{{k\in\Z}\atop{M\le |k|\le N-1}}|\widehat f_k|+\left|
\sum_{{k\in\Z}\atop{|k|\le N-1}}\frac{|k|}{N}\,\widehat f_k\,e^{2\pi ik x}
-\sum_{{k\in\Z}\atop{|k|\le M-1}}\frac{|k|}{M}\,\widehat f_k\,e^{2\pi ik x}
\right|\\&&\qquad\le\epsilon+
\left|\sum_{{k\in\Z}\atop{|k|\le M-1}}\frac{(N-M)\,|k|}{NM}\,\widehat f_k\,e^{2\pi ik x}
\right|+\left|\sum_{{k\in\Z}\atop{M\le|k|\le N-1}}\frac{|k|}{N}\,\widehat f_k\,e^{2\pi ik x}\right|\\&&\qquad\le\epsilon+\frac{N-M}{NM}
\sum_{{k\in\Z}\atop{|k|\le M-1}}|k|\,|\widehat f_k|+\sum_{{k\in\Z}\atop{M\le|k|\le N-1}}|\widehat f_k|\\&&\qquad\le\epsilon+\frac{N-M}{NM}
\sum_{{k\in\Z}\atop{|k|\le L_\epsilon-1}}|k|\,|\widehat f_k|+\frac{N-M}{NM}
\sum_{{k\in\Z}\atop{L_\epsilon\le|k|\le M-1}}|k|\,|\widehat f_k|+\epsilon\\&&\qquad\le\epsilon+
\frac{(N-M)L_\epsilon}{NM}\sum_{{k\in\Z}\atop{|k|\le L_\epsilon-1}}|\widehat f_k|+
\sum_{{k\in\Z}\atop{L_\epsilon\le|k|\le M-1}}|\widehat f_k|+\epsilon\\&&\qquad\le\epsilon+
\frac{L_\epsilon}{M}\sum_{{k\in\Z}}|\widehat f_k|+\epsilon+\epsilon\\&&\qquad\le4\epsilon,
\end{eqnarray*}
as long as~$M$ is sufficiently large (in dependence of~$\epsilon$ and~$L_\epsilon$), from which the desired result follows.

\paragraph{Solution to Exercise~\ref{FEJPOAE1edc4.10}.}
We utilise~\eqref{IFEJ-FORM} and distinguish two cases, when~$|x|\le\frac{1}{4N}$ and when~$|x|\in\left(\frac1{4N},\frac12\right]$.

If~$|x|\le\frac{1}{4N}$, then, for some~$C_1>0$,
\begin{eqnarray*}
F_N(x)\le \frac{(N\pi x)^2}{2N\sin^2(\pi x)}\le\frac{C_1\,N^2}{2N}\le\frac{C_1\,N}{2}\le\frac{C_1\,N}{1+N^2x^2}.
\end{eqnarray*}

If instead~$|x|\in\left(\frac1{4N},\frac12\right]$, then, for some~$C_2>0$,
\begin{eqnarray*}
F_N(x)\le \frac{1}{2N\sin^2(\pi x)}\le\frac{C_2}{2N x^2}=\frac{C_2\,N}{2N^2 x^2}=\frac{C_2\,N}{N^2 x^2+N^2x^2}\le
\frac{C_2\,N}{\frac1{16}+N^2x^2}.
\end{eqnarray*}From these observations, the desired result plainly follows.

\paragraph{Solution to Exercise~\ref{PKSMDu0ojf65bv94.11}.}
We stress that~$f$ is not assumed to be continuous (otherwise the result would follow from Theorem~\ref{FeJ-th.2}).

We deduce from~\eqref{PKSMDu0ojf65bv94.112} that
$$ S_k(0)=\frac{a_0}2+
\sum_{j=1}^k a_j\cos(2\pi j\cdot0)=\frac{a_0}2+
\sum_{j=1}^k a_j\ge0$$
and that, for all~$k\in\{N,\dots,2N\}$, $$S_k(0)\ge
\frac{a_0}2+
\sum_{j=1}^N a_j=S_N(0).$$

Therefore,
\begin{eqnarray*}&& \frac1{2N+1}\sum_{k=0}^{2N}S_k(0)\ge \frac1{2N+1}\sum_{k=N}^{2N}S_k(0)\ge
\frac{N+1}{2N+1} S_N(0)\ge\frac{N+1}{2(2N+1)}\sum_{j=0}^k a_j.
\end{eqnarray*}

From this (see Exercise~\ref{UTSJIFPOAX2e921rhnr}) it follows that
$$\frac{N+1}{2(2N+1)}\sum_{j=0}^k a_j\le\|f\|_{L^\infty(\R)}$$
and accordingly, sending~$N$, $k\to+\infty$,
$$ \sum_{j=0}^{+\infty} a_j\le4\|f\|_{L^\infty(\R)}.$$
From this and Exercise~\ref{UNIQ:ALAmowifj30othgore0-1ILBVIS} the desired results follow.

See also~\cite[Volume~I, Section~2]{MR171116} for related results.

\paragraph{Solution to Exercise~\ref{PKSMDu0ojf65bv94.11124erf54y}.}
Since
$$ \lim_{\theta\to0}\frac{\sin\theta}\theta=1,$$
there exists~$N_0\in\N$, $N_0\ge1$, such that, for all~$\theta\in\left[-\frac1{N_0},\frac1{N_0}\right]$,
\begin{equation}\label{PpmdREVHJKNDybtALSJDMN-1}
\sin\theta\ge\frac\theta2.
\end{equation}

Moreover, on account of~\eqref{jasmx23er} and Exercise~\ref{UTSJIFPOA},
\begin{eqnarray*}&&\frac1N\sum_{k=0}^{N-1}S_{k}(x)=
\sum_{{k\in\Z}\atop{|k|\le N-1}}\left(1-\frac{|k|}{N}\right)\,\widehat f_k\,e^{2\pi ik x}
=i\sum_{k=1}^{N-1}\left(1-\frac{k}{N}\right)\,\big(\widehat f_k-\widehat f_{-k}\big)\,\sin(2\pi kx)\\&&\qquad\qquad\qquad\qquad=\sum_{k=1}^{N-1}\left(1-\frac{k}{N}\right)\,b_k\,\sin(2\pi kx),
\end{eqnarray*}
which, combined with Exercise~\ref{UTSJIFPOAX2e921rhnr}, yields that
\begin{equation}\label{PKSMDu0ojf65bv94.112.112334556}\left|\sum_{k=1}^{N-1}\left(1-\frac{k}{N}\right)\,b_k\,\sin(2\pi kx)\right|\le\|f\|_{L^\infty(\R)}.\end{equation}
Hence, given~$L\in\N$, we choose~$x:=\frac1{4LN_0\pi}$ and~$N:=2L+1$, finding that
\begin{equation}\label{PKSMDu0ojf65bv94.112.1123345}
\|f\|_{L^\infty(\R)}\ge\sum_{k=1}^{2L}\left(1-\frac{k}{2L+1}\right)\,b_k\,\sin\left(\frac{k}{2LN_0}\right).\end{equation}

We also remark that, for all~$k\in\{1,\dots,2L\}$, we have that~$\frac{k}{2LN_0}\le\frac{1}{N_0}$ and thus,
owing to~\eqref{PpmdREVHJKNDybtALSJDMN-1},
$$ \sin\left(\frac{k}{2LN_0}\right)\ge\frac{k}{4LN_0}.$$
Hence, by~\eqref{PKSMDu0ojf65bv94.112.1123} and~\eqref{PKSMDu0ojf65bv94.112.1123345},
\begin{equation}\label{PKSMDu0ojf65bv94.112.112334555}
\|f\|_{L^\infty(\R)}\ge\sum_{k=1}^{2L}\left(1-\frac{k}{2L+1}\right)\,\frac{k\,b_k}{4LN_0}\ge
\sum_{k=1}^{L}\left(1-\frac{k}{2L+1}\right)\,\frac{k\,b_k}{4LN_0}.\end{equation}

Additionally, for all~$k\in\{1,\dots,L\}$,
$$ 1-\frac{k}{2L+1}\ge 1-\frac{L}{2L+1}\ge\frac12$$
and therefore we gather from~\eqref{PKSMDu0ojf65bv94.112.112334555} that
\begin{equation*}
\|f\|_{L^\infty(\R)}\ge\sum_{k=1}^{L}\frac{k\,b_k}{8LN_0}.\end{equation*}

As a result, using again~\eqref{PKSMDu0ojf65bv94.112.112334556},
\begin{eqnarray*}
\left|\sum_{k=1}^{L}b_k\,\sin(2\pi kx)\right|&\le&
\left|\sum_{k=1}^{L}\left(1-\frac{k}{L+1}\right)\,b_k\,\sin(2\pi kx)\right|
+\left|\sum_{k=1}^{L} \frac{k\,b_k}{L+1}\,\sin(2\pi kx)\right|\\
&\le&\|f\|_{L^\infty(\R)}+\sum_{k=1}^{L} \frac{k\,b_k}{L+1}\\&\le&\|f\|_{L^\infty(\R)}+\sum_{k=1}^{L} \frac{k\,b_k}{L}\\
&\le&\|f\|_{L^\infty(\R)}+8N_0\|f\|_{L^\infty(\R)},
\end{eqnarray*}
proving the desired bound.

The requested uniform convergence follows from Exercise~\ref{UNIQ:ALAmowifj30othgore0-1ILBVIS}.

See also~\cite[Volume~I, Section~2]{MR171116} for related results.

\paragraph{Solution to Exercise~\ref{PKSMDu0ojf65bv94.11124erf54y.bi}.}
We focus on the case of
Exercise~\ref{PKSMDu0ojf65bv94.11124erf54y} (the case in Exercise~\ref{PKSMDu0ojf65bv94.11}
is similar, changing the notion of parity and swapping sines and cosines).

Let~$\epsilon>0$. By way of Theorem~\ref{FeJ-th.2}, we find~$N_\epsilon\in\N$ such that, for all~$N\ge N_\epsilon$,
\begin{equation}\label{PKSMDu0ojf65bv94.112.1123.7}\sup_{x\in\R}\left|
\frac1N\sum_{k=0}^{N-1}S_k(x)-f(x)
\right|\le\epsilon.\end{equation}

Let$$ g^{(N)}(x):=f(x)-\frac1N\sum_{k=0}^{N-1}S_k(x)$$
and notice that~$g^{(N)}$ is an odd function (because so are~$f$ and~$S_k$), and that, for every~$k$, $\ell\in\N$,
\begin{eqnarray*}&&2
\int_0^1 S_k(x)\,\sin(2\pi \ell x)\,dx=2
\sum_{k=1}^N b_k\int_0^1\sin(2\pi kx)\,\sin(2\pi \ell x)\,dx={b_k\,\delta_{\ell,k}}
\end{eqnarray*}
and, as a result, the Fourier Series of~$g^{(N)}$ in trigonometric form is
$$ \sum_{\ell=1}^{+\infty} b^{(N)}_\ell\,\sin(2\pi\ell x),$$
where
\begin{eqnarray*}
b^{(N)}_\ell:=b_\ell-\frac1N\sum_{k=0}^{N-1}
{b_k\,\delta_{\ell,k}}=\left(1-\frac1N\right)b_\ell,
\end{eqnarray*}
which is non-negative, due to~\eqref{PKSMDu0ojf65bv94.112.1123}.

Hence, we deduce from Exercise~\ref{PKSMDu0ojf65bv94.11124erf54y} and~\eqref{PKSMDu0ojf65bv94.112.1123.7}
that, for all~$N\ge N_\epsilon$ and all~$M\in\N$,
\begin{equation}\label{PJSMDJPOHDm1u0rjo3ihy3.2r21t5} \sup_{x\in\R}|S_{M,g^{(N)}}(x)|\le C\|g^{(N)}\|_{L^\infty(\R)}\le C\epsilon.\end{equation}

Moreover (see Exercise~\ref{PKS0-3-21}) when~$M\ge N-1\ge k$ we have that
$$ S_{M,S_k} =S_k$$
and thus
$$ S_{M,g^{(N)}}=S_{M,f}-\frac1N\sum_{k=0}^{N-1}S_{M,S_k}=S_{M,f}-\frac1N\sum_{k=0}^{N-1}S_k.$$

From this and~\eqref{PKSMDu0ojf65bv94.112.1123.7} we find that, if~$K$ is large enough,
\begin{equation*}\sup_{x\in\R}\left|S_{M,f}(x)-S_{M,g^{(N)}}(x)-f(x)
\right|\le\epsilon.\end{equation*}
This and~\eqref{PJSMDJPOHDm1u0rjo3ihy3.2r21t5} yield that
\begin{equation*}\sup_{x\in\R}\left|S_{M,f}(x)-f(x)
\right|\le(C+1)\epsilon,\end{equation*}
thus establishing the desired uniform convergence.

\paragraph{Solution to Exercise~\ref{PKSMDu0ojf65bv94.11124erf54y.bi.PY}.}
Let
$$ \phi(x):=\frac{f(x)+f(-x)}{2}\qquad{\mbox{and}}\qquad\psi(x):=\frac{f(x)-f(-x)}{2}.$$
Notice that~$\phi$ and~$\psi$ are continuous and periodic of period~$1$. Also, $\phi$ is even and~$\psi$ is odd.

Moreover, if the Fourier Series of~$f$ in trigonometric form is
$$ \frac{a_0}2+\sum_{k=1}^{+\infty}\Big(a_k\cos(2\pi kx) + b_k\sin(2\pi kx)\Big),$$
with~$a_k$, $b_k\ge0$, we have that, for all~$k\in\N$,
\begin{eqnarray*}
2\int_0^1 \phi(x)\,\cos(2\pi kx)\,dx=\int_0^1 \big(f(x)+f(-x)\big)\,\cos(2\pi kx)\,dx=2\int_0^1 f(x)\,\cos(2\pi kx)\,dx=a_k
\end{eqnarray*}
and
\begin{eqnarray*}&&
2\int_0^1 \phi(x)\,\sin(2\pi kx)\,dx=\int_0^1 \big(f(x)+f(-x)\big)\,\sin(2\pi kx)\,dx\\&&\qquad=\int_0^1 f(x)\,\sin(2\pi kx)\,dx-\int_0^1 f(x)\,\sin(2\pi kx)\,dx=0,
\end{eqnarray*}
giving that the Fourier Series of~$\phi$ in trigonometric form is
$$ \frac{a_0}2+\sum_{k=1}^{+\infty}a_k\cos(2\pi kx) .$$
Similarly, the Fourier Series of~$\psi$ in trigonometric form is
$$\sum_{k=1}^{+\infty}b_k\sin(2\pi kx) .$$

Therefore, we can apply Exercise~\ref{PKSMDu0ojf65bv94.11124erf54y.bi} to the functions~$\phi$ and~$\psi$.

In particular, we know that
$$ \lim_{N\to+\infty} \left(\sup_{x\in\R}\left| S_{N,\phi}(x)-\phi(x)\right|
+\sup_{x\in\R}\left| S_{N,\psi}(x)-\psi(x)\right|\right)=0.$$
Since~$f=\phi+\psi$, we thus conclude that
\begin{eqnarray*}&&
\lim_{N\to+\infty}\sup_{x\in\R}\left| S_{N,f}(x)-f(x)\right|=
\lim_{N\to+\infty}\sup_{x\in\R}\left| S_{N,\phi}(x)+S_{N,\psi}(x)-\phi(x)-\psi(x)\right|
\\&&\qquad\le\lim_{N\to+\infty} \left(\sup_{x\in\R}\left| S_{N,\phi}(x)-\phi(x)\right|+\sup_{x\in\R}\left| S_{N,\psi}(x)-\psi(x)\right|\right)=0
,\end{eqnarray*}
as desired.

See also~\cite[Volume~I, Section~2]{MR171116} and~\cite{MR1563037} for related, and also more general, results.

\paragraph{Solution to Exercise~\ref{FEJPOAE1edc4.1}.}
We know (see e.g.~\cite[Theorem~5.22]{MR3381284}) that~$f$ is finite almost everywhere
and (see e.g.~\cite[Theorem~7.15]{MR3381284}) that almost every point is a Lebesgue density point for~$f$:
more explicitly, for almost every~$x\in\R$ we can assume that~$|f(x)|<+\infty$ and that
$$ \lim_{h\searrow0}\frac1{2h}\int_{-h}^{h}\big|f(x)-f(x+t)\big|\,dt=0.$$
In particular, given~$\epsilon>0$, we pick~$\delta_\epsilon\in\left(0,\frac14\right)$ such that, for all~$h\in(0,\delta_\epsilon]$,
$$ \frac1{2h}\int_{-h}^{h}\big|f(x)-f(x+t)\big|\,dt\le\epsilon.$$

We also define
$$ I(h):=\int_{0}^{h}\big|f(x)-f(x+t)\big|\,dt$$
and we remark that, for all~$h\in[-\delta_\epsilon,\delta_\epsilon]$,
$$ |I(h)|\le2\epsilon h.$$

Hence, after Exercise~\ref{FEJPOAE1edc4.10} and an integration by parts,
\begin{eqnarray*}
&&\int_{\{|t|\le\delta_\epsilon\}} \big|f(x)-f(x+t)\big|\,\left( \frac{\sin(N\pi t)}{\sin(\pi t)}\right)^2\,dt\le
CN^2\int_{\{|t|\le\delta_\epsilon\}} \big|f(x)-f(x+t)\big|\,\frac{dt}{1+N^2t^2}\\&&\qquad
=CN^2\int_{\{|t|\le\delta_\epsilon\}} I'(t)\,\frac{dt}{1+N^2t^2}\\
&&\qquad=\frac{CN^2\big(I(\delta_\epsilon)-I(-\delta_\epsilon)\big)}{1+N^2\delta_\epsilon^2}
+2CN^4\int_{\{|t|\le\delta_\epsilon\}} t\,I(t)\,\frac{dt}{(1+N^2t^2)^2}
\\&&\qquad\le
\frac{4CN^2\epsilon\delta_\epsilon}{1+N^2\delta_\epsilon^2}
+4CN^4\epsilon\int_{\{|t|\le\delta_\epsilon\}} t^2\,\frac{dt}{(1+N^2t^2)^2}\\&&\qquad\le
\frac{4CN^2\epsilon\delta_\epsilon}{1+N^2\delta_\epsilon^2}
+4CN\epsilon\int_{-\infty}^{+\infty} s^2\,\frac{ds}{(1+s^2)^2}\\&&\qquad\le \widetilde C \left(
\frac{N^2\epsilon\delta_\epsilon}{1+N^2\delta_\epsilon^2}+N\epsilon
\right),
\end{eqnarray*}
for some constants~$C$, $\widetilde C>0$.

From this and~\eqref{FAGBS:0-0.2} we arrive at 
\begin{equation*}
\begin{split}&
\left|\int_{(0,1)\cap{\mathcal{D}}_\delta} \big(f(x)-f(y)\big)\,\left( \frac{\sin(N\pi (x-y))}{\sin(\pi (x-y))}\right)^2\,dy
\right|\le  \widetilde C \left(
\frac{N^2\epsilon\delta_\epsilon}{1+N^2\delta_\epsilon^2}+N\epsilon
\right),
\end{split}
\end{equation*}
where
$$ {\mathcal{D}}_\delta:=\bigcup_{\ell\in\Z}(\ell+x-\delta,\,\ell+x+\delta).$$

As a result, recalling~\eqref{IFEJ-FORM}, \eqref{023wefv4567uygfdfgyhuizo0iwhfg0eb5627e01.p1rebi}, and~\eqref{FAGBS:0-0.1},
\begin{eqnarray*}&&\lim_{N\to+\infty}\left|f(x)-\frac1N\sum_{k=0}^{N-1}S_k(x)\right|
\le\lim_{N\to+\infty}
\left(\frac{|f(x)|+\|f\|_{L^1((0,1))}}{2N\sin^2(\pi\delta_\epsilon)}+ \widetilde C \left(
\frac{N\epsilon\delta_\epsilon}{1+N^2\delta_\epsilon^2}+\epsilon
\right)\right)
=\widetilde C\epsilon.
\end{eqnarray*}
Taking now~$\epsilon$ as small as we wish, we obtain the desired result.

\paragraph{Solution to Exercise~\ref{FEJPO-010}.} 
At a first glance, we could proceed as in Exercise~\ref{NONSepCFGDV}, that is, by~\eqref{jasmx23er}, for every~$\ell\in\N$,
\begin{eqnarray*}&&
\int_0^1f(x)\,\cos(2\nu_\ell\pi x)\,dx=\frac{\lambda_\ell}2
\end{eqnarray*}
and therefore
$$ \frac{|\lambda_\ell|}2\le\int_0^1|f(x)|\,|\cos(2\nu_\ell\pi x)|\,dx=\int_0^1f(x)\,|\cos(2\nu_\ell\pi x)|\,dx\le
\int_0^1f(x)\,dx=\frac{\lambda_0}2,$$
giving the desired result -- but only up to a factor~$2$, which is not enough to solve this exercise!

Incidentally, the assumption that none of the~$\nu_k$ divides another has not been utilised: therefore one has to use another, more subtle, strategy. Namely, we integrate the function~$f$ against the Fej\'er Kernel and we use~\eqref{IFEJ-FORM2} and the sign assumption on~$f$ to see that, for all~$\ell\in\N$ and~$x\in\R$,
$$ 0\le\int_0^1 f (y) \,F_N(x-\nu_\ell y)\,dy.$$
This and~\eqref{IFEJ-FORM.bi} give that
\begin{equation}\label{0.m2kmfBSMD-1}\begin{split}
0\le&\sum_{{k\in\Z}\atop{|k|\le N-1}}\left(1-\frac{|k|}{N}\right)\,\int_0^1f(y)\,e^{2\pi ik(x-\nu_\ell y)}\,dy
\\=&\sum_{{k\in\Z}\atop{|k|\le N-1}}\left(1-\frac{|k|}{N}\right)\,\widehat f_{k\nu_\ell}\,e^{2\pi ikx}.
\end{split}\end{equation}
Also, by~\eqref{jasmx23er} and~\eqref{A0kmalseb.a}, we have that~$\lambda_0=2\Re(\widehat f_0)=2\widehat f_0$,
that~$\Im(\widehat f_k)=0$, that~$\widehat f_{-k}=\widehat f_k$,
that~$\widehat f_k=0$ unless~$|k|=\nu_j$ for some~$j$, and that~$\lambda_j=\widehat f_{\nu_j}+\widehat f_{-\nu_j}$.

In particular, since none of the~$\nu_k$ divides another, if~$\widehat f_{k\nu_\ell}\ne0$ then necessarily~$k\in\{0,1,-1\}$.
Accordingly, \eqref{0.m2kmfBSMD-1} boils down to
\begin{equation*}\begin{split}
0\le\,&\widehat f_0+\left(1-\frac{1}{N}\right)\,\big(\widehat f_{\nu_\ell}\,e^{2\pi ix}+\widehat f_{-\nu_\ell}\,e^{-2\pi ix}\big)\\
=\,&\widehat f_0+2\left(1-\frac{1}{N}\right)\,\widehat f_{\nu_\ell}\,\cos(2\pi x)\\=\,&\frac{\lambda_0}2+\left(1-\frac{1}{N}\right)\,\lambda_{\nu_\ell}\,\cos(2\pi x).
\end{split}\end{equation*}
Sending~$N\to+\infty$, we obtain that
$$ 0\le\frac{\lambda_0}2+\lambda_{\nu_\ell}\,\cos(2\pi x)$$
and the desired result thus follows by choosing~$x:=0$ and~$x:=\frac12$.

\paragraph{Solution to Exercise~\ref{kpoqsld0O3.1P2P37jsd1}.}
Yes, they are.

For the optimality in Exercise~\ref{NONSepCFGDV}
suppose, for a contradiction, that~\eqref{12wewiatf} can be sharpened in the form
\begin{equation}\label{9uhgvgh.mjhyt5nmeucpfz0o}
|\widehat f_k|\le\mu\widehat f_0,\end{equation}
for some~$\mu\in(0,1)$.

Then, given~$m\in\R$, consider the function
\begin{equation}\label{bNijmdcvm6tgr6tclb-01} f(x):=\big(2\cos(4\pi x)\big)^{2m}\ge0.\end{equation}
We remark that
\begin{eqnarray*}&&
f(x)=\left( e^{4\pi ix}+e^{-4\pi ix}\right)^{2m}
=\sum_{j=0}^{2m} \left({2m}\atop{j}\right) e^{4\pi i jx}\,e^{-4\pi i(2m-j)x}
=\sum_{k=-m}^{m} \left({2m}\atop{m+k}\right) e^{8\pi i kx}.
\end{eqnarray*}
Hence, we infer from~\eqref{9uhgvgh.mjhyt5nmeucpfz0o} that
$$ \mu \left({2m}\atop{m}\right)=
\mu\widehat f_0\ge|\widehat f_4|=\left({2m}\atop{m+1}\right) $$
and therefore
$$ 1>\mu\ge\lim_{m\to+\infty}\displaystyle\frac{\displaystyle\left({2m}\atop{m+1}\right)}{\displaystyle\left({2m}\atop{m}\right)}=\lim_{m\to+\infty}\frac{m}{m+1}=1,$$
which is a contradiction, showing that the estimate in Exercise~\ref{NONSepCFGDV} is sharp.

\begin{figure}[h]
\includegraphics[height=2.8cm]{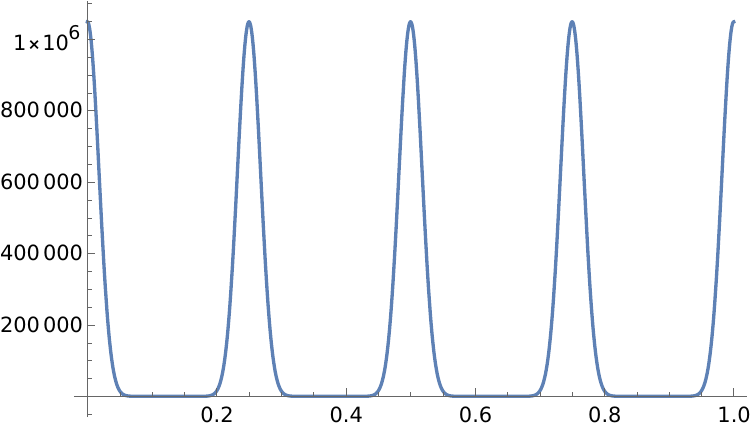}$\,\;\quad$\includegraphics[height=2.8cm]{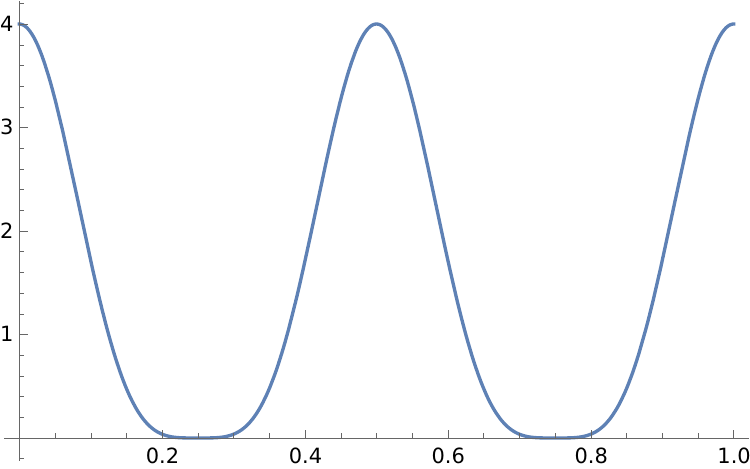}$\,\;\quad$
\includegraphics[height=2.8cm]{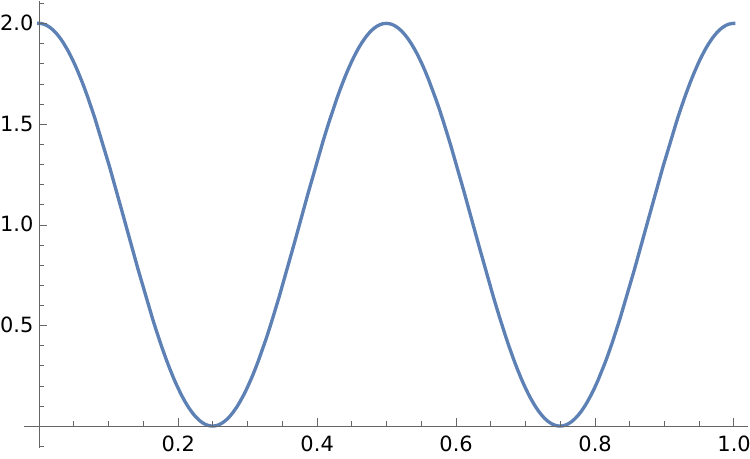}
\centering
\caption{Left: plot of the function in~\eqref{bNijmdcvm6tgr6tclb-01} with~$m:=10$.
Center: plot of the function in~\eqref{bNijmdcvm6tgr6tclb-02}.
Right: plot of the function in~\eqref{bNijmdcvm6tgr6tclb-03}.}\label{9ikm9400ntrtgvcl8b9mb8n.ikdnsEU}
\end{figure}

Now we highlight the necessity of assuming that the frequencies in Exercise~\ref{FEJPO-010} do not divide one another.
For this, we can consider the example
\begin{equation}\label{bNijmdcvm6tgr6tclb-02} f(x):=\big(1+\cos(4\pi x)\big)^2\ge0.\end{equation}
We observe that
$$ f(x)=1+2\cos(4\pi x)+\cos^2(4\pi x)=\frac{3}2+2\cos(4\pi x)+\frac12\cos (8\pi x),$$
which is of the form~\eqref{A0kmalseb.a} with~$\lambda_0:=3$, $\lambda_1:=2$, $\lambda_2:=\frac12$, $\lambda_k:=0$ for all~$k\ge3$, $\nu_1:=2$ and~$\nu_2:=4$. Notice that in this case~$\nu_1$ divides~$\nu_2$ and~\eqref{A0kmalseb.a0-ln} is not satisfied, since
\begin{equation*} |\lambda_1 |=2 >\frac32= \frac{\lambda_0}2.\end{equation*}
This says that for~\eqref{A0kmalseb.a0-ln} to be true, the assumption that the frequencies do not divide one another
cannot be dropped.

For the optimality in Exercise~\ref{FEJPO-010} when the frequencies do not divide one another, 
consider the case of the function
\begin{equation}\label{bNijmdcvm6tgr6tclb-03}f(x)=1+\cos(4\pi x),\end{equation}
which is non-negative and can be written as in~\eqref{A0kmalseb.a}, with~$\lambda_0:=2$, $\lambda_1:=1$, $\nu_1:=2$
(and, for instance, for all~$k\ge2$, $\lambda_k:=0$ and~$\nu_k$ being the $(k-1)$th prime number larger than~$2$).
In this situation, we have that~$|\lambda_1|=\frac{\lambda_0}2$, which attains the optimality in~\eqref{A0kmalseb.a0-ln}.

See Figure~\ref{9ikm9400ntrtgvcl8b9mb8n.ikdnsEU} for the diagrams of the functions utilised in this exercise.

\paragraph{Solution to Exercise~\ref{MSCOCOS}.}
Yes, the series in~\eqref{CSPPCSjKS} converges to a continuous function periodic of period~$1$,
but no, the Fourier Series of this function does not converge uniformly (and in particular it does not converge uniformly to the original function). This is an interesting example, because it says that Theorem~\ref{C1uni} would not hold if one replaced the uniform Dini's Condition~\eqref{DINI-UNO} with the assumption that~$f$ is continuous. 

Also, this example reveals how useful Theorem~\ref{FeJ-th.2} can be in the presence of continuous functions, since while Fourier Series can fail to uniformly converge to such functions, the averaged Fourier Sum always converges.

Moreover, this example is even stronger than what is required here, since, as we will see below, not only
the Fourier Series fails to converge, but also it does not converge at the origin.
In this way, this example shows also that
there exist continuous functions of period~$1$ whose Fourier Series does not converge at a point, a fact that will be retaken in Theorem~\ref{COMADIVESKMDWD} and further expanded in Section~\ref{EXCE}.

To check that the given series is the Fourier Series of a continuous function periodic of period~$1$, we recall
Exercise~\ref{SPDCD-0.02} (used here with~$N:=2\ell$) to see that
\begin{equation*} Q(x,\ell)=2\sin(4\pi \ell x)\sum_{j=1}^{\ell}\frac{\sin(2\pi jx)}j.\end{equation*}
That being so, recalling Exercise~\ref{PIPIO} (used here with~$\gamma_k:=\frac1k$) we obtain that
$$ |Q(x,\ell)|\le 2\big|\sin(4\pi \ell x)\big|\, (\pi+1)\le 2(\pi+1).$$

On this account, given~$M\in\N$, we define
\begin{equation}\label{91wueihdncos2pimx}
f_M(x):=\sum_{k=1}^{M} \alpha_k \,Q(x,n_k)\end{equation}
and we observe that, for all~$M> M'$,
$$ \sup_{x\in\R}|f_M(x)-f_{M'}(x)|\le\sum_{k=M'+1}^{M} \alpha_k  \sup_{x\in\R}|Q(x,n_k)|\le 2(\pi+1)\sum_{k=M'+1}^{M} \alpha_k ,$$
which, in view of~\eqref{ZLAM-01}, is the tail of a convergent series.

Therefore
\begin{equation}\label{OJSN.0o9ijn.plkjhgvbnHFDRTYUJSNBVYujhndcmowkv-91}
{\mbox{$f_M$ converges uniformly to a certain function~$f$.}}\end{equation} Since~$f_M$ is a continuous function, periodic of period~$1$, we infer that so is~$f$. This answers the first question of Exercise~\ref{MSCOCOS}.

To answer the second question of Exercise~\ref{MSCOCOS}, we use the change of index~$m:=2n_k-j$ to deduce that
\begin{equation}\label{OJSN.0o9ijn.plkjhgvbnHFDRTYUJSNBVYujhndcmowkv-1} \begin{split}
f_M(x)&=\sum_{k=1}^{M} \alpha_k\,\left(\sum_{{j\in\Z\setminus\{0\}}\atop{|j|\le n_k}}\frac{\cos(2\pi(2n_k-j)x)}{j}\right)\\
&=\sum_{k=1}^{M} \alpha_k\,\left(\sum_{{m\in\N\cap[n_k,3n_k]}\atop{m\neq 2n_k}}\frac{\cos(2\pi mx)}{2n_k-m}\right).
\end{split}
\end{equation}

Now we use the monotonicity in~\eqref{ZLAM-01a} to find a strictly increasing function, with inverse function~$\psi:\R\to\R$
strictly increasing, that describes the map~$\N\ni k\mapsto n_k$, that is, given~$p\in\R$, the condition~$p\le n_k$ is equivalent to~$\psi(p)\le k$. 

With this notation, we can swap the order of the finite summations in~\eqref{OJSN.0o9ijn.plkjhgvbnHFDRTYUJSNBVYujhndcmowkv-1} and conclude that
\begin{equation} \label{OJSN.0o9ijn.plkjhgvbnHFDRTYUJSNBVYujhndcmowkv-919}\begin{split}
f_M(x)&=\sum_{m=n_1}^{n_M}
\sum_{{k\in\N}\atop{{\psi(m/3)\le k\le\min\{\psi(m),M\}}\atop{k\ne\psi(m/2)}}}
\frac{\alpha_k\,\cos(2\pi mx)}{2n_k-m}\\&=\sum_{m=n_1}^{n_M}a_m\,\cos(2\pi mx)\\&=\sum_{m=1}^{n_M}a_m\,\cos(2\pi mx),
\end{split}
\end{equation}
where
$$ a_m:=\begin{dcases}\displaystyle\sum_{{k\in\N}\atop{{\psi(m/3)\le k\le\min\{\psi(m),M\}}\atop{k\ne\psi(m/2)}}}
\frac{\alpha_k}{2n_k-m} &{\mbox{ if }}m\ge n_1,\\
0&{\mbox{ if }}m< n_1.\end{dcases}
$$

Now we check that (in spite of their rather complicated expression),
\begin{equation}\label{SONLOE}
{\mbox{the coefficients~$a_m$ are the Fourier coefficients of~$f$ in trigonometric form.}}
\end{equation}
To this end, we observe that, as a consequence of~\eqref{ZLAM-01a},
\begin{equation}\label{SONLOE.22}
{\mbox{$n_M\to+\infty$, as~$M\to+\infty$.}}\end{equation}

We also use~\eqref{OJSN.0o9ijn.plkjhgvbnHFDRTYUJSNBVYujhndcmowkv-91} and~\eqref{OJSN.0o9ijn.plkjhgvbnHFDRTYUJSNBVYujhndcmowkv-919} to see that, for all~$p\in\N$,
\begin{eqnarray*}&&
2\int_0^1 f(x)\,\sin(2\pi px)\,dx=2\lim_{M\to+\infty}\int_0^1 f_M(x)\,\sin(2\pi px)\,dx\\
&&\qquad=2\lim_{M\to+\infty}\sum_{m=1}^{n_M} a_m\int_0^1\cos(2\pi mx)\,\sin(2\pi px)\,dx=0
\end{eqnarray*}
and that
\begin{eqnarray*}&&
2\int_0^1 f(x)\,\cos(2\pi px)\,dx=2\lim_{M\to+\infty}\int_0^1 f_M(x)\,\cos(2\pi px)\,dx\\
&&\qquad=2\lim_{M\to+\infty}\sum_{m=1}^{n_M} a_m\int_0^1\cos(2\pi mx)\,\cos(2\pi px)\,dx=
\lim_{M\to+\infty} \sum_{m=1}^{n_M} a_m\,\delta_{mp}=a_p,
\end{eqnarray*}
where~\eqref{SONLOE.22} has been used in the last step.

These observations and the setting of Fourier coefficients in trigonometric form~\eqref{jasmx23er}
establish~\eqref{SONLOE}, as desired.

We now prove that
\begin{equation}\label{SONLOE.2134r2fg2}
{\mbox{the Fourier Series of~$f$ does not converge at the origin.}}\end{equation}
This is rather surprising, especially when we compare this statement to~\eqref{OJSN.0o9ijn.plkjhgvbnHFDRTYUJSNBVYujhndcmowkv-91}.

To prove~\eqref{SONLOE.2134r2fg2}, we suppose, by contradiction, that, given~$\epsilon>0$, there exists~$N_\epsilon\in\N$ such that
\begin{equation}\label{SONLOE.2134r2fg201p2we}
{\mbox{for all~$N_1$, $N_2\ge N_\epsilon$}}\end{equation}
we have that
\begin{equation}\label{SONLOE.2134r2fg201p2we12rt35v5dw9} \epsilon\ge |S_{N_1}(0)-S_{N_2}(0)|.\end{equation}
We choose~$N_1:=n_M$ and~$N_2:=n_{M-1}$, with~$M$ large enough such that~\eqref{SONLOE.2134r2fg201p2we} is satisfied (and this is possible, thanks to~\eqref{SONLOE.22}).

Moreover, it follows from~\eqref{OJSN.0o9ijn.plkjhgvbnHFDRTYUJSNBVYujhndcmowkv-919} and~\eqref{SONLOE},
that, for all~$M\in\N$,
$$ S_{n_M,f}(x)=\sum_{m=1}^{n_M}a_m\,\cos(2\pi mx)=f_M(x).$$
Combining this with~\eqref{91wueihdncos2pimx.01oewjdflnv} and~\eqref{91wueihdncos2pimx} we arrive at
$$ S_{n_M,f}(0)=f_M(0)=\sum_{k=1}^{M} \alpha_k \,Q(0,n_k)
=\sum_{k=1}^{M}\sum_{{j\in\Z\setminus\{0\}}\atop{|j|\le n_k}}\frac{\alpha_k}{j}.$$

From this and~\eqref{SONLOE.2134r2fg201p2we12rt35v5dw9} we arrive at
\begin{eqnarray*}
\epsilon&\ge& \left|\sum_{k=1}^{M}\sum_{{j\in\Z\setminus\{0\}}\atop{|j|\le n_k}}\frac{\alpha_k}{j}-\sum_{k=1}^{M-1}\sum_{{j\in\Z\setminus\{0\}}\atop{|j|\le n_k}}\frac{\alpha_k}{j}\right|\\&=&
\sum_{{j\in\Z\setminus\{0\}}\atop{|j|\le n_M}}\frac{\alpha_M}{j}\\&\ge&
c\alpha_M\,\ln n_M,
\end{eqnarray*}
for some~$c>0$.

Hence, if we call~$\sigma>0$ the left-hand side of~\eqref{ZLAM-02}, we can pick a suitable sub-sequence in~$M$ such that
$$\epsilon\ge c\lim_{M\to+\infty}\alpha_M\,\ln n_M=c\sigma,$$
which is absurd if~$\epsilon$ is small enough.

This establishes~\eqref{SONLOE.2134r2fg2}, thus completing (and even exceeding) the analysis required in Exercise~\ref{MSCOCOS}.

See~\cite[Theorems~(1$\cdot$3) and~(1$\cdot$13) on pages~299--300]{MR1963498}
and~\cite[Theorem~12.35]{MR3381284}
for even more general and deeper results.

\paragraph{Solution to Exercise~\ref{POISSONKERN-ex1}.} By Exercise~\ref{FBA},
\begin{equation}\label{oqpkdwgmintthcrz0frazz0dzarz}
\left| \sum_{{k\in\Z}} r^{|k|}\,\widehat f_k\,e^{2\pi ikx}\right|\le
\sum_{{k\in\Z}} r^{|k|}\,|\widehat f_k|
\le\|f\|_{L^1((0,1))}\sum_{{k\in\Z}} r^{|k|}<+\infty,
\end{equation}
yielding the desired result.

Notice that for~$r=1$ the Abel mean would reduce to the Fourier Series of~$f$ (however, the convergence
of the Fourier Series of an integrable function may fail, see Section~\ref{EXCE}), therefore Abel means are useful tools to ``regularise'' possibly divergent series.
See~\cite[Section~12.7]{MR3381284} for further information about Abel means.

\paragraph{Solution to Exercise~\ref{POISSONKERN-ex2}.}
We have that $$\left|\sum_{k=1}^{+\infty}r^k\cos(2\pi k\theta)\right|\le\sum_{k=1}^{+\infty}r^k<+\infty,$$
since~$r\in[0,1)$, showing the convergence of the series in~\eqref{POISSONKERN-ex2-form1}.

What is more, if~$z\in\C$ with~$|z|<1$,
$$ 1+2\sum_{k=1}^{+\infty}z^k=-1+2\sum_{k=0}^{+\infty}z^k=-1+\frac2{1-z}=\frac{1+z}{1-z}$$
and thus, choosing~$z:=re^{2\pi i\theta}$ with~$r\in(0,1)$ and~$\theta\in\R$, we see that
\begin{equation}\label{MSLM-12345} 1+2\sum_{k=1}^{+\infty}r^ke^{2\pi ik\theta}
=\frac{1+re^{2\pi i\theta}}{1-re^{2\pi i\theta}}.\end{equation}
Taking the real part of this identity, we arrive at
\begin{eqnarray*}
{\mathcal{P}}(r,\theta)=1+2\sum_{k=1}^{+\infty}r^k\cos(2\pi k\theta)=\Re\left( 1+2\sum_{k=1}^{+\infty}r^ke^{2\pi ik\theta}\right)=\Re\left(\frac{1+re^{2\pi i\theta}}{1-re^{2\pi i\theta}}\right).
\end{eqnarray*}
Multiplying both the numerator and the denominator by~$1-re^{-2\pi i\theta}$, we obtain
\begin{eqnarray*}
{\mathcal{P}}(r,\theta)=\Re\left(\frac{1-r^2+re^{2\pi i\theta}-re^{-2\pi i\theta}}{1+r^2-re^{2\pi i\theta}-re^{-2\pi i\theta}}\right)=\frac{1-r^2}{1+r^2-2r\cos(2\pi \theta)},
\end{eqnarray*}
yielding the desired result in~\eqref{POISSONKERN-ex2-form2}.

Furthermore, by Exercise~\ref{POISSONKERN-ex1}, and recalling~\eqref{fasv},
the $r$-Abel mean of~$f$ can be written as
\begin{eqnarray*} &&\sum_{{k\in\Z}} r^{|k|}\,\widehat f_k\,e^{2\pi ikx}=
\widehat f_0+\sum_{k=1}^{+\infty} r^k\,\big(\widehat f_k\,e^{2\pi ikx}
+\widehat f_{-k}\,e^{-2\pi ikx}\big)=
\widehat f_0+\sum_{k=1}^{+\infty} r^k\,\big(\widehat f_k\,e^{2\pi ikx}
+\overline{ \widehat f_k\,e^{2\pi ikx}}\big)\\&&\qquad=
\widehat f_0+2\Re\left(\sum_{k=1}^{+\infty} r^k\,\widehat f_k\,e^{2\pi ikx}\right),
\end{eqnarray*}
which, en passant, proves~\eqref{XIUJHGV9iuhg-2.1pjolsnENASMDP}, and is equal to
\begin{eqnarray*}&&
\Re\left(\int_0^1f(y)\,dy+2
\sum_{k=1}^{+\infty} r^k\,
\int_0^1f(y)\,e^{2\pi ik(x-y)}\,dy\right)\\
&&\qquad=\int_0^1f(y)\,\Re\left(1+2
\sum_{k=1}^{+\infty} r^k\,e^{2\pi ik(x-y)}\right)\,dy.
\end{eqnarray*}

\begin{figure}[h]
\includegraphics[height=3cm]{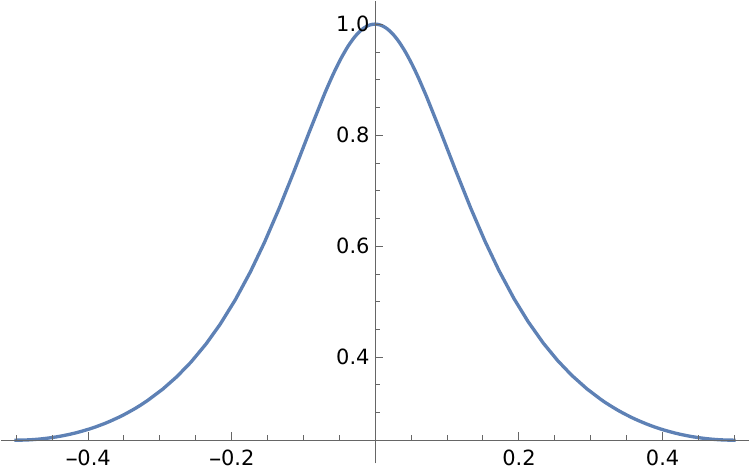}$\quad$\includegraphics[height=3cm]{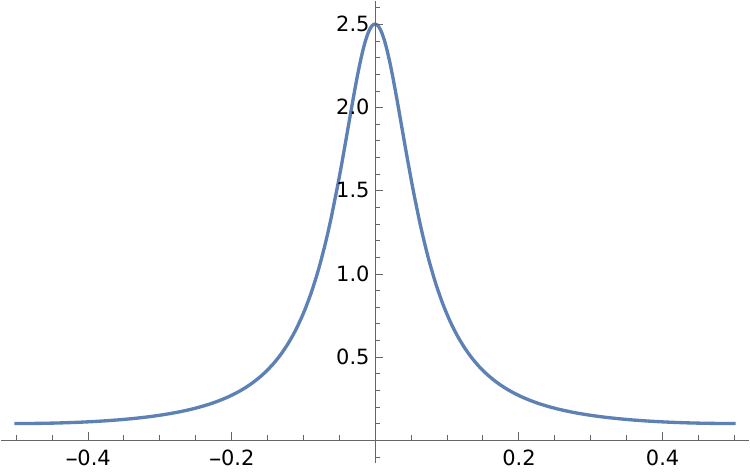}$\quad$\includegraphics[height=3cm]{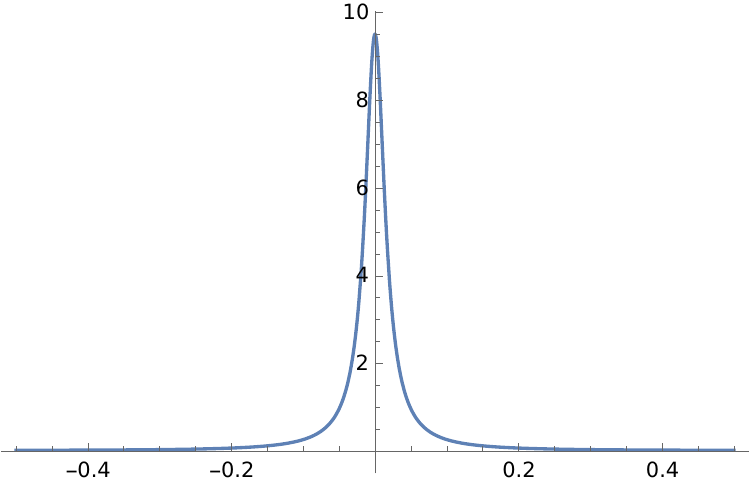}
\centering
\caption{Plot of the periodic Poisson Kernel with~$r\in\left\{\frac13,\frac23,\frac{9}{10}\right\}$.}\label{CBDVMVAMRNSLMD-o1kemrf}
\end{figure}

From this and~\eqref{POISSONKERN-ex2-form1} we obtain~\eqref{POISSONKERN-ex2-form3}, as desired.
See Figure~\ref{CBDVMVAMRNSLMD-o1kemrf} for a visual sketch of the periodic Poisson Kernel. 
Compared with Figure~\ref{amdi.fe}, the reader may appreciate some similarity in the concentration phenomenon with the Fej\'er Kernel, but notice that here ``less oscillations'' are present.

\paragraph{Solution to Exercise~\ref{POISSONKERN-ex3}.}
We use the notation~$\epsilon:=1-r\searrow0$.

Integrating~\eqref{POISSONKERN-ex2-form1}, we see that, for all~$r\in(0,1)$,
$$ \int_0^1{\mathcal{P}}(r,y)\,dy=1.$$
Thus, in light of Exercise~\ref{POISSONKERN-ex2}, the difference between the $r$-Abel mean of~$f$ and~$f$ itself at the point~$x$ can be written as
\begin{equation}\label{POISSONLAMsxCOjsdR1}
\begin{split}&{\mathcal{V}}(x):=
\int_0^1 f(x-y)\,{\mathcal{P}}(r,y)\,dy-f(x)=
\int_0^1 \big(f(x-y)-f(x)\big)\,{\mathcal{P}}(r,y)\,dy\\&\qquad\qquad
=\int_0^1 \frac{(1-r^2)\,\big(f(x-y)-f(x)\big)}{(1-r)^2+4r\sin^2(\pi y)}\,dy=\int_0^1\frac{(2\epsilon-\epsilon^2)\big(f(x-y)-f(x)\big)}{ \epsilon^2+4(1-\epsilon)\sin^2(\pi y)}\,dy.
\end{split}\end{equation}

Now, given~$\eta>0$, to be taken as small as we wish in what follows, we use the continuity of~$f$ at the point~$x$ to find some~$\delta>0$ (depending on~$\eta$ and~$x$, but independent of~$\epsilon$)
such that, whenever~$y\in(-\delta,\delta)$, we have that~$|f(x-y)-f(x)|\le\eta$. In this way, substituting for~$t:=y-1$ and then for~$Y:=\frac{y}\epsilon$, we find that
\begin{eqnarray*}
&& \int_0^\delta\frac{(2\epsilon-\epsilon^2)\big|f(x-y)-f(x)\big|}{ \epsilon^2+4(1-\epsilon)\sin^2(\pi y)}\,dy+
\int_{1-\delta}^1\frac{(2\epsilon-\epsilon^2)\big|f(x-y)-f(x)\big|}{ \epsilon^2+4(1-\epsilon)\sin^2(\pi y)}\,dy\\
&&\qquad\le\int_0^\delta\frac{3\epsilon\big|f(x-y)-f(x)\big|}{ \epsilon^2+\sin^2(\pi y)}\,dy+
\int_{-\delta}^0\frac{3\epsilon\big|f(x-t)-f(x)\big|}{ \epsilon^2+\sin^2(\pi (t+1))}\,dt\\
&&\qquad\le\int_{-\delta}^\delta\frac{3\epsilon\big|f(x-y)-f(x)\big|}{ \epsilon^2+\sin^2(\pi y)}\,dy\\
&&\qquad\le\int_{-\delta}^\delta\frac{C\epsilon\eta}{ \epsilon^2+y^2}\,dy\\
&&\qquad=\int_{-\delta/\epsilon}^{\delta/\epsilon}\frac{C\eta}{ 1+Y^2}\,dY\\&&\qquad\le C\eta,
\end{eqnarray*}
for some~$C>0$ possibly varying from line to line.

On this account and~\eqref{POISSONLAMsxCOjsdR1}, we gather that
\begin{equation*}
|{\mathcal{V}}(x)|\le\int_\delta^{1-\delta}\frac{6\epsilon\|f\|_{L^\infty(\R)}}{ \epsilon^2+\sin^2(\pi \delta)}\,dy+C\eta\le
\frac{6\epsilon\|f\|_{L^\infty(\R)}}{ \epsilon^2+\sin^2(\pi \delta)}\,dy+C\eta.
\end{equation*}
Therefore, sending~$r\nearrow1$ (i.e., $\epsilon\searrow0$),
\begin{equation*}
\lim_{r\nearrow1}|{\mathcal{V}}(x)|\le C\eta.
\end{equation*}
We now take~$\eta$ as small as we wish and we obtain the desired result.

\paragraph{Solution to Exercise~\ref{LPmPOIKE01}.} In virtue of~\eqref{POISSONLAMsxCOjsdR1}
and Minkowski's Integral Inequality (see e.g.~\cite{GENT})
the norm in~$L^p((0,1))$ of the difference between the $r$-Abel mean of~$f$ and~$f$ can be written as
\begin{equation}\label{fetyuVVSnbgfUH3381284} \begin{split}&{\mathcal{W}}:=\left(\int_0^{1}\big|{\mathcal{V}}(x)\big|^p \,dx\right)^{\frac1p}\le\left(
\int_0^{1}\left( 
\int_0^1\frac{3\epsilon\big|f(x-y)-f(x)\big|}{ \epsilon^2+\sin^2(\pi y)}\,dy
\right)^p \,dx\right)^{\frac1p}\\&\qquad\qquad\qquad\le
\int_0^1\left( \int_0^1 \left(\frac{3\epsilon\big|f(x-y)-f(x)\big|}{ \epsilon^2+\sin^2(\pi y)}\right)^p
\,dx\right)^{\frac1p}\,dy\\&\qquad\qquad\qquad=\int_0^1\frac{3\epsilon}{ \epsilon^2+\sin^2(\pi y)}\left( \int_0^1 \big|f(x-y)-f(x)\big|^p\,dx\right)^{\frac1p}\,dy.
\end{split}\end{equation}

Now we pick~$\eta>0$, to be taken arbitrarily small in what follows and we use the continuity of the translations in Lebesgue spaces (see e.g.~\cite[Theorem~8.19]{MR3381284}) to find~$\delta>0$ (depending on~$\eta$, but independent of~$\epsilon$) such that when~$y\in(-\delta,\delta)$
$$ \int_0^1 \big|f(x-y)-f(x)\big|^p\,dx\le\eta^p.$$
On this account,
\begin{eqnarray*}&&
\int_{(0,\delta)\cup(1-\delta,1)}\frac{3\epsilon}{ \epsilon^2+\sin^2(\pi y)}\left( \int_0^1 \big|f(x-y)-f(x)\big|^p
\,dx\right)^{\frac1p}\,dy\\
&&\qquad=\int_{-\delta}^\delta\frac{3\epsilon}{ \epsilon^2+\sin^2(\pi y)}\left( \int_0^1 \big|f(x-y)-f(x)\big|^p
\,dx\right)^{\frac1p}\,dy\\&&\qquad\le\int_{-\delta}^\delta\frac{3\epsilon\eta}{ \epsilon^2+\sin^2(\pi y)}\,dy\\
&&\qquad\le C\eta\int_{-\infty}^{+\infty}\frac{dY}{1+Y^2}\\
&&\qquad\le C\eta,
\end{eqnarray*}
for some~$C>0$ that we freely rename line after line.

Plugging this information into~\eqref{fetyuVVSnbgfUH3381284}, we conclude that
\begin{eqnarray*}{\mathcal{W}}&\le&\int_\delta^{1-\delta}\frac{3\epsilon}{ \epsilon^2+\sin^2(\pi y)}\left( \int_0^1 \big|f(x-y)-f(x)\big|^p\,dx\right)^{\frac1p}\,dy+C\eta
\\&\le&\frac{6\epsilon\,\|f\|_{L^p((0,1))}}{\sin^2(\pi\delta)}+C\eta.
\end{eqnarray*}
Hence, taking the limit as~$r\nearrow1$ (that is~$\epsilon\searrow0$) we gather that
$$ \lim_{r\nearrow1}{\mathcal{W}}\le C\eta.$$
From this, the desired result follows by sending~$\eta\searrow0$.

\paragraph{Solution to Exercise~\ref{LPmPOIKE01.002l42er21fygnb9mnvf}.}
The strategy is similar to that of Exercise~\ref{POISSONKERN-ex2}.
We have that $$\left|\sum_{k=1}^{+\infty}r^k\sin(2\pi k\theta)\right|\le\sum_{k=1}^{+\infty}r^k<+\infty,$$
since~$r\in[0,1)$, showing the convergence of the series in~\eqref{POISSONKERN-ex2-form1-ILCO}.

Also, by~\eqref{MSLM-12345},
\begin{equation*}
1+2\sum_{k=1}^{+\infty}r^k \cos(2\pi ik\theta)+2i\sum_{k=1}^{+\infty}r^k \sin(2\pi ik\theta)=
1+2\sum_{k=1}^{+\infty}r^ke^{2\pi ik\theta}
=\frac{1+re^{2\pi i\theta}}{1-re^{2\pi i\theta}},\end{equation*}
from which we obtain~\eqref{LAMCIKCONJA0-2}.

Furthermore, multiplying both the numerator and the denominator by~$1-re^{-2\pi i\theta}$, we find that
\begin{eqnarray*}
{\mathcal{Q}}(r,\theta)=\Im\left(\frac{1+re^{2\pi i\theta}}{1-re^{2\pi i\theta}}\right)=
\Im\left(\frac{1-r^2+re^{2\pi i\theta}-re^{-2\pi i\theta}}{1+r^2-re^{2\pi i\theta}-re^{-2\pi i\theta}}\right)=\frac{2r\sin(2\pi\theta)}{ (1-r)^2+4r\sin^2(\pi\theta)},
\end{eqnarray*}
yielding the desired result in~\eqref{POISSONKERN-ex2-form2medl213mrpo21t}.

\paragraph{Solution to Exercise~\ref{POISSONKERN-exsadsdghj3Xs23954868.0-1}.}
In view of the uniform convergence in~\eqref{POISSONKERN-ex2-form1-ILCO},
we calculate that
\begin{eqnarray*}&&i
\int_0^1 f(y)\,{\mathcal{Q}}(r,x-y)\,dy=2i
\sum_{k=1}^{+\infty}r^k\int_0^1 f(y)\,\sin(2\pi k(x-y))\,dy\\&&\qquad=
\sum_{k=1}^{+\infty}r^k\int_0^1 f(y)\,e^{2\pi ik(x-y)}\,dy
-\sum_{k=1}^{+\infty}r^k\int_0^1 f(y)\,e^{-2\pi ik(x-y)}\,dy\\&&\qquad=
\sum_{k=1}^{+\infty}r^k\widehat f_k\,e^{2\pi ikx}-
\sum_{k=1}^{+\infty}r^k\widehat f_{-k}\,e^{-2\pi ikx}\\&&\qquad=\sum_{k\in\Z}r^{|k|}\,{\operatorname{sign}}(k)\,\widehat f_{k}\,e^{2\pi ikx},
\end{eqnarray*}
as desired. 

\paragraph{Solution to Exercise~\ref{LPmPOIKE01ijdqsp21rj90jtm}.}
By~\eqref{lmq-21iDgbnmpeTArhmewxi}, if~$r\in(0,1)$ and~$R\in(r,1)$, then
\begin{eqnarray*}g_{R,r}(x):=
f^\Diamond_R(x)-f^\Diamond_r(x)=-i\sum_{k\in\Z}\big(R^{|k|}-r^{|k|}\big)\,{\operatorname{sign}}(k)\,\widehat f_{k}\,e^{2\pi ikx}=
\sum_{k\in\Z}c_{R,r,k}\,e^{2\pi ikx}
\end{eqnarray*}
with
$$ c_{R,r,k}:=-i\big(R^{|k|}-r^{|k|}\big)\,{\operatorname{sign}}(k)\,\widehat f_{k}.$$
We stress that~$\overline{c_{R,r,-k}}=c_{R,r,k}$ and
$$ \sum_{k\in\Z}|c_{R,r,k}|^2 \le\sum_{k\in\Z} \big(R^{|k|}-r^{|k|}\big)^2\,|\widehat f_{k}|^2 \le\sum_{k\in\Z} R^{2|k|}\,|\widehat f_{k}|^2\le\sum_{k\in\Z} |\widehat f_{k}|^2<+\infty,$$
thanks to~\eqref{L2THM.0-02}, and accordingly, in light of Theorem~\ref{THCOL2FB}(ii),
we have that~$c_{R,r,k}$ is the $k$th Fourier coefficient of~$g_{R,r}\in L^2((0,1))$.

Therefore, by~\eqref{L2THM.0-02},
\begin{equation}\label{CegjlesumsO}
\|f^\Diamond_R-f^\Diamond_r\|^2_{L^2((0,1))}= \|g_{R,r}\|^2_{L^2((0,1))}=\sum_{k\in\Z} |c_{k,R,r}|^2  \le\sum_{k\in\Z} \big(R^{|k|}-r^{|k|}\big)^2\,|\widehat f_{k}|^2.
\end{equation}
Hence, given~$\epsilon>0$, we pick~$N_\epsilon\in\N$ such that
$$ \sum_{{k\in\Z} \atop{|k|\ge N_\epsilon+1}}|\widehat f_{k}|^2\le\epsilon$$
and we deduce from~\eqref{CegjlesumsO} that
\begin{eqnarray*}
\|f^\Diamond_R-f^\Diamond_r\|^2_{L^2((0,1))}&\le&\sum_{{k\in\Z} \atop{|k|\le N_\epsilon}}\big(R^{|k|}-r^{|k|}\big)^2\,|\widehat f_{k}|^2+
\epsilon\\&\le&C_\epsilon (R-r)^2\sum_{{k\in\Z} \atop{|k|\le N_\epsilon}}|\widehat f_{k}|^2+
\epsilon
\end{eqnarray*}
and thus we can find~$R_\epsilon\in(0,1)$ such that if~$R\ge r\ge R_\epsilon$ then~$\|f^\Diamond_R-f^\Diamond_r\|^2_{L^2((0,1))}\le2\epsilon$. This gives that~$f^\Diamond_r$ is a Cauchy sequence in~$L^2((0,1))$, from which the desired result follows.

\section{Solutions to selected exercises of Section~\ref{SEC:GIBBS-PH}}

\paragraph{Solution to Exercise~\ref{9k:INTEHigii.0}.}
First off, we see that if~$x\in(-\infty,0]$ then~$Nx-1\le-1$, whence we see that~$ \varphi_N(x)=0$ and thus~\eqref{NOAKsm0.1mef1} holds true.

Let now~$x\in(0,+\infty)$. When~$N\ge\frac2x$, we have that~$Nx-1\ge1$ and therefore~\eqref{NOAKsm0.1mef1} follows in this case as well.

But the limit is not uniform, because
$$ \lim_{N\to+\infty}\sup_{x\in\R}\varphi_N(x)\ge \lim_{N\to+\infty}\varphi_N\left(\frac1N\right)=\varphi(0)=1.$$

\paragraph{Solution to Exercise~\ref{ALLSFAMNLS}.} We have that, for every~$m\in\N\setminus\{0\}$,
\begin{equation}\label{ALLSFAMNLS3}
\begin{split}\Xi_m&:=
\left|\int_0^{\pi m} \frac{\sin \tau }{\pi(N+1)\sin\left(\frac{\tau}{2(N+1)}\right)}\,d\tau
-\frac2\pi\int_0^{\pi m} \frac{\sin \tau }{\tau}\,d\tau\right|
\\&\qquad=
\frac1\pi\left|\int_0^{\pi m} \sin \tau \left(\frac{ 1}{(N+1)\sin\left(\frac{\tau}{2(N+1)}\right)}
-\frac{2}{\tau}\right)\,d\tau\right|\\&\qquad=
\frac2\pi\left|\int_0^{\pi m} \frac{\sin \tau}\tau\cdot\frac{ \frac{\tau}{2(N+1)}-\sin\left(\frac{\tau}{2(N+1)}\right)}{\sin\left(\frac{\tau}{2(N+1)}\right)}
\,d\tau\right|\\&\qquad=
\frac2\pi\left|\int_0^{\pi m} \frac{\sin \tau}\tau\;
\phi\left(\frac{\tau}{2(N+1)}\right)\,d\tau\right|,
\end{split}\end{equation}
where, for all~$x\in(0,2\delta\pi)\subset(0,\pi)$, we set
$$ \phi(x):=\frac{ x-\sin x}{\sin x}.$$

Notice that~$\phi$ is continuous at~$0$, with $$\lim_{x\searrow0}\phi(x)=0$$and, since~$\delta\in\left(0,\frac12\right)$,
\begin{equation}\label{USMdfvk02oekjfv023erf.1qms.1} \sup_{x\in(0,2\delta\pi)}|\phi'(x)|\le C,\end{equation}
for some~$C>0$.

Hence, up to renaming~$C$ line after line, we have that, for all~$\tau\in(0,2\delta\pi)$,
\begin{equation}\label{USMdfvk02oekjfv023erf.1qms.2} \left|\phi\left(\frac{\tau}{2(N+1)}\right)\right|\le\frac{C\tau}N\end{equation}
and therefore
$$\left|\int_0^{\pi} \frac{\sin \tau}\tau\;\phi\left(\frac{\tau}{2(N+1)}\right)\,d\tau\right|\le \frac{C}N.$$
This and~\eqref{ALLSFAMNLS3} give that
$$ \Xi_1\le\frac{C}N$$
and, for all~$m\in\N$ with~$m\ge2$,
\begin{eqnarray*}\Xi_m
&\le&\frac2\pi\left|\int_{\pi}^{\pi m} \frac{\sin \tau}\tau\;
\phi\left(\frac{\tau}{2(N+1)}\right)\,d\tau\right|+\frac{C}N.
\end{eqnarray*}
This, with an integration by parts in the sine function, can be estimated by
\begin{eqnarray*}\Xi_m
&\le&\frac2\pi\Bigg(\frac1{\pi m}\left|\phi\left(\frac{\pi m}{2(N+1)}\right)\right|+
\frac1{\pi}\left|\phi\left(\frac{\pi}{2(N+1)}\right)\right|+
\left|\int_{\pi}^{\pi m} \frac{\cos \tau}{\tau^2}\;
\phi\left(\frac{\tau}{2(N+1)}\right)\,d\tau\right|\\&&\qquad\qquad\qquad
+\frac{1}{2(N+1)}\left|\int_{\pi}^{\pi m} \frac{\cos \tau}\tau\;
\phi'\left(\frac{\tau}{2(N+1)}\right)\,d\tau\right|\Bigg)+\frac{C}N.
\\&\le& C
\left|\int_{\pi}^{\pi m} \frac{\cos \tau}{\tau^2}\;
\phi\left(\frac{\tau}{2(N+1)}\right)\,d\tau\right|+\frac{C}{N}\left|\int_{\pi}^{\pi m} \frac{\cos \tau}\tau\;
\phi'\left(\frac{\tau}{2(N+1)}\right)\,d\tau\right|+\frac{C}N.
\end{eqnarray*}

Hence, using again~\eqref{USMdfvk02oekjfv023erf.1qms.1} and~\eqref{USMdfvk02oekjfv023erf.1qms.2},
\begin{eqnarray*}
\Xi_m\le \frac{C}N
\int_{\pi}^{\pi m} \frac{1}{\tau}\,d\tau+\frac{C}N\le \frac{C\ln N}N,
\end{eqnarray*}
which entails the desired result.

\paragraph{Solution to Exercise~\ref{9k:INTEHigii}.} Suppose that~$c^{(1)}$ and~$c^{(2)}$ belong to the Gibbs Set of~$f$ at~$p$, with~$c^{(1)}<c^{(2)}$.
Let~$c\in(c^{(1)},c^{(2)})$. We need to show that~$c$ belongs to the Gibbs Set of~$f$ at~$p$ too.

For this, for~$j\in\{1,2\}$, let~$p_N^{(j)}$ be sequences such that~$p_N^{(j)}\to p$ as~$N\to+\infty$ and
$$ \lim_{N\to+\infty} S_{N}(p_N^{(j)})=c^{(j)}.$$
Let~$N_\star\in\N$ be large enough such that, for all~$N\ge N_\star$,
$$ S_N(p_N^{(1)})<c<S_N(p_N^{(2)}).$$
We consider the continuous function
$$ \Phi_N(t):=S_N\big((1-t)p_N^{(1)}+tp_N^{(2)}\big).$$
By construction,
$$ \Phi_N(0)<c<\Phi_N(1)$$
and therefore, by the Intermediate Value Theorem, there exists~$t_N\in(0,1)$ such that~$\Phi_N(t_N)=c$.

We define~$p_N:=(1-t_N)p_N^{(1)}+t_N\,p_N^{(2)}$ and we stress that~$p_N$ belongs to the segment joining~$p_N^{(1)}$ and~$p_N^{(2)}$.

As a result,
$$ \lim_{N\to+\infty}p_N=p.$$
Moreover,
$$ S_N(p_N)=\Phi_N(t_N)=c,$$
whence~$c$ belongs to the Gibbs Set of~$f$ at~$p$.

\paragraph{Solution to Exercise~\ref{KMSX:0oikjnhG-1sq0wuohgi}.}
First off, we perform a horizontal translation to set~$x_0:=0$. Then, we apply Theorem~\ref{IGIDE}
to the function
$$g(x):=-f(-x)+\ell_++\ell_-,$$ noticing that
$$ \lim_{x\searrow0}g(x)=\ell_+
\qquad{\mbox{and}}\qquad \lim_{x\nearrow0}g(x)=\ell_-.$$
Hence,
we infer from~\eqref{GIB-a1} that, for all~$\delta\in(0,\mu)$,
\begin{equation}\label{GIB-a1.qe} \begin{split}&
\lim_{N\to+\infty}\inf_{x\in\left(-\delta,0\right)\cup\left(0,\delta\right)}\big(S_{N,f}(x)-f(x)\big)
=-\lim_{N\to+\infty}\sup_{x\in\left(-\delta,0\right)\cup\left(0,\delta\right)}
\big(S_{N,g}(x)-g(x)\big)\\&\qquad
=-\frac{\ell_+-\ell_-}\pi\int_0^{\pi } \frac{\sin \tau }{\tau}\,d\tau+\frac{\ell_+-\ell_-}{2}=-\lambda(\ell_+-\ell_-)\end{split}\end{equation}
and from~\eqref{GIB-a1.las} that
\begin{equation}\label{GIB-a1.las.qe} \begin{split}&\lim_{N\to+\infty}
S_{N,f}\left(-\frac1{2N}\right)=-\lim_{N\to+\infty}S_{N,g}\left(\frac1{2N}\right)+\ell_++\ell_-
\\&\qquad
=-\frac{\ell_{+}-\ell_{-}}\pi\int_0^{\pi } \frac{\sin \tau }{\tau}\,d\tau
+\frac{\ell_++\ell_-}2
=-\lambda(\ell_+-\ell_-)+\ell_-.\end{split}\end{equation}

Now we denote by~${\mathcal{G}}$ the Gibbs Set of~$f$ at the origin and we claim that
\begin{equation}\label{UNIMNIAMDRFHJDEWQCNHUI-01}
{\mathcal{G}}\subseteq \big[\ell_- - \lambda(\ell_+-\ell_-),\,\ell_+ +\lambda(\ell_+-\ell_-)\big].
\end{equation}
To this end, pick~$c\in{\mathcal{G}}$. Then, there exists an infinitesimal sequence~$p_N$ with
$$ \lim_{N\to+\infty} S_{N,f}(p_N)=c.$$
Thus, since, for large~$N$,
$$ S_{N,f}(p_N)-f(p_N)\le
\sup_{x\in\left(-\delta,0\right)\cup\left(0,\delta\right)}\big(S_{N,f}(x)-f(x)\big),$$
we deduce from~\eqref{GIB-a1} that
\begin{equation}\label{8j.aEv7sjncsHV}\begin{split}&
c-\ell_+\le
c-\limsup_{N\to+\infty} f(p_N)=
\liminf_{N\to+\infty} \big(S_{N,f}(p_N)-f(p_N)\big)\\&\qquad\le\lim_{N\to+\infty}
\sup_{x\in\left(-\delta,\right)\cup\left(0,\delta\right)}\big(S_{N,f}(x)-f(x)\big)
=\lambda(\ell_+-\ell_-).
\end{split}\end{equation}
Similarly, since, for large~$N$,
$$ S_{N,f}(p_N)-f(p_N)\ge
\inf_{x\in\left(-\delta,0\right)\cup\left(0,\delta\right)}\big(S_{N,f}(x)-f(x)\big),$$
we deduce from~\eqref{GIB-a1.qe} that
\begin{equation*}\begin{split}&
c-\ell_-\ge
c-\liminf_{N\to+\infty} f(p_N)=
\limsup_{N\to+\infty} \big(S_{N,f}(p_N)-f(p_N)\big)\\&\qquad\ge\lim_{N\to+\infty}
\inf_{x\in\left(-\delta,0\right)\cup\left(0,\delta\right)}\big(S_{N,f}(x)-f(x)\big)
=-\lambda(\ell_+-\ell_-).
\end{split}\end{equation*}
This and~\eqref{8j.aEv7sjncsHV} yield that
$$ c\in\big[\ell_- - \lambda(\ell_+-\ell_-),\,\ell_+ +\lambda(\ell_+-\ell_-)\big],$$
which proves~\eqref{UNIMNIAMDRFHJDEWQCNHUI-01}.

Now we claim that
\begin{equation}\label{UNIMNIAMDRFHJDEWQCNHUI-02}
{\mathcal{G}}\supseteq \big[\ell_- - \lambda(\ell_+-\ell_-),\,\ell_+ +\lambda(\ell_+-\ell_-)\big].
\end{equation}
To establish this, in light of Exercise~\ref{9k:INTEHigii}, it suffices to check that
both~$\ell_- - \lambda(\ell_+-\ell_-)$ and~$\ell_+ +\lambda(\ell_+-\ell_-)$ belong to~${\mathcal{G}}$.
But this holds true, thanks to~\eqref{GIB-a1.las} and~\eqref{GIB-a1.las.qe}, thus completing the proof of~\eqref{UNIMNIAMDRFHJDEWQCNHUI-02}.

The desired result now follows by combining~\eqref{UNIMNIAMDRFHJDEWQCNHUI-01} and~\eqref{UNIMNIAMDRFHJDEWQCNHUI-02}.

\paragraph{Solution to Exercise~\ref{0-o2epkdjfm.019k4-9k:INTEHigii}.} 
We observe that~$|\phi_N(x)|\le1$, therefore
\begin{equation}\label{smdkdKSjd02owdkHAScvkfdgdahJ-3}\begin{split}& \lim_{N\to+\infty}\sup_{x\in(0,\delta)}\big(\phi_N(x)-w(x)\big)=
\lim_{N\to+\infty}\sup_{x\in(0,\delta)}\big(\phi_N(x)-1\big)\le0.\end{split}\end{equation}
Moreover, if~$x\in\R\setminus(\Z/2)$, we have that~$\sin(2\pi x)\ne0$
and consequently
$$ \lim_{N\to+\infty}\big( \sin^2(2\pi x)\big)^{\frac{1}{N}}=\lim_{N\to+\infty}
\big|\sin(2\pi x)\big|^{\frac{2}{N}}=1.$$
That being so, we find that, if~$x\in\R\setminus(\Z/2)$,
\begin{eqnarray*}&&\lim_{N\to+\infty}\big(\phi_N(x)-w(x)\big)=
\lim_{N\to+\infty}\big[
(-1)^{\lfloor2 x\rfloor} \big( \sin^2(2\pi x)\big)^{\frac{1}{N}}-(-1)^{\lfloor2 x\rfloor}\big]\\&&\qquad
=(-1)^{\lfloor2 x\rfloor}\lim_{N\to+\infty}\big[
 \big( \sin^2(2\pi x)\big)^{\frac{1}{N}}-1\big]=0.
\end{eqnarray*}
This proves~\eqref{smdkdKSjd02owdkHAScvkfdgdahJ-1}, as desired, and it also entails that
$$ 0=\lim_{N\to+\infty}\phi_N\left(\frac\delta2\right)-w\left(\frac\delta2\right)\le\lim_{N\to+\infty}\sup_{x\in(0,\delta)}\big(\phi_N(x)-w(x)\big).$$
Combining this inequality and~\eqref{smdkdKSjd02owdkHAScvkfdgdahJ-3}, we obtain~\eqref{smdkdKSjd02owdkHAScvkfdgdahJ-2}.

\begin{figure}[h]
\includegraphics[height=2.8cm]{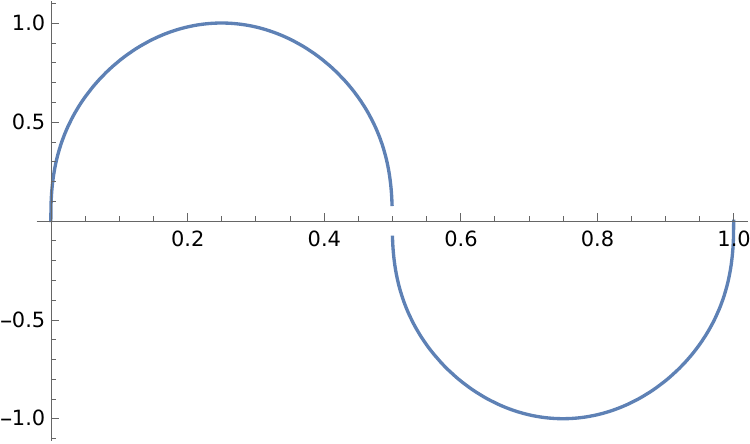}
$\quad$\includegraphics[height=2.8cm]{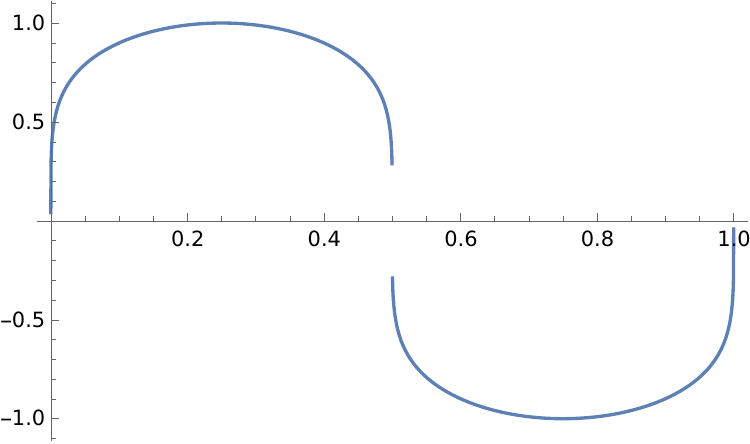}$\quad$\includegraphics[height=2.8cm]{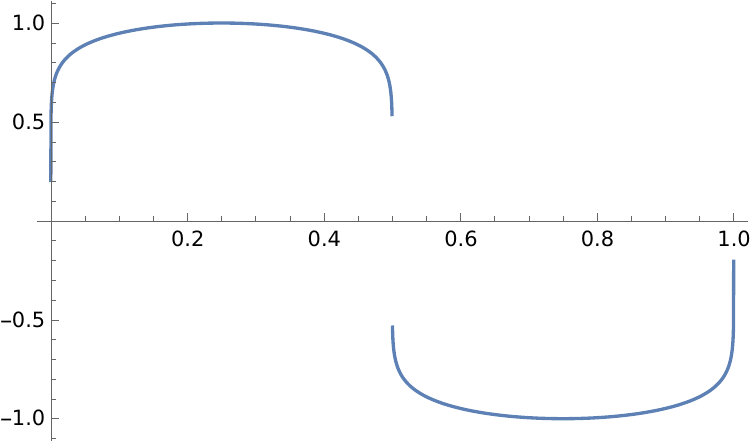}
\centering
\caption{Plot of the function~$(-1)^{\lfloor2 x\rfloor} \big( \sin^2(2\pi x)\big)^{\frac{1}{N}}$ for~$N\in\{5,10,20\}$.}\label{0FSfgpfcnlmsjmofkwlqp-24.fdfghjjnn2meI202}
\end{figure}

See Figure~\ref{0FSfgpfcnlmsjmofkwlqp-24.fdfghjjnn2meI202} to visualise the approximation of the square wave function dealt with in this exercise. Notice the absence of the Gibbs phenomenon in comparison with the Fourier Series approximation sketched in Figure~\ref{09iuygfy876trfdfghjjnn2meI20}.

\section{Solutions to selected exercises of Section~\ref{EXCE}}

\paragraph{Solution to Exercise~\ref{SENFYUJNSBKSxD}.}
Let~$\epsilon>0$. On account of~\eqref{SENFYUJNSBKSxD0}, for all~$k\in\N$ there exists~$L_{\epsilon,k}$ such that, for all~$L\ge L_{\epsilon,k}$, one has that
$$ \left| \sum_{\ell=L+1}^{+\infty}\mu_{k,\ell}\right|\le\epsilon.$$

With greater reason, setting
$$ L_\epsilon:=\max_{k\in\{0,\dots,N\}}L_{\epsilon,k},$$
it follows that, for all~$L\ge L_{\epsilon}$ and all~$k\in\{0,\dots,N\}$, one has that
$$ \left| \sum_{\ell=L+1}^{+\infty}\mu_{k,\ell}\right|\le\epsilon.$$

Hence, if~$M>L\ge L_{\epsilon}$, for all~$k\in\{0,\dots,N\}$, one has that
$$ \left| \sum_{\ell=L+1}^{M}\mu_{k,\ell}\right|=
\left| \sum_{\ell=L+1}^{+\infty}\mu_{k,\ell}- \sum_{\ell=M+1}^{+\infty}\mu_{k,\ell}\right|\le\left| \sum_{\ell=L+1}^{+\infty}\mu_{k,\ell}\right|+\left| \sum_{\ell=M+1}^{+\infty}\mu_{k,\ell}\right|\le2\epsilon.$$

These observations yield that, if~$M>L_{\epsilon}$,
\begin{eqnarray*}&&
\left|\sum_{k=0}^N\sum_{\ell=0}^{+\infty}\mu_{k,\ell}-\sum_{\ell=0}^{M}\sum_{k=0}^N\mu_{k,\ell}\right|\le
\left|\sum_{k=0}^N\sum_{\ell=0}^{L_\epsilon}\mu_{k,\ell}-\sum_{\ell=0}^{M}\sum_{k=0}^N\mu_{k,\ell}\right|+
\sum_{k=0}^N\left|\sum_{\ell=L_\epsilon+1}^{+\infty}\mu_{k,\ell}\right|\\&&\qquad\le
\left|\sum_{k=0}^N\sum_{\ell=0}^{L_\epsilon}\mu_{k,\ell}-\sum_{k=0}^N\sum_{\ell=0}^{M}\mu_{k,\ell}\right|+
(N+1)\epsilon\\&&\qquad=
\left|\sum_{k=0}^N\;\sum_{\ell=L_\epsilon+1}^{M}\mu_{k,\ell}\right|+
(N+1)\epsilon\\&&\qquad\le
\sum_{k=0}^N\left|\sum_{\ell=L_\epsilon+1}^{M}\mu_{k,\ell}\right|+
(N+1)\epsilon\\&&\qquad\le2(N+1)\epsilon+(N+1)\epsilon
\end{eqnarray*}
and then, sending~$M\to+\infty$,
$$ \left|\sum_{k=0}^N\sum_{\ell=0}^{+\infty}\mu_{k,\ell}-\sum_{\ell=0}^{+\infty}\sum_{k=0}^N\mu_{k,\ell}\right|\le
3(N+1)\epsilon.$$
Since~$\epsilon$ is arbitrary, the desired result is established.

An alternative solution goes as follows:
this is a short version of the previous strategy, relying more on the limit notation and on the fact that 
\begin{equation}\label{PROVSPI}\begin{split}&
{\mbox{the limit of a finite sum is the finite sum of the limit,}}
\\&{\mbox{\em{provided that these limits exist and are finite,}}}\end{split}\end{equation}
see e.g.~\cite[Theorem~2, Chapter~5]{zbMATH05043443} -- so, the proof is shorter but one has not to be sloppy with calculus. 

The full argument goes like this. We set
$$ \sigma_{L,k}:=\sum_{\ell=0}^{L}\mu_{k,\ell}.$$
By~\eqref{SENFYUJNSBKSxD0}, we know that, for all~$k\in\{0,\dots,N\}$ the limits
$$ \lambda_k:=\lim_{L\to+\infty}\sum_{\ell=0}^{L}\mu_{k,\ell}=\lim_{L\to+\infty}\sigma_{L,k}$$
exist and are finite.

Therefore,
\begin{equation}\label{PROVSPI2}
\sum_{k=0}^N\sum_{\ell=0}^{+\infty}\mu_{k,\ell}= \sum_{k=0}^N \lambda_k= \sum_{k=0}^N\lim_{L\to+\infty}\sigma_{L,k}=
\lim_{L\to+\infty}\sum_{k=0}^N\sigma_{L,k},\end{equation}
where the last step is due to~\eqref{PROVSPI}.

Also, since finite sums can be freely swapped,
$$ \sum_{k=0}^N\sigma_{L,k}=\sum_{k=0}^N\sum_{\ell=0}^{L}\mu_{k,\ell}=\sum_{\ell=0}^{L}\sum_{k=0}^N\mu_{k,\ell}.$$
One thus inserts this information into~\eqref{PROVSPI2} and obtains the desired result.

\paragraph{Solution to Exercise~\ref{SRTAMAYBHJAOPLKAZXAQNBKNEMINSNC}.}
By way of Theorem~\ref{SMXC22} (used here with~$m:=0$), we have that, for all~$k\in\Z$,
\begin{equation*} |\widehat f_k|\le \|f\|_{L^\infty(\R)} \end{equation*}
and therefore
$$ |S_{N,f}(x_0)|\le\sum_{{k\in\Z}\atop{|k|\le N}} |\widehat f_k|\le2N+1,$$
establishing~\eqref{TACLCSIMOMNDADEKPD.1}.

Now, to check~\eqref{TACLCSIMOMNDADEKPD.2}, the idea would be to choose a continuous function that coincides with the sign of the Dirichlet Kernel.
Strictly speaking, this is not possible, since the Dirichlet Kernel does have zeros (hence the sign function is not continuous there).
To get around with this difficulty, one can ``remove'' small intervals around the zeros of the Dirichlet Kernel, define a function as the sign
of the Dirichlet Kernel outside this pathological set, and then extend it by continuity.

As usual, one needs a bit of analysis to justify tricks of this sort. To this end, by~\eqref{PAKSw-L4}, we recall that the Dirichlet Kernel $$D_N(x)=
\frac{\sin\big((2N+1)\pi x\big)}{\sin(\pi x)}$$
vanishes in~$[0,1]$ if and only if~$x=x_k:=\frac{k}{2N+1}$, with~$k\in\N\cap[0,2N+1]$.

Hence, we consider the intervals~$I_{k,N}:=\left(x_k-e^{-N},x_k+e^{-N}\right)$ and we set
$$ Z_N:=\bigcup_{k\in\N\cap[0,2N+1]}\big( [0,1]\cap I_{k,N}\big).$$
We observe that the Lebesgue measure of~$Z_N$ is bounded by~$4(N+1)e^{-N}$.

We define
$$ [0,1]\setminus Z_N\ni x\mapsto \widetilde f_N(x):=\begin{dcases} 1&{\mbox{ if }}D_N(x)>0,\\
-1&{\mbox{ if }}D_N(x)<0
\end{dcases}$$
and we then extend continuously (e.g., by linear interpolation inside the intervals~$I_{k,N}$) $\widetilde f_N$ to be a continuous
function in~$[0,1]$ with~$\widetilde f_N(0)=\widetilde  f_N(1)$ (and thus, by periodic extension, to a continuous function in all~$\R$, periodic of period~$1$). Notice that, by this construction, $\|\widetilde  f_N\|_{L^\infty(\R)}=1$.

We now let~$f_N(x):=\widetilde f_N(x_0-x)$ and~$Z^\star_N:=x_0-Z_N=\{ x\in\R$ s.t. $x_0-x\in Z_N\}$.
Recalling~\eqref{023wefv4567uygfdfgyhuizo0iwhfg0eb5627e01.pre}, as well as
Exercises~\ref{fr12}, \ref{LGEGMCTLFDCD5D325E28A9N}, and~\ref{K-3PIO},
we see that
\begin{eqnarray*}
|S_{N,f_N}(x_0)|&=&\left|\int_0^1 f_N(y)\,D_N(x_0-y)\,dy\right|
\\&\ge&\left|\int_{[0,1]\setminus Z_N^\star} f_N(y)\,D_N(x_0-y)\,dy\right|- \int_{Z_N^\star} |f_N(y)|\,|D_N(x_0-y)|\,dy\\&=&
\left|\int_{[0,1]\setminus Z_N^\star} \widetilde f_N(x_0-y)\,D_N(x_0-y)\,dy\right|- \int_{Z_N^\star} |f_N(y)|\,|D_N(x_0-y)|\,dy
\\&\ge&\int_{[0,1]\setminus Z_N^\star} |D_N(x_0-y)|\,dy- \int_{Z_N^\star}|D_N(x_0-y)|\,dy
\\&\ge&\int_0^1 |D_N(x_0-y)|\,dy- 2(2N+1)|Z_N|
\\&\ge&c\ln N- 8(N+1)(2N+1)e^{-N}.
\end{eqnarray*} Taking the limit as~$N\to+\infty$, the claim in~\eqref{TACLCSIMOMNDADEKPD.2} follows.

\paragraph{Solution to Exercise~\ref{SRTAMAYBHJAOPLKAZXAQNBKNEMINSNC.L1pe}.}
This is similar to, but slightly different than, Exercise~\ref{SRTAMAYBHJAOPLKAZXAQNBKNEMINSNC}.

We have that, for every~$f\in L^1_{\text{per}}$ with~$\|f\|_{L^1((0,1))}=1$ and any~$k\in\Z$,
$$ |\widehat f_k|\le \int_0^1|f(x)|\,dx=1$$
and therefore
$$ \|S_{N,f}\|_{L^1((0,1))}=\int_0^1\left|\sum_{{k\in\Z}\atop{|k|\le N}}\widehat f_k\,e^{2\pi ikx}\right|\,dx
\le\int_0^1\sum_{{k\in\Z}\atop{|k|\le N}}|\widehat f_k|\,dx\le 2N+1,$$
which establishes~\eqref{TACLCSIMOMNDADEKPD.1pe}.

Now we pick~$\epsilon\in\left(0,\frac14\right)$ and a function~$\psi\in C^\infty_0\left(\left(-\frac14,\frac14\right),\,[0,+\infty)\right)$
with
$$\int_\R\psi(x)\,dx=1.$$
We also let
$$ \psi_\epsilon(x):=\frac1\epsilon\psi\left(\frac{x}\epsilon\right).$$
We extend~$\psi_\epsilon$ outside~$\left(-\frac12,\frac12\right)$ to have a function periodic of period~$1$.
In this way, $\psi_\epsilon\in L^1_{\text{per}}$ and~$\|\psi_\epsilon\|_{L^1((0,1))}=1$.

Also, in view of~\eqref{ERMCBVSUEFASDG}, we remark that the $k$th Fourier coefficient of~$\psi_\epsilon$ is
$$\widehat\psi_{\epsilon,k}=\frac1\epsilon\int_{-1/2}^{1/2} \psi\left(\frac{x}\epsilon\right)\,e^{-2\pi ikx}\,dx
=\int_{-1/4}^{1/4} \psi(y)\,e^{-2\pi ik\epsilon y}\,dy$$ and consequently, by the Dominated Convergence Theorem,
$$ \lim_{\epsilon\searrow0}\widehat\psi_{\epsilon,k}=\int_{-1/4}^{1/4} \psi(y)\,dy=1.$$

As a result,
\begin{eqnarray*}
\sup_{{f\in L^1_{\text{per}}}\atop{\|f\|_{L^1((0,1))}=1}}\|S_{N,f}\|_{L^1((0,1))}&\ge&
\|S_{N,\psi_\epsilon}\|_{L^1((0,1))}
=\int_{-1/2}^{1/2}\left|\sum_{{k\in\Z}\atop{|k|\le N}}\widehat \psi_{\epsilon,k}\,e^{2\pi ikx} \right|\,dx.
\end{eqnarray*}

Hence, by Fatou's Lemma,
\begin{eqnarray*}&&
\sup_{{f\in L^1_{\text{per}}}\atop{\|f\|_{L^1((0,1))}=1}}\|S_{N,f}\|_{L^1((0,1))}\ge\lim_{\epsilon\searrow0}\int_{-1/2}^{1/2}\left|\sum_{{k\in\Z}\atop{|k|\le N}}\widehat \psi_{\epsilon,k}\,e^{2\pi ikx} \right|\,dx\\&&\qquad\ge\int_{-1/2}^{1/2}\left|\sum_{{k\in\Z}\atop{|k|\le N}}\lim_{\epsilon\searrow0}\widehat \psi_{\epsilon,k}\,e^{2\pi ikx} \right|\,dx=\int_{-1/2}^{1/2}\left|\sum_{{k\in\Z}\atop{|k|\le N}}e^{2\pi ikx} \right|\,dx.
\end{eqnarray*}
We now recognise the Dirichlet Kernel presented on page~\pageref{PAKSw-L4} and infer from
Exercise~\ref{K-3PIO} that
$$ \sup_{{f\in L^1_{\text{per}}}\atop{\|f\|_{L^1((0,1))}=1}}\|S_{N,f}\|_{L^1((0,1))}\ge
\int_{-1/2}^{1/2} |D_N(x)|\,dx\ge c\ln N,$$
from which one obtains~\eqref{TACLCSIMOMNDADEKPD.2pe}.

\paragraph{Solution to Exercise~\ref{KDKAMSDW:SDCMILCSKMc}.} We will establish a stronger result, namely that the
set $$ {\mathcal{Z}}:=\left\{ f\in L^1_{\text{per}} {\mbox{ s.t. }} \sup_{N\in\N}\|S_{N,f}\|_{L^1((0,1))}<+\infty\right\}$$
has empty interior in~$L^1_{\text{per}}$.

To prove this, we modify the argument in Theorem~\ref{BaireCategoryTh}.
Namely, for every~$j\in\N$, we let
$$ {\mathcal{Z}}_j:=\left\{ f\in L^1_{\text{per}} {\mbox{ s.t. }} \sup_{N\in\N}\|S_{N,f}\|_{L^1((0,1))}\le j\right\}.$$
Since
$$ {\mathcal{Z}}=\bigcup_{j\in\N}{\mathcal{Z}}_j,$$
in light of the Baire Category Theorem (see e.g.~\cite[Theorem~5.6]{MR924157}) it suffices to show that 
\begin{equation}\label{CSDVAF.pe}
{\mbox{each set~${\mathcal{Z}}_j$ is closed and has empty interior in~$L^1_{\text{per}}$.}}\end{equation}

For this, let
$$ {\mathcal{Z}}_{j,N}:=\left\{ f\in C_{\text{per}} {\mbox{ s.t. }}\|S_{N,f}\|_{L^1((0,1))}\le j\right\}.$$
Let~$f\in L^1_{\text{per}}\setminus {\mathcal{Z}}_{j,N}$ and~$\epsilon>0$. Then, $ \|S_{N,f}\|_{L^1((0,1))}>j$
and if~$g\in L^1_{\text{per}}$ with~$0<\|f-g\|_{L^1((0,1))}\le\epsilon$, setting~$h:=\frac{f-g}{\|f-g\|_{L^1((0,1))}}$ we have 
that
\begin{eqnarray*}
&&\|S_{N,g}\|_{L^1((0,1))}\ge\|S_{N,f}\|_{L^1((0,1))}-\|S_{N,f-g}\|_{L^1((0,1))}
\\&&\qquad=\|S_{N,f}\|_{L^1((0,1))}-\|f-g\|_{L^1((0,1))}\,\left\|S_{N,h}\right\|_{L^1((0,1))}
\ge \|S_{N,f}\|_{L^1((0,1))}-\epsilon\sigma_N,
\end{eqnarray*}
where
$$\sigma_N:=\sup_{{\phi\in L^1_{\text{per}}}\atop{\|\phi\|_{L^1((0,10)}=1}} \|S_{N,\phi}\|_{L^1((0,1))}.$$
We stress that~$\sigma_N<+\infty$, thanks to~\eqref{TACLCSIMOMNDADEKPD.1pe}.

Therefore, if~$\epsilon\in\left(0,\frac{\|S_{N,f}\|_{L^1((0,1))}-j}{\sigma_N}\right)$ we have that~$\|S_{N,g}\|_{L^1((0,1))}>j$, whence~$g\in L^1_{\text{per}}\setminus{\mathcal{Z}}_{j,N}$.

This yields that~$L^1_{\text{per}}\setminus{\mathcal{Z}}_{j,N}$ is open and accordingly~${\mathcal{Z}}_{j,N}$ is closed in~$L^1_{\text{per}}$.

On this account, the set~${\mathcal{Z}}_j$ (which is the countable intersection of the closed sets~${\mathcal{Z}}_{j,N}$
over~$N\in\N$) is also closed.

Thus, to complete the proof of~\eqref{CSDVAF.pe}, we need to show that~${\mathcal{Z}}_j$ has empty interior in~$L^1_{\text{per}}$.
To this end, we argue by contradiction and assume that there exist~$j\in\N$,
$f_j\in {\mathcal{Z}}_j$ and~$\epsilon_j>0$ such that if~$\psi\in L^1_{\text{per}}$ with~$\|f_j-\psi\|_{L^1((0,1))}\le\epsilon_j$, then~$\psi\in{\mathcal{Z}}_j$. 

Hence, for every~$f\in L^1_{\text{per}}$ with~$\|f\|_{L^1((0,1))}=1$, we take~$\psi_j:=f_j-\epsilon_j f$. In this way, we have that
$$\|f_j-\psi_j\|_{L^1((0,1))}=\epsilon_j\|f\|_{L^1((0,1))}=\epsilon_j$$
and therefore~$\psi_j\in{\mathcal{Z}}_j$.

In virtue of this fact,
$$ j\ge\sup_{N\in\N}\|S_{N,\psi_j}\|_{L^1((0,1))}\ge \sup_{N\in\N}\|S_{N,\epsilon_j f}\|_{L^1((0,1))}- \sup_{N\in\N}\|S_{N,f_j}\|_{L^1((0,1))}
\ge\epsilon_j \sup_{N\in\N}\|S_{N,f}\|_{L^1((0,1))}-j$$
and thus, for all~$N\in\N$,
$$ \|S_{N,f}\|_{L^1((0,1))}\le\frac{2j}{\epsilon_j}.$$
Since this is valid for every~$f\in L^1_{\text{per}}$ with~$\|f\|_{L^1((0,1))}=1$, we conclude that
$$ \sup_{{f\in L^1_{\text{per}}}\atop{\|f\|_{L^1((0,1))}=1}} \|S_{N,f}\|_{L^1((0,1))}\le\frac{2j}{\epsilon_j},$$
but this is in contradiction with~\eqref{TACLCSIMOMNDADEKPD.2pe}.

\paragraph{Solution to Exercise~\ref{CLURVB}.}
Suppose not. Then there exists a sequence~$\{\alpha_k\}_{k\in\N}$ such that~$\alpha_k\ge\beta$ and
$$\lim_{k\to+\infty}\left|\int_{\alpha_k}^1\frac{\cos(\pi t)}{t}\,dt\right|=+\infty.$$
Now, if~$\alpha_k$ were bounded, we could extract a sub-sequence~$\alpha_{k_j}$ such that~$\alpha_{k_j}\to\alpha_\star\in[\beta,+\infty)$ as~$j\to+\infty$, but this would entail that
$$ +\infty=\lim_{j\to+\infty}\left|\int_{\alpha_{k_j}}^1\frac{\cos(\pi t)}{t}\,dt\right|
=\left|\int_{\alpha_\star}^1\frac{\cos(\pi t)}{t}\,dt\right|,$$
which is a contradiction.

Hence, $\alpha_k$ must be necessarily unbounded. This gives that there exists a divergent sub-sequence~$\alpha_{k_j}$
and consequently
\begin{eqnarray*}
+\infty&=&\lim_{j\to+\infty}\left|\int_{\alpha_{k_j}}^1\frac{\cos(\pi t)}{t}\,dt\right|\\&
=&\left|\int^{+\infty}_1\frac{\cos(\pi t)}{t}\,dt\right|\\
&=&\left| \int^{+\infty}_1\frac{\sin(\pi t)}{\pi t^2}\,dt\right|\\
&\le&\int^{+\infty}_1\frac{dt}{\pi t^2},
\end{eqnarray*}
which is finite, and this is a contradiction too.

\section{Solutions to selected exercises of Section~\ref{ANY}}

\paragraph{Solution to Exercise~\ref{013p24rtg.01324r5t4y56.0}.}
Let~$q\in\Q\cap(0,+\infty)$.
We observe that~$x+q\in\Q$ if and only if~$x\in\Q$, from which the desired result follows.

\paragraph{Solution to Exercise~\ref{013p24rtg.01324r5t4y56}.}
We denote by~${\mathcal{P}}$ the set in~\eqref{013p24rtg.01324r5t4y56.e} and we stress that, since $f$ is periodic, ${\mathcal{P}}\ne\varnothing$. As a result, we can define
$$T_0:=\inf{\mathcal{P}}$$
and we need to show that this infimum is attained.

For this, we take a sequence of periods~$T_j\in{\mathcal{P}}$ such that~$T_j\to T_0$ as~$j\to+\infty$. Since~$f$ is continuous, for every~$x\in\R$ we have that
$$ f(x+T_0)=\lim_{j\to+\infty}f(x+T_j)=\lim_{j\to+\infty}f(x)=f(x),$$
hence, to show that~$T_0\in{\mathcal{P}}$ we only need to check that~$T_0>0$.

Suppose, by contradiction, that\begin{equation}\label{q83v6301740wofjhgrnIP}T_0=0.\end{equation} Pick an arbitrary point~$x\in\R$. Let~$m_{j,x}\in\Z$ be such that
$$ m_{j,x}\,T_j\le x<(m_{j,x}+1)\,T_j.$$
Then,
$$ |x-m_{j,x}\,T_j|=x-m_{j,x}\,T_j\le (m_{j,x}+1)\,T_j-m_{j,x}\,T_j=T_j,$$
which is infinitesimal as~$j\to+\infty$, due to~\eqref{q83v6301740wofjhgrnIP}.

As a consequence, 
$$\lim_{j\to+\infty}m_{j,x}\,T_j=x$$
and the continuity\footnote{If~$f$ is not continuous, this argument cannot work
and, in fact, a discontinuous periodic function may have arbitrarily small periods, as showcased in Exercise~\ref{013p24rtg.01324r5t4y56.0}. See however Exercise~\ref{013p24rtg.01324r5t4y56.bb} for a counterpart in Lebesgue spaces.} of~$f$ entails that
$$ f(x)=\lim_{j\to+\infty}f(m_{j,x}\,T_j)=f(0).$$
This says that~$f$ is constant, in contradiction\footnote{The assumption that~$f$ is nonconstant is essential, since constant functions have arbitrarily small periods.}
with our assumption.

\paragraph{Solution to Exercise~\ref{013p24rtg.01324r5t4y56.b}.}
If~$T>0$ is any period of~$f$, by the minimality of~$T_0$ we know that~$T\ge T_0$. Hence, we take~$k\in\N\cap[1,+\infty)$ such that
$$ kT_0\le T<(k+1)T_0$$
and we set~$\tau:= T-kT_0$.

Notice that~$\tau\ge0$ and, for all~$x\in\R$,
$$ f(x+\tau)=f(x+T-kT_0)=f(x-kT_0)=f(x),$$
showing that~$\tau$ is a period for~$f$. The minimality of~$T_0$ thereby implies that~$\tau\ge T_0$, unless~$\tau=0$.

But~$\tau<(k+1)T_0-kT_0=T_0$ and therefore~$\tau=0$. Consequently~$T=kT_0$, as advertised.

\paragraph{Solution to Exercise~\ref{013p24rtg.01324r5t4y56.bb}.} Let~$g:=f\chi_{[-3\tau,3\tau]}$ and note that~$g\in L^1(\R)$.
Hence, by the continuity of the translations in~$L^1(\R)$ (see e.g.~\cite[Lemma 4.3]{MR2759829}),
\begin{equation}\label{CTL1iqekrjgt0-23p4tyhalpd3n39fvw-3}
\lim_{h\to0}\int_\R |g(t+h)-g(t)|\,dt=0.\end{equation}

We now consider the set~${\mathcal{L}}$ of Lebesgue points for the function~$f$, namely the collection of points~$p\in\R$ for which
\begin{equation}\label{nrrjE45:2-ekrjgt0-23p4tyhalpd3n39fvw-3.21}
\lim_{r\searrow0}\frac1r\int_{p-r}^{p+r}| f(t)-f(p)|\,dt=0.\end{equation}
By~\cite[Theorem 7.15]{MR3381284}, we have that~$\R\setminus{\mathcal{L}}$ has null Lebesgue measure.

Up to a horizontal translation, we may suppose that~$0\in{\mathcal{L}}$. Then, the desired result will be proved by showing that either the set in~\eqref{013p24rtg.01324r5t4y56.e} 
admits a minimum or, for every~$x\in{\mathcal{L}}$,
\begin{equation}\label{CTL1iqekrjgt0-23p4tyhalpd3n39fvw-3.2}
f(x)=f(0)=:c.
\end{equation}
To check this, we let~$T_0$ be the infimum in~\eqref{013p24rtg.01324r5t4y56.e} and
take a sequence of periods~$T_j\in{\mathcal{P}}$ such that~$T_j\to T_0$ as~$j\to+\infty$. By construction~$T_0\in[0,\tau]$.
Also, by means of~\eqref{CTL1iqekrjgt0-23p4tyhalpd3n39fvw-3} and the fact that~$f$ is periodic of period~$T_j$,
\begin{eqnarray*}0&=&\lim_{j\to+\infty}\int_\R |g(t+T_0-T_j)-g(t)|\,dt\\&=&
\lim_{j\to+\infty}\int_{\R} \big|f(t+T_0-T_j)\chi_{[-3\tau,3\tau]}(y+T_0-T_j)-f(t)\chi_{[-3\tau,3\tau]}(t)\big|\,dt\\&\geq&
\lim_{j\to+\infty}\int_{-2\tau}^{2\tau} \big|f(t+T_0-T_j)\chi_{[-3\tau,3\tau]}(y+T_0-T_j)-f(t)\chi_{[-3\tau,3\tau]}(t)\big|\,dt\\&=&
\lim_{j\to+\infty}\int_{-2\tau}^{2\tau} |f(t+T_0-T_j)-f(t)|\,dt\\&=&
\lim_{j\to+\infty}\int_{-2\tau}^{2\tau} |f(t+T_0)-f(t)|\,dt.\end{eqnarray*}
This yields that~$f(t+T_0)=f(t)$ for a.e.~$t\in[-2\tau,2\tau]$ and consequently, by the fact that~$f$ is also
periodic of period~$\tau$,
for a.e.~$t\in\R$, showing that~$T_0$ is a period for~$f$ as well.

Now we consider two scenarios. First, if~$T_0>0$, then~$T_0\in{\mathcal{P}}$. Thus,~$T_0$ is the minimal element of~$
{\mathcal{P}}$ and we are done.

We can thereby focus on the casee in which
\begin{equation}\label{23d4q83v6301740wofjhgrnIP}
T_0=0.\end{equation}
We pick an arbitrary point~$x\in{\mathcal{L}}$. Let~$m_{j,x}\in\Z$ be such that
$$ m_{j,x}\,T_j\le x<(m_{j,x}+1)\,T_j.$$
Then,
$$ |x-m_{j,x}\,T_j|=x-m_{j,x}\,T_j\le (m_{j,x}+1)\,T_j-m_{j,x}\,T_j=T_j,$$
which is infinitesimal as~$j\to+\infty$, due to~\eqref{23d4q83v6301740wofjhgrnIP}.

As a consequence, recalling~\eqref{CTL1iqekrjgt0-23p4tyhalpd3n39fvw-3} and that~$f$ is periodic of period~$T_j$,
\begin{eqnarray*}
0&=&\lim_{j\to+\infty}\int_\R |g(t-x+m_{j,x}\,T_j)-g(t)|\,dt\\&=&\lim_{j\to+\infty}\int_{\R} 
\big|f(t-x+m_{j,x}\,T_j)\chi_{[-3\tau,3\tau]}(t-x+m_{j,x}\,T_j)-f(t)\chi_{[-3\tau,3\tau]}(t)\big|\,dt
\\&\ge&\lim_{j\to+\infty}\int_{-2\tau}^{2\tau} 
\big|f(t-x+m_{j,x}\,T_j)\chi_{[-3\tau,3\tau]}(t-x+m_{j,x}\,T_j)-f(t)\chi_{[-3\tau,3\tau]}(t)\big|\,dt
\\&=&\lim_{j\to+\infty}\int_{-2\tau}^{2\tau} 
|f(t-x+m_{j,x}\,T_j)-f(t)|\,dt
\\&=&\lim_{j\to+\infty}\int_{-2\tau}^{2\tau} 
|f(t-x)-f(t)|\,dt.\end{eqnarray*}
On this account, we infer that~$f(t-x)=f(t)$ for a.e.~$t\in[-2\tau,2\tau]$ and consequently, by the periodicity of~$f$ with period~$\tau$,
for a.e.~$t\in\R$.

Therefore, since both~$0$ and~$x$ belong to~${\mathcal{L}}$, we can employ~\eqref{nrrjE45:2-ekrjgt0-23p4tyhalpd3n39fvw-3.21}
and conclude that
$$ f(x)=\lim_{r\searrow0}\frac1r\int_{x-r}^{x+r}f(t)\,dt=
\lim_{r\searrow0}\frac1r\int_{x-r}^{x+r}f(t-x)\,dt=\lim_{r\searrow0}\frac1r\int_{-r}^{r}f(s)\,ds=f(0),$$
which proves~\eqref{CTL1iqekrjgt0-23p4tyhalpd3n39fvw-3.2}, as desired.

\paragraph{Solution to Exercise~\ref{QUASDFG-qaudwe4y6.PREC}.}
Suppose, for the sake of contradiction, that~$\phi$ is nonconstant. Then,
due to Exercise~\ref{013p24rtg.01324r5t4y56}, we can consider its fundamental period~$T_0$.
We also know from Exercise~\ref{013p24rtg.01324r5t4y56.b}
that~$T_1=N_1 \,T_0$ and~$T_2=N_2 \,T_0$ for some~$N_1$, $N_2\in\N\cap[1,+\infty)$.

As a result, $ \frac{T_1}{T_2}=\frac{N_1 \,T_0}{N_2 \,T_0}=\frac{N_1 }{N_2 }\in\Q$, in contradiction with our assumption.

\paragraph{Solution to Exercise~\ref{QUASDFG-qaudwe4y6}.}
We suppose, for the sake of contradiction, that~$g$ is periodic and we pick any period~$T>0$. Then, we define
$$F(x):=f_1(x+T)-f_1(x).$$
We stress that
\begin{equation}\label{LP:02rt-peert0-2}
{\mbox{$F$ is continuous}}
\end{equation} and
\begin{equation}\label{LP:02rt-peert0-1}
{\mbox{$F$ is periodic of period~$T_1$.}}
\end{equation}

Also, for all~$x\in\R$,
\begin{eqnarray*}&&
F(x)=g(x+T)-g(x)+f_2(x)-f_2(x+T)=f_2(x)-f_2(x+T),
\end{eqnarray*}
showing that
\begin{equation*}
{\mbox{$F$ is periodic of period~$T_2$ as well.}}
\end{equation*}

This, \eqref{LP:02rt-peert0-2}, and~\eqref{LP:02rt-peert0-1} allow us to use
Exercise~\ref{QUASDFG-qaudwe4y6.PREC} and gather that~$F$ is constant,
namely there exists~$C\in\R$ such that~$F(x)=C$ for all~$x\in\R$. 

We then observe that, by the periodicity of~$f_1$,
\begin{equation*}\begin{split}&
C=\frac1{T_1}\int_0^{T_1}F(x)\,dx=\frac1{T_1}
\left( \int_0^{T_1}f_1(x+T)\,dx-\int_0^{T_1}f_1(x)\,dx\right)\\&\qquad\qquad\quad
=\frac1{T_1}\left( \int_{T}^{T+T_1}f_1(x)\,dx-\int_0^{T_1}f_1(x)\,dx\right)
=0.\end{split}\end{equation*}

{F}rom this, we conclude that, for all~$x\in\R$,
$$f_1(x+T)=f_1(x)+F(x)=f_1(x)$$
and also
$$ f_2(x+T)=f_2(x)-F(x)=f_2(x),$$
showing that
\begin{equation}\label{LP:02rt-peert0-1asVfJ45t-c}
{\mbox{$f_1$ and~$f_2$ are also periodic of period~$T$.}}
\end{equation}

Now, by Exercise~\ref{013p24rtg.01324r5t4y56}, we take the fundamental periods~$T_1^\star$ and~$T_2^\star$
for~$f_1$ and~$f_2$, respectively.

By Exercise~\ref{013p24rtg.01324r5t4y56.b}, we know that
that~$T_1=N_1 \,T_1^\star$ and~$T_2=N_2 \,T_2^\star$ for some~$N_1$, $N_2\in\N\cap[1,+\infty)$ and,
moreover, by~\eqref{LP:02rt-peert0-1asVfJ45t-c}, it holds that~$T=M_1 \,T_1^\star$ and~$T=M_2 \,T_2^\star$ for some~$M_1$, $M_2\in\N\cap[1,+\infty)$.

All in all,
\begin{eqnarray*}
\frac{T_1}{T_2}=\frac{N_1 \,T_1^\star}{N_2 \,T_2^\star}=
\frac{N_1 \,M_1\,M_2\,T_1^\star}{N_2 \,M_1\,M_2\,T_2^\star}=\frac{N_1 \,M_2\,T}{N_2 \,M_1 \,T}=\frac{N_1 \,M_2}{N_2 \,M_1}\in\Q,
\end{eqnarray*}
reaching the desired contradiction.

\begin{figure}[h]
\includegraphics[width=0.3\textwidth]{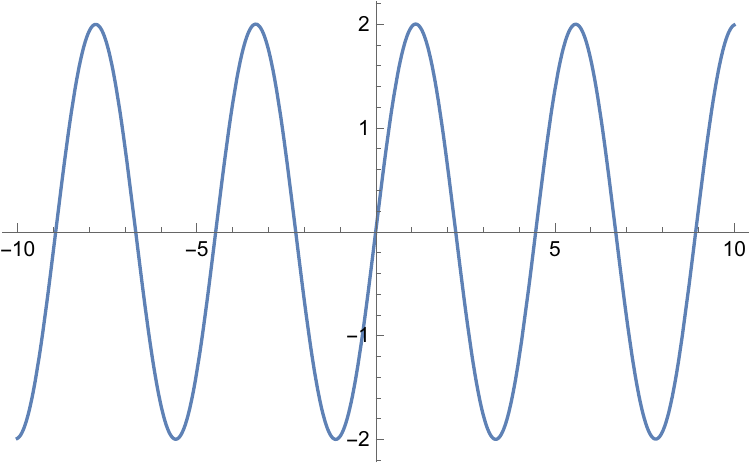}$\quad$
\includegraphics[width=0.3\textwidth]{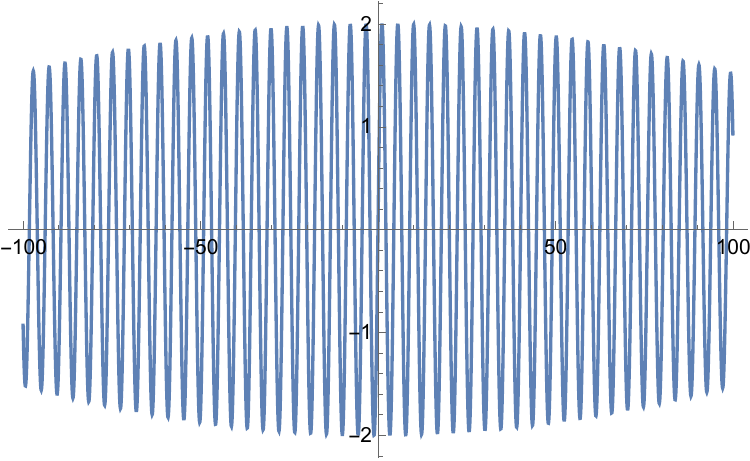}$\quad$
\includegraphics[width=0.3\textwidth]{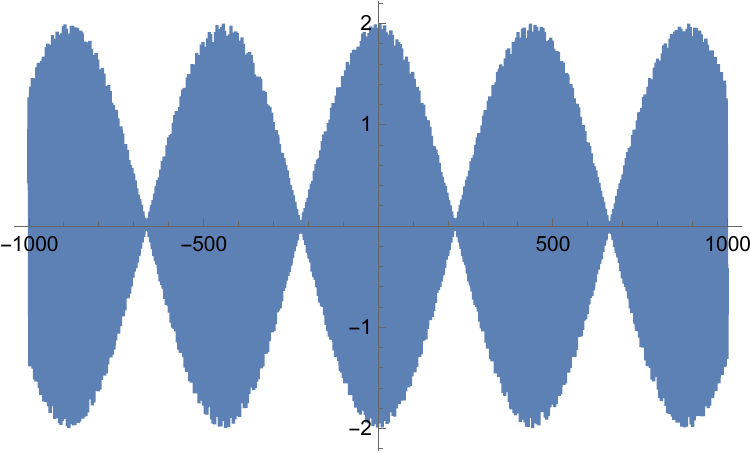}
\centering
\caption{Not a periodic function:
plot of the function~$\sin(\sqrt{2}x)+\sin\left(\frac{7x}5\right)$
with~$x\in[-\ell,\ell]$, where~$\ell\in\{10,\,100,\,1000\}$.}\label{0283ir9fq0fkne4kpjne.023oriSTQEILSR8hu63k.lop1-wdCNMKtk}
\end{figure}

Parenthetically, we point out that addressing the periodicity of a function just by looking at a picture may be tricky,
since different scales and patterns may pop up, see e.g. Fiigure~\ref{0283ir9fq0fkne4kpjne.023oriSTQEILSR8hu63k.lop1-wdCNMKtk}.

For more information about this topic  and related ones, see the discussion on almost periodic functions in~\cite{NEXT}.

\section{Solutions to selected exercises of Section~\ref{SEC:ANUYDIMEAGa}}

\paragraph{Solution to Exercise~\ref{PRECON56789}.}
Since~$E$ is bounded, there exists~$R>0$ such that~$E\subseteq B_R$. As a result, if~$x\in E$ and~$|m|\ge2R$,
$$ |2(x+m)|\ge 2(|m|-|x|)\ge2(|m|-R)\ge2R,$$
giving that~$2(x+m)\not\in E$ and thus~$\chi_E(2(x+m))=0$.

This observation yields that only a finite number of terms (namely, the integers~$m$ with~$|m|<2R$) can contribute to the series defining~$f$, as desired.

In addition, for every~$\ell\in\Z^n$, setting~$j:=m+\ell$ we have that
$$ f(x+\ell)=\sum_{m\in\Z^n}\chi_E(2(x+\ell+m))=\sum_{j\in\Z^n}\chi_E(2(x+j))=f(x),$$ establishing that~$f$ is $\Z^n$-periodic.

It remains to calculate the Fourier coefficients of~$f$. For this, we observe that
$$ \bigcup_{m\in\Z^n}\Big((0,1)^n+m\Big)=\R^n,$$
up to a set of null Lebesgue measure and we proceed as follows:
\begin{eqnarray*}\widehat f_k&=&\sum_{m\in\Z^n}\int_{(0,1)^n}\chi_E(2(x+m))\,e^{-2\pi ik\cdot x}
\,dx\\&=&
\sum_{m\in\Z^n}\int_{(0,1)^n+m}\chi_E(2y)\,e^{-2\pi ik\cdot (y-m)}\,dy
\\&=&\int_{\R^n}\chi_E(2y)\,e^{-2\pi ik\cdot y}\,dy
\\&=&\frac1{2^n}\int_{\R^n}\chi_E(x)\,e^{-\pi ik\cdot x}\,dx
\\&=&\frac1{2^n}\int_{E}e^{-\pi ik\cdot x}\,dx.\end{eqnarray*}

\paragraph{Solution to Exercise~\ref{PRECON56789-BIS}.} We have that
\begin{eqnarray*}
\int_{(0,1)^n}|f(x)|^2\,dx&=&\sum_{m,h\in\Z^n}\int_{(0,1)^n}\chi_E(2(x+m))\,\chi_E(2(x+h))\,dx\\&=&\sum_{m,h\in\Z^n}\int_{(0,1)^n+m}\chi_E(2y)\,\chi_E(2(y+h-m))\,dy\\&=&\sum_{m,j\in\Z^n}\int_{(0,1)^n+m}\chi_E(2y)\,\chi_E(2(y+j))\,dy\\&=&\sum_{j\in\Z^n}\int_{\R^n}\chi_E(2y)\,\chi_E(2(y+j))\,dy\\&=&\frac1{2^n}\sum_{j\in\Z^n}\int_{\R^n}\chi_E(x)\,\chi_E(x+2j)\,dx,
\end{eqnarray*}
where we have used that the above series is a finite sum to freely swap the integral and summation signs.

As a result,
$$ \|f\|_{L^2((0,1)^n)}=\sqrt{\frac1{2^n}\sum_{j\in\Z^n}\int_{\R^n}\chi_E(x)\,\chi_E(x+2j)\,dx}.$$

\paragraph{Solution to Exercise~\ref{PRECON56789-BIS.lo}.}
Let~$j\in\{1,\dots,n\}$ be such that
$$ |k_j|=\max\big\{|k_1|,\dots,|k_n|\big\}.$$
Pick~$\mu_m:=\ell\delta_{j,m}$ and note that
\begin{equation} \label{PRECON56789-BIS.lo2}\prod_{m=1}^n |k_m|^{\mu_m}=|k_j|^\ell.\end{equation}
Besides,
$$ \sum_{m=1}^nk_1^2\le n k_j^2.$$
This and~\eqref{PRECON56789-BIS.lo2} yield the desired result.

\section{Solutions to selected exercises of Section~\ref{ODCABRBAGN}}

\paragraph{Solution to Exercise~\ref{GABNMSOIPEDCFQMw0-1}.}
If~$k=h$ the result is obvious, so we suppose that~$k-h\ne0$.
We let~$\omega:=e^{\frac{2\pi i (k-h)}N}$ and we observe that~$0<|k-h|<N$,
whence~$\omega\ne1$.

Moreover,
$$ \omega^N=e^{{2\pi i (k-h)}}=1$$
and, as a result,
$$ \sum_{j=0}^{N-1} e^{\frac{2\pi ij(k-h)}N}=\sum_{j=0}^{N-1} \omega^j=
\frac{1-\omega^N}{1-\omega}=0.$$

\section{Solutions to selected exercises of Section~\ref{TIDES}}

\paragraph{Solution to Exercise~\ref{TIDEX}.} For all~$m\in\Z$ with~$|m|\le N$, we have that
\begin{eqnarray*}
\int_{t_0}^{t_0+T}h(t)\,e^{-i\omega_m t}\,dt
&=&\sum_{{{k\in\Z}\atop{|k|\le N}}} \int_{t_0}^{t_0+T}\widehat h_{k}\,e^{i(\omega_k-\omega_m)t}\,dt\\&=&\widehat h_{m}\,T+
\sum_{{{k\in\Z}\atop{|k|\le N}}\atop{k\ne m}} \frac{\widehat h_{k}\,\big(e^{i(\omega_k-\omega_m)(t_0+T)}-e^{i(\omega_k-\omega_m)t_0}\big)}{i(\omega_k-\omega_m)}.
\end{eqnarray*}
The desired result follows, since~$|e^{i(\omega_k -\omega_m)(t_0+T)}-e^{i(\omega_k -\omega_m)t_0}|\le
|e^{i(\omega_k -\omega_m)(t_0+T)}|+|e^{i(\omega_k -\omega_m)t_0}|=1+1=2$.

\section{Solutions to selected exercises of Section~\ref{PDESE}}

\paragraph{Solution to Exercise~\ref{UNIZER}.} Suppose, by contradiction that~$w$ does not vanish identically.
Therefore, there exists~$(x_\star,t_\star)\in[0,\ell]\times[0,T)$ such that~$w(x_\star,t_\star)\ne0$.

Let us suppose that~$w(x_\star,t_\star)>0$, the case~$w(x_\star,t_\star)<0$ being similar. Then, we pick~$\epsilon>0$, to be chosen conveniently small, and define
$$ w_\epsilon(x,t):=w(x,t)-\epsilon\,\left(1-e^{-\frac{\pi^2 t}{\ell^2}}\right)\,\sin\left(\frac{\pi x}\ell\right).$$
We notice that
$$ w_\epsilon(x_\star,t_\star)\ge w(x_\star,t_\star)-\epsilon\,\left(1-e^{-\frac{\pi^2 t_\star}{\ell^2}}\right)\ge\frac{w(x_\star,t_\star)}2>0,$$
as long as~$\epsilon$ is sufficiently small, and therefore there exists a point~$(x_\epsilon,t_\epsilon)\in[0,\ell]\times[0,t_\star]$ such that
$$ w_\epsilon(x_\epsilon,t_\epsilon)=\max_{[0,\ell]\times[0,t_\star]} w_\epsilon>0.$$

We stress that~$w_\epsilon(0,t_\epsilon)=w(0,t_\epsilon)=0=w(\ell,t_\epsilon)=w_\epsilon(\ell,t_\epsilon)$ and therefore~$x_\epsilon\in(0,\ell)$. 
Additionally, $w_\epsilon(x_\epsilon,0)=w(x_\epsilon,0)=0$ and consequently~$t_\epsilon\in(0,t_\star]$.

As a result, by the above maximal property, we have that$$\partial_{xx}w_\epsilon(x_\epsilon,t_\epsilon)\le0\qquad{\mbox{and}}\qquad\partial_tw_\epsilon(x_\epsilon,t_\epsilon)\ge0.$$
For this reason,
\begin{eqnarray*}
0&\le&\partial_tw_\epsilon(x_\epsilon,t_\epsilon)-\partial_{xx}w_\epsilon(x_\epsilon,t_\epsilon)\\
&=&\partial_tw(x_\epsilon,t_\epsilon)
-\frac{\pi^2\epsilon}{\ell^2}\,e^{-\frac{\pi^2 t_\epsilon}{\ell^2}}\,\sin\left(\frac{\pi x_\epsilon}\ell\right)
-\partial_{xx}w(x_\epsilon,t_\epsilon)
-\frac{\pi^2\epsilon}{\ell^2}\,\left(1-e^{-\frac{\pi^2 t_\epsilon}{\ell^2}}\right)\,\sin\left(\frac{\pi x_\epsilon}\ell\right)\\
&=&-\frac{\pi^2\epsilon}{\ell^2}\sin\left(\frac{\pi x_\epsilon}\ell\right),
\end{eqnarray*}
which is strictly negative and we have thereby reached the desired contradiction.

\paragraph{Solution to Exercise~\ref{NSSEEV}.}
In the notation of Section~\ref{SEC:ANUYDIMEAGa}, we consider formal Fourier Series expansions
in the space variable~$x$ for the velocity~$v(x,t)$ and the pressure~$p(x,t)$, at any given time~$t>0$, namely
\begin{equation}\label{WTTGSRIOMGFR}
\begin{split}
&v(x,t)=\sum_{k\in\Z^3} \widehat v_k(t)\,e^{2\pi ik\cdot x}\\
{\mbox{and }}\qquad&p(x,t)=\sum_{k\in\Z^3} \widehat p_k(t)\,e^{2\pi ik\cdot x}.
\end{split}\end{equation}
Since~$v$ takes values in~$\R^3$, we have that~$\widehat v_k$ is a function from~$t$ to~$\R^3$,
say~$(\widehat v_{1,k}(t),\widehat v_{2,k}(t),\widehat v_{3,k}(t))$,
but we stress that~$p$ is a scalar function, therefore~$\widehat p_k$  is a function from~$t$ to~$\R$.

Also, since each component of~$v$ and~$p$ is real valued, by~\eqref{fasv} we have that
\begin{equation}\label{WTTGSRIOMGFR.mjuT21e3s} \widehat v_{-k}=\overline{\widehat v_k}\qquad{\mbox{and}}\qquad
\widehat p_{-k}=\overline{\widehat p_k}.\end{equation}
In this notation, without any attempt of discussing convergence issues,
we show that the Navier-Stokes equation~\eqref{NS6THTgBplaISw} is formally equivalent to
\begin{equation}\label{NS6THTgBplaISw.FOU}
\begin{dcases}
\displaystyle \frac{d\widehat v_k}{dt}+4\pi^2\nu|k|^2\,\widehat v_k=2\pi i\sum_{h\in\Z^3} (h\cdot\widehat v_{k-h})\left(
\frac{(k\cdot\widehat v_h)k}{|k|^2}-\widehat v_h
\right),\\
k\cdot\widehat v_k=0.
\end{dcases}\end{equation}

To see this, we first observe that the solenoidal velocity condition in~\eqref{NS6THTgBplaISw} formally corresponds to
\begin{eqnarray*}
0=\sum_{m=1}^3\frac{ \widehat v_m}{\partial x_m}=
\sum_{{k\in\Z^3}\atop{1\le m\le3}} \frac{\partial }{\partial x_m}\big(\widehat v_{m,k}(t)\,e^{2\pi ik\cdot x}\big)=
\sum_{{k\in\Z^3} \atop{1\le m\le3}}2\pi i k_m\,\widehat v_{m,k}(t)\,e^{2\pi ik\cdot x}
\end{eqnarray*}
and accordingly
$$\sum_{m=1}^3k_m\,\widehat v_{m,k}(t)=0,$$
which is the second equation in~\eqref{NS6THTgBplaISw.FOU}.

Moreover, the momentum conservation equation in~\eqref{NS6THTgBplaISw} formally corresponds to
the following identity for each~$j\in\{1,2,3\}$:
\begin{eqnarray*}
0&=&\frac{\partial v_j }{\partial t}+\sum_{m=1}^3 v_m\frac{\partial v_j}{\partial x_m} -\nu \,\Delta v_j+\frac{\partial p}{\partial x_j}\\&=&\sum_{k\in\Z^3} \frac{d\widehat v_{j,k}(t)}{dt}\,e^{2\pi ik\cdot x}
+\sum_{{q,h\in\Z^3}\atop{1\le m\le3}} 2\pi ih_m\,\widehat v_{m,q}(t)\,\widehat v_{j,h}(t)\,e^{2\pi i(q+h)\cdot x}
\\&&\qquad+4\pi^2\nu\sum_{k\in\Z^3}|k|^2\, \widehat v_{j,k}(t)\,e^{2\pi ik\cdot x}
+\sum_{k\in\Z^3}2\pi ik_j \widehat p_k(t)\,e^{2\pi ik\cdot x}
\end{eqnarray*}
and therefore, taking the $k$th harmonic,
\begin{equation*}
\frac{d\widehat v_{j,k}}{dt}+\sum_{{h\in\Z^3}\atop{1\le m\le3}} 2\pi ih_m\,\widehat v_{m,k-h}\,\widehat v_{j,h}
+4\pi^2\nu\,|k|^2\, \widehat v_{j,k}+2\pi ik_j \widehat p_k=0,
\end{equation*}
or, equivalently, in vector notation,
\begin{equation}\label{NS6THTgBplaISw.FOU.2}
\frac{d\widehat v_{k}}{dt}+2\pi i\sum_{{h\in\Z^3}} h\cdot\widehat v_{k-h}\,\widehat v_{h}
+4\pi^2\nu\,|k|^2\, \widehat v_{k}+2\pi ik \widehat p_k=0.
\end{equation}

To find the first equation in~\eqref{NS6THTgBplaISw.FOU} it is now convenient to use again
the incompressibility condition to get rid of the pressure dependence in~\eqref{NS6THTgBplaISw.FOU.2}.
Namely, the Navier-Stokes equation~\eqref{NS6THTgBplaISw} also formally implies that
\begin{eqnarray*}
0&=&\nabla\cdot\left(
\frac{\partial v }{\partial t}+(v \cdot \nabla )v -\nu \,\Delta v +\nabla p
\right)\\&=&\sum_{m,r=1}^3\frac{\partial v_m}{\partial x_r}\,
\frac{\partial v_r}{\partial x_m}+\Delta p
\end{eqnarray*}
and then, in Fourier modes,
\begin{eqnarray*}
0&=&-4\pi^2
\sum_{{q,h\in\Z^3}\atop{1\le m,r\le3}} q_r\,h_m\widehat v_{m,q}(t)\,\widehat v_{r,h}(t)
\,e^{2\pi i(q+h)\cdot x}
-4\pi^2\sum_{k\in\Z^3} |k|^2\,\widehat p_k(t)\,e^{2\pi ik\cdot x}.
\end{eqnarray*}
As a result, taking the $k$th mode,
\begin{eqnarray*} 0&=&\sum_{{h\in\Z^3}\atop{1\le m,r\le3}} (k-h)_r\,h_m\,\widehat v_{m,k-h}\,\widehat v_{r,h}
+|k|^2\,\widehat p_k\\&=&
\sum_{{h\in\Z^3}}\big( h\cdot \widehat v_{k-h}\big)\,\big((k-h)\cdot\widehat v_{h}\big)
+|k|^2\,\widehat p_k.
\end{eqnarray*}

From this and the second equation in~\eqref{NS6THTgBplaISw.FOU}, it follows that
\begin{eqnarray*}
\widehat p_k&=&-\frac1{|k|^2}\sum_{{h\in\Z^3}}\big( h\cdot \widehat v_{k-h}\big)\,\big((k-h)\cdot\widehat v_{h}\big)
\\&=&-\frac1{|k|^2}\sum_{{h\in\Z^3}}\big( h\cdot \widehat v_{k-h}\big)\,\big(k\cdot\widehat v_{h}\big).
\end{eqnarray*}
This and~\eqref{NS6THTgBplaISw.FOU.2} yield that
\begin{eqnarray*}0&=&
\frac{d\widehat v_{k}}{dt}
+4\pi^2\nu\,|k|^2\, \widehat v_{k}+2\pi i \left(\sum_{{h\in\Z^3}} h\cdot\widehat v_{k-h}\,\widehat v_{h}
+k\widehat p_k\right)\\&=&
\frac{d\widehat v_{k}}{dt}
+4\pi^2\nu\,|k|^2\, \widehat v_{k}+2\pi i \sum_{{h\in\Z^3}} \left( h\cdot\widehat v_{k-h}\,\widehat v_{h}
-\frac{k}{|k|^2}\big( h\cdot \widehat v_{k-h}\big)\,\big(k\cdot\widehat v_{h}\big)
\right)\\&=&
\frac{d\widehat v_{k}}{dt}
+4\pi^2\nu\,|k|^2\, \widehat v_{k}+2\pi i \sum_{{h\in\Z^3}}\big( h\cdot \widehat v_{k-h}\big) \left( \widehat v_{h}
-\frac{k(k\cdot\widehat v_{h})}{|k|^2}
\right)
\end{eqnarray*} 
and this is the first equation in~\eqref{NS6THTgBplaISw.FOU}.

From the physical point of view the first equation in~\eqref{NS6THTgBplaISw.FOU}
is interesting since it reveals the different scales at which different phenomena take place.
Indeed, very roughly speaking,
the viscosity term contributes in~\eqref{NS6THTgBplaISw.FOU} as a friction of the order of~$\nu|k|^2\,|\widehat v_k|$,
while the inertial contribution coming from the transport term is of the order of~$|k|\,|\widehat v_k|^2$.
It is therefore tempting to conjecture that the viscous terms prevail over the inertial ones for large frequencies,
say when~$\nu|k|\gtrapprox |\widehat v_k|$, and that the influence of viscosity can be somewhat disregarded
up to a sufficiently high frequency scale.
This may suggest, in some sense, that the energy gets dissipated only at high frequency space,
flowing, especially in the absence of viscosity, as a ``cascade'' from some modes to others (this is however a delicate business,
deeply depending on the dimension too, which cannot be resolved simply at a heuristic level).
See also Exercise~\ref{2NSSEEV.bis}, as well as~\cite{MR1428905, MR1872661, MR3929468, MR4475666}
and the references therein for more information about these phenomena.

\paragraph{Solution to Exercise~\ref{2NSSEEV}.}
We will show that
\begin{equation}\label{ksmD9cvbHA034or-1204-11}
{\mathcal{V}}(t)=4\pi^2\sum_{k\in\Z^3}
\sum_{\ell=1}^3|k|^2\, |\widehat v_{\ell,k}(t)|^2={\mathcal{D}}(t).
\end{equation}
Interestingly, one of the consequences of~\eqref{ksmD9cvbHA034or-1204-11}
is that the enstrophy (i.e. a measurement of the squared total vorticity of a fluid) and the overall rate of deformation coincide (very roughly speaking, showing that, in some sense,
a fluid motion is the superposition of an infinitesimal linear displacement, which does not alter
significantly the mutual position of fluid particles, and a swirling).

Let us proceed with the computations needed, arguing at a formal level, disregarding any possible convergence issue.
We use the notation in~\eqref{WTTGSRIOMGFR} and we calculate the first component of the curl of~$v$ as
\begin{eqnarray*}
\frac{\partial v_3}{\partial x_2}(x,t)-\frac{\partial v_2}{\partial x_3}(x,t)=
2\pi i\sum_{k\in\Z^3} \big(k_2\, \widehat v_{3,k}(t)-k_3\, \widehat v_{2,k}(t)\big)\,e^{2\pi ik\cdot x}
\end{eqnarray*}
and similarly one can compute the second and third components of the curl of~$v$ just by permuting indices.

Hence, by~\eqref{WTTGfdgtchofrTLA:0-2} and Parseval's Identity~\eqref{PARS}, {\footnotesize
\begin{equation}\label{ksmD9034or-1204}
\begin{split}
{\mathcal{V}}(t)&=
4\pi^2\sum_{k\in\Z^3} \Big(\big|k_2\, \widehat v_{3,k}(t)-k_3\, \widehat v_{2,k}(t)\big|^2
+\big|k_1\, \widehat v_{3,k}(t)-k_3\, \widehat v_{1,k}(t)\big|^2+
\big|k_2\, \widehat v_{1,k}(t)-k_1\, \widehat v_{2,k}(t)\big|^2\Big)\\
&=4\pi^2\sum_{k\in\Z^3}\Bigg[
\sum_{\ell=1}^3 \left(|k|^2-k_\ell^2\right)\, |\widehat v_{\ell,k}(t)|^2
-k_1k_2\,\big(\widehat v_{1,k}\,\overline{\widehat v_{2,k}}+\overline{\widehat v_{1,k}}\,\widehat v_{2,k}\big)\\&\qquad\qquad\qquad\qquad
-k_1k_3\,\big(\widehat v_{1,k}\,\overline{\widehat v_{3,k}}+\overline{\widehat v_{1,k}}\,\widehat v_{3,k}\big)
-k_2k_3\,\big(\widehat v_{2,k}\,\overline{\widehat v_{3,k}}+\overline{\widehat v_{2,k}}\,\widehat v_{3,k}\big)
\Bigg].
\end{split}
\end{equation}}

Besides, from the second relation in~\eqref{NS6THTgBplaISw.FOU}, {\footnotesize
\begin{eqnarray*}
0&=& \sum_{k\in\Z^3}\big| k_1\, \widehat v_{1,k}+k_2\, \widehat v_{2,k}+k_3\, \widehat v_{3,k}\big|^2\\&=& \sum_{k\in\Z^3}
\sum_{\ell=1}^3k_\ell^2\, |\widehat v_{\ell,k}(t)|^2\\&&\quad
+\sum_{k\in\Z^3}\Big(k_1k_2\,\big(\widehat v_{1,k}\,\overline{\widehat v_{2,k}}+\overline{\widehat v_{1,k}}\,\widehat v_{2,k}\big)
+k_1k_3\,\big(\widehat v_{1,k}\,\overline{\widehat v_{3,k}}+\overline{\widehat v_{1,k}}\,\widehat v_{3,k}\big)
+k_2k_3\,\big(\widehat v_{2,k}\,\overline{\widehat v_{3,k}}+\overline{\widehat v_{2,k}}\,\widehat v_{3,k}\big)\Big).
\end{eqnarray*}}
Adding~$4\pi^2$ times this quantity to~\eqref{ksmD9034or-1204}, we obtain
the first identity in~\eqref{ksmD9cvbHA034or-1204-11}, as claimed.

In addition, since
$$ \frac{\partial v_m}{\partial x_j}(x,t)=2\pi i\sum_{k\in\Z^3} k_j\,\widehat v_{m,k}(t)\,e^{2\pi ik\cdot x},$$
recalling~\eqref{WTTGfdgtchofrTLA:0-2.9rjHNX} and Parseval's Identity~\eqref{PARS} we see that
\begin{equation*}
{\mathcal{D}}(t)=4\pi^2\sum_{{k\in\Z^3}\atop{1\le j,m\le3}}k_j^2\,|\widehat v_{m,k}(t)|^2,
\end{equation*}
leading to the second identity in~\eqref{ksmD9cvbHA034or-1204-11}.

\paragraph{Solution to Exercise~\ref{2NSSEEV.bis}.} We formally show that
\begin{equation}\label{RMANMIASOW98}
{\mathcal{E}}'(t)=-\nu{\mathcal{V}}(t)+
\int_{(0,1)^3}\left( p(x,t)+\frac{|v(x,t)|^2}2\right)
\nabla\cdot v(x,t)\,dx.
\end{equation}

To check this, we proceed as follows.
In view of~\eqref{WTTGfdgtchofrTLA:0-1}, some integrations by parts formally give that
\begin{equation}\label{MCRDCJCSSTR120}\begin{split}
{\mathcal{E}}'(t)&=\int_{(0,1)^3} v(x,t)\cdot \frac{\partial v}{\partial t}(x,t)\,dx\\&=
\int_{(0,1)^3} v(x,t)\cdot \left( \nu \,\Delta v (x,t)-\nabla p(x,t)-(v(x,t) \cdot \nabla )v(x,t) \right)\,dx\\&=
-\nu\int_{(0,1)^3}|\nabla v(x,t)|^2\,dx +
\int_{(0,1)^3} p(x,t)\nabla\cdot v(x,t)\,dx\\&\qquad\qquad\qquad\quad-\int_{(0,1)^3}
v(x,t)\cdot\Big((v(x,t) \cdot \nabla )v(x,t)\Big)\,dx.\end{split}
\end{equation}

Actually, one more integration by parts shows that
\begin{eqnarray*}
&&\int_{(0,1)^3}v(x,t)\cdot\Big((v(x,t) \cdot \nabla )v(x,t)\Big) \,dx\\&&\qquad=
-\int_{(0,1)^3}\left( | v(x,t)|^2\,\nabla\cdot v(x,t)+v(x,t)\cdot\Big(
(v(x,t) \cdot \nabla )v(x,t) \Big)\right)\,dx,
\end{eqnarray*}
that is
\begin{eqnarray*}
&&\int_{(0,1)^3}v(x,t)\cdot\Big((v(x,t) \cdot \nabla )v(x,t) \Big)\,dx=
-\frac12\int_{(0,1)^3} |v(x,t)|^2\,\nabla\cdot v(x,t)\,dx.
\end{eqnarray*}
We plug this identity into~\eqref{MCRDCJCSSTR120}, recall~\eqref{WTTGfdgtchofrTLA:0-2} and~\eqref{ksmD9cvbHA034or-1204-11},
and obtain~\eqref{RMANMIASOW98}, as desired.

Now, on the one hand, recalling the solenoidal velocity condition in~\eqref{NS6THTgBplaISw}, one would
formally deduce from~\eqref{RMANMIASOW98} that
\begin{equation}\label{RMANMIASOW981}
{\mathcal{E}}'(t)=-\nu{\mathcal{V}}(t),
\end{equation}
which would elegantly suggest that the energy dissipation is proportional to viscosity and enstrophy
(and this might be close to what we might have expected, since usually energy dissipation is a byproduct of friction,
and, in light of~\eqref{ksmD9cvbHA034or-1204-11}, this viscous dissipation is expected to
consume the energy stronger in the presence of high ranges of frequencies).

On the other hand, to make sense of the additional term in~\eqref{RMANMIASOW98} one needs to require that the velocity field~$v$ is sufficiently regular (roughly speaking, one wants to integrate the square of the velocity against the derivative of the velocity itself, so, dimensionally, it would require some control on the ratio between~$v$ and a third of the space variable).
The energy balance of the Navier-Stokes and Euler equations are indeed a very delicate matter
(related to the so-called \index{Onsager's conjecture} Onsager's conjecture, now an established result,
see~\cite[Theorem~6.24]{MR4475666}).

We stress that this is not just a mathematical artifact or an unnecessary complication of an abstract theory:
the fact that~\eqref{RMANMIASOW981} provides a too naive expectation
is highlighted by physical and numerical observations, and
several experiments confirmed this anomalous dissipation of energy, approaching constant, nontrivial
values in the small viscosity regime, especially in the presence of turbulence, see e.g.~\cite{MR5851, MR36116, MR1428905, MR1872661, MR3866888, MR3929468, MR4475666} and the references therein.

\section{Solutions to selected exercises of Section~\ref{POISSONKERN-ex4}}

\paragraph{Solution to Exercise~\ref{SPKEMLEMFCf01}.} We prove the second claim (the first one being similar). Arguing for a contradiction, we suppose that
\begin{equation}\label{CABNScmcmcs.1D} M:=\max_{B_1} u>0,\end{equation}
where we used the short notation~$u:=u_f$.

By Sard's Lemma (see e.g.~\cite[Theorem~1.3 on page~69]{zbMATH03555096} and the references therein), for almost all~$m\in(0,M)$ the level set~$\{u=m\}$ is noncritical, namely
$$ \big\{{\mbox{$x\in B_1$ s.t. $u(x)=m$ and~$\nabla u(x)=0$}}\big\}=\varnothing$$
and therefore~$\{u=m\}$ is a curve of class at least~$C^1$.
For this reason, we can pick one of these noncritical values~$m\in(0,M)$, use the Divergence Theorem and find that
\begin{eqnarray*}
0=\int_{ \{u>m\} }\Delta u(x)\,\big(u(x)-m\big)\,dx=-\int_{ \{u>m\} }|\nabla u(x)|^2\,dx.
\end{eqnarray*} 
Therefore, we conclude that~$\nabla u(x)=0$ for all~$x\in\{u>m\}$.
On this account, $u$ is necessarily constant, and constantly equal to~$m$, in~$\{u>m\}$, in contradiction with~\eqref{CABNScmcmcs.1D}.

\paragraph{Solution to Exercise~\ref{U3fNIDIRIY.1}.} Suppose that there are two solutions, say~$u$ and~$v$, and consider their difference~$w:=u-v$. Then, $w\in C^2(B_1)\cap C(\overline{B_1})$ and we have that
$$\begin{dcases}\Delta w=0&{\mbox{ in }}B_1,\\w=0&{\mbox{ on }}\partial B_1.\end{dcases}$$
By Exercise~\ref{SPKEMLEMFCf01}, we conclude that~$w$ vanishes identically.

\paragraph{Solution to Exercise~\ref{U3fNIDIRIY.1bisa}.} We observe that~$B_r(p_0)\Subset B_1$.
We also provide a general remark about real and complex integrals. Namely, writing~$p_0=(x_0,y_0)\in\R^2$,
given~$v:\R^2\to\R$, we can parameterise the circle~$\partial B_r(p_0)$ by
$$\gamma(\phi):=\big(x_0+r\cos(2\pi\phi),y_0+r\sin(2\pi\phi)\big)$$ with~$\phi\in[0,1)$,
and obtain that
\begin{equation}\label{U3fNIDIRIY.1bisabgt67}
\int_{\partial B_r(p_0)} v=\int_0^1 v(\gamma(\phi))\,|\dot\gamma(\phi)|\,d\phi=
2\pi r\int_0^1 v\big(x_0+r\cos(2\pi\phi),y_0+r\sin(2\pi\phi)\big)\,d\phi.
\end{equation}

Moreover, if~$V:\C\to\R$ is such that~$v(x,y)=V(x+iy)$ for all~$(x,y)\in\R^2$, we can parameterise the circle~${\mathcal{C}}_r(z_0)$ of radius~$r$ centred at~$z_0:=x_0+iy_0$ and travelled anticlockwise by~$z(\phi)=z_0+re^{2\pi i\phi}$, with~$\phi\in[0,1)$, finding that
\begin{eqnarray*}
\frac1{2\pi i}\oint_{{\mathcal{C}}_r(z_0)}\frac{V(z)}{z-z_0}\,dz=\int_0^1 V\big(z_0+re^{2\pi i\phi}\big)\,d\phi=
\int_0^1 v\big(x_0+r\cos(2\pi\phi),y_0+r\sin(2\pi\phi)\big)\,d\phi.
\end{eqnarray*}

Comparing this with~\eqref{U3fNIDIRIY.1bisabgt67}, we have that
\begin{equation}\label{U3fNIDIRIY.1bisabgt6765.012}
\frac1{2\pi r}\int_{\partial B_r(p_0)} v=
\frac1{2\pi i}\oint_{{\mathcal{C}}_r(z_0)}\frac{V(z)}{z-z_0}\,dz.
\end{equation}

Now we define
$$F(z):=\widehat f_0+2\sum_{k=1}^{+\infty} \widehat f_k\,z^k$$
and we stress that~$F$ is holomorphic in the unit disc (recall~\eqref{oqpkdwgmintthcrz0frazz0dzarz}).

We also use~\eqref{MjAPLSvbiyKamdfRbr} to see that
$$ u_f=\Re F,$$
and similarly we consider~$v_f:=\Re F$.

Moreover, by Cauchy Integral Formula in complex analysis (see e.g.~\cite[Theorem~4.1 on page~45]{MR1976398}),
$$ \frac1{2\pi i}\oint_{{\mathcal{C}}_r(z_0)}\frac{F(z)}{z-z_0}\,dz= F(z_0).$$
We can thus employ~\eqref{U3fNIDIRIY.1bisabgt6765.012} and conclude that
\begin{equation*}\begin{split}&
u_f(p_0)=\Re F(z_0)=\Re\left( \frac1{2\pi i}\oint_{{\mathcal{C}}_r(z_0)}\frac{F(z)}{z-z_0}\,dz\right)\\&\quad=
\Re\left( \frac1{2\pi i}\oint_{{\mathcal{C}}_r(z_0)}\frac{\Re F(z)}{z-z_0}\,dz+i\frac1{2\pi i}\oint_{{\mathcal{C}}_r(z_0)}\frac{\Im F(z)}{z-z_0}\,dz\right)\\&\quad=
\Re\left( \frac1{2\pi r}\int_{\partial B_r(p_0)}u_f+i\frac1{2\pi r}\int_{\partial B_r(p_0)}v_f\right)= \frac1{2\pi r}\int_{\partial B_r(p_0)}u_f.
\end{split}
\end{equation*}
The proof of~\eqref{U3fNIDIRIY.1bisa-EQ2} is thereby complete.

Then, to prove~\eqref{U3fNIDIRIY.1bisa-EQ1} we utilise~\eqref{U3fNIDIRIY.1bisa-EQ2}, \eqref{U3fNIDIRIY.1bisabgt67}, and polar coordinate centred at~$p_0$ to find that
\begin{eqnarray*}&&
\int_{B_r(p_0)}u_f
=\int_0^r\left(
2\pi \rho\int_0^1 u_f\big(x_0+r\cos(2\pi\phi),y_0+r\sin(2\pi\phi)\big)\,d\phi\right)\,d\rho\\&&\qquad
=\int_0^r\left(
\int_{\partial B_\rho(p_0)} u_f\right)\,d\rho=
\int_0^r\left(2\pi \rho\,u_f(p_0)\right)\,d\rho=
\pi r^2 u_f(p_0),
\end{eqnarray*}
as desired.

\section{Solutions to selected exercises of Section~\ref{TEojemy04i7u:Asdweg54uk:028tuyjh}}

\paragraph{Solution to Exercise~\ref{jlwmeotrjho590iukjiqdhfpirnHSNdlm0ujt1}.}
It is easily seen that the series in~\eqref{u:funcorners}
converges uniformly in every compact subset of~$(0,1)\times(0,1)$, together with its derivatives.

As a result, one can perform differentiations under the summation sign
to see that~$\Delta u=0$. Also, since~$\sin((2k+1)\pi x)=0$ when~$x\in\N$, we have that~$u(0,y)=u(1,y)=0$. Along the same lines, since~$\sinh((2k+1)\pi y)=0$ when~$y=0$, it follows that~$u(x,0)=0$.

It remains to show that~$u$ attains the boundary datum along~$(0,1)\times\{1\}$. In order to achieve this,
we recall that, by Exercise~\ref{ojld03-12d}, the square wave
$$f(x)=\begin{dcases}1&{\mbox{ if }}x\in\left[0,\frac12\right),\\
-1&{\mbox{ if }}x\in\left[\frac12,1\right),\end{dcases}$$
can be written in the form
$$f(x)=\sum_{j=0}^{+\infty} \frac4{\pi(2j+1)}\sin(2\pi(2j+1)x),$$
where this pointwise identity is valid for all~$x\in(0,1)\setminus\left\{\frac12\right\}$,
due to Theorem~\ref{C1uni}.

Hence, for all~$x\in(0,1)$, $$1=f\left(\frac{x}2\right)=
\sum_{j=0}^{+\infty} \frac4{\pi(2j+1)}\sin(\pi(2j+1)x)=u(x,1),$$
as desired.

\paragraph{Solution to Exercise~\ref{jlwmeotrjho590iukjiqdhfpirnHSNdlm0ujt1BIS}.}
Since
\begin{equation}\label{joqfw0t7u0439867-iukj}
\frac{\sinh((2k+1)\pi y)}{\sinh((2k+1)\pi)}-e^{(2k+1)\pi(y-1)}
=\frac{e^{(2k+1)\pi(y-1)} -e^{(2k+1)\pi(1-y)}}{
(e^{(2k+1)\pi}-e^{-(2k+1)\pi})\,e^{(2k+1)\pi}},\end{equation}
for all~$y\in[0,1]$ we have that
$$\left|\frac{\sinh((2k+1)\pi y)}{\sinh((2k+1)\pi)}
-e^{(2k+1)\pi(y-1)}\right|\le
\frac{2e^{(2k+1)\pi} }{
(e^{(2k+1)\pi}-e^{-(2k+1)\pi})\,e^{(2k+1)\pi}}=
\frac{2}{e^{(2k+1)\pi}-e^{-(2k+1)\pi}}.$$
As a consequence, setting
$$ v(x,y):=\frac4\pi\sum_{k=0}^{+\infty}
\frac1{2k+1}\left(\frac{\sinh((2k+1)\pi y)}{\sinh((2k+1)\pi)}-e^{(2k+1)\pi(y-1)}\right)
\,\sin((2k+1)\pi x),$$we have that~\eqref{u:funcorners:W} holds true
and the above series converges uniformly in~$[0,1]\times[0,1]$,
thus entailing the desired continuity of~$v$.

It also follows from~\eqref{joqfw0t7u0439867-iukj} that the derivatives
of the approximating sum of~$v$ converge uniformly in~$[0,1]\times[0,1]$,
ensuring that~$v$ is differentiable
infinitely many times and its derivatives are continuous in~$[0,1]\times[0,1]$.

\paragraph{Solution to Exercise~\ref{jlwmeotrjho590iukjiqdhfpirnHSNdlm0ujt1BIS.04}.}
We have that
\begin{eqnarray*} \frac{1-Z^t}{1-Z}&=&\frac{1-r^t e^{\pi it\phi}}{1-r e^{\pi i\phi}}\\&=&
\frac{1-r^t \cos(\pi t\phi)-ir^t \sin(\pi t\phi)}{1-r\cos(\pi \phi)-ir\sin(\pi \phi)}\\
&=&\frac{\big(1-r^t \cos(\pi t\phi)-ir^t \sin(\pi t\phi)\big)\big(1-r\cos(\pi \phi)+ir\sin(\pi \phi)\big)}{\big(1-r\cos(\pi \phi)-ir\sin(\pi \phi)\big)
\big(1-r\cos(\pi \phi)+ir\sin(\pi \phi)\big)}\\&=&\frac{\big(1-r^t \cos(\pi t\phi)-ir^t \sin(\pi t\phi)\big)\big(1-r\cos(\pi \phi)+ir\sin(\pi \phi)\big)}{\big(1-r\cos(\pi \phi)\big)^2 +r^2\sin^2(\pi \phi)}.
\end{eqnarray*}

Notice also that
$$r^2\sin^2(\pi \phi)\ge\delta^2\sin^2(\pi\delta).$$
As a consequence,
\begin{eqnarray*} &&\left|\frac{1-Z^t}{1-Z}\right|\le
\frac{9}{\big(1-r\cos(\pi \phi)\big)^2 +r^2\sin^2(\pi \phi)}\le\frac{9}{r^2\sin^2(\pi \phi)}\le
\frac9{\delta^2\sin^2(\pi\delta)}.
\end{eqnarray*}

\paragraph{Solution to Exercise~\ref{jlwmeotrjho590iukjiqdhfpirnHSNdlm0ujt2}.}
In~$[0,1]\times[0,1-\rho]$, the advertised uniform convergence follows
from the fact that
$$\left|
\frac{\sinh((2k+1)\pi y)}{\sinh((2k+1)\pi)}
\right|\le\frac{e^{(2k+1)\pi y}+1}{e^{(2k+1)\pi}-1}\le
\frac{e^{(2k+1)\pi (1-\rho)}+1}{e^{(2k+1)\pi}-1},$$
which, for large~$k$, is bounded from above by~$4e^{-(2k+1)\pi\rho}$.

Hence, we can focus on the uniform convergence in~$[\rho,1-\rho]\times[1-\rho,1]$.
Equivalently, we check the convergence properties of the series in~\eqref{u:funcorners:Z}.
In the above domain, by Exercise~\ref{jlwmeotrjho590iukjiqdhfpirnHSNdlm0ujt1BIS.04} we have that, for all~$t\in(0,+\infty)$,
$$ \left|\frac{1-Z^t}{Z(1-Z)}\right|\le C,$$
where~$C>0$ depends only on~$\rho$.

We also employ the summation by parts formula in Exercise~\ref{SBPF}, used here with
$$ \alpha_k:=\frac{1}{2k+1}\qquad{\mbox{and}}\qquad\beta_k:=\sum_{j=0}^k Z^{2j-1}.$$
Thus, we find that
$$ |\beta_k|=\left|
\frac1Z\sum_{j=0}^k Z^{2j}
\right|=\left|
\frac{1-Z^{2(k+1)}}{Z(1-Z)}
\right|\le C$$
and that, for all~$n\ge m$,
\begin{eqnarray*}\sum_{k=m}^n \frac{Z^{2k+1}}{2k+1}&=&
\sum_{k=m}^{n}\alpha_{k}(\beta_{k+1}-\beta_{k})\\&=&\left(\alpha_{n}\beta_{n+1}-\alpha_{m}\beta_{m}\right)-\sum_{k=m+1}^{n}\beta_{k}(\alpha_{k}-\alpha_{k-1})\\&=&\left(\alpha_{n}\beta_{n+1}-\alpha_{m}\beta_{m}\right)+\sum_{k=m+1}^{n}\frac{2\beta_{k}}{4k^2-1}.\end{eqnarray*}

Consequently, up to freely renaming~$C$, when~$(x,y)\in[\rho,1-\rho]\times[1-\rho,1]$ we have that
\begin{eqnarray*}\left|\sum_{k=m}^n \frac{Z^{2k+1}}{2k+1}\right|\le
\frac{C}n+\frac{C}m+\sum_{k=m+1}^{n}\frac{C}{4k^2-1}
\le\frac{C}m,\end{eqnarray*}
yielding the desired uniform convergence.

\paragraph{Solution to Exercise~\ref{jlwmeotrjho590iukjiqdhfpirnHSNdlm0ujt2hgd}.}
By~\eqref{SDMLSPMCAUSDM:D}, we have that
$$ 1+Z=2 - \pi r(\cos\theta -i \sin\theta)+O(r^2)$$
and therefore
\begin{equation}\label{REWSDV:qwdfgvbweEQscFsdRWEGRBN-2}\begin{split}
\alpha_+&=\arctan\left(\frac{\pi r \sin\theta+O(r^2)
}{2 - \pi r\cos\theta+O(r^2)}\right)\\&=
\frac{\pi r \sin\theta+O(r^2)
}{2 - \pi r\cos\theta+O(r^2)}+O(r^2)\\&=
\frac{\pi r \sin\theta}2+O(r^2).
\end{split}\end{equation}

The computation of~$\alpha_-$ is a bit more involved,
since we would like to obtain an estimate on the remainders
in terms of powers of~$r$ which are uniform
with respect to the polar angle~$\theta\in\left(0,\frac\pi2\right)$, and this is somewhat delicate since~$s:=\sin\theta$
and~$c:=\cos\theta$ are positive quantities for a given~$\theta\in\left(0,\frac\pi2\right)$, but they vanish at one endpoint
of this interval.

To overcome this difficulty,
we deduce from~\eqref{SDMLSPMCAUSDM:D} that
\begin{equation*}
1-Z= 1-e^{-\pi r(c- is)}=1-e^{-\pi cr}\cos(\pi sr)
-ie^{-\pi cr}\sin(\pi sr),
\end{equation*}
whence
\begin{equation}\label{aldjwSDFGHNergpoehgnbitor:024rotikght}
\begin{split}
\alpha_-&=-\arccot\left(\frac{1-e^{-\pi cr}\cos(\pi sr)}{e^{-\pi cr}\sin(\pi sr)}\right)\\
&=-\arccot\left(\frac{e^{\pi cr}-\cos(\pi sr)}{\sin(\pi sr)}\right)\\
&=-\arccot\left(\frac{
1+\pi cr+\frac{\pi^2c^2r^2}2+O(c^3r^3)-\left( 1-\frac{ \pi^2s^2 r^2}2 +O(s^4 r^4)\right)}{\pi sr+O(s^3r^3)}\right)\\&=-\arccot\left(\frac{
\pi cr+\frac{\pi^2r^2}2+O(r^3)
}{\pi sr+O(s^3r^3)}\right)\\&=-\arccot\left(\frac{
c+\frac{\pi r}{2}+O( r^2)
}{s+O(s^3r^2)}\right).\end{split}
\end{equation}

Furthermore,
\begin{equation}\label{qwd0fjovD34ERfFwm12}\begin{split}
\Upsilon&:=
\frac{
c+\frac{\pi r}{2}+O( r^2)
}{s+O(s^3r^2)}-\cot\theta\\&=
\frac{
c+\frac{\pi r}{2}+O(r^2)
}{s+O(s^3r^2)}-\frac{c}s\\
&=
\frac{\pi r}{2 s}\big(1+O(r)\big).\end{split}
\end{equation}

We claim that
\begin{equation}\label{qowdjfvn0238e3josjndDD}
\alpha_-=-\theta+\frac{\pi rs}2+O(r^2)
\end{equation}
and we prove this by distinguishing two cases, namely
when~$|s|\ge\frac{\sqrt2}2$
and when~$|s|<\frac{\sqrt2}2$.

If~$|s|\ge\frac{\sqrt2}2$,
in view of~\eqref{qwd0fjovD34ERfFwm12} we have that~$\Upsilon=O(r)$.
Thus, we perform a Taylor expansion and we see that
\begin{eqnarray*}
\arccot(\cot\theta+\Upsilon)&=&
\arccot(\cot\theta)-\frac{\Upsilon}{1+\cot^2\theta}+O(\Upsilon^2)\\
&=&\theta-\frac{\frac{\pi r}{2 s}\big(1+O(r)\big)}{1+\frac{c^2}{s^2}}+O(r^2)\\&=&\theta-
\frac{\pi rs}{2 }+O(r^2).
\end{eqnarray*}
This and~\eqref{aldjwSDFGHNergpoehgnbitor:024rotikght} yield~\eqref{qowdjfvn0238e3josjndDD} in this case.

Suppose now that~$|s|<\frac{\sqrt2}2$. Then,
we have that~$|c|\ge\frac{\sqrt2}2$ and consequently
\begin{eqnarray*}
\arccot(\cot\theta+\Upsilon)&=&
\arctan\left( \frac1{\cot\theta+\Upsilon}\right)\\&=&
\arctan\left( \frac1{\frac{c}s+
\frac{\pi r}{2 s}\big(1+O(r)\big)}\right)\\&=&
\arctan\left( \frac{s}{c+
\frac{\pi r}{2}\big(1+O(r)\big)}\right)\\&=&\arctan\left( \frac{s}{c}\left(1-\frac{\pi r}{2c}+O(r^2)\right)\right)\\&=&
\arctan\left( \frac{s}{c}-\frac{\pi rs}{2c^2}+O(sr^2)\right)\\&=&\theta-\frac{\frac{\pi rs}{2c^2}}{1+\frac{s^2}{c^2}}+O(r^2)\\&=&\theta-
\frac{\pi rs}{2 }+O(r^2).
\end{eqnarray*}
Combining this information with~\eqref{aldjwSDFGHNergpoehgnbitor:024rotikght},
we complete the proof of~\eqref{qowdjfvn0238e3josjndDD}.

\section{Solutions to selected exercises of Section~\ref{Pkjqn90-902}}

\paragraph{Solution to Exercise~\ref{TELES}.}
For all~$k\ge2$,
$$ \frac{1}{k^2}\le \frac{1}{k(k-1)}=\frac1{k-1}-\frac{1}{k}.$$
As a result, for all~$N\in\N$, $N\ge2$,
$$ \sum _{k=2}^{N}\frac{1}{k^2}\le\sum _{k=2}^{N}\left(\frac1{k-1}-\frac{1}{k}\right)=
1-\frac1N.$$
Therefore,
$$ \sum _{k=1}^{+\infty}\frac{1}{k^2}=1+\lim_{N\to+\infty}
\sum _{k=2}^{N}\frac{1}{k^2}\le 1+\lim_{N\to+\infty}\left(
1-\frac1N\right)=2.$$

\paragraph{Solution to Exercise~\ref{TELES-INTR}.}
We have that \begin{eqnarray*}&&\sum _{k=1}^{+\infty }\frac{1}{k^2}=1+\sum _{k=2}^{+\infty }
\int_{k-1}^{k}\frac{dx}{k^2}\le1+\sum _{k=2}^{+\infty }
\int_{k-1}^{k}\frac{dx}{x^2}=1+\int_{1}^{+\infty}\frac{dx}{x^2}=1+1=2.\end{eqnarray*}

\paragraph{Solution to Exercise~\ref{AOJsnkwe9rgeotg2-34rt-p0}.}
We consider the sawtooth waveform, i.e. the periodic extension~$f_{\text{per}}$ of period~$1$, as defined in~\eqref{FPER}, of the function
$$ [0,1)\ni x\longmapsto f(x):=x-\frac12.$$
We know (see Exercise~\ref{SA:W}) that the~$k$th Fourier coefficient of~$f_{\text{per}}$ equals~$\frac{i}{2\pi k }$ when~$k\ne0$ and~$0$ when~$k=0$.

Moreover, $f_{\text{per}}\in L^2((0,1))$ and therefore (see Theorem~\ref{THCOL2FB}) Parseval's Identity holds true.

As a result,
$$ \frac12\sum_{k=1}^{+\infty}\frac{1}{\pi^2 k^2}=
\sum_{k\in\Z} |\widehat f_k|^2 = \|f\|^2_{L^2((0,1))}=\int_0^1\left(x-\frac12\right)^2\,dx=\frac1{12},$$
from which we obtain~\eqref{BAFO}.

\paragraph{Solution to Exercise~\ref{AOJsnkwe9rgeotg2-34rt}.} The gist is to develop both sides of~\eqref{EUJNC-P} in power series near the origin. To this end, we take~$x\in\left(-\frac12,\frac12\right)$ and we observe that, for all~$N\in\N$,
\begin{eqnarray*}
\ln\left(\prod_{j=1}^{N}\left(1-\frac{x^2}{ j^2}\right)\right)=
\sum_{j=1}^{N}\ln\left(1-\frac{x^2}{ j^2}\right)=
-\sum_{j=1}^{N}\sum_{\ell=1}^{+\infty}\frac{x^{2\ell}}{ j^{2\ell}\,\ell}.
\end{eqnarray*}
By the absolute convergence of the series, we can swap the order of summation and find that
\begin{eqnarray*}&&
\ln\left(\prod_{j=1}^{N}\left(1-\frac{x^2}{ j^2}\right)\right)=
-\sum_{\ell=1}^{+\infty}\sum_{j=1}^{N}\frac{x^{2\ell}}{ j^{2\ell}\,\ell}\\&&\qquad=-\sum_{\ell=1}^{+\infty}\sum_{j=1}^{+\infty}\frac{x^{2\ell}}{ j^{2\ell}\,\ell}+\epsilon_N=-\sum_{\ell=1}^{+\infty}\frac{\zeta(2\ell)\,x^{2\ell}}{\ell}+\epsilon_N,
\end{eqnarray*}
with~$\epsilon_N$ infinitesimal as~$N\to+\infty$.

From this and Euler's sine product formula~\eqref{EUJNC-P}, we gather that
\begin{equation}\label{L-2pkjfr320-7t439hgonb89}\begin{split}&
\ln\left(\frac{\sin (\pi x)}{\pi x}\right)
=\ln\left(\prod_{j=1}^{+\infty}\left(1-\frac{x^2}{ j^2}\right)\right)=
\lim_{N\to+\infty}\ln\left(\prod_{j=1}^{N}\left(1-\frac{x^2}{ j^2}\right)\right)\\&\qquad=-
\lim_{N\to+\infty}\left(\sum_{\ell=1}^{+\infty}\frac{\zeta(2\ell)\,x^{2\ell}}{\ell}+\epsilon_N\right)=-
\sum_{\ell=1}^{+\infty}\frac{\zeta(2\ell)\,x^{2\ell}}{\ell}.\end{split}
\end{equation}

Also, since, near the origin,
$$ \frac{\sin (\pi x)}{\pi x}=\frac{1}{\pi x}\sum_{\ell=0}^{+\infty}\frac{(-1)^\ell\,(\pi x)^{2\ell+1}}{(2\ell+1)!}
=\sum_{\ell=0}^{+\infty}\frac{(-1)^\ell\,(\pi x)^{2\ell}}{(2\ell+1)!}=
1-\frac{(\pi x)^{2}}{6}+o(x^2),
$$
we have that
$$ \ln\left(\frac{\sin (\pi x)}{\pi x}\right)=\ln\left(1-\frac{(\pi x)^{2}}{6}+o(x^2)\right)=
-\frac{(\pi x)^{2}}{6}+o(x^2).$$

We combine this and~\eqref{L-2pkjfr320-7t439hgonb89} and we conclude that
$$ -\frac{(\pi x)^{2}}{6}+o(x^2)=- \zeta(2)\,x^{2} +o(x^2),$$
leading to~\eqref{BAFO}.

\paragraph{Solution to Exercise~\ref{AOJsnkwe9rgeotg2-34rt-DOU}.}
This method was invented in~\cite{MR737691}. On the one hand, we consider the change of variable
$$ (X,Y):=\left(\frac{x+y}2,\frac{x-y}2\right)$$
giving that
\begin{equation*}
I:=\int_0^1\int_0^1 \frac{dx\,dy}{1-xy}=2\iint_S\frac{dX\,dY}{1-X^2-Y^2},
\end{equation*}
where~$S$ is the square with vertexes in~$(0,0)$, $\left(\frac12,\frac12\right)$, $\left(\frac12,-\frac12\right)$ and~$(1,0)$.

Hence,
\begin{equation*}
I=4\int_{0}^{1/2}\left(\int_0^X\frac{dY}{1-X^2-Y^2}\right)\,dX+4\int_{1/2}^{1}\left(\int_0^{1-X}\frac{dY}{1-X^2-Y^2}\right)\,dX.
\end{equation*}
Using that, for all~$a>0$ and~$\xi\in\R$,
$$ \int_0^\xi \frac{dt}{a^2+t^2}=\frac{\arctan (\xi/a)}a,$$
and substituting with trigonometric functions,
we obtain that
\begin{equation}\label{PFQDF87}I=\frac{\pi^2}{18}+\frac{\pi^2}9=
\frac{\pi^2}6.\end{equation}

On the other hand, using that, for all~$r\in(-1,1)$,
$$ \frac1{1-r}=\sum_{k=0}^{+\infty} r^k,$$ taking~$r:=xy$
by the Monotone Convergence Theorem we see that
\begin{eqnarray*}
I=\int_0^1\int_0^1 \sum_{k=0}^{+\infty} (xy)^k\,dx\,dy=\sum_{k=0}^{+\infty} 
\int_0^1\int_0^1(xy)^k\,dx\,dy=\sum_{k=0}^{+\infty} \frac1{(k+1)^2}=\sum_{k=1}^{+\infty}\frac1{k^2}. 
\end{eqnarray*} Comparing this with~\eqref{PFQDF87} we obtain the desired result.

\paragraph{Solution to Exercise~\ref{SOMMS1}.}
We employ the summation by parts formula in Exercise~\ref{SBPF}, with
$$ \alpha_k:=\sum_{j=1}^k\frac1j=H_k\qquad{\mbox{and}}\qquad\beta_k:=-\frac1k.$$
Thus,
$$\beta_{k+1}-\beta_k=\frac1k-\frac1{k+1}=\frac1{k(k+1)}$$
and we thereby infer from~\eqref{BYPA} that \begin{eqnarray*}
\sum_{k=1}^{n}\frac{H_k}{k(k+1)}=\left(-\frac{H_n}{n+1}+1\right)+\sum_{k=2}^{n}\frac1{k^2}=
-\frac{H_n}{n+1}+\sum_{k=1}^{n}\frac1{k^2}
.\end{eqnarray*}
Since (see e.g. Exercise~\ref{CESA})
$$ \lim_{n\to+\infty}\frac{H_n}{n+1}=0,$$
we thus find that
$$ \sum_{k=1}^{+\infty}\frac{H_k}{k(k+1)}=\sum_{k=1}^{+\infty}\frac1{k^2}.$$
The desired result now follows from Proposition~\ref{BASEL}.

\paragraph{Solution to Exercise~\ref{SOMMS2}.}
We have that
$$ 1+2+\dots+n=\frac{n(n+1)}2,$$
hence the desired result follows from Exercise~\ref{SOMMS1}.

\paragraph{Solution to Exercise~\ref{EUFO}.}
Let us consider the first~$N$ primes~$1<p_1<p_2<\dots<p_N$. By the Prime Factorization Theorem, every integer greater than~$1$ can be represented uniquely (up to the order of the factors) as a product of prime numbers (and these prime numbers are less than or equal to the given integer, otherwise their product would overcome it).

As a result, for every
$n\in\N\cap [1,p_N]$ there exist
$$j_1(n),\dots,j_{M(n)}(n)\in\{1,\dots,N\}\qquad{\mbox{ and }}\qquad
\ell_1(n),\dots,\ell_{M(n)}(n)\in\N$$ such that
$$ n=p_{j_1(n)}^{\ell_1(n)}\dots p_{\ell_{M(n)}(n)}^{\ell_{M(n)}(n)}.$$
In particular, if~$n\in\N$ and~$n$ cannot be written as~$p_1^{\ell_1}\dots p_N^{\ell_N}$ then necessarily~$n\ge p_N+1$.

That is, if
$$ {\mathcal{F}}:=\big\{
{\mbox{$n\in\N$ s.t. $n$ cannot be written as $p_1^{\ell_1}\dots p_N^{\ell_N}$}}\big\},$$
we have that \begin{equation}\label{019o2iwje0239o8ytygvSoqmL} {\mathcal{F}}\subseteq [p_N+1,+\infty).\end{equation}

Actually, the observation in~\eqref{019o2iwje0239o8ytygvSoqmL} is useful if we know already that there are infinitely many primes (which is certainly the case, see Exercise~\ref{EUCL}), but we do not really need this information to make this argument work: simply, if there were only finitely many primes, we take~$N$ so large that~$\{p_1,\dots,p_N\}$ collects all the primes and then in this case we can replace~\eqref{019o2iwje0239o8ytygvSoqmL} by
\begin{equation}\label{019o2iwje0239o8ytygvSoqmL-al} {\mathcal{F}}=\varnothing.\end{equation}

Moreover, if~$x\in(0,+\infty)$, setting~$y:=\ln x$ we see that, for all~$s\in\C$,
\begin{eqnarray*}
x^s=x^{\Re s+i\Im s}=x^{\Re s} \,e^{i y\Im s}
\end{eqnarray*}
and therefore
\begin{equation}\label{019o2iwje0239o8ytygvSoqmL2}
|x^s|=x^{\Re s}.\end{equation}

Consequently, by~\eqref{019o2iwje0239o8ytygvSoqmL}, \eqref{019o2iwje0239o8ytygvSoqmL-al} and~\eqref{019o2iwje0239o8ytygvSoqmL2}, if~$\Re s>1$,
\begin{equation}\label{0iuhgvf8yf85868g5fP-2}\begin{split}&
\left|\sum_{n=1}^{+\infty}\frac1{n^s}-
\sum_{\ell_1,\dots,\ell_N=0}^{+\infty}\frac1{(p_1^{\ell_1}\dots p_N^{\ell_N})^s}\right|=\left|\sum_{n\in{\mathcal{F}}}\frac1{n^s}\right|\le
\sum_{n\in{\mathcal{F}}}\left|\frac1{n^s}\right|\\&\qquad\le\begin{dcases}\displaystyle
\sum_{n=p_N+1}^{+\infty}\frac1{n^{\Re s}}&{\mbox{ if there exist infinitely many primes,}}\\ 0&{\mbox{ if there exist only finitely many primes.}}\end{dcases}\end{split}
\end{equation}

Now we recall the geometric series for~$z\in\C$ with~$|z|<1$, namely
\begin{equation}\label{019o2iwje0239o8ytygvSoqmL3} \sum_{\ell=0}^{+\infty} z^\ell=\frac1{1-z}.\end{equation}
Since, by~\eqref{019o2iwje0239o8ytygvSoqmL2}, for all~$j\in\{1,\dots,N\}$,
$$ |p_j^{-s}|=p_j^{-\Re s}<1,$$
we can apply~\eqref{019o2iwje0239o8ytygvSoqmL3} with~$z:=p_j^{-s}$ and we find that
$$ \sum_{\ell_j=0}^{+\infty} p_j^{-\ell_j s}=\frac1{1-p_j^{-s}}.$$
Therefore,
$$ \prod_{j=1}^N \frac1{1-p_j^{-s}}=
\sum_{\ell_1,\dots,\ell_N=0}^{+\infty} p_1^{-\ell_1 s}\dots p_N^{-\ell_N s}=\sum_{\ell_1,\dots,\ell_N=0}^{+\infty}\frac1{(p_1^{\ell_1}\dots p_N^{\ell_N})^s}.$$
Plugging this information into~\eqref{0iuhgvf8yf85868g5fP-2}, we infer that
\begin{equation}\label{09879iujh7ygf5tgb234rfsdfghj284P1} \begin{split}&\left|\sum_{n=1}^{+\infty}\frac1{n^s}-\prod_{j=1}^N \frac1{1-p_j^{-s}}\right|\\&\qquad
\le\begin{dcases}\displaystyle
\sum_{n=p_N+1}^{+\infty}\frac1{n^{\Re s}}&{\mbox{ if there exist infinitely many primes,}}\\ 0&{\mbox{ if there exist only finitely many primes.}}\end{dcases}\end{split}\end{equation}
Now, if there exist only finitely many primes, the observation in~\eqref{09879iujh7ygf5tgb234rfsdfghj284P1} proves that
\begin{equation}\label{02i0rjf94yjcxjcjn2nrgt1ds}
\sum_{n=1}^{+\infty}\frac1{n^s}=\prod_{{{p{\text{ prime}}}}}\frac1{1-p^{-s}}.
\end{equation}
If instead there exist infinitely many primes (which, again, is the real case, see Exercise~\ref{EUCL}), then we still obtain~\eqref{02i0rjf94yjcxjcjn2nrgt1ds} by taking the limit as~$N\to+\infty$ in~\eqref{09879iujh7ygf5tgb234rfsdfghj284P1}.

From~\eqref{02i0rjf94yjcxjcjn2nrgt1ds} and the definition of the Riemann zeta function in~\eqref{Zzetas} we obtain~\eqref{87uye80217yrg2353fhv5h6h434125}, as desired.

See~\cite[Chapter~7]{MR2063737} for more information about Euler's product formula, which is named there ``the Golden Key'', another proof of it that exploits nice ``cancellations'', its original statement by Euler (in Latin), and its link with the Riemann Hypothesis.

\paragraph{Solution to Exercise~\ref{EUFOanc}.}
This follows from Proposition~\ref{BASEL}, equation~\eqref{Zzetas} and Exercise~\ref{EUFO}.

\paragraph{Solution to Exercise~\ref{EUCL}.} For example, one can give one proof by following Euclid's beautiful line of ideas,
one proof by elementary trigonometry, two proofs using Euler's product formula,
and one proof by estimating the prime-counting function.

First proof: suppose by contradiction that there are only finitely many primes
$$p_1,\dots,p_N>1.$$
Let~$P := p_1p_2\dots p_N+1$. Since~$P$ is strictly larger than any of the above primes, $P$ cannot be a prime
and therefore there exists~$j\in\{1,\dots,N\}$ such that~$p_j$ divides~$P$.
Since~$p_j$ also divides~$p_1p_2\dots p_N$, necessarily~$p_j$ divides~$1$, contradiction.

Second proof: suppose by contradiction that there are only finitely many primes
$$p_1,\dots,p_N>1.$$
Then, for every~$m\in\{1,\dots,N\}$
$$ K_m:=\frac{1}{p_m}\prod_{j=1}^N p_j=p_1\dots p_{m-1}p_{m+1}\dots p_N\in\N$$
and therefore
$$ \sin\left( \frac\pi{p_m}\left(1+2\prod_{j=1}^N p_j\right)\right)=\sin\left( \frac\pi{p_m}+2 K_m\pi\right)=
\sin\left( \frac\pi{p_m}\right).
$$
Using this and the fact that~$\sin\left( \frac\pi{p_m}\right)>0$, we find that
\begin{equation}\label{k0-ur9032ythv6b824cvnc45v7j2wedfghjk}
\prod_{m=1}^N\sin\left( \frac\pi{p_m}\left(1+2\prod_{j=1}^N p_j\right)\right)>0.
\end{equation}

Now we claim that there exists~$m_\star\in\{1,\dots,N\}$ such that
\begin{equation}\label{KSD1limHBSlnerp34}
\frac1{p_{m_\star}}\left(1+2\prod_{j=1}^N p_j\right)\in\N.
\end{equation}
Indeed, either the number
$$ 1+2\prod_{j=1}^N p_j$$
is prime, hence it must be one among~$\{p_1,\dots,p_N\}$, for which the quantity in~\eqref{KSD1limHBSlnerp34} becomes~$1$, or it is not prime, and in this case it can be divided by one of the primes among~$\{p_1,\dots,p_N\}$,
for which the quantity in~\eqref{KSD1limHBSlnerp34} becomes an integer.

In any case, the claim in~\eqref{KSD1limHBSlnerp34} is established, and this yields that
$$ \sin\left( \frac\pi{p_{m_\star}}\left(1+2\prod_{j=1}^N p_j\right)\right)=0.$$

As a result,
$$\prod_{m=1}^N\sin\left( \frac\pi{p_m}\left(1+2\prod_{j=1}^N p_j\right)\right)=0,$$
in contradiction with~\eqref{k0-ur9032ythv6b824cvnc45v7j2wedfghjk}.

This proof was presented in one single line in~\cite{MR3352806}.

Third proof: suppose by contradiction that there are only finitely many primes
$$p_1,\dots,p_N>1.$$
By Exercise~\ref{EUFO} (and we stress that the argument presented to
establish Exercise~\ref{EUFO} does not use that there are infinitely many primes!), for every~$L\in\N$,
$$ \prod_{j=1}^N {\frac{1}{1-p_j^{-1}}}=\lim_{\epsilon\searrow0}
\prod_{j=1}^N {\frac{1}{1-p_j^{-1-\epsilon}}}=
\lim_{\epsilon\searrow0}\zeta(1+\epsilon)\ge
\lim_{\epsilon\searrow0}\sum_{k=1}^L\frac{1}{k^{1+\epsilon}}=\sum_{k=1}^L\frac{1}{k}
.$$
Hence, sending~$L\to+\infty$,
$$ +\infty>\prod_{j=1}^N {\frac{1}{1-p_j^{-1}}}\ge\lim_{L\to+\infty}\sum_{k=1}^L\frac{1}{k}=\sum_{k=1}^{+\infty}\frac{1}{k}=+\infty,$$
contradiction.

Fourth proof: suppose by contradiction that there are only finitely many primes, then the right-hand side of~\eqref{mnhgreg-190r2f439uut9032yf936hgv6b} would be a rational number, but the left-hand side of~\eqref{mnhgreg-190r2f439uut9032yf936hgv6b} is irrational (since~$\pi$ is irrational, see~\cite[Section~16]{zbMATH05043443}), contradiction.

Fifth proof: Let~$\pi(x)$ be the prime-counting function \index{prime-counting function}, i.e. the number of prime numbers less than or equal to some given~$x\in[1,+\infty)$.

Given~$n\in\N$, by the Prime Factorization Theorem we can write it as
$$ n= p_1^{\ell_1}\dots p_{\pi(n)}^{\ell_{\pi(n)}},$$
for suitable primes~$p_1,\dots,p_{\pi(n)}$.

We write~$\ell_j=2\mu_j+\theta_j$, with~$\mu_j\in\N$ and~$\theta_j\in\{0,1\}$ and $$m:=p_1^{\mu_1}\dots p_{\pi(n)}^{\mu_{\pi(n)}}.$$
Notice that~$ m^2\le n$ and
$$ n= m^2\,p_1^{\theta_1}\dots p_{\pi(n)}^{\theta_{\pi(n)}}.$$
Since there are~$m 2^{\pi(n)}\le \sqrt{n} 2^{\pi(n)}$ possible different numbers of this form, we have that
$$ n\le \sqrt{n} 2^{\pi(n)}$$
and therefore
$$ \pi(n)\ge \frac{\log_2 n}{2},$$
which diverges as~$n\to+\infty$, showing that there are infinitely many primes.

See e.g.~\cite{MR1377060} for other proofs.

\paragraph{Solution to Exercise~\ref{HAPRO-0}.} The gist is to pass to the limit Euler's product formula.
Here are the technical details. Suppose not. Then, there exists~$N$ large enough such that
\begin{equation}\label{RI} \prod_{{p{\text{ prime}}}\atop{p>N}}\frac1{1-p^{-1}}\le 1.\end{equation}
Let now~$s\in(1,+\infty)$. Since
$$\frac1{1-p^{-s}}\le\frac1{1-p^{-1}},$$
we deduce from~\eqref{RI} that
\begin{equation*} \prod_{{p{\text{ prime}}}\atop{p>N}}\frac1{1-p^{-s}}\le1.\end{equation*}
Then, by Euler's product formula~\eqref{87uye80217yrg2353fhv5h6h434125}, for all~$s\in(1,+\infty)$ and~$M\in\N$,
$$\prod_{{p{\text{ prime}}}\atop{1<p\le N}}{\frac{1}{1-p^{-s}}}+1\ge
\prod_{{p{\text{ prime}}}\atop{p>1}}{\frac{1}{1-p^{-s}}}
=\sum_{k=1}^{+\infty}\frac1{k^s}\ge\sum_{k=1}^{M}\frac1{k^s}.$$
We thus send~$s\searrow1$ and we get that
$$\prod_{{p{\text{ prime}}}\atop{1<p\le N}}{\frac{1}{1-p^{-1}}}+1\ge\sum_{k=1}^{M}\frac1{k}.$$
We now send~$N\to+\infty$ and we find that
$$\prod_{{p{\text{ prime}}}\atop{p>1}}{\frac{1}{1-p^{-1}}}+1\ge\sum_{k=1}^{M}\frac1{k}$$
and therefore, sending now~$M\to+\infty$,
$$\prod_{{p{\text{ prime}}}\atop{p>1}}{\frac{1}{1-p^{-1}}}+1\ge+\infty,$$
in contradiction with our initial assumption.

\paragraph{Solution to Exercise~\ref{HAPRO}.} This is a refinement of an argument already presented in the solution to Exercise~\ref{AOJsnkwe9rgeotg2-34rt} in which we use Euler's product formula and the Taylor expansion of the logarithm.

Specifically, by a Taylor expansion we have that
\begin{equation}\label{02jf90ewnv4578-21366}
-\sum_{{p{\text{ prime}}}\atop{p>1}} \ln\left(1-\frac1p\right)
=\sum_{{p{\text{ prime}}}\atop{p>1}}\sum_{k=1}^{+\infty}\frac{1}{kp^k}.
\end{equation}

Moreover, since, for all~$k\ge2$,
$$ \sum_{{p{\text{ prime}}}\atop{p>1}}\frac{1}{p^k}\le\sum_{j=2}^{+\infty}\frac1{j^k}\le\int_1^{+\infty}\frac{dx}{x^k}=\frac1{k-1},$$
we obtain that
\begin{eqnarray*}
&& \sum_{{p{\text{ prime}}}\atop{p>1}}\sum_{k=2}^{+\infty}\frac{1}{kp^k}
\le\sum_{k=2}^{+\infty}\frac{1}{k(k-1)},
\end{eqnarray*}
which is a convergent series.

This and~\eqref{02jf90ewnv4578-21366} give that the series in~\eqref{02jf90ewnv4578-213660} and the left-hand side of~\eqref{02jf90ewnv4578-21366} only differ by a constant:
more specifically, we have proved that
$$ S:=\sum_{{p{\text{ prime}}}\atop{p>1}}\sum_{k=2}^{+\infty}\frac{1}{kp^k}<+\infty$$
and therefore, by~\eqref{02jf90ewnv4578-21366},
$$ -\sum_{{p{\text{ prime}}}\atop{p>1}} \ln\left(1-\frac1p\right)-S
=\sum_{{p{\text{ prime}}}\atop{p>1}}\frac{1}{p}.$$
Hence, the desired result is equivalent to
\begin{equation}\label{02jf90ewnv4578-213660b} -\sum_{{p{\text{ prime}}}\atop{p>1}} \ln\left(1-\frac1p\right)=+\infty.\end{equation}
To check this, we notice that, for every~$N\in\N$,
\begin{eqnarray*}&&
-\sum_{{p{\text{ prime}}}\atop{p>1}} \ln\left(1-\frac1p\right)=\sum_{{p{\text{ prime}}}\atop{p>1}} \ln\frac1{1-p^{-1}}\ge
\sum_{{p{\text{ prime}}}\atop{1<p\le N}} \ln\frac1{1-p^{-1}}=
\ln\left(\prod_{{p{\text{ prime}}}\atop{1<p\le N}}\frac1{1-p^{-1}}\right).
\end{eqnarray*}
Hence, taking the limit as~$N\to+\infty$,
$$ -\sum_{{p{\text{ prime}}}\atop{p>1}} \ln\left(1-\frac1p\right)\ge\ln\left(\prod_{{p{\text{ prime}}}\atop{p>1}}\frac1{1-p^{-1}}\right),$$ and the latter is infinite (see Exercise~\ref{HAPRO-0}), which completes the proof of~\eqref{02jf90ewnv4578-213660b}.

Incidentally, this also gives yet another (slightly different) proof that there are infinitely many primes, to be compared with Exercise~\ref{EUCL}.

\paragraph{Solution to Exercise~\ref{x1kdo021ps6lungajkd}.}
By Exercise~\ref{0ei02r3290234pfp42} (taking~$x:=\frac12$ here and noticing that~$\cos(k\pi)=(-1)^k$) we have that
$$ \frac1{2^4}=\frac1{80}+\sum_{k=1}^{+\infty}\frac{\pi^2 k^2-6}{2\pi^4 k^4},$$
from which one obtains the desired result. 

\paragraph{Solution to Exercise~\ref{KALMS0owk3klxslkfw-e}}
By Exercise~\ref{ka-2} and Theorem~\ref{BASw}, we know that, for all~$x\in[0,1]$,
$$x(1-x)=\frac16-\sum_{k\in\Z\setminus\{0\}}\frac{1}{2\pi^2 k^2}\,e^{2\pi ikx},$$
with uniform convergence of the series.

Then, we can use Proposition~\ref{AJSaa} and obtain that Parseval's Identity in~\eqref{PARS} holds true,
yielding that
$$ \frac1{36}+\sum_{k\in\Z\setminus\{0\}}\frac{1}{4\pi^4 k^4}=
\int_0^1 x^2(1-x)^2\,dx=\frac1{30},$$
which gives the desired result.

\paragraph{Solution to Exercise~\ref{aujsoJs893405}}

The required series equals
\begin{equation}\label{mci12} \frac{\pi^2}2\left(\frac56-\frac{2e}{e^2-1}\right).\end{equation}
To check this, we extend~$f$ to a periodic function of period~$1$ and we observe that this is a continuous, even function.

Thus, we can consider its Fourier Series in trigonometric form
(recall~\eqref{jasmx23er} and
Exercise~\ref{smc203e2rf436b.2900rjm4on}). In this context, since
\[ 4\int_0^{1/2} \cosh(2x)\,\cos(2\pi kx)\,dx=
\frac{ (-1)^k ( e^2-1)}{e( \pi^2 k^2 + 1)},
\]
the Fourier Series of~$f$ takes the form
\begin{equation}\label{90OJSldfw0eojrf2} \frac{ e^2-1}{2e}
+\sum_{k=1}^{+\infty}\frac{ (-1)^k ( e^2-1)}{e( \pi^2 k^2 + 1)}\,\cos(2\pi kx).\end{equation}
This series resembles that in~\eqref{90OJSldfw0eojrf}, the main difference being the absence of a factor~$k^2$ in the denominator: this suggests that we need to integrate the above expression twice to obtain the desired result.
More precisely, by virtue of Theorem~\ref{THCOL2FB}(i), we know that the series in~\eqref{90OJSldfw0eojrf2}
converges to~$f$ in~$L^2((0,1))$ and, as a consequence, for every~$x\in\left[0,\frac12\right]$,
\begin{eqnarray*}
\sinh(2x)&=&2\int_0^x \cosh(2\eta)\,d\eta
\\&=&\int_0^x\left[
\frac{ e^2-1}{e}
+2\sum_{k=1}^{+\infty}\frac{ (-1)^k ( e^2-1)}{e( \pi^2 k^2 + 1)}\,\cos(2\pi k\eta)
\right]\,d\eta\\&=&
\frac{ (e^2-1)x}{e}
+\sum_{k=1}^{+\infty}\frac{ (-1)^k ( e^2-1)}{e\pi k( \pi^2 k^2 + 1)}\,\sin(2\pi kx).
\end{eqnarray*}
We note the above series converges uniformly in~$\left[0,\frac12\right]$ and thus, for all~$x\in\left[0,\frac12\right]$
(and, in fact, for all~$x\in\left[-\frac12,\frac12\right]$ by even symmetry),
\begin{equation}\label{90OJSldfw0eojrf26}
\begin{split}
\cosh(2x)-1&=
2\int_0^x \sinh(2\eta) d\eta\\ &=
\int_0^x\left[
\frac{2 (e^2-1)\eta}{e}
+2\sum_{k=1}^{+\infty}\frac{ (-1)^k ( e^2-1)}{e\pi k( \pi^2 k^2 + 1)}\,\sin(2\pi k\eta)\right]\,d\eta\\&=
\frac{(e^2-1)x^2}{e}
-\sum_{k=1}^{+\infty}\frac{ (-1)^k ( e^2-1)}{e\pi^2 k^2( \pi^2 k^2 + 1)}\,
\big(\cos(2\pi kx)-1\big).
\end{split}\end{equation}

In addition, by means of
Exercise~\ref{ka-2} and Theorem~\ref{THCOL2FB}(i), we know that
$$x(1-x)=\frac16+\sum_{k\in\Z\setminus\{0\}}
\frac{e^{2\pi ikx}}{2\pi^2 k^2},$$
with convergence in~$L^2((0,1))$ and therefore, replacing~$x$ with~$x+\frac12$,
$$\frac14-x^2=\frac16+\sum_{k\in\Z\setminus\{0\}}
\frac{e^{2\pi ik\left(x+\frac12\right)}}{2\pi^2 k^2}=
\frac16+\sum_{k=1}^{+\infty}
\frac{\cos\left(2\pi k x+\pi k\right)}{\pi^2 k^2}=\frac16+\sum_{k=1}^{+\infty}
\frac{(-1)^k\cos (2\pi k x)}{\pi^2 k^2}
,$$
with convergence in~$L^2\left(\left(-\frac12,\frac12\right)\right)$.

Combining together this observation, \eqref{90OJSldfw0eojrf2} and~\eqref{90OJSldfw0eojrf26}, we conclude that, in~$L^2\left(\left(-\frac12,\frac12\right)\right)$,
\begin{eqnarray*}&&
\frac{ e^2-1}{2e}
+\sum_{k=1}^{+\infty}\frac{ (-1)^k ( e^2-1)}{e( \pi^2 k^2 + 1)}\,\cos(2\pi kx)=\cosh(2x)\\&&\qquad=
\frac{(e^2-1)x^2}{e}
-\sum_{k=1}^{+\infty}\frac{ (-1)^k ( e^2-1)}{e\pi^2 k^2( \pi^2 k^2 + 1)}\,
\big(\cos(2\pi kx)-1\big)\\&&\qquad=
\frac{e^2-1}{e}
\left(
\frac1{12}-\sum_{k=1}
\frac{(-1)^k\cos (2\pi k x)}{\pi^2 k^2}\right)
-\sum_{k=1}^{+\infty}\frac{ (-1)^k ( e^2-1)\,\cos(2\pi kx)}{e\pi^2 k^2( \pi^2 k^2 + 1)}+1
\\&&\qquad\qquad
+\sum_{k=1}^{+\infty}\frac{ (-1)^k ( e^2-1)}{e\pi^2 k^2( \pi^2 k^2 + 1)}+1.
\end{eqnarray*}
Hence, considering the zero-order Fourier coefficients of the above expansions,
\begin{eqnarray*}
\frac{ e^2-1}{2e}=
\frac{e^2-1}{12e}
+\sum_{k=1}^{+\infty}\frac{ (-1)^k ( e^2-1)}{e\pi^2 k^2( \pi^2 k^2 + 1)}+1,
\end{eqnarray*}
yielding the desired result in~\eqref{mci12}.

\section{Solutions to selected exercises of Section~\ref{CSLE:A:SEZZ0-1}}

\paragraph{Solution to Exercise~\ref{BNCDVCuPXasjmdc.1e}.} Equality in~\eqref{lk2340-012.cnh} takes place if and only if equality takes place in~\eqref{BNCDVCuPXasjmdc.1}, that is if and only if~$|\widehat f_k|=|k\,\widehat f_k|$ for all~$k\in\Z$.
As a consequence, $\widehat f_k=0$ for all~$k\in\Z\setminus\{-1,\,1\}$ and thus equality in~\eqref{lk2340-012.cnh} is equivalent to~$f(x)=a\cos(2\pi x)+b\sin(2\pi x)$ for some~$a$, $b\in\R$.

We also remark that the existence of these functions attaining equality in~\eqref{lk2340-012.cnh} shows that the constant
in the right-hand side of~\eqref{lk2340-012.cnh} is optimal.

\paragraph{Solution to Exercise~\ref{BNCDVCuPXasjmdc.1e.102i}.} 
The idea is (up to scaling) to repeat~$f$ in an odd symmetric fashion and periodise, then apply~\eqref{lk2340-012.cnh} and obtain the desired estimates. The technical details go as follows.
We define
\begin{equation}\label{12.BNCDVCuPXasjmdc.1e.102i.qw2edfmv013}
g(x):=f(2x).\end{equation}
In this way, $g\in C^1\left(\left[ 0,\frac12\right]\right)$, with
\begin{equation}\label{BNCDVCuPXasjmdc.1e.102i.qw2edfmv013}
g(0)=g\left(\frac12\right)=0.\end{equation}

Thus, we extend~$g$ in~$\left(-\frac12,0\right)$ via an odd reflection.
Then, we extend~$g$ to a periodic function of period~$1$ (and, with a slight abuse this notation, we denote this extension by~$g$ as well). Notice that such a periodic function~$g$ is, by construction, odd symmetric.

Since~$g$ may develop corners at points of~$\frac{\Z}2$, we consider a smooth, even, periodic function~$\tau_\epsilon:\R\to[0,1]$ of period~$\frac12$ such that~$\tau_\epsilon=0$ in~$[-\epsilon,\epsilon]$, $\tau_\epsilon=1$ in~$\left[2\epsilon,\frac12-2\epsilon\right]$ and~$|\tau_\epsilon'|\le\frac2\epsilon$. We also consider the function~$g_\epsilon:=g\tau_\epsilon$ which is now periodic of period~$1$, belongs to~$C^1(\R)$ and, owing to its odd symmetry,
$$ \int_0^1g_\epsilon(x)\,dx=0.$$
On this account, we can use~\eqref{lk2340-012.cnh} and conclude that
\begin{equation}\label{lk2340-012.cnh.09-1.p2}
\int_0^1 |g_\epsilon(x)|^2\,dx\le\frac1{4\pi^2}\int_0^1 |g_\epsilon'(x)|^2\,dx.
\end{equation}

Moreover, since~$g_\epsilon(x)\to g(x)$ for all~$x\in(0,1)$ as~$\epsilon\searrow0$, we deduce from Fatou's Lemma that
\begin{equation}\label{lk2340-012.cnh.09-1.p1}\liminf_{\epsilon\searrow0}
\int_0^1 |g_\epsilon(x)|^2\,dx\ge\int_0^1 |g(x)|^2\,dx=2\int_0^{\frac12} |g(x)|^2\,dx=\int_0^1 |f(y)|^2\,dy.
\end{equation}

In addition, if
$$I_\epsilon:=(0,2\epsilon)\cup\left(\frac12-2\epsilon,\frac12+2\epsilon\right)\cup(1-2\epsilon,1),$$
we have that~$\tau_\epsilon'=0$ in~$(0,1)\setminus I_\epsilon$ and, as a consequence of~\eqref{BNCDVCuPXasjmdc.1e.102i.qw2edfmv013}, for all~$x\in I_\epsilon$,
$$|g(x)|\le 4\epsilon\|g'\|_{L^\infty((0,1/2))}=
8\epsilon\|f'\|_{L^\infty((0,1))}.$$
Consequently,
\begin{eqnarray*}&&
\| g'-g_\epsilon' \|_{L^2((0,1))}\\&&\qquad\le\| g'-g'\tau_\epsilon \|_{L^2((0,1))}+\| g\tau_\epsilon' \|_{L^2((0,1))}\\
&&\qquad\le\sqrt{\|g'\|^2_{L^\infty((0,1/2))}\int_0^1 (1-\tau_\epsilon(x))^2\,dx}+\sqrt{
\sup_{x\in I_\epsilon}|g(x)|^2
\int_{I_\epsilon} |\tau_\epsilon'(x)|^2 \,dx}\\
&&\qquad\le\sqrt{4\|f'\|^2_{L^\infty((0,1))}\int_0^1 (1-\tau_\epsilon(x))^2\,dx}+16\sqrt{8\epsilon}\,\|f'\|_{L^\infty((0,1))}.\end{eqnarray*}
The Dominated Convergence Theorem thereby entails that
\begin{eqnarray*}
\lim_{\epsilon\searrow0}\| g'-g_\epsilon' \|_{L^2((0,1))}=0\end{eqnarray*}
and correspondingly
\begin{eqnarray*}
\liminf_{\epsilon\searrow0}\| g'_\epsilon \|_{L^2((0,1))}
\le\liminf_{\epsilon\searrow0}\Big(\| g'\|_{L^2((0,1))}+\| g'-g_\epsilon' \|_{L^2((0,1))}\Big)=\| g'\|_{L^2((0,1))}.
\end{eqnarray*}

After this, \eqref{lk2340-012.cnh.09-1.p2}, and~\eqref{lk2340-012.cnh.09-1.p1}, we gather that
\begin{eqnarray*}&&
\int_0^1 |f(y)|^2\,dy\le
\liminf_{\epsilon\searrow0}\int_0^1 |g_\epsilon(x)|^2\,dx\le\frac1{4\pi^2}\liminf_{\epsilon\searrow0}\int_0^1 |g_\epsilon'(x)|^2\,dx\\&&\qquad\le\frac1{4\pi^2}\int_0^1 |g'(x)|^2\,dx=\frac1{2\pi^2}\int_0^{\frac12} |g'(x)|^2\,dx=
\frac2{\pi^2}\int_0^{\frac12} |f'(2x)|^2\,dx\\&&\qquad=\frac1{\pi^2}\int_0^1 |f'(y)|^2\,dy,
\end{eqnarray*}
as desired.

\paragraph{Solution to Exercise~\ref{BNCDVCuPXasjmdc.1e.102}.} For equality to hold in~\eqref{lk2340-012.cnh.102},
we have that equality must hold for~$g_\epsilon$ in~\eqref{lk2340-012.cnh.09-1.p2}
(and thus, up to taking limits, for~$g$). Therefore, by Exercise~\ref{BNCDVCuPXasjmdc.1e}, $g(x)=
a\cos(2\pi x)+b\sin(2\pi x)$ for some~$a$, $b\in\R$. But since~$g$ is odd, necessarily~$a=0$.

Hence, after~\eqref{12.BNCDVCuPXasjmdc.1e.102i.qw2edfmv013}, we have that~$f(x)=g\left(\frac{x}2\right)=
b\sin(\pi x)$ for some~$b\in\R$.

\paragraph{Solution to Exercise~\ref{CSLE:A}.} The gist is to argue for a contradiction and construct a trigonometric polynomial with ``too many zeros'' (with some care to avoid multiplicity issues). The technical details go as follows.

First of all, 
\begin{equation}\label{SUDFIS}
{\mbox{it suffices to prove~\eqref{lk2340-012} for all~$y\in\left[0,\frac{1}{2N}\right]$,}}\end{equation}
because, if that holds true, then we can consider the trigonometric polynomial~$\widetilde f(x):=f(x_0-x)$, which satisfies
$$ \widetilde f(0)=f(x_0)=\max_{x\in\R}|f(x)|=\max_{x\in\R}|\widetilde f(x)|$$
and write that, for all~$y\in\left[-\frac{1}{2N},0\right)$,
$$ f(x_0+y)=\widetilde f(-y)\ge \widetilde f(0)\,\cos(-2\pi Ny)=f(x_0)\,\cos(2\pi Ny).$$

Moreover, without loss of generality, one can assume that, for all~$y\in(0,1)$,
\begin{equation}\label{plk2340-012}
f(x_0+y)<f(x_0)
\end{equation}
Indeed, for all~$\epsilon>0$, one can look at
$$ f_\epsilon(x):=f(x)+\epsilon\cos(2\pi(x-x_0)),$$
notice that~$f_\epsilon$ is also a trigonometric polynomial of degree~$N\ge1$ and that, for all~$y\in(0,1)$,
$$ f_\epsilon(x_0+y)=f(x_0+y)+\epsilon\cos(2\pi y)<f(x_0+y)+\epsilon\le f(x_0)+\epsilon=f_\epsilon(x_0).$$
Hence, if the desired result holds true under the additional assumption~\eqref{plk2340-012}, thus in particular for~$f_\epsilon$,
we have that, for all~$y\in\left[0,\frac{1}{2N}\right]$,
$$ f(x_0+y)+\epsilon\cos(2\pi y)=f_\epsilon(x_0+y)\ge f_\epsilon(x_0)\,\cos(2\pi Ny)=
\big(f(x_0)+\epsilon\big)\,\cos(2\pi Ny)$$
and so the desired claim for the original~$f$ follows by sending~$\epsilon\searrow0$.

After these preliminary observations, we deal with the core argument.
To prove~\eqref{lk2340-012}, suppose not. Then, by~\eqref{SUDFIS},
there exists~$y_*\in\left[0,\frac{1}{2N}\right]$ such that
$$ f(x_0+y_*)< f(x_0)\,\cos(2\pi Ny_*).$$
Thus, by continuity, we can find
\begin{equation}\label{y0-01} y_0\in\left[0,\frac{1}{2N}\right)\end{equation} such that
\begin{equation}\label{y0-02} f(x_0+y_0)< f(x_0)\,\cos(2\pi Ny_0).\end{equation}
Notice that~$y_0\ne0$ (otherwise~$f(x_0)< f(x_0)\,\cos(0)= f(x_0)$, contradiction).
Hence, we consider the points~$\zeta_k:=\frac{k}{2N}$, with~$k\in\{1,\dots,2N\}$.

We define the trigonometric polynomial of degree~$N$
$$ g(y):=f(x_0+y)- f(x_0) \,\cos(2\pi Ny).$$
We point out that, for all~$k\in\{0,\dots,2N\}$,
$$ g(\zeta_k)=f(x_0+\zeta_k)- f(x_0) \,\cos(k\pi)=f(x_0+\zeta_k)-(-1)^k f(x_0)$$
and so, in light of~\eqref{plk2340-012},
\begin{equation} \begin{split}\label{y0-01-2}
&{\mbox{if $k$ is odd, then~$g(\zeta_k)>0$,}}\\
&{\mbox{if $k$ is even, then~$g(\zeta_k)<0$.}}\end{split}\end{equation}

Consequently, 
\begin{equation}\label{vv9Hn5tCedfTgRcpwsd-1}\begin{split}&
{\mbox{$g$ possesses at least one zero in each}}\\&{\mbox{of the open intervals~$(\zeta_k,\zeta_{k+1})$,
for~$k\in\{1,\dots,2N-1\}$.}}\end{split}\end{equation}

Additionally, $g(y_0)<0$, due to~\eqref{y0-02}, and~$g(\zeta_1)>0$, due to~\eqref{y0-01-2}.
Thus, in view of~\eqref{y0-01}, we have that
\begin{equation}\label{vv9Hn5tCedfTgRcpwsd-2}
{\mbox{$g$ possesses at least one zero in the open interval~$(y_0,\zeta_1)$.}}\end{equation}

Furthermore, $g(0)=f(x_0)- f(x_0) \,\cos(0)=0$.
As a result, owing to~\eqref{vv9Hn5tCedfTgRcpwsd-1} and~\eqref{vv9Hn5tCedfTgRcpwsd-2}, we have that~$g$ possesses at least~$(2N-1)+1+1=2N+1$ zeros,
in contradiction with Exercise~\ref{ILTRI}.

\paragraph{Solution to Exercise~\ref{CSLE:B}.} We can assume that~$N\ge1$, otherwise~$f$ is constant and the result is obvious.

Suppose first that~$f:\R\to\R$.
Let~$x_0\in\R$ be such that
\begin{equation}\label{FDEX0} |f'(x_0)|=\max_{x\in\R}|f'(x)|.\end{equation}
Without loss of generality, up to replacing~$f$ with~$-f$, we can assume that~$f'(x_0)\ge0$.
We can thereby use Exercise~\ref{CSLE:A} applied to the trigonometric polynomial~$f'$ and find that, for every~$y\in\left[-\frac{1}{2N},\frac{1}{2N}\right]$,
$$ f'(x_0+y)\ge f'(x_0)\,\cos(2\pi Ny).$$
We integrate this inequality over~$y\in\left(-\frac{1}{4N},\frac{1}{4N}\right)$ and we conclude that
\begin{eqnarray*}&&
2\max_{x\in\R}|f(x)|\ge
f\left(x_0+\frac{1}{4N}\right)-f\left(x_0-\frac{1}{4N}\right)=
\int_{-\frac{1}{4N}}^{\frac{1}{4N}}f'(x_0+y)\,dy\\&&\qquad\ge f'(x_0)\,\int_{-\frac{1}{4N}}^{\frac{1}{4N}}\cos(2\pi Ny)\,dy=
\frac{f'(x_0)}{\pi N}=\frac{1}{\pi N}\max_{x\in\R}|f'(x)|.
\end{eqnarray*}
This proves~\eqref{BVELE} when~$f:\R\to\R$.

If~$f:\R\to\C$, we use the following rotation trick: given~$x_0$ as in~\eqref{FDEX0}, we
pick~$\alpha\in[0,2\pi)$ such that
$$e^{i\alpha} f'(x_0)\in[0,+\infty).$$ Then, we
apply the result that we have just proven to the real-valued trigonometric polynomial~$g(x):=\Re\big(e^{i\alpha}f(x)\big)$.

In this way, 
$$ |g(x)|\le|e^{i\alpha}f(x)|=|f(x)|$$
and
\begin{eqnarray*}&& |g'(x)|=\big| \Re\big(e^{i\alpha}f'(x)\big)\big|\le \big|e^{i\alpha}f'(x)\big|=|f'(x)|\le
|f'(x_0)|=e^{ i\alpha} f'(x_0)\\&&\qquad\qquad\qquad\qquad=
\Re\big(e^{i\alpha}f'(x_0)\big)=g'(x_0),
\end{eqnarray*}
yielding that
$$ \max_{x\in\R}|g'(x)|=g'(x_0)=|f'(x_0)|.$$
Hence, applying~\eqref{BVELE} to the real-valued trigonometric polynomial~$g$,
\begin{eqnarray*}|f'(x_0)|=
\max_{x\in\R}|g'(x)|\le 2\pi N\max_{x\in\R}|g(x)|\le2\pi N\max_{x\in\R}|f(x)|,
\end{eqnarray*}
as desired.

\paragraph{Solution to Exercise~\ref{BVELEb}.}
One can take, for instance, $\cos(2\pi Nx)$.

\section{Solutions to selected exercises of Section~\ref{CSLE:A:GEO}}

\paragraph{Solution to Exercise~\ref{kmsdcJAx.iton4.143.e}.}
We use the assumptions and notations displayed in the proof of Theorem~\ref{kmsdcJAx.iton4.143}.
Besides, up to a translation, we can suppose that the barycenter of the curve under consideration lies at the origin, namely
$$\int_0^1x(t)\,dt=0=\int_0^1y(t)\,dt.$$
We also notice that, given~$\sigma>0$, for all~$a$, $b\ge0$, it follows from the Cauchy-Schwarz Inequality that
$$ 2ab\le \sigma a^2+\frac{b^2}\sigma.$$

Then, by the Poincar\'e-Wirtinger Inequality of Theorem~\ref{CSLE:A:SEZZ0-1-pwin},
\begin{eqnarray*}
2\int_0^1 \big|x(t)y'(t)\big|\,dt\le
\int_0^1 
\left( \sigma|x(t)|^2+\frac{|y'(t)|^2}\sigma\right)\,dt\le\int_0^1 
\left( \frac{\sigma\,|x'(t)|^2}{4\pi^2}+\frac{|y'(t)|^2}\sigma\right)\,dt
\end{eqnarray*}
and similarly
\begin{eqnarray*}
2\int_0^1 \big|x'(t)y(t)\big|\,dt\le\int_0^1 
\left( \frac{|x'(t)|^2}\sigma+\frac{\sigma\,|y'(t)|^2}{4\pi^2}\right)\,dt.
\end{eqnarray*}

These observations, combined with~\eqref{SMXC22.12wedmMKAx}, entail that
\begin{eqnarray*}&&4\pi{\mathcal{A}}=2\pi\int_0^1 
\big(x(t)y'(t)-y(t)x'(t)\big)\,dt\le\pi\left( \frac\sigma{4\pi^2}+\frac1\sigma\right)
\int_0^1 
\left( |x'(t)|^2+|y'(t)|^2\right)\,dt\\&&\qquad\qquad=\pi\left( \frac\sigma{4\pi^2}+\frac1\sigma\right)
\int_0^1 |\gamma'(t)|^2\,dt=\pi\left( \frac\sigma{4\pi^2}+\frac1\sigma\right).\end{eqnarray*}
Choosing~$\sigma:=2\pi$, we conclude that~$4\pi{\mathcal{A}}\le1=4\pi{\mathcal{L}}$, as desired.

We also stress that the equality sign is attained if and only if equality holds in the
Cauchy-Schwarz Inequality and in the
Poincar\'e-Wirtinger Inequality of Theorem~\ref{CSLE:A:SEZZ0-1-pwin}, see Exercise~\ref{BNCDVCuPXasjmdc.1e}, that is
\begin{equation*}\begin{dcases}&|y'(t)|=\sigma |x(t)|=2\pi |x(t)|,\\ &
|x'(t)|=\sigma |y(t)|=2\pi| y(t)|,\\ &
x(t)=a\cos(2\pi t)+b\sin(2\pi t), \\ &y(t)=c\cos(2\pi t)+d\sin(2\pi t),\end{dcases}\end{equation*}
for some~$a$, $b$, $c$, $d\in\R$.

Since
$$\begin{dcases}\displaystyle& |c\cos(2\pi t)+ d\sin(2\pi t)|=| y(t)|=\frac{|x'(t)|}{2\pi }=|-a\sin(2\pi t)+ b\cos(2\pi t)|,\\
\displaystyle& |a\cos(2\pi t)+ b\sin(2\pi t)|=| x(t)|=\frac{|y'(t)|}{2\pi }=|-c\sin(2\pi t)+ d\cos(2\pi t)|,
\end{dcases}$$
we find (picking~$t=0$) that~$c\in\{b,-b\}$ and~$d\in\{a,-a\}$, as well as
(picking~$t=1/8$) that~$ab=-cd$.

As a result,
\begin{eqnarray*}&&
(x(t))^2+(y(t))^2=(a^2+c^2)\cos^2(2\pi t)+(b^2+d^2)\sin^2(2\pi t)+2(ab+cd)\sin(2\pi t)\cos(2\pi t)\\
&&\qquad=a^2+c^2,
\end{eqnarray*}which is a circle.

\paragraph{Solution to Exercise~\ref{GEOEX:per}.}
The parameterization in~\eqref{RTSm} will be obtained from the
``normal angle'' (also known as ``turning angle''). That is, given~$t\in\R$, we consider the vector~$\nu_t$ identified by the complex number~$e^{it}$ and the halfplane with exterior normal equal to~$\nu_t$. By parallel sliding such an halfplane in the direction of~$\nu_t$, we can touch the curve at a point~$p_t$ at which the exterior unit normal of the curve coincides with~$\nu_t$ (which coincides with the complex number~$e^{it}$). Correspondingly, the unit tangent direction of the curve coincides with a $(\pi/2)$-rotation of~$\nu_t$, i.e., up to choosing an orientation of the curve, with~$ie^{it}$.

We thus consider the parameterization of the curve given by~$\gamma(t):=p_t$.
By construction, the tangent vector~$ \gamma'(t)$ is proportional to~$ie^{it}$,
which is~\eqref{RTSm}.

Actually, by~\eqref{RTSm} it also follows that this parameterization is of period~$2\pi$, in view of the uniqueness of the normal (hence the tangent) vector
for strictly convex curves. 

\paragraph{Solution to Exercise~\ref{GEOEX}.} We identify the plane~$\R^2$ with~$\C$
and we write the closed curve~$\gamma$ as a Fourier Series, say
$$ \gamma(t)=\sum_{k=-\infty}^{+\infty} \gamma_k \,e^{i\omega kt},$$
for some~$\omega>0$ and coefficients~$\gamma_k\in\C$, and we have that
\begin{eqnarray*}
\gamma(t)-\gamma\left(t+\frac\pi{\omega n}\right)&=&
\sum_{k=-\infty}^{+\infty} \gamma_k \left(e^{i\omega kt}-e^{i\omega k\left(t+\frac\pi{\omega n}\right)}\right)\\&=&\sum_{k=-\infty}^{+\infty} \gamma_k\, e^{i\omega k t}\left(1-e^{\frac{i\pi k}{n}}\right).
\end{eqnarray*}

Since the curve is smooth, the decay of the Fourier coefficients (see Theorem~\ref{0CINFITH1})
and the Fubini-Tonelli Theorem guarantee that
one can swap the order of integration and sum and obtain that
\begin{eqnarray*}\int_0^{2\pi/\omega}\left(
\gamma(t)-\gamma\left(t+\frac\pi{\omega n}\right)\right)\,e^{-i\omega n t}\,dt&=&\sum_{k=-\infty}^{+\infty}\int_0^{2\pi/\omega} \gamma_k \,e^{i\omega (k-n)t}\left(1-e^{\frac{i\pi k}{n}}\right)\,dt.
\end{eqnarray*}
Now,
$$ \int_0^{2\pi/\omega} e^{i\omega(k-n) t}\,dt=\begin{cases}
2\pi/\omega & {\mbox{ if }}k=n,\\0&{\mbox{ otherwise,}}
\end{cases}$$yielding that
\begin{equation}\label{INT}\int_0^{2\pi/\omega}\left(
\gamma(t)-\gamma\left(t+\frac\pi{\omega n}\right)\right)\,e^{-i\omega n t}\,dt=
\frac{2\pi \gamma_n}\omega \,\left(1-e^{\frac{i\pi n}{n}}\right)=\frac{4\pi \gamma_n}\omega.
\end{equation}

So far, the convexity of the curve has not been used, and the calculations have
been developed for any parameterization of the curve.
We now use the convexity assumption (in fact, a slightly stronger
convexity assumption)
to obtain a parameterization
in which the tangent vector has the particularly simple expression described in Exercise~\ref{GEOEX:per}.
That is, to establish~\eqref{GEOEX:INEQ}, we assume that the curve is strictly convex and we will exploit~\eqref{RTSm} (once
strictly convex curves are settled, one can use an approximation argument,
shadow a convex curve by a strictly convex one, obtain the final inequality for this
sequence of approximating curves,
pass the final inequality to the limit,
and obtain its validity for convex curves).

The periodicity of the parameterization in Exercise~\ref{GEOEX:per}
also says that we can use~\eqref{INT} with~$\omega=1$
in this setting. In particular, choosing~$n=1$ in~\eqref{INT}, we obtain that
\begin{eqnarray*}
4\pi \gamma_1=
\int_0^{2\pi}\left(
\gamma(t)-\gamma\left(t+\pi\right)\right)\,e^{-i t}\,dt\le
2\pi{\mathcal{D}}.
\end{eqnarray*}

Now, by~\eqref{RTSm},
\begin{eqnarray*}
\mu(t)=
-i\gamma'(t)\,e^{-it}=
\sum_{k=-\infty}^{+\infty} k \gamma_k \,e^{i k t}\,e^{-it}=\sum_{k=-\infty}^{+\infty} k \gamma_k \,e^{i (k-1) t}
\end{eqnarray*}
and therefore
\begin{eqnarray*}
\int_0^{2\pi}\mu(t)\,dt=
\int_0^{2\pi}\sum_{k=-\infty}^{+\infty} k \gamma_k \,e^{i (k-1) t}\,dt=
2\pi\gamma_1\le\pi{\mathcal{D}}.
\end{eqnarray*}

But~\eqref{RTSm} also gives that
$$ |\gamma'(t)|=|i\mu(t)\,e^{it}|=|\mu(t)|=\mu(t),$$
and consequently
$$ {\mathcal{L}}=\int_0^{2\pi}\mu(t)\,dt\le\pi{\mathcal{D}},$$
as desired.

\paragraph{Solution to Exercise~\ref{NOAKMSDceikfMKsCfDtgbdcAtVkmdcYI-1}.}
No. For example, one can look at a smooth curve obtained by joining as many concentric circles as we want.

\begin{figure}[h]
\includegraphics[height=4.1cm]{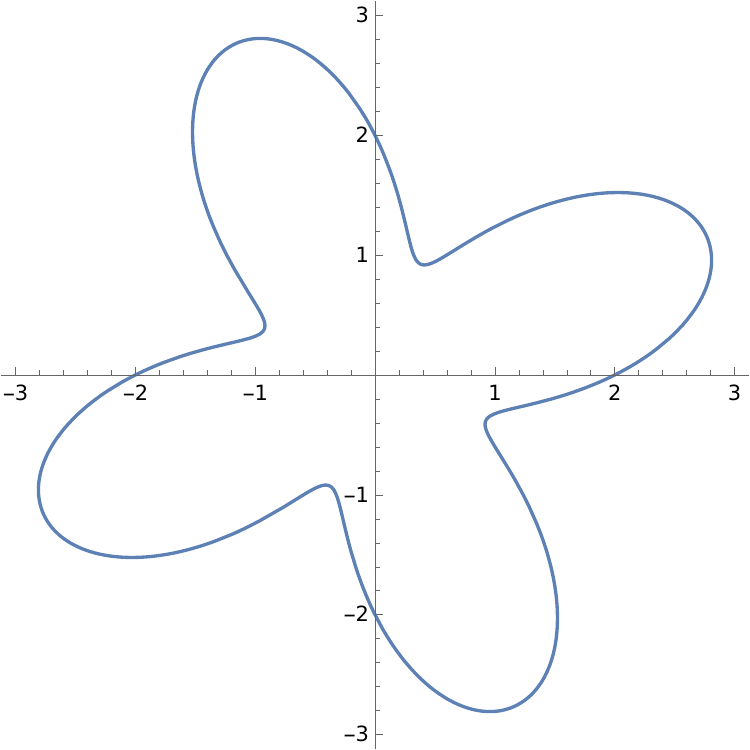}$\,\;\qquad\quad$\includegraphics[height=4.1cm]{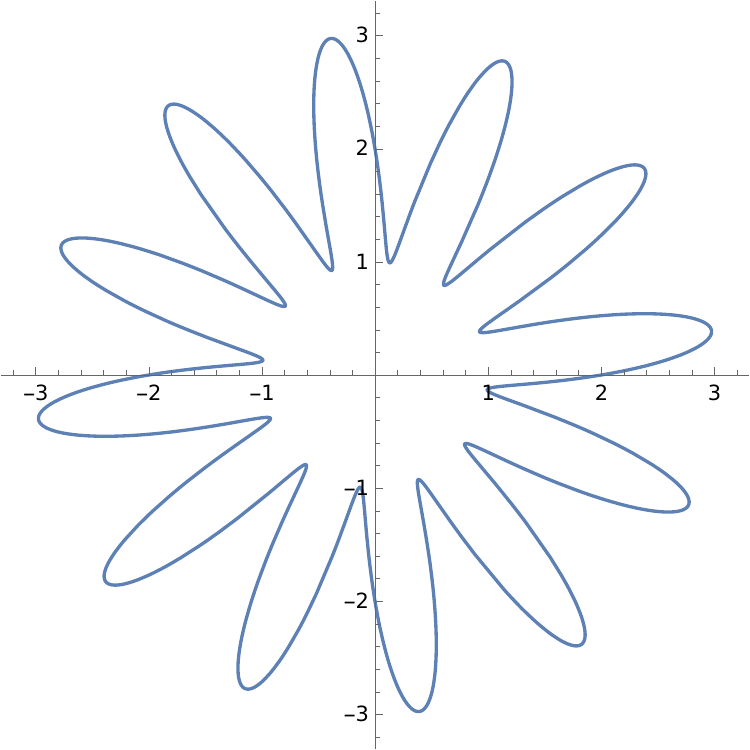} $\,\;\qquad\quad$\includegraphics[height=4.1cm]{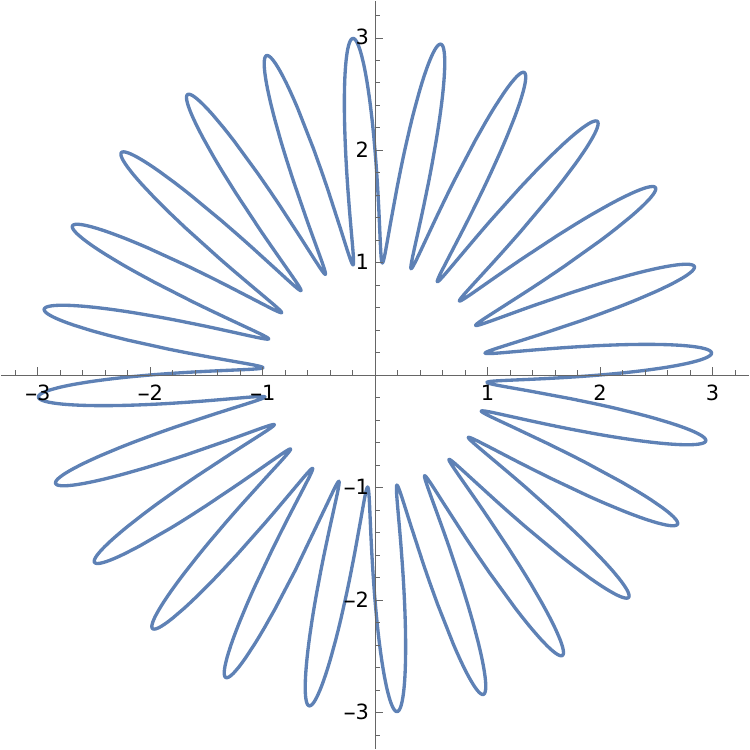}
\centering
\caption{The curve in~\eqref{TCuddBESGINMoTHASqer-1} with~$N\in\{4,12,24\}$.}\label{RET12P.91o2eujrfne.023oriSTQEILSRCNMKtk}
\end{figure}

\paragraph{Solution to Exercise~\ref{NOAKMSDceikfMKsCfDtgbdcAtVkmdcYI-2}.}
No, starshapedness is not enough either. As a counterexample, one can
look at
\begin{equation}\label{TCuddBESGINMoTHASqer-1}
\gamma(t)=(2+\sin(Nt))e^{it}\end{equation} for~$N\in\N$ large,
see Figure~\ref{RET12P.91o2eujrfne.023oriSTQEILSRCNMKtk}.

\paragraph{Solution to Exercise~\ref{NOAKMSDceikfMKsCfDtgbdcAtVkmdcYI-2b}.}
Yes, indeed the inequality in~\eqref{GEOEX:INEQ}
is optimal, since it is attained (among other curves) by the circle.

\paragraph{Solution to Exercise~\ref{NOAKMSDceikfMKsCfDtgbdcAtVkmdcYI-3}.}
No, convex curves do not necessarily satisfy the equality in~\eqref{GEOEX:INEQ}: for example, one can consider an approximation of a square.

\section{Solutions to selected exercises of Section~\ref{CSLE:A:WEIE}}

\paragraph{Solution to Exercise~\ref{LMSTNKNOLZOJOqdU}.}
For instance, we can use complex exponentials:
\begin{eqnarray*}&&
\cos(m\theta)+\cos((m-2)\theta)-2\cos\theta\,\cos((m-1)\theta)\\&&\qquad=
\Re\big(e^{im\theta}\big)+\Re\big(e^{i(m-2)\theta}\big)-\Re\big(2\cos\theta\,\cos((m-1)\theta)\big)\\&&\qquad=
\Re\left(e^{im\theta}+e^{i(m-2)\theta}-\frac{1}{2}\big(e^{i\theta}+e^{-i\theta}\big)\big(e^{i(m-1)\theta}+e^{-i(m-1)\theta}\big)\right)\\&&\qquad=\frac12\Re\Big(
\big( e^{im\theta}-e^{-im\theta}\big)+\big(e^{i(m-2)\theta}-e^{-i(m-2)\theta}\big)\Big)\\&&\qquad=0.
\end{eqnarray*}

\paragraph{Solution to Exercise~\ref{LMSTNKNOLZOJOqdU2}.} We argue by induction.
When~$N=0$, we can take~$T_0(t):=1$. When~$N=1$, we can take~$T_1(t):=t$.

Now, suppose that the claim holds true for the indices~$\{0,\dots,N-1\}$, with~$N\ge1$. We employ
Exercise~\ref{LMSTNKNOLZOJOqdU} and, from the inductive hypothesis, we obtain that
\begin{eqnarray*}
\cos(N\theta)&=&2\cos\theta\,\cos((N-1)\theta)-\cos((N-2)\theta)\\
&=&2\cos\theta\,T_{N-1}(\cos\theta)-T_{N-2}(\cos\theta).
\end{eqnarray*}
The desired result thus follows by setting
$$ T_N(t):=2t\,T_{N-1}(t)-T_{N-2}(t).$$
We stress that the coefficient in front of the term of higher degree for~$T_N$ is twice
the one for~$T_{N-1}$, which also proves~\eqref{eTGBSq01ijm43g-01i24}.

\section{Solutions to selected exercises of Section~\ref{RADOSE}}

\paragraph{Solution to Exercise~\ref{RADOEXILE}.}
The answer is~$\sqrt\pi\,e^{-p^2}$.

To get this result, we proceed as follows. We use the notation in~\eqref{RAD:TRA:DE1}.
Since the vectors~$(-\sin\phi,\cos\phi)$ and~$(\cos\phi,\sin\phi)$ are orthogonal and of unit length,
we see that~$|r_{p,\phi}(t)|^2=t^2+p^2$. Hence, if~$f(X):=e^{-|X|^2}$,
we deduce from~\eqref{RAD:TRA:DE} that
$${\mathcal{R}}_f(p,\phi)=\int_{-\infty}^{+\infty} e^{-t^2-p^2}\,dt=\sqrt\pi\,e^{-p^2}.$$

\paragraph{Solution to Exercise~\ref{RADOEXILE.2}.} In a sense yes, even if we focused only
on the case of ``nice'' functions when we presented the
Radon Transfrom in~\eqref{RAD:TRA:DE}.

Indeed, formally one can apply~\eqref{RAD:TRA:DE} with~$f:=\delta_0$ and obtain that
\begin{eqnarray*}
{\mathcal{R}}_{\delta_0}(p,\phi)=\int_{-\infty}^{+\infty}\delta_0\big(r_{p,\phi}(t)\big)\,dt=
\begin{dcases}
1 & {\mbox{ if the straight line~$r_{p,\phi}$ passes through the origin,}}\\
0 & {\mbox{ otherwise.}}
\end{dcases}\end{eqnarray*}
That is,
$$ {\mathcal{R}}_{\delta_0}(p,\phi)=\begin{dcases}
1 & {\mbox{ if $p=0$,}}\\
0 & {\mbox{ otherwise.}}
\end{dcases}$$

\paragraph{Solution to Exercise~\ref{ANSBEHINGFVA6tYTRAhJAS}.}
Since the Laplacian is invariant under rotation (see e.g.~\cite{ELEM2}), given two mutually orthogonal unit vectors~$\varpi_1$ and~$\varpi_2$ we have that $$\Delta u=\partial_{\varpi_1}^2 u+\partial_{\varpi_2}^2 u,$$
where, for all~$j\in\{1,2\}$, we used the standard notation~$\partial_{\varpi_j}:=\varpi_j\cdot\nabla$ for directional derivatives.

In particular, since the vectors~$\varpi_1:=(-\sin\phi,\cos\phi)$ and~$\varpi_2:=(\cos\phi,\sin\phi)$ are orthogonal and of unit length, we have that
\begin{eqnarray*} \Delta u\big(r_{p,\phi}(t)\big)&=&\partial_{\varpi_1}^2 u\big(\varpi_1t+\varpi_2p\big)+\partial_{\varpi_2}^2 u\big(\varpi_1t+\varpi_2p\big)\\&=&\partial_t^2\Big( u\big(\varpi_1t+\varpi_2p\big)\Big)+\partial^2_p\Big(u\big(\varpi_1t+\varpi_2p\big)\Big).\end{eqnarray*}

As a result, by~\eqref{RAD:TRA:DE} and the fact that~$u$ is compactly supported (and so are its derivatives),
\begin{eqnarray*}
{\mathcal{R}}_{\Delta u}(p,\phi)&=&\int_{-\infty}^{+\infty}\Delta u\big(r_{p,\phi}(t)\big)\,dt\\
&=&\int_{-\infty}^{+\infty}\partial_{t}^2\Big( u\big(\varpi_1t+\varpi_2p\big)\Big)\,dt
+\int_{-\infty}^{+\infty}\partial^2_p\Big(u\big(\varpi_1t+\varpi_2p\big)\Big)\,dt\\&=&0
+\partial^2_p\left(\int_{-\infty}^{+\infty}u\big(\varpi_1t+\varpi_2p\big)\,dt\right)\\&=&\partial^2_p{\mathcal{R}}_u(p,\phi),
\end{eqnarray*}as desired.

The observation pointed out in this exercise is interesting because it suggests that the Radon Transform is useful in partial differential equations, for instance, because it reduces operators acting on several variables to derivatives in a single direction.

\paragraph{Solution to Exercise~\ref{RADOEX1}.} By~\eqref{RADODE1},
we know that~${\mathcal{T}}_0(x)=1$ and~${\mathcal{T}}_1(x)=x$.

Consequently,
$$ I_0(a)=\int_1^a \frac{a\,dx}{x\,\sqrt{(a^2-x^2)(x^2-1)}}$$
and, substituting for~$y:=\frac{a}x$,
$$ I_1(a)=\int_1^a \frac{x\,dx}{\sqrt{(a^2-x^2)(x^2-1)}}=
\int_1^a \frac{a\,dy}{y\,\sqrt{(y^2-1)(a^2-y^2)}}=I_0(a).$$

Instead, substituting above for~$t:=\sqrt{1-\frac{x^2}{a^2}}$ and recognising the derivative of the inverse of the sine function,
$$ I_1(a)=\int_0^{\sqrt{1-\frac1{a^2}}}\frac{dt}{\sqrt{1-\frac1{a^2}-t^2}}\,dt=\frac\pi2.$$

\paragraph{Solution to Exercise~\ref{RADOEX2}.} We recall that, for all~$\xi$, $\eta\in\R$,
$$\cosh(\xi+\eta)=\cosh \xi\cosh \eta+\sinh \xi\sinh \eta.$$
Particularly, for all~$n\in\N\cap[1,+\infty)$ and~$\xi\in\R$,
$$\cosh((n+1)\xi)=\cosh (n\xi)\cosh \xi+\sinh(n \xi)\sinh\xi$$
and
$$\cosh((n-1)\xi)=\cosh (n\xi)\cosh\xi-\sinh(n \xi)\sinh \xi.$$
Thus, for all~$x\ge1$, we infer from~\eqref{RADODE1} that
\begin{eqnarray*} {\mathcal{T}}_{n+1}(x)&=&
\cosh\big( (n+1)\arccosh x\big)\\&=&
\cosh \big( n\arccosh x\big)\cosh\big( \arccosh x\big)+\sinh\big( n\arccosh x\big)\sinh\big( \arccosh x\big)\\&=&
{\mathcal{A}}+{\mathcal{B}}
,\end{eqnarray*}
where we used the short notation
\begin{eqnarray*}&& {\mathcal{A}}:=\cosh \big( n\arccosh x\big)\cosh\big( \arccosh x\big)=
x\cosh \big( n\arccosh x\big)\\
{\mbox{and}}\qquad&&
{\mathcal{B}}:=\sinh\big( n\arccosh x\big)\sinh\big( \arccosh x\big)=
\sqrt{x^2-1}\,\sinh\big( n\arccosh x\big)
,\end{eqnarray*}
and similarly
$${\mathcal{T}}_{n-1}(x)={\mathcal{A}}-{\mathcal{B}}.$$

We also recall that, for all~$\xi$, $\eta\in\R$,
$$\cos(\xi+\eta)=\cos \xi\cos \eta-\sin \xi\sin \eta.$$
In particular, for all~$n\in\N\cap[1,+\infty)$ and~$\xi\in\R$,
$$\cos((n+1)\xi)=\cos (n\xi)\cos\xi-\sin(n \xi)\sin \xi$$
and
$$\cos((n-1)\xi)=\cos (n\xi)\cos\xi+\sin(n \xi)\sin \xi.$$
Then, if~$x\in(1,a)$, we deduce from~\eqref{RADODE1} that
\begin{eqnarray*} {\mathcal{T}}_{n+1}\left(\frac{x}a\right)&=&
\cos\left( ( n+1)\arccos\frac{x}a\right)\\&=&
\cos\left( n\arccos\frac{x}a\right)\cos\left( \arccos\frac{x}a\right)-\sin\left( n\arccos\frac{x}a\right)\sin\left( \arccos\frac{x}a\right)\\&=&{\mathcal{C}}-{\mathcal{D}},
\end{eqnarray*}
where
\begin{eqnarray*}&& {\mathcal{C}}:=\cos\left( n\arccos\frac{x}a\right)\cos\left( \arccos\frac{x}a\right)=
\frac{x}a\,\cos\left( n\arccos\frac{x}a\right)
\\
{\mbox{and}}\qquad&&
{\mathcal{D}}:=\sin\left( n\arccos\frac{x}a\right)\sin\left( \arccos\frac{x}a\right)=
\sqrt{1-\frac{x^2}{a^2}}\,\sin\left( n\arccos\frac{x}a\right),\end{eqnarray*}
and, in the same vein,
$$ {\mathcal{T}}_{n-1}\left(\frac{x}a\right)={\mathcal{C}}+{\mathcal{D}}.$$

Consequently, we find the following nice simplification:
\begin{eqnarray*}
&&{\mathcal{T}}_{n+1}(x)\,{\mathcal{T}}_{n+1}\left(\frac{x}a\right)-{\mathcal{T}}_{n-1}(x)\,{\mathcal{T}}_{n-1}\left(\frac{x}a\right)=({\mathcal{A}}+{\mathcal{B}})({\mathcal{C}}-{\mathcal{D}})-
({\mathcal{A}}-{\mathcal{B}})({\mathcal{C}}+{\mathcal{D}})\\
&&\qquad=({\mathcal{A}}{\mathcal{C}}-{\mathcal{A}}{\mathcal{D}}+{\mathcal{B}}{\mathcal{C}}-{\mathcal{B}}{\mathcal{D}})-
({\mathcal{A}}{\mathcal{C}}+{\mathcal{A}}{\mathcal{D}}-{\mathcal{B}}{\mathcal{C}}-{\mathcal{B}}{\mathcal{D}})\\&&\qquad=
2({\mathcal{B}}{\mathcal{C}}-{\mathcal{A}}{\mathcal{D}}).
\end{eqnarray*}

On this account and~\eqref{RADODE1.gfbv}, we gather that
\begin{equation}\label{NS-kxc-RADODE1.gfbv}\begin{split}
I_{n+1}(a)-I_{n-1}(a)&=\int_1^a \frac{a\,
\left[ {\mathcal{T}}_{n+1}(x)\,{\mathcal{T}}_{n+1}\left(\frac{x}a\right)-
{\mathcal{T}}_{n-1}(x)\,{\mathcal{T}}_{n-1}\left(\frac{x}a\right)
\right]}{x\,\sqrt{(a^2-x^2)(x^2-1)}}\,dx\\
&=\int_1^a \frac{2a\,
\left( {\mathcal{B}}{\mathcal{C}}-{\mathcal{A}}{\mathcal{D}}
\right)}{x\,\sqrt{(a^2-x^2)(x^2-1)}}\,dx
.\end{split}
\end{equation}

It is also useful to observe that{\footnotesize{
\begin{eqnarray*}
&&\frac{a\,\left( {\mathcal{B}}{\mathcal{C}}-{\mathcal{A}}{\mathcal{D}}
\right)}{x\,\sqrt{(a^2-x^2)(x^2-1)}}\\
&&\quad=\frac{a\,
\left(\frac{x}a\, \sqrt{x^2-1}\,\sinh\big( n\arccosh x\big)\cos\left( n\arccos\frac{x}a\right)
-x\,\sqrt{1-\frac{x^2}{a^2}}\,\cosh \big( n\arccosh x\big)\,\sin\left( n\arccos\frac{x}a\right)
\right)}{x\,\sqrt{(a^2-x^2)(x^2-1)}}\\
&&\quad=\frac{\sinh\big( n\arccosh x\big)\cos\left( n\arccos\frac{x}a\right)
}{\sqrt{a^2-x^2}}
-\frac{\cosh \big( n\arccosh x\big)\,\sin\left( n\arccos\frac{x}a\right)}{\sqrt{x^2-1}}.
\end{eqnarray*}}}

Now, to spot additional cancellations, we observe that
$$ \partial_x\left( \sinh\big( n\arccosh x\big)\right)=\frac{n \cosh\big(n \arccosh x\big)}{\sqrt{x^2-1}}
$$
and
$$ \partial_x\left(\sin\left( n\arccos\frac{x}a\right)\right)=-\frac{
n\,\cos\left( n\arccos\frac{x}a\right)}{\sqrt{a^2-x^2}},
$$
leading to
\begin{eqnarray*}
&&-\frac{na\,\left( {\mathcal{B}}{\mathcal{C}}-{\mathcal{A}}{\mathcal{D}}\right)}{x\,\sqrt{(a^2-x^2)(x^2-1)}}\\
&&\quad=\sinh\big( n\arccosh x\big)\,\partial_x\left(\sin\left( n\arccos\frac{x}a\right)\right)
+\partial_x\left( \sinh\big( n\arccosh x\big)\right)\,\sin\left( n\arccos\frac{x}a\right)\\
&&\quad=\partial_x\left( \sinh\big( n\arccosh x\big)\,\sin\left( n\arccos\frac{x}a\right)\right).
\end{eqnarray*}

This and~\eqref{NS-kxc-RADODE1.gfbv} return that
\begin{eqnarray*}&&
-\frac{n\big(I_{n+1}(a)-I_{n-1}(a)\big)}2\\&&\quad=-
\int_1^a \frac{na\,
\left( {\mathcal{B}}{\mathcal{C}}-{\mathcal{A}}{\mathcal{D}}
\right)}{x\,\sqrt{(a^2-x^2)(x^2-1)}}\,dx\\&&\quad=
\int_1^a \partial_x\left( \sinh\big( n\arccosh x\big)\,\sin\left( n\arccos\frac{x}a\right)\right)\,dx\\&&\quad=
\sinh\big( n\arccosh a\big)\,\sin\left( n\arccos 1\right)
-\sinh\big( n\arccosh 1\big)\,\sin\left( n\arccos\frac{1}a\right)\\&&\quad=
\sinh\big( n\arccosh a\big)\,\sin(0)
-\sinh(0)\,\sin\left( n\arccos\frac{1}a\right)\\&&\quad=0.
\end{eqnarray*}
Accordingly, we obtain that~$I_{n+1}(a)=I_{n-1}(a)$, which is the desired result.

\paragraph{Solution to Exercise~\ref{RADOEXX}.} We argue by induction over~$m$. When~$m\in\{0,1\}$,
the desired claim follows from Exercise~\ref{RADOEX1}. 

Let us now take~$m\ge2$. By virtue of Exercise~\ref{RADOEX2} and the inductive hypothesis, we see that~$I_m(a)=I_{m-2}(a)=\frac\pi2$, as desired.

\paragraph{Solution to Exercise~\ref{RADOEXXCR}.} Several proofs are possible.
The proof that we present here
relies on calculus, expressing the area element on the space of lines in terms of the arclength parameter and the angle made by the line with the curve.

Let~${\mathcal{L}}$ be the length of~$\gamma$. We parameterise~$\gamma$ by its arclength, denoting this parameterization by~$(x(t),y(t))$, with~$t\in[0,{\mathcal{L}})$. 

Since, by the arclength setting, the vector~$(x'(t),y'(t))$ has unit length, we write~$x'(t)=\cos\tau(t)$ and~$y'(t)=\sin\tau(t)$, for some angular function~$\tau(t)$
(from the geometric point of view~$\tau(t)$ represents
the angle between the tangent to the curve at the point~$(x(t),y(t))$ and the horizontal coordinate axis).

For each~$\theta\in[0,2\pi)$ and~$t\in[0,{\mathcal{L}})$, we define
\begin{equation}\label{EplacGHnEG.1}
\begin{split}&p(\theta,t):=
x(t)\cos\left(\theta+\tau(t)-\frac\pi2\right) +y(t)\sin\left(\theta+\tau(t)-\frac\pi2\right),\\&
\phi(\theta,t):= \theta+\tau(t)-\frac\pi2,\\ {\mbox{and }}\;\quad&
f(\theta,t):=\left(p(\theta,t),\;
\phi(\theta,t)\right).\end{split}\end{equation}

We see that
\begin{equation*}\begin{split}
\left| \det\frac{\partial f}{\partial(\theta,t)}\right|&=
\left| \det\left( \begin{matrix}
-x\sin\phi +y\cos\phi&\,
&x'\,\cos\phi +y'\,\sin\phi
-x\tau'\sin\phi +y\tau'\cos\phi
\\ 1&\,& \tau'
\end{matrix}\right)\right|\\&=\big|
x'\,\cos\phi +y'\,\sin\phi\big|\\&=\big|
x'\,\sin(\theta+\tau)-y'\,\cos(\theta+\tau)\big|\\&=\big|
x'\,(\sin\theta\cos\tau+\cos\theta\sin\tau)-y'\,(\cos\theta\cos\tau-\sin\theta\sin\tau)\big|\\&=\big|
\sin\theta\cos^2\tau+\cos\theta\sin\tau\cos\tau-\cos\theta\sin\tau\cos\tau+\sin\theta\sin^2\tau|\\&=|\sin\theta|.\end{split}
\end{equation*}

For this reason, employing the Change of Variables Formula (see e.g.~\cite[Theorem~2 on page~99]{MR3409135}) we obtain that
\begin{equation}\label{EplacGHnEG.14}
\iint_{[0,2\pi)\times[0,{\mathcal{L}})}|\sin\theta|\,d\theta\,dt=\iint_{\R^2}
N_\star(p,\phi)\,dp\,d\phi,
\end{equation}
where~$N_\star(p,\phi)$ is the number of elements~$(\theta,t)$ for which~$f(\theta,t)=(p,\phi)$.

\begin{figure}[h]
\includegraphics[height=5.99cm]{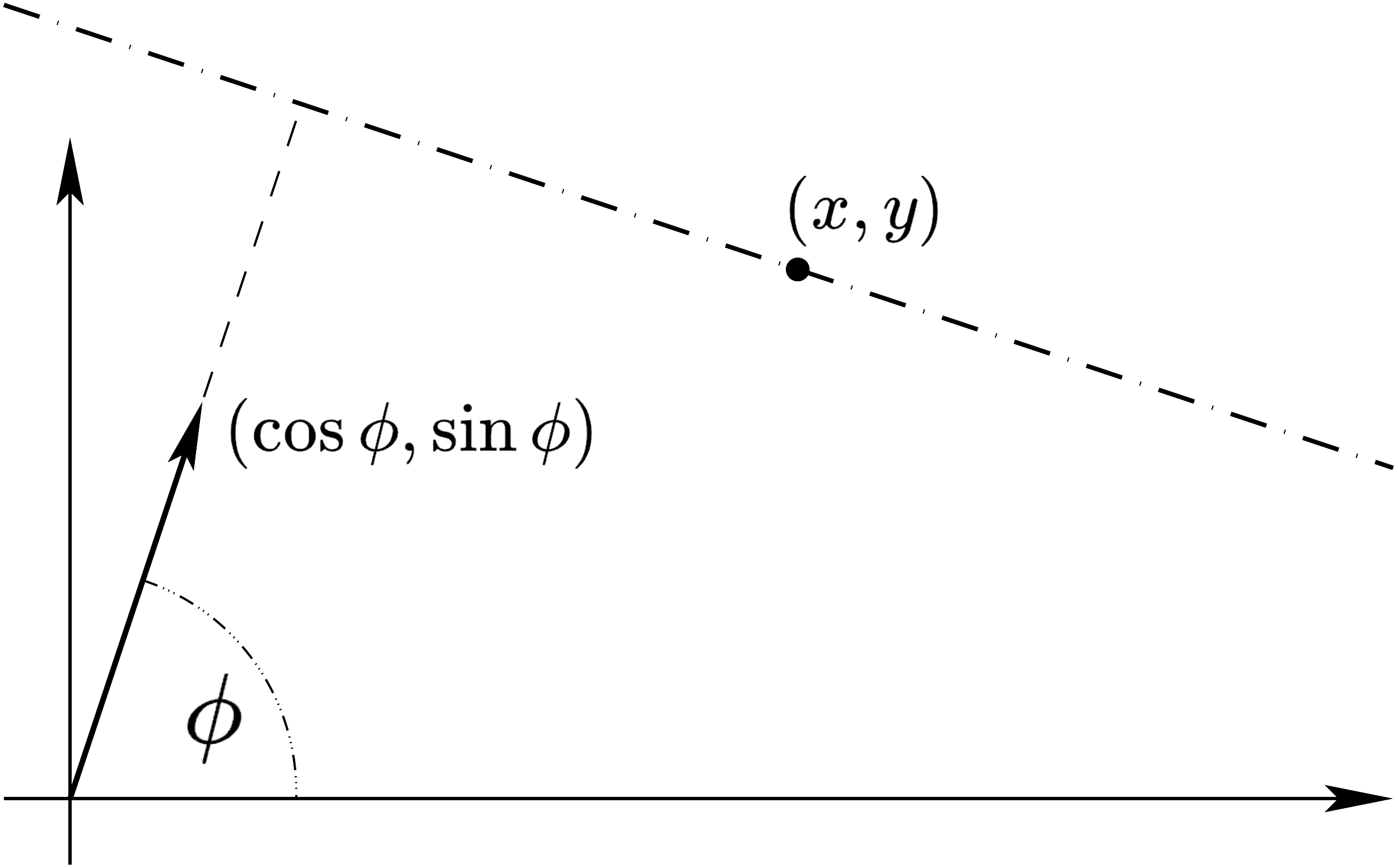}
\centering
\caption{A point~$(x,y)$ belongs to the straight line~$r_{p,\phi}$ if and only if~$p=(\cos\phi,\sin\phi)\cdot(x(t),y(t))$.}\label{RnaBba90.qpwd.BVA-1mnlocET12P.91o2eujrfne.023oritk}
\end{figure}

Now, in light of~\eqref{EplacGHnEG.1}, we observe that~$f(\theta,t)=(p,\phi)$ if and only if~$p=x(t)\cos\phi +y(t)\sin\phi=
(\cos\phi,\sin\phi)\cdot(x(t),y(t))$, that is if any only if the point~$(x(t),y(t))$ of the curve belongs\footnote{If one wants an analytical proof of this fact, rather than relying on Figure~\ref{RnaBba90.qpwd.BVA-1mnlocET12P.91o2eujrfne.023oritk}, it is possible to argue as follows.

Suppose, on the one hand that a point~$(x,y)\in\R^2$ lies on the straight line~$r_{p,\phi}$.
Then, by~\eqref{RAD:TRA:DE1}, there exists~$T\in\R$ such that~$(x,y)=
(-\sin\phi,\cos\phi)T+(\cos\phi,\sin\phi)p$. This gives that~$
(\cos\phi,\sin\phi)\cdot(x,y)=(\cos\phi,\sin\phi)\cdot(-\sin\phi,\cos\phi)T+(\cos\phi,\sin\phi)\cdot(\cos\phi,\sin\phi)p=p$.

Now suppose, on the other hand, that~$p=(\cos\phi,\sin\phi)\cdot(x,y)=x\cos\phi+y\sin\phi$.
Let~$T:=y\cos\phi-x\sin\phi$. Then,
\begin{eqnarray*}
&&(-\sin\phi,\cos\phi)T+(\cos\phi,\sin\phi)p=
(-\sin\phi,\cos\phi)(y\cos\phi-x\sin\phi)+(\cos\phi,\sin\phi)p\\&&\qquad=
( -y\sin\phi\cos\phi+x\sin^2\phi+p\cos\phi,\,
y\cos^2\phi-x\sin\phi\cos\phi+p\sin\phi)\\&&\qquad=
\big( -y\sin\phi\cos\phi+x\sin^2\phi+(x\cos\phi+y\sin\phi)\cos\phi,\,
y\cos^2\phi-x\sin\phi\cos\phi+(x\cos\phi+y\sin\phi)\sin\phi\big)=(x,y),
\end{eqnarray*}
whence~$(x,y)\in\R^2$ lies on the straight line~$r_{p,\phi}$.}
to the straight line~$r_{p,\phi}$ in the notation of~\eqref{RAD:TRA:DE1}, see Figure~\ref{RnaBba90.qpwd.BVA-1mnlocET12P.91o2eujrfne.023oritk} 

There is one caveat, though: here~$p$ can also be negative,
while in~\eqref{RAD:TRA:DE1} one restricted to positive quantities, so we are counting
the straight lines twice (one with positive~$p$, one with negative~$p$, which coincides with the first up to rotating~$\phi$ by~$\pi$). These observations give that~$N_\star(p,\phi)=2N_\gamma(p,\phi)$ and accordingly
\begin{equation}\label{EplacGHnEG.15}
\iint_{\R\times[0,2\pi)}
N_\star(p,\phi)\,dp\,d\phi=2\iint_{[0,+\infty)\times[0,2\pi)}
N_\gamma(p,\phi)\,dp\,d\phi.
\end{equation}

Furthermore, by a direct calculation,
$$ \int_0^{2\pi}|\sin\theta|\,d\theta=4.$$

Combining this, \eqref{EplacGHnEG.14}, and~\eqref{EplacGHnEG.15}, we conclude that
$$ 4{\mathcal{L}}=2\iint_{[0,+\infty)\times[0,2\pi)}
N_\gamma(p,\phi)\,dp\,d\phi,$$
whence
$$ {\mathcal{L}}=\frac12\iint_{[0,+\infty)\times[0,2\pi)}
N_\gamma(p,\phi)\,dp\,d\phi,$$
which is the desired Crofton Formula.

See~\cite[Chapter~3, Section~2]{MR2162874} for more details on this type of arguments and their links to differential forms,
and~\cite[Chapter~3]{MR2168892} for a different method of proof.

See also~\cite{MR60840, MR488215, MR2162874, MR3469669} for an introduction to integral geometry
and~\cite{MR1836759} for a historical perspective.

\paragraph{Solution to Exercise~\ref{RADOEXXISP.02.234}.}
Since the length and the intersection number are invariant under rigid motions, we may suppose that the circle is centred at the origin.
Thus, to compute the quantity in~\eqref{LA97529SFOISMOsmSjA}, pick any~$\phi\in[0,2\pi)$ and note that the straight line~$r_{p,\phi}$ in~\eqref{RAD:TRA:DE1} intersects the unit circle twice when~$p\in[0,1)$ and once when~$p=1$ (also, no intersections take place when~$p>1$). This gives that, for all~$\phi\in[0,2\pi)$,
$$ N_\gamma(p,\phi)=\begin{dcases} 2 & {\mbox{ if }}p\in[0,1),
\\ 1 & {\mbox{ if }}p=1,\\0 & {\mbox{ if }}p\in[1,+\infty),
\end{dcases}$$
and accordingly
$$ \frac12\iint_{[0,+\infty)\times[0,2\pi)} N_\gamma(p,\phi)\,dp\,d\phi=
\frac12\int_0^{2\pi} \left( \int_0^12\,dp\right)\,d\phi=2\pi.$$

\paragraph{Solution to Exercise~\ref{RADOEXXISP.02.234.RSG}.}
Since the notions of both length and intersections are invariant under rigid motions,
without loss of generality we may suppose that this segment has the form~$\{0\}\times[0,\ell]$
(the cases~$\{0\}\times[0,\ell)$, $\{0\}\times(0,\ell]$, and~$\{0\}\times(0,\ell)$ being similar).

In this setting, we have that a straight line~$r_{p,\phi}$ as in~\eqref{RAD:TRA:DE1}
intersects this segment if and only if there exists~$t\in\R$ for which~$t\sin\phi-p\cos\phi=0$
and~$t\cos\phi+p\sin\phi\in[0,\ell]$, that is, if and only\footnote{Strictly speaking, one should distinguish here the case in which~$\sin\phi=0$, but we can ignore this caveat, since the corresponding values of~$\phi$ lie in a set of null Lebesgue measure and therefore will not contribute to the integral formulas.} if~$\frac{p}{\sin\phi}\in[0,\ell]$,
and thus, since~$p\ge0$, if and only if~$\phi\in[0,\pi]$ and~$p\in[0,\ell\sin\phi]$.

Since segments and straight lines cannot intersect more than once, this says that
$$ N_\gamma(p,\phi)=\begin{dcases}
1&{\mbox{ if~$\phi\in[0,\pi]$ and~$p\in[0,\ell\sin\phi]$,}}\\
0&{\mbox{ otherwise}}
\end{dcases}$$
and consequently
$$ \frac12\iint_{[0,+\infty)\times[0,2\pi)}
N_\gamma(p,\phi)\,dp\,d\phi
=\frac12\int_0^\pi\left(\int_0^{\ell\sin\phi}\,dp\right)\,d\phi=\frac\ell2\int_0^\pi \sin\phi\,d\phi=\ell
.$$

\paragraph{Solution to Exercise~\ref{RADOEXXISP.dp}.}
Let~$\phi\in[0,2\pi)$ and consider the straight line~$r_{p,\phi}$ in~\eqref{RAD:TRA:DE1}.
By convexity, either~$r_{p,\phi}$ is tangent to~$\gamma_1$ (which occurs at most at two values of~$p$),
or~$N_{\gamma_1}(p,\phi)=2$, because a line cannot cut transversely a convex curve at three points.

Also, since~$\gamma_2$ contains~$\gamma_1$, whenever~$N_{\gamma_1}(p,\phi)\ne0$,
we have that~$N_{\gamma_2}(p,\phi)\ge2$.

Consequently, if~${\mathcal{L}}_j$ denotes the length of~$\gamma_j$, with~$j\in\{1,2\}$, we infer
from the Crofton Formula in Exercise~\ref{RADOEXXCR} that \begin{eqnarray*}&&{\mathcal{L}}_1
=\frac12\iint_{[0,+\infty)\times[0,2\pi)} N_{\gamma_1}(p,\phi)\,dp\,d\phi
=\frac12\int_0^{2\pi}\left(\int_{{p\ge0}\atop{N_{\gamma_2}(p,\phi)\ne0}}2 \,dp\right)\,d\phi\\&&\qquad\qquad
\le
\frac12\iint_{[0,+\infty)\times[0,2\pi)} N_\gamma(p,\phi)\,dp\,d\phi=
{\mathcal{L}}_2.\end{eqnarray*}

\paragraph{Solution to Exercise~\ref{RADOEXXISP.dpVAR}.}
In the notation of Exercise~\ref{RADOEXXCR}, suppose that~$N_{\gamma_1}(p,\phi)\le N_\star$, for some~$N_\star\in\N$.

Then, by the Crofton Formula,
\begin{equation}\label{QEILMAMDOEXXCR}{\mathcal{L}}_1=
\frac12\iint_{[0,+\infty)\times[0,2\pi)} N_{\gamma_1}(p,\phi)\,dp\,d\phi
\le 
\frac{N_\star}2
\iint_{{[0,+\infty)\times[0,2\pi)}\atop{ N_{\gamma_1}(p,\phi)\ne0}}\,dp\,d\phi.
\end{equation}

Additionally, since~$\gamma_1$ is encircled by~$\gamma_2$, anytime a straight line intersects~$\gamma_1$ it has to intersect~$\gamma_2$ at least twice (once to enter and once to exit the region encompassed by~$\gamma_2$).
This gives that if~$N_{\gamma_1}(p,\phi)\ne0$ then~$N_{\gamma_2}(p,\phi)\ge2$.

Consequently,
$$ \iint_{{[0,+\infty)\times[0,2\pi)}\atop{ N_{\gamma_1}(p,\phi)\ne0}}\,dp\,d\phi\le
\frac12\iint_{{[0,+\infty)\times[0,2\pi)}\atop{ N_{\gamma_1}(p,\phi)\ne0}}
N_{\gamma_2}(p,\phi)\,dp\,d\phi\le{\mathcal{L}}_2.$$

Combining this and~\eqref{QEILMAMDOEXXCR} we conclude that
$$ {\mathcal{L}}_1\le\frac{N_\star\,{\mathcal{L}}_2}2$$
and therefore~$N_\star\ge\frac{2 {\mathcal{L}}_1}{ {\mathcal{L}}_2}$. 

Since~$N_\star\in\N$, this gives that~$N_\star\ge\left\lfloor\frac{2{\mathcal{L}}_1}{{\mathcal{L}}_2}\right\rfloor$,
which is the desired result.

\paragraph{Solution to Exercise~\ref{GECOB}.} 
It is more convenient to think at the equivalent problem in which the needle stays still, e.g. as the segment~$\{0\}\times[0,\ell]$ in the plane, and ``the floor'' is tossed over the needle'', i.e. in the notation of~\eqref{RAD:TRA:DE1},we consider a family of parallel lines~$\big\{r_{p+k,\phi}$, with $k\in\Z\big\}$, where~$p$ and~$\phi$ are randomly chosen (with uniform distribution) in~$[0,+\infty)\times[0,2\pi)$.

Actually, since the configuration obtained would be the same up to exchanging~$\phi$ and~$\pi-\phi$, we can suppose that~$p$ and~$\phi$ are randomly chosen in~$[0,+\infty)\times[0,\pi)$. 

\begin{figure}[h]
\includegraphics[height=5.99cm]{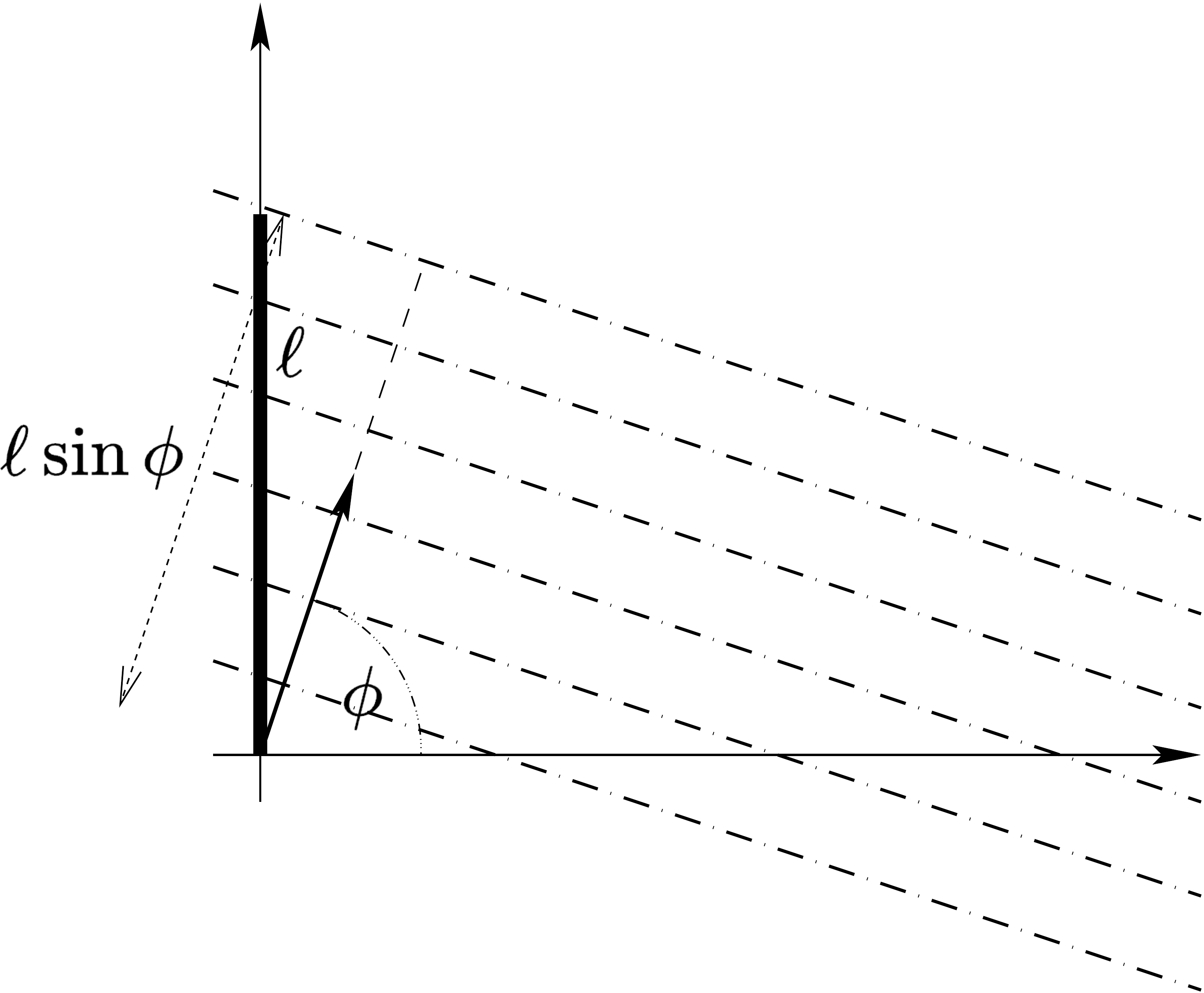}
\centering
\caption{Intersecting a needle with a family of parallel straight lines.}\label{RnaBb2aincC90.qpwd.BVA-1mnlocET12P.91o2eujrfne.023oritk}
\end{figure}

Now, since~$\ell\le1$, we see that if~$r_{p_0,\phi}$ intersects the needle, then~$r_{p,\phi}$ also intersects the needle
for~$p$ in an interval of length~$\ell\sin\phi$, see Figure~\ref{RnaBb2aincC90.qpwd.BVA-1mnlocET12P.91o2eujrfne.023oritk}.

More explicitly,
if~$P_\phi$ collects all the~$p_0\in[0,1)$ for which~$r_{p_0,\phi}$ intersects the needle, then~$P_\phi$ has length~$\ell\sin\phi$. 

Now, we are asked to find the probability that the random family of straight lines intersects the needle, which is the number of occurrences in which~$r_{p_0,\phi}$ intersects the needle, integrated over~$p_0\in[0,1)$ and~$\phi\in[0,\pi)$, divided by the total measure of this space: in view of the above considerations,
this quantity is $$\frac{\displaystyle\int_0^\pi\left(\int_{P_\phi}dp\right)\,d\phi}{\displaystyle\iint_{[0,1)\times[0,\pi)}dp\,d\phi}=\frac{\displaystyle\int_0^\pi\ell\sin\phi\,d\phi}\pi=\frac{2\ell}\pi,$$
as desired.

As a side remark, let us show that this result is compatible with
the Crofton Formula of Exercise~\ref{RADOEXXCR}, with~$\gamma:=\{0\}\times[0,\ell]$ representing the needle. To this end,
in the notation of Exercise~\ref{RADOEXXCR}, we know from the previous considerations that~$N_\gamma(p,\phi)=1$ for all~$p\in P_\phi$, as well as~$N_\gamma(p,\phi)=0$ for all~$p\not\in P_\phi$, and therefore
\begin{eqnarray*}&& \frac12\iint_{[0,+\infty)\times[0,2\pi)} N_\gamma(p,\phi)\,dp\,d\phi=
\iint_{[0,+\infty)\times[0,\pi)} N_\gamma(p,\phi)\,dp\,d\phi\\&&\qquad=
\int_0^\pi\left(\int_{P_\phi}dp\right)\,d\phi=\int_0^\pi\ell\sin\phi\,d\phi=\ell,
\end{eqnarray*}
which is precisely the length of the needle, in agreement with Exercise~\ref{RADOEXXCR}.

\section{Solutions to selected exercises of Section~\ref{DEDESECTDET}}

\paragraph{Solution to Exercise~\ref{EDG:GWHEQmpec.2.DE4}.} We let~$x_0:=0$ and~$x_{\ell+1}:=1$ and
we know that~$f$ is differentiable in the intervals~$
(x_0,x_1)$, $(x_1,x_2)$, $\dots$, $(x_{\ell-1},x_\ell)$, $(x_\ell,x_{\ell+1})$. With this notation, integrating by parts we have that
\begin{eqnarray*}
E_N(x)&=&\frac{i\pi}{N}\sum_{{k\in\Z}\atop{|k|\le N}} k\int_0^1 f(y)\,e^{2\pi ik( x-y)}\,dy\\
&=&\frac{i\pi}{N}\sum_{{k\in\Z}\atop{|k|\le N}}\sum_{j=0}^\ell k\int_{x_j}^{x_{j+1}} f(y)\,e^{2\pi ik( x-y)}\,dy\\
&=&-\frac{1}{2N}\sum_{{k\in\Z}\atop{|k|\le N}}\sum_{j=0}^\ell \int_{x_j}^{x_{j+1}} f(y)\,\frac{d}{dy}\big(e^{2\pi ik( x-y)}\big)\,dy\\
&=&-\frac{1}{2N}\sum_{{k\in\Z}\atop{|k|\le N}}\sum_{j=0}^\ell \Big(f(x_{j+1}^-)\,e^{2\pi ik( x-x_{j+1})}-f(x_j^+)\,e^{2\pi ik( x-x_j)}\Big)\\&&\qquad+\frac{1}{2N}\sum_{{k\in\Z}\atop{|k|\le N}}\sum_{j=0}^\ell \int_{x_j}^{x_{j+1}} f'(y)\,e^{2\pi ik( x-y)}\,dy,
\end{eqnarray*}
where
$$ f(x_j^+):=\lim_{x\searrow x_j}f(x)\qquad{\mbox{and}}\qquad
f(x_j^-):=\lim_{x\nearrow x_j}f(x).$$

We point out that, by~\eqref{EDG:GWHEQmpec.2.DE9},
\begin{eqnarray*}&&
\sum_{j=0}^\ell \Big(f(x_{j+1}^-)\,e^{2\pi ik( x-x_{j+1})}-f(x_j^+)\,e^{2\pi ik( x-x_j)}\Big)\\&&\qquad=
\sum_{j=1}^{\ell+1} f(x_{j}^-)\,e^{2\pi ik( x-x_{j})}
-\sum_{j=0}^\ell f(x_j^+)\,e^{2\pi ik( x-x_j)}\\&&\qquad=
f(1)e^{2\pi ik( x-1)}-f(0)e^{2\pi ik x}+\sum_{j=1}^{\ell} \Big(f(x_{j}^-)-f(x_j^+)\Big)\,e^{2\pi ik( x-x_{j})}\\&&\qquad=
0-\sum_{j=1}^{\ell} \sigma_j\,e^{2\pi ik( x-x_{j})}.
\end{eqnarray*}

From these observations we arrive at
\begin{equation}\label{EDG:GWHEQmpec.2.DE47}\begin{split}
E_N(x)=\frac{1}{2N}\sum_{{k\in\Z}\atop{|k|\le N}}\sum_{j=1}^\ell \sigma_j\,e^{2\pi ik( x-x_{j})}+\frac{1}{2N}\sum_{{k\in\Z}\atop{|k|\le N}}\sum_{j=0}^\ell \int_{x_j}^{x_{j+1}} f'(y)\,e^{2\pi ik( x-y)}\,dy.
\end{split}
\end{equation}
Now, in light of Lemma~\ref{KASMqwdfed123erDKLI},
\begin{eqnarray*}&&
\sup_{x\in[0,1)}\left|\sum_{{k\in\Z}\atop{|k|\le N}}\sum_{j=0}^\ell \int_{x_j}^{x_{j+1}} f'(y)\,e^{2\pi ik( x-y)}\,dy\right|\\&&\qquad=
\sup_{x\in[0,1)}\left|\sum_{j=0}^\ell \int_{x_j}^{x_{j+1}}
\frac{f'(y)\,\sin\big((2N+1)\pi (x-y)\big)}{\sin(\pi( x-y))}
\,dy\right|\\&&\qquad\le M\sup_{x\in[0,1)}\sum_{j=0}^\ell \int_{x_j}^{x_{j+1}}\left|
\frac{\sin\big((2N+1)\pi (x-y)\big)}{\sin(\pi( x-y))}\right|\,dy\\&&\qquad= M\sup_{x\in[0,1)} \int_{0}^{1}\left|
\frac{\sin\big((2N+1)\pi (x-y)\big)}{\sin(\pi( x-y))}\right|\,dy\\&&\qquad= 2M \int_{0}^{1/2}\left|
\frac{\sin\big((2N+1)\pi \theta\big)}{\sin(\pi\theta)}\right|\,d\theta.
\end{eqnarray*}
As a result (see Exercise~\ref{K-3PIO.eDhnZmasKvcTYhFA.1}) we conclude that
\begin{equation}\label{EDG:GWHEQmpec.2.DE47.oo23}
\sup_{x\in[0,1)}\left|\sum_{{k\in\Z}\atop{|k|\le N}}\sum_{j=0}^\ell \int_{x_j}^{x_{j+1}} f'(y)\,e^{2\pi ik( x-y)}\,dy\right|\le C\ln N,
\end{equation}
for some~$C>0$.

Employing again Lemma~\ref{KASMqwdfed123erDKLI}, we also see that, for all~$j_0\in\{1,\dots,\ell\}$ and~$x\in[0,1)\setminus\{x_{j_0}\}$,
\begin{equation}\label{Mefscy-CYEDG:GWHEQmpec.2.DE47}
\left|\sum_{{k\in\Z}\atop{|k|\le N}} \sigma_{j_0}\,e^{2\pi ik( x-x_{j_0})}\right|=
\left|\frac{\sigma_{j_0}\,\sin\big((2N+1)\pi ( x-x_{j_0})\big)}{\sin\big(\pi( x-x_{j_0})\big)}\right|\le
\frac{C}{\big|\sin\big(\pi( x-x_{j_0})\big)\big|},
\end{equation}
up to renaming~$C>0$.

From this, \eqref{EDG:GWHEQmpec.2.DE47}, and~\eqref{EDG:GWHEQmpec.2.DE47.oo23}, we deduce that, if~$j_0\in\{1,\dots,\ell\}$,
\begin{equation*}\begin{split}
\big|E_N(x_{j_0})-\sigma_{j_0}\big|&\le\frac{C}{N}\sum_{{1\le j\le\ell}\atop{j\ne j_0}}\frac1{\big|\sin\big(\pi( x-x_{j_0})\big)\big|}+\frac{C\ln N}{N},\end{split}
\end{equation*}
up to renaming~$C$, from which the desired result in~\eqref{EDG:GWHEQmpec.2.DE3} follows.

Furthermore, by~\eqref{EDG:GWHEQmpec.2.DE47}, \eqref{EDG:GWHEQmpec.2.DE47.oo23}, and~\eqref{Mefscy-CYEDG:GWHEQmpec.2.DE47}, we also find that if~$x\in[0,1)\setminus\{x_1,\dots,x_\ell\}$, then
\begin{eqnarray*}
\big|E_N(x)\big|\le\frac{|\sigma_{j_0}|}{2N}+\frac{C}{N\,\big|\sin\big(\pi( x-x_{j_0})\big)\big|}+\frac{C\ln N}{N},
\end{eqnarray*} from which we obtain~\eqref{EDG:GWHEQmpec.2.DE1}.

\paragraph{Solution to Exercise~\ref{SASCaCmnmsPcplstqE}.}
\begin{equation}\label{SASCaCmnmsPcplstq}
E_N(x)=\frac{\pi}{N}\sum_{k=1}^N k\big(  b_k \cos(2 \pi k x)-
a_k \sin(2 \pi k x) 
\big).\end{equation}
By means of~\eqref{fasv} and~\eqref{jasmx23er}, we can write~\eqref{EnaENSIADCac}
in the form
\begin{eqnarray*}
&& \frac{i\pi}{N}\sum_{{k\in\Z}\atop{0<|k|\le N}} k\widehat f_k\,e^{2\pi i kx}
=\frac{i\pi}{N}\sum_{k=1}^N k\Big(\widehat f_k\,e^{2\pi i kx}-
\widehat f_{-k}\,e^{-2\pi i kx}
\Big)\\&&\qquad
=\frac{i\pi}{N}\sum_{k=1}^N k\Big(\widehat f_k\,e^{2\pi i kx}-
\overline{\widehat f_{k}\,e^{2\pi i kx}}
\Big)
=-\frac{2\pi}{N}\sum_{k=1}^N k\, \Im\Big(\widehat f_k\,e^{2\pi i kx}
\Big)\\&&\qquad
=-\frac{\pi}{N}\sum_{k=1}^N k\,\Im\Big(\big(a_k-ib_k\big)
\big(\cos(2\pi kx)+i\sin(2\pi kx)\big)\Big),\end{eqnarray*}
leading to~\eqref{SASCaCmnmsPcplstq}.

\begin{figure}[h]
\includegraphics[height=2.99cm]{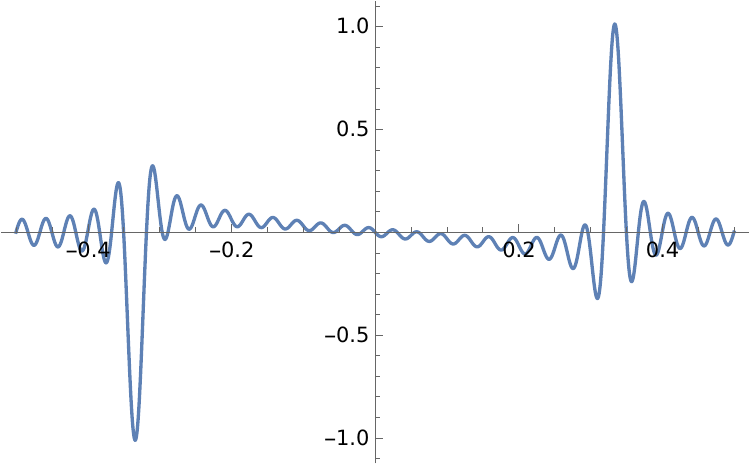}$\,\;\quad$\includegraphics[height=2.99cm]{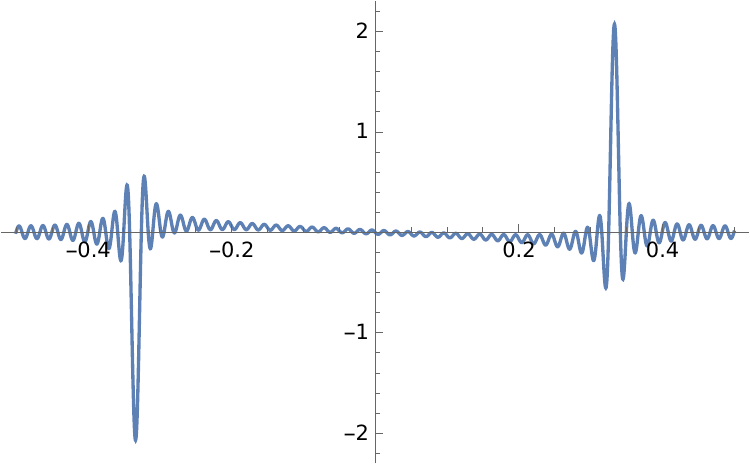} $\,\;\quad$\includegraphics[height=2.99cm]{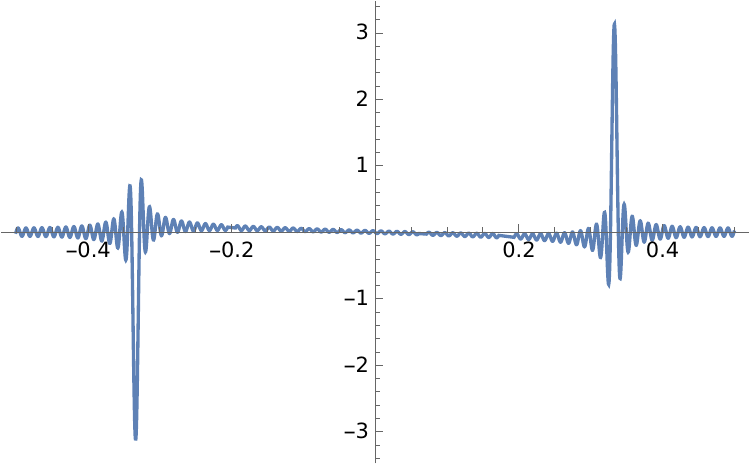}
\centering
\caption{Diagrams of the function in~\eqref{SASCaCmnmVEQ}
for~$M\in\{10,\,20,\,30\}$.}\label{RnaBbamnlocET12P.91o2eujrfne.023oritk}
\end{figure}

\paragraph{Solution to Exercise~\ref{SASCaCmnmsPcplstqE-vis}.}
We make use of Exercises~\ref{NU90i3orjf:12oeihfnvZMAS} and~\ref{SASCaCmnmsPcplstqE}.
Namely, off the back of~\eqref{02938urfjiikmsiiklo0PSnm} and~\eqref{SASCaCmnmsPcplstq}, choosing~$N:=3M$ for simplicity, we see that
\begin{equation}\label{SASCaCmnmVEQ}\begin{split}
-\frac{3M \,E_{3M}(x)}\pi&=
\sum_{j=1}^M \frac{2}{9\pi^2 j} \sin(6 \pi j x) \\&\qquad
+
\sum_{j=0}^{M-1}
\frac{(2\pi^2(3j+1)^2-9)\sqrt{3}-6\pi (3j+1)}{18\pi^3 (3j+1)^2}
\sin(2 \pi (3j+1) x)\\ &\qquad
-
\sum_{j=0}^{M-1}
\frac{(2\pi^2(3j+2)^2-9)\sqrt{3}+6\pi (3j+2)}{18\pi^3 (3j+2)^2}
\sin(2 \pi (3j+2) x),
\end{split}\end{equation}
whose graph is sketeched in Figure~\ref{RnaBbamnlocET12P.91o2eujrfne.023oritk} (notice the vanishing
away from~$x=\pm\frac13$ and the rising singularities at those points, in agreement with~\eqref{EDG:GWHEQmpec.2.DE1} and~\eqref{EDG:GWHEQmpec.2.DE3};
to be compared with Figure~\ref{C7n7okjntghbrdGHJKFRO903.lk90fv0-12920045-23-5283012P}).

\section{Solutions to selected exercises of Section~\ref{TWHOLIKSMDRA902}}

\paragraph{Solution to Exercise~\ref{BESSEL-FC-EX1}.} For all~$k\in\N$, $a>0$, and~$N$, $M\in\N$,
$$ \sup_{x\in[-a,a]} \sum_{j=N}^{M}\frac{| x|^{2j+k}}{2^{2j+k}\,j!\,(j+k)!}\le\sum_{j=N}^{M}\frac{a^{2j+k}}{2^{2j+k}\,j!\,(j+k)!},$$
showing that the series defining~$J_k$ converges locally uniformly.

Actually, a similar argument applies to the derivatives of~$J_k$, therefore, for all~$k\in\N$,
\begin{eqnarray*}
&& x^2 J''_k(x)+x\,J'_k(x)+(x^2-k^2)\,J_k(x)\\&&\qquad=
\sum_{j=0}^{+\infty}\frac{(-1)^j \,(2j+k)(2j+k-1)x^{2j+k}}{2^{2j+k}\,j!\,(j+k)!}
+\sum_{j=0}^{+\infty}\frac{(-1)^j\,(2j+k) x^{2j+k}}{2^{2j+k}\,j!\,(j+k)!}\\&&\qquad\qquad\qquad
+(x^2-k^2)
\sum_{j=0}^{+\infty}\frac{(-1)^j x^{2j+k}}{2^{2j+k}\,j!\,(j+k)!}\\&&\qquad=
\sum_{j=0}^{+\infty}\frac{(-1)^j \,x^{2j+k}}{2^{2j+k}\,j!\,(j+k)!}\Big(
(2j+k)(2j+k-1)+(2j+k)-k^2\Big)\\&&\qquad\qquad\qquad+\sum_{j=0}^{+\infty}\frac{(-1)^j\,x^{2(j+1)+k}}{2^{2j+k}\,j!\,(j+k)!}
\\&&\qquad=
4\sum_{j=1}^{+\infty}\frac{(-1)^j \,j\,x^{2j+k}}{2^{2j+k}\,j!\,(j+k-1)!}+\sum_{j=0}^{+\infty}\frac{(-1)^j\,x^{2(j+1)+k}}{2^{2j+k}\,j!\,(j+k)!}\\
\\&&\qquad=
-4\sum_{J=0}^{+\infty}\frac{(-1)^J \,x^{2(J+1)+k}}{2^{2(J+1)+k}\,J!\,(J+k)!}+\sum_{j=0}^{+\infty}\frac{(-1)^j\,x^{2(j+1)+k}}{2^{2j+k}\,j!\,(j+k)!}\\&&\qquad=0,
\end{eqnarray*}
which establishes~\eqref{BEDIEQG} when~$k\in\N$.

Now, if~$k\in\Z\cap(-\infty,-1]$, we use the previous observation and~\eqref{BEDIEQGnega} to see that
$$ x^2 J''_k(x)+x\,J'_k(x)+(x^2-k^2)\,J_k(x)=(-1)^k\big(x^2 J''_{-k}(x)+x\,J'_{-k}(x)+(x^2-k^2)\,J_{-k}(x)\big)=0,$$
which completes the proof of~\eqref{BEDIEQG}.

\paragraph{Solution to Exercise~\ref{BESSEL-FC-EX2fac}.} For~$k\in\N$, the desired claim follows from~\eqref{BEDIEQG.010}.
This and~\eqref{BEDIEQGnega} also give the desired result when~$k$ is a negative integer.

\paragraph{Solution to Exercise~\ref{BESSEL-FC-EX2}.} In light of~\eqref{BEDIEQG.010} and~\eqref{BEDIEQGnega} we have that, when~$m\in\N\cap[1,+\infty)$,
\begin{equation}\label{BEDIEQG.010.0122395}
J_{-m}(x):=(-1)^{m}J_{m}(x)=\sum_{j=0}^{+\infty}\frac{(-1)^{j+m} x^{2j+m}}{2^{2j+m}\,j!\,(j+m)!}.
\end{equation}

Also, since, for all~$\zeta\in\C$,
$$ e^\zeta=\sum_{h=0}^{+\infty} \frac{\zeta^h}{h!},$$
with locally uniform convergence of the above series, we have that
\begin{eqnarray*}&&\exp\left(\frac{x}2\left(z-\frac1z\right)\right)=
e^{\frac{xz}2}\,e^{-\frac{x}{2z}}=
\sum_{h,j=0}^{+\infty} \frac{(-1)^j\,x^{h+j}z^{h-j}}{2^{h+j}\,j!\,h!}\\&&\qquad=
\sum_{{j\ge0}\atop{h\ge j}} \frac{(-1)^j\,x^{h+j}z^{h-j}}{2^{h+j}\,j!\,h!}+
\sum_{{j\ge0}\atop{0\le h\le j-1}} \frac{(-1)^j\,x^{h+j}z^{h-j}}{2^{h+j}\,j!\,h!}
.\end{eqnarray*}
Hence, the index substitutions~$k:=h-j$ and~$m:=j-h$, together with~\eqref{BEDIEQG.010}, lead to
\begin{equation}\label{BEDIEQG.010.0122395-924}\begin{split}\exp\left(\frac{x}2\left(z-\frac1z\right)\right)&=
\sum_{{j\ge0}\atop{k\ge 0}} \frac{(-1)^j\,x^{2j+k}z^{k}}{2^{2j+k}\,j!\,(j+k)!}+
\sum_{{j\ge0}\atop{1\le m\le j}} \frac{(-1)^j\,x^{2j-m}z^{-m}}{2^{2j-m}\,j!\,(j-m)!}\\&=\sum_{k=0}^{+\infty}J_k(x)\,z^k+
\sum_{{m\ge1}\atop{j\ge m}} \frac{(-1)^j\,x^{2j-m}z^{-m}}{2^{2j-m}\,j!\,(j-m)!}.\end{split}\end{equation}

Moreover, comparing with~\eqref{BEDIEQG.010.0122395}, for all~$m\in\N\cap[1,+\infty)$, the index substitution~$\ell:=j-m$ gives that
\begin{eqnarray*}
\sum_{j=m}^{+\infty} \frac{(-1)^j\,x^{2j-m}}{2^{2j-m}\,j!\,(j-m)!}=\sum_{\ell=0}^{+\infty} \frac{(-1)^{\ell+m}\,x^{2\ell+m}}{2^{2\ell+m}\,(\ell+m)!\,\ell!}=J_{-m}(x).
\end{eqnarray*}
We plug this information into~\eqref{BEDIEQG.010.0122395-924} and we conclude that
$$ \exp\left(\frac{x}2\left(z-\frac1z\right)\right)=\sum_{k=0}^{+\infty}J_k(x)\,z^k+\sum_{m=1}^{+\infty}J_{-m}(x)\,z^{-m},$$
which yields~\eqref{GENEBESS}.

\begin{figure}[h]
\includegraphics[height=2.8cm]{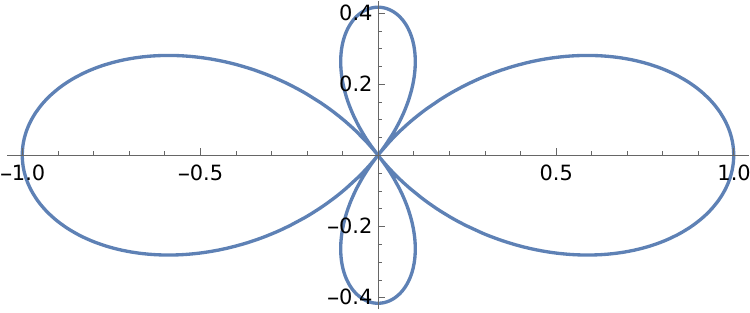}$\,\;\qquad\quad$\includegraphics[height=2.8cm]{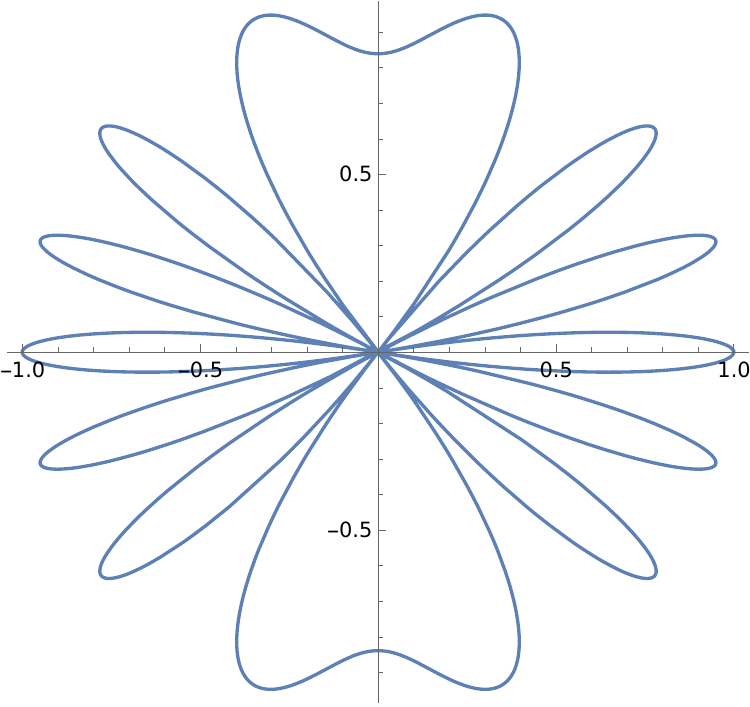} $\,\;\qquad\quad$\includegraphics[height=2.8cm]{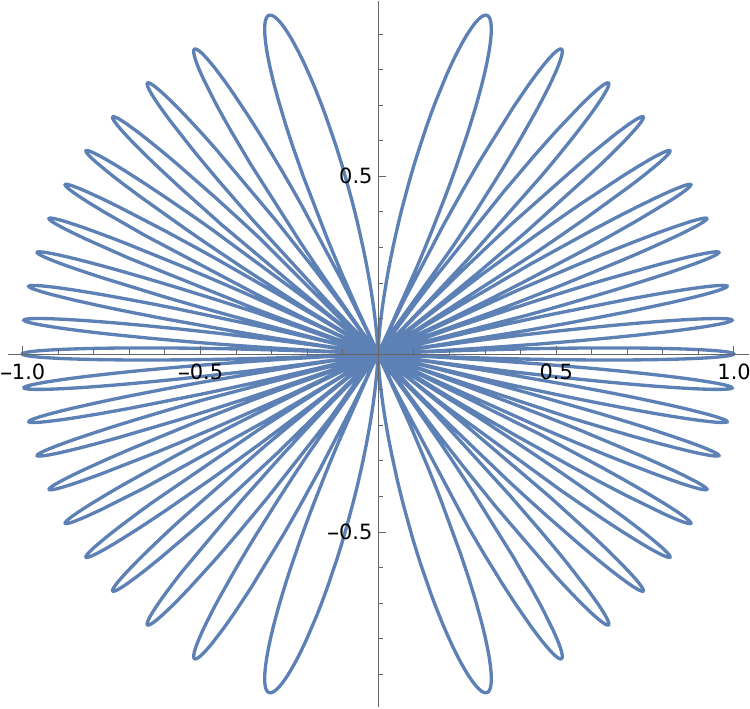}
\centering
\caption{Diagrams of the curves parametrised in polar coordinates by~$r=\cos( \beta\sin\theta)$, with~$\beta\in\{2,\,10,\,33\}$.}\label{RET12P.91o2eujrfne.023oritk}
\end{figure}

\paragraph{Solution to Exercise~\ref{BESSEL-FC-EX3}.} On account of Exercise~\ref{BESSEL-FC-EX2}, using the notation~$z:=e^{i\theta}$, and recalling~\eqref{BEDIEQGnega}, we see that
\begin{eqnarray*}&&
\cos\left(x\sin\theta\right)+i\sin\left(x\sin\theta\right)
=\exp\left(ix\sin\theta\right)=
\exp\left(\frac{x}2\left(e^{i\theta}-e^{-i\theta}\right)\right)\\&&\qquad=\sum_{{k\in\Z}}J_k(x)\,e^{ik\theta}=
\sum_{{k\in\Z}}J_k(x)\,\cos(k\theta)+i\sum_{{k\in\Z}}J_k(x)\,\sin(k\theta)\\&&\qquad=
J_0(x)+\sum_{k=1}^{+\infty} \big(J_k(x)+J_{-k}(x)\big)\,\cos(k\theta)+i\sum_{k=1}^{+\infty}\big(J_k(x)-J_{-k}(x)\big)\,\sin(k\theta)\\&&\qquad=J_0(x)+\sum_{k=1}^{+\infty} \big(J_k(x)+(-1)^kJ_{k}(x)\big)\,\cos(k\theta)+i\sum_{k=1}^{+\infty}\big(J_k(x)-(-1)^kJ_k(x)\big)\,\sin(k\theta)\\&&\qquad=J_0(x)+\sum_{m=1}^{+\infty} 2J_{2m}(x)\,\cos(2m\theta)+i\sum_{m=0}^{+\infty}2J_{2m+1}(x)\,\sin((2m+1)\theta).
\end{eqnarray*}
The desired result follows by taking the real and imaginary parts.

\begin{figure}[h]
\includegraphics[height=2.8cm]{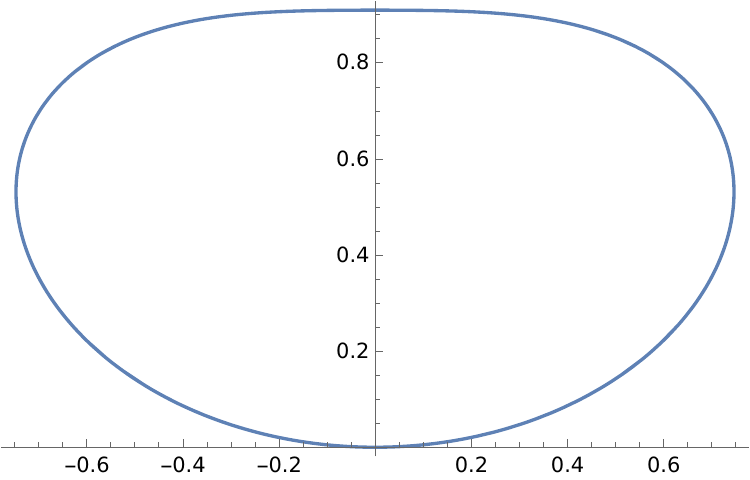}$\,\;\qquad\quad$\includegraphics[height=2.8cm]{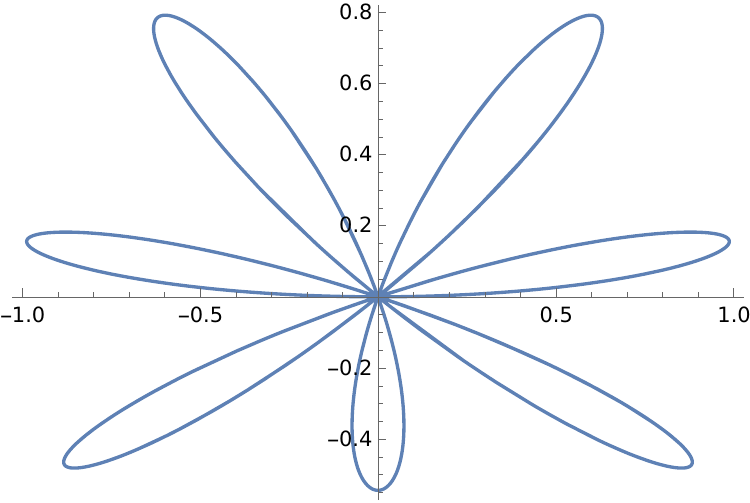}
$\,\;\qquad\quad$\includegraphics[height=2.8cm]{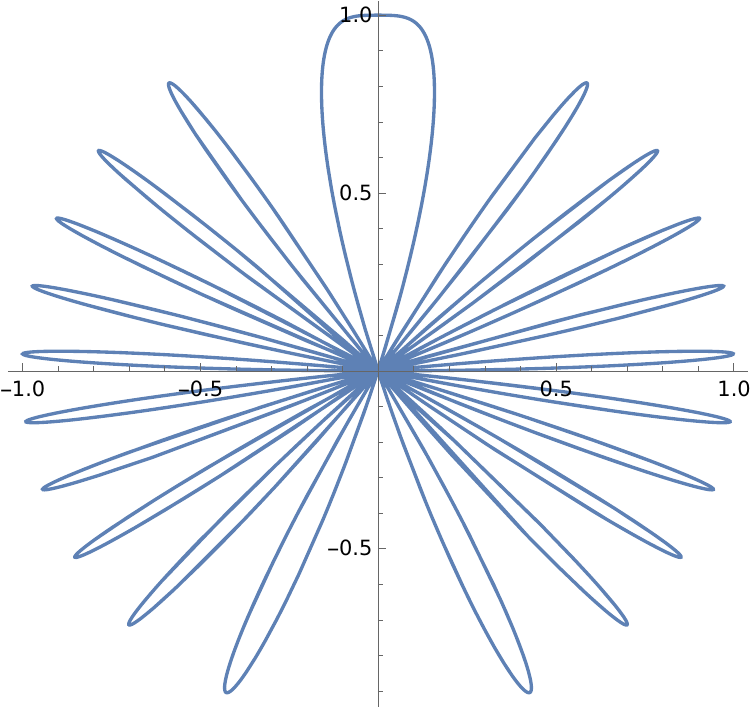}
\centering
\caption{Diagrams of the curves parametrised in polar coordinates by~$r=\sin( \beta\sin\theta)$, with~$\beta\in\{2,\,10,\,33\}$.}\label{RET12P.91o2eujrfne.023oritks}
\end{figure}

On a possibly unrelated topic, let us mention that the functions showing up in this exercise can be used to produce cute pictures, see Figures~\ref{RET12P.91o2eujrfne.023oritk} and~\ref{RET12P.91o2eujrfne.023oritks}.

\paragraph{Solution to Exercise~\ref{BESSEL-FC-EX4}.} After Exercise~\ref{BESSEL-FC-EX3},
used here with~$\beta:=\epsilon$ and~$\theta:=2\pi\omega' t$,
we have that
\begin{eqnarray*}&&
\sin\big(2\pi\omega t+\epsilon \sin(2\pi\omega' t)\big)\\&&\qquad=
\sin\big(2\pi\omega t\big)\,\cos\big(\epsilon \sin(2\pi\omega' t)\big)+\cos\big(2\pi\omega t\big)\,\sin\big(\epsilon \sin(2\pi\omega' t)\big)\\&&\qquad=\sin\big(2\pi\omega t\big)\,\left(J_0(\epsilon) + 2\sum_{m=1}^{+\infty}J_{2m}(\epsilon)\,\cos(4m\pi\omega' t)\right)\\&&\qquad\qquad\qquad+2\cos\big(2\pi\omega t\big)\,\sum_{m=0}^{+\infty}J_{2m+1}(\epsilon)\,\sin(2(2m+1)\pi\omega' t)
\\&&\qquad=J_0(\epsilon) \,\sin\big(2\pi\omega t\big)+ 2\sum_{m=1}^{+\infty}J_{2m}(\epsilon)\,\sin\big(2\pi\omega t\big)\,\cos(4m\pi\omega' t) \\&&\qquad\qquad\qquad+2\sum_{m=0}^{+\infty}J_{2m+1}(\epsilon)\,\sin(2(2m+1)\pi\omega' t)\,\cos\big(2\pi\omega t\big).\end{eqnarray*}

Hence, by Prosthaphaeresis,
\begin{eqnarray*}&&
\sin\big(2\pi\omega t+\epsilon \sin(2\pi\omega' t)\big)\\&&\qquad=J_0(\epsilon) \,\sin\big(2\pi\omega t\big)+ \sum_{m=1}^{+\infty}J_{2m}(\epsilon)\,
\Big(\sin\big(2\pi(\omega t+2m\omega' t)\big)+\sin\big(2\pi(\omega t-2m\omega' t)\big)\Big) \\&&\qquad\qquad+\sum_{m=0}^{+\infty}J_{2m+1}(\epsilon)\,\Big(\sin\big(2\pi ((2m+1)\omega' t+\omega t)\big)+\sin\big(2\pi ((2m+1)\omega' t-\omega t)\big)\Big)
\\&&
\qquad=J_0(\epsilon) \,\sin\big(2\pi\omega t\big)+ 
\sum_{k=1}^{+\infty} J_k(\epsilon)\,\sin\big(2\pi(\omega+k\omega')t\big)+\sum_{k=1}^{+\infty} (-1)^kJ_k(\epsilon)\,\sin\big(2\pi(\omega-k\omega')t\big).
\end{eqnarray*}
The desired result now follows from~\eqref{BEDIEQGnega}.

\paragraph{Solution to Exercise~\ref{BESSEL-FC-EX4-DIFFB.imp}.}
We employ Exercise~\ref{BESSEL-FC-EX3}
(and the uniform convergence of the series involved in order to swap the summation and integral signs) and we see that
\begin{equation}\label{MANSDODKJGrGBskkipalISNDGHoyo1} \frac1\pi\int_0^\pi \cos(x\sin\theta)\,d\theta=
J_0(x) + \frac2\pi\sum_{m=1}^{+\infty}J_{2m}(x)\,\int_0^\pi\cos(2m\theta)\,d\theta=J_0(x).\end{equation}

We also remark that, using the substitution~$\alpha:=\pi-\theta$,
\begin{eqnarray*}&&
\int_0^\pi \cos(x\sin\theta)\,d\theta=\int_0^{\pi/2} \cos(x\sin\theta)\,d\theta+
\int_{\pi/2}^\pi \cos(x\sin\theta)\,d\theta\\&&\qquad=\int_0^{\pi/2} \cos(x\sin\theta)\,d\theta+
\int_0^{\pi/2} \cos(x\sin(\pi-\alpha))\,d\alpha=2\int_0^{\pi/2} \cos(x\sin\theta)\,d\theta.
\end{eqnarray*}
Hence, comparing with~\eqref{MANSDODKJGrGBskkipalISNDGHoyo1} and substituting for~$\sin\theta=:t$,
\begin{equation*} J_0(x)=\frac2\pi\int_0^1\frac{\cos(xt)}{\sqrt{1-t^2}}\,dt.\end{equation*}
The desired result is thereby established.

\paragraph{Solution to Exercise~\ref{BESSEL-FC-EX4-DIFFBmahUAHSNd01}.}
After Exercise~\ref{BESSEL-FC-EX4-DIFFB.imp} we see that, for all~$m\in\N\setminus\{0\}$,
\begin{eqnarray*}J_0(m\pi)=\frac2\pi\int_0^1\frac{\cos(m\pi t)}{\sqrt{1-t^2}}\,dt
.\end{eqnarray*}
Substituting for~$\tau:=m t$, we find that
\begin{equation}\label{GCgCpaqdmdINcpFK01}\begin{split}\frac{\pi J_0(m\pi)}2=
\int_0^{m}\frac{\cos(\pi \tau)}{\sqrt{m^2-\tau^2}}\,d\tau=\sum_{j=0}^{m-1}
\int_{j}^{j+1}\frac{\cos(\pi \tau)}{\sqrt{m^2-\tau^2}}\,d\tau
.\end{split}\end{equation}

Also, setting~$\sigma:=\tau-j$, we see that
\begin{eqnarray*}
&&\int_{j}^{j+1}\frac{\cos(\pi \tau)}{\sqrt{m^2-\tau^2}}\,d\tau=
\int_{0}^{1}\frac{\cos(\pi \sigma+\pi j)}{\sqrt{m^2-(\sigma+j)^2}}\,d\sigma=\int_{0}^{1}\frac{(-1)^j\,\cos(\pi \sigma)}{\sqrt{m^2-(\sigma+j)^2}}\,d\sigma.
\end{eqnarray*}
Hence, we define
$$ \Phi(\sigma):=\sum_{j=0}^{m-1}\frac{(-1)^j}{\sqrt{m^2-(\sigma+j)^2}}$$
and we gather that
\begin{equation}\label{Gc3245cfI-2j1-93ijtg} \sum_{j=0}^{m-1}\int_{j}^{j+1}\frac{\cos(\pi \tau)}{\sqrt{m^2-\tau^2}}\,d\tau=
\int_{0}^{1} \cos(\pi \sigma)\,\Phi(\sigma)\,d\sigma.\end{equation}

Now, we use an integration by parts and the fact that~$\sin0=\sin\pi=0$ to see that
\begin{eqnarray*}&&
\int_{0}^{1} \cos(\pi \sigma)\,\Phi(\sigma)\,d\sigma=\frac1\pi
\int_{0}^{1}\frac{d}{d\sigma}\big( \sin(\pi \sigma)\big)\,\Phi(\sigma)\,d\sigma\\&&\qquad=-
\frac1\pi
\int_{0}^{1}\sin(\pi \sigma)\,\Phi'(\sigma)\,d\sigma=\sum_{j=0}^{m-1}(-1)^{j+1}
\int_{0}^{1}\frac{\sin(\pi \sigma)\, (j + \sigma)}{\pi\,(m^2 - (j + \sigma)^2)^{\frac32}}
\,d\sigma\\&&\qquad=\sum_{j=0}^{m-1}(-1)^{j+1} c_j(m),
\end{eqnarray*}
where, for all~$m\in\N$ and~$j\in\{0,\dots,m-1\}$, we set
$$ c_j(m):=\int_{0}^{1}\frac{\sin(\pi \sigma)\, (j + \sigma)}{\pi\,(m^2 - (j + \sigma)^2)^{\frac32}}
\,d\sigma.$$

For this reason, we write~\eqref{Gc3245cfI-2j1-93ijtg} as
\[ \sum_{j=0}^{m-1}\int_{j}^{j+1}\frac{\cos(\pi \tau)}{\sqrt{m^2-\tau^2}}\,d\tau=\sum_{j=0}^{m-1}(-1)^{j+1} c_j(m)\]
and then, by~\eqref{GCgCpaqdmdINcpFK01},
\begin{equation}\label{GCgCpaqdmdINcpFK01sp2}\frac{\pi J_0(m\pi)}2=\sum_{j=0}^{m-1}(-1)^j c_j(m).\end{equation}

Now we remark that, for all~$\sigma\in(0,1)$ and~$s\in(0,m-1)$,
\begin{equation}\label{GCgCpaqdmdINcpFK01sp} \frac{d}{ds} \left(\frac{s + \sigma}{(m^2 - (s + \sigma)^2)^{\frac32}}\right)=\frac{m^2 + 2 (s+ \sigma)^2}{(m^2 - (s + \sigma)^2)^{\frac52}}>0.\end{equation}

We also stress that~$c_0(m)>0$ and that, when~$j\in\{0,\dots,m-2\}$,
\begin{eqnarray*}
c_{j+1}(m)-c_j(m)&=&\int_{0}^{1}\frac{\sin(\pi \sigma)}{\pi}\left(\frac{j +1+ \sigma}{(m^2 - (j +1+ \sigma)^2)^{\frac32}}-\frac{j + \sigma}{(m^2 - (j + \sigma)^2)^{\frac32}}\right)
\,d\sigma>0,
\end{eqnarray*}
thanks to~\eqref{GCgCpaqdmdINcpFK01sp}.

From this and~\eqref{GCgCpaqdmdINcpFK01sp2} we infer that~$J_0(m\pi)>0$ when~$m$ is even and~$J_0(m\pi)<0$ when~$m$ is odd.

\begin{figure}[h]
\includegraphics[height=3.94cm]{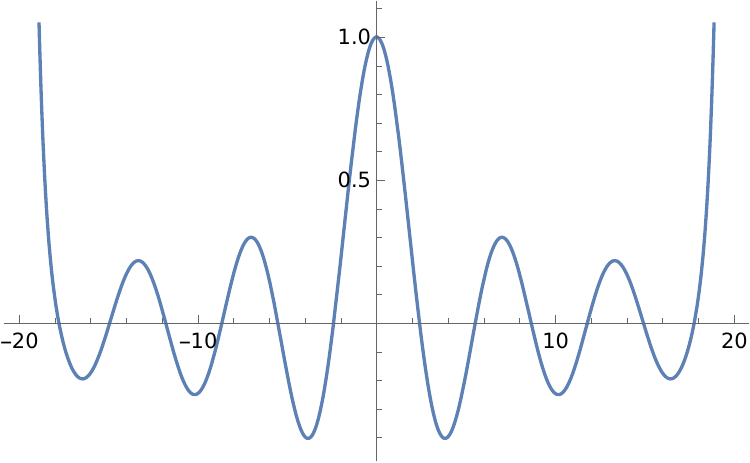}$\,\;\quad\quad$\includegraphics[height=3.94cm]{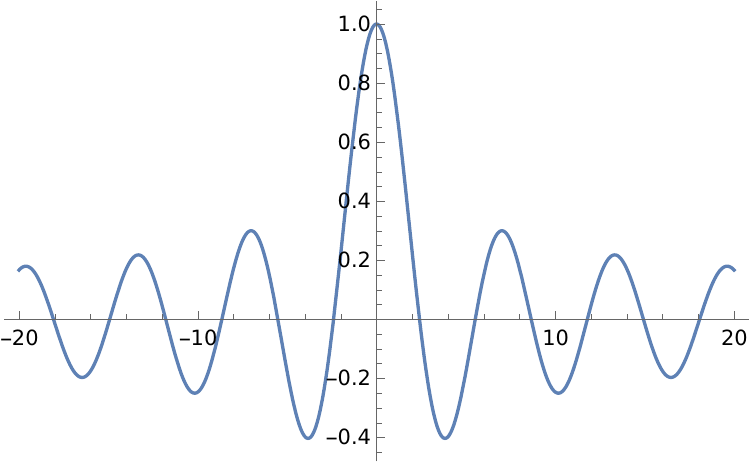}
\centering
\caption{Left: Plot of~$\displaystyle\sum_{j=0}^{22}\frac{(-1)^j x^{2j}}{2^{2j}\,(j!)^2}$.
Right: plot of~$J_0$.}\label{4DET24OikdnsEUxXqewr0}
\end{figure}

See Figure~\ref{4DET24OikdnsEUxXqewr0} for a sketch of~$J_0$ and of an approximation of it.

\paragraph{Solution to Exercise~\ref{BESSEL-FC-EX4-DIFFB}.}
The fact that~$J_0$ possesses infinitely many zeros follows from Exercise~\ref{BESSEL-FC-EX4-DIFFBmahUAHSNd01}, so we only need to show that these zeros are discrete.

We give two proofs of this. First proof:
since~$J_0$ is real analytic, by the uniform convergence of the series in Exercise~\ref{BESSEL-FC-EX1},
its zeros cannot accumulate on each other (see e.g.~\cite[Corollaries~1.2.4 and~1.2.7]{MR1916029}).

Second proof: by~\eqref{BEDIEQG.010} we know that
\begin{equation}\label{BEDIEQG.0100k2Zewdf.1}
J_0(0)=1.\end{equation}
Now, by contradiction, we suppose that there is a sequence~$\{x_\ell\}_{\ell\in\N}$ such that~$J_0(x_\ell)=0$ for each~$\ell\in\N$ and~$x_\ell$ converges to some~$x_\star$ as~$\ell\to+\infty$.

Hence, we have that~$J_0(x_\star)=0$ and therefore
\begin{equation}\label{BEDIEQG.0100k2Zewdf.2}
x_\star\ne0,\end{equation} owing to~\eqref{BEDIEQG.0100k2Zewdf.1}.

Moreover, by Rolle's Theorem, we can find~$y_\ell$ on the segment joining~$x_\ell$ to~$x_{\ell+1}$ such that~$J_0'(y_\ell)=0$. We observe that~$y_\ell$ also converges to~$x_\star$ as~$\ell\to+\infty$, thus~$J'_0(x_\star)=0$.

Since~\eqref{BEDIEQG} is a Cauchy problem in the vicinity of~$x_\star$, due to~\eqref{BEDIEQG.0100k2Zewdf.2},
the uniqueness of the solution of ordinary differential equations gives that~$J_0$ must vanish identically, in contradiction with~\eqref{BEDIEQG.0100k2Zewdf.1}.

See~\cite[Chapter~V]{MR193286} and the references therein for more information on the zeros of Bessel functions.

\paragraph{Solution to Exercise~\ref{BESSEL-FC-EX4-DIFFB-HYPB}.} We stress that the existence of~$\lambda$ is guaranteed by Exercise~\ref{BESSEL-FC-EX4-DIFFB}.

Also
$$ U(1,t)=J_0(\lambda )\,\cos(\sqrt{c}\,\lambda t)=0,$$
by definition of~$\lambda$, and
$$ \partial_t U(r,t)=-\sqrt{c}\,\lambda\,J_0(\lambda r)\,\sin(\sqrt{c}\,\lambda t),$$
yielding that~$\partial_t U(r,0)=0$.

Besides, one can rewrite the equation in~\eqref{BESSEL-FC-EX4-DIFFB-HYPB.e}
in polar coordinates as
\begin{equation}\label{BESSEL-FC-EX4-DIFFB-HYPB.ef}
\partial^2_t U=c\partial^2_r U+\frac{c\partial_r U}r.\end{equation}

In this spirit, recalling Bessel's Differential Equation~\eqref{BEDIEQG} with~$k=0$, we find that
\begin{eqnarray*}&&
c\partial^2_r U(r,t)+\frac{c\partial_r U(r,t)}r-\partial^2_t U(r,t)\\&&\qquad=
c\lambda^2 J_0''(\lambda r)\,\cos(\sqrt{c}\,\lambda t)+\frac{c\lambda J_0'(\lambda r)\,\cos(\sqrt{c}\,\lambda t)}{r}
-c\lambda^2 J_0(\lambda r)\,\cos(\sqrt{c}\,\lambda t)
\\&&\qquad=\frac{c\,\cos(\sqrt{c}\,\lambda t)}{r^2}\Big(
\lambda^2 r^2\,J_0''(\lambda r)+\lambda r\, J_0'(\lambda r)
- (\lambda r)^2\,J_0(\lambda r)\Big)\\&&\qquad=0.
\end{eqnarray*}
This gives that~\eqref{BESSEL-FC-EX4-DIFFB-HYPB.ef} is satisfied for all~$r\in(0,1)$
and therefore~\eqref{BESSEL-FC-EX4-DIFFB-HYPB.e} is satisfied in the unit disk, possibly outside the origin.

Hence, to complete the requested task, we need to show that~\eqref{BESSEL-FC-EX4-DIFFB-HYPB.e} is satisfied at the origin too. For this, we utilise~\eqref{BEDIEQG.010} to see that, near the origin,
$$ J_0(\xi)=1-\frac{ \xi^{2}}{4}+O(\xi^4)$$
and accordingly~\eqref{BESSEL-FC-EX4-DIFFB-HYPB.ui} can be written in Cartesian coordinates near the origin as
$$ u(x,y,t)=J_0\left(\lambda \sqrt{x^2+y^2}\right)\,\cos(\sqrt{c}\,\lambda t)
=\left( 1-\frac{ \lambda^{2}(x^2+y^2)}{4}+O\big((x^2+y^2)^2\big)\right)\,\cos(\sqrt{c}\,\lambda t).$$
As a result, $u$ is of class~$C^2$ in the vicinity of the origin, hence we can pass to the limit~\eqref{BESSEL-FC-EX4-DIFFB-HYPB.e} from values outside the origin and obtain that it is satisfied at the origin too, as desired.

We notice that Exercise~\ref{BESSEL-FC-EX4-DIFFB-HYPB} can be seen as a suitable ``hyperbolic'' version
of the ``elliptic'' problem studied in Section~\ref{POISSONKERN-ex4}.

See e.g.~\cite[Chapter~X]{MR193286} for more information
about the use of Bessel Functions to describe the displacement of a membrane and related problems.

\section{Solutions to selected exercises of Section~\ref{LIMOTORSE}}

\paragraph{Solution to Exercise~\ref{PelapeammilonDF}.}
We revisit the proof of Theorem~\ref{ERGOB}. Namely, let~$\epsilon>0$ and (recalling~\eqref{ERGOB.eqRoscjh})
pick a trigonometric polynomial~$S_\epsilon$ such that~$\|f-S_{\epsilon}\|_{L^\infty(\R^n)}\le2\epsilon$.

Then,
\begin{equation}\label{GLDShatqwPelapeammilonDF-1}
\left|\int_0^1f(x)\,dx-\int_0^1S_\epsilon(x)\,dx\right|\le2\epsilon\end{equation}
and
\begin{equation}\label{GLDShatqwPelapeammilonDF-2}
\left|\frac1K\sum_{k=0}^{K-1} f(\gamma k)-\frac1K\sum_{k=0}^{K-1} S_\epsilon(\gamma k)\right|\le2\epsilon.\end{equation}

Moreover, writing
$$ S_{\epsilon}(x)=\sum_{{j\in\Z}\atop{|j|\le N_\epsilon}}c_{\epsilon,j}\,e^{2\pi ijx},$$
for suitable~$c_{\epsilon,j}$ and~$N_\epsilon$, we see that
\begin{eqnarray*}
&&\sum_{k=0}^{K-1} S_\epsilon(\gamma k)=
\sum_{{0\le k\le K-1}\atop{{j\in\Z}\atop{|j|\le N_\epsilon}}}c_{\epsilon,j}\,e^{2\pi ij\gamma k}=Kc_{\epsilon,0}+
\sum_{{{j\in\Z}\atop{1\le|j|\le N_\epsilon}}}\frac{c_{\epsilon,j}\,\big(1-e^{2\pi ij\gamma K}\big)}{1-e^{2\pi ij\gamma}},
\end{eqnarray*}
and we stress that the latter denominator does not vanish, since~$\gamma$ is irrational.

As a result,
\begin{eqnarray*}&&
\lim_{K\to+\infty}\frac1K\sum_{k=0}^{K-1} S_\epsilon(\gamma k)=c_{\epsilon,0}+\lim_{K\to+\infty}\frac1K
\sum_{{{j\in\Z}\atop{1\le|j|\le N_\epsilon}}}\frac{c_{\epsilon,j}\,\big(1-e^{2\pi ij\gamma K}\big)}{1-e^{2\pi ij\gamma}}\\&&\qquad=c_{\epsilon,0}=\int_0^1 \left(\sum_{{j\in\Z}\atop{|j|\le N_\epsilon}}c_{\epsilon,j}\,e^{2\pi ijx}\right)\,dx=\int_0^1S_\epsilon(x)\,dx.
\end{eqnarray*}
From this, \eqref{GLDShatqwPelapeammilonDF-1}, and~\eqref{GLDShatqwPelapeammilonDF-2} we gather that
$$ \limsup_{K\to+\infty}
\left|\frac1K\sum_{k=0}^{K-1} f(\gamma k)-\int_0^1f(x)\,dx\right|\le
4\epsilon+\limsup_{K\to+\infty}
\left|\frac1K\sum_{k=0}^{K-1} S_\epsilon(\gamma k)-\int_0^1S_\epsilon(x)\,dx\right|=4\epsilon.$$
The desired result now follows by sending~$\epsilon\searrow0$.

\paragraph{Solution to Exercise~\ref{PelapeammilonDF2}.} No. For instance, take~$\gamma:=\sqrt2$
and consider the set~${\mathcal{Z}}:=\{ \sqrt{2}\,k+h,$ with~$h$, $k\in\Z\}$.

We stress that~${\mathcal{Z}}$ is countable, whence of null Lebesgue measure, and therefore, if
$$ f(x):=\begin{dcases} 1&{\mbox{ if }}x\in {\mathcal{Z}},\\ 0&{\mbox{ if }}x\in\R\setminus {\mathcal{Z}},
\end{dcases}$$
then~$f$ is periodic of period~$1$ and
$$ \int_0^1f(x)\,dx=0.$$
However,
$$ \sum_{k=0}^{K-1} f(\gamma k)=\sum_{k=0}^{K-1} f(\sqrt{2}\, k)\sum_{k=0}^{K-1} 1=K.$$
As a consequence, the result in Exercise~\ref{PelapeammilonDF} does not hold true for this discontinuous function~$f$.

\section{Solutions to selected exercises of Section~\ref{PelapeammilonS}}

\paragraph{Solution to Exercise~\ref{Pelapeammilon-14.1}.}
If~$I=[a,b]$ we are precisely in the setting of Theorem~\ref{PelapeammilonT}, hence we can suppose that~$I$ is one of the intervals~$[a,b)$, $(a,b]$ or~$(a,b)$. In particular, we can assume that~$a<b$, otherwise~$I=\varnothing$ and the desired result is obvious.

Hence, we can pick~$\epsilon\in\left(0,\frac{b-a}2\right)$, to be taken as small as we wish in what follows, and we observe that~$[a+\epsilon,b-\epsilon]\subseteq I\subseteq[a,b]$ and thus~${\mathcal{N}}_{[a+\epsilon,b-\epsilon]}(K)\le {\mathcal{N}}_{I}(K)\le{\mathcal{N}}_{[a,b]}(K)$.

We thereby infer from Theorem~\ref{PelapeammilonT} that
$$ b-a-2\epsilon=
\lim_{K\to+\infty}\frac{{\mathcal{N}}_{[a+\epsilon,b-\epsilon]}(K)}{K}\le
\lim_{K\to+\infty}\frac{{\mathcal{N}}_{I}(K)}{K}\le\lim_{K\to+\infty}\frac{{\mathcal{N}}_{I}(K)}{K}=b-a,$$
and the desired result follows by sending~$\epsilon\searrow0$.

\paragraph{Solution to Exercise~\ref{Pelapeammilon4}.} If~$a=0$ and~$b=1$ the result is obvious, so we can suppose that~$[a,b]$ does not exhaust~$[0,1]$.
Let
$$ T:=\bigcup_{\ell\in\N}[\ell+a,\ell+b].$$
Let also~$\epsilon>0$ small enough such that
$$ T^+_\epsilon:=\bigcup_{p\in T} (p-\epsilon,p+\epsilon)$$
does not exhaust the whole of~$\R$.

Let also
$$ T^-_\epsilon:=\big\{p\in\R{\mbox{ s.t. }}(p-\epsilon,p+\epsilon)\in T\big\}.$$

Then, pick functions~$f_\epsilon\in C^\infty_0( T^+_\epsilon,[0,1])$ with~$f_\epsilon=1$ in~$T$
and~$g_\epsilon\in C^\infty_0( T,[0,1])$ such that~$g_\epsilon=1$ in~$T^-_\epsilon$.

In virtue of Exercise~\ref{PelapeammilonDF}, we know that
$$\limsup_{K\to+\infty}\frac1K\sum_{k=0}^{K-1} f_\epsilon(\gamma k)=\int_0^1f_\epsilon(x)\,dx\le b-a+2\epsilon.$$
and
$$\liminf_{K\to+\infty}\frac1K\sum_{k=0}^{K-1} g_\epsilon(\gamma k)=\int_0^1g_\epsilon(x)\,dx\ge b-a-2\epsilon.$$

Since~$f_\epsilon(x)=1$ when~$\{x\}\in[a,b]$, we also have that $$
\limsup_{K\to+\infty}\frac1K\sum_{k=0}^{K-1} f_\epsilon(\gamma k)\ge\limsup_{K\to+\infty}\frac{{\mathcal{N}}_{[a,b]}(K)}K.$$
Moreover, since~$g_\epsilon(x)=0$ when~$\{x\}\not\in[a,b]$,
$$
\liminf_{K\to+\infty}\frac1K\sum_{k=0}^{K-1} g_\epsilon(\gamma k)\le\liminf_{K\to+\infty}\frac{{\mathcal{N}}_{[a,b]}(K)}K.$$

All in all,
$$ b-a-2\epsilon\le\liminf_{K\to+\infty}\frac{{\mathcal{N}}_{[a,b]}(K)}K\le\limsup_{K\to+\infty}\frac{{\mathcal{N}}_{[a,b]}(K)}K\le
b-a+2\epsilon$$
and the desired result follows by sending~$\epsilon\searrow0$.

\paragraph{Solution to Exercise~\ref{Pelapeammilon5}.} Yes. Let~$L\in\N$, to be taken as large as we wish, and set
$$ [0,1]\ni x\longmapsto f_L(x):=\sum_{j=0}^{L-1} f\left(\frac{j}L \right)\chi_{\left[\frac{j}L ,\frac{j+1}L\right)}(x),$$
where the notation for the characteristic function of a set has been used (see~\eqref{CHARAFA}). 

We extend~$f_L$ as a periodic function of period~$1$. Also, we observe that, for all~$x\in[0,1)$,
taking~$j_x\in\{0,\dots,L-1\}$ such that~$x\in\left[\frac{j_x}L ,\frac{j_x+1}L\right)$, we have that
\begin{eqnarray*}
|f_L(x)-f(x)|=\left| f\left(\frac{j_x}L \right)-f(x)\right|\le\sup_{{\xi,\eta\in[0,1]}\atop{|\xi-\eta|\le 1/L}}|f(\xi)-f(\eta)|.
\end{eqnarray*}
From this and the continuity of~$f$ assumed in Exercise~\ref{PelapeammilonDF}, we infer that
\begin{equation}\label{Pelapeammilon-14.1.22}
\lim_{L\to+\infty}\sup_{x\in\R}|f_L(x)-f(x)|=0.
\end{equation}

Moreover,
\begin{eqnarray*}
\left|\frac1K\sum_{k=0}^{K-1} f_L(\gamma k)-\int_0^1f_L(x)\,dx\right|
\le\sum_{j=0}^{L-1}\left|f\left(\frac{j}L \right)\right|\;\left|\frac{1}K\sum_{k=0}^{K-1} \chi_{\left[\frac{j}L ,\frac{j+1}L\right)}(\gamma k)-
\frac{1}L\right|.
\end{eqnarray*}
Hence, by Theorem~\ref{PelapeammilonT} (actually, by the version of Theorem~\ref{PelapeammilonT}
given in Exercise~\ref{Pelapeammilon-14.1}),
\begin{eqnarray*}&&\limsup_{K\to+\infty}
\left|\frac1K\sum_{k=0}^{K-1} f_L(\gamma k)-\int_0^1f_L(x)\,dx\right|\leq\limsup_{K\to+\infty}
\sum_{j=0}^{L-1}\left|f\left(\frac{j}L \right)\right|\;\left|
\frac{{\mathcal{N}}_{{\left[\frac{j}L ,\frac{j+1}L\right)}}(K)}K-
\frac{1}L\right|=0.
\end{eqnarray*}

This and~\eqref{Pelapeammilon-14.1.22} yield that
\begin{eqnarray*}&&\limsup_{K\to+\infty}
\left|\frac1K\sum_{k=0}^{K-1} f(\gamma k)-\int_0^1f(x)\,dx\right|\\&&\qquad\leq\limsup_{L\to+\infty}
\limsup_{K\to+\infty}
\left[
\frac1K\sum_{k=0}^{K-1} |f(\gamma k)-f_L(\gamma k)|+\int_0^1|f(x)-f_L(x)|\,dx\right]=0,
\end{eqnarray*}as desired.

\paragraph{Solution to Exercise~\ref{DIpFA:0}.} We start by proving~\eqref{DIpFA:0ladj01-211}. To this end, we argue as follows.
Given~$\alpha\in\R^{n-1}$ and~$r\in\N^{n-1}$, we denote by~$ \{\alpha\cdot r\}$ the fractional part of~$\alpha\cdot r$, according to the fractional part notation introduced on page~\pageref{FRaPARTADE}.

Thus, given~$\alpha\in\R^{n-1}$ and~$M\in\N\setminus\{0\}$, we let~$ I_M$ be the set containing all the numbers of the form~$ \{\alpha\cdot r\}$, for all~$r=(r_1,\dots,r_{n-1})\in\N^{n-1}$
such that~$r_m\in[1, M]$ for all~$m\in\{1,\dots,n-1\}$.

We now suppose that the vector~$(\alpha,1)$ is rationally independent (according to the notation in~\eqref{NORESA}), i.e. for every~$k\in\Z^{n}\setminus\{0\}$ we have that~$(\alpha,1)\cdot k\ne0$. In this situation, all the elements~$\{\alpha\cdot r\}$ in the definition of~$I_M$ are necessarily different.

As a result, $I_M$ contains~$M^{n-1}$ points (since this is the cardinality of the sets made of~$r=(r_1,\dots,r_{n-1})\in\N^{n-1}$ such that~$r_m\in[1, M]$ for all~$m\in\{1,\dots,n-1\}$).

Moreover,~$I_M\subseteq(0,1)$. 

Therefore, if~$L\in\N\cap(1,M^{n-1})$ and we divide~$(0,1)$ into intervals~$J_1,\dots,J_L$ of length~$\frac1L$, we find (e.g., by the Pigeonhole Principle) that at least two points of~$I_M$ must lie in the same interval, say~$J_j$ for some~$j\in\{1,\dots,L\}$.

On this account, there exist~$r^\star$ and~$r^\sharp\in\N^{n-1}$ with~$r^\star_m$, $r^\sharp_m\in[1,M]$ for all~$m\in\{1,\dots,n-1\}$ such that
$$ \big|\{\alpha\cdot r^\star\}- \{\alpha\cdot r^\sharp\}\big|\le\frac1L.$$
In this way, we have found~$q=(q_1,\dots,q_{n-1})\in\Z^{n-1}$ and~$p\in\Z$ such that
\begin{equation}\label{CUYd0cYmoLCY2tgcy2}
|q_m|\le M \,{\mbox{ for all $m\in\{1,\dots,n-1\}$ and }}\,
|\alpha\cdot q-p|\le\frac1L.
\end{equation}

Now, let~$\omega\in\R^n$ be rationally independent. In particular, we have that~$\omega\ne0$ and therefore, up to reordering variables, we can suppose that
\begin{equation}\label{CMn-1jmdf02o3lefVTvio92hjaL}
0\ne|\omega_n|=\max_{m\in\{1,\dots,n\}}|\omega_m|.
\end{equation}
Accordingly, we can define~$\alpha\in\R^{n-1}$ by~$\alpha_m:=\frac{\omega_m}{\omega_{n}}$ for all~$m\in\{1,\dots,n-1\}$.

We stress that, for all~$k\in\Z^{n}\setminus\{0\}$, setting~$k':=(k_1,\dots,k_{n-1})\in\Z^{n-1}$,
\begin{equation}\label{ifohggbrghiytu836y837630}
(\alpha,1)\cdot k=\alpha\cdot k'+k_n=\frac{1}{\omega_n} (\omega\cdot k)\ne0,\end{equation}
and therefore the vector~$(\alpha,1)$ is rationally independent.

We can thereby employ~\eqref{CUYd0cYmoLCY2tgcy2} with~$L:=\frac{M^{n-1}}{2}$ and, for each~$M\in\N\setminus\{0\}$, we find~$q^{[M]}\in\Z^{n-1}$ and~$p^{[M]}\in\Z$ such that
\begin{equation}\label{kXapsd-c-1}
\big|q^{[M]}\big|\le CM \;{\mbox{ and }}\;
\big|\alpha\cdot q^{[M]}-p^{[M]}\big|\le\frac{C}{M^{n-1}},
\end{equation}
for some constant~$C>0$.

We also remark that~$|\alpha|\le C$, up to renaming~$C$, due to~\eqref{CMn-1jmdf02o3lefVTvio92hjaL} and therefore
\begin{equation}\label{kXapsd-c-2} \big|p^{[M]}\big|\le \big|\alpha\cdot q^{[M]}\big|+
\big|\alpha\cdot q^{[M]}-p^{[M]}\big|\le |\alpha|\,\big|q^{[M]}\big|+
\frac{C}{M^{n-1}}\le C|\alpha|M+\frac{C}{M^{n-1}}\le C M,\end{equation}
up to renaming~$C$ over and over.

Now we set~$k^{[M]}:=\big( q^{[M]},-p^{[M]}\big)\in\Z^n$ and we deduce from~\eqref{ifohggbrghiytu836y837630} and~\eqref{kXapsd-c-1}
that
\begin{equation}\label{KPD:MfcovaBUn}\begin{split}
&\big|\omega\cdot k^{[M]}\big|=
\big|\omega\cdot\big( q^{[M]},-p^{[M]}\big)\big|=|\omega_n|\,\big|(\alpha,1)\cdot\big( q^{[M]},-p^{[M]}\big)\big|
\le\frac{C\,|\omega|}{M^{n-1}}
,
\end{split}\end{equation}
up to keeping renaming~$C$.

Now we dive into the details of the proof of~\eqref{DIpFA:0ladj01-211}.
We assume that~$n\ge2$. For the sake of contradiction, we suppose that~\eqref{DIpFA:0ladj01-211} is violated, hence for all~$K>0$ there exists~$\omega^{[K]}\in\R^n$ for which the set~${\mathcal{E}}_K$ collecting all~$k\in\Z^n\setminus\{0\}$ such that~$|\omega^{[K]}\cdot k|\le\frac{K\,|\omega^{[K]}|}{|k|^{n-1}}$ consists only of finitely many elements (in particular, we have that necessarily~$\omega^{[K]}\ne0$ and, even more, that~$\omega^{[K]}$ is rationally independent).

Thus, we can define
\begin{equation}\label{OJSNLavax2} \mu_K:=\min_{k\in{\mathcal{E}}_K}\big|\omega^{[K]}\cdot k\big|.\end{equation}
Now, in view of~\eqref{KPD:MfcovaBUn}, for all~$M\in\N\setminus\{0\}$, we can find~$k^{[M,K]}\in\Z^{n}$ such that
\begin{equation*}\big|\omega^{[K]}\cdot k^{[M,K]}\big|\le
\frac{C\,\big|\omega^{[K]}\big|}{M^{n-1}},
\end{equation*}
for some constant~$C>0$.

In particular, if~$K$ is sufficiently large, we have that~$k^{[M,K]}
\in{\mathcal{E}}_K$ and therefore, by~\eqref{OJSNLavax2},
\begin{equation} \label{kXapsd-c-7}\mu_K\le
\frac{C\,\big|\omega^{[K]}\big|}{M^{n-1}}.\end{equation}

Accordingly, we can take the limit as~$M\to+\infty$ and find that~$\mu_K=0$.
As a consequence, there exists~$k^\star\in\Z^n\setminus\{0\}$ such that~$ \omega^{[K]}\cdot k_\star=0$, against the fact that~$\omega^{[K]}$ is rationally independent.
This contradiction establishes~\eqref{DIpFA:0ladj01-211}, as desired.

Now we prove~(i). We notice that if~$n=1$ the claim is trivial, hence we suppose that~$n\ge2$. Assume, for a contradiction, that there exist~$\tau\in(0,n-1)$,
$\omega\in\R^n$ and~$\gamma>0$ such that, for every $k\in\Z^n\setminus\{0\}$,
$$ |\omega\cdot k|\ge\frac{\gamma}{|k|^\tau}.$$
In the wake of~\eqref{DIpFA:0ladj01-211} we take a diverging sequence of integers~$k_\ell\in\Z^n$ fulfilling~\eqref{DIpFA:0ladj01-211}, namely~$|k_\ell|\to+\infty$ as~$\ell\to+\infty$ and
$$|\omega\cdot k_\ell|\le\frac{C\,|\omega|}{|k_\ell|^{n-1}}.$$
As a result,
$$ \gamma\le \lim_{\ell\to+\infty}|\omega\cdot k_\ell|\;|k_\ell|^{\tau}\le\lim_{\ell\to+\infty}
\frac{C\,|\omega|}{|k_\ell|^{n-1-\tau}}=0,$$
which is a contradiction, completing the proof of~(i).

We now focus on the proof of~(ii). To this end we observe that, to establish~(ii), it suffices to prove that, for all~$R>0$, the set~$ \big(\R^n\setminus{\mathcal{D}}_{n,\tau}\big)\cap B_R$ has null Lebesgue measure.

We also note that
\begin{eqnarray*}&& \big(\R^n\setminus{\mathcal{D}}_{n,\tau}\big)\cap B_R\\&&\,=\left\{ \omega\in B_R{\mbox{ s.t. for all $\gamma>0$ there exists $k_\gamma \in\Z^n\setminus\{0\}$ for which }}
|\omega\cdot k_\gamma|<\frac{\gamma}{|k_\gamma|^\tau}\right\}\\&&\,=\bigcap_{\gamma>0}
{\mathcal{A}}_\gamma,\end{eqnarray*}
where
$$ {\mathcal{A}}_\gamma:=
\left\{ \omega\in B_R{\mbox{ s.t. there exists $k\in\Z^n\setminus\{0\}$ for which }}
|\omega\cdot k|<\frac{\gamma}{|k|^\tau}\right\}.$$

Hence, the desired result in~(ii) would follow if we can prove that
\begin{equation}\label{TBCGTBRTGTT630TM}
{\mbox{the Lebesgue measure of ${\mathcal{A}}_\gamma$ is less than or equal to $C\gamma R^{n-1}$,}}
\end{equation}
for some constant~$C>0$.

To check~\eqref{TBCGTBRTGTT630TM}, we notice that
\begin{equation}\label{TBCGTBRTGTT630TM2} {\mathcal{A}}_\gamma=\bigcup_{k\in\Z^n\setminus\{0\}} {\mathcal{A}}_{\gamma,k},\end{equation}
where
$$ {\mathcal{A}}_{\gamma,k}:=
\left\{ \omega\in B_R{\mbox{ s.t. }} |\omega\cdot k|<\frac{\gamma}{|k|^\tau}\right\}.$$
Now, the set~${\mathcal{A}}_{\gamma,k}$ is a slab normal to~$\frac{k}{|k|}$ of height~$\frac{2\gamma}{|k|^{\tau+1}}$, intersected~$B_R$, therefore the measure of~${\mathcal{A}}_{\gamma,k}$ is controlled from above by~$\frac{C\gamma R^{n-1}}{|k|^{\tau+1}}$.

Because of this, our assumption on~$\tau$ in~(ii), and~\eqref{TBCGTBRTGTT630TM2}, we find that the measure of~${\mathcal{A}}_{\gamma}$ is controlled from above by
$$ \sum_{k\in\Z^n\setminus\{0\}}\frac{C\gamma R^{n-1}}{|k|^{\tau+1}}\le C\gamma R^{n-1},$$
up to renaming~$C$.

This proves~\eqref{TBCGTBRTGTT630TM} and thus the claim in~(ii).

\paragraph{Solution to Exercise~\ref{DIpFA:1}.} The arguments will rely on Exercise~\ref{DIpFA:0}.

Suppose first that~$r\in(-\infty,2n]$.
Given~$x\in\R^n$, we let~$\omega:=(x,2\pi)\in\R^{n+1}$.
We recall~\eqref{DIpFA:0ladj01-211} (used here in dimension~$n+1$ rather than~$n$)
to find infinitely many~$(k,h)\in(\Z^n\times\Z)\setminus\{0\}$ such that
\begin{equation}\label{YhuijdIpFA:0}|x\cdot k+2\pi h|=|\omega\cdot(k,h)|\le\frac{C\,|\omega|}{(|k|+|h|)^{n}}.\end{equation}

We can thus pick a sequence~$(k_j,h_j)\in(\Z^n\times\Z)\setminus\{0\}$ such that~\eqref{YhuijdIpFA:0} is satisfied and$$|(k_j,h_j)|\to+\infty$$ as~$j\to+\infty$.

We claim that
\begin{equation}\label{YhuijdIpFA:1}
\limsup_{j\to+\infty}|k_j|=+\infty.
\end{equation}
Indeed, suppose not. Then there exists~$K$ such that~$|k_j|\le K$ and therefore, in light of~\eqref{YhuijdIpFA:0},
$$ 2\pi|h_j| -K|x|\le
2\pi|h_j| -|x\cdot k_j|\le|x\cdot k_j+2\pi h_j|\le\frac{C\,|\omega|}{(|k_j|+|h_j|)^{n}}\le C|\omega|.
$$This yields that~$|h_j|$ is also bounded, in contradiction with our construction, and the proof of~\eqref{YhuijdIpFA:1} is thereby complete.

In this way, we have found~$k_j\in\Z^n$ with~$|k_j|\to+\infty$ as~$j\to+\infty$ such that
there exists~$h_j\in\Z$ such that
\begin{equation*}|x\cdot k_j+2\pi h_j|\le\frac{C\,|\omega|}{(|k_j|+|h_j|)^{n}}\le\frac{C\,|\omega|}{|k_j|^{n}}.
\end{equation*}
This gives that~$x\cdot k_j$ is very close to a multiple of~$2\pi$ for~$j$ sufficiently large.

Hence, up to renaming~$C$,
\begin{equation*}
\cos(x \cdot k_j)\ge 1-\frac{C|\omega|^2}{|k_j|^{2n}}.\end{equation*}
Accordingly, using that, for all~$y\in\R$,
\begin{equation}\label{PK-2GHNA-naDFlms} \lim_{\eta\to+\infty}\left(1-\frac{y}{\eta}\right)^\eta=\frac1{e^y},\end{equation}
we find that, for~$j$ sufficiently large,
$$ \big(\cos( x\cdot k_j)\big)^{|k_j|^{2n}}\ge 
\frac{1}{2e^{C|\omega|^2}}=:c\in(0,1)$$
and thus 
$$\big(\cos( x\cdot k_j)\big)^{|k_j|^r}\ge c^{|k_j|^{r-2n}}\ge c.$$
This violates the necessary condition for the series in~\eqref{ER:jDi01} (no matter what invasion of~$\Z^n$ one chooses, in the notation of page~\pageref{INVA-01}) and~(i) is thereby established.

We now focus on the proof of~(ii).
To this end, suppose that~$r\in(2n,+\infty)$, 
and we can thus pick~$\tau\in\left(n,\frac{r}2\right)$.
We recall the notation in
Exercise~\ref{DIpFA:0} (with~$n+1$ in the place of~$n$)
and stress that~$\R^{n+1}\setminus {\mathcal{D}}_{n+1,\tau}$
has null Lebesgue measure in~$\R^{n+1}$ (see in particular item~(ii) in Exercise~\ref{DIpFA:0}).

In this situation, by Fubini's Theorem (see e.g.~\cite[Theorem~6.1]{MR3381284}),
for almost every~$\xi\in\R$ we have that~$\big(\R^{n+1}\setminus {\mathcal{D}}_{n+1,\tau}\big)\cap\{x_{n+1}=\xi\}$ has null measure in~$\R^n$.

In particular, we can pick~$\xi\in(0,+\infty)$ and find a subset~${\mathcal{Z}}$ of~$\R^n$ with null Lebesgue measure such that for every~$y\in\R^n\setminus {\mathcal{Z}}$ we have that~$(y,\xi)\in{\mathcal{D}}_{n+1,\tau}$, that is, by~\eqref{DEFIDIOF}, for all~$(k,h)\in(\Z^n\times\Z)\setminus\{0\}$,
\begin{equation}\label{josd-2we} |y\cdot k+\xi h|\ge\frac{\gamma}{(|k|+|h|)^\tau},\end{equation}
for some~$\gamma>0$.

Let now consider a suitable dilation of~${\mathcal{Z}}$. Namely, we look at
$$ {\mathcal{Z}}_\star:=\frac{\pi}{\xi}{\mathcal{Z}}=\left\{
\frac{\pi \varsigma}{\xi}{\mbox{ with }}\varsigma\in{\mathcal{Z}}
\right\}$$
and we stress that~$ {\mathcal{Z}}_\star$ has also null Lebesgue measure in~$\R^n$.

Thus, for every~$x\in\R^n\setminus{\mathcal{Z}}_\star$,
every~$k\in\Z^n\setminus\{0\}$, and every~$h\in\Z$,
applying~\eqref{josd-2we} to~$y:=\frac{\xi x}{\pi}\in\R^n\setminus{\mathcal{Z}}$
we find that
\begin{equation}\label{DIpF0podk-02woieufjhBISHF2D5t78-1qWD}
|x\cdot k+\pi h|=\left|\frac{\pi y}{\xi}\cdot k+\pi h \right|
=\frac{\pi}{\xi}|y\cdot k+\xi h|\ge\frac{\pi\gamma}{\xi\,(|k|+|h|)^\tau}.\end{equation}

We claim that there exists~$\lambda\in(0,1]$ such that for all~$x\in\R^n\setminus{\mathcal{Z}}_\star$,
$k\in\Z^n\setminus\{0\}$, and~$h\in\Z$ we have that
\begin{equation}\label{DIpF0podk-02woieufjhBISHF2D5t78-1qWDn}
|x\cdot k+\pi h|\ge\frac{\lambda}{ (2+|x|)^\tau\,|k|^\tau}.\end{equation}
To check this, we first observe that when~$|h|\ge(1+|x|)|k|$ we have that
$$|x\cdot k+\pi h|\ge\pi|h|-|x|\,| k|\ge|k|\ge1,$$
hence~\eqref{DIpF0podk-02woieufjhBISHF2D5t78-1qWDn} is satisfied.

If instead~$|h|<(1+|x|)|k|$ we use~\eqref{DIpF0podk-02woieufjhBISHF2D5t78-1qWD} to see that $$
|x\cdot k+\pi h|\ge\frac{\pi\gamma}{\xi\,(|k|+|h|)^\tau}\ge\frac{\pi\gamma}{\xi\,\Big(|k|+(1+|x|)|k|\Big)^\tau},
$$
from which~\eqref{DIpF0podk-02woieufjhBISHF2D5t78-1qWDn} follows
by taking~$\lambda:=\min\left\{1,\frac{\pi\gamma}{\xi}\right\}$.

On the wake of~\eqref{DIpF0podk-02woieufjhBISHF2D5t78-1qWDn},
we observe that~$x\cdot k$ is separated a bit from multiples of~$\pi$
(where the cosine takes values~$\pm1$),
and thus we find that,
for all~$x\in\R^n\setminus{\mathcal{Z}}_\star$ and all~$k\in\Z^n\setminus\{0\}$,
$$\big|\cos(k\cdot x)\big|\le 1-\frac{1}{C (2+|x|)^{2\tau}\,|k|^{2\tau}},$$
for some~$C>0$.

Hence, recalling~\eqref{PK-2GHNA-naDFlms}, for large~$|k|$ we infer that
$$\big|\cos(k\cdot x)\big|^{|k|^{2\tau}}\le c,
$$
for some~$c\in(0,1)$, possibly depending on the given~$x\in\R^n\setminus{\mathcal{Z}}_\star$.

Consequently, for large~$|k|$,
$$\big|\cos(k\cdot x)\big|^{|k|^{r}}\le c^{|k|^{r-2\tau}},
$$ which is a term of a convergent series, thanks to our assumption on~$\tau$.

This gives that when~$r\in(2n,+\infty)$ the series in~\eqref{ER:jDi01}
converges absolutely for almost every~$x\in\R^n$.

\begin{figure}[h]
\includegraphics[height=2.97cm]{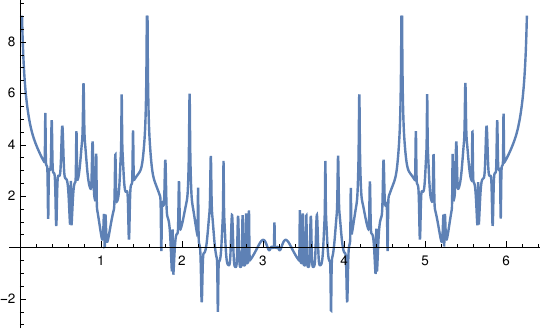}$\,\;\quad$\includegraphics[height=2.97cm]{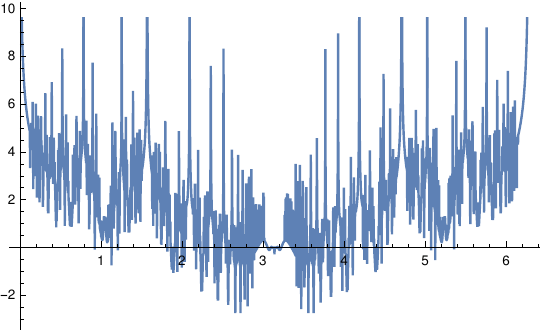}$\,\;\quad$
\includegraphics[height=2.97cm]{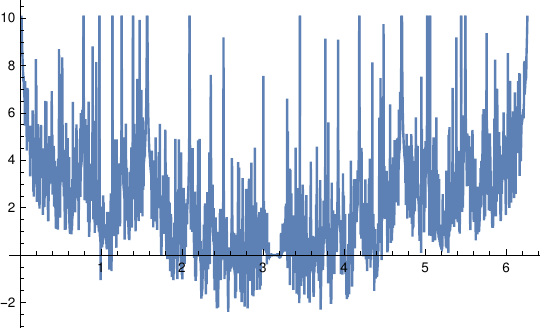}
\centering
\caption{Sketch of the plot of~$\displaystyle\sum_{{k\in\Z}\atop{|k|\le N}} \big(\cos(k x)\big)^{|k|^3}$ with~$N\in\{10,\,25,\,100\}$.}\label{RET12ojwdfldvmfSCNALP}
\end{figure}

However, that series diverges on a dense set in~$\R^n$, namely on the set of rationally dependent~$x\in\R^n$. Indeed, for these values of~$x$ there exists~$k_\star\in\Z^n\setminus\{0\}$ such that~$k_\star\cdot x=0$. For this reason,
for all~$\ell\in\N$, we have that~$\ell k_\star\cdot x=0$ and thus~$\cos(\ell k_\star\cdot x)=1$, showing that the necessary condition for the series in~\eqref{ER:jDi01}
to converge is violated (regardless of the possible invasion of~$\Z^n$ that one can choose).
The proof of~(ii) is thus complete.

See Figure~\ref{RET12ojwdfldvmfSCNALP} for a rough sketch of the situation analyzed here.

See also~\cite{MR1543717, MR1543937} for other approaches to Exercise~\ref{DIpFA:1}.

\section{Solutions to selected exercises of Section~\ref{SE:MINK}}

\paragraph{Solution to Exercise~\ref{MINKE}.} Take~$(-1,1)^n$.

\section{Solutions to selected exercises of Section~\ref{CLRTNSSECMAS}}

\paragraph{Solution to Exercise~\ref{MNVCEC-aIN}.} Let~$\epsilon>0$. We split the circle in small disjoint intervals~$I_1,\dots,I_{K_\epsilon}$ such that, for all~$j\in\{1,\dots,K_\epsilon\}$,
$$ \sup_{I_j}f-\inf_{I_j}f\le\epsilon$$
We also define
$$ f_\epsilon(x):=\sum_{j=1}^{K_\epsilon} \left(\inf_{I_j}f\right)\,\chi_{I_j}(x).$$
and we remark that~$|f(x)-f_\epsilon(x)|\le\epsilon$.

As a result, for all~$m\in\{1,2\}$,
\begin{equation}\label{DSTAFIPMAkcHivvba-1}
\left|\int_\Omega f\big( X_m(\omega)\big)\,d{\mathbb{P}}_\omega-
\int_\Omega f_\epsilon\big( X_m(\omega)\big)\,d{\mathbb{P}}_\omega\right|\le
\int_\Omega \left|f\big( X_m(\omega)\big)- f_\epsilon\big( X_m(\omega)\big)\right|\,d{\mathbb{P}}_\omega\le\epsilon\,{\mathbb{P}}(\Omega)=\epsilon.
\end{equation}

Similarly, since
$$|f(x)f(y)-f_\epsilon(x)f_\epsilon(y)|\le
|f(x)|\,|f(y)-f_\epsilon(y)|
+|f(x)-f_\epsilon(x)|\,|f_\epsilon(y)|\le2\|f\|_{L^\infty(\R/\Z)}\,\epsilon,$$
we find that
\begin{equation}\label{DSTAFIPMAkcHivvba-2}
\left|\int_\Omega f\big( X_1(\omega)\big)\,f\big( X_2(\omega)\big)
\,d{\mathbb{P}}_\omega-
\int_\Omega f_\epsilon\big( X_1(\omega)\big)\,f_\epsilon\big( X_2(\omega)\big)\,d{\mathbb{P}}_\omega\right|\le2\|f\|_{L^\infty(\R/\Z)}\,\epsilon.\end{equation}

Moreover, for all~$m\in\{1,2\}$,
\begin{equation}\label{DSTAFIPMAkcHi:7}
\int_\Omega f_\epsilon\big( X_m(\omega)\big)\,d{\mathbb{P}}_\omega=
\sum_{j=1}^{K_\epsilon} \left(\inf_{I_j}f\right)\,
\int_\Omega\chi_{I_j}\big( X_m(\omega)\big)\,d{\mathbb{P}}_\omega=
\sum_{j=1}^{K_\epsilon}\left(\inf_{I_j}f\right)\,{\mathbb{P}}(X_m\in I_j)
\end{equation}
and therefore, by~\eqref{DSTAFIPMAkcHi}, 
\begin{eqnarray*}
&& \int_\Omega f_\epsilon\big( X_1(\omega)\big)\,d{\mathbb{P}}_\omega\,\int_\Omega f_\epsilon\big( X_2(\omega)\big)\,d{\mathbb{P}}_\omega=\sum_{j,h=1}^{K_\epsilon}
\left(\inf_{I_j}f\right)\,\left(\inf_{I_h}f\right)\,
{\mathbb{P}}(X_1\in I_j){\mathbb{P}}(X_2\in I_h)\\&&\qquad=
\sum_{j,h=1}^{K_\epsilon}\left(\inf_{I_j}f\right)\,\left(\inf_{I_h}f\right)\,{\mathbb{P}}(X_1\in I_j,\,X_2\in I_h)\\&&\qquad=\int_\Omega
\sum_{j,h=1}^{K_\epsilon} \left(\inf_{I_j}f\right)\,\left(\inf_{I_h}f\right)\,
\chi_{I_j}\big( X_1(\omega)\big)\,\chi_{I_h}\big( X_2(\omega)\big)\,d{\mathbb{P}}_\omega\\
&&\qquad=\int_\Omega f_\epsilon\big( X_1(\omega)\big)\,f_\epsilon\big( X_2(\omega)\big)\,d{\mathbb{P}}_\omega.
\end{eqnarray*}

Gathering this information, \eqref{DSTAFIPMAkcHivvba-1} and~\eqref{DSTAFIPMAkcHivvba-2}, we conclude that {\footnotesize
\begin{eqnarray*}&&
\left|\int_\Omega f\big( X_1(\omega)\big)\,f\big( X_2(\omega)\big)\,d{\mathbb{P}}_\omega-
\int_\Omega f\big( X_1(\omega)\big)\,d{\mathbb{P}}_\omega\,
\int_\Omega f\big( X_2(\omega)\big)\,d{\mathbb{P}}_\omega\right|\\&&\qquad\le
\left|\int_\Omega f_\epsilon\big( X_1(\omega)\big)\,f_\epsilon\big( X_2(\omega)\big)\,d{\mathbb{P}}_\omega-
\int_\Omega f_\epsilon\big( X_1(\omega)\big)\,d{\mathbb{P}}_\omega\,
\int_\Omega f_\epsilon\big( X_2(\omega)\big)\,d{\mathbb{P}}_\omega\right|+4\epsilon\big(1+\|f\|_{L^\infty(\R/\Z)}\big)\\&&\qquad=4\epsilon\big(1+\|f\|_{L^\infty(\R/\Z)}\big).
\end{eqnarray*}}
The desired result now follows by taking the limit as~$\epsilon\searrow0$.

\paragraph{Solution to Exercise~\ref{MNVCEC-aIN:2}.} This is a variant of Exercise~\ref{MNVCEC-aIN}, from which we borrow the notation. From~\eqref{DSTAFIPMAkcHienGHaBe}, \eqref{DSTAFIPMAkcHivvba-1}, and~\eqref{DSTAFIPMAkcHi:7},
\begin{eqnarray*}
&&\left|\int_\Omega f\big( X_1(\omega)\big)\,d{\mathbb{P}}_\omega-
\int_\Omega f\big( X_2(\omega)\big)\,d{\mathbb{P}}_\omega\right|\\&&\qquad\le
\left|\int_\Omega f_\epsilon\big( X_1(\omega)\big)\,d{\mathbb{P}}_\omega-
\int_\Omega f_\epsilon\big( X_2(\omega)\big)\,d{\mathbb{P}}_\omega\right|+2\epsilon\\&&\qquad=
\left|\sum_{j=1}^{K_\epsilon}\left(\inf_{I_j}f\right)\,\big({\mathbb{P}}(X_1\in I_j)-{\mathbb{P}}(X_2\in I_j)\big)\right|
+2\epsilon\\&&\qquad=2\epsilon.
\end{eqnarray*}
The desired result now follows by taking the limit as~$\epsilon\searrow0$.

\paragraph{Solution to Exercise~\ref{MNVCEC-aIN:3}.} This is a classic in the theory of Fourier Transforms, but it is often useful in many contexts. To prove the desired result, we define
$$ f(\zeta):=\int_{-\infty}^{+\infty} e^{-ax^2-2\pi i\zeta x}\,dx
\qquad{\mbox{and}}\qquad g(\zeta):=\sqrt{\frac{\pi}a}\,e^{-\frac{\pi^2\zeta^2}{a}}.$$
We notice that
$$ f(0)=\int_{-\infty}^{+\infty} e^{-ax^2}\,dx=\sqrt{\frac{\pi}a}=g(0),$$
that, by the Dominated Convergence Theorem and an integration by parts,
\begin{eqnarray*}
&& f'(\xi)=-\int_{-\infty}^{+\infty}2\pi ix\, e^{-ax^2-2\pi i\zeta x}\,dx=
\frac{\pi i}a\int_{-\infty}^{+\infty}\frac{d}{dx}\left( e^{-ax^2}\right)\,e^{-2\pi i\zeta x}\,dx\\
&&\qquad\qquad\qquad=-\frac{2\pi^2 \zeta}a\int_{-\infty}^{+\infty}\ e^{-ax^2}\,e^{-2\pi i\zeta x}\,dx=-\frac{2\pi^2 \zeta}a\,f(\zeta)
\end{eqnarray*}
and that
$$ g'(\zeta)=-\frac{2\pi^2\zeta}a\,\sqrt{\frac{\pi}a}\,e^{-\frac{\pi^2\zeta^2}{a}}=-\frac{2\pi^2 \zeta}a\,g(\zeta).$$
This gives that both~$f$ and~$g$ satisfy the same ordinary differential equation, with the same initial condition,
and accordingly they must be equal by the uniqueness of the solution of ordinary differential equations.

\paragraph{Solution to Exercise~\ref{MNVCEC-aIN:4}.} We remark that
$$ \sup_{x\in[0,1]} \sum_{{k\in\Z}\atop{|k|>N}} e^{-\frac{(x+k)^2}2}\le\sum_{{k\in\Z}\atop{|k|>N}} e^{\frac{2 |k|- k^2}2},
$$which is the tail of a convergent series.

This gives that the series in~\eqref{WRAPGA} converges uniformly and therefore~${\mathcal{W}}$ is continuous
and periodic of period~$1$.
 
We thus calculate its Fourier coefficients and we observe that, by the uniform convergence of the series
and the substitution~$y:=x+k$, for each~$m\in\Z$ we have that
\begin{eqnarray*}&&\sqrt{2\pi}\,\widehat{\mathcal{W}}_m= \sum_{k\in\Z} \int_0^1 e^{-\frac{(x+k)^2}2-2\pi imx}\,dx=
\sum_{k\in\Z} \int_k^{k+1} e^{-\frac{y^2}2-2\pi im(y-k)}\,dy\\&&\qquad\qquad\qquad= \sum_{k\in\Z} \int_k^{k+1} e^{-\frac{y^2}2-2\pi imy}\,dy= \int_{-\infty}^{+\infty} e^{-\frac{y^2}2-2\pi imy}\,dy.
\end{eqnarray*}

Also, Exercise~\ref{MNVCEC-aIN:3} used with~$a:=\frac12$ and~$\zeta:=m$ returns that
$$\int_{-\infty}^{+\infty} e^{-\frac{ x^2}2-2\pi imx}\,dx=\sqrt{2\pi}\,e^{-2\pi^2 m^2}$$
and therefore
$$ \widehat{\mathcal{W}}_m=e^{-2\pi^2 m^2}.$$

The Fourier Series of~${\mathcal{W}}$ is thereby
$$ \sum_{k\in\Z}e^{-2\pi ikx-2\pi^2 k^2}$$
and the desired result follows by Theorem~\ref{BASw}.

\paragraph{Solution to Exercise~\ref{02ktml-PMSdPR-2lXC135.30568emahbP}.}
The bottom line of this exercise is that we are basically 
``flipping coins'':
any digit in the binary expansion has 50\% chance of being either~$0$ or~$1$, and the value of the~$k$th digit does not influence the value of the~$j$th digit if~$k\ne j$. To formalize the mathematical arguments we proceed as follows.

Given~$j\in\N$, we observe that
$$\Omega=\Omega_{j,0}\cup\Omega_{j,1},$$
with~$\Omega_{j,0}:=
\{\omega{\mbox{ s.t. }}\omega_j=0\}$ and~$\Omega_{j,1}:=\{\omega{\mbox{ s.t. }}\omega_j=1\}$,
and that the translation of~$\Omega_{j,0}$ by~$\frac1{2^j}$ gives~$\Omega_{j,1}$.

Consequently, we can compute the corresponding Lebesgue measures and find that $$1=|\Omega|=|\Omega_{j,0}|+|\Omega_{j,1}|=2|\Omega_{j,0}|,$$
returning that \begin{equation}\label{MCAMSOLofnOLMSM0fmg-04o-ktgmDXM.2345}
|\Omega_{j,0}|=\frac12=|\Omega_{j,1}|.\end{equation}

For this reason,
$$ \int_\Omega X_j(\omega)\,d\omega=\frac12|\Omega_{j,1}|=\frac14$$
and
$$ \int_\Omega X_j^2(\omega)\,d\omega=\frac14|\Omega_{j,1}|=\frac18.$$

Moreover, for all~$\omega\in\Omega$ we have that~$X_j(\omega)\in\left\{0,\frac12\right\}$ and accordingly, for all~$A\subseteq\R$,
$$  {\mathbb{P}}(X_j\in A)={\mathbb{P}}\left(X_j\in \left(A\cap\left\{0,\frac12\right\}\right)\right)=
\begin{dcases}
0 &{\mbox{ if }}\displaystyle A\cap\left\{0,\frac12\right\}=\varnothing,\\
|\Omega_{j,0}|&{\mbox{ if }}\displaystyle A\cap\left\{0,\frac12\right\}=\{0\},\\
|\Omega_{j,1}|&{\mbox{ if }}\displaystyle A\cap\left\{0,\frac12\right\}=\left\{\frac12\right\},\\
1 &{\mbox{ if }}\displaystyle A\supseteq\left\{0,\frac12\right\}
.\end{dcases}$$
It follows from this and~\eqref{MCAMSOLofnOLMSM0fmg-04o-ktgmDXM.2345} that
$$  {\mathbb{P}}(X_j\in A)=
\begin{dcases}
0 &{\mbox{ if }}\displaystyle A\cap\left\{0,\frac12\right\}=\varnothing,\\
\displaystyle\frac12&{\mbox{ if either $\displaystyle A\cap\left\{0,\frac12\right\}=\{0\}$ or $\displaystyle A\cap\left\{0,\frac12\right\}=\left\{\frac12\right\}$}},\\
1 &{\mbox{ if }}\displaystyle A\supseteq\left\{0,\frac12\right\}
.\end{dcases}$$
Since the latter does not depend on~$j$, we have that all~$X_j$'s are identically distributed.

Now, to check that they are mutually independent, we want to show that, for any~$N\in\N$ and any sets~$A_0,\dots,A_N$,
\begin{equation}\label{MCAMSOLofnOLMSM0fmg-04o-ktgmDXM.2345b}
{\mathbb{P}}(X_0\in A_0,\dots,\,X_N\in A_N)= {\mathbb{P}}(X_0\in A_0)\dots {\mathbb{P}}(X_N\in A_N).\end{equation}
To check this, we can suppose that, for each~$j\in\{0,\dots,N\}$, either~$A_j\cap\left\{0,\frac12\right\}=\{0\}$ or~$A_j\cap\left\{0,\frac12\right\}=\left\{\frac12\right\}$: indeed, if~$A_j\cap\left\{0,\frac12\right\}=\varnothing$ one obtains~$0$ in both sides of~\eqref{MCAMSOLofnOLMSM0fmg-04o-ktgmDXM.2345b}, while if~$A_j\cap\left\{0,\frac12\right\}=\left\{0,\frac12\right\}$ then certainly~$X_j\in A_j$ and the corresponding term can be removed in both sides of~\eqref{MCAMSOLofnOLMSM0fmg-04o-ktgmDXM.2345b}.

Thus, we denote by~$\kappa_j$ the element in~$A_j\cap\left\{0,\frac12\right\}$ and by~$\eta_j$ the element in~$\left\{0,\frac12\right\}\setminus\{\kappa_j\}$. Let also~$T_j:\Omega\to\Omega$ be the map altering the~$j$th digit, namely
\begin{eqnarray*}&& T_j(0.\omega_0\dots\omega_{j-1}\kappa_j\omega_{j+1}\dots)=(0.\omega_0\dots\omega_{j-1}\eta_j\omega_{j+1}\dots)
\\{\mbox{and}}\qquad&&
T_j(0.\omega_0\dots\omega_{j-1}\eta_j\omega_{j+1}\dots)=(0.\omega_0\dots\omega_{j-1}\kappa_j\omega_{j+1}\dots).\end{eqnarray*}
We note that~$T_j$ is a translation (of magnitude~$\pm \frac1{2^{j+1}}$) and~$T_j\circ T_j$ is the identity.

Therefore,
$$ {\mathbb{P}}(X_j\in A_j)=\int_{\{X_j(\omega)\in A_j\} }d\omega=\int_{\{X_j(T(\varpi))\in A_j\} }d\varpi=
\int_{\{X_j(\varpi)\not\in A_j\}} d\varpi={\mathbb{P}}(X_j\not\in A_j),$$
returning that
$$ {\mathbb{P}}(X_j\in A_j)=\frac{{\mathbb{P}}(X_j\in A_j)+{\mathbb{P}}(X_j\not\in A_j)}2=\frac12$$
and, as a result,
\begin{equation}\label{MCAMSOLofnOLMSM0fmg-04o-ktgmDXM.2345c3}
{\mathbb{P}}(X_0\in A_0)\dots {\mathbb{P}}(X_N\in A_N)=\frac1{2^{N+1}}.\end{equation}

Along the same lines, one could take~$S_j$ to be either~$T_j$ or the identity and consider the translation~$S:=S_0\circ\dots\circ S_N$. Since there are two possible choices for each~$S_j$, there are~$2^{N+1}$ translations of these types, leading to
$$ {\mathbb{P}}(X_0\in A_0,\dots,\,X_N\in A_N)=\frac1{2^{N+1}}.$$
This, together with~\eqref{MCAMSOLofnOLMSM0fmg-04o-ktgmDXM.2345c3}, proves~\eqref{MCAMSOLofnOLMSM0fmg-04o-ktgmDXM.2345b}, as desired.

See e.g.~\cite[pages~3--5]{MR1700749},
\cite[Exercise~7.4]{MR2247694},
and~\cite[Chapter~2, Section~11]{MR2977961} for further details on this type of constructions.



\begin{appendix}

  \chapter{One more example confirming Theorem~\ref{COMADIVESKMDWD}}\label{ANEXDIHALA}

Here we provide yet another example to confirm that the Fourier Series of a continuous function does not, in general, converge pointwise:

\begin{proof}[Another proof of Theorem~\ref{COMADIVESKMDWD}]
We define~$f$ in the interval~$\left[0,\frac12\right]$, with~$f(0)=f\left(\frac12\right)=0$, and we extend it by an even reflection to define it in the interval~$\left[-\frac12,\frac12\right]$. Then, we extend it by periodicity so to have it periodic of period~$1$.

To define the function in~$\left[0,\frac12\right]$, we set~$f(0):=0$ and 
$$ n_j:=3^{\sigma_j},\qquad{\mbox{where}}\qquad\sigma_j:={\sum_{m=0}^j m^4}.$$
Notice that~$n_0=1$, that~$n_{j+1}>n_j$ and that~$n_j\to+\infty$ as~$j\to+\infty$.

Hence, for all~$x\in\left(0,\frac12\right]$ we can identify a unique~$j_x\in\N$ such that
$$ x\in \left(\frac1{2n_{j_x+1}},\frac1{2n_{j_x}} \right].$$
Thus, for every~$x\in\left(0,\frac12\right]$ we can define
\begin{equation}\label{DEFEFFEDISKMAHOEKER} f(x):=\frac{\sin (4\pi n_{j_x+1}\,x)}{(j_x+2)^2}.\end{equation}
We remark that
\begin{equation} \label{HanSkdmP}\begin{split}&\lim_{x\searrow\frac{1}{2n_j}}f(x)=\frac{\sin (2\pi n_{j}/n_{j})}{(j+1)^2}=0
\\ {\mbox{and}}\qquad&
\lim_{x\nearrow\frac{1}{2n_j}}f(x)=\frac{\sin (2\pi n_{j+1}/n_j)}{(j+2)^2}=\frac{\sin (2\pi 3^{(j+1)^4})}{(j+2)^2}=0.\end{split}
\end{equation}
Besides,
$$ \lim_{x\searrow0}j_x=+\infty$$
and therefore
$$ \lim_{x\searrow0}|f(x)|\le\lim_{x\searrow0}\frac{1}{j_x^2}=0.$$
This and~\eqref{HanSkdmP} entail that~$f$ (and actually the extension of~$f$) is a continuous function.

Thus, to complete the proof of Theorem~\ref{COMADIVESKMDWD}, we need to show that~$S_{N}(0)$ does not converge.
To this end, we use Corollary~\ref{9.0km2df0.2wnedfo.1-2x} (and we will translate the integral over the period,
according to Exercise~\ref{fr12})
to see that
\begin{equation*}\begin{split} S_N(x)-f(x)=
\int_{-1/2}^{1/2} \big(f(y)-f(x)\big)\,\frac{\sin\big((2N+1)\pi (y-x)\big)}{\sin(\pi(y- x))}\,dy
\end{split}\end{equation*}
and therefore, since~$f$ is even,
\begin{equation}\label{0polwdf0erf-2}\begin{split} S_N(0)&=
\int_{-1/2}^{1/2} f(y)\,\frac{\sin\big((2N+1)\pi y\big)}{\sin(\pi y)}\,dy\\&=
2\int_{0}^{1/2} f(y)\,\frac{\sin\big((2N+1)\pi y\big)}{\sin(\pi y)}\,dy.
\end{split}\end{equation}

Moreover,
\begin{eqnarray*}
&&\left|\int_{0}^{1/2} f(y)\,\frac{\sin\big((2N+1)\pi y\big)}{\sin(\pi y)}\,dy-\int_{0}^{1/2} f(y)\,\frac{\sin(2N\pi y)}{\sin(\pi y)}\,dy\right|
\\&&\qquad\le\int_{0}^{1/2} |f(y)|\,\frac{\big|\sin(2N\pi y)\big|\big(1-\cos(\pi y)\big)+\big|\cos(2N\pi y)\big|\sin(\pi y)}{\sin(\pi y)}\,dy\\&&\qquad\le\|f\|_{L^\infty(\R)}\int_{0}^{1/2}\frac{\big(1-\cos(\pi y)\big)+\sin(\pi y)}{\sin(\pi y)}\,dy,
\end{eqnarray*}
which is a finite quantity.

Hence, if we denote quantities that are bounded independently of~$N$ by the symbol~{\mbox{$\Xi$}}, we can write~\eqref{0polwdf0erf-2} in the form
\begin{equation}\label{0polwdf0erf-3}\begin{split} S_N(0)=
2\int_{0}^{1/2} f(y)\,\frac{\sin(2N \pi y)}{\sin(\pi y)}\,dy+{\mbox{{\mbox{$\Xi$}}}}.
\end{split}\end{equation}

Similarly,
\begin{eqnarray*}
&&\left|
\int_{0}^{1/2} f(y)\,\frac{\sin(2N \pi y)}{\sin(\pi y)}\,dy-\int_{0}^{1/2} f(y)\,\frac{\sin(2N \pi y)}{\pi y}\,dy
\right|
\\&&\qquad\le\|f\|_{L^\infty(\R)}\int_{0}^{1/2} f(y)\,\frac{|\sin(\pi y)-\pi y|}{\pi y\sin(\pi y)}\,dy,
\end{eqnarray*}
which is a finite quantity and we thereby rewrite~\eqref{0polwdf0erf-3} in the form
\begin{equation}\label{0polwdf0erf-4}\begin{split} \pi S_N(0)=
2\int_{0}^{1/2} f(y)\,\frac{\sin(2N \pi y)}{y}\,dy+{\mbox{{\mbox{$\Xi$}}}}.
\end{split}\end{equation}

Now we use the trigonometric formula in Exercise~\ref{SPDCD-0.01}
to see that
\begin{eqnarray*}&& 2\sin (4\pi n_{j+1}\,y)\,\sin(4\pi n_{k}\,y)\\&&\qquad=
\cos\big( 2\pi(n_{j+1}-n_k)y\big)-\cos\big( 2\pi(n_{j+1}+n_k)y\big)
\\&&\qquad=
\cos\big( 2\pi|n_{j+1}-n_k|y\big)-\cos\big( 2\pi|n_{j+1}+n_k|y\big)
\end{eqnarray*}
and correspondingly, using the substitutions~$t:=2|n_{j+1}\pm n_k|y$,
\begin{equation}\label{0poel2f-1kmdm93.11}\begin{split}&
\left|2\int_{1/(2n_{j+1})}^{1/(2n_{j+1})} \frac{\sin (4\pi n_{j+1}\,y)\,\sin(4\pi n_{k}\,y)}{y}\,dy\right|\\
&\qquad\le\left|\int_{1/(2n_{j+1})}^{1/(2n_j)}\frac{\cos\big( 2\pi|n_{j+1}-n_k|y\big)}{y}\,dy\right|+
\left|\int_{1/(2n_{j+1})}^{1/(2n_j)}\frac{\cos\big( 2\pi|n_{j+1}+n_k|y\big)}{y}\,dy\right|\\
&\qquad=\left|\int_{|n_{j+1}-n_k|/n_{j+1}}^{|n_{j+1}-n_k|/n_{j}}\frac{\cos(\pi t)}{t}\,dt\right|+
\left|\int_{|n_{j+1}+n_k|/n_{j+1}}^{|n_{j+1}+n_k|/n_{j}}\frac{\cos(\pi t)}{t}\,dt\right|\\&\qquad\le
\left|\int_{|n_{j+1}-n_k|/n_{j}}^1\frac{\cos(\pi t)}{t}\,dt\right|+\left|\int_{|n_{j+1}+n_k|/n_{j}}^1\frac{\cos(\pi t)}{t}\,dt\right|\\&\qquad\quad\qquad+
\left|\int_{|n_{j+1}-n_k|/n_{j+1}}^1\frac{\cos(\pi t)}{t}\,dt\right|+
\left|\int_{|n_{j+1}+n_k|/n_{j+1}}^1\frac{\cos(\pi t)}{t}\,dt\right|.
\end{split}\end{equation}

Furthermore, if~$p\in\N$ and~$q\in\N\setminus\{0\}$,
\begin{eqnarray*}&& n_{p+q}-n_p=3^{\sigma_p+(p+1)^4+\dots+(p+q)^4}-3^{\sigma_p}=
3^{\sigma_p}\big( 3^{(p+1)^4+\dots+(p+q)^4}-1\big)\\&&\qquad\ge3^{\sigma_p}(3-1)\ge2n_p.\end{eqnarray*}
As a result, if~$k>j+1$ then
\begin{eqnarray*}&&
\min\left\{ \frac{|n_{j+1}-n_k|}{n_j},\;\frac{|n_{j+1}+n_k|}{n_j},\;\frac{|n_{j+1}-n_k|}{n_{j+1}},\;\frac{|n_{j+1}+n_k|}{n_{j+1}}\right\}
=\frac{|n_{j+1}-n_k|}{n_{j+1}}\\&&\qquad=\frac{n_k-n_{j+1}}{n_{j+1}}\ge2
\end{eqnarray*}
and if instead~$k<j+1$ then
\begin{eqnarray*}&&
\min\left\{ \frac{|n_{j+1}-n_k|}{n_j},\;\frac{|n_{j+1}+n_k|}{n_j},\;\frac{|n_{j+1}-n_k|}{n_{j+1}},\;\frac{|n_{j+1}+n_k|}{n_{j+1}}\right\}
=\frac{|n_{j+1}-n_k|}{n_{j+1}}\\&&\qquad=\frac{n_{j+1}-n_{k}}{n_{j+1}}\ge\frac{n_{j+1}-\frac{n_{j+1}}3}{n_{j+1}}=\frac23.
\end{eqnarray*}

All in all, if~$k\ne j+1$,
$$ \min\left\{ \frac{|n_{j+1}-n_k|}{n_j},\;\frac{|n_{j+1}+n_k|}{n_j},\;\frac{|n_{j+1}-n_k|}{n_{j+1}},\;\frac{|n_{j+1}+n_k|}{n_{j+1}}\right\}\ge\frac23.$$

From this and~\eqref{0poel2f-1kmdm93.11} we infer that, if~$k\ne j+1$,
\begin{eqnarray*}&&
\left|2\int_{1/(2n_{j+1})}^{1/(2n_{j+1})} \frac{\sin (4\pi n_{j+1}\,y)\,\sin(4\pi n_{k}\,y)}{y}\,dy\right|
\\&&\qquad\le4\sup_{\alpha\ge2/3}\left|\int_{\alpha}^1\frac{\cos(\pi t)}{t}\,dt\right|\le C,
\end{eqnarray*} for some~$C>0$ (see Exercise~\ref{CLURVB}).

Consequently,
$$ \left|2\sum_{j\in\N\setminus\{k-1\}}\int_{1/(2n_{j+1})}^{1/(2n_j)}\frac{\sin (4\pi n_{j+1}\,x)\,\sin(4n_k \pi y)}{(j+2)^2 y}\,dy\right|\le\sum_{j\in\N}\frac{C}{(j+1)^2},$$
which is a finite quantity.

This, \eqref{DEFEFFEDISKMAHOEKER} and~\eqref{0polwdf0erf-4} entail that, choosing~$N:=2n_k$,
\begin{equation}\label{kaHAfT83.1}\begin{split} \pi S_{2n_k}(0)&=
2\int_{0}^{1/2} f(y)\,\frac{\sin(4n_k \pi y)}{y}\,dy+{\mbox{{\mbox{$\Xi$}}}}\\
&=2\sum_{j=0}^{+\infty}\int_{1/(2n_{j+1})}^{1/(2n_j)}\frac{\sin (4\pi n_{j+1}\,x)\,\sin(4n_k \pi y)}{(j+2)^2 y}\,dy+{\mbox{{\mbox{$\Xi$}}}}\\
&=2\int_{1/(2n_{k})}^{1/(2n_{k-1})}\frac{\sin^2(4n_k \pi y)}{(k+1)^2 y}\,dy+{\mbox{{\mbox{$\Xi$}}}}.
\end{split}\end{equation}

Now we use the substitution~$t:=4n_k \pi y$ and observe that
\begin{eqnarray*} && 2\int_{1/(2n_{k})}^{1/(2n_{k-1})}\frac{\sin^2(4n_k \pi y)}{(k+1)^2y}\,dy=
\frac{2}{(k+1)^2}\int_{2\pi}^{2 n_k\pi/n_{k-1}}\frac{\sin^2 t}{t}\,dt\\&&\qquad=
\frac{2}{(k+1)^2}\int_{2\pi}^{2\pi \,3^{k^4}}\frac{\sin^2 t}{t}\,dt\ge \frac{2}{(k+1)^2}\sum_{\ell=1}^{3^{k^4}-1}
\int_{2\ell\pi+\frac\pi4}^{2\ell\pi+\frac\pi2}\frac{\sin^2 t}{t}\,dt\\&&\qquad\ge\frac{\sqrt2}{(k+1)^2}\sum_{\ell=1}^{3^{k^4}-1}
\int_{2\ell\pi+\frac\pi4}^{2\ell\pi+\frac\pi2}\frac{dt}{2(\ell+1)\pi}=\frac{\sqrt2}{8(k+1)^2}\sum_{\ell=1}^{3^{k^4}-1}
\frac{1}{\ell+1}\\&&\qquad\ge c\,k^2,
\end{eqnarray*}
for some~$c>0$.

Hence, by~\eqref{kaHAfT83.1},
\begin{eqnarray*} \pi S_{2n_k}(0)\ge c\,k^2
+{\mbox{{\mbox{$\Xi$}}}},\end{eqnarray*}
which diverges.
\end{proof}

\end{appendix}


\backmatter

\bibliography{biblio}
\bibliographystyle{alpha}

\printindex

\end{document}